%% file: qpe.tex
\documentclass{tran-l}
\usepackage{amssymb}
\usepackage{amsmath}
\usepackage{comment}
\usepackage{graphicx,psfrag}
\usepackage{subfigure}
\usepackage{enumerate}
\usepackage[flushright]{caption2}

\newcommand{\NN}{{\mathbb N}}

\newcommand{\ZZ}{{\mathbb Z}}

\newcommand{\RR}{{\mathbb R}}
\newcommand{\nul}{\emptyset }
\newcommand{\D}{\Delta }
\newcommand{\de}{\delta }
\newcommand{\bt}{\beta }
\newcommand{\z}{\zeta }
\newcommand{\io}{\iota }
\newcommand{\ka}{\kappa }
\newcommand{\lm}{\lambda }
\newcommand{\x}{\chi }
\newcommand{\Lm}{\Lambda }
\newcommand{\e}{\varepsilon }

\newcommand{\G}{\Gamma }
\newcommand{\g}{\gamma }
\newcommand{\al}{\alpha }
\renewcommand{\t}{\tau }
\renewcommand{\th}{\theta}
\newcommand{\s}{\sigma }
\renewcommand{\S}{\Sigma }
\newcommand{\W}{\Omega}
\newcommand{\w}{\omega }
\newcommand{\pd}{\partial}
\newcommand{\dcup}{{\ \dot{\cup}\ }}
\newcommand{\bdcup}{\dot{\bigcup}}
\newcommand{\mbf}{\mathbf}
\newcommand{\wt}{\widetilde}
\renewcommand{\^}{\wedge}
\newcommand{\cPQ}{{\mathcal{QE}}}
\newcommand{\cSPQ}{{\mathcal{SQE}}}
\newcommand{\cPG}{{\mathcal{GE}}}
\newcommand{\cSPG}{{\mathcal{SGE}}}

\newcommand{\maps}{\longrightarrow}
\newcommand{\bl}{\backslash}
\newcommand{\+}{\times}

\newcommand{\cA}{{\mathcal A}}
\newcommand{\cB}{{\mathcal B}}
\newcommand{\cC}{{\mathcal C}}
\newcommand{\cD}{{\mathcal D}}
\newcommand{\cE}{{\mathcal E}}
\newcommand{\cF}{{\mathcal F}}

\newcommand{\cH}{{\mathcal H}}
\newcommand{\cI}{{\mathcal I}}
\newcommand{\cL}{{\mathcal L}}
\newcommand{\cM}{{\mathcal M}}
\newcommand{\cN}{{\mathcal N}}
\newcommand{\cO}{{\mathcal O}}
\newcommand{\cP}{{\mathcal P}}
\newcommand{\cQ}{{\mathcal Q}}
\newcommand{\cR}{{\mathcal R}}
\newcommand{\cV}{{\mathcal V}}
\newcommand{\cal}{\mathcal}
\newcommand{\inr}{\textrm{int}}
\newcommand{\pprime}{{\prime\prime}}
\newcommand{\Baa}{Baumslag60}
\newcommand{\Bab}{Baumslag65}
\newcommand{\BMR}{BaumslagMyasnikovRemeslennikov99}
\newcommand{\BP}{BogleyPride92}
\newcommand{\CEa}{ComerfordEdmunds84}
\newcommand{\CEb}{ComerfordEdmunds81}
\newcommand{\Cu}{Culler81}

\newcommand{\DHa}{DuncanHowie91a}
\newcommand{\DHb}{DuncanHowie91b}
\newcommand{\DHc}{DuncanHowie93}
\newcommand{\ET}{EdjvetThomas97}
\newcommand{\Fe}{Fenn83}
\newcommand{\Fr}{Friel01}
\newcommand{\Ge}{Gersten87}
\newcommand{\GL}{GrigorchukLysionok92}
\newcommand{\GT}{GoldsteinTurner79}
\newcommand{\Ha}{Howie89}
\newcommand{\Hb}{Howie90}
\newcommand{\Hc}{Howie91}
\newcommand{\HT}{HowieThomas93}

\newcommand{\KMb}{KharlampovichMyasnikov982}
\newcommand{\KMc}{KharlampovichMyasnikov983}
\newcommand{\KMd}{KharlampovichMyasnikov984}
\newcommand{\Lya}{Lyndon60a}
\newcommand{\Lyb}{Lyndon60b}
\newcommand{\LS}{LyndonSchutzenberger62}
\newcommand{\LM}{LipschutzMiller71}
\newcommand{\Mak}{Makanin82}
\newcommand{\Mat}{Matiyasevic93}
\renewcommand{\MR}{MyasnikovRemeslennikov96}
\newcommand{\Pri}{Pride91}
\newcommand{\Ro}{Rourke79}
\newcommand{\RS}{RipsSela95}
\newcommand{\Sea}{Sela01}
\newcommand{\Seb}{Sela02}
\newcommand{\Sh}{Short84}
\newcommand{\ovr}[1]{\overline{ #1 }} 
\newcommand{\supstar}[1]{#1^{*}}
\renewcommand{\labelenumi}{(\theenumi)}
\renewcommand{\theenumi}{\roman{enumi}}

\numberwithin{equation}{section}
\title{Exponential Genus Problems in one--relator products of groups}
\author{A.J.\ Duncan}
\address{School of Mathematics and Statistics\\
Merz Court\\
University of Newcastle upon Tyne\\
Newcastle upon Tyne\\
NE1 7RU\\
UK}
\email{a.duncan@ncl.ac.uk}
\urladdr{http://www.ncl.ac.uk/math/}
\subjclass[2000]{Primary 20F65 20F05 20F10; Secondary 20F06 57M07} 
\date{\today}
\thanks{Research partially supported by EPSRC grants B/91/RFH/1148
and GR/N28313}
\begin{document}
\newcommand{\be}{\begin{enumerate}}
\newcommand{\ee}{\end{enumerate}}
\newcommand{\bd}{\begin{description}}
\newcommand{\ed}{\end{description}}
\newtheorem{theorem}{Theorem}[section]
\newtheorem{lemma}[theorem]{Lemma}
\newtheorem{defn0}[theorem]{Definition}
\newenvironment{defn}{\begin{defn0} \rm}{\end{defn0}}
\newtheorem{assum}[theorem]{Assumption}
\newtheorem{corol}[theorem]{Corollary}
\newtheorem{prop}[theorem]{Proposition}
\newtheorem{nte}[theorem]{Note}
\newenvironment{note}{\begin{nte} \rm}{\end{nte}}
\newtheorem{exam}[theorem]{Example}
\newenvironment{exx}{\begin{exam} \rm}{\end{exam}}
\newtheorem{rmk}{Remark}
\newenvironment{remark}{\begin{rmk} \rm}{\end{rmk}}
\newlength{\storemargini}
\newlength{\storemarginii}
\bibliographystyle{plain}
\input in1
\input in2
\input in3
\input in4
\bibliography{/home/najd2/tex/bib/qpe,/home/najd2/tex/bib/algebra}
\end{document}

%% file: in1.tex
\maketitle


\begin{abstract}
Exponential equations in free groups were studied initially 
by Lyndon and Sch\"utzenberger \cite{\LS}
and then by Comerford and Edmunds \cite{\CEa}. Comerford and Edmunds
showed that the problem of determining
whether or not the class of quadratic exponential equations have solution is decidable, in
finitely generated free groups.  In this paper we show that for finite systems of quadratic exponential
equations decidability passes, under certain hypotheses, from the factor groups to 
free products and one--relator products.
\end{abstract}
\section{Introduction}
\label{int}
In \cite{\LS} Lyndon and Sch\"utzenberger studied the equation \[a^M=b^Nc^P,\]
where $a, b$ and $c$ are fixed elements of a  free group and $M,N$ and $P$ are variables which 
may take integer variables. This is an example of an exponential equation. In fact it is shown in 
\cite{\LS} that this equation has no solution with $M$, $N$ and $P$ all greater than $1$, unless
$a$, $b$ and $c$ are all powers of a common element.
Comerford and Edmunds
\cite{\CEa} studied more general exponential equations and showed that the problem of determining
whether or not the class of quadratic exponential equations have solution is decidable, in
finitely generated free groups. The definition of Comerford and Edmunds is that adopted 
here.\footnote{In \cite{\CEa} the word ``parametric'' is used instead of
``exponential''.
}
In fact Lyndon also introduced exponential groups (on which an associative ring acts by exponention) which
he used to characterise the solution sets of equations in $1$ variable in free groups \cite{\Lya}, \cite{\Lyb}.
Lyndon considered the word problem in an exponential group and this can also be regarded as an exponential
equation, but using a definition which allows more general values of exponents than that used here and
in \cite{\CEa}.  Exponential groups have been widely studied (see for example \cite{\Baa}, \cite{\Bab},
\cite{\BMR}, \cite{\MR}, \cite{\KMb}, \cite{\KMc}, \cite{\KMd}). In particular these groups play an important
role in work on Tarskii's conjectures on the elementary theory of free groups 
(see \cite{\KMb}, \cite{\KMc}, \cite{\Sea} and \cite{\Seb}).
In this paper we consider exponential equations in the
setting of one--relator products of groups. We show that for finite systems of quadratic exponential
equations decidability is, under some extra hypotheses, inherited by free products and one--relator products.

Let $G$ be a group given by a recursive presentation $\langle C|R\rangle$, write $\ovr u$ for the 
element of $G$ represented by $u\in F(C)$ and  let $X$ be a recursive set.
An equation over $G$ with variables in $X$ is an expression of the form 
\begin{equation}\label{int_eqn}
w=1,
\end{equation} 
where $w$ is an element of $F(C)$. 
A solution to (\ref{int_eqn}) is a homomorphism $\phi:F(C\dcup X) \maps G$ such that
$\phi(w)=1$ and $\phi(c)=\ovr{c}$, for all $c\in C$.
A class of equations over $G$ is said to be solvable (or decidable) if there exists an algorithm
which will determine whether or not a given equation of the class has a solution. As it is well known
that the word problem is not solvable in certain finitely presented groups, in general equations are
not solvable. However in suitably restricted classes of groups large classes of equations are known
to be solvable. One of the most striking results in this direction is Makanin's: that finite systems
of equations are solvable in finitely generated free groups \cite{\Mak}. 
Furthermore
Rips and Sela \cite{\RS} have
generalised this result to
torsion--free hyperbolic groups.

An equation is quadratic if every variable which occurs at all 
does so exactly twice. For example, if $u,v\in F(C)$,
and 
and $x,y,z\in X$, for then the equations $u=1$,  $ux^{-1}vxy^2=1$ and $ux^{-1}vx[y,z]=1$ are all quadratic.
These examples are all of the form
\[ 
\prod_{i=1}^{n}x_i^{-1}u_ix_i\,
\prod_{i=1}^{t}\left [x_{i+n},x_{i+n+t}\right ]
\prod_{i=1}^{p}x_{i+n+2t}^2=1,
\]
with $u_i\in F(C)$ and $x_i\in X$.
In fact all quadratic equations are equivalent to at least one such equation: where equations with the
same solution sets are deemed equivalent. The {\em genus} 
of an $n$--tuple $\ovr u_1,\ldots ,\ovr u_n$ of elements
of $G$ is the set of pairs of integers $(t,p)$ such that the equation above has solution.
The genus problem is that of determining whether or not a given pair $(t,p)$ belongs to the 
genus of a given $n$--tuple of elements of $G$.
It is known (see Example \ref{exx_classic_genus} 
and \cite{\DHa}, \cite{\DHb} and \cite{\DHc}) that decidability of the genus
problem is inherited from the factor groups by one--relator products, under suitable hypotheses on the
factors and the relator.  

In this paper we generalise some of the results of Example \ref{exx_classic_genus} mentioned above
to quadratic exponential equations. If $\lm$ and $\mu$ are unknowns, $u,v \in F(C)$ and $x,y,z\in X$ then
the expressions $u^\lm=1$,  $u^\lm x^{-1}v^\mu xy^2=1$ and $u^\lm x^{-1}v^\mu x[y,z]=1$
are all examples of quadratic exponential equations. A solution to one of these equations consists of 
a map $\al:\{\lm,\mu\}\maps \ZZ$, say $\al(\lm)=l$, $\al(\mu)=m$, 
and a homomorphism $\phi:F(C)\maps G$ such that $\phi$ is a solution to
the quadratic equation $u^l=1$,  $u^lx^{-1}v^m xy^2=1$ or  $u^l x^{-1}v^m x[y,z]=1$, as appropriate.
Usually some conditions are imposed on the values that $\lm$ and $\mu$ can take. For instance in the exponential
equations above it is reasonable to insist that $\lm,\mu\neq 0$. More generally 
solutions with $\lm>l$ and $\mu>m$ may be of interest.
We leave the formal definition 
to Section \ref{qpequations} but roughly speaking
a quadratic exponential equation consists of an expression, similar to one of those above, 
together with a finite system of  equations, congruences and inequalities on $\lm,\mu,\ldots$, which we
call a parameter system. 
We restrict to linear parameter systems, since  as shown by Matiyasevi\u{c}, Hilbert's tenth problem is unsolvable 
for non--linear systems \cite{\Mat}.
In \cite{\CEa} an algorithm is constructed which will determine whether or not  a given quadratic exponential
equation has a solution in the case where $G$ is a finitely generated free group (see Example \ref{exx_CE}).


Let $A$ and $B$ be groups.  A {\em one--relator product} of $A$ and $B$ is
a group of the form $G=(A*B)/N(s)$, where $A*B$ is the free product,
$s\in A*B$ is a cyclically reduced word of length at least $2$, and $N(s)$ 
is the normal closure of $s$ in $A*B$. 
The exponential genus problem is a natural generalisation of
the genus problem, described above, from quadratic to quadratic exponential equations. We show that
under suitable conditions the decidability of the exponential genus problem is inherited by $G$ from
the factor groups $A$ and $B$. 
In order to prove this result a stronger hypothesis on $A$ and $B$ has been imposed: namely 
that simultaneous systems of quadratic exponential equations are solvable in $A$ and
$B$, in the following sense. Let $Q_1,\ldots ,Q_m, Q_{m+1},\ldots,Q_{m+n}$ be quadratic exponential 
equations subject to a system of parameters $\cL$. If $Q_i$ is an exponential equation over $A$, for
$i=1,\ldots,m$ and over $B$, for $i=m+1,\ldots,m+n$ then we regard these equations as simultaneous. 
For example take $a\in A$ and $b\in B$, where $Q_1$ is  $a^\lm=1$ and $Q_2$ is $b^\mu=1$, subject to the parameter
system consisting of $\lm>0$, $\mu>0$ and 
$\lm=3\mu-7$.
If there
is an algorithm to decide whether or not such systems of simultaneous  equations have 
a simultaneous solution then we say that
quadratic exponential equations are solvable over $(A,B)$. This gives rise to  a corresponding simultaneous 
definition of exponential genus (Definition \ref{L_genus}).
We prove that if quadratic exponential equations are solvable over $(A,B)$ then
quadratic exponential equations defined over $G$ are also solvable, as long as the relator is a sufficiently
high power and not of a particularly pathogenic type (see Corollary \ref{plain}). 
In fact we show that simultaneous systems of 
quadratic equations defined over one--relator groups of the required form are
also solvable (see Corollary  \ref{main_eqn}).  This follows from the  main result of the paper, Theorem \ref{main},
which is that if the exponential genus problem is solvable in a class of groups
then it is also solvable in one--relator products, of the required form, of these groups. 
All these results hold if we replace one--relator products with free products. This means in particular that
we obtain the  result of \cite{\CEa}, that quadratic exponential equations are solvable in finitely generated
free groups (see Corollary \ref{cyclic_corol}).

The paper is structured as follows. In Sections \ref{qwords} and \ref{qpequations} we describe exponential
quadratic equations and exponential genus. In Section \ref{resolution} we describe how, given an exponential
quadratic equation we can form a finite set, called a resolution, of corresponding equations in a normal form, 
which we call special. The original equation has a solution if and only if some element of the
resolution has a solution so this allows us to restrict attention to special equations.
In Section \ref{decisions} we describe in more detail the decision problems under consideration.
In Section \ref{qpe_and_genus} we state the main result of the paper, Theorem \ref{main}, and show
how it is a consequence of the more technical Theorem \ref{mainc}. The remainder of the paper is devoted to
proving Theorem \ref{mainc}. We use pictures which are 
described here following, mainly, Howie \cite{\Ha}, \cite{\Hb} and
\cite{\Hc}. The use of pictures in group theory goes back to Rourke \cite{\Ro} (see
also Fenn \cite{\Fe}) and was developed to apply to one--relator products by
Short \cite{\Sh}. The theory has been further developed and applied to various aspects of geometric 
group theory in, for example, \cite{\BP}, \cite{\ET}, \cite{\HT}
and \cite{\Pri}. 
Here we extend use of pictures made in \cite{\DHa}, \cite{\DHb} and \cite{\DHc}. 
 We introduce pictures and establish some of their basic properties in Section \ref{pictures}. In particular
we describe the link between solutions of equations and pictures. A picture can be thought of as a graph on
a compact surface $\S$ with boundary. Each  boundary component is assigned a label and these labels correspond
to the group elements which occur in a quadratic equation. A solution to the quadratic equation arises 
from a picture and vice--versa. In Section \ref{corridor} we define corridors which
are later used to bound the degree of exponents occuring in certain minimal solutions to exponential equations. 
In Sections \ref{angle} and 
\ref{curvature} we assign angle and then curvature to vertices and regions of pictures. These are used together
with corridors to 
bound the size of pictures that correspond to solutions to equations. To obtain this bound we analyse the curvature
of regions of our pictures and show, in Sections \ref{config_C}, \ref{config_D} and \ref{angle_final} that non--negative
curvature occurs only at the boundary and in certain well defined configurations, see Figures \ref{F1} and \ref{F2}
 for example. This allows us, in Section \ref{isoperimetry}, to bound the number of edges of a picture corresponding to a special equation, suitably
minimised. We complete the proof of Theorem \ref{mainc} in Section \ref{proof_mainc}. 

It seems likely that the results proved here should hold with at least some of the restrictions on
the relator  removed. 
In fact it would be surprising if the conditions on the factor groups and relators of Example \ref{exx_classic_genus}
were not sufficient to prove  Theorem \ref{main}. Indeed I believe that  
it is possible to push the proof given here to obtain results more like those of Example \ref{exx_classic_genus}
without very much more effort. 
However, given the length of the proof, 
I do not feel that any increase in its complexity is desirable, especially for such minor gains.
A more systematic approach to the arguments of Section
\ref{pictures} onwards might considerably shorten the argument and allow some unnecessary hypotheses to be
removed.
\section {Quadratic words}\label{qwords}

Let $X$ and $D$ be sets indexed by $\NN$. A {\em word} with
{\em coefficients} in $D$ and {\em variables} in $X$ is a finite
(possibly empty) sequence of elements of $(D\dcup X)^{\pm 1}$.
Suppose $w=w_1\cdots w_n$ is a word, where $w_i\in (D\dcup X)^{\pm
1}$, and let $a\in D\dcup X$. If $w_i=a^\e$, where $\e=\pm 1$,
then $w_i$ is called an {\em occurrence} of $a$ in $w$, and $a$ is
said to {\em occur} in $w$.

Define 
\begin{gather*}
L_X(w)=\{x\in X\,:\, x \textrm{ occurs in } w\},\\
L_D(w)=\{d\in D\,:\, d \textrm{ occurs in } w\}\textrm{ and}\\
L(w)=L_D(w)\dcup L_X(w).
\end{gather*} 
Given $a\in D\dcup X$ define $O_a(w)$
to be the number of occurrences of $a$ in $w$.
The word $w$ with coefficients in $D$ and variables in $X$ is $quadratic$ if
$O_x(w)=2$, for all $x\in L_X(w)$, and $O_d(w)=1$, for all $d\in L_D(w)$.
A list of $k
\ge 1$ words $\mbf w=(w_1,\dots w_k)$ with coefficients in $D$ and variables
in $X$ is called a {\em system} of words, of {\em dimension} $k$,
if  $L(w_i)\cap L(w_j)=\nul$, whenever $i\neq j$.
If $\mbf w=(w_1,\ldots ,w_k)$ is a system of words then we define
$L_D(\mbf w)=\cup^k_{i=1}L_D(w_i)$, $L_X(\mbf w)=\cup^k_{i=1}L_X(w_i)$ and
$L(\mbf w)=L_D(\mbf w)\dcup L_X(\mbf w)$.
An element $a\in L(\mbf w)$ is said to {\em occur} in $\mbf w$.
The {\em length} of $\mbf w$ is $|\,\mbf w|=\sum_{i=1}^k|w_i|$ and $\mbf w$ is
{\em quadratic} if $w_i$ is quadratic, for $i=1,\ldots k$.
%

A word $w$  in an alphabet $A\dcup A^{-1}$ is {\em cyclically
reduced} if no cyclic permutation of $w$ has a subword of the form
$aa^{-1}$ or $a^{-1}a$, with $a\in A$.
In \cite{\GL} Grigorchuk and Lysionok take a cyclically reduced quadratic word
$w=a_1^{\e_1}\cdots a_n^{\e_n}$, with coefficients in $D$ and variables in $X$, where
$a_i\in D\dcup X$ and $\e_i=\pm 1$,
and describe a surface $\S(w)$ associated to $w$. This is done as
follows. 
Let $I$ denote the unit interval
$I=\left [0,1\right ]$ oriented from $0$ to $1$. Choose a partition
$0=p_0,\ldots ,p_n=1$ of  $I$ into $n$ sub--intervals of equal length. Label the sub--interval
$\left [p_{i-1},p_{i}\right ]$
with $a_i$ and orient it from $p_{i-1}$ to $p_i$ if $\e_i=1$ and in
the opposite sense otherwise. Identify $p_0$ and $p_n$ to form a labelled oriented complex $B$
homeomorphic to $S^1$. Let $\D$ be a disk with boundary $\de \D$
and identify $\de \D$ with $B$.
Then $\S(w)$ is the quotient space of $\D$ formed
by identifying sub--intervals of $B$ which have the same label, respecting the orientation
of sub--intervals. The result is a compact surface $\S(w)$, which 
has boundary if and only if $L_D(w)$ is
non--empty. Furthermore $\S(w)$ is orientable if and only if every element of $L_X(w)$ occurs
once with positive exponent and once with negative exponent. We denote the number of boundary
components of $\S(w)$ by $\bt(w)$. Let $\S^\prime(w)$ be the surface obtained from $\S(w)$ by
capping off each boundary component with a  disk. Define the {\em genus} of $w$,
genus$(w)=\{(t,p)\in \NN^2\,:\,\S(w)$ is homeomorphic to the connected sum of $t$ torii
and $p$ projective planes$\}$.

If $\S(w)$ has boundary component $\bt$ then $\bt$ is divided into oriented sub--intervals
labelled by elements of $D$. Take an end point of one of these intervals as base point
of $\bt$. If $\S(w)$ is orientable we may choose an orientation on $\bt$ such that, travelling
round $\bt$ in the direction of orientation every sub--interval $[p_{i-1},p_i]$ that is
encountered is traversed from $p_{i-1}$ to $p_i$. Otherwise we fix
one of the two possible orientations of $\bt$. Then traversing $\bt$ once, 
in the direction of orientation, 
starting from the base point, 
and recording the labels as they occur, with positive exponent if the
sub--interval orientation agrees with that of $\bt$ and negative exponent otherwise, we obtain
a word which we call a {\em boundary label} of $\bt$. If $\S(w)$ has boundary components
$\bt_1,\ldots ,\bt_n$ with boundary labels $b_1,\ldots ,b_n$, respectively, then
$(b_1, \ldots ,b_n)$ is called a {\em boundary labels list} for $w$ of {\em length} $n$,
and $(w;b_1, \ldots ,b_n)$
is called a {\em labelled quadratic word}. Note that it is implicit in the definition of
labelled quadratic word that $w$ is cyclically reduced.

Let $F=F(D\dcup X)$ denote the free group on the set $D\dcup X$
and let Aut$_fF$ denote the subgroup of
the automorphism group of $F$ which consists of those  elements
which fix all but finitely many elements of $D\dcup X$.
Given an automorphism $\th \in \textrm{Aut}_f(F)$ let $\widetilde{u\th}$ denote the reduced
word representing $u\th\in F$. Let $(w;b_1,\ldots ,b_n)$ be a labelled quadratic word. If
$\wt{w \th}$ is a quadratic word, $\S(w)$ is homeomorphic to $\S(\wt{w\th})$ and
$\wt{w \th}$ has a boundary labels list $(\wt{b_1\th},\ldots ,\wt{b_n\th})$, then
$\th$ is called an {\em admissible transformation} of $(w;b_1,\ldots ,b_n)$ to
$(\wt{w \th};\wt{b_1\th},\ldots ,\wt{b_n\th})$.

Given non-negative integers $\xi,n,t,p$
define the quadratic word
\[q(\xi,n,t,p)=\prod_{i=1+\xi}^{n+\xi}x_i^{-1}d_ix_i\,
\prod_{i=1+\xi}^{t+\xi}\left [x_{i+n},x_{i+n+t}\right ]
\prod_{i=1+\xi}^{p+\xi}x_{i+n+2t}^2\,.\]
Let $n$ be a non--negative integer. A {\em partition} of $n$ of {\em
length} $k$ is a $k$--tuple $(n_1,\ldots,n_k)$ of integers such that
$n_i\ge 0$ and $\sum^k_{i=1}n_i=n$. 
Let  $\mbf n=(n_1,\ldots, n_k)$, $\mbf t=(t_1,\ldots t_k)$ 
and $\mbf p=(p_1,\ldots,p_k)$ be partitions of length $k>0$, of
integers $n$, $t$ and $p$ respectively, such that
$n_j+t_j+p_j>0$, for $1\le j\le k$. Then we say that $(\mbf n,\mbf t,\mbf p)$ is
a {\em positive} $3${\em --partition} of $(n,t,p)$ (of {\em length} $k$).
If $\mbf h$ is any partition, of length $k$ of a non--negative integer $h$
and  $(\mbf n,\mbf t,\mbf p)$ is a positive $3$--partition then we say that
$(\mbf h,\mbf n,\mbf t,\mbf p)$ is a {\em positive} $4${\em --partition} of $(h,n,t,p)$.
Given partitions 
$\mbf{h^\prime}=(h^\prime_1,\ldots,h^\prime_l)$, 
$\mbf{n^\prime}=(n^\prime_1,\ldots,n^\prime_l)$, 
$\mbf{t^\prime}=(t^\prime_1,\ldots,t^\prime_l)$ and 
$\mbf{p^\prime}=(p^\prime_1,\ldots,p^\prime_l)$ such that 
$(\mbf h^\prime,\mbf n^\prime,\mbf t^\prime,\mbf p^\prime)$ is a positive 
$4$--partition, we
say that  
$(\mbf h^\prime,\mbf n^\prime,\mbf t^\prime,\mbf p^\prime)
\le (\mbf h,\mbf n,\mbf t,\mbf p)$
if there exist integers
$j_1,\ldots ,j_l$ such that $1\le j_1<\cdots<j_l\le k$ with 
$h^\prime_i\le h_{j_i}$, $n^\prime_i\le n_{j_i}$,
$t^\prime_i\le t_{j_i}$ and $p^\prime_i\le p_{j_i}$,
for $i=1,\ldots, l$. If, in addition, $t^\prime_i= t_{j_i}$ and $p^\prime_i= p_{j_i}$,
for $i=1,\ldots, l$, we say that 
$(\mbf h^\prime,\mbf n^\prime,\mbf t^\prime,\mbf p^\prime)$ is 
a {\em sub--partition} of $(\mbf h,\mbf n,\mbf t,\mbf p)$ and write
$(\mbf h^\prime,\mbf n^\prime,\mbf t^\prime,\mbf p^\prime)
\le_g (\mbf h,\mbf n,\mbf t,\mbf p)$.
 
Given a positive $3$--partition $(\mbf n,\mbf t,\mbf p)$, of length $k$ as above, define
\begin{align}
\xi_1 &=\xi_1(\mbf n,\mbf t,\mbf p) =0 \textrm{ and }\nonumber\\ 
\xi_{j+1} &= \xi_{j+1}(\mbf n,\mbf t,\mbf p) =\xi_{j}+n_j+2t_j+p_j, \textrm{ for } j=1,\ldots ,k-1.\label{xij}
\end{align}
Next define the system
of quadratic words
\begin{equation}\label{qntp}
\mbf q (\mbf n,\mbf t,\mbf p)=(q(\xi_j,n_j,t_j,p_j)\,:\,j=1,\ldots ,k).
\end{equation}
Note that $q(\xi_j,n_j,t_j,p_j)$ has boundary labels list $(d_{1+\xi_j},\ldots ,d_{n_j+\xi_j})$.
A system $\mbf w=(w_1,\ldots ,w_l)$ of  quadratic words is in {\em standard form} if
$\mbf w = \mbf q(\mbf n,\mbf t,\mbf p)$, for some positive $3$--partition 
$(\mbf n,\mbf t,\mbf p)$.

Suppose that $\mbf w=(w_1,\ldots ,w_k)$ is a system of  quadratic words. We define the 
{surface associated to} $\mbf w$ to be the compact surface
\[\S(\mbf w)=\bdcup^k_{i=1}\S(w_i).\] 
Assume that $w_i$
has boundary labels list $\mbf b^i=(b^i_1,\ldots ,b^i_{n_i})$ of length $n_i$.
We denote $((w_1;\mbf b^1),\ldots ,(w_k;\mbf b^k))$ by $(\mbf w; \mbf b)$. We call
$(\mbf w;\mbf b)$ a {\em system of labelled quadratic words} which we say is in
{\em standard form} if $\mbf w$ is in standard form.
Given $\th\in \textrm{Aut}_f(F)$, let $\wt{\th(\mbf b^i})$ denote
$(\wt{\th(b^i_1 )},\ldots ,\wt{\th(b^i_{n_i})})$ and $(\wt{\th(\mbf w)};\wt{\th(\mbf b)})$ denote
$((\wt{\th(w_1)};\wt{\th(\mbf b^1)}),\ldots ,(
\wt{\th(w_k )};\wt{\th(\mbf b^k )}))$. Then
$\th$ is called an
{\em admissible transformation} of $(\mbf w;\mbf b)$ to $(\wt{\th(\mbf w )};\wt{\th(\mbf b )})$ if
$\th$ is an admissible transformation of $(w_i;\mbf b^i)$ to $(\wt{\th(w_i)};\wt{\th(\mbf b^i )})$,
for $i=1,\ldots ,k$, and $L(\th(x))\cap L(\th(y))=\nul$ whenever $x\in
L(w_i)$, $y\in L(w_j)$ and $i\neq j$.

The following proposition can be proved exactly as the case $k=1$ of \cite[Proposition 2.1]{\GL}.

\begin{prop}[cf. {\cite[Proposition 2.1]{\GL}}]\label{stdform0}
Using the notation above, let $(\mbf w;\mbf b)$ be a system of labelled quadratic words,
of dimension $k$,
let $n_i=\bt(w_i)$ and let $(t_i,p_i)\in$ genus$(w_i)$,
for $i=1,\ldots ,k$. Let $\mbf n=(n_1,\ldots ,n_k)$, $\mbf t=(t_1,\ldots t_k)$ and
$\mbf p=(p_1,\ldots ,p_k)$. Then there exists an automorphism
$\eta\in \textrm{Aut}_f(F(D\dcup X))$ such that $\eta$ is an admissible transformation
of $(\mbf w;\mbf b)$ to $(\mbf q(\mbf n,\mbf t,\mbf p);\mbf d)$ and $\eta(F(D))=F(D)$.
Furthermore $\eta$ can be effectively constructed.
\end{prop}

\section{Quadratic Exponential Equations}\label{qpequations}

Let $I$ be a recursive set and, for all $i\in I$, let
$\langle C_i|R_i \rangle$ be a recursive presentation for a group
$H_i$ and let $H=*_{i\in I}H_i$. Let $J$ be a recursively enumerable set and, for all $j\in J$,
let $s_j\in H$, of free--product length at least $2$,  and let $N$ be the normal closure of
$\{s_j\}_{j\in J}$ in $H$.

Let $C=\dcup _{i\in I} C_i$, $R=\dcup _{i\in I}R_i$ and, for each $j\in J$,
choose $\wt{s_j}\in F(C)$ such that $\wt{s_j}$ represents $s_j\in H$. Let
$\mbf{\wt{s}}=\{\wt{s_j}\}_{j\in J}$. Then $H$ has a presentation
$\langle C|R\rangle$ and $G=H/N$ has a presentation
$\langle C|R\cup \mbf{\wt{s}} \rangle$, which we call the {\em
  natural} presentations of $H$ and  $G$, respectively. Given groups $H_i$ with fixed
presentations as above and groups $H$ and $G$ as described we shall always assume
$H$ and $G$ are given by fixed natural presentations. 
By a {\em word} in $H$ we mean an element of $H$ written in
free product normal form. The {\em length}, $l(h)$, of an element $h$ of $H$ refers to the 
free product
length of a word in the normal form for $H$ representing $h$
(as opposed to the length as an element
of $F(C)$). 
If $u$ and $v$ are words in $H$ and $u$ is a subword of some cyclic permutation
of $v$ then we say that $u$ is a {\em cyclic subword} of $v$.

A {\em quadratic equation over} $G$ is a pair of the form $(w=1,\bt)$, where
$w$ is a quadratic word with coefficients in $D$ and variables in $X$, and 
$\bt$ is a map $\bt:L_D(w)\maps H$. A {\em solution} to the quadratic equation $(w=1,\bt)$
consists of a homomorphism $\phi:F(L(w))\maps H$ such that $\phi(w)\in N$ and 
$\phi(d)=\bt(d)$, for all $d\in L_D(w)$. In this section we define a generalisation
of quadratic equations to what we shall call quadratic exponential equations. 

We assume that $H_i$ has solvable word problem, for all $i\in I$, and
consequently $H$ has solvable word problem: that is, that there exists
an algorithm which, given $w\in F(C)$ will determine whether or
not $w=1$ in $H$ and, if so, will output a product
$\sum^n_{k=1}u_k^{-1}r_k^{\e_k}u_k=w$
in $F(C)$, where $u_k\in F(C)$, $r_k\in R$ and $\e_k=\pm 1$.
Given this assumption we may formulate equations and decision problems for 
$\langle C|R\cup\mbf{\wt{s}}\rangle$ in terms of words in the free product $H$ instead of 
words in $F(C)$.

Let $\Lm=\{\lm_1,\ldots \}$ be a set indexed by $\NN$, which we call a
set of {\em parameters}. Let $N$ be the free $\ZZ$--module on $\Lm$ and let $M$ be the free
$\ZZ$--module $\ZZ\oplus N$. Elements of $M$ are regarded as linear polynomials over $\Lm$.
An ordered pair $(h,m)$, $h\in H$, $m\in M$, $m\notin \ZZ$, $h$ cyclically reduced is
called
a {\em proper exponential} $H${\em --letter} and the set of all such ordered
pairs is denoted $(H,\Lm)_0$. An ordered pair $(h,n)$, $h\in H$, $l(h)=n\in
\ZZ$ is called a {\em degenerate exponential} $H${\em --letter} and the
set of all such ordered pairs is denoted $(H,\Lm)_1$.
The set $(H,\Lm)_0\cup(H,\Lm)_1$ is called the
set of {\em exponential} $H${\em --letters} and denoted $H^\Lm$.
We regard elements of $F(H^\Lm)$ as reduced words on $H^\Lm\dcup (H^{\Lm})^{-1}$.
Given $h\in H$ we identify $h$ with $(h,l(h))$, so as a set $H\subset H^\Lm$.
An element $(h,f)\in H^\Lm$ such that $h\in H_i$, for some $i\in I$, is
called a {\em minor} element of $H^\Lm$

If $h\in H$ and $r\in \ZZ$ with $l(h)\ge r\ge 0$, define $\io(h,r)$ to
be the maximal initial subword of $h$ of length $r$ and $\t(h,r)$ to be
the maximal terminal subword of $h$ of length $r$. Given $a\in H$, with
$l(a)=m$, and $\al\in \ZZ$, let $q,r\in \ZZ$,
such that $|\,\al |=qm+r$ and $0\le r< m$, and  define
\[a\^\al=\left\{\begin{array}{ll}
a^q\io(a,r) & \textrm{if } \al\ge 0\\
a^{-q}\io(a^{-1},r) & \textrm{if } \al<0
\end{array} \right. .\]

We have $a\^ql(a)=a^q$ so, if $l(a)>1$ and $\al\neq 0$, then $a^\al\neq
a\^\al$. Also if $a$ is cyclically reduced
\[\begin{array}{lcl}
(a\^\al)^{-1}&=&  \left\{\begin{array}{ll}
(\io(a,r)^{-1}\t(a,m-r)^{-1})\^\al & \textrm{if } \al\ge 0\\
(\t(a,r)\io(a,m-r))\^(-\al) & \textrm{if } \al< 0
\end{array} \right. ,\\ \\
&=& \left\{\begin{array}{ll}
(\t(a,m-r)\io(a,r))\^(-\al) & \textrm{if } \al\ge 0\\
(\t(a^{-1},m-r)\io(a^{-1},r))\^\al & \textrm{if } \al<0
\end{array} \right. .
\end{array}\]

Suppose $a$ is cyclically reduced, $|\,\al|=ql(a)+r$, $0\le r < l(a)$,
$a=a_0a_1$, with $l(a_0)=r$, $l(a_1)=l(a)-r$,
and $a=b_0b_1$ with $l(b_1)=r$, $l(b_0)=l(a)-r$. Then
\[a\^\al=\left\{\begin{array}{ll}
a^qa_0 & \textrm{if } \al\ge 0\\
a^{-q}b^{-1}_1 & \textrm{if } \al < 0
\end{array} \right. , \]
and
\[a^{-1}\^\al=\left\{\begin{array}{ll}
a^{-q}b^{-1}_1 & \textrm{if } \al\ge 0\\
a^qa_0 & \textrm{if } \al<0
\end{array} \right., \]
so $a\^(-\al)=a^{-1}\^\al$.

A {\em parameter system} over $\Lm$ is a finite system of linear
diophantine equations, congruences and inequalities in the variables
$\Lm$. A {\em solution} to the parameter system $\cL$ is a
retraction $\al$ from $M$ to $\ZZ$ (i.e. a $\ZZ$--module homomorphism from
$M$ onto $\ZZ$, fixing $\ZZ$ pointwise) such that substitution of
$\al(\lm)$ for $\lm$, for each $\lm\in \Lm$, solves $\cL$. 
A parameter system is {\em consistent} if it has a solution.
Let $\al$ be a $\ZZ$--module retraction from $M$ onto $\ZZ$. We define
the {\em parameter evaluation homomorphism} $\hat\al:F(H^\Lm)\maps H$
corresponding  to $\al$, by setting $\hat\al(h,f)=h\^\al(f)$, for $(h,f)\in
H^\Lm$.

Next we generalise the notation for presentations of $H$ and $G$ to
collections of such groups and presentations. Let $I$, $H_i$, $\langle C_i|R_i\rangle$,
$\langle C|R\rangle$ and $H$ be defined as above. Let $K$
be a recursive subset of $\NN$ and let $X_k\subseteq I$ be a recursive set, for
all $k\in K$. Define 
\begin{equation}\label{Hsys}
H_{*,k}=\ast_{j\in X_k}H_j,\textrm{ for } k\in K.
\end{equation} 

Let $\mbf w=(w_1,\ldots, w_k)$ be  a system of  cyclically reduced quadratic words,
with coefficients in $D$ and variables in $X$, let $\bt:F(L_D(\mbf
w))\maps F(H^\Lm)$ be a homomorphism and let $a_1, \ldots ,a_k$ be
elements of $K$. If $\bt(d)\in F(H_{*,a_i}^\Lm)$, for all $d\in
L_D(w_i)$, $i=1,\ldots ,k$, then the ordered pair
$(\mbf w,\bt)$ is called a {\em quadratic exponential word} over $((X_k)_{k\in K};(H_i)_{i\in I})^\Lm$ 
with {\em basis} $(a_1,\ldots ,a_k)$. 

We define a {\em quadratic} {\em exponential} {\em equation}   over
$((X_k)_{k\in K};(H_i)_{i\in I}))^\Lm$ 
with {\em basis}  $(a_1,\ldots ,a_k)$
to be a pair $(\mbf w=1,\bt)$, where $(\mbf w,\bt)$ is a
quadratic exponential word  over 
$((X_k)_{k\in K};(H_i)_{i\in I})^\Lm$ with basis $(a_1,\ldots ,a_k)$. 

For each $k\in K$ let $\mbf{s}_k$ be a recursively enumerable subset of $H_{*,k}$,
such that each element $s\in \mbf{s}_k$ satisfies $l(s)\ge 2$,
and define 
\begin{equation}\label{Gsys}
G_k=H_{*,k}/N_k,
\end{equation}
where
$N_k$ is the normal closure of  $\mbf{s}_k$ in $H_{*,k}$.
We allow the possibility that $\mbf s_k$ is empty, so it may be that 
$G_k=H_{*,k}$. 
Given a quadratic exponential equation $W=(\mbf w=1,\bt)$ 
over $((X_k)_{k\in K};(H_i)_{i\in I})^\Lm$ with basis $(a_1,\ldots ,a_k)$ a
{\em solution}  to $W$
in $(G_k)_{k\in K}$ is a 
pair $(\phi,\al)$ such that
$\al$ is a retraction from $M$ to $\ZZ$, $\phi :F(L(\mbf w))\maps H$ is
a homomorphism, $\phi(a)\in H_{*,a_i}$, for all $a\in L(w_i)$, 
$\phi(w_i)\in N_{a_i}$, for $i=1,\ldots ,k$, and
$\phi(d)=\hat{\al}(\bt(d))$, for all $d\in L_D(\mbf w)$.

If $W$ is an equation over $((X_k)_{k\in K};(H_i)_{i\in I})^\Lm$ for which solutions are
sought in
$(G_{k})_{k\in K}$ then we say that $W$ has {\em environment} 
\[((X_k)_{k\in K};(H_i)_{i\in I} ;(G_k)_{k\in K})^\Lm.\] To specify an equation fully its environment
and
basis must be given, but we omit these when the meaning is clear.
If, in the definition of the
$H_{*,k}$, the set $K$ has a unique element $s$ and $X_s=I$ then we refer to quadratic equations with
environment
$((X_k)_{k\in K};(H_i)_{i\in I} ;(G_k)_{k\in K})^\Lm$
as equations over $H^\Lm$ with solutions in $G$.

Let $u \in F(H^\Lm)$ and let $\cL$ be a consistent parameter system. We say
$u$ is {\em constrained} by $\cL$ if $\cL$ implies that $f\equiv k (\mod
l(a))$, for some $k\in \ZZ$, for all proper exponential $H$--letters
$(a,f)\in H^\Lm$ occurring in $u$, with $l(a)\ge 2$.
A {\em constrained quadratic exponential equation} consists of a triple
$(\mbf w=1,\bt,\cL)$, where $(\mbf w=1,\bt)$ is a quadratic exponential
equation, $\cL$ is a consistent parameter system and $\bt(d)$ is constrained
by $\cL$, for all $d\in L_D(\mbf w)$. A quadratic exponential
equation  $(\mbf w=1,\bt)$ or a constrained quadratic exponential equation  
$(\mbf w=1,\bt,\cL)$
is in {\em standard form} if $\mbf w$ is in standard form.
A {\em solution}
to the constrained quadratic exponential equation $(\mbf w=1,\bt,\cL)$ is
a solution $(\phi,\al)$ to $(\mbf w=1,\bt)$ such that $\al$ is a solution
to $\cL$. 

An element $(a,f)\in H^\Lm$ is said to {\em occur} in 
the system of quadratic words $(\mbf w,\bt)$ 
if $(a,f)$ occurs in $\bt(d)$ for some
$d\in L_D(\mbf w)$.
A polynomial $f\in M$ is said to {\em occur}  in 
the quadratic exponential
equation $(\mbf w=1,\bt)$, or in the constrained quadratic exponential equation $(\mbf w=1,\bt,\cL)$, 
if $(a,f)\in H^\Lm$ occurs in  $(\mbf w,\bt)$,
for some $a\in H$.
Also a polynomial $f$ is said to {\em occur}
in the parameter system $\cL$ if $\cL$ contains an equation $f=0$ or a
congruence $f\equiv 0 (\mod k)$, for some $k\in \ZZ$.
A polynomial $f$ is said to {\em occur} in the constrained quadratic
exponential equation $(\mbf w=1,\bt,\cL)$ if $f$ occurs either in $(\mbf
w=1,\bt)$ or in $\cL$.
An element $\lm\in \Lm$ is said to {\em occur} in $f\in M$ if the coefficient
of $\lm$ in $f$ is non--zero. Also $\lm$ is said to {\em occur} in a quadratic
exponential equation $(\mbf w=1,\bt,\cL)$ if $\lm$ occurs in an element $f\in M$
which occurs in $(\mbf w=1,\bt,\cL)$.

The $H${\em --length} of a 
exponential $H$--letter $(a,f)$ is
$l_H(a,f)=l(a\^ f)$.
The $H${\em --length} of a reduced word 
$u=(a_1,f_1)^{\e_1}\ldots (a_n,f_n)^{\e_n}\in F(H^\Lm)$, 
where $\e_i=\pm 1$, is defined to be
\begin{equation}
\label{H-length}
l_H(u)=\sum_{i=1}^n l_H(a_i,f_i).
\end{equation} 
The {\em exponential length} of $u$ is defined to
be 
\begin{equation}
\label{exponential-length}
|u|_\Lm=n-|\,\{i:1\le i \le n \textrm{ and } (a_i,f_i) \textrm{ is
  degenerate}\}|.
\end{equation} 
If $\al:M\maps \ZZ$ is a retraction the {\em exponent length} of $\al$ on $u$ is 
\begin{equation}
\label{exponent-length}
|u|_\al=\sum_{i=1}^n |\,\al(f_i)|.
\end{equation} 
The $H^\Lm${\em --length} of $u$ is defined to be 
\begin{equation}
\label{HL-length}
|u|=n.
\end{equation} 
Let $W=(\mbf w=1,\bt,\cL)$ be a constrained quadratic exponential equation, with
$\mbf w=(w_1,\ldots, w_k)$.
Then the $H${\em --length}, {\em exponential coordinates}, {\em exponential length} and 
$H^\Lm${\em --length} of $W$, are 
\begin{align*}
l_H(W)&=\sum_{d\in L_D(\mbf w)}l_H(\bt(d)),\\[0.5em]
(W)_\Lm&=\left(\sum_{d\in L_D(w_1)}|\,\bt(d)|_\Lm, \ldots, \sum_{d\in L_D(w_k)}|\,\bt(d)|_\Lm\right),\\[0.5em]
|W|_\Lm&=\sum_{d\in L_D(\mbf w)}|\,\bt(d)|_\Lm \textrm{ and }\\[0.5em]
%
|W|&=\sum_{d\in L_D(\mbf w)}|\,\bt(d)|, 
\end{align*}
respectively. 

A solution $(\phi,\al)$ to $W$ is said to be {\em exponent--minimal} if, given any solution 
$(\phi^\prime,\rho)$ to $W$, the inequality
\[\sum_{d\in L_D(\mbf w)}|\,\bt(d)|_\rho\ge\sum_{d\in L_D(\mbf w)}|\,\bt(d)|_\al \]
holds.

Let $z\in \ast_{i\in I} H_i$, with $l(z)=n$, $z=u_1\cdots u_n$, where $u_j\in H_{i_j}$,
for some $i_j\in I$ and $j=1,\ldots ,n$. Define  supp$(z)$ the {\em support} of $z$ to 
be $\{i_1,\ldots, i_n\}$, so $z\in *_{j\in \textrm{ supp}(z)}H_j$. We define supp$(z,f)$
the {\em support}  of $(z,f)\in H^\Lm$ to be supp$(z)$. If $w\in F(H^\Lm)$,
$w=w_1^{\e_1}\cdots w_n^{\e_n}$, with $w_j\in H^\Lm$ and $\e_j=\pm 1$, we define supp$(w)$ the {\em support} 
of $w$ to be $\cup_{j=1}^n\textrm{ supp}(w_j)$. Let $\mbf w=(w_1,\ldots ,w_n)$ be a system of words. 
If $W=(\mbf w =1, \bt,\cL)$ is a constrained quadratic
exponential equation with environment
\begin{equation}\label{big_env}
((X_k)_{k\in K};(H_i)_{i\in I} ;(G_k)_{k\in K})^\Lm
\end{equation}
and basis $(a_1,\ldots ,a_k)$ then we define the {\em support} of $W$
\[\textrm{supp}(W)=(X_{a_1}^\prime,\ldots ,X_{a_k}^\prime),\] where
\[X^\prime_{a_i}= \bigcup \{\textrm{supp}(\bt(d)):d\in L_D(w_i)\}\cup \bigcup\{\textrm{supp}(s):s\in \mbf s_{a_i }\}.\]
Define $I^\prime=\cup_{i=1}^kX_{a_k}^\prime$, $\quad H^\prime_{*,a_i}=*_{j\in X^\prime_{a_i}}H_j$, $\quad N^\prime_{a_i}$ to be the normal closure of 
$\mbf s_{a_k}\in H^\prime_{*,a_i}$ and $G^\prime_{a_i}=H^\prime_{*,a_i}/N^\prime_{a_i}$, for $i=1,\ldots ,k$.
Then $W$ is a constrained quadratic exponential equation with environment 
\begin{equation}\label{small_env}
((X^\prime_{a_j})_{j=1,\ldots, k};(H_i)_{i\in I^\prime} ;(G^\prime_j)_{j=1,\ldots,k})^\Lm
\end{equation}
and basis  $(a_1,\ldots ,a_k)$. Every solution to $W$ with environment (\ref{small_env}) is a 
solution to $W$ with environment (\ref{big_env}). The converse also holds as in the next lemma.
\begin{lemma}\label{basis_reduce}
Every solution to $W$ with environment (\ref{big_env}) determines canonically a solution
to $W$ with environment (\ref{small_env}).
\end{lemma}
{\em Proof.} Let $(\phi, \al)$ be a solution to $W$ with environment (\ref{big_env}). Fix $i\in \{1,\ldots , k\}$,
suppose $X_{a_i}^\prime=\{i_1,\ldots, i_s\}$ and set $X^{\pprime}=X_{a_i}\backslash X_{a_i}^\prime$.
Let \[A=*_{j=1}^s H_{i_j}\textrm{ and } B=*_{j\in X^\pprime} H_{j}.\]
Denote by $f$ the canonical map $f:A*B\maps A$ such that ker$(f)$ is the normal closure of $B$ in $A*B$.
Then $N_{a_i}$ is by definition the normal closure of $\mbf s_{a_i}$ in $A*B$ and $\mbf s_{a_i}\subseteq A$, so
$f$ maps $N_{a_i}$ into $N^\prime_{a_i}$. Define $\phi^\prime:L(w_i)\maps A$ by 
$\phi^\prime(a)=f\circ \phi(a)$, for all $a\in L(w_i)$.
If $d\in L_D(w_i)$ then $\phi(d)=\hat\al\bt(d)$ and supp$(\bt(d))\subseteq X^\prime_{a_i}$ so $\phi(d)\in A$.
Hence $\phi^\prime(d)=\hat\al\bt(d)$, for all $d\in L_D(w_i)$. Also $\phi^\prime(w_i)\in N^\prime_{a_i}$ as
$\phi(w_i)\in N_{a_i}$. As $L(w_i)\cap L(w_j)=\emptyset$, when $i\neq j$, we can define $\phi^\prime$ in this way
simultaneously over $L(w_i)$, for $i=1,\ldots , k$. Then $(\phi^\prime,\al)$ is a solution to $W$ with environment
(\ref{small_env}).
\begin{exx}\label{exx_eqn}
Let $I=\NN$ and in each example below let $H_i=\langle C_i|R_i \rangle$ be trivial unless $C_i$ and $R_i$ are 
given explicitly.
\be
\item\label{exx_eqn_1} Let $K=\{1,2\}$, $X_1=\{2\}$ and $X_2=\{1,2\}$. Let $C_i=\{c_{i1},c_{i2}\}$, for $i=1,2$, 
$R_1=\{c_{11}^3,c_{12}^3\}$, 
$R_2=\{c_{21}^2,c_{22}^2\}$, 
$\mbf s_1=\emptyset$ and 
$\mbf s_2=\{c_{11}c_{22}c_{21}c_{12}c_{22}c_{21}\}$. 
Then $H_{*,1}=H_2$, $H_{*,2}=H_1*H_2$, $G_1=H_2$ and $G_2=H_{*,2}/\langle\langle \mbf s_2\rangle\rangle$. 
Let $w_1=d_1$, 
$w_2=d_2x_3^{-1}d_3x_3$ and 
$\mbf w=(w_1,w_2)$. Let 
$\bt(d_1)= (c_{22},\lm_1)$,  
$\bt(d_2)=(c_{12}c_{22}c_{21},2)(c_{11}c_{22}c_{21},f_1)(c_{11}c_{22},2)$, and 
$\bt(d_3)=(c_{21}c_{22}c_{11}^2,f_2)(c_{22}c_{11}^2,2)$, 
where 
$f_1=\lm_1-2$ and 
$f_2 = \lm_1-\lm_2$. 
Let $\cL_0$ consist of the congruences
$\lm_i \equiv 0 (\mod 2)$ and the inequalities $\lm_i>0$, for $i=1,2$. Let $\cL_1$ consist of the same congruences
and the inequalities $\lm_1>4$, $\lm_2>2$. 
Then $W_i=(\mbf w =1,\bt,\cL_i)$ is  a quadratic 
exponential equation with environment $(X_1,X_2;(H_i)_{i\in I};G_1,G_2)^\Lm$ and basis $(1,2)$. 
These equations also have environment
$(X_1,X_2;H_1,H_2;G_1,G_2)^\Lm$. 
Let $\al_0$ be a retraction of $M$ to $\ZZ$ such that $\al_0(\lm_1)=4$ and $\al_0(\lm_2)=4$. 
Let $\phi_0(x_3)=c_{11}c_{22}$ and set $\phi_0(d_i)=\hat\al_0\bt(d_i)$, for $i=1,2,3$. 
Then $(\phi_0,\al_0)$ is a solution to $W_0$ but not $W_1$.
Let $\al_1$ be such that $\al_1(\lm_1)=6$ and $\al_1(\lm_2)=4$. 
Let $\phi_1(x_3)=c_{11}c_{22}$ and set $\phi_1(d_i)=\hat\al_1\bt(d_i)$, for $i=1,2,3$.
Then $(\phi_1,\al_1)$ is a solution to $W_0$ and $W_1$.
\item Let $C_i,R_i,H_i$ be as in Example \ref{exx_eqn}(\ref{exx_eqn_1}), 
for $i=1,2$ and let $C_j,R_j$ be arbitrary for $j\ge 3$. 
Let $K=\NN$ and $X_1=\{1,2\}$, $X_2=\{1,2,3\}$ and $X_j=\{j,j+1\}$, for $j\ge 3$, so
$H_{*,1}=H_1*H_2$, $H_{*,2}=H_1*H_2*H_3$ and $H_{*,j}=H_j*H_{j+1}$, $j\ge 3$. Let $\mbf s_i$ be as in  
Example \ref{exx_eqn}(\ref{exx_eqn_1}), for 
$i=1,2$ and let $\mbf s_j=\emptyset$, for $j\ge 3$. For $i=0$ or $1$, the equation $W_i$ of 
Example \ref{exx_eqn}(\ref{exx_eqn_1}) is an equation with 
environment $((X_i)_{i\in K};(H_i)_{i\in I}; (G_i)_{i\in K})^\Lm$ and basis $(1,2)$, with $G_i$ defined as usual. 
The support of $W_i$ is
$(X^\prime_1,X^\prime_2)$, where $X_i^\prime$ is the set $X_i$ of 
Example \ref{exx_eqn}(\ref{exx_eqn_1}). Therefore $W_i$ is 
an equation with environment $(X^\prime_1,X^\prime_2;H_1,H_2;G^\prime_1,G^\prime_2)^\Lm$ 
and basis $(1,2)$, for $G^\prime_i$ defined as usual. Moreover $W_i$ has
a solution with respect to the first named environment if and only if it has a solution with respect to the latter.
\ee
\end{exx}
The next example illustrates equations corresponding to some common group theoretic questions.
\begin{exx}\label{exx_eqn_probs} Let $I=\NN$, and let $C_i$, $R_i$ (and $H_i$) be given, for all $i\in I$.
Let $K=\{1\}$, $X_1$ be a recursive subset of $I$ and $\mbf s_1\subseteq H_{*,1}$ and $G_1$ be defined as usual.
Let $\mbf w=(w_1)$ and let $W=(\mbf w=1,\bt,\cL)$, where $w_1$, $\bt$ and $\cL$ are described below. 
The equation $W$ will in each case have environment $(X_1;(H_i)_{i\in I};G_1)^\Lm$ and basis $(1)$.
\be
\item\label{exx_eqn_probs_1} Let $w_1=d_1$ and $\bt(d_1)=(h_1,\lm_1)$, where $h_1\in H_{*,1}$.
Assume $\cL$ contains the subsystem $\cL_0$ consisting of the inequality $\lm_1\neq 0$ and the congruence
$\lm_1\equiv 0 (\mod{l(h_1)})$. 
As no elements of $X$ occur in $W$, given a solution $\al$ to $\cL$ 
we may determine whether there is a solution $(\phi,\al)$ to $W$ by setting $\phi=\hat\al\bt$
and evaluating $\phi(w_1)$. If there is a solution then we have $h_1^{m}=1$ in $G_1$, where
$m l(h_1)=\al(\lm_1)$. When $\cL=\cL_0$ finding a solution is equivalent to showing that $h_1$ has finite order.
\item\label{exx_eqn_probs_2} With $w_1$ and $\cL$ as in the previous case let
$\bt(d_1)=(h_1,\lm_1)(h_2,l(h_2))$, where $h_1,h_2\in H_{*,1}$.
Again, given a solution $\al$ to $\cL$ 
we may determine whether there is a solution $(\phi,\al)$ to $W$ by setting $\phi=\hat\al\bt$
and evaluating $\phi(w_1)$. If there is a solution then we have $h_1^{m}=h_2$ in $G_1$, with $m$ as before a
non--zero integer.
\item\label{exx_eqn_probs_3}  Let $w_1=d_1x_2^{-1}d_2x_2$ and $\bt(d_i)=(h_i,\lm_i)$, where $h_i\in H_{*,1}$,
for $i=1,2$.
Assume $\cL$ contains the subsystem $\cL_0$ consisting of the inequalities $\lm_i\neq 0$ and the congruences
$\lm_i\equiv 0 (\mod{l(h_i)})$, for $i=1,2$. 
 If there is a solution then $h_1^{m}$ is conjugate to $h_2^n$ in $G_1$, for some non--zero integers
$m$ and $n$. When $\cL=\cL_0$ finding a solution is equivalent to showing that some power of $h_1$ 
is conjugate to a power of $h_2$. If $\cL$ is the union of $\cL_0$ with the equation $\lm_2=l(h_2)$ then finding
a solution is equivalent to showing that some non--trivial power of $h_1$ is conjugate to $h_2$.
\ee
\end{exx} 
\subsection{Quadratic equations and genus}\label{qpe_and_genus}
 Let $I$ and $K$
be recursive sets and for each $k\in K$ let $X_k$ be a recursive subset of $I$.
Let $H_{*,k}$ be the free--product defined in (\ref{Hsys}), for all $k\in K$.
Let $\mbf s_k$ be a recursively enumerable subset of $H_{*,k}$ and define $G_k$ as in 
(\ref{Gsys}).
Let $\cL$ be a consistent parameter system. Let 
$\mbf z=(z_1,\ldots,z_n)$ be an $n$--tuple of elements of
$F(H^\Lm)$ such that $z_i$ is  constrained by $\cL$ and supp$(z_i)\subseteq X_{a_i}$, for some $a_i\in K$, for 
$i=1,\ldots ,n$. 
The set of all such $\mbf z$ is denoted by 
\[(\dcup_{k\in K}H_{*,k}^\Lm,\cL).\]

Let 
$\mbf n=(n_1,\ldots, n_k)$, $\mbf t=(t_1,\ldots t_k)$ and $\mbf p=(p_1,\ldots,p_k)$
be $k$--tuples of integers such that $(\mbf n,\mbf t,\mbf p)$ is a positive $3$--partition
of $(n,t,p)$, for some integers $t$ and $p$.
Define 
$\xi_j=\xi_j(\mbf n,\mbf t,\mbf p)$, for $1\le j\le k$, as in (\ref{xij}) and define
\begin{align}\label{muj} 
\mu_1&=\mu_1(\mbf n)=0 \textrm{ and }\nonumber\\
\mu_{j+1}&=\mu_{j+1}(\mbf n)=\mu_j+n_j,
\textrm{ for }j=1,\ldots ,k.
\end{align}
If $\mbf z\in (\dcup_{k\in K}H_{*,k}^\Lm,\cL)$ and there exist
$a_1,\ldots, a_k\in K$ such that supp$(z_i)\subseteq X_{a_j}$, for all $i$ with 
$1+\mu_j\le i\le \mu_{j+1}$ and for $j=1,\ldots ,k$,  then we say that $\mbf z$
is {\em supported} over $\mbf n$ with {\em basis}
$(a_1,\ldots ,a_k)$. 

Assume that $\mbf z$ is supported
over
$\mbf n$ with basis $(a_1,\ldots ,a_k)$.  
Define a map $\bt=\bt(\mbf z,\mbf n,\mbf t,\mbf p)$ from the subset
$\cup_{j=1}^{k}\{d_{\xi_j+1},\ldots,d_{\xi_j+n_j}\}$ of $D$ to   
$F(H^\Lm)$  
by
\begin{equation}\label{bjz} 
\bt(d_{\xi_j+i})= z_{\mu_j+i}, \textrm { for } j=1,\ldots,k, 
\end{equation}
and let 
$\mbf q=\mbf q(\mbf n,\mbf t,\mbf p)$.  The 
{\em quadratic exponential equation associated to} 
$(\mbf z,\cL,\mbf n,\mbf t,\mbf p)$ is defined to be the equation 
$(\mbf q=1,\bt,\cL)$ which we denote by
$Q=Q(\mbf z,\cL,\mbf n,\mbf t,\mbf p)$. Note that
$Q$ is an equation with environment
$((X_k)_{k\in K};(H_i)_{i\in I};(G_k)_{k\in K})^\Lm$
and basis $(a_1,\ldots ,a_k)$.
\begin{defn}\label{L_genus}
The $\cL${\em --genus} of 
$\mbf z$ in $(G_k:k\in K)$, denoted $\cL${--genus}$(\mbf z)$ is the set of
positive $3$--partitions $(\mbf n,\mbf t,\mbf p)$ such that
$\mbf z$ is supported over $\mbf n$ and 
$Q$ 
 has a solution.
\end{defn}

Now let $(\mbf h,\mbf n,\mbf t,\mbf p)$ be a positive $4$--partition. 
We shall denote by 
\begin{equation}\label{eqn_form}
(\dcup_{k\ge 1}H_{*,k}^\Lm,\cL,\mbf n,\mbf h)
\end{equation}
the set
of $n$--tuples $\mbf z\in (\dcup_{k\in K}H_{*,k}^\Lm,\cL)$ such that $\mbf z$ is supported over
$\mbf n$,  and $\sum_{i=1}^{n_j}|z_{\mu_j+i}|_\Lm\le h_j$, for $j=1,\ldots ,k$.
(Note that given the latter condition the $j$th exponential coordinate
$h_j^\prime$ of $Q$ satisfies
$h_j^\prime\le h_j$, for all $j$.) If the set $K$ has a unique element, $s$ say, and $X_s=I$ 
then we write  \[(H^\Lm,\cL,\mbf n,\mbf h)\]
instead of  ({\ref{eqn_form}) and to simplify notation we sometimes use the latter even though $K$ has more
than one element.

Given 
$z=\prod_{i=1}^{\eta}(h_i,f_i)\in F(H^\Lm)$ define
\begin{equation}\label{double-h}
\w(z)= \sum_{i=1}^{\eta}l(h_i)(l(h_i)-1).
\end{equation}
Define 
\begin{align}
|\,\mbf{s}| & = \textrm{max}\{l(s):s\in \mbf{{s}}_j, j\in K \},\label{mbfs}\\
W_0(\mbf z) &=  \sum_{i=1}^{n}\w(z_i),\textrm{ where }\w(z_i)\textrm{
  is defined in } (\ref{double-h}),\label{W0} \\
W_1(\mbf z) &=  \sum_{i=1}^{n}|\,z_i|,\textrm{ where }|\,z_i|\textrm{
  is defined in } (\ref{HL-length}),\label{W1}\\
W_2(\mbf z) &= \max\{l(h) : (h,f) \textrm{ occurs in } z_i, 
\textrm{ for some } i\}\textrm{ and}\label{W2}\\
W_3(\mbf z) &= \max \{|\,\mbf{{s}}| ,W_2^2(\mbf z)\}.\label{W3}
\end{align}
We use $W_i$ instead of $W_i(\mbf z)$ when appropriate.

Next define
\begin{equation}\label{M_j}
M(j)=\left\{
\begin{array}{ll}
W_2^2+1, & \mbox{if $j=0$}\\
W_3W_2^2(W_2^2+2|\,\mbf s|^2(|\,\mbf{{s}}|/2+1)^2), & \mbox{if $j=1$}\\
8|\,\mbf{{s}}|^{12}W_2^2+|\,\mbf{{s}}|, & \mbox{if $j=2$}
\end{array}
\right.
\end{equation}
and 
\begin{equation}\label{M_z}
M(\mbf z)=\max\{M(j):j=0,1,2\}.
\end{equation}

Let the $H^\Lm$--length of $z_i$ be $|\,z_i|=\eta_i$,  for $i=1,\ldots
,n$. Define
\begin{align}
\wt{\mu_1} & =0 \textrm{ and }\nonumber\\
\wt{\mu_{j+1}}& =\wt{\mu_j}+\eta_j, \textrm{ for } j=1,\ldots, n 
\end{align}
(so $\wt{\mu_{n+1}}=W_1$). Write
\[z_j=\prod_{i=1}^{\eta_j}(h_{i+\wt{\mu_j}},f_{i+\wt{\mu_j}}),\] for $j=1,\ldots
,n$. Consider the list \[(h_1,f_1),\ldots,(h_{W_1},f_{W_1})\] of
letters occuring in the above expressions for $z_1,\ldots, z_n$.
We call the $p$th element $(h_p,f_p)$ of this list the $p$th {\em  letter} 
of $\mbf z$, for $1\le p\le W_1$. 

Let $\lm_1,\ldots, \lm_{|\,W_1|}$ be distinct elements of $\Lm$. For
each $p$ such that $(h_p,f_p)$ is a proper exponential $H$--letter define
\[\ovr{(h_p,f_p)}=(h_p,\lm_p)\in H^\Lm.\]
For
each $p$ such that $(h_p,f_p)$ is a degenerate exponential $H$--letter define
\[\ovr{(h_p,f_p)}=(h_p\^f_p,f_p)\in H^\Lm.\]
For $j=1,\ldots n$, let 
\[\ovr{z_j}=\prod_{i=1}^{\eta_j}\ovr{(h_{i+\mu_j},f_{i+\mu_j})},\]
and define $\ovr{\mbf z}=(\ovr{z_1},\ldots,\ovr{z_n})$. 
Given $p$ such that $(h_p,f_p)$ is a proper exponential $H^\Lm$--letter
let $k_p$ be the unique integer such that $\cL$ implies $f_p\equiv k_p
(\mod l(h_p))$ and $0\le k_p<l(h_p)$. Define the {\em homogeneous}
system of parameters $\cH=\cH(\mbf z)$ {\em associated} to $\mbf z$
and $\cL$ to be the
parameter system:
\[\lm_p-1\ge 0,\]
\[\lm_p-k_p\equiv 0 (\mod l(h_p)),\]
for all $p$ with $1\le p\le W_1$ such that $(h_p,f_p)$ is a proper
exponential $H$--letter and $l(h_p)>1$.
Then $\ovr{\mbf z}$ is constrained by $\cH(\mbf z)$ and supported over $\mbf n$ and
\[\sum_{i=1}^n|\,z_i|_{\Lm}=\sum_{i=1}^n|\,\ovr{z_i}|_{\Lm},\]
so $\ovr{\mbf z}\in (\dcup_{k\in K}H_{*,k}^\Lm,\cH,\mbf n,\mbf h)$.
\begin{defn}\label{homog_eqn}
Let $\mbf z\in (\dcup_{k\in K}H_{*,k}^\Lm,\cL,\mbf n,\mbf h)$. The {\em homogeneous} equation
associated to $\mbf z$ is the quadratic exponential equation $Q_{\cH}=Q_{\cH}(\ovr{\mbf
  z},\cH, \mbf n,\mbf t,\mbf p)$ associated to $(\ovr{\mbf z},\cH,
\mbf n,\mbf t,\mbf p)$. 
\end{defn}
\begin{exx}\label{exx_z}
With 
$I$, 
$C_i$, $R_i$, $K$, $X_k$,  $\mbf s_k$, $H_{*,k}$, and $G_k$ as in 
Example \ref{exx_eqn}(\ref{exx_eqn_1}) let
\begin{align*}
z_1 & =(c_{22},\lm_1),\\
z_2 & =(c_{12}c_{22}c_{21},2)(c_{11}c_{22}c_{21},f_1)(c_{11}c_{22},2),
\textrm{ and }\\
 z_3 & =(c_{21}c_{22}c_{11}^2,f_2)(c_{22}c_{11}^2,2),
\end{align*}
where 
$f_1=\lm_1-2$ and 
$f_2 = \lm_1-\lm_2$. 
Let $\cL_0$ and $\cL_1$ be as in Example \ref{exx_eqn}(\ref{exx_eqn_1}). Then
$\mbf z=(z_1,z_2,z_3)$ is constrained by $\cL_i$, $i=1,2$. Let $n=3$, $t=p=0$ and
$\mbf n=(1,2)$. As supp$(z_1)=X_1$ and supp$(z_i)=X_2$, $i=1,2$, $\mbf z$ is supported 
over $\mbf n$ with basis $(1,2)$. Since $t=p=0$ we have $\mbf t=\mbf p=(0,0)$, 
$\mu_j=\xi_j$, for $j=0,1,2$ and 
$\mu_1=0,$ $\mu_2=1$ and $\mu_3=3$. The map $\bt$ from $\{d_1,d_2,d_3\}$ to $F(H^\Lm)$ satisfies
$\bt(d_i)=z_i$, for $i=1,2,3$. We have $\mbf q=(q_1,q_2)$, where $q_1=x_1^{-1}d_1x_1$ and 
$q_2=x_2^{-1}d_2x_2x_3^{-1}d_3x_3$.
Let $Q_i=Q(\mbf z,\cL_i,\mbf n,\mbf t,\mbf p)$, for $i=1,2$ and let $W_i$ be
as given in  Example \ref{exx_eqn}(\ref{exx_eqn_1}). A solution to $Q_i$ is also a solution to $W_i$.
A solution $(\phi,\al)$ to $W_i$ gives rise to solutions for $Q_i$: a solution is obtained by 
setting $\phi(x_1)$ equal to any element of $H_{*,1}$ and $\phi(x_2)=1\in H_{*,2}$. Since we 
found solutions for $W_1$ and $W_2$ in Example \ref{exx_eqn}(\ref{exx_eqn_1}) it follows that
$(\mbf n, \mbf t, \mbf p)\in \cL_i$--genus$(\mbf z)$, for $i=1,2$.

We have $\w(z_1)=0$, $\w(z_2)=6$ and $\w(z_3)=4$ and $|\,\mbf s|=4$ so  $W_0=10$. We have
$|\,z_1|=1$, $|\,z_2|=3$ and $|\,z_3|=2$ so  $W_1=6$. As $W_2=2$ we have $W_3=4$.
We have $M(0)=4$, $M(1)=576$ and $M(2)=131074$ so $M(\mbf z)=131074$.

We have $\wt{\mu_1}=0$, $\wt{\mu_2}=1$, $\wt{\mu_3}=4$ and $\wt{\mu_4}=6$ so
\begin{align*}
\ovr{z}_1 & =(c_{22},\lm_1),\\
\ovr{z}_2 & =(c_{12}c_{22}c_{21},2)(c_{11}c_{22}c_{21},\lm_3)(c_{11}c_{22},2),
\textrm{ and }\\
\ovr{z}_3 & =(c_{21}c_{22}c_{11}^2,\lm_6)(c_{22}c_{11}^2,2),
\end{align*} 
Then $\cH$ consists of the inequalities $\lm_p-1\ge 0$, for $p=1,3,6$, and the congruences
$\lm_p\equiv 0 (\mod 2)$, for $p=3,6$. The homogeneous equation
associated to $\mbf z$ is $Q_{\cH}=(\mbf q=1,\bt_{\cH},\cH)$, where $\bt_{\cH}(d_i)=\ovr{z}_i$.
Define $\phi$ from $L(Q_{\cH})$ to $H_{*,1}\dcup H_{*,2}$ by $\phi(x_i)=1$, $i=1,2$ and 
$\phi(x_3)=c_{11}c_{22}$. Define $\al$ to be the retraction such that $\al(\lm_1)=6$, $\al(\lm_3)=4$ and 
$\al(\lm_6)=2$. Then $(\phi,\al)$ is a solution to $Q_{\cH}$. Compare this solution to the solution
$(\phi_1,\al_1)$ to $W_i$ given in Example \ref{exx_eqn}(\ref{exx_eqn_1}). We have obtained $\phi$ from
$\phi_1$ and set $\al(\lm_p)$ equal to the value $\al_1(g_p)$, where $g_p$ is the polynomial occuring in the 
$p$th letter of $\mbf z$. The same construction fails using $\al_0$ from Example \ref{exx_eqn}(\ref{exx_eqn_1})
because we obtain $\al(\lm_6)=0$.
\end{exx}
\begin{exx}\label{exx_z_probs}
Let
$I$, 
$C_i$, $R_i$, $K$, $X_k$,  $\mbf s_k$, $H_{*,k}$, and $G_k$ as in 
Example \ref{exx_eqn_probs}. In this case equations have environment $(X_1;(H_i)_{i\in I};G_1)^\Lm$.
\be
\item\label{exx_z_probs_1} Let $\mbf n=\mbf h=(1)$. The set $(H_{*,1}^\Lm,\cL,\mbf n,\mbf h)$ consists of elements
of the form $\mbf z=(u_0(h,f)u_1)$, where $u_i$ is a product of degenerate exponential letters, that is
$u_i\in F((H_{*,1},\Lm)_0)$, for $i=1,2$.
The equation associated to such $\mbf z$ is $Q=((x_1^{-1}d_1x_1)=1,\bt,\cL)$, where $\bt(d_1)=u_0(h,f)u_1$.
This is equivalent to the equation $Q^\prime=((d_1)=1,\bt,\cL)$, in that solutions to $Q$ are solutions to 
$Q^\prime$ and, given a solution to $Q^\prime$ any value of $\phi(x_1)$ gives a solution to $Q$. Suppose that
$\mbf z=((h,f))$, 
$f=\lm_1$ and  $\cL$ contains the subsystem $\cL_0$ consisting of $\lm_1\neq 0$ and $\lm_1\equiv 0 (\mod l(h))$. 
Then $Q^\prime$ is the equation
of Example \ref{exx_eqn_probs}(\ref{exx_eqn_probs_1}). 
Hence
$(\mbf n, (0), (0))\in \cL_0$--genus$(\mbf z)$
if and only if $h$ has finite order in $G_1$.

On the other hand if $\mbf z=((h_1,\lm_1)(h_2,l(h_2))$ then, as in Example \ref{exx_eqn_probs}(\ref{exx_eqn_probs_2}),
$(\mbf n, (0), (0))\in \cL_0$--genus$(\mbf z)$
if and 
only if $h_1^m=h_2$, for some non--zero integer $m$.
\item\label{exx_z_probs_2}  Let $\mbf n=(1)$ and $\mbf h=(2)$. 
The set $(H_{*,1}^\Lm,\cL,\mbf n,\mbf h)$ consists of elements
of the form $\mbf z=(z_1)$, where $z_1\in H_{*,1}^\Lm$, $z_1=u_0w_1u_1w_2u_2$, with $u_i,\in  F((H_{*,1},\Lm)_0)$
and 
$w_i\in H^\Lm_{*,1}$.
The equation associated to such $\mbf z$ is $Q=((x_1^{-1}d_1x_1)=1,\bt,\cL)$, 
where $\bt(d_1)=z_1$.
This is equivalent to the equation $Q^\prime=((d_1)=1,\bt,\cL)$, as in the previous case.
With $z_1=(h_1,\lm_1)(h_2,\lm_2)$ and $\cL_0$ consisting of the inequalities $\lm_i\neq 0$ and the 
congruences $\lm_i\equiv 0 (\mod l(h_i))$ 
it follows that
$(\mbf n, (0), (0))\in \cL_0$--genus$(\mbf z)$
if and 
only if a non--zero power of $h_1$ is equal to a non--zero power of $h_2$.
\item\label{exx_z_probs_3} Let $\mbf n=\mbf h=(2)$. The set $(H_{*,1}^\Lm,\cL,\mbf n,\mbf h)$ consists of elements
of the form $\mbf z=(z_1,z_2)$, where $z_i\in H_{*,1}^\Lm$, $z_i=u_iw_iv_i$, with $u_i,v_i \in  F((H_{*,1},\Lm)_0)$
and either 
$w_i=(h_{1},f_1)(h_2,f_2)$, with $f_1,f_2\in M\backslash \ZZ$ and $w_{i-1}\in F((H_{*,1},\Lm)_0)$ or
$w_1=(h_{1},f_1)$, $w_2=(h_2,f_2)$, with $f_i\in M$ (and $w_i$ may be trivial). 
The equation associated to such $\mbf z$ is $Q=((x_1^{-1}d_1x_1x_2^{-1}d_2x_2)=1,\bt,\cL)$, 
where $\bt(d_i)=z_i$.
This is equivalent to the equation $Q^\prime=((d_1x_2^{-1}d_2x_2)=1,\bt,\cL)$, in the sense of the previous 
case of the Example. With $\mbf z=((h_1,\lm_1),(h_2,\lm_2))$ and $\cL_0$ as in 
Example \ref{exx_eqn_probs}(\ref{exx_eqn_probs_3}) it follows that
$(\mbf n, (0), (0))\in \cL_0$--genus$(\mbf z)$
if and 
only if a non--zero power of $h_1$ is conjugate to a non--zero power of $h_2$.
\ee
\end{exx}

\section{Resolutions of quadratic equations}\label{resolution}
Throughout this section the sets $I$, $K$ and $X_k$ and the groups $H_i$, $H_{*,k}$ and $G_k$
are defined as in Section \ref{qpequations}, (\ref{Hsys}) and (\ref{Gsys}). In particular we
assume that $H_i$ has solvable word problem, for all $i\in I$.

Given a $k$--tuple $(a_1,\ldots ,a_k)$ of integers we say that the $s$--tuple
$(b_{i_1},\ldots ,b_{i_s})$ of integers is a {\em sub--basis} of $(a_1,\ldots ,a_k)$
if $(b_1,\ldots , b_k)$ is  a permutation of $(a_1,\ldots ,a_k)$ and $1\le i_1<i_2< \cdots<i_s\le k$. 
\begin{defn}
Let $W=(\mbf w=1,\bt,\cL)$ be a constrained quadratic exponential equation 
with environment
$((X_k)_{k\in K};(H_i)_{i\in I} ;(G_k)_{k\in K})^\Lm$ and
basis $(a_1,\ldots ,a_k)$. 
Let $R$  be a set $R=\{W_i=(\mbf{w_i}=1,\g_i,{\cL}_i)\,:\, i=1,\ldots ,n\}$ 
of constrained quadratic exponential equations such that $W_i$ has environment
$((X_k)_{k\in K};(H_i)_{i\in I} ;(G_k)_{k\in K})^\Lm$ and basis equal to a sub--basis of $(a_1,\ldots ,a_k)$,
for $i=1,\ldots ,n$.
Then $R$ 
is called a {\em resolution} of $W$ if the following
conditions hold.
\be
\item If $W$ has a   solution then $W_i$
has a    solution, for some $i$.
\item There is an algorithm which, given a
solution to $W_i$, for some $i$, outputs a
 solution to $W$ in $(G_{k})_{k\in K}$.
\ee
\end{defn}
A resolution $R$ of $W$ is said to have property $\mathcal \mbf P$ if every element of $R$ has
property $\mathcal \mbf P$.

\begin{theorem}[cf. {\cite[Proposition 2.2]{\GL}}]\label{stdform}
Let $W=(\mbf w=1,\bt,\cL)$ be a  constrained quadratic exponential
equation.
Then there exists an algorithm which, given input $W$, outputs a resolution
of $W$ in standard form.
\end{theorem}

{\em Proof.}  Let $\mbf w=(w_1,\ldots ,w_k)$ and let $(b^i_1,\ldots ,b^i_{n_i})$ be a boundary labels list for $w_i$,
$i=1,\ldots ,k$. Choose $\eta\in$Aut$_f(F)$ using Proposition \ref{stdform0}, such
that $\eta$ is an admissible transformation of $(\mbf w,\mbf b)$
to $(\mbf v,\mbf d)$, where $\mbf v=\mbf q(\mbf n,\mbf t,\mbf p)$, with $\mbf n$, $\mbf t$
and $\mbf p$ as described in Proposition \ref{stdform0}, and
$\eta(F(D))=F(D)$.

Let $T=L(\mbf w)$ and define $U=\cup_{a\in T}L(\eta(a))$ so $L(\mbf v)\subset U$. Since
$\eta$ fixes all but finitely many elements of $D\dcup X$ there is a
finite set $V\subset D\dcup X$ such that $\eta(F(V))=F(V)$ and $T\cup
U\subset V$. Then $D\cap T\subset D\cap V$ so we may define a map
$\bt^*:F(D\cap V)\maps F(H^\Lm)$ with $\bt^*(d)=\bt(d)$, for all $d\in
D\cap T$ and  $\bt^*(d)=1$, for all $d\notin
D\cap T$.

As $\eta(F(D))=F(D)$ we have $\eta(F(D\cap V))=F(D\cap V)$, so we may
define a map $\de^*:F(D\cap V)\maps F(H^\Lm)$ by
$\de^*=\bt^*\circ\eta^{-1}|_{F(D\cap V)}$. Define $\de:F(L_D(\mbf v))\maps
F(H^\Lm)$ by $\de=\de^*|_{F(L_D(\mbf v))}$. If $d\in L_D(\mbf v)$ then
$d=\wt{\eta(b^i_j)}$, for appropriate $i$ and $j$, so
$\de(d)=\de^*(d)=\bt^*\eta^{-1}(d)=\bt^*(b^i_j)=\bt(b^i_j)$.
It follows that 
$(\mbf v=1,\de,\cL)$ is a constrained quadratic exponential equation with the same
environment 
and basis 
as $W$.

Suppose $(\mu,\al)$ is a solution to $(\mbf
v=1,\de,\cL)$. 
Let $T_i=L_X(w_i)$ and $U_i=\cup_{x\in T_i} L(\eta(x))$, for $i=1,\ldots ,k$.
Let $S_i=L_X(v_i)$ and $V_i=\cup_{x\in S_i} L(\eta^{-1}(x))$, for $i=1,\ldots ,k$.
Define a map $\mu^*:V\maps H$  by
\[\mu^*(a)=\left\{\begin{array}{l}
\mu(a),\textrm{ if } a\in L(\mbf v)\\
\hat\al\de^*(a),\textrm{ if } a\in D\cap V\backslash L(\mbf v)\\
\textrm{any element of }H_{*,a_i} \textrm{ if } a\in U_i\backslash L(\mbf v)\\
\textrm{any element of }H \textrm{ otherwise}
\end{array}\right. ,\]
and extend $\mu^*$ to a homomorphism from $F(V)$ to $H$.
(Thus if $a\in D\cap V$ then $\mu^*(a)=\hat\al\de^*(a)$.) If $d\in
D\cap T$ then $\eta(d)\in F(D\cap U)$ so
$\mu^*(\eta(d))=\hat\al\de^*(\eta(d))= \hat\al\bt^*(d)=\hat\al\bt(d)$.
Also $\mu^*\eta(a)\in H_{*,a_i}$, for all $a\in L(w_i)$, and  $\mu^*\eta(w_i)$ 
is conjugate in $H_{*,a_i}$ to
$\mu^*\wt{\eta(w_i)}=\mu\wt{\eta(w_i)}=\mu(v_i)\in N_{a_i}$, for
$i=1,\ldots ,k$. It follows that, setting $\phi
=\mu^*\circ\eta|_{F(T)}$, we have  a solution $(\phi,\al)$ to $W$.

Conversely suppose $(\phi,\al)$ is a
solution to $W$. Define a map $\phi^*:V\maps H$ by
\[\phi^*(a)=\left\{\begin{array}{l}
\phi(a),\textrm{ if } a\in T\\
\hat\al\bt^*(a),\textrm{ if } a\in D\cap V\backslash T\\
\textrm{any element of }H_{*,a_i} \textrm{ if } a\in V_i\backslash T\\
\textrm{any element of }H \textrm{ otherwise}
\end{array}\right. ,\]
and extend $\phi^*$ to a homomorphism $\phi^*:F(V)\maps H$. Define
$\mu=\phi^*\circ\eta^{-1}|_{L(\mbf v)}$. Then
$\mu(a)\in H_{*,a_i}$, for all $a\in L(v_i)$, and 
$\mu(v_i)=\phi^*\eta^{-1}(v_i)$ is conjugate to $\phi^*(w_i)\in
N_{a_i}$, for $i=1,\ldots ,k$. If  $d\in L_D(\mbf v)$, we have
$\mu(d)=\phi^*\eta^{-1}(d)=\hat\al\bt^*\eta^{-1}(d)=\hat\al\de^*(d)=\hat\al\de(d)$.
Hence $(\mu,\al)$ is a
solution to
$(\mbf v=1,\de,\cL)$. It follows that $R=\{(\mbf v=1,\de,\cL)\}$ is a resolution of
$W$ in standard form. Since $\eta, U_i, V_i, \bt^*, \de^*, L(\mbf w)$ and $L(\mbf v)$ may all
be effectively constructed there exists an algorithm as claimed.\\

The process of finding  $\eta, \mbf v, U, V, \bt^*, \de^*, L(\mbf w)$ and $L(\mbf v)$ is called {\em
reduction to standard form}.

Let $u$ be an element of $F(H^\Lm)$ such that $u$ is constrained by a 
parameter system $\cL$. We say $(u,{\cL})$ has $H^\Lm(i)${\em --redundancy} if
condition $i$ below holds and that $(u,{\cL})$ has $H^\Lm${\em
--redundancy} if $(u,{\cL})$ has $H^\Lm(i)$--redundancy for some $i$, with
$1\le i\le 7$.
\renewcommand{\theenumi}{\arabic{enumi}}
\renewcommand{\labelenumi}{$i=$\theenumi:}
\be
\item\label{h=1} A generator $(h,f)\in H^\Lm$ occurs in $u$ with $h=1\in H$.
\item\label{f=0} A generator $(h,f)\in H^\Lm$ occurs in $u$ and $\cL$ implies that
$f=0$.
\item\label{a=pp} A proper exponential $H$--letter $(a,f)$ occurs in $u$ with $a=a_0^s$, where
$a_0$ is a cyclically reduced element of $H$, $s\in \ZZ$, $|\,s|>1$ and
$l(a)=m>1$.
\item\label{negexp} A generator $(h,f)\in H^\Lm$ occurs in $u$ with exponent $-1$.
\item\label{constpoly} The word $u$ contains a subword $(a_1,f_1)\cdots (a_n,f_n)$, where
$(a_i,f_i)\in H^\Lm$ and either 
\be
\item 
$n>1$ and $f_1,\ldots ,f_n\in \ZZ$ or 
\item
there is at least one $i$ such that
$(a_i,f_i)$ is a proper exponential $H$--letter and, for each $i=1,\ldots ,n$, 
if
$(a_i,f_i)$ is a proper exponential $H$--letter then $\cL$ implies that $f_i=m_i\in\ZZ$.
\ee
\item\label{sameword} The word $u$ contains a subword $(a,f)(b,g)$, where
$(a,f),(b,g)\in H^\Lm$ at least one of which is a proper exponential $H$--letter,
$l(a)=l(b)=m$,
$\cL$ implies that $f\equiv k(\mod m)$, for some $k$ with $0\le k<m$,
$a=g_1\cdots g_m$,
$g_i\in H_{k_i}$ and $b^\e=g_{k+1}\cdots g_mg_1\cdots g_k$, with  $\e=\pm 1$.
\item\label{cancelword} The word $u$ contains a subword $(a,f)(b,g)$, where
$(a,f),(b,g)\in H^\Lm$ at least one of which is a proper exponential $H$--letter,
$l(a)=m$, $l(b)=n$, $a=a_1\cdots
a_m$, $b=b_1\cdots b_n$, $a_i\in H_{k_i}$, $b_j\in H_{l_j}$, $\cL$
implies that $f\equiv k(\mod m)$, for some $k$ with $0\le k<m$, and
$a_k=b_1^{-1}\in H$.
\ee
\renewcommand{\labelenumi}{(\theenumi)}
\renewcommand{\theenumi}{\roman{enumi}}
If $(u,{\cL})$ has no $H^\Lm$--redundancy it is said to be $H^\Lm${\em
--irredundant}. The equation $(\mbf w=1,\bt,{\cL})$ is said to be
$H^\Lm${\em --redundant} if $(\bt(d),{\cL})$ is $H^\Lm$--redundant, for some
$d\in L_D(\mbf w)$, and $H^\Lm${\em --irredundant} otherwise.

If $s\in\ZZ$, with $1\le s\le 7$, and $(u,{\cL})$ is $H^\Lm(i)$--irredundant for
$i=1,\ldots ,s$ we say
$(u,{\cL})$ is $H^\Lm(1,s)${\em --irredundant}. The equation $(\mbf
w=1,\bt,{\cL})$ is $H^\Lm(1,s)${\em --irredundant} if $(\bt(d),{\cL})$ is
$H^\Lm(1,s)$--irredundant, for all $d\in L_D(\mbf w)$.

Given a  constrained quadratic exponential equation $W=(\mbf w=1,\bt,\cL)$ and $d\in L_D(\mbf
w)$, suppose $(\bt(d),\cL)$ has $H^\Lm(i)$--redundancy. For each $i$
we define a set
$R(W,i)$ of constrained quadratic exponential equations. When $i\neq 4$ the set $R(W,i)=\{W^\prime\}$, where
$W^\prime=(\mbf w=1,\bt^\prime,\cL)$ with $\bt^\prime$ obtained from
$\bt$ by setting $\bt^\prime(d^\prime)=\bt(d^\prime)$, if
$d^\prime\neq d$ and $\bt^\prime(d)$ equal to an element of
$F(H^\Lm)$ dependent on $i$ as follows.
\be
\item[$i=$\ref{h=1} or \ref{f=0}:] Delete, from $\bt(d)$, all occurrences of the generator
$(h,f)$ at which redundancy occurs and freely reduce the resulting word, to give $\bt^\prime(d)$.
\item[$i=$\ref{a=pp}:] If $H^\Lm(\ref{a=pp})$--redundancy occurs at a
subword $(a,f)$ of $\bt(d)$, where $a=a_0^s$, as in the description of
$H^\Lm(\ref{a=pp})$--redundancy above, then replace all occurrences of $(a,f)$ in
$\bt(d)$ by $(a_0,f)$, to give $\bt^\prime(d)$. Note
that, if $\cL$ implies that $f\equiv k(\mod l(a))$ then
$\cL$ implies   that $f\equiv k(\mod l(a_0))$.
\item[$i=$\ref{constpoly}:] Suppose $H^\Lm(\ref{constpoly})$--redundancy occurs at a
subword 
$(a_1,f_1)\cdots (a_n,f_n)$
of $\bt(d)$. For $i$ such that $(a_i,f_i)$ is a proper exponential
$H$--letter set $m_i\in\ZZ$ such that $\cL$ implies $f_i=m_i$. For all
other $i$ set $m_i=f_i$. Set $c_i=a_i\^m_i$, for $i=1,\ldots ,n$, to give $\bt^\prime(d)$. 
Replace $(a_1,f_1)\cdots (a_n,f_n)$ in
$\bt(d)$ by $(c,h)$, where $c=c_1\cdots c_n$ and $h=l(c)\in \ZZ$.
\item[$i=$\ref{sameword}:] If $H^\Lm(\ref{sameword})$--redundancy occurs at a
subword $(a,f)(b,g)$ of $\bt(d)$ then replace all occurrences of $(a,f)(b,g)$ in
$\bt(d)$ by $(a,f+\e g)$, to give $\bt^\prime(d)$. Note that, if $\cL$ implies that $f\equiv k(\mod
m)$ and $g\equiv l(\mod m)$, then $\cL$ implies that $f+\e g\equiv
k+\e l(\mod m)$.
\item[$i=$\ref{cancelword}:] Suppose $H^\Lm(\ref{cancelword})$--redundancy occurs at a
subword $(a,f)(b,g)$ of $\bt(d)$, as described in the case
$i=$\ref{cancelword} above. We obtain $\bt^\prime(d)$ from $\bt(d)$ as follows.
If $(a,f)$ and $(b,g)$ are both proper exponential
$H$--letters then replace all occurrences of $(a,f)(b,g)$ in
$\bt(d)$ by $(a,f-1)(b_2\cdots b_nb_1,g-1)$.  
If $(a,f)$ is a
proper exponential $H$--letter but $(b,g)$ is not then replace all occurrences of
$(a,f)(b,g)$ in $\bt(d)$ by $(a,f-1)(\t(b,g-1),g-1)$.
If $(b,g)$ is a
proper exponential $H$--letter but $(a,f)$ is not then replace all occurrences of
$(a,f)(b,g)$ in $\bt(d)$ by $(\io(a,f-1),f-1)(b_2\cdots b_nb_1,g-1)$. 
\ee 
We consider finally the case
where $W$ has occurrence of $H^\Lm(\ref{negexp})$--redundancy.
\be
\item[$i=$\ref{negexp}:] Suppose $\bt(d)=w_0(h,f)^{-1}w_1$, where
$w_i\in F(H^\Lm)$, $(h,f)\in H^\Lm$ and $\bt(d)$ is reduced as written.
If $(h,f)$ is a degenerate exponential $H$--letter then $f=l(h)$ and we 
define $\bt_1:F(L_D(\mbf w))\maps H$ by
$\bt_1(d)=w_0(h^{-1},f)w_1$
and $\bt_1(d^\prime)=\bt(d^\prime)$, for $d^\prime\neq d$.
Then $R(W,\ref{negexp})=\{(\mbf w=1,\bt_1,\cL)\}$. 

On the 
other hand,
if $(h,f)$ is a proper exponential $H$--letter,
let $m=l(h)$ and let $k\in \ZZ$, such that $0\le k <m$ and $\cL$ implies
that $f\equiv k(\mod m)$. 
For $i=1,2$ define $\bt_i:F(L_D(\mbf w))\maps F(H^\Lm)$ by $\bt_i(v)=\bt(v)$,
if $v\neq d$,
\[\bt_1(d)=w_0([\io(h,k)^{-1}\t(h,m-k)^{-1}],f)w_1\] and
\[\bt_2(d)=w_0([\t(h,k)\io(h,m-k)],-f)w_1.\]
Let $\cL_1=\cL\cup\{f\ge 0\}$, let $\cL_2=\cL\cup\{-f>0\}$.
Then any solution $\al$ to $\cL$ is a solution to either $\cL_1$ or $\cL_2$ and, if
$\al$ is a solution to $\cL_i$, we have
\[\hat\al\bt(d^\prime)=\hat\al\bt_i(d^\prime),\textrm{ for all }d^\prime\in
L_D(\mbf w), i=1,2.\]
Define $W_i=(\mbf w=1,\bt_i,\cL_i)$, for $i=1,2$ and set
$R(W,\ref{negexp})=\{W_1,W_2\}$.
\ee

\begin{lemma}\label{Red} Let $W=(\mbf w=1,\bt,\cL)$ be a
constrained quadratic exponential equation. Then there exists an algorithm which, given $W$,
outputs an $H^\Lm(1,\ref{negexp})$--irredundant resolution $R$
of $W$ in standard form.
\end{lemma}

{\em Proof.}  We may assume, using Theorem \ref{stdform}, that $W$ is in
standard form. The solvability of the word problem for $H$ allows us to
determine  whether $W$ contains an occurrence of $H^\Lm(i)$--redundancy,
for some fixed $i$, $1\le i\le 3$, and to find an occurrence if one exists. Replacing $W$
by the unique equation of the set $R(W,i)$ we obtain an equation $W^\prime$ with fewer
occurrences of $H^\Lm(i)$--redundancy than $W$. Repeating this process we
eventually obtain an equation $W^\prime$ which is
$H^\Lm(i)$--irredundant, for $i=1,2$ and $3$. We may now eliminate an
occurrence of $H^\Lm(\ref{negexp})$--redundancy by replacing $\{W^\prime\}$ with
$R(W^\prime,\ref{negexp})$,
if necessary. Note that, by construction,
no element of $R(W^\prime,\ref{negexp})$ contains an occurrence of
$H^\Lm(i)$--redundancy, with $i=1,2$ or $3$. We may repeat the removal
of occurrences of $H^\Lm(\ref{negexp})$--redundancy on
elements of $R(W^\prime,\ref{negexp})$ 
to
obtain a set $R$ of
$H^\Lm(1,\ref{negexp})$--irredundant,  constrained, quadratic exponential equations.

It remains to show that if $W$ is a  constrained quadratic exponential equation then
$R(W,i)$ is a resolution of $W$, for $i=1,2,3,4$. This is easy to see if
$i=1$ or $2$.

Suppose we have $H^\Lm(\ref{a=pp})$--redundancy at the occurrence of
$(a,f)$ in $\bt(d)=w_0(a,f)^\e w_1$, where $\e=\pm 1$ and $d\in L_D(\mbf
w)$. As $W$ is constrained, $\cL$ implies that $f\equiv k(\mod m)$,
for some integer $k$, with $0\le k< m$. Let $a_0$ be a cyclically
reduced word and $s>1$ an integer such that $a=a_0^s$. By definition of
$H^\Lm(\ref{a=pp})$--redundancy $l(a_0)>0$. Set $n=l(a_0)$ and let $p,l\in \ZZ$
such that $k=pn+l$, with $0\le l<n$. Assume $\al$ is a solution to
$\cL$. Then, if $\al(f)\ge 0$, say $\al(f)=qsn+k=qm+k$, we have
$\hat\al(a,f)=a^{q}\io(a,k) = a_0^{sq}a_0^{p}\io(a_0,l) =
a_0\^(qsn+pn+l) = \hat\al(a_0,f)$. On the other hand, if $\al(f)<0$, say
$-\al(f)= qm+k$, then $\hat\al(a,f)=a^{-q}\io(a^{-1},k) = a_0^{-sq}a_0^{-p}\io(a_0^{-1},l) =
a_0\^-(qsn+pn+l) = \hat\al(a_0,f)$. Hence, if $\al$ is a solution
to $\cL$ then $\hat\al\bt(d)=\hat\al\bt^\prime(d)$,
for all $d\in L_D(\mbf w)$.
Therefore solutions of $W$ and $W^\prime$ coincide and $R(W,3)$ is a resolution of $W$.

When $H^\Lm(\ref{negexp})$--redundancy occurs at a exponential $H$--letter which is
not proper, then the resulting set $R(W,\ref{negexp})$ is clearly a resolution
of $W$. In the remaining case a solution to $W$ is also a solution to $W_i$, for
$i=1$ or $2$. Conversely given a solution to $W_i$ it is a solution to
$W$. It follows that $R(W,i)$ is a resolution of $W$, for $i=1,2,3$ and
$4$. 

As $\mbf w$ is the quadratic word appearing in every element of each
 of the resolutions used above to replace $W$, the 
resolution obtained is also in standard form.\\

The application of the algorithm of Lemma \ref{Red}  to produce an $H^\Lm(1,\ref{negexp})$--irredundant resolution
of an equation is called $H^\Lm(1,\ref{negexp})${\em --reduction}.

The
{\em dual} of a exponential $H$--letter $(a,f)$
is  \[\supstar{(a,f)} =(a^{-1},-f).\]
Note that $\supstar{(\supstar{(a,f)})}=(a,f)$. There is
a well--defined map $\nu:H^\Lm\maps F(H^\Lm)$ given by $\nu(a,f)
=\supstar{(a,f)}$ 
and
$\nu(b,g)=(b,g)$, if $(b,g)\neq (a,f)$.
Such
a map induces an endomorphism $\nu$ of  $F(H^\Lm)$, called the
{\em dualization} at $(a,f)$.

We define the {\em deletion} at $(a,f)\in H^\Lm$
to be the endomorphism of $F(H^\Lm)$ defined by
$\de(a,f)=\de(\supstar{(a,f)}) = 1$ and $\de(b,g)=(b,g)$, for all
$(b,g)\not\in\{(a,f),\supstar{(a,f)}\}$.

Let $W=(\mbf w=1,\bt,\cL)$ be a  constrained quadratic exponential equation.
Dualization at $(a,f)$ is {\em
admissible} for $W$ if $(a,f)$ occurs in $\bt(d)$, for some $d\in
L_D(\mbf w)$ and if $\cL\cup\{-f>0\}$ is consistent. Deletion
at  $(a,f)$ is {\em admissible} for $W$ if $(a,f)$ occurs in $\bt(d)$, for some $d\in
L_D(\mbf w)$ and if $\cL\cup\{f=0\}$ is consistent.

Define $\cM(\mbf
w,\bt)=\{f\in M\,:\,(a,f)$ occurs in $(\mbf w,\bt)\}$.
Then $W$ is called {\em positive} if it is $H^\Lm(1,\ref{negexp})$--irredundant
and $\cL$ implies that $f>0$, for all $f\in \cM(\mbf w,\bt)$.

\begin{lemma}\label{Posi}
Let $W=(\mbf w =1,\bt,\cL)$ be an $H^\Lm(1,\ref{negexp})$--irredundant, consistent,  quadratic exponential
equation. Then there exists an
algorithm which, when given $W$, outputs
a positive resolution of $W$ in standard form.
\end{lemma}

\noindent{\em Proof.}
If $(a,f)$ is a degenerate exponential $H$--letter occuring in $W$ and $f<0$ then we may replace
$(a,f)$ by $(a^{-1},-f)$. We may therefore restrict attention to proper exponential $H$--letters.
Given $f\in \cM(\mbf w,\bt)$  define
$H^\Lm(\mbf w,\bt,f)=\{(a,f)\in H^\Lm\,:\,(a,f)$ occurs in $(\mbf w,\bt)\}$.
Suppose  $f\in \cM(\mbf w,\bt)$  and assume that $H^\Lm(\mbf w,\bt,f)=\{(a_1,f),\ldots ,(a_n,f)\}$ and that
$H^\Lm(\mbf w,\bt,-f)=\{(a_{n+1},-f),\ldots ,(a_{n+p},-f)\}$.
%
If
$\cL\cup\{f<0\}$ is consistent,
let $\nu_i$ be the dualization at $(a_i,f)$, $i=1,\ldots ,n+p$, 
define $\nu^f=\nu_n\circ \cdots \circ \nu_1$ and set
$W_-(f)=(\mbf w=1,\nu^f\circ \bt,\cL\cup\{-f>0\})$.
 If
$\cL\cup\{f=0\}$ is consistent,
let $\de_i$ be the deletion at $(a_i,f)$.
Define $\de^f=\de_{n+p}\circ\cdots \circ\de_n\circ \cdots \circ \de_1$ and set
$W_0(f)=(\mbf w=1,\de^f\circ \bt,\cL\cup\{f=0\})$.
Finally, if $\cL\cup \{f>0\}$
is consistent, define $\nu^{-f}=\nu_{n+p}\circ \cdots \circ \nu_{n+1}$ and set $W_+(f)=(\mbf w=1,\bt,\cL\cup\{f>0\})$.
Let $R_0$ be the set consisting of whichever of $W_-(f)$, $W_0(f)$ and $W_+(f)$ are defined.
As $\cL$
is consistent $R_0$ is non--empty.

Note that for any solution $\al$ to $\cL$ we have
\[\begin{array}{lclcl}
\hat\al(a,f) &=& a\^\al(f) &=& a^{-1}\^ -\al(f)\\
            &=& a^{-1}\^\al(-f) &=& \hat\al(a^{-1},-f)\\
            &=& \hat\al(\supstar{(a,f)}).
\end{array}\]
Hence $(\phi,\al)$ is a solution to $W$ with $\al(f)<0$ if and only if $(\phi,\al)$ is a
solution to $W_-(f)$. Furthermore $(\phi,\al)$ is a solution to $W$ with
$\al(f)>0$ or $\al(f)=0$ if and only if $(\phi,\al)$ is a
solution to $W_+(f)$ or $W_0(f)$, respectively. Thus $R_0$ is a resolution of $W$ and as $\cL$ is
linear there is an algorithm which constructs $R_0$ from input $W$.

Let $\cN (W)=\{f\in \cM(\mbf w,\bt)\,:\, \cL\cup\{f\le 0\}$ is consistent $\}$ and
suppose $f\in \cN(W)$.
By definition, $f\notin \cN(W_+(f))$. Also $f\notin \cN(W_0(f))\cup \cN(W_-(f))$ as
$H^\Lm(\mbf w,\de^f \circ \bt,f)=H^\Lm(\mbf w,\nu^f \circ \bt,f)=\nul$.
Similarly $-f\notin \cN(W_0(f))\cup \cN(W_-(f))\cup \cN(W_+(f))$.
It follows that, for $W^\prime\in R_0$,
we have $|\,\cN(W^\prime)|<|\,\cN(W)|$.

No $H^\Lm(i)$--redundancy, with $i=1,3$ or $4$, is introduced in
the process of forming the elements of $R_0$. If
$H^\Lm(\ref{f=0})$--redundancy occurs in some element $W^\prime$ of
$R_0$, then its subsequent removal does not increase
$|\,\cN(W^\prime)|$. Thus we may assume that $R_0$ is
$H^\Lm(1,\ref{negexp})$--irredundant and repeat the above process on
elements of $R_0$. Eventually a positive resolution of $W$ is
obtained.  As in the proof of Lemma \ref{Red} we may assume that $W$ 
and this resolution of $W$ are both in standard form.\\

The process of forming a positive resolution is called {\em positive} $H^\Lm(1,4)${\em --reduction}.

\begin{lemma}\label{Irr}
Let $W=(\mbf w=1,\bt,\cL)$ be a  constrained, quadratic exponential equation.
Then there exists an algorithm which, when given $W$, outputs
an  $H^\Lm$--irredundant, positive resolution of $W$ in standard form.
\end{lemma}

\noindent{\em Proof.} Lemma \ref{Posi} implies that we may assume
that $W$ is positive. If $W$ contains an occurrence of
$H^\Lm(\ref{constpoly})$--redundancy then the unique equation
$W^\prime$ of $R(W,\ref{constpoly})$, described above, is
positive. Clearly $R(W,\ref{constpoly})$ is a resolution of $W$
and $W^\prime$ has fewer occurrences of
$H^\Lm(\ref{constpoly})$--redundancy than $W$. Furthermore the sum of exponential length and $H^\Lm$--length 
decreases in passing from $W$ to $W^\prime$.
Then $W^\prime$ may be $H^\Lm(1,\ref{negexp})$--reduced, if
necessary, to produce a resolution of $W$ every element of which
is  
positive and has smaller exponential length or $H^\Lm$--length than $W$. 
We may continue this process to arrive at a resolution
which is positive and $H^\Lm(1,\ref{constpoly})$--irredundant. We call
the process of forming this resolution {\em positive} 
$H^\Lm(1,\ref{constpoly})${\em --reduction}.

Note that the solvability of the word problem in $H$ allows the
identification of occurrences of $H^\Lm(\ref{sameword})$-- and
$H^\Lm(\ref{cancelword})$--redundancy. Assume that $W$ is an
element of the positive $H^\Lm(1,\ref{constpoly})$--irredundant 
resolution obtained in the previous paragraph.
Given an occurrence of $H^\Lm(\ref{sameword})$--redundancy we note
that, as $W$ is positive, the set $R(W,\ref{sameword})$ defined
above is a resolution of $W$. Also, the unique element of
$R(W,\ref{sameword})$ has smaller $H^\Lm$--length than $W$.
Hence we may form the resolution $R(W,\ref{sameword})$,
to each element of which we apply positive $H^\Lm(1,\ref{constpoly})$--reduction, %
if necessary, and repeat the process.
As none of the operations involved increases $H^\Lm$--length, this process terminates and  the
resulting resolution is positive and $H^\Lm(1,\ref{sameword})$--irredundant.
We call the process of forming this resolution
{\em positive} $H^\Lm(1,\ref{sameword})${\em --reduction}.

Given $W$ in the above resolution, suppose $W$ contains an
occurrence of $H^\Lm(\ref{cancelword})$--redundancy at
$(a,f)(b,g)$, as described in the definition of
$H^\Lm$--redundancy. The set $R(W,\ref{cancelword})$ is a
resolution of $W$, as $W$ is positive. Let $W^\prime$ be the
unique element of this set. As above, we may form a positive
$H^\Lm(1,\ref{sameword})$--irredundant resolution $R^\prime$, of
$W^\prime$,
which is also a resolution of $W$. 
Elements of $R(W,i)$, where $i=\ref{h=1},\ref{f=0}$ or $\ref{sameword}$ have smaller $H^\Lm$--length than $W$. The
element of $R(W,\ref{a=pp})$ has smaller $H$--length than $W$. In the case of $W^\prime\in R(W,\ref{constpoly})$
the sum of the $H^\Lm$--length and exponential length is smaller than that of $W$. No operation involved in positive 
$R(1,\ref{sameword})$--reduction increases any of $H^\Lm$--length, $H$--length or exponential length. Furthermore 
once $H^\Lm(\ref{negexp})$-redundancy has been removed it is not reintroduced by any subsequent 
removal of $H^\Lm(i)$--redundancy, for any other $i$. Hence we may assume that in this application of positive 
$H^\Lm(1,\ref{sameword})$--reduction to form $R^\prime$ we either reduce the sum of $H^\Lm$--length, $H$--length and
exponential length, or do nothing, and so introduce no new occurrence of $H^\Lm(\ref{cancelword})$--redundancy. 
(Note that if $\cL$ implies
$f-1$ or $g-1$ is equal to $0$ then the process of
$H^\Lm$--reduction removes $(a,f)$ or $(b,g)$ from the image of
$L_D(\mbf w)$ under $\bt$.)

In  forming $R^\prime$ we construct $\bt^\prime$ by replacing
occurrences of $u_0=(a,f)(b,g)$ by a word $u_1$. We call this the
$R7${\em --process} at $u_0$. If $(u_1,\cL)$ is
$H^\Lm(7)$--redundant the $R7$--process can be repeated at $u_1$,
replacing $u_1$ with $u_2$. If this can be repeated $n$ times we
obtain a sequence $u_0,\ldots ,u_n$ and we say that the
$R7$--process can be repeated $n$ times at $u_0$. Suppose that
there exists an integer $K$ (dependent only on $(a,f)$ and
$(b,g)$) such that the  $R\ref{cancelword}$--process can be
repeated only $K$ times, at $(a,f)(b,g)$. Then, provided this
holds at all occurrences of $H^\Lm(\ref{cancelword})$--redundancy,
the $R\ref{cancelword}$--process can be repeated at any further
occurrences of $H^\Lm(\ref{cancelword})$-redundancy to form an
$H^\Lm$--irredundant, positive resolution of $W$.

If $(a,f)$ is a degenerate exponential $H$--letter then $f\in \ZZ$ and the
$R\ref{cancelword}$--process
can be repeated at most $f$ times at $(a,f)(b,g)$.
A similar statement holds if $(b,g)$ is degenerate.
Suppose then that $(a,f)$ and $(b,g)$ are proper exponential $H$--letters.
Let $l(a)=m$, $l(b)=n$ and
assume first that  $\cL$ implies that
$f\equiv 0 (\mod m)$. Suppose also that the $R\ref{cancelword}$--process can be repeated
$m+n$ times at $(a,f)(b,g)$.
Then, for all $r$ with $1\le r\le m+n$, we have $a^{-1}\^ r=b\^ r$.
This implies in particular that $m=1$ if and only if $n=1$. If $m=n$ then $W$ has
$H^\Lm(\ref{sameword})$--redundancy, so we may assume $m\neq n$ and $m,n>1$. Let $r=m+n$ and
$u=a^{-1}\^ r=b\^ r$. Then $u$ has periods $m$ and $n$ so, from \cite{\Ha}, $u$ has period
$e=\textrm{gcd}(m,n)$. As $e|m$ and $e|n$, both $a$ and $b$ have period $e$. Thus $W$ also
has $H^\Lm(\ref{a=pp})$--redundancy, a contradiction. Hence, when $\cL$ implies $f\equiv 0(\mod m)$,
the $R\ref{cancelword}$--process cannot continue for more than
$m+n-1$ iterations at $(a,f)(b,g)$.

Now suppose that $\cL$ implies $f\equiv k(\mod m)$, with $0\le k\le m-1$. After $k$ iterations
of the $R\ref{cancelword}$--process
$(a,f)(b,g)$ is replaced by $(a,f-k)(b^\prime,g-k)$, for some cyclic permutation $b^\prime$ of
$b$.  As $\cL$ implies $f-k\equiv 0 (\mod m)$, it follows that at most $m+n-1$ further
iterations are possible. Hence the $R\ref{cancelword}$--process
cannot be repeated more than $m+n+k-1<2m+n-1$ times at $(a,f)(b,g)$.

Hence in all cases there is an integer $K$ computable from $(a,f)(b,g)$ such that
the $R\ref{cancelword}$--process can be repeated only $K$ times at $(a,f)(b,g)$.
Thus,
we can effectively  form an
$H^\Lm$--irredundant, positive resolution of $W$.  
As in the proof of Lemma \ref{Red} we may assume that $W$ 
and this resolution of $W$ are both in standard form. 
This
completes the proof of the Lemma.\\

Let $u\in F(H^\Lm)$ constrained by a system of parameters $\cL$. We say $(u,\cL)$ is
{\em cyclically} $H^\Lm${\em --redundant} if $(u^\prime,\cL)$ is $H^\Lm$--redundant, for
some cyclic permutation $u^\prime$ of $u$.
If $(u,\cL)$ is $H^\Lm$--irredundant and $(u^\prime,\cL)$ is cyclically
$H^\Lm$--redundant for some cyclic permutation of $u$ we say $(u,\cL)$
is {\em strictly cyclically} $H^\Lm${\em --redundant}. In this case if
$u=p_1\cdots p_n$, with $p_i\in (H^\Lm)^{\pm 1}$ then the word $p_np_1$
must be $H^\Lm(i)$--redundant, for $i=5,6$ or $7$.

The equation $(\mbf w=1,\bt,\cL)$ is {\em cyclically} $H^\Lm${\em
--redundant} if $\bt(d)$ is cyclically $H^\Lm$--redundant for some $d\in
L_D(\mbf w)$ and {\em cyclically} $H^\Lm${\em --irredundant} otherwise.

Let $W=(\mbf w=1,\bt,\cL)$ be $H^\Lm$--irredundant. If $d\in
L_D(\mbf w)$ and $\bt(d)$ is cyclically $H^\Lm$--redundant then we
have $\bt(d)=(b,g)p(a,f)$, where $(b,g)$, $(a,f)\in H^\Lm$, $p\in
F(H^\Lm)$, and $(a,f)(b,g)$ has $H^\Lm(i)$--redundancy, with $i=5,6$
or $7$. Given this situation we define a set $R^\prime
(W,i)=\{(\mbf w=1,\bt_i,\cL) \}$, where
$\bt_i(d^\prime)=\bt(d^\prime)$, for $d^\prime\neq d$ and $\bt_i(d)$
is dependent on $i$ as follows.
\be
\item[$i=$\ref{constpoly}:] $\bt_{\ref{constpoly}}(d)=p(c,h)$,
where $(c,h)$ is as defined in $R(W,\ref{constpoly})$.
\item[$i=$\ref{sameword}:] $\bt_{\ref{sameword}}(d)=p(a,f+\e g)$, where $(a,f)$, $(b,g)$ and $\e$
are as described in the definition of $H^\Lm(\ref{sameword})$--redundancy.
\item[$i=$\ref{cancelword}:] Let $(a,f)(b,g)$ satisfy the conditions for $H^\Lm(\ref{cancelword})$-redundancy.
In the case where $(a,f)$ and $(b,g)$ are both proper exponential $H$--letters
then let
$\bt_{\ref{cancelword}}(d)=(b_2\cdots b_nb_1,g-1)p(a,f-1)$. 
If $(a,f)$ is a proper exponential $H$--letter but $(b,g)$ is not then
$\bt_{\ref{cancelword}}(d)=(\t(b,g-1),g-1)p(a,f-1)$.
If $(b,g)$ is a
proper exponential $H$--letter but $(a,f)$ is not then
$\bt_{\ref{cancelword}}(d)=(b_2\cdots b_nb_1,g-1)p(\io(a,f-1),f-1)$. 
\ee

\begin{lemma}\label{Cirr} 
Let $W=(\mbf w=1,\bt,\cL)$ be a constrained quadratic exponential equation. Then
there exists an algorithm which, given $W$, outputs a cyclically $H^\Lm$--irredundant, positive
resolution of $W$ in standard form.
\end{lemma}

\noindent {\em Proof.} Given Lemma \ref{Irr} we may assume $W$ is
$H^\Lm$--irredundant and positive. As in the proof of Lemma \ref{Irr} we
may use the sets $R^\prime(W,i)$ to remove occurrences of
$H^\Lm$-redundancy in $p(a,f)(b,g)$, at each stage replacing $W$ by a
set with elements which are $H^\Lm$-irredundant, positive and such that
either the sum of $H^\Lm$--length, $H$--length and 
exponential length is smaller or there are
fewer occurrences of cyclic $H^\Lm$-redundancy than at the previous
stage. It therefore suffices to show that each of the sets
$R^\prime(W,i)$ is a resolution of $W$.

Suppose we have an occurrence of cyclic $H^\Lm$-redundancy at $\bt(d)$, as
described above, and suppose $d$ occurs in $w_i$. Let
$w_i=u_0x^{-1}dxu_1$ and $\bt(d)=(b,g)p(a,f)$. Suppose $(\phi,\al)$ is a
solution to $W$. Define $\phi^\prime :D\dcup X\maps H$ by
$\phi^\prime(v)=\phi(v)$, if $v\in D\dcup X\backslash\{x,d\}$. Define
$\phi^\prime(d)=\hat\al\bt_i(d)$. If $i=5$ or $6$ define
$\phi^\prime(x)=(b\^\al(g))^{-1}\phi(x)$. If $i=7$ define
$\phi^\prime(x)=b^{-1}_1\phi(x)$. Then clearly $\phi^\prime
|_{L_D(\mbf w)}=\hat\al\bt_i$ and it is easy to check that $\phi^\prime(w_i)=\phi(w_i)$.
Hence $(\phi^\prime,\al)$ is a solution to the unique element $W_i$ of
$R^\prime(W,i)$. Conversely, if $(\phi^\prime,\al)$ is a solution to
$W_i$, for some $i$, then define $\phi :F(D\dcup X) \maps H$ by
$\phi(v)=\phi^\prime(v)$, for $v\in D\dcup X\backslash \{x,d\}$,
$\phi(d)=\hat\al\bt(d)$ and $\phi(x)$ using the same rules as above,
relating $\phi$ and $\phi^\prime$. Then $(\phi,\al)$ is a solution to
$W$.
Hence $R^\prime(W,i)$ is a resolution of $W$ in each case.  
As in the proof of Lemma \ref{Red} we may assume that $W$ 
and these resolutions of $W$ are  in standard form.\\

Let $a$ be a reduced word in $H$ and $\mbf s$ a set of words of $H$. 
Suppose that $a$ has a reduced subword $u$ such that
$uv$ is a cyclic permutation of an element of $\mbf{{s}}{^{\pm 1}}$, 
for some $v\in H$, with
$l(v)< l(u)$. Then we say that $a$ is {\em relator--reducible} with respect
to $(H,\mbf s)$. In this case, if
$a=a_0ua_1$, where $l(a)=l(a_0)+l(u)+l(a_1)$ and $a_i\in H$, we say that the reduced word
representing $a_0v^{-1}a_1\in H$ is an {\em elementary relator--reduction} of $a$,
with respect to $(H,\mbf s)$,
{\em replacing} $u$ by $v^{-1}$. If $a$ is not relator--reducible it is said to be
{\em relator--reduced}, with respect to $(H,\mbf s)$.

Let $N$ be the normal closure of $\mbf s\subset H$. If the word problem is solvable in $H$ then given $a\in H$ we may determine whether or not $a$
is relator--reducible and if so perform an elementary relator--reduction. This decreases
the length of $a$ and so we may perform further elementary relator--reductions until a
relator--reduced word $a^\prime \in H$, with $aN=a^\prime N$, is obtained. Then $a^\prime$
is called a {\em relator--reduction} of $a$, with respect to $(H,\mbf s)$.
A degenerate exponential $H$--letter $(a,f)\in H^\Lm$ is said to be 
{\em relator--reducible}, with respect to $(H,\mbf s)$,
if $a$ is relator--reducible, with respect to $(H,\mbf s)$. 
A non--degenerate exponential
$H$--letter $(a,f)\in H^\Lm$ is said to be {\em relator--reducible}
if some cyclic permutation of $a$ is relator--reducible. A exponential $H$--letter which
is not relator--reducible is called {\em relator--reduced}.
Let $\mbf w=\{w_1,\ldots ,w_k\}$ be a system of quadratic words. 
The equation $W=(\mbf w=1,\bt,\cL)$ with environment 
$((X_k)_{k\in K};(H_i)_{i\in I} ;(G_k)_{k\in K})^\Lm$
and 
basis $(a_1,\ldots, a_k)$
is said to be  {\em relator--reducible} 
if an  element $(a,f)\in H^\Lm_{*,a_i}$, which is  relator--reducible with
respect to $\mbf s_{a_j}$, occurs in $\bt(d)$ for some $d\in L_D(w_i)$.  
Otherwise $W$ is said
to be {\em relator--reduced}.

\begin{lemma}\label{Sred} Let $W=(\mbf w=1,\bt,\cL)$ be a
constrained, quadratic exponential equation. Then there exists an algorithm
which, given $W$, outputs a relator--reduced, cyclically
$H^\Lm$--irredundant, positive resolution of $W$ in standard form.
\end{lemma}

\noindent{\em Proof.} Given Lemma \ref{Cirr} we may assume $W$ is
cyclically $H^\Lm$--irredundant and positive. Suppose a
relator--reducible $H_{*,a_i}$--letter $(a,f)$ occurs in $\bt(d)$, for
some $d\in L_D(w_i)$, with $\bt(d)=p_0(a,f)p_1$. Assume first
that $(a,f)$ is degenerate. Define an equation $W^\prime=(\mbf
w=1,\bt^\prime,\cL)$ as follows. Let
$\bt^\prime(d^\prime)=\bt(d^\prime)$, for all $d^\prime\in L_D(\mbf
w)\bl\{d\}$ and $\bt^\prime(d)=p_0(a^\prime,l(a^\prime))p_1$, where
$a^\prime$ is a relator--reduction of $a$ with respect to 
$(H_{*,a_i},\mbf s_{a_i})$.  Clearly
solutions of $W$ and $W^\prime$ coincide and $W^\prime$ has
smaller $H$--length than $W$.  We may now replace $W^\prime$ by a
positive, cyclically $H^\Lm$--irredundant resolution of $W^\prime$,
observing that none of the operations involved increases
$H$--length. This process may be repeated until a positive,
cyclically $H^\Lm$--irredundant resolution of $W$ is obtained, with
no relator--reducible degenerate exponential $H$--letters.

Assume then that all relator--reducible exponential $H$--letters
occurring in $W$ are proper and that $(a,f)$ is
relator--reducible and occurs in $\bt(d)$, with $d\in L_D(w_i)$
and that $\bt(d)=p_0(a,f)p_1$. 
Define {\em crr}$(W)$ to be the sum of $l_H(a,f)$ over all occurrences 
of relator--reducible proper exponential $H$--letters occuring in $W$. 

Suppose first that $a$ is
relator--reducible and  has relator--reduction $a_0$. Let
$a_0=b^{-1}a_1b$, where $a_1$ is cyclically reduced and $b$ is
reduced. Suppose $\cL$ implies that $f\equiv k (\mod l(a))$, with
$0\le k< l(a)$. Let $\lm\in\Lm$ such that $\lm$ does not occur in $W$
and let $c=\io(a,k)$. Suppose $c$ has relator--reduction $c_0$.
Define equations $W_i=(\mbf w =1,\bt_i,\cL_i)$, for $i=0$ and $1$,
as follows. Let $\bt_i(d^\prime)=\bt(d^\prime)$, for all
$d^\prime\neq d$, let $\bt_0(d)=p_0(c_0,l(c_0))p_1$ and let
$\bt_1(d)=p_0(b^{-1},l(b))(a_1,\lm l(a_1))(bc_0,l(bc_0))p_1$. Set
$\cL_0=\cL\cup\{f=k\}$ and $\cL_1=\cL\cup\{g=0,\lm>0\}$, where
$g=f-k-\lm l(a)$.

The pair $(\phi,\al)$ is a solution to $W_0$ if and only
$(\phi,\al)$ is a solution to $W$ with $\al(f)=k$. Suppose
$(\phi,\al)$ is a solution to $W$. If $\al(f)=ql(a)+k$, then we may
assume that $q\ge 0$, since $W$ is positive, and in the light of the previous
remark we assume that $q>0$. Given that
$q>0$, define a retraction $\al^\prime :M\maps \ZZ$ by
$\al^\prime(\mu)=\al(\mu)$, for all $\mu\in \Lm\backslash\{\lm\}$, and
$\al^\prime(\lm)=q$. We have
$\hat\al^\prime\bt_1(d)=\hat\al(p_0)b^{-1}a_1^qbc_0\hat\al(p_1)$, so
$\hat\al^\prime\bt_1(d)N_{a_i}=\hat\al\bt(d)N_{a_i}$. Thus, if $\phi^\prime
:F(L_D(\mbf w))\maps H$ is defined by $\phi^\prime(a)=\phi(a)$,
for $a\neq d$ and $\phi^\prime(d)=\hat\al^\prime\bt_1(d)$ we have
$\phi^\prime|_{L_D(\mbf w)}=\hat\al^\prime\bt_1$ and
$\phi^\prime(w_j)N_{a_j}=\phi(w_j)N_{a_j}=N_{a_j}$, for all $j$. Hence
$(\phi^\prime,\al^\prime)$ is a solution to $W_1$.

Conversely, given  a solution $(\phi^\prime,\al^\prime)$ to $W_1$ we
have, since $\cL_1$ implies that $g=0$, that
$l(a)\al^\prime(\lm)=\al^\prime(f)-k=ql(a)$
so $\hat\al^\prime\bt_1(d)N_{a_i}=\hat\al^\prime\bt(d)N_{a_i}$.
Therefore, if we define $\phi : F(L(\mbf w))\maps H$ by
$\phi|_{L(\mbf w)\backslash\{d\}} = \phi ^\prime|_{L(\mbf
w)\backslash\{d\}}$ and $\phi(d)=\hat\al^\prime\bt(d)$ it follows,
since $\phi(w_j)N_{a_j}=\phi^\prime(w_j)N_{a_j}$, for all $w_j$, that
$(\phi,\al^\prime)$ is a solution to $W$. Hence $R=\{W_0,W_1\}$ may be effectively constructed, given $W$,
and is
a resolution of $W$. Note that although $a_1$ is relator--reduced it may have a cyclic permutation
which is relator--reducible. However $l(a_1)<l(a)$ so crr$(W_1)<$ crr$(W)$, in any case.

Now suppose  $\bt(d)=p_0(a,f)p_1$, $\cL$ implies that $f\equiv k
(\mod l(a))$, $a$ is relator--reduced but that $a=a^\prime
a^\pprime$ and the cyclic permutation $a^\pprime a^\prime$ of $a$
is relator--reducible. Let $a^\pprime a^\prime$ have
relator--reduction $a_0$ and let $a_0=b^{-1}a_1b$, where $a_1$ is
cyclically reduced. Define the equations $W_i=(\mbf
w=1,\bt_i,\cL_i)$, for $i=0,2$ and $3$ as follows. Let $W_0$ be as
above and, for $i=2$ and $3$, let $\bt_i|_{L_D(\mbf
w)\bl\{d\}}=\bt|_{L_D(\mbf w)\bl\{d\}}$. Define
\[\bt_2(d)=p_0(a_2,l(a_2))p_1,\] where $a_2$ is the
relator--reduction of $a\^(l(a)+k)$. Define
\[\bt_3(d)=p_0(d_0,l(d_0))(a_1,\lm l(a_1))(d_1,l(d_1))p_1,\] where $\lm\in
\Lm$ but $\lm$ does not occur in $W$, $d_0$ is the
relator--reduction of $a^\prime b^{-1}$ and $d_1$ is the
relator--reduction of $ba^\pprime \io(a,k)$. Let
$\cL_2=\cL\cup\{f-l(a)+k=0\}$ and $\cL_3=\cL\cup\{g=0,\lm>0\}$,
where $g=f-\lm l(a)-(k+l(a))$. Let $R^\prime$ be the set consisting
of those of $W_0,W_2$ and $W_3$ that are defined. Then as in the
case where $a$ is reducible it follows that $R^\prime$ is a
resolution of $W$ and again crr$(W_i)<$ crr$(W)$, for $i=0,2,3$.

Note that, as $W$ is positive, elements of $R$ or $R^\prime$ are
positive. We may now replace each element of $R$ or $R^\prime$ 
by a cyclically
$H^\Lm$--irredundant, positive resolution, observing that this
process introduces no new occurrences of relator--reducibility at
proper exponential $H$--letters. After removal of any
relator--reducibility which may have been introduced at
degenerate exponential $H$--letters we obtain a resolution of $W$
every element of which is cyclically $H^\Lm$--irredundant,
positive and has no relator--reducible degenerate exponential
$H$--letters.
Furthermore elements $W^\prime$ of this resolution satisfy crr$(W^\prime)<$
crr$(W)$ since this is true of $R$ and $R^\prime$.
Continuing this
process we obtain, eventually, the required resolution (when crr$(W^\prime)=0$). 
As in the proof of Lemma \ref{Red} we may assume that $W$ 
and this resolution are both   in standard form.\\

As before, let $W=(\mbf w=1,\bt,\cL)$ be a constrained quadratic exponential
equation with basis $(a_1,\ldots, a_k)$. Let $s$ be a cyclic permutation of an element of
$\mbf{{s}}_{a_i}^{\pm 1}$ such that $s=s^m_0$, where $s_0$ is
cyclically reduced, $m\ge 1$ and $s_0\in H_{*,a_i}$. Suppose a proper exponential letter
$(s_0,f)\in
H^\Lm_{*,a_i}$ occurs in the $i$th component of $W$: 
that is there exists $d\in L_D(w_i)$ such
that $\bt(d)=p_0(s_0,f)p_1$, for some $p_0,p_1\in F(H^\Lm_{*,a_i})$. Then we say
that 
$(s_0,f)$ 
is a  
{\em relator--unconstrained} exponential $H$--letter occuring in $\bt(d)$, 
with respect to $(H_{*,a_i},\mbf s_{a_i})$. 
The 
equation $W$ is called {\em relator--unconstrained} if a
relator--unconstrained exponential $H$--letter occurs in $\bt(d)$, for 
some $d\in L_D(\mbf w)$ and
{\em relator--constrained} otherwise.

\begin{lemma}\label{Scon}
Let $W=(\mbf w=1,\bt,\cL)$ be a  constrained, quadratic exponential
equation. Then there exists an algorithm which, given $W$, outputs
an relator--constrained, relator--reduced, cyclically
$H^\Lm$--irredundant, positive resolution of $W$ in standard form.
\end{lemma}

\noindent{\em Proof.} Given Lemma \ref{Sred} we may assume $W$ is
relator--reduced, cyclically $H^\Lm$--irredundant and positive.
Suppose that a relator--unconstrained exponential $H$--letter
$(s_0,f)$ occurs in $W$, say $\bt(d)=p_0(s_0,f)p_1$, for some $d\in
L_D(w_u)$, where $s=s_0^m$ is a cyclic permutation of an element of
$\mbf{{s}}_{a_u}^{\pm 1}$, $m\ge 1$, $1\le u\le k$. 

Let $\{j_1,\ldots, j_p\}$ be the set of integers such that
$0<j_i<l(s)/2$ and $\cL_i=\cL\cup\{f\equiv -j_i(\mod(l(s))\}$
is consistent, for $i=1,\ldots ,p$. 
Let $\{j_{p+1},\ldots, j_{p+q}\}$ be the set of integers such that
$0\le j_i\le l(s)/2$ and $\cL_i=\cL\cup\{f\equiv j_i(\mod(l(s))\}$
is consistent, for $i=p+1,\ldots ,p+q$. Let $T=L_D(\mbf w)$ and, for $i=1,\ldots ,p$, define
$\bt_i:L_D(\mbf w)\maps F(H^\Lm)$ by 
$\bt_i|_{T\backslash \{d\}}= \bt|_{T\backslash \{d\}}$
and $\bt_i(d)=p_0(s_0^{-1}\^j_i,j_i)p_1$. For $i=p+1,\ldots ,p+q$,
define 
$\bt_i:F(H^\Lm)\maps H$ by 
$\bt_i|_{T\backslash \{d\}}= \bt|_{T\backslash \{d\}}$
and $\bt_i(d)=p_0(s_0\^j_i,j_i)p_1$.
Let $W_i=(\mbf w=1,\bt_i,\cL_i)$, for $i=1,\ldots ,p+q$ and 
$R=\{W_1,\ldots ,W_{p+q}\}$. 

For all $i$ it follows that $W_i$ is a positive, constrained,
quadratic exponential equation with smaller exponential length than $W$. 
To see that $R$ is a resolution of $W$, suppose first that $(\phi,\al)$ 
is a solution to $W$.  Let $\al(f)=bl(s)+c$, where $0\le c<l(s)$, 
$b,c\in \ZZ$.  There is either 
an $i$ such that $c= j_i$, with $p+1\le i\le p+q$, or there is an $i$ such that
$c=l(s)-j_i$, with $1\le i\le q$. In either case $\al$ is a solution to $\cL_i$.  
Fix such an
$i$ and 
define $\phi^\prime :F(L(\mbf w))\maps H$ by
$\phi^\prime |_{L(\mbf w)\bl \{d\}}=\phi|_{L(\mbf w)\bl \{d\}}$ and
$\phi^\prime(d)=\hat\al\bt_i(d)$.
Assume first that $p+1\le i\le p+q$ so that $0\le j_i\le l(s)/2$. 
 Then 
\begin{align*}
\phi^\prime(d)&=\hat\al(p_0)(s_0\^j_i)\hat\al(p_1)\\
 & = \hat\al(p_0)s^{-b}(s_0\^\al(f))\hat\al(p_1).
\end{align*}
Thus $\phi^\prime(d)N_{a_u}=\phi(d)N_{a_u}$ and $(\phi^\prime,\al)$ is a solution to $W_i$.
Similarly, if $1\le i\le q$, so that $0< j_i< l(s)/2$, then
\begin{align*}
\phi^\prime(d)&=\hat\al(p_0)(s_0^{-1}\^j_i)\hat\al(p_1)\\
 &=\hat\al(p_0)s^{-b-1}s^{b+1}(s_0^{-1}\^(l(s)-c))\hat\al(p_1)\\
 &=\hat\al(p_0)s^{-b-1}(s_0\^(l(s)b+c))\hat\al(p_1)\\
 & = \hat\al(p_0)s^{-b-1}(s_0\^\al(f))\hat\al(p_1).
\end{align*}
Thus $\phi^\prime(d)N_{a_u}=\phi(d)N_{a_u}$ and $(\phi^\prime,\al)$ is a solution to $W_i$.

Conversely, if $(\phi^\prime,\al^\prime)$ is a solution to $W_i\in
R$ then $\al^\prime$ is a solution to $\cL$, as $\cL\subset \cL_i$.
Define $\phi : F(L(\mbf w))\maps H$ by $\phi|_{L(\mbf w)\backslash
\{d\}}=\phi^\prime |_{L(\mbf w)\backslash \{d\}}$ and
$\phi(d)=\hat\al^\prime\bt(d)=\hat\al^\prime(p_0)(s_0\^\al^\prime(f))\hat\al^\prime (p_1)$. 
If $1+p\le i\le p+q$ then $\al^\prime(f)=bl(s)+j_i$, for some $b\in\ZZ$, so as above
it follows that $\phi^\prime(d)N_{a_u}=\phi(d)N_{a_u}$ and that $(\phi,\al^\prime)$ is a
 solution to $W$. Similarly, if $1\le i\le q$ then $\al^\prime(f)=bl(s)-j_i$, for some $b\in\ZZ$, 
so
\begin{align*}
\phi^\prime(d)&=\hat\al^\prime(p_0)(s_0^{-1}\^j_i)\hat\al^\prime(p_1)\\
 & = \hat\al^\prime(p_0)s^{-b}(s_0\^\al^\prime(f))\hat\al^\prime(p_1),
\end{align*}
and $(\phi,\al^\prime)$ is a solution to $W$.

Therefore $R$ is a resolution of $W$. Using Lemma \ref{Sred}, we may replace $R$ by a resolution
which is relator--reduced, cyclically
$H^\Lm$--irredundant, positive  and in standard form.  None
of the operations involved increase exponential length and so 
repetition of this process gives the required resolution.\\

Let $W=(\mbf w=1,\bt,\cL)$ be a constrained quadratic exponential
equation. Then $W$ is called {\em singular} if 
$\bt(d)=1\in F(H^\Lm)$, for some $d\in L_D(\mbf w)$ and {\em non--singular} otherwise. The set
$S(W)=\{d\in L_D(\mbf w)\,:\, \bt(d)=1\}$ is called the set of 
{\em singularities} of $W$.
 
\begin{lemma}\label{sing}
Let $W=(\mbf w=1,\bt,\cL)$ be a  constrained, quadratic exponential
equation. Then there exists an algorithm which, given $W$, outputs
a non--singular, relator--constrained, relator--reduced, cyclically
$H^\Lm$--irredundant, positive resolution of $W$ in standard form.
\end{lemma}

\noindent{\em Proof.} We may assume that $W$ is in standard form so that
$\mbf w=\mbf q(\mbf n,\mbf t,\mbf p)$, where $(\mbf n,\mbf t, \mbf p)$ is 
some positive 3--partition of $(n,t,p)\in \ZZ^3$.
With the above notation assume $S(W)=\{d_{i_1},\ldots,d_{i_s}\}$, where 
$0<i_1 < \cdots < i_s$. Define $i_0=0$, $i_{s+1}=n+2t+p+1$ and 
$T=\{i_1,\ldots ,i_{s}\}$  and define
a map \[\s:\{1,\ldots ,n+2t+p\}\backslash T \maps \{1,\ldots ,n+2t+p-s\},\]
by \[\s(i)=i-j, \textrm{ when } i_j<i <i_{j+1}.\]

Suppose $\mbf w=w_1,\ldots , w_k$, with $w_j=q(\xi_j,n_j,t_j,p_j)$, for $j=1,\ldots ,k$.
Let $S_j=T\cap \{\xi_j+1,\ldots ,\xi_j+n_j\}$, for $j=1,\ldots ,k$. Form a new system
of equations $W^\prime=(\mbf w^\prime=1,\bt^\prime,\cL)$ from $W$ as follows. 
If $n_j-|\,S_j|=t_j=p_j=0$ then
delete $w_j$ from the system. Otherwise replace $w_j$ with a new equation $w_j^\prime$ 
formed by 
\be
\item deleting $x_i^{-1}d_ix_i$, from $w_j$ when $i\in S_j$,
\item replacing $x_i^{-1}d_ix_i$ with 
$x_{\s(i)}^{-1}d_{\s(i)}x_{\s(i)}$, when $i\notin S_j$,
\item replacing $[x_i,x_{i+t_j}]$ with $[x_{\s(i)},x_{\s(i+t_j)}]$, when 
$\xi_j+n_j+1\le i\le \xi_j+n_j+t_j$,
\item replacing $x_i^2$ with $x_{\s(i)}^2$, when 
$\xi_j+n_j+2t_j\le i\le \xi_j+n_j+2t_j+p_j$.
\ee
Define $\bt^\prime:L_D(\mbf w^\prime)\maps F(H^\Lm)$ by 
\[\bt^\prime(d_i)=\bt(d_{\s^{-1}(i)}).\]

Now suppose that $(\phi,\al)$ is a solution to $W$. Define 
\[\phi^\prime(a_i)=\phi(a_{\s^{-1}(i)}),\]
for $a_i=x_i\in L_X(\mbf w)$ and $a_i=d_i\in L_D(\mbf w).$ 
By construction $(\phi^\prime,\al)$
is then a solution to $W^\prime$. Conversely, let $(\phi^\prime,\al)$ be a solution to
$W^\prime$. For $i\notin T$,  define $\phi(a_i)=\phi^\prime(a_{\s(i)})$, where 
$a_i=x_i\in L_X(\mbf w)$ 
or $a_i=d_i\in L_D(\mbf w)$. Define $\phi(d_i)=\phi(x_i)=1$, for $i\in T$.
Then $(\phi,\al)$ is a solution to $W$. Hence $\{W^\prime\}$ is a resolution of $W$.

Given Lemma \ref{Scon} we may assume that $W$, and so also $W^\prime$, is
relator--constrained, relator--reduced, cyclically
$H^\Lm$--irredundant and positive.\\

The process of forming the resolution of Lemma \ref{sing} is
called {\em reduction of singularities}.
%

A parameter system $\cL$ is said to be {\em normalised} if
$\cL=\cC\dcup\cE\dcup\cI$, where: 
\bd
\item[$\cP1$] $\cC$ is a finite set of congruences of the form 
$f\equiv 0(\mod m)$, where $f\in M$ and $m> 0$, $m\in \ZZ$;
\item[$\cP2$] $\cE$ is a finite set of diophantine equations of the form
$g=0$, $g\in M$;
\item[$\cP3$] $\cI$ is a finite set of diophantine inequalities of the
form $h\ge 0$, for $h\in M$;
\item[$\cP4$] if $c\in \ZZ$ and $c$ occurs in a congruence of $\cC$
with modulus $m$ then  $0\le c<m$. 
\ed
Henceforward we shall
assume that if $(\mbf w=1,\bt,\cL)$ is a constrained quadratic exponential equation, then
$\cL$ is normalised.

A resolution of a  constrained, quadratic exponential equation which is in
standard form, cyclically $H^\Lm$--irredundant, positive, relator--reduced,
relator--constrained, non--singular
and normalised 
is called {\em special}.

\begin{theorem}\label{special}
Let $W=(\mbf w=1,\g,\cL)$ be a constrained,   quadratic exponential equation with environment 
$((X_k)_{k\in K};(H_i)_{i\in I} ;(G_k)_{k\in K})^\Lm$. Then
there exists an algorithm which, given $W$ outputs a special resolution of $W$.
\end{theorem}

Now let $(\mbf h,\mbf n,\mbf t,\mbf p)$ be a positive 4--partition and  
$\mbf z\in (\dcup_{k\ge 1}H_{*,k}^\Lm,\cL,\mbf n,\mbf h).$  
If $Q(\mbf z,\cL,\mbf n,\mbf t,\mbf p)$
is special we call $\mbf z$ a {\em special} element of 
$(\dcup_{k\ge 1}H_{*,k}^\Lm,\cL,\mbf n,\mbf h)$.
\begin{lemma}\label{homog-special}
If $\mbf z$ is special then the homogeneous equation 
\[Q_{\cH}=Q(\ovr{\mbf z},\cH,\mbf n,\mbf t,\mbf p)\] 
associated to $Q(\mbf z,\cL,\mbf n,\mbf t,\mbf p)$ is special.
\end{lemma}

\section{Decision problems}\label{decisions}
Let $I$ be a recursive set and, for all $i\in I$, let $C_i$ be a recursive 
set and $R_i$ a recursively enumerable  
subset of $F_i=F(C_i)$. Let $H_i$ be the group with presentation
$\langle C_i|R_i\rangle$, for all $i\in I$. Let $K=I$ and set $X_k=\{k\}$, for
all $k\in K$. Then $F_{*,i}=F_i$, for all $i\in I$.

The {\em quadratic exponential equation problem} or 
${\cPQ}${\em --problem} in the indexed set $(\langle C_i|R_i\rangle : i\in I)$ 
is that  of determining whether or not, given a
constrained quadratic equation $W=({\mbf w}=1,\bt,\cL)$   with environment
$((\{i\})_{i\in I};(F_i)_{i\in I};(H_i)_{i\in I})^\Lm$, the equation $W$
has a solution and, if so, of finding one. In fact, given that the word problem is
solvable in $H_i$, for all $i$, once the existence of a solution to $W$ has been established
then a solution $(\phi,\al)$ may be found by enumeration of all possible values of $\phi(x_j)$
and $\al(\lm_j)$, for $x_j$ and $\lm_j$ occuring in $W$. However the algorithms described here
for the ${\cPQ}$--problem will determine whether a solution exists by constructing one.
A similar comment applies to all the other decision problems considered herein.  

Since we wish to consider the extension of solvability from indexed families of groups
such as $(H_i)_{i\in I}$ to free--products and one--relator products of such groups we
reformulate this problem in terms of free--products. 
Take $I$, $C_i$, $R_i$, $F_i$ and $H_i$ to be defined
as at the beginning of this section. As in Section \ref{qpequations}  
let $K$ be a recursive set and $X_k\subset I$ be a
recursive subset of $I$, for all $k\in K$. Define $H_{*,k}$ as in (\ref{Hsys}),
let $\mbf s_k$ be a recursive subset of $H_{*,k}$, all elements of which have length at least $2$,  
and define $N_k$ and $G_k$ as in 
(\ref{Gsys}). The groups $H_{*,k}$ and $G_k$ are equipped with natural 
presentations as defined in Section \ref{qpequations}.

The ${\cPQ}${\em --problem} in the indexed set $(G_k: k\in K)$ is to determine, given a
constrained quadratic equation $W$  with environment 
\[((X_k)_{k\in K};(H_i)_{i\in I} ;(G_k)_{k\in K})^\Lm,\] 
whether
or not $W$ has a solution and, if so, to find one.
The {\em special quadratic exponential equation problem} or 
$\cSPQ${\em --problem} in $(G_k: k\in K)$ is defined in the same way as the
$\cPQ$--problem with the additional condition that $W$ is special.

We say that the ${\cPQ}$--problem in $(G_k: k\in K)$ is {\em solvable} if there is
an algorithm which, given input such a constrained quadratic exponential 
equation outputs a solution if one exists and otherwise indicates the
absence of such a solution.
{\em Solvability} of 
the ${\cSPQ}$--problem in $(G_k: k\in K)$ is defined analogously.

Now let $(\mbf h,\mbf n,\mbf t,\mbf p)$ be a positive 4--partition. 
In the notation above
the {\em exponential genus problem at} 
$(\mbf h,\mbf n,\mbf t,\mbf p)$, or
$\cPG(\mbf h,\mbf n,\mbf t,\mbf p)${\em --problem}, in $(G_k: k\in K)$ is that of
determining whether or not, given a consistent parameter system $\cL$
and an $n$--tuple $\mbf z$ of elements of $(\dcup_{k\ge 1}H_{*,k}^\Lm,\cL,\mbf n,\mbf h)$, 
the triple $(\mbf n,\mbf t,\mbf p)\in \cL$--genus$(\mbf z)$ and
if so of finding a solution to 
$Q(\mbf z,\cL,\mbf n,\mbf t,\mbf p)$. 
The {\em special exponential genus problem at} $(\mbf h,\mbf n,\mbf t,\mbf p)$, or
$\cSPG(\mbf h,\mbf n,\mbf t,\mbf p)${\em --problem}, in $(G_k: k\in K)$ is defined in
the same way as
the exponential  genus problem with the additional constraint that
$\mbf z$ must be special. 
Note that these are decision problems for
presentations of the groups $G_k$ given by some fixed environment. 
In particular the $\cPG(\mbf h,\mbf n,\mbf t,\mbf p)$--problem in 
$(\langle C_i|R_i\rangle , i\in I)$ 
is to 
determine whether or not, given 
an $n$--tuple $\mbf z$ of elements of $(\dcup_{i\in  I}F_{i}^\Lm,\cL,\mbf n,\mbf h)$, 
the triple $(\mbf n,\mbf t,\mbf p)\in \cL$--genus$(\mbf z)$ and
if so to find a solution to 
$Q(\mbf z,\cL,\mbf n,\mbf t,\mbf p)$, an equation
with environment $((\{i\})_{i\in I};(F_i)_{i\in I};(H_i)_{i\in I})^\Lm$.

The $\cPG(\mbf h,\mbf n,\mbf t,\mbf p)$--problem in $(G_k: k\in K)$
is {\em solvable} if there
exists an algorithm which, given a consistent parameter system $\cL$
and an $n$--tuple $\mbf z$ of elements of 
$(\dcup_{k\in  K}H_{*,k}^\Lm,\cL,\mbf n,\mbf h)$, 
decides whether or not $(\mbf n,\mbf t,\mbf p)\in
\cL$--genus$(\mbf z)$ and if so outputs a solution to 
$Q(\mbf z,\cL,\mbf n,\mbf t,\mbf p)$. The 
$\cSPG(\mbf h,\mbf n,\mbf t,\mbf p)$--problem is {\em solvable} if the
obvious analogous conditions are satisfied.
\begin{exx}\label{exx_CE}
Let $I=K=X_1=\{1\}$, let $C_1=\{c_1,\ldots,c_n\}$ and let $R_1=\emptyset$.
Then $H_1=F_1$, the free group of rank $n$. The $\cPQ$--problem in $F_1$ is to 
determine, given a quadratic exponential equation with environment $(X_1;F_1;F_1)^\Lm$, whether or not the equation
has solution and if so of finding one. This problem is shown to be solvable by Comerford
and Edmunds in \cite{\CEa}. 
\end{exx}
\begin{exx}\label{exx_LM}
In \cite{\LM} Lipschutz and Miller describe a family of decision problems  several of which
may be given in terms of quadratic exponential equations. 
We repeat the descriptions of \cite{\LM} here for such cases.
Given any two of these problems
Lipschutz and Miller construct a group presentation in which the first problem is solvable and the 
second is not. 
\be
\item The {\em order problem} is to determine the order of a given group element. With the 
notation of Example \ref{exx_z_probs} 
we see (loc. cit. (\ref{exx_z_probs_1})) that if the $\cPG((1),(1),(0),(0))$--problem is solvable in $G_1$ then we can determine whether or not 
the element $h_1\in H_{*,1}$ has
finite order in $G_1$. If so then we obtain a non--zero integer $n$ such that $h_1^n=1$ in $G_1$.
Define equations \[W_j=((w)=1,\bt_j,\cL_0), \textrm{ for } j=1,\ldots ,n,\]
where $w=d_1$, $\bt_j(d_1)=(h_1^j,l(h_1^j))$ and $\cL_0=\emptyset$. The order of $h_1$ is the smallest
$j$ such that $W_j$ has a solution. The $W_j$'s are instances of the $\cPG((0),(1),(0),(0))$--problem in
$G_1$ and so are also instances 
of the $\cPG((0),(1),(0),(0))$--problem. Hence if the $\cPG((0),(1),(0),(0))$--problem is
solvable in $G_1$ we may solve the order problem in $G_1$.
\item The {\em power problem} is to determine, given two elements of a group, whether or not the first is 
a power of the second. Again Example \ref{exx_z_probs}(\ref{exx_z_probs_1}) shows that we can solve
this problem in $G_1$ if the $\cPG((0),(1),(0),(0))$--problem is
solvable in $G_1$. 
\item The {\em generalised power problem} is to determine, given two elements of a group, 
whether or not a power of the first is 
a power of the second.  Example \ref{exx_z_probs}(\ref{exx_z_probs_2}) shows that we can solve
this problem in $G_1$ if the $\cPG((2),(1),(0),(0))$--problem is
solvable in $G_1$. 
\ee
\end{exx}
\begin{exx}\label{exx_power_conj}
Continuing the development of the previous example, the {\em power conjugacy problem} is to 
determine, given two elements of a group, 
whether or not a power of the first is conjugate to 
a power of the second. From Example \ref{exx_z_probs}(\ref{exx_z_probs_3}) it follows that that we can solve
this problem in $G_1$ if the $\cPG((2),(2),(0),(0))$--problem is
solvable in $G_1$. 
\end{exx}

Now let $\cF$ be a collection of recursive group presentations (elements of which we'll
refer to as groups).
The $\cPQ${\em --problem} in $\cF$ is the union of the $\cPQ$--problem over all  
environments $((\{i\})_{i\in I};F(C_i)_{i\in I};(H_i)_{i\in I})^\Lm$, where 
$I$ is a recursive set, $(\langle C_i:R_i\rangle)_{i\in I}$ is a family of presentations indexed by $I$, 
such that $\langle C_i:R_i\rangle$ is an element of $\cF$ and a  presentation for $H_i$, for all $i\in I$. 
The $\cSPQ$--problem, the $\cPG(\mbf h,\mbf n,\mbf t,\mbf p)$--problem, 
the $\cSPG(\mbf h,\mbf n,\mbf t,\mbf p)$--problem, the $\cPG$--problem and 
 the $\cPG$--problem
in $\cF$ are all  defined analogously.

Given the collection $\cF$ define $\cF^*$ to be the collection of groups of the 
form $\ast_{i\in I} H_i$, where $I$ is a recursive set and $H_i$ belongs to $\cF$, 
for all $i\in I$. Each element of $\cF^*$ is given by its natural presentation.
Note that every element of $\cF$ belongs to $\cF^*$.
Next define $\cF^*_{\sim}$ 
to be the collection consisting of all groups of $\cF^*$ together
with groups of the  form $H/N$, where $H=\ast_{i\in I} H_i$ is an element of $\cF^*$,
 as above, with $I$ a set of size at least $2$,  and  $N$ is the normal closure of 
a recursively enumerable subset of elements, of length at least $2$, of $H$. 
Again each element of $\cF^*_{\sim}$ is given
by its natural presentation.

Given a subcollection $\cC$ of groups of $\cF^*_{\sim}$ the $\cPQ$--problem 
in $\cC$ can therefore
be described as follows. Let $I$ and $K$ be a recursive sets 
such that $\langle C_i|R_i\rangle$
is a presentation of a group $H_i$ of $\cF$, for all $i\in I$. Let $X_k$ be 
a recursive subset of $I$ and $H_{*,k}$ be defined by (\ref{Hsys}), for each $k\in K$.
Let $\mbf s_k$ be a recursively enumerable  subset of elements, of length at least $2$, 
of  $H_{*,k}$ and let $N_k$ and $G_k$ be defined
by (\ref{Gsys}), for each $k\in K$. Finally assume that these sets have been chosen
so that $G_k$ belongs to $\cC$, for all $k\in K$. The $\cPQ$--problem in $\cC$ is to 
determine, given such data and 
a constrained quadratic exponential equation $W$ with environment
$((X_k)_{k\in K};(H_i)_{i\in I} ;(G_k)_{k\in K})^\Lm$,
whether or not $W$ has a solution and, if so, to find one.
Similarly the $\cPG(\mbf h,\mbf n,\mbf t,\mbf p)$--problem 
in $\cC$ is to determine, given such data, a consistent parameter system $\cL$ and 
an element $\mbf z\in (\dcup_{k\in K} H^\Lm_{*,k},\cL,\mbf n,\mbf h)$ whether
or not $(\mbf n,\mbf t,\mbf p)\in \cL\textrm{--genus}(\mbf z)$ and, if so, to 
find a solution to $Q(\mbf z,\cL,\mbf n,\mbf t,\mbf p)$.
The $\cSPQ$ and $\cSPG(\mbf h,\mbf n,\mbf t,\mbf p)$--problems 
in $\cC$ have similar descriptions, differing only in
that $W$ and $\mbf z$ must be special. The $\cPG$ and $\cSPG$--problems
also have similar descriptions.

\begin{theorem}\label{PQtoSPQ}  Let $\cF$ be a collection of groups with solvable word problem
and let $\cC$ be a subcollection of
$\cF^*_\sim$.
The  ${\cPQ}$--problem in $\cC$ is solvable 
if and only if the
${\cSPQ}$--problem in $\cC$ is solvable.
\end{theorem}

\noindent
{\em Proof.}
If the  ${\cPQ}$--problem is solvable then the  $\cSPQ$--problem 
is obviously also solvable.  Conversely suppose the $\cSPQ$--problem
in $\cC$ is  solvable. Let $W=(\mbf w=1,\bt,\cL)$ be a constrained
quadratic exponential equation. Using the algorithm of Theorem
\ref{special} we obtain a special resolution $R$ of $W$. As the
$\cSPQ$--problem is solvable we may determine, for each element $Y\in R$   
whether or not a
solution to $Y$  exists. If there is such a solution
then, by definition of resolution, we have an algorithm which outputs
a solution of $W$. Otherwise $W$ has no solution.\\[1em]
\medskip

\begin{lemma}\label{specialg}
Let $\cF$ be a collection of groups and let $\cC$ be a subcollection of
$\cF^*_\sim$.
Assume all the groups in $\cF$ have solvable word problem
and let $(\mbf h,\mbf n,\mbf t,\mbf p)$ be a positive 4--partition of $(h,n,t,p)$. 
Then the
following are equivalent.
\be 
\item\label{specialgi} The $\cPG(\mbf h,\mbf n,\mbf t,\mbf p)$--problem 
in $\cC$ is solvable.
\item\label{specialgii} The
$\cSPG(\mbf h^\prime,\mbf{n^\prime},\mbf t^\prime,\mbf p^\prime)$--problem 
in $\cC$ is solvable, for 
all positive 4--partitions 
$(\mbf h^\prime,\mbf n^\prime,\mbf t^\prime,\mbf p^\prime)$ of 
$(h,n^\prime,t,p)$ such that $n^\prime \in \ZZ$ and  
$(\mbf h^\prime,\mbf n^\prime,\mbf t^\prime,\mbf p^\prime)
\le (\mbf h,\mbf n,\mbf t,\mbf p)$.
\item\label{specialgiii} The
$\cSPG(\mbf h^\prime,\mbf{n^\prime},\mbf t^\prime,\mbf p^\prime)$--problem 
in $\cC$ is solvable, for 
all  positive 4--partitions 
$(\mbf h^\prime,\mbf n^\prime,\mbf t^\prime,\mbf p^\prime)
\le_g (\mbf h,\mbf n,\mbf t,\mbf p)$.
\ee
\end{lemma}
\par\noindent{\em Proof.} 
For the duration of the proof choose recursive sets
$I$ and $K$ 
such that $\langle C_i|R_i\rangle$
is a presentation of a group $H_i$ of $\cF$, for all $i\in I$. Choose $X_k$ to be 
a recursive subset of $I$ and let $H_{*,k}$ be defined by (\ref{Hsys}), for each $k\in K$.
Choose a recursively enumerable set $\mbf s_k$ of $H_{*,k}$ and let $N_k$ and $G_k$ be defined
by (\ref{Gsys}), for each $k\in K$. Finally assume that these sets have been chosen
so that $G_k$ belongs to $\cC$, for all $k\in K$.

We show first that (\ref{specialgi}) implies 
(\ref{specialgiii}). Suppose that the 
$\cPG(\mbf h,\mbf n,\mbf t,\mbf p)$--problem in $\cC$ is solvable. Then  the 
 $\cSPG(\mbf h,\mbf n,\mbf t,\mbf p)$--problem is obviously also
 solvable. Assume that
$\mbf n=(n_1,\ldots,n_k)$ is a partition of an integer $n$ and  define
$\mu_j=\mu_{j}(\mbf n)$, for $j=1,\ldots ,k+1$, as in (\ref{muj}). 
Let $(\mbf h^\prime,\mbf n^\prime,\mbf t^\prime,\mbf p^\prime)$ be a positive 4--partition 
with
$\mbf{h^\prime}=(h_1^\prime,\ldots,h_l^\prime)$,
$\mbf{n^\prime}=(n_1^\prime,\ldots,n_l^\prime)$,
$\mbf{t^\prime}=(t_1^\prime,\ldots,t_l^\prime)$, 
$\mbf{p^\prime}=(p_1^\prime,\ldots,p_l^\prime)$
and  $(\mbf h^\prime,\mbf n^\prime,\mbf t^\prime,\mbf p^\prime) \le_g (\mbf h,\mbf n,\mbf t,\mbf p)$. 
Then there exist integers
$j_1,\ldots ,j_{l}$ such that $1\le j_1<\cdots<j_{l}\le k$ with
$h^\prime_i\le h_{j_i}$,
$n^\prime_i\le n_{j_i}$,
$t^\prime_i = t_{j_i}$,
and $p^\prime_i = p_{j_i}$,
for $i=1,\ldots, l$.  Define
$\mu^\prime_j=\mu_{j}(\mbf n^\prime)$, for $j=1,\ldots ,l+1$.
Let $\cL$ be a system of parameters and let  
$\mbf z^\prime=(z_1^\prime,\ldots ,z_{n^\prime}^\prime)$ be an
element of $(\dcup_{k\in K} H^\Lm_{*,k},\cL,\mbf n^\prime,\mbf h^\prime)$ with
basis $(a^\prime_1,\ldots, a^\prime_l)$. 
Define an element
$\mbf z=(z_1,\ldots ,z_n)$ of $(\dcup_{k\in K} H^\Lm_{*,k},\cL,\mbf n,\mbf h)$ as follows. 
For $s$ such that $\mu_{j_i}+1\le s\le \mu_{j_i}+n^\prime_i$ define
$z_s=z^\prime_{s^\prime}$, where $s^\prime=s- \mu_{j_i}+\mu^\prime_{i}$. For all
other $s\in\{1,\ldots ,n\}$ define $z_s=(1,0)\in H^\Lm$. Then 
$\mbf z$ has basis $(a_1,\ldots ,a_k)$, where $a_{j_i}=a^\prime_i$, and $a_j=1$, if
$j\notin \{j_1,\ldots ,j_l\}$.

Recalling (\ref{xij}) and  (\ref{qntp}) let
$\xi_j=\xi_j(\mbf n,\mbf t,\mbf p)$
for $j=1,\ldots k$, $\mbf q=\mbf q(\mbf n,\mbf t,\mbf p)$,
$\xi^\prime_j=\xi_j(\mbf n^\prime,\mbf t^\prime,\mbf p^\prime)$
for $j=1,\ldots l$ and   
$\mbf q^\prime=\mbf q(\mbf n^\prime,\mbf t^\prime,\mbf p^\prime)$. 
Recalling (\ref{bjz}) let $\bt= \bt(\mbf z, \mbf n,\mbf t,\mbf p)$, 
$\bt^\prime= \bt(\mbf z^\prime, \mbf n^\prime,\mbf t^\prime,\mbf p^\prime)$,
$Q=Q(\mbf z,\cL,\mbf n,\mbf t,\mbf p)=(\mbf q=1,\bt,\cL)$
and $Q^\prime=Q(\mbf z^\prime,\cL,\mbf n^\prime,\mbf t^\prime,\mbf p^\prime)=
(\mbf q^\prime=1,\bt^\prime,\cL)$.
%
Given a homomorphism $\phi:F(L(\mbf q))\maps H$ we may construct a homomorphism
$\phi^\prime:F(L(\mbf q^\prime))\maps H$ as follows. For each $i\in\{1,\ldots ,l\}$ and
$s^\prime$ with 
$1+\xi^\prime_i\le s^\prime\le n_i^\prime+\xi^\prime_i$, define
\[\phi^\prime(x_{s^\prime})=\phi(x_s) \textrm{ and } 
\phi^\prime(d_{s^\prime})=\phi(d_s),\]
where $s=s^\prime-\xi_i^\prime+\xi_{j_i}$.
For $s^\prime$ with 
$n^\prime_i+\xi^\prime_i< s^\prime\le n_i^\prime+2t^\prime_i+p^\prime_i+\xi^\prime_i$, define
\[\phi^\prime(x_{s^\prime})=\phi(x_s),\]
where $s=s^\prime-\xi_i^\prime-n^\prime_i+\xi_{j_i}+n_{j_i}$.
If $(\phi,\al)$ is a solution to $Q$ then $(\phi^\prime,\al)$ is a solution to $Q^\prime$.

Similarly, given a homomorphism $\phi^\prime:F(L(\mbf q^\prime))\maps H$ 
we may construct a homomorphism
$\phi^{\prime\prime}:F(L(\mbf q))\maps H$ as follows.
For each $j_i\in\{j_1,\ldots ,j_l\}$ and
$s$ with 
$1+\xi_{j_i}\le s\le n_i^\prime+\xi_{j_i}$, define
\[\phi^{\prime\prime}(x_{s})=\phi^\prime(x_{s^\prime}) \textrm{ and } 
\phi^{\prime\prime}(d_{s})=\phi^\prime(d_{s^\prime}),\]
where $s=s^\prime-\xi_i^\prime+\xi_{j_i}$.
For each $j_i\in\{j_1,\ldots ,j_l\}$ and
$s$ with 
$n_{j_i}+\xi_{j_i}< s\le n_{j_i}+2t_{j_i}+p_{j_i}+\xi_{j_i}$, define
\[\phi^{\prime\prime}(x_{s})=\phi^\prime(x_{s^\prime}),\] where 
$s=s^\prime-\xi_i^\prime-n^\prime_i+\xi_{j_i}+n_{j_i}$. 
For all other $i$, in an appropriate
range,  define $\phi(d_i)=1$
and $\phi(x_i)=1$.
If $(\phi^\prime,\al)$ is a solution to $Q$ then $(\phi^{\prime\prime},\al)$ 
is a solution to $Q^\prime$.
Therefore solutions of
$Q$ correspond to solutions of $Q^\prime$. The algorithm for the solution
of the $\cPG(\mbf h,\mbf n,\mbf t,\mbf p)$--problem can thus be used to
solve the $\cSPG(\mbf h^\prime,\mbf{n^\prime},\mbf t^\prime,\mbf p^\prime)$--problem.
Therefore (\ref{specialgi}) implies (\ref{specialgiii}).

Suppose now that (\ref{specialgii}) holds, so the 
$\cSPG(\mbf h^\prime,\mbf n^\prime,\mbf t^\prime,\mbf p^\prime)$--problem is
solvable for all 
all positive 4--partitions $(\mbf h^\prime,\mbf n^\prime,\mbf t^\prime,\mbf p^\prime)$ 
of 
$(h,n^\prime,t,p)$ such that $n^\prime \in \ZZ$ and  
$(\mbf h^\prime,\mbf n^\prime,\mbf t^\prime,\mbf p^\prime)\le (\mbf h,\mbf n,\mbf t,\mbf p)$.
Given 
$\mbf z\in (\dcup_{k\in K} H_{*,k}^\Lm,\cL,\mbf n,\mbf h)$, where $\cL$ is some  parameter system, the
equation
$Q=Q(\mbf z,\cL,\mbf n,\mbf t,\mbf p)$ is in standard form. From Lemma
\ref{special} it follows that there is an algorithm which outputs a
special resolution $R$ of $Q$. Furthermore the only operation involved
in producing $R$ which changes any of the $n_i$, $t_i$, $p_i$ or the length $k$ of the
partition, is reduction of
singularities. Reduction of singularities never reduces $t_i$, $p_i$, $t$ or $p$  and never 
increases $n_i$ or $k$. 
As none of the operations involved in producing $R$ increases the exponential length of any 
of the equations involved 
it follows that if 
$(\mbf q(\mbf{n^\prime},\mbf t^\prime,\mbf p^\prime)=1,\cL^\prime,\bt^\prime)$ 
is an element of $R$ with boundary labels list $(z_1^\prime,\ldots z_{n^\prime}^\prime)$
then there is a partition $\mbf h^\prime$ of $h$ such that 
$\mbf z^\prime\in (\dcup_{k\in K} H_{*,k}^\Lm,\cL,\mbf n^\prime,\mbf h^\prime)$,  
$(\mbf h^\prime,\mbf n^\prime,\mbf t^\prime,\mbf p^\prime)
\le (\mbf h,\mbf n,\mbf t,\mbf p)$ and 
$(\mbf h^\prime,\mbf n^\prime,\mbf t^\prime,\mbf p^\prime)$ is a partition of 
$(h,n^\prime,t,p)$, for
some $n^\prime \in \ZZ$. The algorithm for the  
$\cSPG(\mbf h^\prime,\mbf{n^\prime},\mbf t^\prime,\mbf p^\prime)$--problem 
may therefore be used to determine whether this
element of $R$ has a solution and if so to find one. From the
definition of resolution it now follows that an algorithm for the
$\cPG(\mbf h,\mbf n,\mbf t,\mbf p)$--problem
exists. Therefore (\ref{specialgii}) implies (\ref{specialgi}). Finally, it follows
directly from the definition of $\le_g$ that 
(\ref{specialgiii}) implies (\ref{specialgii}).\\

We define the $\cPG${\em --problem} and $\cSPG${\em --problem} in $\cC$ to be
the unions of the $\cPG(\mbf h,\mbf n,\mbf t,\mbf p)$--problems and the
$\cSPG(\mbf h,\mbf n,\mbf t,\mbf p)$--problems, respectively, over all positive 4--partitions
$(\mbf h,\mbf n,\mbf t,\mbf p)$ of length $k$, for all $k>0$. 
\begin{lemma}\label{eqn-gen}
Let $\cF$ be a collection of groups and let $\cC$ be a subcollection of
$\cF^*_\sim$.
Assume all the groups in $\cF$ have solvable word problem.
Then the
following are equivalent.
\be
\item\label{eqn-geni} The $\cPQ$--problem in $\cC$ is solvable.
\item\label{eqn-genii} The $\cPG$--problem  in $\cC$ is solvable.
\item\label{eqn-geniii} The $\cSPG$--problem  in $\cC$ is solvable.
\ee
\end{lemma}
\par\noindent{\em Proof.} That assertions (\ref{eqn-genii}) and  (\ref{eqn-geniii}) 
are equivalent
follows from Lemma \ref{specialg}. Suppose the $\cPQ$--problem is solvable. 
In the notation of the proof of Lemma \ref{specialg} Let 
$(\mbf h,\mbf n,\mbf t,\mbf p)$ be a positive 4--partition of 
$(h,n,t,p)$, $\cL$ a parameter system and $\mbf z$
an element of $(\dcup_{k\in K} H_{*,k}^\Lm,\cL,\mbf n,\mbf h)$. The algorithm for the 
$\cPQ$--problem may then be used to determine a solution to the equation 
$Q(\mbf z,\cL, \mbf n,\mbf t,\mbf p)$. 
The $\cPG$--problem is therefore
solvable.

Conversely suppose that the $\cPG$--problem is 
solvable. Using the notation of the proof of  Lemma \ref{specialg} again
let $W=(\mbf w=1,\cL,\bt)$ be 
a constrained quadratic exponential equation with 
environment 
$((X_k)_{k\in K};(H_i)_{i\in I} ;(G_k)_{k\in K})^\Lm$.
Using Theorem \ref{stdform} we may assume that
$\mbf w$ is in standard form and so $\mbf w=\mbf q(\mbf n,\mbf t,\mbf p)$, for some 
$\mbf n,\mbf t$ and $\mbf p$. Let $\mbf z$ be a boundary labels list for $(\mbf w, \bt)$
and let
$Q=Q(\mbf z,\cL,\mbf n,\mbf t,\mbf p)$. Then $W$ has a solution if and only if 
$(\mbf n,\mbf t,\mbf p)\in \cL$--genus$(\mbf z)$. Therefore, 
if $Q$ has exponential coordinates $\mbf h$,
the algorithm for the 
$\cPG(\mbf h, \mbf n,\mbf t,\mbf p)$--problem may be used to determine a solution to $W$. 
It follows
that the $\cPQ$--problem is solvable.

\subsection{Decision problems in free products and one--relator products}
\label{dec_in_op}
We shall consider  exponential genus problems for collections of groups consisting
of free--products and 
 one--relator products of groups in which such problems are already solvable.

Let $\cF$ be a collection of group presentations and define $\cF^*_\sim$ as in the 
previous section. Define the subcollection $\cC_m$ of $\cF^*_\sim$ to consist of
$\cF^*$ together with groups of the form $H/N$, where $H=\ast_{i\in I} H_i$, $H_i$ belongs to
$\cF$, 
$I$ is a recursive set of size at least $2$ 
and $N$ is the normal closure of an element $r^d$
of $H$ where $r$ is cyclically reduced, not a proper power in $H$, $l(r)\ge 2$ and
$d\ge m$. We call such groups $d${\em --torsion} one--relator products with factors in $\cF$. 
If $H/N$ is a $d$--torsion group then we may always write $H=A*B$, where $A$ and $B$ belong to 
$\cF^*$ and the conditions on $r$ still hold: that is $H/N$ is a $d$--torsion group with exactly two
factors in $\cF^*$.
 Let $s$ denote the element $r^d$ of $H$.  Then 
$l(s)=dl(r)$ is at least $2d$ and $s$ is cyclically reduced as a word in $H$.

A $d$--torsion group $A*B/N$,  
where $N$ is the normal closure of $r^d$, is said
to be {\em exceptional} of {\em type} $E(p,*,d)$ if $r=aUbU^{-1}$, with $a,b\in A\dcup B$
and $a^p=1$ and of {\em type} $E(p,q,d)$ if in addition $b^q=1$. 
A $d$--torsion group may be exceptional in more than one way, involving
different $p,q,d$: see \cite{\DHc} for more details. We define the collection of groups
$\cC_m^+$ to consist of $\cF^*$ together with those groups of $\cC_m\backslash\cF^*$ 
which are not (in any way) exceptional of type $E(2,*,m)$.
\begin{exx}\label{exx_classic_genus}
Let $\cF$ be a collection of group presentations and consider the  
$\cPG(\mbf h,\mbf n,\mbf t,\mbf p)$--problem in $\cF$, where $\mbf h$ is a partition of $0$.
Suppose \[\mbf z=(z_1,\ldots ,z_n)\in (\dcup_{i\in I}H_i^\Lm,\cL,\mbf n,\mbf h),\]
for some recursive set $I$, groups $H_i\in\cF$, and consistent parameter system $\cL$.
Since no parameters
occur in $\mbf z$, for all $j$,  we may regard $z_j$ an element
of $H_i$, for some $i\in I$. The parameter system $\cL$ has no bearing on the solutions to 
the equation associated to $\mbf z$ and we can replace it by the empty set. We write {\em genus} instead
of $\emptyset$--genus. If $\mbf n=(n_1,\ldots, n_k)$, $\mbf t=(t_1,\ldots ,t_k)$ and 
$\mbf p=(p_1,\ldots ,p_k)$ then 
\[(\mbf n, \mbf t,\mbf p)\in \textrm{genus}(\mbf z)
\textrm{ if and only if }(n_i, t_i,p_i)\in \textrm{genus}(z_i),\textrm{ for }i=1,\ldots ,k.\]
Therefore the  $\cPG(\mbf h,\mbf n,\mbf t,\mbf p)$--problem reduces to the union of the 
\[\cPG((0),(n_i),(t_i),(p_i))\textrm{--problems.}\] 
That is, we need not consider systems of equations of size 
greater than $1$ if  $\mbf h$ is a partition of zero. When $\mbf h$ is a partition of zero we refer to 
the $\cPQ$-problem as the $\cQ$--problem.
It is known (see \cite{\CEb}, \cite{\Cu}, \cite{\Ge}, \cite{\GT}) 
that the $\cQ$--problem is solvable in free groups and that 
if the $\cQ$--problem is solvable in $\cF$ then 
the $\cQ$--problem is solvable in $\cF^*$. 

Assume now that $n,t,p$ are non--negative integers such that $n+t+p>0$ and 
that the $\cPG((0),(n^\prime),(t^\prime),(p^\prime))$--problem is 
solvable in $\cF$, for all $n^\prime$, $t^\prime$, $p^\prime$ such that
\begin{gather} 
t^\prime+p^\prime/2 \le t+p/2 \textrm{ and }\nonumber\\ 
n^\prime \le n+2(t-t^\prime)+(p-p^\prime).\label{genus_cond}
\end{gather}
Now let $m$ be a fixed integer and let 
$\cC$ be a collection of 
one--relator products with factors
in $\cF$ such that every element of $\cC$ is $d$--torsion, for some $d\ge m$. It is known 
(\cite{\DHa}, \cite{\DHb}, \cite{\DHc}) that if (\ref{genus_cond}) and 
one of the following conditions (\ref{genus_cond_1})--(\ref{genus_cond_3}) holds, 
for all groups $G$ of $\cC$, then the 
$\cPG((0),(n),(t),(p))$--problem is 
solvable in $\cC$. 
\be
\item\label{genus_cond_1} $m\ge 5$ and $G$ is not in any way exceptional of type $E(2,3,5)$ or $E(2,3,6)$.
\item $m\ge 4$ and no letter occuring in $r$ has order $2$.
\item\label{genus_cond_3} $m\ge 2$ and the groups of $\cF$ are locally indicable.
\ee 
\end{exx}

The previous example motivates the following definition.
Let $(\mbf h,\mbf n, \mbf t,\mbf p)$ 
be a positive 4--partition of length $k$ of $(h,n,t,p)$, for some $h,n,t,p\in \ZZ$. 
Define $\cN(\mbf h,\mbf n,\mbf t,\mbf p)$ to 
be the set of positive 4--partitions 
$(\mbf{h^\prime},\mbf{n^\prime},\mbf{p^\prime},\mbf{t^\prime})$, of length $k^\prime$,
satisfying the conditions (\ref{N1})--(\ref{N4}) below. Assume
$\mbf h^\prime=(h_1^\prime,\ldots,h_{k^\prime}^\prime)$, 
$\mbf n^\prime=(n_1^\prime,\ldots,n_{k^\prime}^\prime)$, 
$\mbf t^\prime=(t_1^\prime,\ldots,t_{k^\prime}^\prime)$ and  
$\mbf p^\prime=(p_1^\prime,\ldots,p_{k^\prime}^\prime)$.
\be
\item\label{N1} There exists a partition $\mbf l=(l_1,\ldots,l_k)$ such that
$0\le l_j\le\max\{1,h_j\}$, for $j=1,\ldots, k$ and 
\[k^\prime=\sum_{j=1}^k l_j.\] 
\item Setting $q_0=0$ and $q_j=\sum_{s=1}^j l_s$, for $j=1,\ldots ,k$, we have
\[\sum_{i=q_{j-1}+1}^{q_j} h^\prime_i\le h_j, \textrm{ for } j=1,\ldots ,k.\]
\item\label{N4} For $q_{j-1}+1\le i\le q_j$ and $1\le j\le k$
\be
\item $t^\prime_i+p^\prime_i/2 \le t_j+p_j/2$ and 
\item $n^\prime_i \le n_j+2(t_j-t^\prime_i)+(p_j-p_i^\prime)$.
\ee
\ee

The main result proved here is the following.
\begin{theorem}
\label{main}
Let $(\mbf h,\mbf n,\mbf t,\mbf p)$ be a positive $4$--partition and let $\cF$ be
a collection of group presentations with solvable 
$\cPG(\mbf h^\prime,\mbf n^\prime,\mbf t^\prime,\mbf p^\prime)$--problem, for 
all 
$(\mbf h^\prime,\mbf n^\prime,\mbf t^\prime,\mbf p^\prime)\in 
\cN(\mbf h,\mbf n,\mbf t,\mbf p)$. Then the 
$\cPG(\mbf h^\prime,\mbf n^\prime,\mbf t^\prime,\mbf p^\prime)$--problem is solvable
in the subcollection $\cC^+_6$ of $\cF^*_\sim$ , for all 
$(\mbf h^\prime,\mbf n^\prime,\mbf t^\prime,\mbf p^\prime)\in 
\cN(\mbf h,\mbf n,\mbf t,\mbf p)$.
\end{theorem}
\medskip

\noindent
Most of the remainder of the paper is given to proving a restricted version
of this theorem. Here we state the restricted version 
and show how Theorem \ref{main}
follows from it.
\begin{theorem}
\label{mainc}
Let $(\mbf h,\mbf n,\mbf t,\mbf p)$ be a positive $4$--partition, let $I$ be
a recursive set  and let 
$A, B$, $H_i, i\in I$, groups, with fixed presentations given, such that the 
$\cPG(\mbf h^\prime,\mbf n^\prime,\mbf t^\prime,\mbf p^\prime)$--problem 
is solvable in the indexed set of groups
$(A, B)\cup(H_i:i\in I)$, 
for 
all 
$(\mbf h^\prime,\mbf n^\prime,\mbf t^\prime,\mbf p^\prime)\in 
\cN(\mbf h,\mbf n,\mbf t,\mbf p)$.
Let $H=A*B$ and  let $s=r^m$, 
where $r\in H$ is cyclically reduced, not a proper power in $H$, $l(r)\ge 2$ and
$m\ge 6$.
Set $G=H$ or $G=H/N$, where $N$ is the normal closure of $s$ in $H$. In the 
latter case assume
further that $G$ is not in any way exceptional of type $E(2,*,m)$.
Then the 
$\cPG(\mbf h^\prime,\mbf n^\prime,\mbf t^\prime,\mbf p^\prime)$--problem is solvable
in $(G, A, B)\cup( H_i:i\in  I )$, for 
all 
$(\mbf h^\prime,\mbf n^\prime,\mbf t^\prime,\mbf p^\prime)\in 
\cN(\mbf h,\mbf n,\mbf t,\mbf p)$.
\end{theorem}
\medskip

\noindent
{\em Proof of Theorem \ref{main}.} 
Assume that Theorem \ref{mainc} holds and set
$\cN=\cN(\mbf h,\mbf n,\mbf t,\mbf p)$. Suppose that  $\cC$ is a subcollection
of the  collection $\cF^*_\sim$ for which it has been shown 
that if   the hypotheses of Theorem \ref{main} hold then the 
$\cPG(\mbf h,\mbf n,\mbf t,\mbf p)$--problem is solvable in $\cC$. Given
any $(\mbf h^\prime,\mbf n^\prime,\mbf t^\prime,\mbf p^\prime)\in 
\cN$ we have $\cN(\mbf h^\prime,\mbf n^\prime,\mbf t^\prime,\mbf p^\prime)
\subseteq \cN$, so the hypotheses of Theorem \ref{main} hold for 
$(\mbf h^\prime,\mbf n^\prime,\mbf t^\prime,\mbf p^\prime)$ instead of 
$(\mbf h,\mbf n,\mbf t,\mbf p)$. Hence the 
$\cPG(\mbf h^\prime,\mbf n^\prime,\mbf t^\prime,\mbf p^\prime)$--problem is solvable
is solvable in $\cC$. In particular it suffices to show that the 
$\cPG(\mbf h,\mbf n,\mbf t,\mbf p)$--problem is solvable in $\cC^+_6$.

Assume then that the hypotheses of Theorem \ref{main} hold.
First we shall prove the following. Let $l$ be a positive integer and 
$(A_1,\ldots ,A_l,B_1,\ldots ,B_l) \cup(
H_i:i\in  I )$ be an indexed set of groups all of which belong to $\cF$. For
$i=1,\ldots ,l$ let $H_{*,i}=A_i * B_i$ and let $N_i$ be a subgroup of $H_{*,i}$ 
satisfying either
\be[{$\cO$}1.]
\item\label{Gfree} $N_i=1\in H_{*,i}$ or 
\item\label{Gone} $N_i$ is the normal closure of $r_i^{m_i}$ in $H_{*,i}$, 
where $r_i\in H_{*,i}$ and $m_i\ge 6$ is an integer such that 
$H_{*,i}/N_i$ is an $m_i$--torsion one--relator product which is not exceptional of 
type $E(2,*,m)$.
\ee
Set $G_i=H_{*,i}/N_i$, $i=1,\ldots ,l$.
Then the $\cPG(\mbf h^\prime,\mbf n^\prime,\mbf t^\prime,\mbf p^\prime)$--problem 
is solvable in \[(G_1,\ldots ,G_l,A_1,\ldots ,A_l,B_1,\ldots ,
B_l)\cup(H_i:i\in  I ),\] for all $(\mbf h^\prime,\mbf n^\prime,\mbf t^\prime,\mbf p^\prime)
\in \cN$.  

The case $l=1$ follows from Theorem \ref{mainc}. Assume inductively that
the result holds for some $l-1\ge 1$. Then the 
$\cPG(\mbf h^\prime,\mbf n^\prime,\mbf t^\prime,\mbf p^\prime)$--problem is solvable
in \[(G_1,\ldots ,G_{l-1},A_1,\ldots ,A_l,B_1,\ldots ,B_l)\cup(H_i:i\in  I ),\] 
for all $(\mbf h^\prime,\mbf n^\prime,\mbf t^\prime,\mbf p^\prime)
\in \cN$. Using Theorem \ref{mainc} again the 
$\cPG(\mbf h^\prime,\mbf n^\prime,\mbf t^\prime,\mbf p^\prime)$--problem is solvable
in \[(G_1,\ldots ,G_{l},A_1,\ldots ,A_l,B_1,\ldots ,B_l)\cup(H_i:i\in  I ),\] as required. 

Now let $l_1,\ldots ,l_n$ be integers, $l_i\ge 2$ and let $A_{i,1},\ldots,A_{i,l_i}$
belong to $\cF$, for $i=1,\ldots ,n$. Define $H_{*,i}=*\{A_{i,j}:j=1,\ldots ,l_i\}$ 
and let $N_i$ 
be a subgroup of $H_{*,i}$ satisfying $\cO$\ref{Gfree} or $\cO$\ref{Gone}, 
for $i=1,\ldots n$. We shall show that, if $G_i=H_{*,i}/N_i$,  then the 
$\cPG(\mbf h^\prime,\mbf n^\prime,\mbf t^\prime,\mbf p^\prime)$--problem is solvable
in \[(G_1,\ldots ,G_{l},A_{1,1},\ldots ,A_{n,l_n})\cup( B_i:i\in  I ),\] for any recursive
set I and groups $B_i$ of $\cF$, and for all 
$(\mbf h^\prime,\mbf n^\prime,\mbf t^\prime,\mbf p^\prime)
\in \cN$. 

The case $l_i=2$ for all $i$ is covered by the argument above. Using induction
on $\sum_{i=1}^n(l_i-2)$ we assume that $l_i>2$, for some $i$. Without loss of 
generality we may assume that $i=1$ and set 
\[H_{*,1}^\prime=A_{1,2}*\cdots *A_{1,l_1}.\] 
The inductive hypothesis implies that the 
$\cPG(\mbf h^\prime,\mbf n^\prime,\mbf t^\prime,\mbf p^\prime)$--problem is solvable
in \[(H_{*,1}^\prime,G_2,\ldots ,G_{n},A_{1,1},\ldots ,A_{n,l_n})\cup( B_i:i\in  I ),\]
for all 
$(\mbf h^\prime,\mbf n^\prime,\mbf t^\prime,\mbf p^\prime)
\in \cN$.
As $G_1=(A_{1,1}*H_{*,1}^\prime)/N_1$, it follows from Theorem \ref{mainc} that
the 
$\cPG(\mbf h^\prime,\mbf n^\prime,\mbf t^\prime,\mbf p^\prime)$--problem is solvable
in \[(G_1,H_{*,1}^\prime,G_2,\ldots ,G_{n},A_{1,1},\ldots ,A_{n,l_n})\cup( B_i:i\in  I )\]
and so in particular in 
\[(G_1,G_2,\ldots ,G_{n},A_{1,1},\ldots ,A_{n,l_n})\cup( B_i:i\in  I ),\] 
for all 
$(\mbf h^\prime,\mbf n^\prime,\mbf t^\prime,\mbf p^\prime)
\in \cN$,
as required.

Now let 
$R$ be a recursive set and let $\cI$ be a set of groups, indexed by $R$, all of which belong to
$\cC^+_6$. Then we can write
$R =J\cup K$, where $J$ and $K$ are recursive sets,  and there exists an element $B_j$ of $\cF$, for all
$j\in J$,  and  $G_k$ of  $\cC^+_6$ but not to $\cF$, for all $k\in K$, such that 
$\cI=(B_j)_{j\in J}\cup {G_k}_{k\in K}$. Thus there exist
recursive sets $I$, $X_k\subseteq I$ 
and groups $H_i$, $H_{*,k}$ and $N_k$, such that $H_{*,k}=*_{i\in X_k}H_i$, with $k\in K$, 
$H_i\in \cF$,  $G_k=H_{*,k}/N_k$, where  $H_{*,k}$, $N_k$ and $G_k$ satisfy
$\cO$\ref{Gfree} and $\cO$\ref{Gone}.
Suppose that
\[\mbf z\in \left(\bdcup_{k\in K}H_{*,k}^\Lm\dcup \bdcup_{j\in J}B^\Lm_j,\cL,\mbf n,\mbf h\right).\]
The equation $Q$ associated to $\mbf z$ has environment 
\begin{equation}\label{big_env_1}
((X_k)_{k\in K},(\{j\})_{j\in J};(H_i)_{i\in I},(B_j)_{j\in J};(G_k)_{k\in K},(B_j)_{j\in J})^\Lm
\end{equation}
and 
we may assume, after reordering if necessary, that $\mbf z$ and $Q$ have basis
$(k_1,\ldots, k_l,j_1,\ldots ,j_{l^\prime})$, where $k_i\in K$ and $j_i\in J$.
From Lemma \ref{basis_reduce} a solution to $Q$ with environment (\ref{big_env_1}) determines a solution
to $Q$ with environment 
\begin{align}\label{small_env_1}
&(X^\prime_{k_1},\ldots, X^\prime_{k_l},\{j_1\}\ldots,\{j_{l^\prime}\};\nonumber\\
&(H_i)_{i\in I^\prime},B_{j_1},\ldots,B_{j_{l^\prime}};\nonumber\\
&G^\prime_{k_1},\ldots,G^\prime_{k_l},B_{j_1},\ldots,B_{j_{l^\prime}})^\Lm
\end{align}
and the same basis,
where supp$(Q)=(X^\prime_{k_1},\ldots, X^\prime_{k_l},\{j_1\}\ldots,\{j_{l^\prime}\})$,  
$I^\prime=X^\prime_{k_1}\cup \cdots \cup X^\prime_{k_l}$ and $G^\prime_k$ is defined in (\ref{small_env}). 
Conversely 
a solution to $Q$ with environment (\ref{small_env_1}) is a solution to $Q$ with environment
(\ref{big_env_1}).
By construction $G^\prime_{k_i}=H^\prime_{*,k_i}/N_i^\prime$, where $H^\prime_{*,i}$, $N^\prime_i$ and 
$G^\prime_i$, satisfy $\cO$\ref{Gfree} and $\cO$\ref{Gone}. 
Hence $(\mbf n,\mbf t,\mbf p)\in \cL$--genus$(\mbf z)$ if  and only if there exists a solution to $Q$ with
environment (\ref{small_env_1}). We may determine whether or not $Q$ has a solution  with this environment if
we
can solve the $\cPG(\mbf h,\mbf n,\mbf t,\mbf p)$--problem 
in 
\begin{equation}\label{last_env}
(G^\prime_{k_1},\ldots,G^\prime_{k_l})\cup(H_i:i\in  I^\prime)\cup(B_{j_1},\ldots,B_{j_{l^\prime}}).
\end{equation}
Since $N^\prime_{k_s}$ is either trivial or the normal closure of a single element, it follows that $X^\prime_{k_s}$
is finite, for $s=1,\ldots ,l$. Hence, from the previous case, 
the $\cPG(\mbf h,\mbf n,\mbf t,\mbf p)$--problem is solvable in
(\ref{last_env})
 and the Theorem follows.
\begin{corol}\label{main_eqn}
Let $(\mbf h,\mbf n,\mbf t,\mbf p)$ be a positive $4$--partition and let $\cF$ be
a collection of group presentations with solvable 
$\cPQ$--problem.
Then the 
$\cPQ$--problem is solvable
in the subcollection $\cC^+_6$ of $\cF^*_\sim$.
\end{corol}
\medskip

\noindent
{\em Proof.} This follows from Theorem \ref{main} and Lemma \ref{eqn-gen}.
\begin{corol}\label{cyclic_corol}
Let $\cF$ be the cyclic groups. 
The $\cPQ$--problem is solvable in $\cC_6^+$. That is in free--products 
and $d$--torsion one--relator products which belong to $\cC^+_6$, with cyclic factors.
In particular the $\cPQ$--problem is solvable in  recursively generated free groups. 
\end{corol}

\noindent{\em Proof. } The $\cPQ$--problem is solvable in $\cF$ since each instance reduces to a system
of linear Diophantine equations, congruences and inequalities. The result follows from Corollary \ref{main_eqn}.

\begin{corol}\label{plain}
Let $A$ and $B$ be groups such that the $\cPG$--problem is solvable in $(A,B)$.
Let $G=A*B/N$ be a $d$--torsion one--relator product of $A$ and $B$, for $d\ge 6$, which is not in
any way special of type $E(2,*,d)$. Then the $\cPG$--problem is solvable in $(G,A,B)$. In particular,
if $z_1,\ldots ,z_k$ are elements of $A*B$ and $\cL$ is a consistent parameter system in 
$\{\lm_1,\ldots \}$ then there is an algorithm to decide whether or not there are integers $n_1,\ldots,n_k$
and elements $x_i\in A*B$ such that
\be
\item
\[\prod_{i=1}^{k}x_i^{-1}z_i^{n_i}x_i\prod_{i=k+1}^{n+t}[x_i,x_{i+t}]\prod_{i=k+2t+1}^{k+2t+p}x_i^2=1,
\textrm{ in $G$ and }
\]
\item there  is a solution to $\cL$ with $\lm_i=n_i$, $i=1,\ldots ,k.$ 
\ee
\end{corol}

%% file: in2.tex
\section{Pictures }\label{pictures}

Let $A$ and $B$ be groups and $m\geq 2$ an integer and set $H=A*B$.  
 Let $X^\prime_A$ 
and $X_B^\prime$ be CW--complexes with $\pi_1(X_A^\prime,x_A)=A$ and
$\pi_1(X_B^\prime,x_B)=B$, for some $x_A\in X_A^\prime$ and $x_B\in X_B^\prime$. 
Let $O=A$ or  $O=B$ and  let $e_O$
be a copy of the 
unit interval $[0,1]$.
Form the quotient space $X_O$ of $X^\prime_O\dcup e_O$ 
by identifying $x_O$ and $1\in e_O$.  
Form the quotient space $Y$ of $X_A\dcup X_B$ 
by identifying $0\in e_A$ and $0\in e_B$ to a point $p$.
Then $\pi_1(Y,p)=H$. 

Fix a word $r\in H$ with 
$l (r)\geq 2$, such that $r$ is not a proper power in $H$.  Let
$s=r^m$, $N(s)=\langle\langle s\rangle\rangle^H$  and 
$G=H/N(s).$ Suppose that $s=a_1b_1\cdots a_nb_n$, when written in normal form for the free product.
Let $\D_s$ be a closed unit disk with oriented boundary 
divided into $2n$ intervals. Label the intervals consecutively $a_1,b_1,\ldots, a_n,b_n$,
so that $\pd\D_s$, read from a suitable starting point in the direction of orientation, is
labelled $s$. Form the quotient space
$Z$ of $Y\dcup \D_s$ by identifying $\pd\D_s$ with a path in $(Y,p)$ representing 
$s\in H$, according to the label of $\pd\D_s$. Then $\pi_1(Z,p)=G$.
Alternatively, let $G=H$ and define $\D_s=p\in Y$, so that $Z=Y$ and $\pi_1(Z,p)=G$ again.

Let $\S$ be a compact surface and let $f:\S\maps Z$ be a continuous map satisfying the
following properties.
\be[$S$1.]
\item \label{sketch1} $f$ is transverse to the midpoint $c$ of $\D_s$ and 
the closure of $f^{-1}(\inr(\D_s))$ is a regular neighbourhood of $f^{-1}(c)$.
\item \label{sketch2} $f$ is transverse to $p$ on 
$\S_0=\S\bl f^{-1}(\inr(\D_s))$ and $f^{-1}(e_A\cup e_B)$ is  a
regular neighbourhood of $f^{-1}(p)$ in $\S_0$. 
\ee

Let $\G= f^{-1}(\D_s)$. Then $(f,\G)$ is a {\em sketch} on $\S$ over $G$ (with respect
to $Z$) which we refer to merely as $\G$, when the meaning is clear.  A connected 
component of the closure of $f^{-1}(\inr(\D_s))$ is called a {\em vertex}
of $\G$ and the set of all vertices is denoted $\cV(\G)$ or just $\cV$.
Condition
$S$\ref{sketch1} implies that $\cV$ 
consists of a collection of 
embedded disks in $\inr(\S)$, each of which is a regular neighbourhood of a point
of $f^{-1}(c)$.
If $G=H$ then $\cV=\emptyset$. 
A connected component of $f^{-1}(p)$ is called an {\em arc} of $\G$ and the 
set of all arcs is denoted $\cA=\cA(\G)$. Condition $S$\ref{sketch2} implies
that $\cA$ consists of a set of properly embedded $1$--submanifolds of $\S_0$.

A connected component of $\pd v\bl \{\cup e: e\in \cA\}$, where $v\in \cV$,  is called a {\em vertex corner} 
of $\G$. 
A component of $\pd\S\bl \G$ is called a {\em boundary corner} of $\G$. 
The closure of 
a connected component of $\S\bl \G$ is called  a {\em region} of $\G$ and the 
set of regions is denoted $\cR(\G)$ or $\cR$.
Every vertex and boundary corner belongs to a unique region of $\G$.

Now suppose that $\S^\prime$ is a connected component of $\S$ 
and define $\G^\prime=\S^\prime \cap \G$. If $\S^\prime$ is 
orientable we may choose consistent 
orientations of each region and vertex of $\S^\prime$. We call such a
choice an {\em orientation} of $\G^\prime$. If $\bt$ is a boundary component
of $\S^\prime$ then an orientation $\z$ of $\G^\prime$ gives
rise to a unique orientation $\z(\bt)$ of $\bt$. Similarly the orientation of
regions induces an orientation on $\pd v$, which we denote
$\z(\pd v)$, for all $v\in \cV$. Thus $\z(\pd v)$ is the opposite of
the orientation induced on $\pd v$ by the orientation of $v$.  
If $l$ is a sub--interval of a corner $c$ of $\G$ then there is a unique
orientation $\z(l)$ of $l$ induced from the orientation of the
boundary component (of $\S^\prime$ or $v\in \cV$) to which $c$ belongs.
 
If $\S^\prime$ is non--orientable then we  choose an 
orientation 
for each connected component of $\pd\S^\prime$, 
for each connected component of $\pd \D$, for all regions
$\D$ of $\G^\prime$, 
and for  $\pd v$, for all vertices $v$ of $\G^\prime$.
To each vertex $v$ we assign the orientation opposite to that chosen for 
$\pd v$.
Call the collection of these orientations $\z$
an {\em orientation} of $\G^\prime$.
Let $l$ be a sub--interval of a corner $c$ of $\G^\prime$. Then 
$c$ belongs to a boundary component $\bt$ of $\S^\prime$ or to the boundary $\bt$ of 
a vertex. The orientation $\zeta(\bt)$  induces an 
orientation of $l$ which we call $\zeta(l)$.

We now fix an orientation for $\G^\prime$ for all connected
components $\S^\prime$  of $\S$. We call these orientations an {\em orientation}
of $\G$ and denote this by $\z$. We refer
to the orientation induced on $y$ by $\z$ as $\z(y)$, whenever $y$ is an appropriate
subset of $\S$. 

Let $(f,\G)$ be a 
sketch on $\S$ 
with orientation $\z$. We define a labelling function $f^l$ from corners of $\G$ to $H$
as follows. Let $c$ be a corner of $\G$  with
orientation $\z(c)$ and let 
$q$ be the closed  path $f(c)$ in $Y$. 
If $c$ meets an arc of $\G$ then $q$ is a directed path based at $p$.
In this case we 
define the {\em label} $f^l(c)$ {\em of} $c$ {\em read in the direction} $\zeta(c)$ 
to be the element of $H$ 
represented by the based directed path $q$.
If $c$ does not meet an arc of $\G$ then we choose a base point
$x$ on $c$ and regard $q$ as a directed path based at $f(x)$.
Since $c$ meets no arc of $\G$, the path $q$ lies in $Y\backslash p$ and we may
choose a path $\al$ from $p$ to $x$ such that int$(\al)$ lies in $X_A$ or $X_B$. 
Then the directed path $\al q\al^{-1}$ represents an element $g$ of $A\dcup B$. 
The label of $c$, read in the direction of orientation $\z(c)$, is defined to 
be $g$.  
If $c$ read in the direction $\zeta(c)$ 
has label $w$ then 
the label of $c$ {\em read in the direction opposite to} $\zeta(c)$ is $w^{-1}$.
\begin{defn}\label{pic_def}
The sketch $\G$ endowed with the orientation $\z$ and the labelling function $f^l$ is
called a {\em picture} on $\S$ over $G$.
\end{defn}

\begin{figure}
\psfrag{a}{$a$}
\psfrag{b}{$b$}
\psfrag{c}{$c$}
\psfrag{d}{$d$}
\psfrag{a^-1}{$a^{-1}$}
\psfrag{b^-1}{$b^{-1}$}
\psfrag{c^-1}{$c^{-1}$}
\psfrag{d^-1}{$d^{-1}$}
\begin{center}
  \mbox{
\subfigure[Positive\label{vertex-pos}]
{ \includegraphics[scale= 0.4,clip]{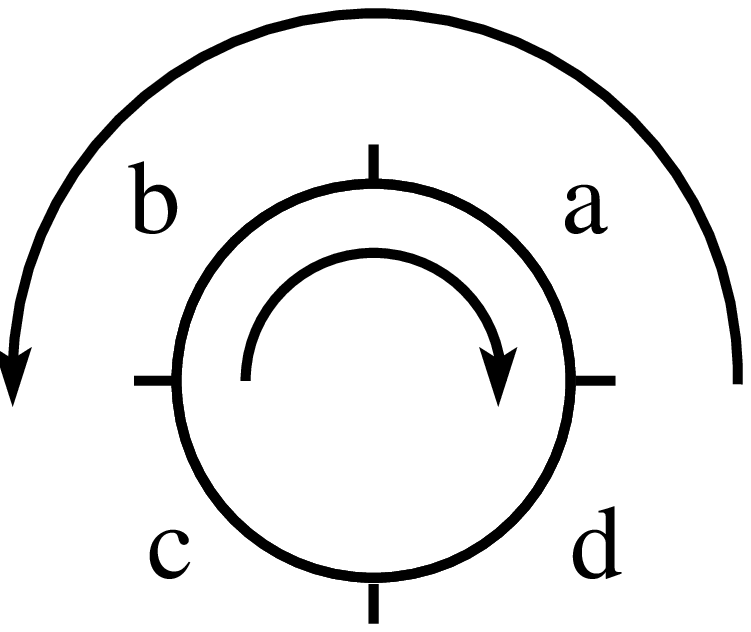} }\qquad
\subfigure[Negative\label{vertex-neg} ]  
{ \includegraphics[scale=0.4,clip]{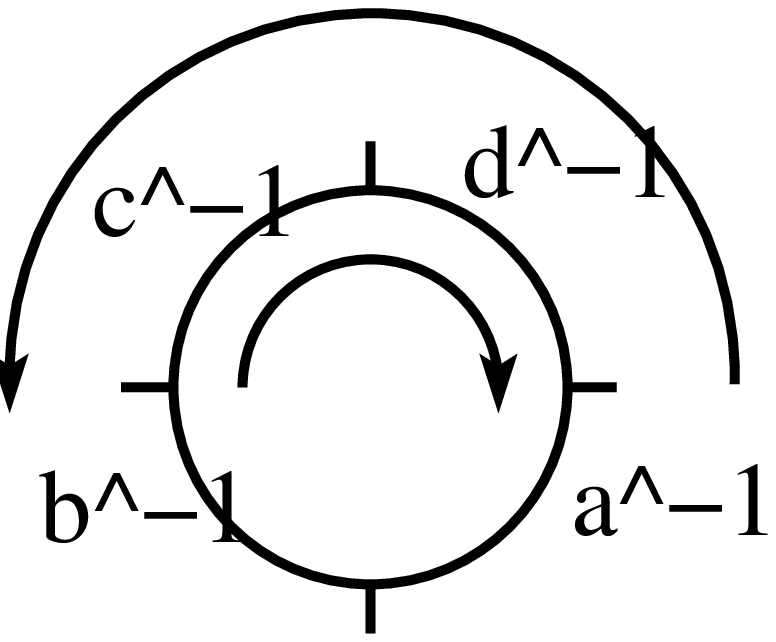} } 
}
\caption{The orientation of the boundary is shown outside and the orientation of the 
vertex inside: $s=abcd$.}\label{vertex-pos-neg}
\end{center}
\end{figure}

A picture over $G=H$ may be regarded as a picture over $G=H/N(s)$, having no vertices. Conversely
a picture over $G=H/N(s)$ with no vertices may be regarded as a picture over $H$.
If $v$ is a vertex of $\G$ and $x$ is a point of 
$\pd v\cap \{\cup e:e\in \cA\}$ then the {\em label} of $v$ with {\em base point} $x$
is the word of $H$ obtained
by reading off the labels of corners of $v$ in the order they occur whilst travelling
once round $\pd v$, from $x$ in the direction $\z(\pd v)$: where the label of 
each corner $c$ 
is read in the direction of orientation induced from that of  
$\zeta(\pd v)$, namely $\zeta(c)$.
Since $f(\pd v)=\pd \D_s$ we have:
\be[$P$1.]
\item \label{pic_cond0} the label of $v$ (with any base point) 
is a cyclic permutation of $s^{\e}$, where
$\e=\pm 1$. 
\ee
The orientation $\z$ on $\pd v$ induces an orientation on 
$f(\pd v)=\pd \D_s$. If this is the same as the fixed orientation of $\pd \D_s$ 
then we may take $\e = 1$ and we
say that $v$ is a {\em positive} vertex (with respect to $\z$). Otherwise we take $\e=-1$ and
say $v$  is a {\em negative} vertex (with respect to $\z$), as shown in Figure \ref{vertex-pos-neg}. 
Similarly, if $\bt$ is a connected component of
$\pd\S$ or a connected component $\pd\D$, for some region
$\D$ of $\G$, let $x\in \bt\cap \{\cup e:e\in \cA(\G)\}$ or, if $\bt$ consists of 
a boundary corner meeting no arc of $\G$, let $x$ be the chosen base point of this
corner. Then
the {\em label} of $\bt$ with {\em base point} $x$ is the word obtained by reading    
once round $\bt$, from $x$ in the direction of $\z(\bt)$. When reference
is made to the label of such $\bt$ it is to be assumed that a base
point has been chosen. The {\em boundary labels} of
a region $\D$ are the labels of the connected components of $\pd\D$.

 As a consequence of properties $S$\ref{sketch1} and $S$\ref{sketch2} 
and the definitions above a 
picture $\G$ on $\S$ satisfies the following properties 
$P$\ref{pic_cond1}--$P$\ref{pic_cond3}. 
\be[$P$1.]\label{piccond}
\addtocounter{enumi}{1}
\item \label{pic_cond1} Let $\D$ be a region of $\G$ with boundary labels 
$d_1,\ldots, d_n$. Then either $d_i\in A$, for $i=1,\ldots , n$, in which case
$\D$ is called an $A${\em --region} or
$d_i\in B$, for $i=1,\ldots , n$, in which case
$\D$ is called a $B${\em --region}.
\item  \label{pic_cond2} 
Let $O=A$ or $B$ and let $\D$ be an $O$--region, with boundary labels
$d_1,\ldots , d_n$, which is homeomorphic to the connected sum of an $n$-punctured
sphere with $t$ torii and $p$ projective planes. Then there exist elements 
$w_1,\ldots , w_n, x_1,\ldots, x_t,y_1,\ldots, y_t,z_1,\ldots,z_p\in O$ such that
\begin{equation}\label{region_cond}
\prod_{i=1}^{n} w_i^{-1}d_iw_i\prod_{i=1}^t[x_i,y_i]\prod_{i=1}^p z_i^2=1,
\end{equation}
in $O$.
\item  \label{pic_cond3} 
If $a$ is an arc of $\G$ then $a$ separates an $A$--region from a $B$--region.
\ee

Conversely properties $P$\ref{pic_cond0} to $P$\ref{pic_cond3} characterise pictures.
\begin{lemma}\label{pic_comb}
Let $\S$ be a compact surface and let $\G$ consist of
\be
\item   a finite collection $\cV$ of disjoint
closed disks in int$(\S)$ (called vertices);
\item a 
finite collection  $\cA$ of disjoint connected
properly embedded $1$--submanifolds of 
$\S_0=\S\backslash\cup_{v\in \cV}\textrm{int}(v)$ (called arcs);
\item 
regions and corners and 
an orientation $\z$ defined as
in the case of the sketch above;  
\item a labelling function $f^l$
from corners of $\G$ to $A\cup B$, which respects orientation, and a collection
of labels defined
as for pictures above.
\ee
If $\G$ satisfies $P$\ref{pic_cond0} to $P$\ref{pic_cond3} then $\G$  is a picture on $\S$
over $G$.
\end{lemma}
We define the {\em genus} of a compact connected surface 
$\S$, with $k$ boundary components,  to be 
genus$(\S)=(2-\chi(\S)-k)/2$. The {\em genus} of a disconnected compact surface
is defined to be the sum of the genera of its connected components.

An arc of a picture, which is not a closed curve, both ends of which
meet $\pd \S$ is called an arc of {\em type} $I$. 
Let $\D$ be a region of $\G$. Then we define
\begin{align*}
t(\D)&= \mbox{ number of arcs of } \G \mbox{ in }\pd\D,\\
\bt(\D)&= \mbox{ number of components of }[\pd\S\backslash \{\mbox{arcs of }\G\}]\cap\pd\D\\
&=\textrm{ number of boundary corners of } \D,\\
\rho(\D)&= 
\mbox{ number of components of }[\G\backslash \{ \mbox{arcs of } \G\}]\cap\pd\D\\
&=\textrm{ the number of vertex corners of }\D,\\
\g(\D)&= \mbox{ number of components } \pd\D \textrm{ meeting no 
arc of }\G\\
\e(\D)&= \textrm{ number of arcs of $\G$ in }\pd\D \mbox{ meeting no vertex of }\G\\
&= \textrm{the number of arcs which are closed curves or of type } I \textrm{ in }\pd\D.
\end{align*}
If $\x(\D )=1$ and $t(\D )=2$ then we call $\D$ {\em collapsible}.

\subsection{Boundary intervals}
Let $\bt$ be a connected component of $\pd\S$, with orientation $\z$, and let
$b_1,\ldots ,b_n$ be $n$ points occuring consecutively around
$\bt$ (read in the direction of $\z(\bt)$). We call $(b_1,\ldots ,b_n)$ a 
{\em boundary partition} of $\bt$. The closed intervals $[b_{i},b_{i+1}]$ of $\bt$,
for $i=1,\ldots ,n$ (subscripts modulo $n$) are called {\em boundary intervals} of 
$\bt$. The boundary interval $[b_{i},b_{i+1}]$ has {\em length} $|[b_{i},b_{i+1}]|$ equal
to the number of connected components of $[b_{i},b_{i+1}]\bl \G$.
(The $b_i$'s do not have to belong to $\G\cap \bt$.)

Let $\cL$ be a consistent system of parameters
and $z=\prod_{i=1}^n(h_i,m_i)\in F(H^\Lm)$. The boundary partition $(b_1,\ldots ,b_n)$ is said to
$\cL${\em --admit} $z$ if (with subscripts modulo $n$)
\be
\item  $|[b_{i},b_{i+1}]|=1$, for all $i$ such that $l(h_i)=1$,
\item $|[b_{i},b_{i+1}]|=l(h_i\^m_i)$, for all $i$ such that $(h_i,m_i)$ is degenerate and
\item $|[b_{i},b_{i+1}]|=q_il(h_i)+k_i$, for some $k_i,q_i\in\ZZ$ such that  $q_i\ge 0$ and $\cL$ implies that
$m_i\equiv k_i (\mod l(h_i))$, for all $i$ such that $l(h_i)>1$ and $(h_i,m_i)$ is non--degenerate.
\ee
If $\al$ is a solution to $\cL$, the boundary partition  $(b_1,\ldots ,b_n)$ $\cL$--admits $z$ 
and $[b_{i},b_{i+1}]$ maps to a path in $Y$ representing
$\hat\al(h_i,m_i)$, 
 then we say that 
 $[b_{i},b_{i+1}]$ has {\em label} $\hat\al(h_i,m_i)$ and 
{\em prime label} $(h_i,m_i)$, 
for $i=1,\ldots ,n$.
In this case the label of $\bt$ is $\hat\al(z)$ and we say 
that $\bt$ has {\em prime label} $z$ with {\em partition}
$(b_1,\ldots ,b_n)$.

Let $(\mbf n,\mbf t,\mbf p)$
be a positive 3--partition of $(n,t,p)$ of 
length $k$, with $\mbf n=(n_1,\ldots, n_k)$ etc.. 
We shall say that a surface has {\em type} 
$(\mbf n,\mbf t,\mbf p)$
if it is a disjoint union of compact 
surfaces $\S_i$, $i=1,\ldots ,k$, such that $\S_i$ has
$n_{i}$ boundary components, genus $t_{i}+p_{i}/2$ and is orientable if $p_{i}=0$.
In this case, if $\G$ is a picture on $\S_1\dcup \cdots \dcup \S_k$ we say that $\G$ has 
{\em partitioned boundary labels list} 
$(u_1,\ldots,u_n)$ 
if 
the boundary labels of $\G\cap \S_i$ are $u_{\mu_i+1},\ldots ,u_{\mu_i+n_{i}}$, for
$i=1,\ldots , k$, where
$\mu_i$ is defined as above.

Now with $H=A*B$ and $G=H/N(s)$ or $G=H$ as above, let $(\mbf h,\mbf n,\mbf t,\mbf p)$
be a positive 4--partition, 
$\mbf z$  a special
element of $(H^\Lm,\cL,\mbf n,\mbf h)$ and $\al$ a solution  to $\cL$. Suppose that $\G$ is a picture 
on a surface 
of type $(\mbf n,\mbf t,\mbf p)$ with partitioned boundary labels list
$(\hat\al(z_1),\ldots,\hat\al(z_n))$.  
Then we may choose a boundary partition $\mbf b_i$ for the boundary component $\bt_i$
labelled $\hat\al(z_i)$  
such that $\mbf b_i$ $\cL$--admits $z_i$.
Then $\bt_i$ has prime label $z_i$ with partition $\mbf b_i$, for all $i$.
\begin{defn}\label{boundary_partition}
Given such boundary partitions we set $\mbf b= \mbf b_1,\ldots,\mbf b_n$ and say that 
$\G$ has {\em boundary partition} $\mbf b$ and {\em prime labels} $\mbf z$  with
{\em partition} $\mbf b$. For ease of notation we write   $\hat\al(\mbf z)$ for
$(\hat\al(z_1),\ldots,\hat\al(z_n))$.
\end{defn}

\begin{defn}
\label{h-vertex}
Let $a$ be a  boundary arc and let $x$ be an endpoint of $a$ on $\pd\S$.
Then $x$ meets 2 distinct boundary corners: denote the closures of these corners
by $c_1$ and $c_2$. If int$(c_1\cup c_2)$ contains a point of a boundary partition
then $a$ is called an $H^\Lm${\em --arc}. 
A vertex incident to an $H^\Lm$--arc is
called an $H^\Lm${\em --vertex}. The sets of $H^\Lm$-arcs and $H^\Lm$--vertices of
$\G$ are denoted by $\cA_B$ and $\cV_B$, respectively. 
\end{defn}

\begin{defn}
A boundary interval  with prime label a minor $H^\Lm$ letter
 is called a {\em minor} boundary interval. A minor boundary interval with 
label of exponential length at least $1$ is called a {\em partisan} boundary
interval.
The boundary components of $\S$ 
containing partisan boundary intervals
are called {\em partisan} 
boundary components. A region of $\G$ containing a partisan boundary interval in
its boundary is called a {\em partisan} region.
\end{defn}
%
For future use we record the following result.
\begin{lemma}\label{HLbound}
The number of $H^\Lm$--vertices of $\G$ is at most $2W_1$ (see (\ref{W1})).
\end{lemma}

\noindent
{\em Proof.} Traverse each boundary component of $\pd\S$ once in the
direction of orientation $\z$ counting $H^\Lm$--vertices as they occur.
\medskip

%
\subsection{Bridge moves and vertex cancellation}
\begin{figure}
\psfrag{e}{$e$}
\psfrag{u}{$u$}
\psfrag{N}{$N$}
\begin{center}
  \mbox{
\subfigure[$N$ non--orientable\label{twisted_loop}]
{ \includegraphics[scale= 0.3,clip]{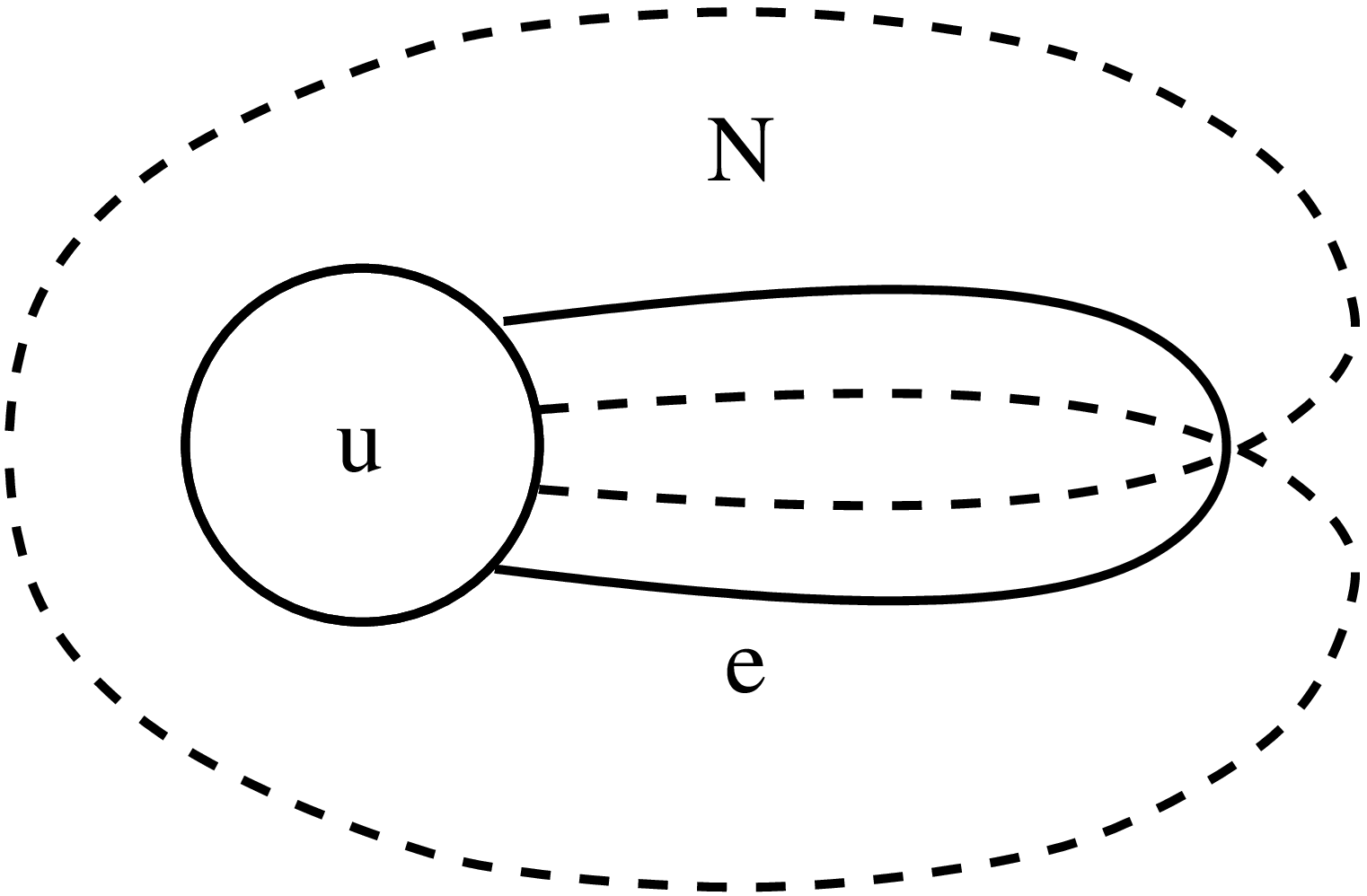} }\qquad
\subfigure[$N$ orientable\label{flat_loop}]  
{ \includegraphics[scale=0.3,clip]{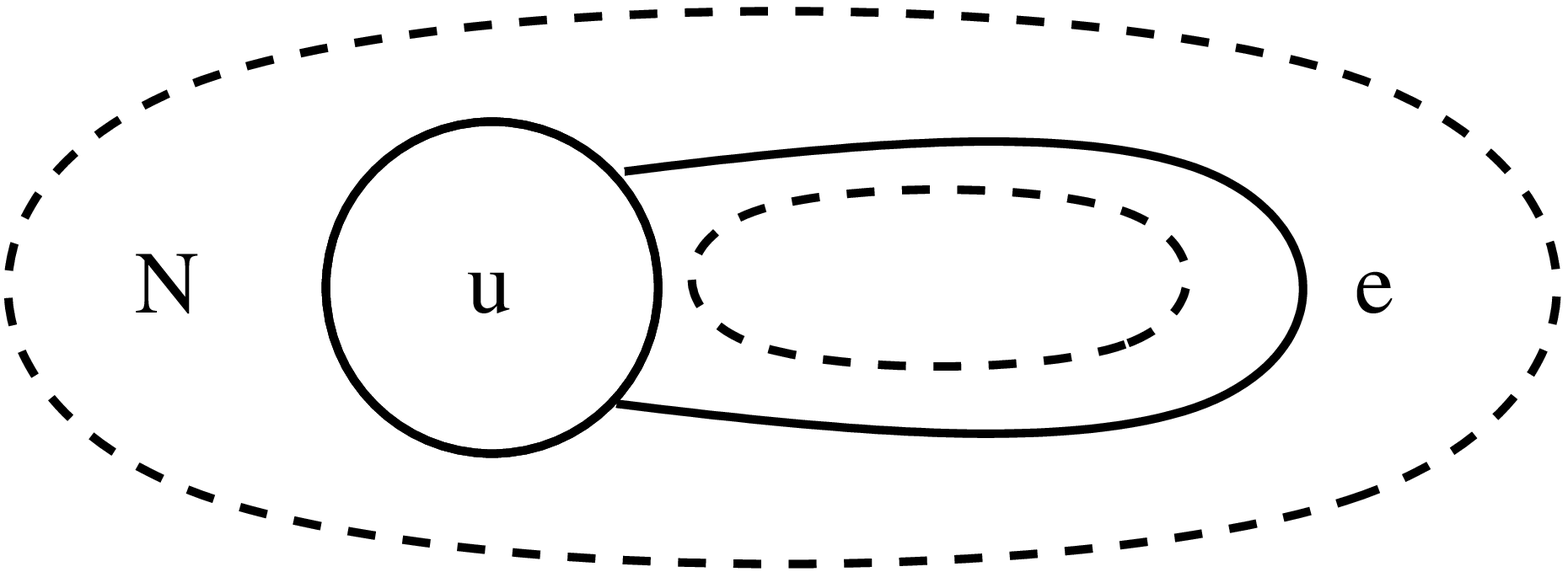} } 
}
\caption{}\label{loop}
\end{center}
\end{figure}

\begin{figure}
\psfrag{e}{$e$}
\psfrag{r^m}{$r^m$}
\psfrag{~r^m}{$\sim r^{-m}$}
\begin{center}
  \mbox{
\subfigure[Confluent\label{e-confluent-y}]
{ \includegraphics[scale= 0.37,clip]{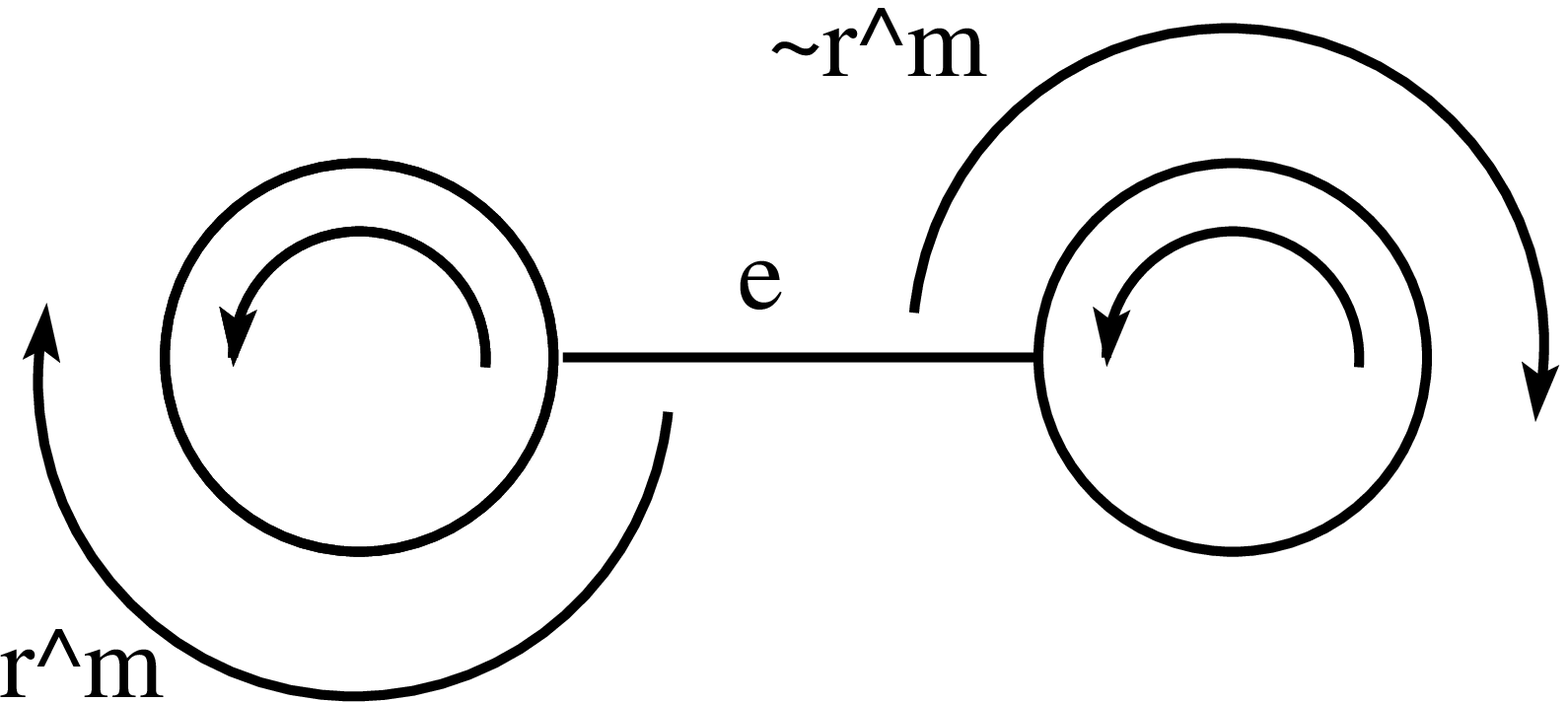} }\qquad
\subfigure[Not confluent\label{e-confluent-n} ]  
{ \includegraphics[scale=0.37,clip]{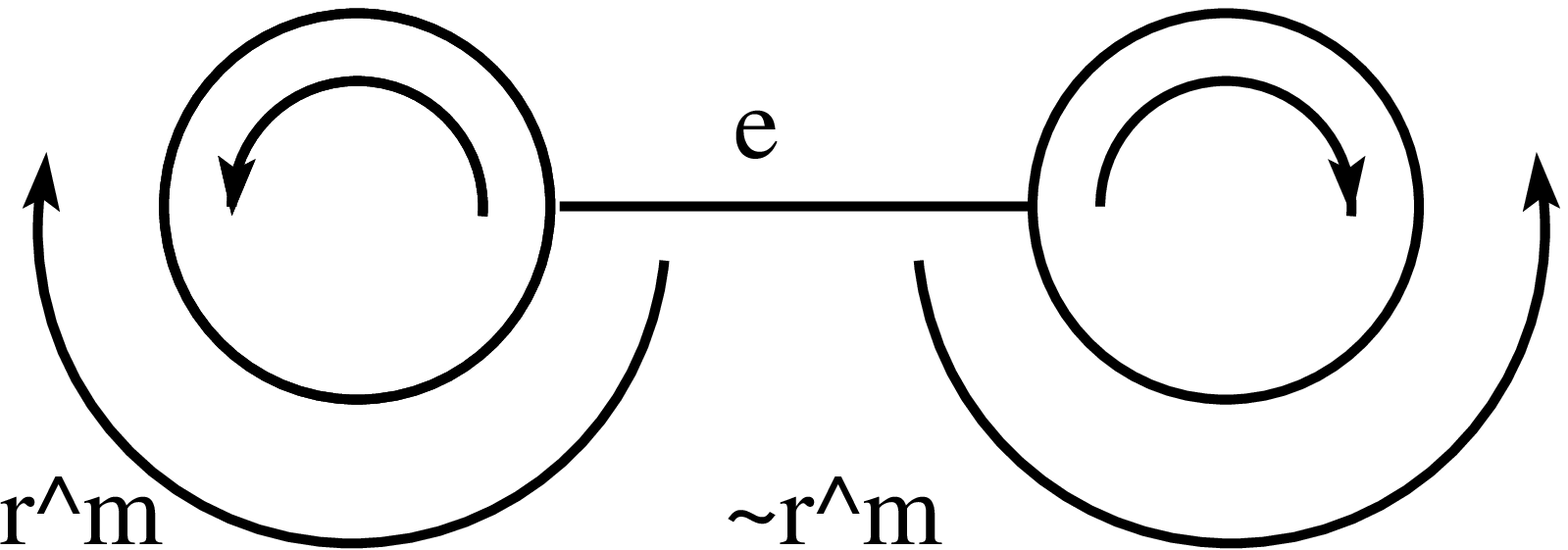} } 
}
\caption{$\sim r^{-m}$ denotes a cyclic permutation of $r^{-m}$}\label{e-confluent}
\end{center}
\end{figure}

Let $u_1$ and $u_2$ be vertices of a picture $\G$, 
with orientation $\z$, and let $e$ be an arc with end--points $x_1$ and $x_2$
on $u_1$ and $u_2$, respectively. 
Let $N$ be a regular
neighbourhood of $u_1\cup u_2\cup e$. If $N$ is orientable
let $\z_e$ be an orientation of $N$.
If $\z_e(u_i)=\z(u_i)$, for both $i=1$ and $2$, or $\z_e(u_i)\neq\z(u_i)$, for both $i=1$ and $2$,
then set $\de=1$. Otherwise set $\de=-1$. If $N$ is non--orientable set $\de=1$.
The latter occurs if $u_1=u_2$ and $N$ is
homeomorphic to a M\"obius band (see Figure \ref{twisted_loop}).
For $i=1,2$ set $\e_i=1$ if $u_i$ is positive and $\e_i=-1$ if $u_i$ is negative.
\begin{defn}\label{e_confluent}
If $\de\e_1\e_2=-1$ we say that $u_1$ and $u_2$ are $e${\em --confluent}.
\end{defn}
Thus $u_1$ and $u_2$ are $e$--confluent if $f$ preserves orientation $\z_e$ 
on precisely one of $\pd u_1$ and
$\pd u_2$. 
Figure \ref{e-confluent-y} shows a pair of $e$--confluent vertices, with $\e_1=\de=1$ and $\e_2=-1$.
Figure \ref{e-confluent-n} shows a pair of vertices which are not $e$--confluent, 
with $\e_2=\de=-1$ and $\e_1=1$.

If $u_1=u_2$ and a regular neighbourhood $N$ of $u_1\cup e$ is non--orientable then $u_1$ is $e$--confluent
to itself, as in Figure \ref{twisted_loop}. On the other hand if $N$ is orientable, as in Figure \ref{flat_loop},
then $u_1$ is not $e$--confluent to itself.

\begin{figure}
\psfrag{a}{$a$}
\psfrag{b}{$b$}
\psfrag{c}{$c$}
\psfrag{d}{$d$}
\psfrag{a^-1}{$a^{-1}$}
\psfrag{b^-1}{$b^{-1}$}
\psfrag{c^-1}{$c^{-1}$}
\psfrag{d^-1}{$d^{-1}$}
\psfrag{e}{$e$}
\psfrag{w1}{$w_1$}
\psfrag{w1-1}{$w_1^{-1}$}
\begin{center}
{\mbox{
\subfigure[$\de=1$\label{vertex_cancel-cy}]
{ \includegraphics[scale=0.37]{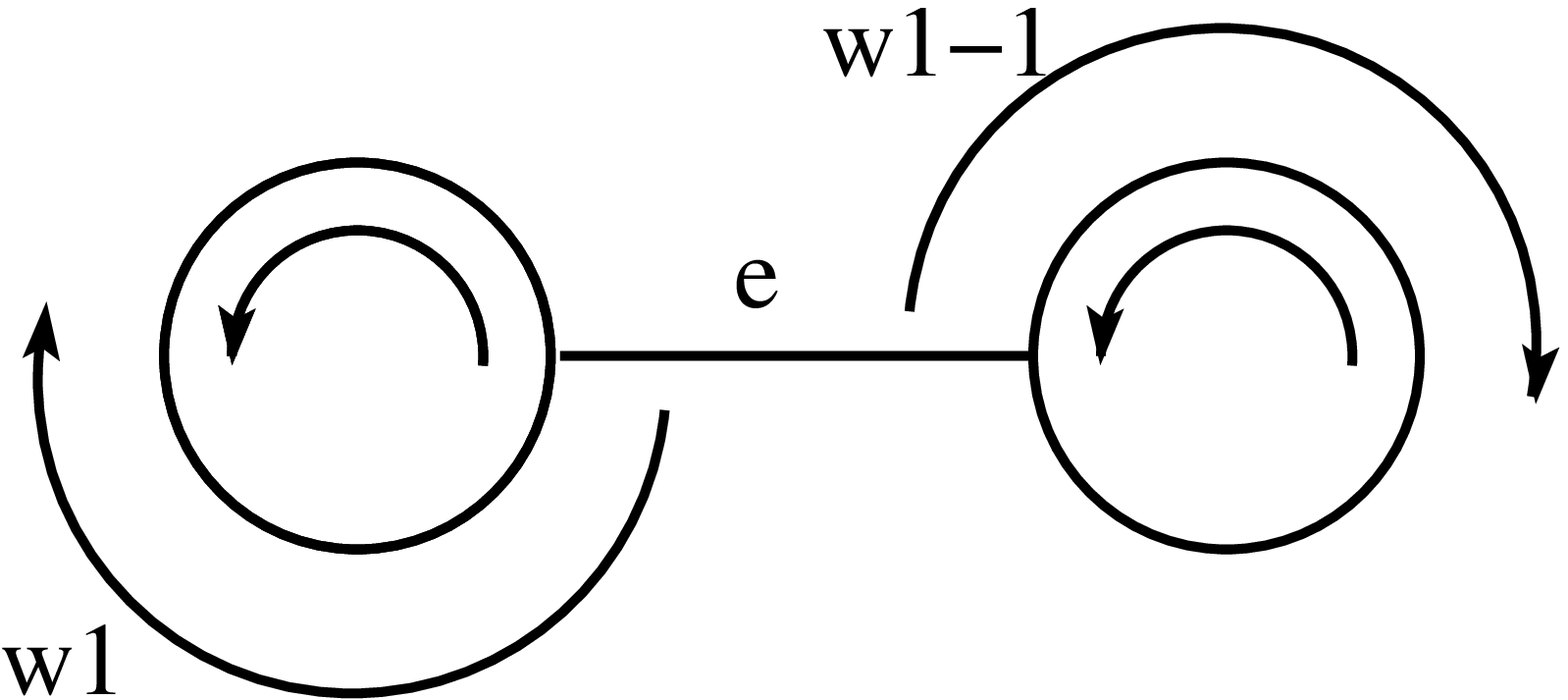} }\qquad
\subfigure[$\de=-1$\label{vertex_cancel-cn} ]  
{\includegraphics[scale=0.37]{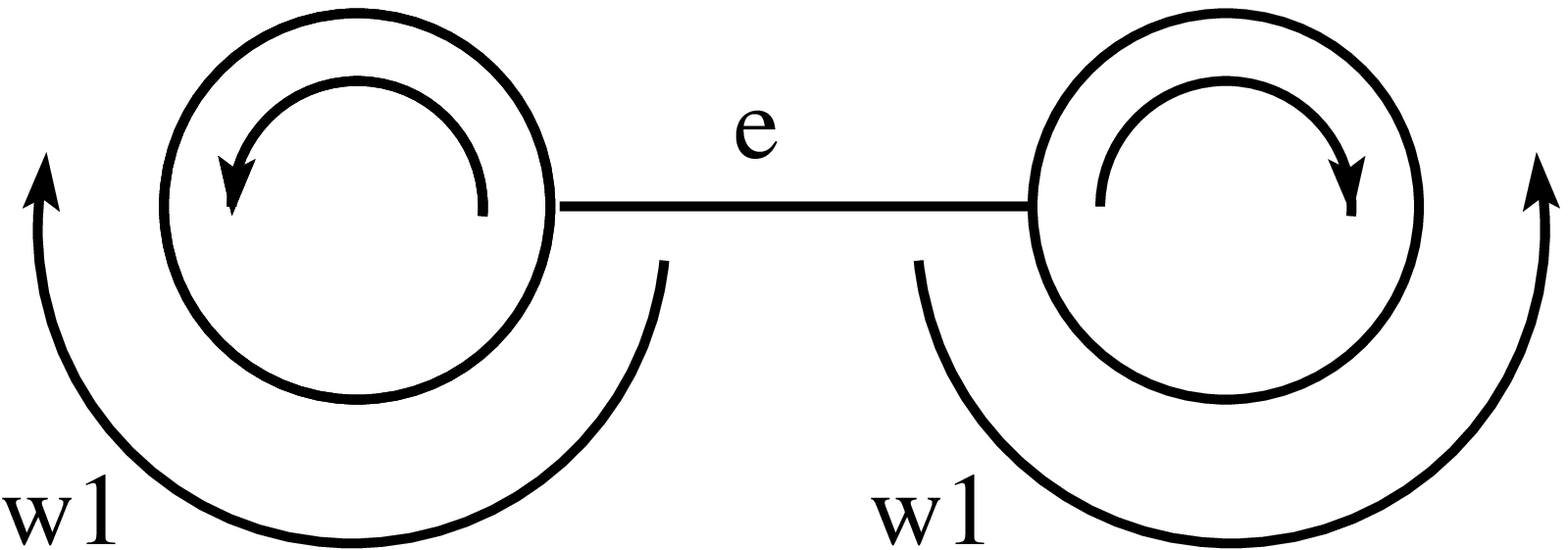} } 
}}
\caption{Cancelling vertices}\label{vertex_cancel}
\end{center}
\end{figure}

In the above setting let $w_i$ be the label of $u_i$ with base point $x_i$ (so
$w_i$ is obtained by reading in the direction opposite to $\z(u_i)$). 
If $w_1w_2^{\de}=1$
in $H$ then $u_1$ and $u_2$ are said to be vertices which {\em cancel}
(see Figure \ref{vertex_cancel}). 
A picture
with no pair of distinct vertices which cancel is said to be {\em reduced}. 

Suppose that $u_1$ and $u_2$ cancel. If $u_1$ and $u_2$ are $e$--confluent then we may
assume, without loss of generality, that $\e_1=\de=1$ and $\e_2=-1$, that $w_1=r^m$ and $w_2$ is 
a cyclic permutation of $r^{-m}$, so $w_2^{-1}$ is a cyclic permutation of $r^m$. 
Suppose that $r=r_0r_1$ and $w_2^{-1}=(r_1r_0)^m$. As $u_1$ and $u_2$  cancel and $\de=1$ we have
\[(r_0r_1)^m=w_1=w_2^{-1}=(r_1r_0)^m.\]
It follows that $r^m$ has periods $\lm=\textrm{min}\{|r_0|,|r_1|\}$ and $l$. From
\cite[Section 3, Proposition 1]{\Ha} it follows that $r^m$, and hence $r$, has period $d=\gcd(\lm,l)$.
By assumption $r$ is not a proper power so this cannot happen unless $d=l$ and $\lm=0$, in which case
$w_2^{-1}=r^m$.

Now suppose that $u_1$ and $u_2$ cancel but are not $e$--confluent. Without loss of generality
we may assume that $\e_1=\e_2=\de=1$, that $w_1=r^m$ and that $w_2=r_1r^{m-1}r_0$, where $r=r_0r_1$.
We have 
\[(r^{-1}_1r_0^{-1})^m=w_1^{-1}=w_2=r_1r^{m-1}r_0,\]
so $r_1^2=r_0^2=1$. It follows that $r$ has a cyclic permutation $xUyU^{-1}$, where
$x,y\in A\dcup B$, with $x^2=y^2=1$. That is $G$ is of type $E(2,2,m)$.

\begin{figure}
\psfrag{a}{$a$}
\psfrag{b}{$b$}
\psfrag{'g}{$\g$}
\begin{center}
{
\includegraphics[scale=0.4]{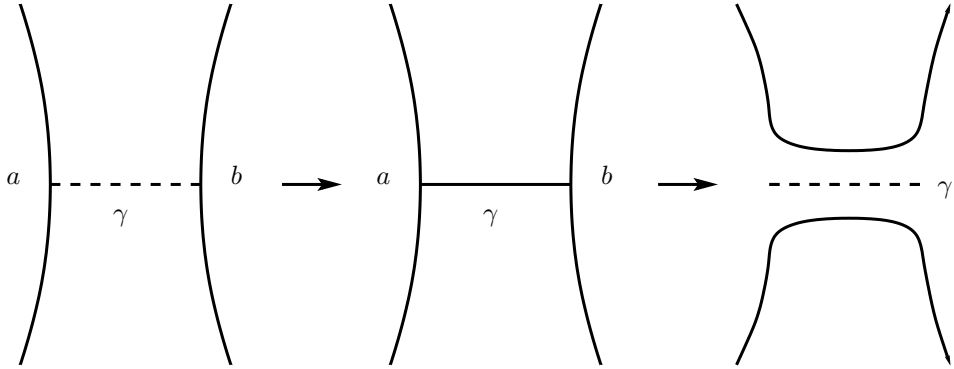}
}\\
\caption{A bridge move}\label{bridge_move}
\end{center}
\end{figure}

Now suppose that $a$ and $b$ are distinct arcs of $\G$ and that $\g$ is 
the image of an embedding of $[0,1]$ 
in $\S$ with one end--point on $a$ and the other on $b$ such that $\g\cap \G=\{a,b\}$. Then 
$f(\g)$ is a closed path in $(X_O,p)$, for $O=A$ or $B$. If $f(\g)$ is nullhomotopic, then
after a homotopy  (constant outside a regular neighbourhood of $\g$) we may assume that
$f(\g)=p$. A further homotopy can be effected to make $f(\g)=x_O$, resulting in 
a new picture, as shown in Figure
\ref{bridge_move}. We call the transformation from the original to the
new picture  a {\em bridge move}.
    
\begin{figure}
\psfrag{a}{$a$}
\psfrag{b}{$b$}
\psfrag{u}{$u$}
\psfrag{v}{$v$}
\psfrag{a^-1}{$a^{-1}$}
\psfrag{b^-1}{$b^{-1}$}
\begin{center}
{
\includegraphics[scale=0.4]{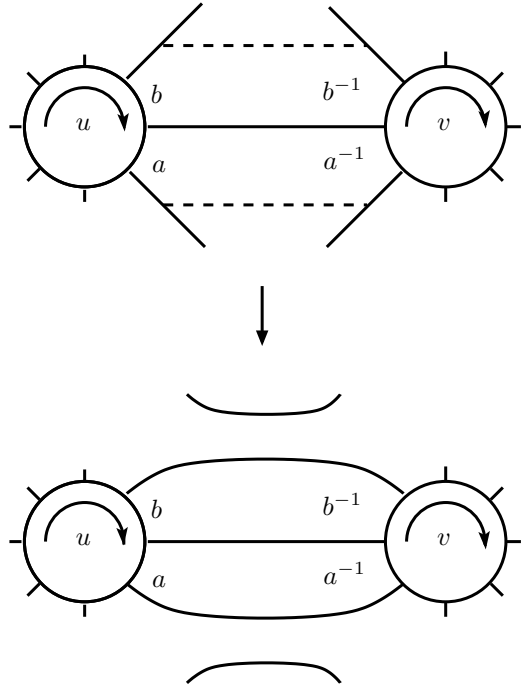}
}\\
\caption{Bridge moves on a pair of cancelling vertices}\label{bridge_pair}
\end{center}
\end{figure}

\begin{figure}
\psfrag{a}{$a$}
\psfrag{b}{$b$}
\psfrag{u}{$u$}
\psfrag{v}{$v$}
\psfrag{a^-1}{$a^{-1}$}
\psfrag{b^-1}{$b^{-1}$}
\begin{center}
{
\includegraphics[scale=0.4]{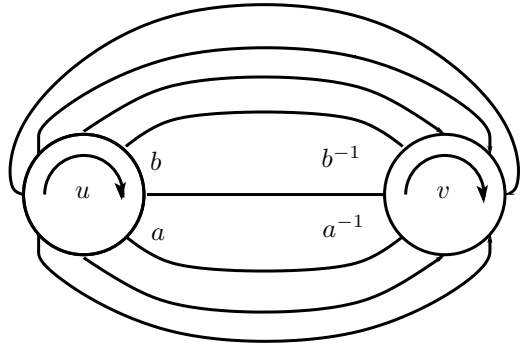}
}\\
\caption{A floating dipole}\label{floating_dipole}
\end{center}
\end{figure}

If $u$ and $v$ are a pair of distinct vertices which cancel then we can perform
bridge moves on arcs incident  $u$ and $v$, as shown in Figure \ref{bridge_pair},
to increase the number of arcs joining $u$ to  $v$. Repeating  such
bridge moves we obtain a picture in which all arcs incident to $u$ are incident to $v$,
in which case we call $u$, $v$ and their incident arcs a {\em floating dipole}: see Figure 
\ref{floating_dipole}.
An arc $c$ of $\G$ such that $c$ is a simple closed curve bounding a
disk $D$ with int$(D)\cap\G=\nul$ is called a {\em floating arc}.

If $C$ is a simple closed arc in $\S\bl \G$ and one component of $\S\bl C$ is 
a disk $C_D$ such that $\G\cap C_D$ is a floating dipole or floating
arc then $f(C)$ is null--homotopic in $Y$. Performing a homotopy, constant outside
$C_D\cup C$, we may replace $f$ by a map such that $f(C_D)=x_O$, where $O=A$ or $B$.
Therefore  we may form a 
new picture on $\S$ by removing all arcs and vertices of $C_D\cap\G$ from $\G$.
We call this operation {\em removal} of a floating dipole or arc.
We consider the equivalence relation $\sim$ generated by the relations 
(\ref{pic_equiv1})--(\ref{pic_equiv3}) below 
on the set of pictures on 
a surface $\S$,
with $n$ boundary components  
with fixed boundary labels $u_1,\ldots ,u_n\in H$.
Suppose $\G$ and $\G^\prime $ are such pictures on $\S$ over $G$. Then 
$\G \sim \G^\prime$ if $\G^\prime$ is obtained from $\G$ by
\be
\item\label{pic_equiv1} an isotopy of $\S$,
\item\label{pic_equiv2} a bridge move or
\item\label{pic_equiv3} removal of a floating dipole or floating arc. 
\ee

\begin{defn}\label{pic_effic}
A picture is {\em efficient} if it has fewest arcs in its equivalence class. 
\end{defn}
It 
follows that an efficient picture is reduced, since otherwise we can perform 
bridge moves on arcs incident to a pair of vertices which cancel until a floating
dipole is produced. Removal of the floating dipole results in a picture with
fewer arcs, contrary to efficiency of the original picture. 

Let $g$ and $n$ be fixed integers and $u_1,\ldots ,u_n$ fixed elements
of $H$. Consider the set $T$ of all pictures
over $G$ on surfaces 
with $n$ boundary components and genus at most $g$, 
with boundary labels $u_1,\ldots ,u_n$.  
\begin{defn}\label{pic_min}
A  
picture in $T$ is said to
be {\em minimal} if no other picture in $T$ has fewer arcs.
\end{defn}
Minimal pictures are efficient and efficient pictures are reduced. 
In general it is not the case that reduced pictures
are efficient or that efficient pictures are minimal.
\begin{prop}[\cite{\DHa}]\label{no_closed_arc}
 Let $\G$ be a minimal picture over $G$ on a 
surface $\S$. Then no arc of $\G$ is closed curve.
\end{prop}
\medskip

\noindent
{\em Proof.}
Suppose that $c$ is an arc of $\G$ and a closed curve. The closure
of $\S\backslash c$ has $2$ boundary components $c_1$ and $c_2$ corresponding to $c$.
Form a new surface $\S_1$ by capping  $c_i$ with a disk, for $i=1,2$. Then $\G\backslash c$
is a picture on $\S_1$, since $f(c)= p\in Z$. The genus of $\S_1$ is no more than that of $\S$
and $\G\backslash c$ has the same boundary labels as $\G$, but fewer arcs. 
This contradicts the minimality of 
$\G$.
\medskip

\begin{prop}[\cite{\DHa}]\label{genreg}
 Let $\G$ be a minimal picture over $G$ on a compact
connected surface $\S$. Then
$\chi(\D)\ge\chi(\S)$, for all regions $\D$ of $\G$.
\end{prop}
\medskip

\noindent
{\em Proof.} Suppose that $\D$ is a region such that
$\chi(\D)<\chi(\S)$. Let $\S_0=\S\bl\textrm{int}(\D)$ and let 
$\Omega_0,\ldots,\Omega_k$ be the connected components of $\S_0$. Then
$\chi(\D)<\chi(\S)=\chi(\S_0)+\chi(\D)=\sum_{i=0}^k\chi(\Omega_i)+\chi(\D).$
Therefore $\chi(\Omega_i)>0$, so $\Omega_i$ is a disk or a sphere, for some $i$. 
As $\S$ is connected this means that $\Omega_i$ is a disk with 
$\Omega_i\cap \D=C$, where $C$ is homeomorphic to $S^1$, for such $i$.
Let $\Omega^\prime$ be a regular neighbourhood of $\Omega_i$ 
and let $\G^\prime=\G\cap \Omega^\prime$. Then $\Omega^\prime$ is a
disk and the boundary of $\Omega^\prime$ lies in $\D$. Hence the
boundary label of $\G^\prime$ is an element $u$ of $A\dcup B$ such
that $u=1$ in $G$. As the Freiheitssatz holds for $G$
(\cite{Howie89}) it follows that $u=1$ in $A$ or $B$, as appropriate. Thus we may
replace $\G^\prime$ with the empty picture on $\Omega^\prime$ without
altering $\S$ or the boundary labels of $\G$. This contradicts the
minimality of $\G$.
\subsection{Parallel arcs and the graph of a picture}
Given two arcs $a_1$ and $a_2$ of a picture $\G$ such that either
\be
\item $a_1$ and $a_2$ bound an annulus $S$ with int$(S)\cap \G=\nul$ or 
\item there are corners $c_1$ and
$c_2$ of $\G$ such that $a_1$, $a_2$, $c_1$ and $c_2$ bound a disk
$\D$ with int$(\D)\cap \G=\nul$,
\ee
then $a_1$ and $a_2$ are said to be {\em parallel} arcs.
There is an  equivalence relation  on the set of arcs of $\G$ generated 
by the parallel relation on arcs. The equivalence class of an arc $a$ is
called an {\em edge} of $\G$, denoted $\bar{a}$. If an edge consists of 
$d$ parallel arcs we say that it has {\em width} $d$. We form the {\em
  graph} $\bar\G$
of $\G$
on $\S$ as follows. Collapse each vertex of $G$ to a point to give a
vertex of $G$.
Collapse each class of arcs $\bar a$ onto $a$: the resulting arcs are
the edges of $\bar \G$, except for any that are closed loops (meeting
no vertex of $\G$) and these are discarded. 
In addition there is one vertex on $\pd\S$ wherever an edge of 
$\bar\G$ meets $\pd\S$. If $e^+$ is one end of an edge $e$ of $\bar \G$ then
the end points of arcs of $\G$ which are identified to form $e^+$ comprise an
{\em end} of the class of arcs corresponding to $e$. If $\bar a$ is a class of arcs with  at least
one end on the boundary then we say that $\bar a$ is a {\em boundary} class. Otherwise
$\bar a$ is an {\em interior} class. 
An arc of $\G$ with both ends on $\pd\S$ is called an {\em arc of type} $I$.
An edge of $\bar\G$ which is the image of a class of arcs of type $I$ 
is called an {\em edge of type} $I$.
We define the {\em degree} of a vertex $v$ of $\G$ to be the degree of the 
corresponding vertex in the graph $\bar\G$, that is
$\deg(v)=$ number of classes of arc incident to $v$.

If $a$ is an arc of an edge $e$ and $a$ is incident to vertices $u$ and $v$ then
$u$ and $v$ are $a$--confluent if and only if they are $b$--confluent for all arcs
$b$ of $e$. In this case we say that $u$ and $v$ are $e${\em --confluent}.

The following facts are established in \cite{\Ha}. A word $x_1\cdots x_n$ is said to have
{\em period} $\lm$ if $x_{i+\lm}=x_i$, for all $i$ such that $1\le i\le n-\lm$. 
\begin{lemma}\label{bound_edge_width}
Assume $G$ is not of type $E(2,*,m)$.
Let vertices $u_1$ and $u_2$ of a picture $\G$ over $G$  be joined by an edge $e$ of 
width $k$. Assume either that $\G$ is reduced  and $u_1$ and $u_2$ are distinct or that 
$\G$ is minimal. Then $k\le l-2$. 
\end{lemma}
\medskip

\begin{figure}
\psfrag{rm}{$r^m$}
\psfrag{r-m}{$r^{-m}$}
\psfrag{s1}{$s_1$}
\psfrag{s2}{$s_2$}
\begin{center}
{
\includegraphics[scale=0.4]{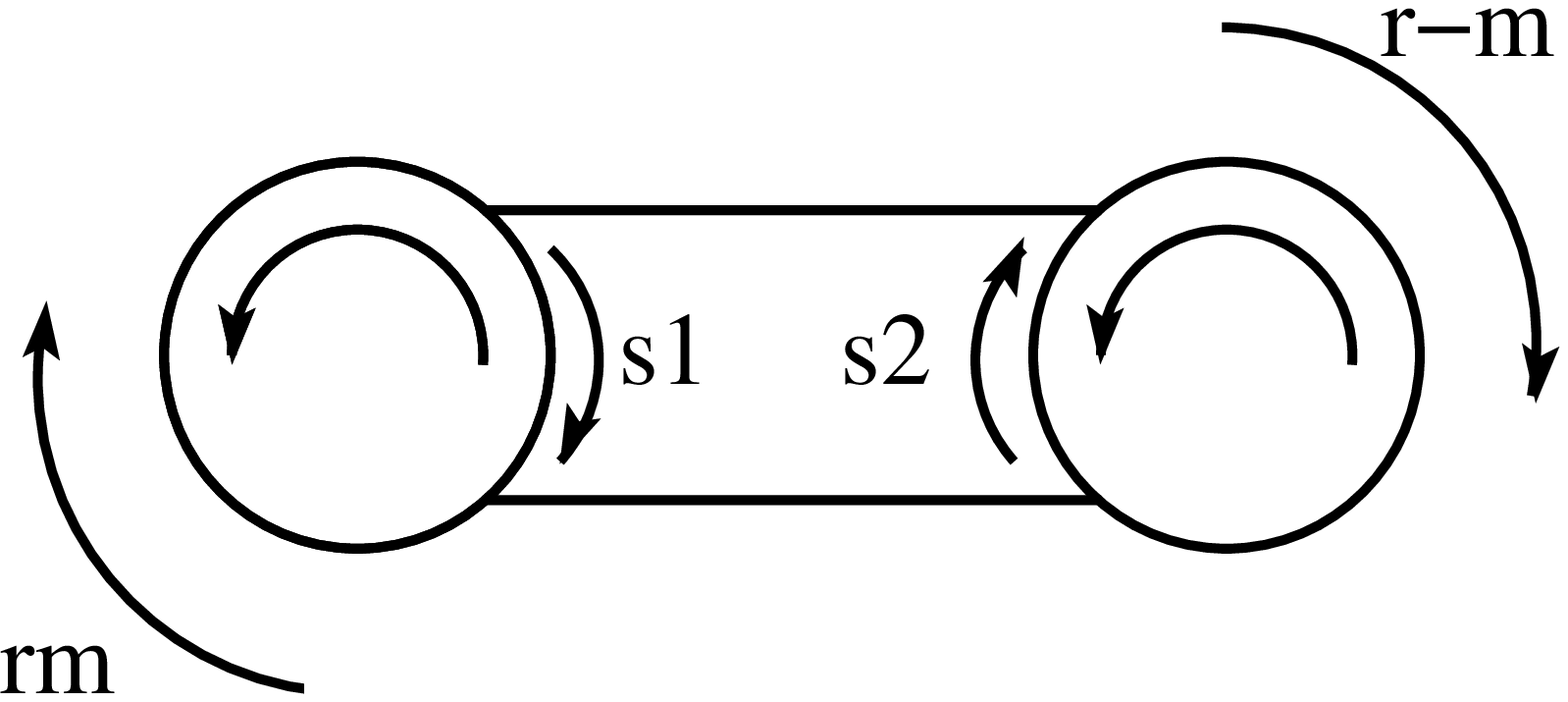}
}\\
\caption{}\label{confluent_edge}
\end{center}
\end{figure}

\noindent{\it Proof.} The ends of $e$ determine subwords $s_1$ and $s_2$ of the boundary labels
$w_1$ and $w_2$ of $u_1$ and $u_2$, respectively, with $l(s_1)=l(s_2)=k-1$. Suppose first that 
$u_1$ and $u_2$ are $e$--confluent.  Without loss of generality we may assume that $\e_1=\de=1$ and
$\e_2=-1$ (with the notation of Definition \ref{e_confluent}). 
Then $s_1$ is a cyclic subword of $r^m$, $s_2$ is a cyclic subword of $r^{-m}$ and $s_2^{-1}=s_1$
(see Figure \ref{confluent_edge}).
Therefore $r^m$ has identical subwords $s_1$ and $s_2^{-1}$. We may assume that $s_1$ is an initial
subword of $r^m$. Then $r^m$ begins $r_0s_2^{-1}$, where $r_0$ is an initial subword of $r$ of length
$l(r_0)=\lm\ge 0$. As $r$ is cyclically reduced in $A*B$ it follows that $\lm$ is even.
We may also assume that $\lm\le l/2$ (otherwise reverse the roles of $s_1$ and $s_2$).
If $k\ge l-1$ this implies that the union of $s_1$ and $s_2$ is an initial subword of $r^m$ of period
$\lm$. If $\lm>0$ then $\lm\ge 2$ so  $l(s_2)+\lm\ge l$ and 
then it follows that $r^m$ has periods $\lm$ and $l$. Hence, from
\cite[Section 3, Proposition 1]{\Ha}, $r^m$ and also $r$ have period $d=\gcd(\lm,l)$. This cannot occur
as $r$ is not a proper power. Therefore $\lm=0$ and $s_1$ and $s_2^{-1}$ coincide as subwords of 
$r^m$.  In this case, if $u_1\neq u_2$ then $u_1$ and $u_2$ cancel.
As $\G$ is reduced this is a contradiction.

\begin{figure}
\psfrag{a}{$a$}
\psfrag{b}{$b$}
\psfrag{c}{$c$}
\psfrag{s1}{$s_1$}
\psfrag{s2}{$s_2$}
\psfrag{N}{$N$}
\begin{center}
{
\includegraphics[scale=0.4]{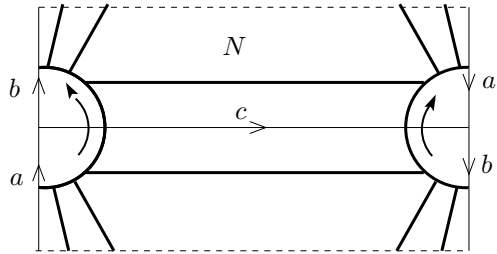}
}\\
\caption{$N$ after cutting along arc $ab$}\label{mobius_band}
\end{center}
\end{figure}

\begin{figure}
\psfrag{a}{$a$}
\psfrag{b}{$b$}
\psfrag{c}{$c$}
\psfrag{N}{$N$}
\begin{center}
{
\includegraphics[scale=0.4]{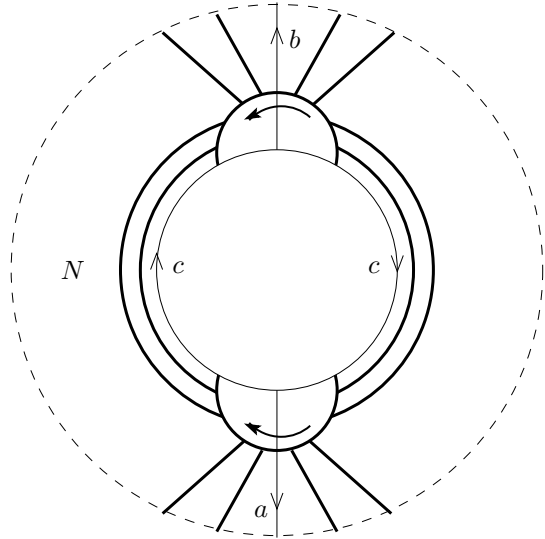}
}\\
\caption{$N$ after cutting along $c$ and pasting $a$ and $b$}\label{cut_mobius_band}
\end{center}
\end{figure}

\begin{figure}
\psfrag{a}{$a$}
\psfrag{b}{$b$}
\psfrag{c}{$c$}
\psfrag{N}{$N$}
\begin{center}
{
\includegraphics[scale=0.4]{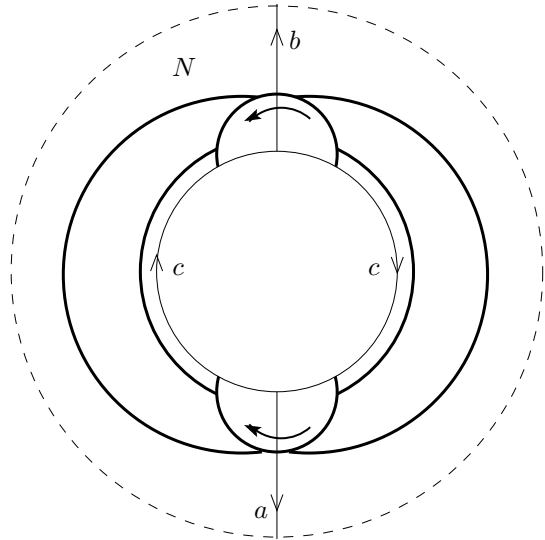}
}\\
\caption{$N$ after performing bridge moves}\label{floating_projective}
\end{center}
\end{figure}

On the other hand, if $\lm=0$, $u_1=u_2$ and $u_1$ is $e$--confluent to itself
then a regular neighbourhood $N$ of $u_1\cup e$ is homeomorphic to a M\"obius band. In this case
we can perform bridge moves on arcs incident to $u_1$ and outside $e$ to  form a class containing
all arcs incident to $u_1$. We illustrate this diagrammatically: suppose $N$ is as shown in 
Figure \ref{mobius_band}. 
Cutting along $c$ and pasting together both copies of $a$ and $b$ we obtain the diagram
of $N$ shown in Figure \ref{cut_mobius_band}.
Clearly the number of arcs incident to the upper and lower halves of the vertex in this diagram
are equal and so we may perform bridge moves as claimed. The result is shown in Figure \ref{floating_projective}.
Since the dotted boundary of this diagram has trivial label we may remove $u_1$ and all its incident
arcs from $N$, leaving a picture on $\S$ with the same boundary labels as $\G$ but fewer arcs. This contradicts
minimality of $\G$. 
Thus if $u_1$ and $u_2$ are $e$--confluent then 
$k\le l-2$.

Now suppose that $u_1$ and $u_2$ are not $e$--confluent and that $k\ge l-1$. Without loss of
generality we may assume that $\e_1=\e_2=\de$ and that $s_1$ is an initial subword of $r^m$.
Then $s_2$ is a cyclic subword of $r^m$ and $s_1$ is identically equal to $s_2^{-1}$, as a word in $A*B$.
Again we may assume that $r^m$ begins $r_0 s_2=U$, with $l(r_0)=\lm\le l/2$. If $s_1$ has even length
then $s_1$ and $s_2$ 
begin with letters of different groups,  
so $\lm$ must be odd. Similarly if $l(s_1)$ is odd then $\lm$ must be even. As $k\ge l-1$ we have
$l(s_i)\ge l-2$ and $l(U)=l(s_2)+\lm\ge l-1$. We may write $s_1=r_0t$, where $s_2=tr_1$, since 
$\lm\le l-2\le l(s_i)$. Hence $s_1=s_2^{-1}$ implies $t=t^{-1}$ and $r_0=r_1^{-1}$. 
Hence $U=r_0s_2=r_0tr_0^{-1}$, where $t=t^{-1}$, so
$t$ is of odd length and its
middle letter has order $2$. Therefore  $U=VyV^{-1}$, for some $V$. If $l(U)=l-1$ then 
$r=VyV^{-1}x$ and $G$ is of type $E(2,*,m)$, a contradiction. Suppose $l(U)\ge l$. We have $r^m=VyV^{-1}W$ and 
we may write $V=V_0xU_0$, where $x\in A\dcup B$, $l(U_0yU_0^{-1})=l-1$ and $l(V_0)\ge 0$.
Hence $r^m$ has a cyclic permutation $xU_0yU_0^{-1}x^{-1}V_0^{-1}WV_0$ which must have period $l$. 
Again a contradiction arises as $x=x^{-1}$ and 
$G$ is of type $E(2,2,m)$. 
\begin{defn}\label{minimalistic}
A picture $\G$ on a surface $\S$ is called {\em minimalistic} if the following conditions hold.
\be
\item No arc of $\G$ is a closed curve.
\item $\chi(\D)\ge \chi(\S_\D)$, whenever $\D$ is a region of $\G$ and $\S_\D$ is the connected 
component of $\S$ containing $\D$.
\item No interior edge of $\G$ has width more than $l$. If $G$ is not of type $E(2,*,m)$ then no
interior edge has width more than $l-2$.
\ee
\end{defn}
From Proposition \ref{no_closed_arc}, Proposition \ref{genreg} 
and Lemma \ref{bound_edge_width}, minimal pictures are minimalistic. The converse
does not hold in general.
\subsection{Routes and distance}
\begin{figure}
\psfrag{v}{$v$}
\psfrag{u}{$u$}
\psfrag{e}{$e$}
\psfrag{'D}{$\D$}
\begin{center}
  \mbox{
\subfigure[paired vertices $u$ and $v$\label{PV}]
{\includegraphics[scale= 0.6,clip]{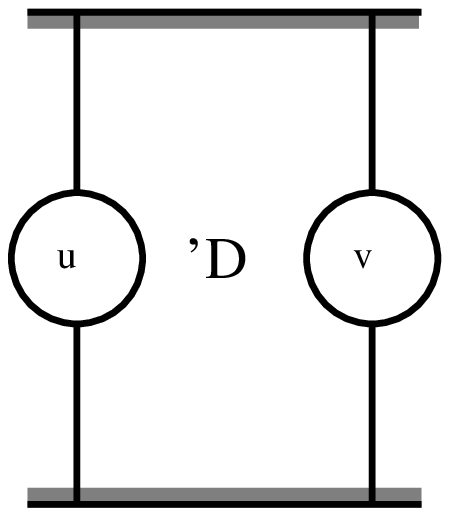} }\qquad
\subfigure[vertex $u$ paired to class of arcs $e$\label{PE}]  
{ \includegraphics[scale=0.6,clip]{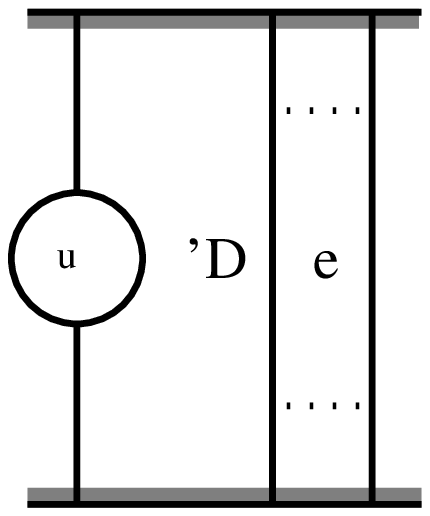} } 
}
\caption{}\label{paired}
\end{center}
\end{figure}

If $u$ and $v$ are vertices incident to a region $\D$ satisfying
$\rho(\D)=\bt(\D)=2$, $\chi(\D)=1$, $t(\D)=4$ and $\e(\D)=0$, as shown in
Figure \ref{PV}, then $u$ and $v$ are said to be {\em paired}. If $e$
is a class of arcs having both end points on $\pd \S$ and $e\cap \pd
\D \neq \nul$,
where $\D$ is a region satisfying $\rho(\D)=1$, $\bt(\D)=2$, $\chi(\D)=1$, $t(\D)=3$ 
and $\e(\D)=0$, as shown in
Figure \ref{PE}, then $u$ and $e$ are said to be {\em paired}. 

A {\em route} in $\G$, from $x_1$ to $x_n$,
is a sequence $x_1,\ldots,x_n$, where $x_i$ is a vertex or
a class of arcs such that the set $\{x_i,x_{i+1}\}$    consists either of
\be
\item two vertices which are joined by an arc or 
\item a vertex and  a boundary class of arcs paired to that vertex or
\item two vertices which are paired, 
\ee
for $i=1,\ldots ,n$.
Let $l_1=\{i: x_i\textrm{ is a boundary class of arcs}\}$ and 
$l_2=\{i:x_i\textrm{ and } x_{i+1}\textrm{ are both vertices}\}$. The {\em length} of the
route is then $|l_1|+|l_2|$. A {\em route} in $\bar\G$ is defined to be the
image of  a route in $\G$ under the quotient map $\G\maps \bar\G$. 
A route is {\em simple} if $x_i\neq x_j$, for
all $i\neq j$.

If $S$ is a subset of $\cV$ then a vertex $v$ is said to be {\em distance}
$d$ from $S$ if the length of a shortest route from $v$ to a vertex 
of $S$ is $d$. If there is no route from 
$v$ to a vertex of $S$ then the distance of $v$ from $S$ is defined to
be $\infty$. If $e$ is a class of arcs that is paired to a vertex then
the {\em distance} of $e$ from $S$ is $d+1$, where $d$ is the minimum of
the distances from $S$ of the vertices paired to $e$.
\subsection{Pictures and equations}

The correspondence between solutions of quadratic equations and pictures is given
in the next three results.
\begin{theorem}\label{fundpic}
Let $u_1,\ldots, u_k$ be elements of $H$ and let $t,p$ be non--negative
integers. Then there exist elements $w_1,\ldots , w_{k+2t+p}\in H$ such that 
\begin{equation}\label{pic_eqn}
w=\prod_{i=1}^k w_i^{-1}u_iw_i 
\prod_{i=1}^t [w_{i+k},w_{i+k+t}]\prod_{i=1}^p w_{i+k+2t}^2=1 \in G
\end{equation}
if and only if there exists a picture $\G$ over $G$ on a surface $\S$,
of genus $t+p/2$, orientable if $p=0$, 
with $k$ boundary components which have 
labels $u_1,\ldots ,u_k$ (with respect to some 
orientation of $\pd \S$). 
\end{theorem}
\medskip

\noindent
{\em Proof.} 
In the notation of Section \ref{qwords}, let $F=F(D\dcup X)$ and
$q=q(0,k,t,p)$.
Then $\S(q)$ is the connected sum of a $k$--punctured sphere with 
$t$ torii and $p$ projective planes, and so has genus $t+p/2$ and is orientable
if $p=0$,  and has boundary components
$\bt_1,\ldots, \bt_k$ with labels $d_1,\ldots, d_k$, respectively.

Suppose a picture $(f,\G$) with orientation $\z$ 
exists on a surface $\S$ of genus $t+p/2$, orientable if $p=0$, with boundary
components $\wt\bt_1,\ldots ,\wt{\bt_k}$ having labels $u_1,\ldots ,u_k$. 
Assume first that $\S$ is connected.
Then $\S$ is homeomorphic to $\S(q)$ and we may choose this
homeomorphism $\theta :\S(q)\maps \S$ so
that the boundary component $\bt_i$ is mapped  to
$\wt\bt_i$, so that the base point of $\bt_i$ is mapped to a base point of $\wt{\bt_i}$ and so that
the orientation  of $\wt{\bt_i}$ is the same as that induced by $\theta(\bt_i)$.
By construction $\S(q)$ is a quotient of the disk $D$ with boundary 
$\pd D$ labelled $q$, so $f\circ \theta$ can be regarded as a map from
$D$ to $Z$. Then $f\circ\theta|_{\pd D}$ maps each labelled interval
of $\pd D$ to a path in $Y$ representing an element of $H$. 
Let $\phi$ be the map from $F$ to $H$ induced by
$f\circ \theta|_{\pd D}$ on the labelled intervals of $\pd D$. 
Then $\phi$ maps the boundary label $q$ of $D$ 
to an element of $H$ representing the identity of $G$.
Furthermore $\phi(d_i)=u_i$, for $i=1,\ldots, k$ and 
$\phi(x_i)=w_i\in H$, for 
$i=1,\ldots,k+2t+p$, and  it follows that $w=1$ in $G$, as required.

Suppose now that  $\S$ has connected components
$\S_1,\ldots, \S_n$, for some $n\ge 2$. We may assume that $\S_j$ is 
the connected sum of a $k_j$--punctured sphere with $t_j$ torii and 
$p_j$ projective planes where $k_j>0$, $\sum_{j=1}^n k_j=k$, 
$\sum_{j=1}^n t_j=t$ and $\sum_{j=1}^n p_j=p$. 
Let the boundary labels of $\G_i=\S_i\cap \G$ be $u_{i_1},\ldots ,u_{i_{k_j}}$,
where $i_1<\cdots <i_{k_j}$. The argument above implies the existence
of elements $w_{ij}$, $x_{ij}$, $y_{ij}$, $z_{ij}\in H$ such that
\[
 \prod_{i=1}^{k_j} w_{ij}^{-1}u_{i_j}w_{ij} 
\prod_{i=1}^{t_j} [x_{ij},y_{ij}]\prod_{i=1}^{p_j} z_{ij}^2=1 \in G,
\]
for $j=1,\ldots ,n$.
The equality (\ref{pic_eqn}) follows after some rearrangement of these equations. 

Conversely, suppose that there exist $w_i$ in $H$ satisfying
(\ref{pic_eqn}). Let $F$, $q$ and $\S(q)$ be defined as above, where 
$\S(q)$ is a quotient of a disk $D$ with boundary label $q$.  Let $\hat d_i$ and 
$\hat x_i$ be
the oriented sub--intervals of $\pd D$ labelled $d_i$ and $x_i$, respectively. 
For $g\in H$ let $\wt g$ be a closed path in $(Y,p)$ representing $g$. Define a map 
$f_\pd:\pd D\maps Z$ by $f_\pd(\hat x_i)= \wt w_i$ and $f_\pd(\hat d_i)=\wt u_i$. 
As $w=1$ in $G$ the map $f_\pd$ maps $\pd D$ to 
a path which is nullhomotopic in $Z$ and so extends to a map 
$f:D\maps Z$. As $f_\pd$ is identical on any two intervals of $\pd D$ which are identified
in forming $\S(q)$, $f$ factors through $\S(q)$ and  we may 
regard $f$ as a map from $\S(q)$ to $Z$. We may make 
$f$ transverse to the mid--point $c$ of $\D_s$ with a homotopy which fixes 
$\pd\S(q)$. A further homotopy ensures that each component of the closure
of $f^{-1}(\inr(\D_s))$ is a regular neighbourhood of $f^{-1}(c)$. Let $\S_0=\S(q)\bl 
f^{-1}(\inr(\D_s))$. Then $f$ restricted to $\pd\S_0$ is transverse to $p$ and so
there is  a homotopy of $f$, relative to $\pd\S_0$, to a map which is transverse to
$p$ on $\S_0$. Therefore we may assume that $f$ gives rise to a sketch $\G$ 
on $\S(q)$.
Extend the orientation of the $\hat{d_i}$'s to an orientation of regions and vertices of
$\S(q)$. (This can be done consistently if $\S(q)$ is orientable.) Let $f^l(d_i)=u_i$ and 
define $f^l$ on vertex corners as in the definition of a picture. Then $(f,\G)$ is a picture
on $\S$ over $G$ with boundary labels $u_1,\ldots ,u_k$, as required. 
\begin{corol}\label{pic_qpe}
Let $(\mbf h,\mbf n,\mbf t,\mbf p)$
be a positive 4--partition, 
$\cL$ a consistent system of parameters and $\mbf z$ an
element of $(H^\Lm,\cL,\mbf n,\mbf h)$. Then
the triple $(\mbf n,\mbf t,\mbf p)\in \cL$--genus$(\mbf z)$ if and only if
there exists a solution $\al$ to $\cL$ and a minimalistic picture $\G$ on a surface 
of type $(\mbf n,\mbf t,\mbf p)$ with partitioned boundary labels list
 $(\hat\al(z_1),\ldots,\hat\al(z_n))$. 
\end{corol}
\medskip

\noindent
{\em Proof.} 
Given a solution $(\phi,\al)$ to $Q=Q(\mbf z,\cL,\mbf n,\mbf t,\mbf p)$
we have $\phi(q(\xi_j,n_j,t_j,p_j))=1$ in $G$, for $j=1,\ldots ,k$.
It follows, from Theorem \ref{fundpic}, that there exists a picture
$\G_j$ over $G$ on a surface $\S_j$ of genus $t_j+p_j/2$, with $n_j$ boundary
components labelled $\phi(d_{\xi_j+1}),\ldots ,\phi(d_{\xi_j+n_j}).$ 
If $\G_j$ exists then there also exists a minimalistic picture satisfying the same
conditions, so we may assume that $\G_j$ is minimalistic.
As $(\phi,\al)$
is a solution to $Q$, $\phi(d_{\xi_j+i})=\hat\al\bt(d_{\xi_j+i})=\hat\al (z_{\mu_j+i})$,
for all $i,j$. The picture $\G=\G_1\cup\cdots\cup\G_k$ on $\S_1\dcup\cdots\dcup\S_k$
therefore satisfies the given requirements.

Conversely, given a solution $\al$ to $\cL$ and a picture $\G$ on a surface 
$\S$ as described, suppose that $\S$ has connected components $\S_1,\ldots,\S_k$, where $\S_j$ has
genus $t_j+p_j/2$ and $n_j$ boundary components labelled 
$\hat\al(z_{\mu_j+1}),\ldots,\hat\al(z_{\mu_j+n_j}).$
Let $\G_j=\G\cap \S_j$ and $q_j=q(\xi_j,n_j,t_j,p_j)$ and
define $\bt_j:F(L_D(q_j))\maps H$ by $\bt_j(d_{\xi_j+i})=\hat\al(z_{\mu_j+i})$, 
for $j=1,\ldots ,k$ and $i=1,\ldots, n_j$. Then,  from Theorem \ref{fundpic},
$\G_j$ gives rise to a solution $\phi_j$ to $(q_j=1,\bt_j)$. Define $\phi :
F(L(\mbf q))\maps H$ by $\phi(x_i)=\phi_j(x_i)$, for 
$\xi_j+1\le i \le \xi_j+n_j+2t_j+p_j$ and $\phi(d_{\xi_j+i})=\hat\al(z_{\mu_j+i})$,
for $1\le i\le n_j$. Then $(\phi,\al)$ is a solution to $Q$. 

\begin{defn}\label{pic_soln}
Let $(\mbf h,\mbf n,\mbf t,\mbf p)$
be a positive 4--partition, 
$\cL$ a consistent system of parameters and $\mbf z$ an
element of $(H^\Lm,\cL,\mbf n,\mbf h)$. 
Let $\al$ be a solution to $\cL$ and $\G$ a picture on a surface 
of type $(\mbf n,\mbf t,\mbf p)$ with partitioned boundary labels list
 $(\hat\al(z_1),\ldots,\hat\al(z_n))$. Then we say that $(\al,\G)$ {\em
   corresponds} to the  solution $(\phi,\al)$ of $Q(\mbf z,\cL,\mbf n,\mbf t,\mbf p)$
constructed in Corollary \ref{pic_qpe}. The pair $(\al,\G)$ is 
called {\em exponent--minimal} if it corresponds to an 
exponent--minimal solution 
$(\phi,\al)$ to $Q(\mbf z,\cL,\mbf n,\mbf t,\mbf p)$.  
\end{defn}

When there is no ambiguity we refer to $\G$, rather than $(\al,\G)$, as exponent--minimal.
Observe that if an exponent--minimal picture with prime labels $\mbf z$ exists then 
a reduced (or efficient or minimal) exponent--minimal picture with the same prime
labels also exists, because all pictures in an equivalence class and in $T$  have the same 
boundary labels.
\section{Corridors}\label{corridor}

Let $(\mbf h,\mbf n,\mbf t,\mbf p)$ be a positive 4--partition, $\cL$
a consistent system of parameters, $\mbf z$ a special element of
$(H^\Lm,\cL,\mbf n,\mbf h)$ and $\al$ a solution to $\cL$.  From now
until further notice we assume that $\G$ is a picture on a surface of
type $(\mbf n,\mbf t,\mbf p)$ with  prime
labels $\mbf z$, labelled by $\hat\al(\mbf z)$ and with boundary partition 
$\mbf b=(b_1,\ldots ,b_{W_1})$.

\begin{figure}
\begin{center}
  { \includegraphics[scale=0.4]{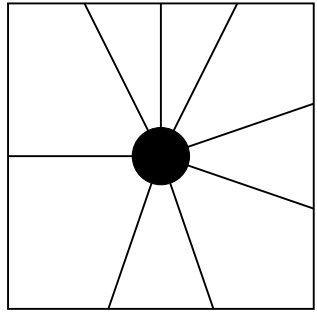}
    }\\
\caption{A picture of type square--$1$ with base $D(2,3,1,2)$}\label{square-1}
\end{center}
\end{figure}

Let $D$ be the disk $[a,a+1]\+[0,1]\subset\RR^2$, where $a\ge 0$ is an  
integer, and let $x_0,y_0,x_1$
and $y_1$ be non--negative integers with $x_i>0$, $i=1,2$. Define a
set $D(x_0,x_1,y_0,y_1)$ of {\em distinguished} points of $D$ to be
\be
\item $(a,i/(y_0+1))$, for $i=1,\ldots,y_0$;
\item $(a+1,i/(y_1+1))$, for $i=1,\ldots,y_1$;
\item $(a+(i/(x_0+1)),0)$, for $i=1,\ldots,x_0$ and
\item $(a+(i/(x_1+1)),1)$, for $i=1,\ldots,x_1$.  \ee A picture $\G_D$ over $G$
  on $D$ is said to have {\em type square--}$0$ (with {\em base}
  $D(1,1,0,0)$) if $\G_D$ has no vertices and one arc
  $\{1/2\}\+[0,1]$.  A picture $\G_D$ on $D$ is said to have {\em type
    square--}$1$ with {\em base} $D(x_0,x_1,y_0,y_1)$ if $x_0,x_1,y_0$
  and $y_1$ are non--negative integers such that $x_0,x_1>0$,
  $x_0+x_1+y_0+y_1=ml$, $\G_D$ has one vertex $v$ and $\G_D$ has $ml$
  arcs each beginning at a point of $\pd v$ and ending at one
  of the distinguished points of $D(x_0,x_1,y_0,y_1)$.  A picture of
  type square--$1$ with base $D(2,3,1,2)$ is shown in Figure
  \ref{square-1}.

\begin{figure}
\begin{center}
  { \includegraphics[scale=0.4]{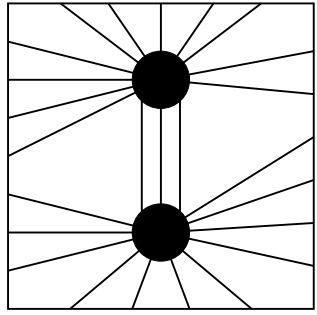}
    }\\
\caption{A picture of type square--$2$ with base $D(4,5,7,6)$}\label{square-2}
\end{center}
\end{figure}

A picture $\G_D$ over $G$ on $D$ is said to have {\em type square--}$2$ with
{\em base} $D(x_0,x_1,y_0,y_1)$ if the following conditions hold.  \be
\item $x_0,x_1,y_0$ and $y_1$ are positive integers such that
  $y_0,y_1\ge 3$ and $x_0+x_1+y_0+y_1\le 2ml-1$.
\item $\G_D$ has $2$ vertices $v_0$ and $v_1$ and
  $ml+(x_0+x_1+y_0+y_1)/2$ arcs. $\G_D$ has one non--empty class of
  arcs joining $v_0$ to $v_1$ and $\G_D$ has $x_i$  arcs
  joining $v_i$ to distinguished points of $D$ on $[a,a+1]\+\{i\}$,
  for $i=1,2$.
\item There are integers $r_i$, $i=0,1$ such that $0<r_i<y_i$ and
  $\G_D$ has $r_i$  arcs joining $v_i$ to distinguished points
  of $\{a+i\}\+[0,1]$ and $y_i-r_i$  arcs joining $v_{1-i}$ to
  distinguished points of $\{a+i\}\+[0,1]$.  
\ee 
A picture of type
  square--$2$ with base $D(4,5,7,6)$, $r_0=3$  and $r_1=2$ is shown in Figure
  \ref{square-2}.
A picture of type square--$j$ is called a {\em
  square--picture} of type $j$, where $j=0,1$ or $2$. An arc $e$ of a square--picture
meeting $[a,a+1]\+\{i\}$, for $i=0$ or $1$ is called a {\em lateral}
arc of the picture. An arc $e$ of a square--picture of type 2 incident
to both vertices is called an {\em internal--lateral} arc.

\begin{figure}
\psfrag{a}{$a$}
\psfrag{1}{$1$}
\psfrag{c0}{$c_0$}
\psfrag{c1}{$c_1$}
\begin{center}
  { \includegraphics[scale=0.4]{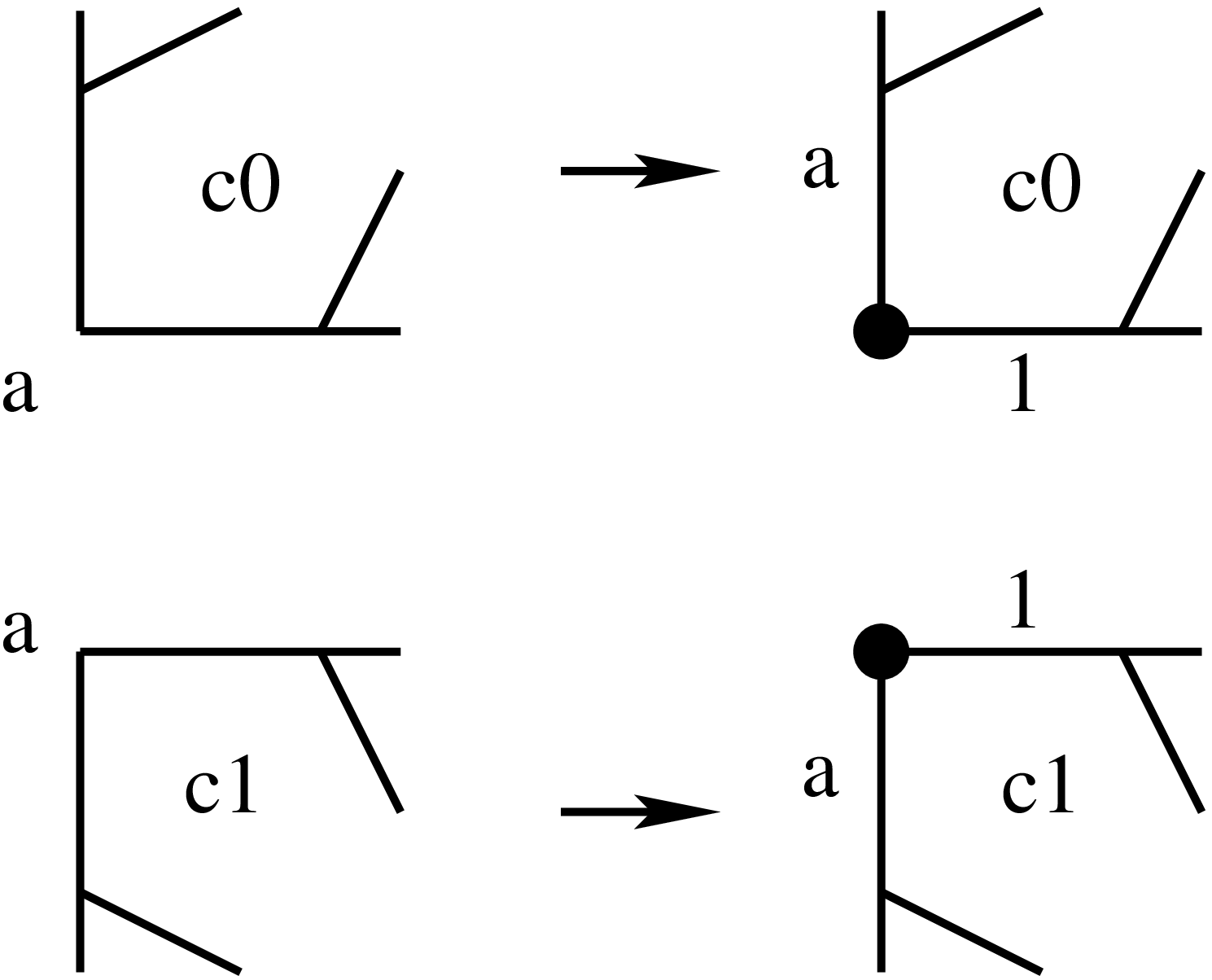}
    }\\
\caption{}\label{left_trivial}
\end{center}
\end{figure}

A square--picture $(f_D,\G_D)$ comes equipped with a labelling
function $f_D^l$, from corners of $\G_D$ to $A\dcup B$. Since we wish
to map square--pictures into arbitrary pictures in such a way that the
intervals $([a+(k/(x_i+1),a+((k+1)/(x_i+1)]\+\{i\})$ correspond to
boundary corners of $\G$ we now refine this labelling function. We
redefine the label of the corner $c_0=(\{a\}\+[0,1/(y_0+1)])\cup
([a,a+(1/(x_0+1))]\+\{0\})$ setting
$f^l_D([a,a+(1/(x_0+1))]\+\{0\})=1$ and $f^l_D(\{a\}\+[0,1/(y_0+1)])=$
the original label of $c_0$. Similarly, we redefine the labelling
function on the corner $c_1=(\{a\}\+[y_0/(y_0+1),1])\cup
([a,a+(1/(x_1+1))]\+\{1\})$ so that the subinterval
$[a,a+(1/(x_1+1))]\+\{1\}$ has trivial label and the defining property
of boundary labels of regions of pictures is preserved (see
Figure \ref{left_trivial}).  
We say such a
labelling function on $\G_D$ is {\em left--trivial}.  If $\G_D$ is of
type $0$ all labels are trivial so any labelling function is
left--trivial. If $\G_D$ is not of type $0$ then, since a regular
neighbourhood of an arc of $\G_D$ maps to a contractible subset of
$Z$, we may achieve a left--trivial labelling by a homotopy of $f_D$
fixing all of $D$ except for the components of $D\bl\G_D$ meeting
$c_0$ and $c_1$.  We shall henceforward assume that square--pictures
have left trivial labelling functions.

If $D$ is the disk $[a,b]\+[0,1]$, where $a,b\in \ZZ$, $0\le a<b$, we define $D_c$ to be
the disk $[c,c+1]\+[0,1]$, for $c\in\ZZ$, $a\le c<b$.  We define a {\em sum} of
pictures of {\em type} 1 to be a picture $\G_D$ on the disk
$D=[a,b]\+[0,1]$, where 
\be
\item $b-a\ge 1$,
\item\label{sum-square} $\G_D\cap D_c$ is a picture of type square--0
  or square--1, for $a\le c\le b-1$, and there is at least one $c$ in
  this range such that $\G_D\cap D_c$ is a picture of type square--1
  and 
\item if $\G_D\cap D_c$ has base $D_c[x_0,x_1,y_0,y_1]$ then $\G_D\cap
  D_{c+1}$ has base $D_{c+1}[x^\prime_0,x^\prime_1,y_1,y^\prime_1]$,
  for $a\le c\le b-2$.  
\ee 

\begin{figure}
\psfrag{a}{$a$}
\psfrag{b}{$b$}
\begin{center}
\includegraphics[scale= 0.4,clip]{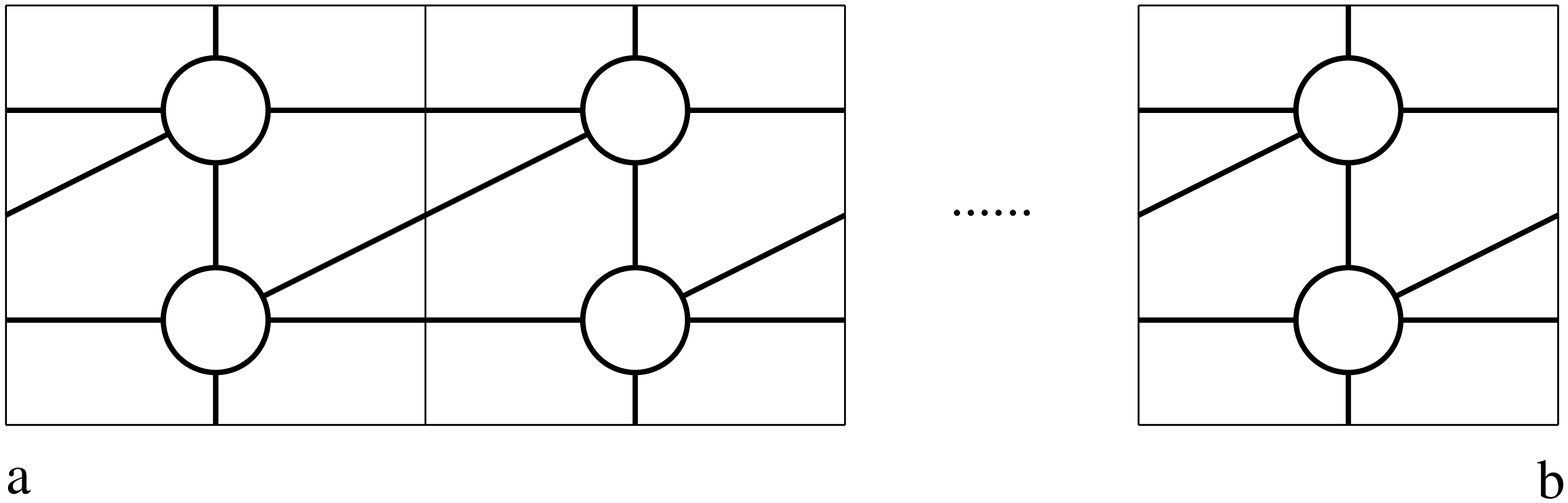} 
\caption{}\label{square_j_2_cond_a}
\end{center}
\end{figure}

\begin{figure}
\psfrag{a}{$a$}
\psfrag{b}{$b$}
\begin{center}
\includegraphics[scale= 0.4,clip]{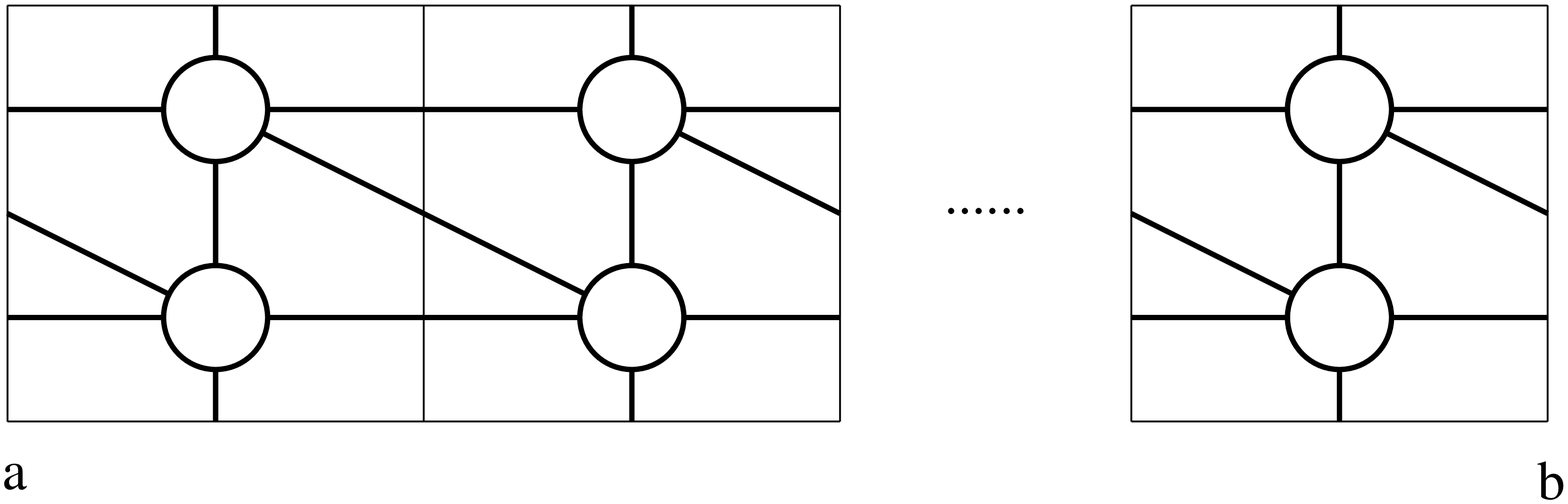} 
\caption{}\label{square_j_2_cond_b}
\end{center}
\end{figure}

For $j=0$ or $2$ we define a sum of pictures of {\em type} $j$
similarly except that in condition (\ref{sum-square}) all pictures
must be of type square--$j$ and if $j=2$ then the following condition
must also hold.  
\be
\setcounter{enumi}{3}
\item\label{square_2_cond} 
If $b-a>1$ then $\G_D$ has the form shown in either Figure \ref{square_j_2_cond_a} or
\ref{square_j_2_cond_b} (where all the edges shown contain at least one arc).
\ee
Given a sum of square--pictures $\G_D$ on
$D=[a,b]\+[0,1]$ we set $\G_D(c)=\G_D\cap D_c$, for $a\le c\le b-1$.  A
sum of square--pictures on $[a,b]\+[0,1]$ is said to have {\em width}
$b-a$. The lateral and internal--lateral arcs of the square--pictures
making up  a
sum of square--pictures are called {\em lateral} and {\em
internal--lateral} arcs, respectively, of the sum.
  
Let $D=[a,a+1]\+[0,1]$ and let $\G_D$ be a square--picture of type $j$
with base $D(x_0,x_1,y_0,y_1)$ on $D$.  Suppose $j>0$.  If $j=1$ let
$v_0$ be the vertex of $\G_D$. If $j=2$ let $v_i$ be the vertex of
$\G_D$ incident to $[a,a+1]\+\{i\}$, for $i=0,1$. Let $a_0$ be the arc
joining $v_0$ to $(a+1/(x_0+1),0)$. Travelling clockwise around $\pd
v_0$, starting from the endpoint of $a_0$, let $a_1$ be the first arc
encountered which does not meet $\{a\}\+[0,1]$. If $j=1$ then $a_1$
connects $v_0$ to $(a+1/(x_1+1),1)$. If $j=2$ then $a_1$ connects $v_0$
to $v_1$. In the latter case, travelling clockwise around $\pd v_1$,
starting from $a_1$, let $a_2$ be the first arc encountered which does
not meet $\{a\}\+[0,1]$: so $a_2$ joins $v_1$ to $(a+1/(x_1+1),1)$.  The
label of $v_i$, read from an appropriate starting point in a clockwise
direction, is $r^{\de_i m}$, where $\de_i=\pm 1$, for $0\le i\le j-1$.
The label $u_i$ on the subinterval of $\pd v_i$, read clockwise, from
$a_i$ to $a_{i+1}$ is then a subword of $r^{\de_i m}$, for $0\le i\le
j-1$.  Let $o_i$, $\rho_i$ and $\s_i$ be
non--negative integers such that $\rho_i<m$, $o_i<l(r)$, $\s_i<l(r)$
and $u_i=\t(r^{\de_i},o_i)(r^{\de_i})^{\rho_i}\io(r^{\de_i},\s_i)$, 
for $0\le i\le j-1$.   In
cases $j=0$ or $1$ define $\de_i=o_i=\rho_i=\s_i=0$, for $j\le i\le 1$.
\begin{defn}
  The {\em left interior marking} of $\G_D$ is the sequence
  \[\de_0,o_0,\rho_0,\s_0,\de_1,o_1,\rho_1,\s_1.\]
\end{defn}

\begin{figure}
\psfrag{a0}{$a_0$}
\psfrag{a1}{$a_1$}
\psfrag{v0}{$v_0$}
\begin{center}
  { \includegraphics[scale=0.4]{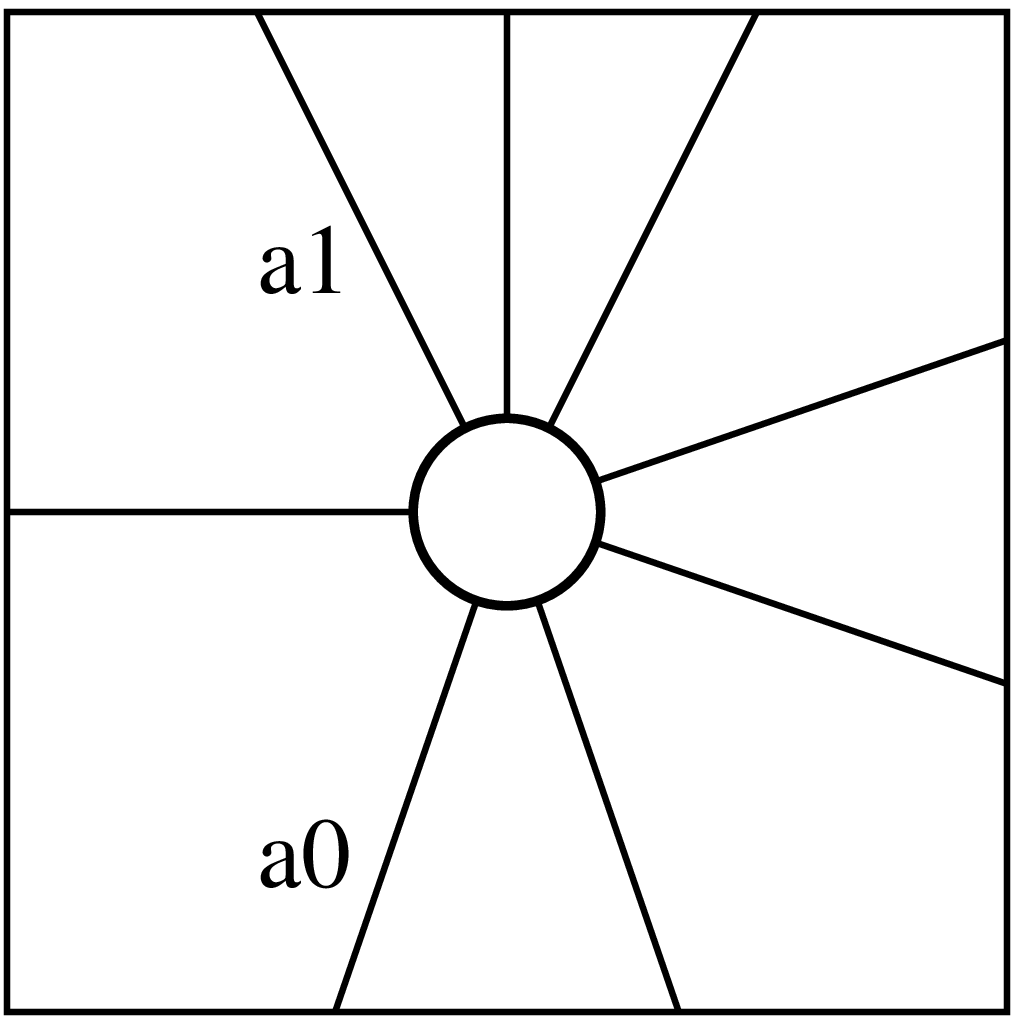}
    }\\
\caption{}\label{square-1_exx}
\end{center}
\end{figure}

\begin{figure}
\psfrag{a0}{$a_0$}
\psfrag{a1}{$a_1$}
\psfrag{v0}{$v_0$}
\psfrag{a2}{$a_2$}
\psfrag{v1}{$v_1$}
\begin{center}
  { \includegraphics[scale=0.4]{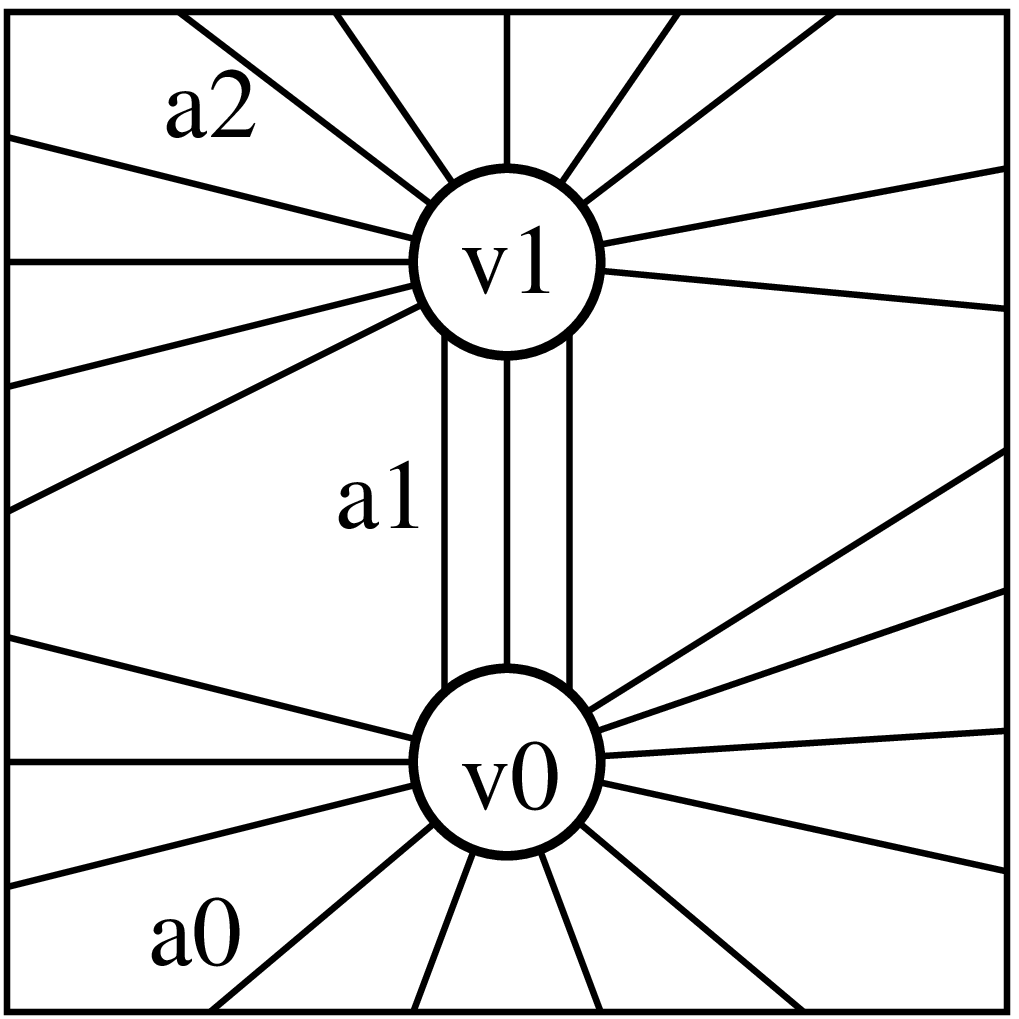}
    }\\
\caption{}\label{square-2_exx}
\end{center}
\end{figure}

For example suppose that, in Figure \ref{square-1_exx},
$r=ab$, $m=4$ and the label on $v_0$ read clockwise from $a_0$ is $(ab)^4$.
Then $\G_D$ has left interior marking $1,0,1,0,0,0,0,0$. 
Again, in Figure \ref{square-2_exx} suppose $r=ab$, $m=7$, the 
label of $v_0$ read from $a_0$ is $(ba)^7$ and the label of 
$v_1$ read from $a_1$ is $(a^{-1}b^{-1})^7$. Then $\G_D$ has left interior
marking $1,1,1,1,-1,1,2,0$.

Let $D=[a,b]\+[0,1]$ and suppose that there exists a map $\mu:D\maps
\S$ and boundary intervals $P_0$ and $P_1$ of $\G$ satisfying the following
properties.  
\be
\item $\mu$ is injective on $D\bl (\{0,1\}\+[0,1])$;
\item either $\mu$ is injective on $\{0,1\}\+[0,1])$ or
  $\mu(a,t)=\mu(b,t)$ or $\mu(a,t)=\mu(b,1-t)$, for all $t\in [0,1]$;
\item $\pd \S\cap \mu(D)=\mu([a,b]\+\{0,1\})$ and
  $\mu([a,b]\+\{i\})\subset P_i$, for $i=0,1$;
\item $\mu(\{c\}\+[0,1])$ is transverse to $\G$ and meets no vertex of
  $\G$, for $c\in\ZZ$ with $a\le c\le b$.  
\ee 
Then we call $\mu$ a {\em binding map}
(on $D$) between $P_0$ and $P_1$.
Now suppose that $P_i$ is not a partisan
boundary interval, for $i=0,1$, and that no arc of $\G$ contained in
$\mu(D)$ is an $H^\Lm$--arc. Let $(f_D,\G_D)$ be a sum of
square--pictures of type $j$ on $D$, with left--trivial labelling
function $f^l_D$, and let $\mu(D)=B$.  Let $c_0,c_1$ be the left--hand
corners of $D$ defined above and in addition let
$c_2=(\{b\}\+[0,1/(y_1+1)])\cup ([b-1+(x_0/(x_0+1)),b]\+\{0\})$ and
$c_3=(\{b\}\+[y_1/(y_1+1),1])\cup ([b-1+(x_1/(x_1+1)),b]\+\{1\})$.
Recall that $f$ is the map from $\S$ to $Z$ described in the
definition of picture.

Assume that $\mu(\G_D)=\G\cap B$ and 
\be
\item $f^l(\mu(c))=f^l_D(c)$, for all corners $c$ of $\G_D$ except
  $c_0,c_1, c_2$ and $c_3$,
\item $f(\mu([a,a+(1/(x_i+1))]\+\{i\}))$ is nullhomotopic, for
  $i=0,1$, and 
\item $d_i=\mu([b-1+(x_i/(x_i+1)),b]\+\{i\})\subset c^\prime_i$, where
  $c^\prime_i$ is a boundary corner of $\G$ and $f(c_i^\prime\bl d_i)$
  is nullhomotopic, for $i=0,1$.
\ee 
Then $C=(D,\G_D,\mu)$ is called a $j${\em --corridor} of {\em
width} $b-a$ with binding map $\mu$ {\em binding} $P_0$ and
$P_1$. In this case $\G_B=\G\cap \mu(D)$ is a picture on $B$ and we
call $B$ the {\em image} of $C$. 
The images $\mu(x)$ of arcs $x$ of the   picture $\G_D$ are called
{\em arcs} of $C$. 
If $\mu(a,t)=\mu(b,t)$ or $\mu(a,t)=\mu(b,1-t)$, for all $t\in [0,1]$,
then all arcs of $C$ are arcs of $\G$. Otherwise, given  an arc 
$x$ of $\G_D$, the arc $\mu(x)$ of
$C$ is an arc of $\G$ if and only if $x$ does not meet
$\{a,b\}\+[0,1]$. Nonetheless it is convenient to refer to $\mu(x)$ as
an arc of $C$, for all arcs $x$ of $\G_D$. 
A
corridor of width 1 is called a {\em compartment}. The {\em left
  interior marking of rank} $s$ of the corridor $C$ is defined to be
the left interior marking of $\G_D(s)$, for all $s$ such that $a\le
s\le b-1$.

If $e$ is a lateral arc, or internal--lateral arc, of $\G_D$
then $\mu(e)$ is called a {\em lateral} arc, or {\em
  internal--lateral} arc, of the corridor $(D,\G_D,\mu)$.  If
$C=(D,\G_D,\mu)$ and $C^\prime=(D^\prime,\G_D^\prime,\mu^\prime)$ are
corridors with identical sets of lateral arcs then we say that $C$ and
$C^\prime$ are {\em equivalent} and this is an equivalence relation on
the set of corridors. If $C$ and $C^\prime$ are equivalent then
$\mu^\prime$ is homotopic to a map $\mu^{\prime\prime}$ such that
$\mu^{\prime\prime}(D^\prime)=\mu(D)$ and
$(D^\prime,\G_D^\prime,\mu^{\prime\prime})$ is a corridor equivalent
to $C$. If $C^\prime=(D^\prime,\G_D^\prime,\mu^\prime)$ is a corridor
such that $\{$lateral arcs of
$C^\prime\}\subseteq\{$lateral arcs of $C\}$
then $C^\prime$ is called a {\em sub--corridor} of $C$.  If $a\le
c<d\le b$ and $D^\prime = [c,d]\+[0,1]$ then the map
$\mu^\prime=\mu|_{D^\prime}$ is called the binding map {\em induced}
from $\mu$ on $D^\prime$ and the corridor $(D^\prime,
\G_{D^\prime},\mu^\prime)$ is called the sub--corridor {\em induced}
from $(D,\G_D,\mu)$, denoted $C[c,d]$ or $C[c]$, if $c+1=d$.  
If $C$ is corridor and is not a sub--corridor of
any corridor of greater width then we say that $C$ is a {\em maximal}
corridor. Obviously two corridors are equivalent if and only
if each one is a sub--corridor of the other.
If $C=(D,\G_D,\mu)$ and 
$C^\prime=(D^\prime,\G_{D^\prime},\mu^\prime)$ are  corridors, with
images $B=\mu(D)$ and $B^\prime=\mu^\prime(D^\prime)$,  we
call $B\cap B^\prime$ the {\em intersection} of $C$
and $C^\prime$, denoted $C\cap C^\prime$, and if
$B\cap B^\prime\neq \nul$ we shall say that $C$ and
$C^\prime$ {\em overlap}.  

\begin{defn}
  Choose one representative from each equivalence class of maximal
  $j$--corridors of $\G$, where $j=0$ or $1$.
  These corridors and the sub--corridors induced from them are called
  {\em designated} corridors of $\G$ (of type 0 and 1).
\end{defn}
(We defer the definition of designated 2--corridor until later.)  
We may always choose the maximal designated corridor $(D,\G_D,\mu)$
such that 
$D=[0,c]\+[0,1]$, for some integer $c$.  In the absence of any
explicit statement to the contrary it is to be assumed that all
corridors are designated corridors.  

\begin{figure}
\begin{center}
  { \includegraphics[scale=0.5]{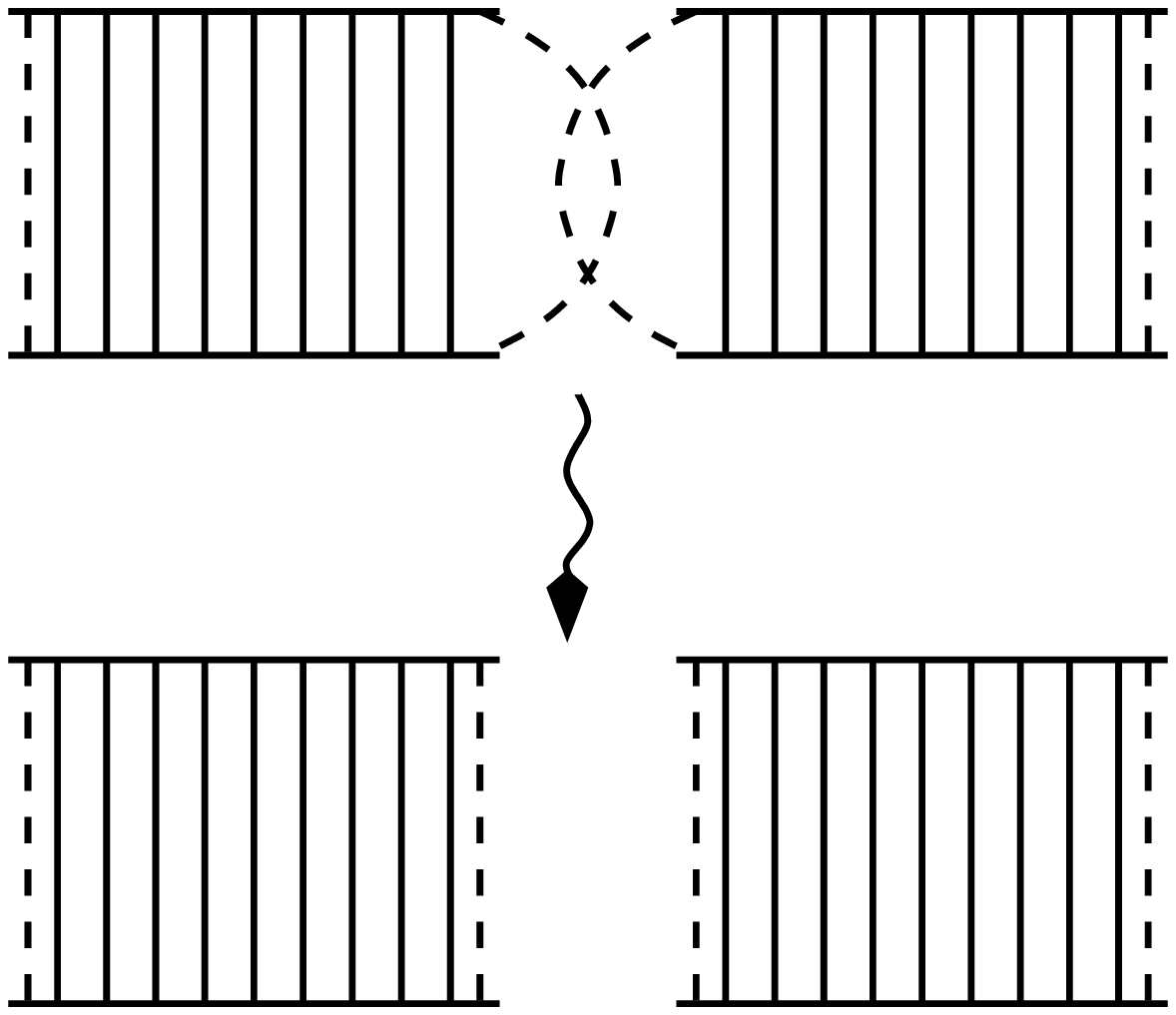}
    }\\
\caption{}\label{Corridor0-i}
\end{center}
\end{figure}

If $(D,\G_D,\mu)$ is a  maximal $0$--corridor then $\G_D$ consists
of arcs from a class of
arcs of type $I$.  If two distinct designated maximal 0--corridors
overlap then we may choose different representatives of the two
equivalence classes of corridors (see Figure \ref{Corridor0-i})
which do not overlap. We assume then that maximal designated
0--corridors do not overlap. If two distinct designated maximal
1--corridors overlap but the intersection contains no arc or vertex
then, as in the type 0 case, we designate different representatives
which do not overlap.  We may assume then that the intersection
contains an arc or vertex. If the intersection contains a vertex then
it contains the lateral arcs incident
to this vertex. After performing a homotopy if necessary we have two
maximal corridors whose intersection contains a compartment. Therefore we may
extend each corridor to contain all compartments of the other. Since
the original corridors are maximal they must be equivalent,
contradicting uniqueness of representatives. If the intersection
contains an arc but no vertex then this arc must be the unique lateral
arc of a compartment of type 0. The same argument applies in this
case.  Hence we may assume that maximal designated 1--corridors do not
overlap.

\begin{figure}
\psfrag{x0}{$x_0$}
\psfrag{x1}{$x_1$}
\psfrag{u0}{$u_0$}
\psfrag{u1}{$u_1$}
\psfrag{e0}{$e_0$}
\psfrag{e1}{$e_1$}
\psfrag{P0}{$P_0$}
\psfrag{P1}{$P_1$}
\begin{center}
  { \includegraphics[scale=0.4]{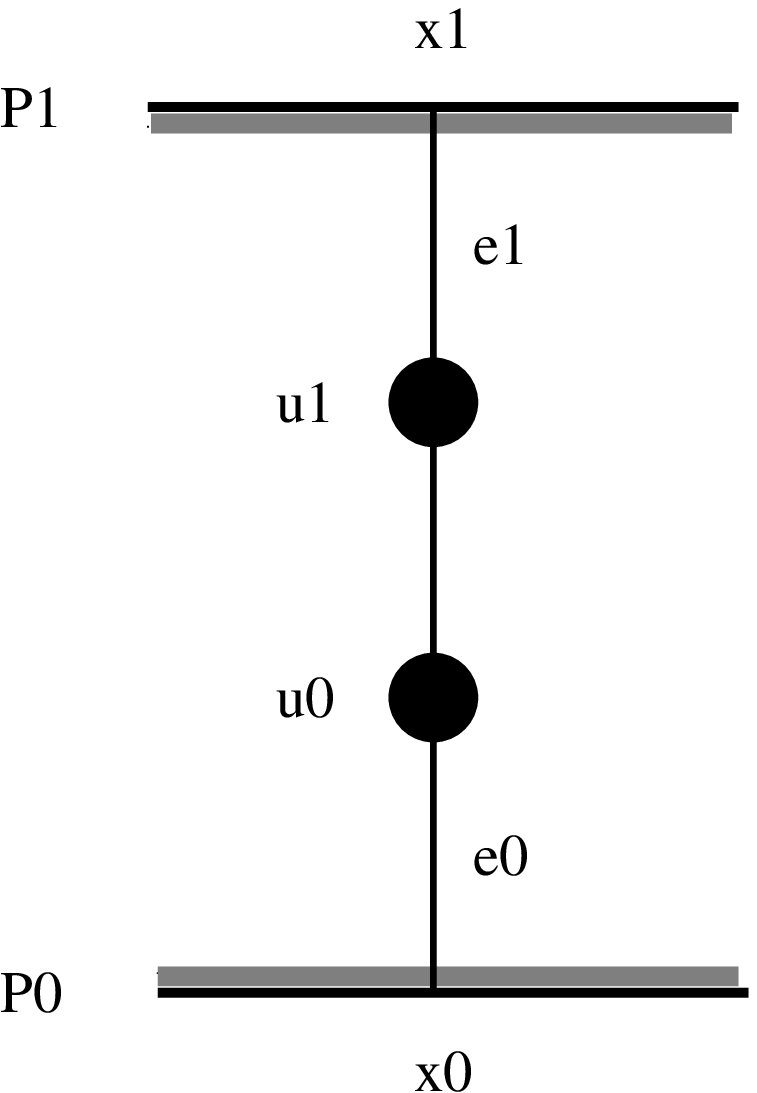}
    }\\
\caption{}\label{2corridor-1-2}
\end{center}
\end{figure}

\begin{figure}
\psfrag{v1}{$v_1$}
\psfrag{u0}{$u_0$}
\psfrag{u1}{$u_1$}
\begin{center}
  \mbox{
\subfigure[$u_1\in C^\prime \lbrack j-1\rbrack$\label{diag1-2-1}]
{ \includegraphics[scale= 0.4,clip]{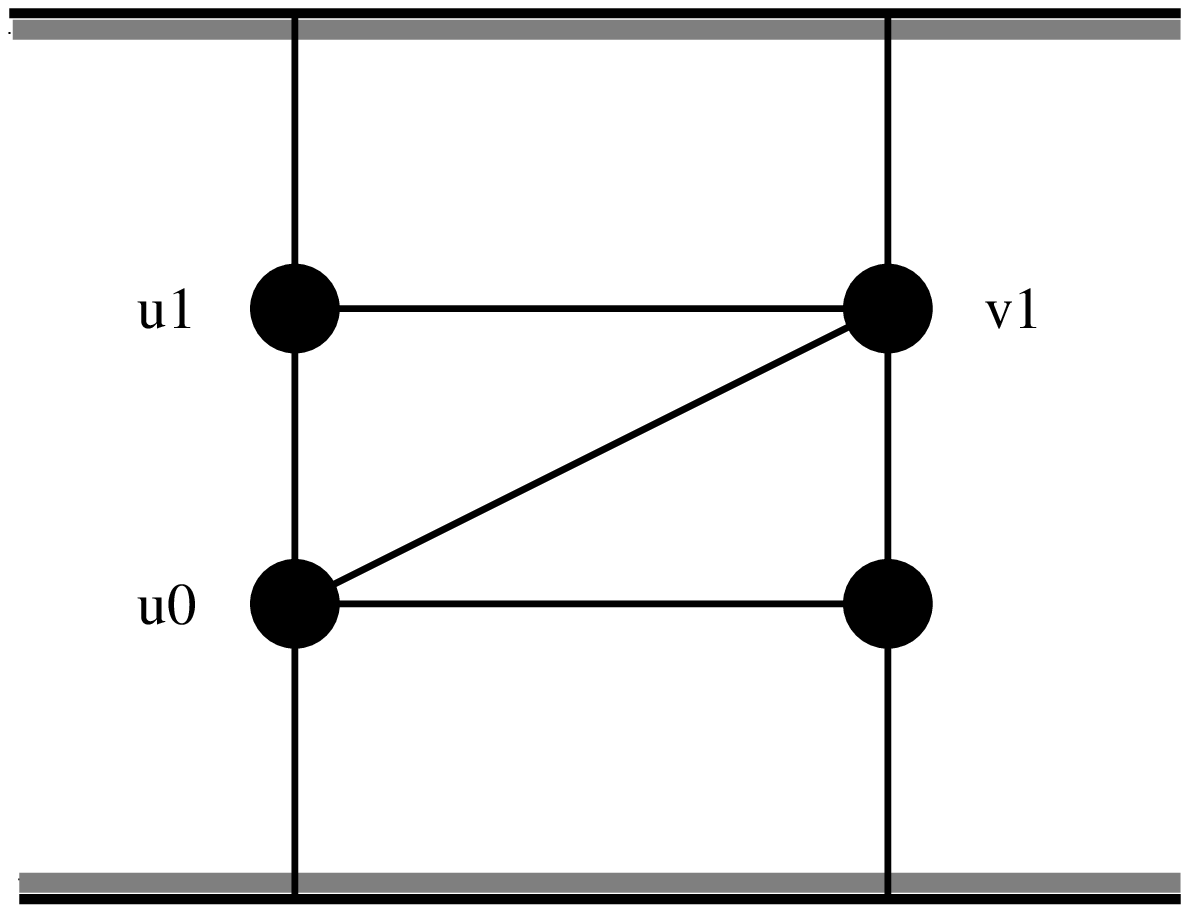} }\qquad
\subfigure[$u_1\in C^\prime  \lbrack j+1 \rbrack$ ]  
{ \includegraphics[scale=0.4,clip]{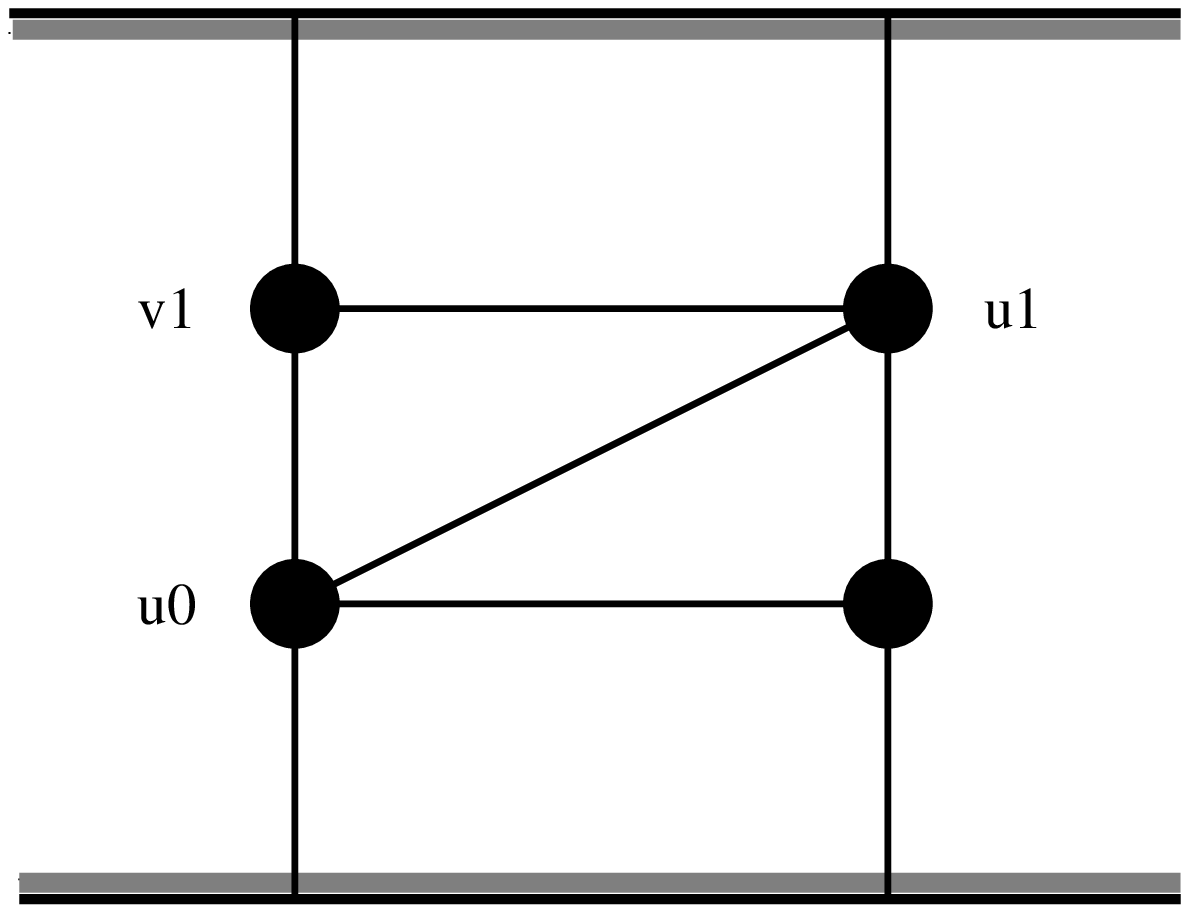} } 
}
\caption{}\label{diag1-2}
\end{center}
\end{figure}

Maximal 2--corridors require further analysis. To begin with we
shall, from now on, call the 2--corridors defined above {\em
  simple}
2--corridors. For ease of reference we shall also refer to corridors
of type 0 and 1 as {\em simple}: so simple $j$--corridors and
$j$--corridors are the same thing, for $j=0,1$. Let $C=(D,\G_D,\mu)$ be
a simple 2--corridor and let $(B,\G_B)=(\mu(D),\G\cap\mu(D))$. Then the
vertices of $\G_B$ together with those arcs of $\G$ whose closure is
contained in $B$ map to a subgraph $\bar C$ of $\bar \G$ on $\S$,
which we call the {\em graph} of $C$. The images of lateral and
internal--lateral arcs of $C$ are called {\em internal} and {\em
  internal--lateral} edges, respectively, of $\bar C$. Let $C^\prime
=(D^\prime,\G_{D^\prime},\mu^\prime)$ be a simple 2--corridor.  Then $C$ is a
sub--corridor of $C^\prime$ if and only if $\bar C\subset \bar
C^\prime$.  Moreover, if $\mu(D)\cap\mu^\prime(D^\prime)$
contains an arc $a$ then every arc equivalent to $a$ also belongs to
$\mu(D)\cap\mu^\prime(D^\prime)$. Hence
$\mu(D)\cap\mu^\prime(D^\prime)$ contains an arc if and only
if $\bar C\cap \bar C^\prime$ contains and edge.  Suppose now that $C$
and $C^\prime$ both bind boundary intervals $P_0$ and $P_1$. If $\bar
C\subset \bar C^\prime$ and both $C$ and $C^\prime$ have width $1$
then clearly $C$ and $C^\prime$ are equivalent. Suppose next that $\bar
C\subset \bar C^\prime$, the
width of $C$ is 1 and the width of $C^\prime$ is at least 2.
Assume further that the graph $\bar C$ meets $P_i$ at a point $x_i$,
let $e_i$ be the edge of $\bar C$ incident to $x_i$ and let $u_i$ be
the vertex of $\bar C$ incident to $e_i$, for $i=0,1$. (See Figure
\ref{2corridor-1-2}.)
As $\bar C\subset \bar C^\prime$ the vertices $u_i$ also belong to
$\bar C^\prime$. There are two possibilities.  \be
\item\label{2corr-align} The internal--lateral edge of $\bar C$
  joining $u_0$ and $u_1$ is an internal--lateral edge of $\bar
  C^\prime$. In this case $C$ is equivalent to a compartment of
  $C^\prime$. 
\item\label{2corr-diag} The internal--lateral edge of $\bar C$ joining
  $u_0$ and $u_1$ is not an internal--lateral edge of $\bar C^\prime$.
  Therefore either $u_0$ or $u_1$ is joined to a vertex of $\bar
  C^\prime$ by an internal--lateral edge of $\bar C^\prime$ which is
  not an internal--lateral edge of $\bar C$.  Suppose that $u_0$ is
  joined to vertex $v_1$ of $\bar C^\prime$ by an internal--lateral
  edge of $\bar C^\prime$ which is the image of a class of arcs in
  a compartment $C^\prime[j]$ of $C^\prime$. Then $u_1$ is a vertex of
  compartment $C^\prime[j-1]$ or $C^\prime[j+1]$ 
  of $C^\prime$
  (see Figure \ref{diag1-2}).  In this case $C$ is not equivalent to
  any compartment of $C^\prime$. The case $u_1$ joined to a vertex of
  $\bar C^\prime$ by an internal--lateral edge of $\bar C^\prime$ is
  similar.  
\ee

\begin{figure}
\psfrag{a}{$a$}
\psfrag{b}{$b$}
\psfrag{c}{$c$}
\psfrag{d}{$d$}
\psfrag{e}{$e$}
\begin{center}
  { \includegraphics[scale=0.6]{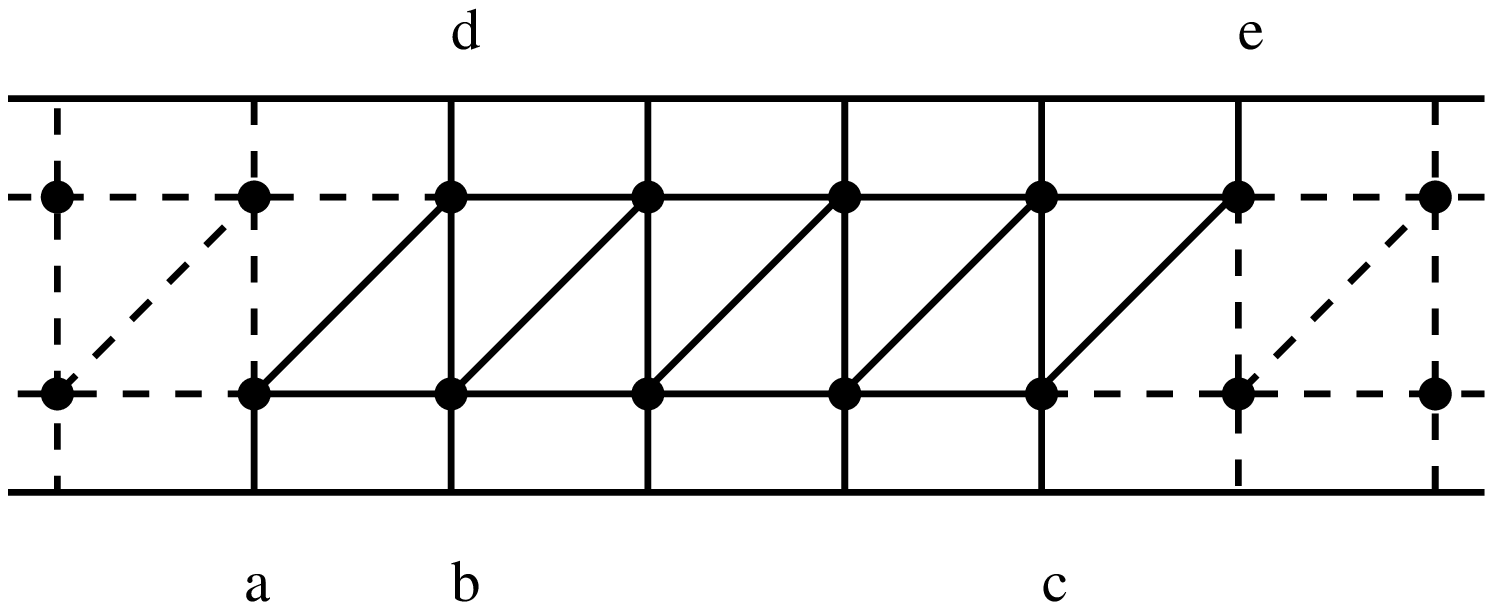}
    }\\
\caption{}\label{sub-diagonal}
\end{center}
\end{figure}

In general if  $\bar C\subset \bar C^\prime$ then  
without loss of generality we may assume that
$D=[0,c]\+[0,1]$ and $D^\prime=[0,c^\prime]\+[0,1]$, for some integers
$c$ and $c^\prime$. 
If the compartment $C[0]$ is equivalent to a compartment of $C^\prime$ (as in
case (\ref{2corr-align}) above) then $C$ is equivalent to a
sub--corridor of $C^\prime$.  If, on the other hand, $C[0]$ is not
equivalent to a compartment of $C^\prime$ (as in case (\ref{2corr-diag})
above) then no compartment of $C$ is equivalent to a compartment of
$C^\prime$. (This follows by induction starting with the case (\ref{2corr-diag}) above.)
In this case $C$ is a sub--corridor of $C^\prime$ but equivalent
to no sub--corridor of $C^\prime$ and we say that $C$ is {\em diagonal} to $C^\prime$. By
definition, if $C$ has width at least 2 and is diagonal to $C^\prime$
then $C^\prime$ has a sub--corridor which is diagonal to $C$.  For
example, consider Figure \ref{sub-diagonal}, which shows the graphs of
two corridors $C\subset C^\prime$: the edges of $\bar C$ are shown solid,
whilst those of $\bar C^\prime\bl\bar C$ are dashed. The
internal--lateral edges of $\bar C$ are shown diagonal whilst those of $\bar
C^\prime$ are vertical. Thus the left hand compartment of $C$ meets
the boundary at $a$ and $d$ whilst the right hand compartment meets it at $c$ and
$e$.  The sub--corridor of $C^\prime$ between $b$ and $c$ is diagonal
to $C$.
We say that simple 2--corridors $C_0$ and $C_1$, both binding boundary intervals $P_0$ and 
$P_1$ and both  of width $c$, 
are {\em diagonally--equivalent} if there is a sub--corridor of $C_i$ of width
$c-1$ which is diagonal to $C_{1-i}$, for $i=0,1$.

If $C_2$ is a corridor of width $2$ then there is a unique equivalence class of 
compartments diagonal to $C_2$. This is a consequence of condition 
(\ref{square_2_cond}) for sums of pictures of type $2$. 
It follows that if $C_0$ is a corridor diagonally--equivalent to a 
corridor $C_1$ and if $C_1$ is 
diagonally--equivalent to a corridor  $C_2$ then $C_0$ and $C_2$ are equivalent.
In particular if $C_0$ 
and $C_1$ are corridors
diagonally--equivalent to a corridor $C_3$ then $C_0$ and $C_1$ are equivalent.
\begin{lemma}\label{2corr-maxint}
Let $C
$ and $C^\prime
$ 
be maximal simple 2--corridors. 
\be
\item\label{2corr-maxint-1} If a simple sub--corridor of $C$ 
is equivalent to a sub--corridor of $C^\prime$
then $C$ and $C^\prime$ are equivalent.
\item\label{2corr-maxint-2} If a simple sub--corridor of $C$ is diagonal
to $C^\prime$ then $C$ and $C^\prime$ are diagonally--equivalent.
\ee
\end{lemma}
\medskip

\noindent
{\em Proof.}
\be
\item
If some sub--corridor of $C$ is 
equivalent to a
sub--corridor of $C^\prime$ then there is a corridor $C^\pprime$
containing all compartments of $C$ and $C^\prime$. 
Then both $C$ and $C^\prime$ are sub--corridors  of $C^\pprime$
and, as $C$ and $C^\prime$ are maximal, $C^\pprime$ is a sub--corridor
of them both. Hence $C$ and $C^\prime$ are equivalent to $C^\pprime$.
\item
Assume that 
some sub--corridor
$C_s$ of $C$ is diagonal to $C^\prime$. Choose $p<q$ such that 
$C_s=C[p,q]$ is  maximal
in this respect and write $C_s=(D_s,\G_s,\mu_s)$. As $C_s$ is maximal among
sub--corridors of $C$ diagonal to $C^\prime$, either there is a vertex
of $\mu(\G_D({p-1}))$ which is not a vertex of $\mu^\prime(\G_{D^\prime})$
or $p=0$.  Similarly, either there is a vertex of $\mu(\G_D(q))$ which is
not a vertex of $\mu^\prime(\G_{D^\prime})$ or $q=c$. If both $p=0$ and
$q=c$ then $C$ is a sub--corridor of $C^\prime$ and the width of $C^\prime$
is necessarily greater than that of $C$, contrary to maximality of $C$.
Thus either $p>0$ or $q<c$. Assume first that $p>0$ and that there is a vertex
of $\mu(\G_D({p-1}))$ which is not a vertex of $\mu^\prime(\G_{D^\prime})$. 
Since the vertices of $\mu(\G_D(p))$ belong to different compartments of $C^\prime$,
it follows that one of the vertices of $\mu(\G_D({p-1}))$ belongs to a compartment
of $C^\prime$ (and is joined to a vertex of $\mu(\G_D(p))$ by an internal--lateral
arc of $C^\prime$).
Hence one vertex of $\mu(\G_D({p-1}))$ is contained in
$\mu^\prime(\G_{D^\prime})$. If $p\ge 2$ it is therefore 
possible to extend $C^\prime$ to a wider corridor by adjoining one vertex of 
$\mu(\G_D({p-1}))$
and one vertex of $\mu(\G_D({p-2}))$. This is contrary to the maximality of $C^\prime$ and
so if $p>0$ then $p=1$. A similar argument shows that if $q<c$ then $q=c-1$.
Now if $p=1$ and $q=c-1$ then it follows that $C^\prime$ is 
a sub--corridor of $C$, so $C$ and $C^\prime$ are equivalent, contrary to hypothesis. 
Hence either
$p=1,q=c$ or $p=0,q=c-1$. It follows that the width of $C_s$ is one less than that
of $C$ and that the compartment of $C$ not in $C_s$ has precisely 
one vertex in $C^\prime$. As $C^\prime$ has a sub--corridor diagonal to $C$ the same
holds on reversing the roles of $C$ and $C^\prime$. Thus 
$C$ and $C^\prime$ are diagonally--equivalent.
 \ee

From the proof of the second part of the above Lemma it follows that if $C$ and
$C^\prime$ are diagonally equivalent then $\bar{C^\prime}$ contains exactly one 
vertex and one lateral edge not in $\bar C$, and vice--versa.

Given that $C=(D,\G_D,\mu)$ and 
$C^\prime=(D^\prime,\G_{D^\prime},\mu^\prime)$ are 
diagonally--equivalent simple 2--corridors we define the sum 
$C + C^\prime$ to be $(D\dcup D^\prime,\G_D\dcup \G_{D^\prime},\mu\dcup\mu^\prime)$
and call this an {\em extended} 2--corridor, {\em extending} $C$ and
$C^\prime$. The set of {\em lateral} arcs of an extended 2--corridor $C+C^\prime$ is the
union $\{$lateral arcs of $C\}\cup \{$lateral arcs of $C^\prime\}$.
We define the set of {\em sub--corridors} of the extended 2--corridor 
$E=C+C^\prime$ to be the union of $\{E\}$ with the sets of 
sub--corridors of $C$ and $C^\prime$.  
(We shall not need to consider extended corridors as sub--corridors of other
extended corridors.)
If $E=C+C^\prime$ is an extended 2--corridor then we define the {\em graph}
of $E$ to be $\bar E=\bar C\cup\bar C^\prime$ and the {\em image} of
$E$ to be $B\cup B^\prime$, where $B$ and $B^\prime$ are the images of
$C$ and $C^\prime$, respectively. By the remark following the Lemma above 
$\bar E$ has one more lateral edge than both $\bar C$ and $\bar C^\prime$. 

We define the set of 2--corridors to be the union of the
set of all simple 2--corridors and the set of all extended 2--corridors.
  We extend the
equivalence relation defined on simple 2--corridors by defining 2--corridors 
$C$ and $C^\prime$ to be {\em equivalent} if they have the same lateral arcs. 
As graphs of 
simple $2$--corridors have an even number of lateral edges, graphs of 
extended $2$--corridors
have an odd number of lateral edges and so no simple $2$--corridor can be 
equivalent to an extended $2$--corridor.
A 
2--corridor $C$ is {\em maximal} in the set of 2--corridors if 
$C$ is equivalent to $C^\prime$, for all 2--corridors $C^\prime$ such that
$\{$lateral edges of $C\}\subset \{$lateral edges of $C^\prime\}$. 
If $E$ and $E^\prime$ are 2--corridors with images $B$ and
$B^\prime$, respectively, then we say that $E$ and $E^\prime$ have  
{\em intersection} $B\cap B^\prime$  and that they {\em
  overlap}
if $B\cap B^\prime\neq\nul$. 

\begin{lemma}\label{maxl-ext}
Let $E=C_0+C_1$ be an extended 2--corridor 
where $C_0$ and $C_1$ are diagonally--equivalent simple $2$--corridors.
Then the following hold. 
\be
\item\label{maxl-ext-1} If $E$ is maximal as a $2$--corridor then 
$C_0$ and $C_1$  are maximal simple 2--corridors.
\item\label{maxl-ext-2} if $E=C_0+C_1^\prime$ then $C_1$ is equivalent to $C_1^\prime$.
\ee
\end{lemma}
\medskip

\noindent
{\em Proof.}
Let $C_0$ and $C_1$ have width $c$.
\be
\item 
Let $C^\prime_i$ be a maximal simple $2$--corridor containing $C_i$, $i=0,1$.
Then $C_0$ is  a sub--corridor of $C_0^\prime$ diagonal to $C^\prime_1$. From 
Lemma \ref{2corr-maxint}.(\ref{2corr-maxint-2}), 
$C_0^\prime$ and $C^\prime_1$ are diagonally--equivalent.
Hence we may form the sum $C_0^\prime+C_1^\prime$. As $E$ is maximal it is 
equivalent to $C_0^\prime+C_1^\prime$. If $e$ is a lateral edge of 
$\bar C^\prime_0\backslash \bar C_0$ this implies that $e$ is a lateral edge of 
$\bar C_1\backslash \bar C_0$. 
As $C_0$ and $C_1$ are diagonally equivalent there is only one
such edge. However, if $\bar C_0^\prime$ is wider than $\bar C_0$ then it must contain
at least $2$ more lateral edges than $\bar C_0$. Hence $C_0^\prime$ and $C_0$ are
equivalent. A similar argument applies to $C_1$.
\item There is a unique vertex $v$ of $\bar E$ which belongs to $\bar C_1$ but not 
$\bar C_0$. Therefore $v$ is the unique vertex of $\bar E$ 
which belongs to  $\bar C_1^\prime$  but not  $\bar C_0$. Hence $C_1$ and 
$C_1^\prime$ must have the same vertices which implies that they are equivalent.
\ee

\begin{lemma}\label{2+corr-maxint}
Let $E$ and $E^\prime$ be maximal 2--corridors. 
\be
\item \label{2+corr-maxint-1}If a simple sub--corridor of $E$ is 
equivalent to a simple sub--corridor of $E^\prime$ or
\item\label{2+corr-maxint-2} a simple sub--corridor of $E$ is diagonal
to a simple sub--corridor of $E^\prime$
\ee
then $E$ and $E^\prime$ are equivalent 2--corridors.
\end{lemma}
\medskip

\noindent
{\em Proof.}
If $E$ and $E^\prime$ are simple 2--corridors then, from Lemma \ref{2corr-maxint},
(\ref{2+corr-maxint-1}) follows immediately and (\ref{2+corr-maxint-2}) is not possible.
If $E$ is a simple 2--corridor and $E^\prime$ an extended 2--corridor then it follows from
Lemma \ref{maxl-ext} and Lemma \ref{2corr-maxint} with either (\ref{2+corr-maxint-1}) or (\ref{2+corr-maxint-2}), 
that $E^\prime$ is equivalent to $E+C^\prime$, where $C^\prime$ is a simple sub--corridor of $E^\prime$, 
contrary to maximality
of $E$. If $E^\prime$ is a simple 2--corridor and $E$ an extended 2--corridor then again a 
contradiction arises.
Hence we assume that both $E$ and $E^\prime$ are extended 2--corridors,
say $E=C_0+C_1$ and $E^\prime=C_0^\prime+C_1^\prime$.
\be
\item If a  simple sub--corridor of $E$ is equivalent to a  simple sub--corridor of $E^\prime$ then,
for some $i,j\in\{0,1\}$ a sub--corridor of $C_i$ is equivalent to a sub--corridor
of $C^\prime_j$. From Lemmas \ref{maxl-ext}(\ref{maxl-ext-1}) and 
\ref{2corr-maxint}(\ref{2corr-maxint-1}) it follows that
$C_i$ is equivalent to $C^\prime_j$. Thus $E=C_i+C_{1-i}$ is equivalent to
$C^\prime_j+C_{1-i}$. 
Hence $C_{1-i}$ is diagonally--equivalent to  $C^\prime_j$ which is 
diagonally--equivalent to $C^\prime_{1-j}$. It follows that $C_{1-i}$ and $C^\prime_{1-j}$ are
equivalent and therefore so are $E$ and $E^\prime$.
\item If a  simple sub--corridor of $E$ is diagonal to a  simple sub--corridor of $E^\prime$ then,
for some $i,j\in\{0,1\}$ a sub--corridor of $C_i$ is diagonal to 
$C^\prime_j$. From Lemma \ref{maxl-ext}(\ref{maxl-ext-1}) and Lemma 
\ref{2corr-maxint}(\ref{2corr-maxint-2}) it follows that $C_i$ and $C^\prime_j$
are diagonally--equivalent. As $C_i$ is diagonally--equivalent to $C_{1-i}$ it follows
that $C^\prime_j$ is equivalent to $C_{1-i}$. As in the previous case   it follows 
that $E$ and $E^\prime$ are equivalent.
\ee

\begin{lemma}\label{ext2corr-int}
Let $E=C_0+C_1$ and $E^\prime=C_0^\prime+C_1^\prime$ 
be maximal extended 2--corridors. Then either $E$ and $E^\prime$ are
equivalent or $\bar E\cap \bar E^\prime\subset \bar C_i\cap \bar C^\prime_j$, for
some $i,j\in\{0,1\}$.
\end{lemma}
\medskip

\noindent
{\em Proof.}
If $\bar{E} \cap \bar{E}^\prime=\emptyset$ then the result holds trivially.
Suppose then that $E$ and $E^\prime$ are not equivalent and that 
$\bar{E} \cap \bar{E}^\prime\neq\emptyset$. Note that this implies that
the images of $E$ and $E^\prime$ are disks (and not annulii or M\"obius bands).
 There is a unique vertex
$u_0$ of $\bar C_0$ not in $\bar C_1$ and
a unique vertex
$u_1$ of $\bar C_1$ not in $\bar C_0$. Similarly,
 there is a unique vertex
$v_0$ of $\bar C^\prime_0$ not in $\bar C^\prime_1$ and
a unique vertex
$v_1$ of $\bar C^\prime_1$ not in $\bar C^\prime_0$. Fix $j=0$ or $1$: if the 
vertex $v_j$ is not in $\bar E$ then 
$\bar E\cap \bar E^\prime\subset \bar E  \cap \bar C^\prime_{1-j}$. 
Now if $u_i$ is not in $\bar C^\prime_{1-j}$, for 
$i=0,1$, then  
$\bar E\cap \bar E^\prime\subset \bar C_{1-i} \cap \bar C^\prime_{1-j}$. 
On the other hand if $u_i\in \bar C^\prime_{1-j}$, for $i=0$ and $1$, then 
$\bar C^\prime_{1-j}$ contains $\bar E$, contrary to maximality of $E$. Finally suppose
that $v_j$ belongs to $\bar E$, for $j=0$ and 1. As $\bar E$ and $\bar E^\prime$ are
both graphs of extended 2--corridors on a disk this implies that $\bar E^\prime$ 
is a subgraph of $\bar E$. As $E$ and $E^\prime$ are not equivalent it follows that
$\bar E^\prime$ is a proper subgraph of $\bar E$ and so $E^\prime$ is not maximal,
a contradiction.

\begin{defn}
Choose one representative $E$ from each equivalence class of maximal
2--corridors of $\G$.
If $E$ is the representative of an equivalence class of maximal
extended 2--corridors then we call $E$ a {\em designated} extended
2--corridor. Corresponding to each designated extended 2--corridor
$E$ choose an ordered pair
$(C,C^\prime)$ of diagonally--equivalent simple 2--corridors such that
$E=C+C^\prime$ and say that $E$ has {\em foundation} $C$ and
{\em extension} $C^\prime$. If $C$ is the chosen representative of an
equivalence class of maximal simple 2--corridors or $C$ is the foundation of 
a designated extended 2--corridor
then we call $C$ and all its induced sub--corridors {\em designated} 
simple 2--corridors. We say that a 2--corridor is {\em designated} if it
is either a designated extended 2--corridor or a designated simple 2--corridor.
\end{defn}
We assume from now on that all 2--corridors are designated. In particular if 
$C_s$ is a proper sub--corridor of the extended 2--corridor $E$ with foundation
$C$ then, by definition, $C_s$ is a sub--corridor induced from $C$. 
(Note that if $E$ is a designated extended 2--corridor with
foundation $C$ and extension $C^\prime$ then $C$ is a designated
simple 2--corridor but $C^\prime$ is not.) 

\begin{figure}
\psfrag{C}{$C$}
\psfrag{C'}{$C^\prime$}
\begin{center}
  { \includegraphics[scale=0.6]{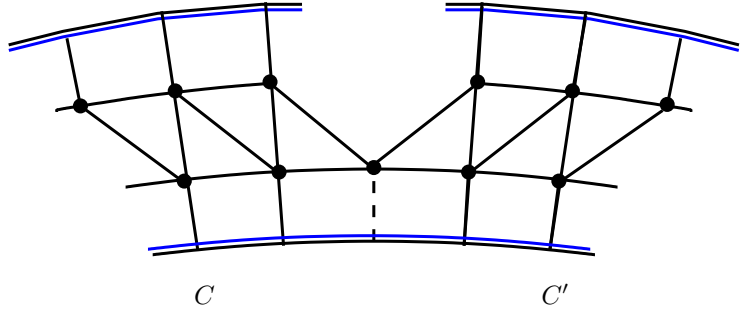}
   }\\
\caption{The unique edge of the intersection is shown as a dashed line}\label{keel}
\end{center}
\end{figure}

Suppose now that $E=C_0+C_1$ and $E^\prime=C_0^\prime+C_1^\prime$ 
are  maximal designated 
2--corridors, with foundations $C_0$ and $C_0^\prime$,
respectively. Assume further that  $E$ and $E^\prime$ are not
equivalent. If $E$ and $E^\prime$ overlap then we may assume
that the images of $E$ and $E^\prime$ are not annulii or M\"obius bands and, as
in the case of simple 2--corridors,
that the graphs $\bar E$ and $\bar E^\prime$ have at least an edge in common.
Furthermore, from Lemma \ref{2+corr-maxint},  no sub--corridor of $E$ is equivalent or 
diagonal to a sub--corridor of $E^\prime$ and, from Lemma
\ref{ext2corr-int}, $\bar E\cap \bar E^\prime\subset \bar C_i\cap \bar
C_j^\prime$, for some $i,j\in \{0,1\}$.
It follows that no compartment of $C_i$ is equivalent or diagonal to
a sub--corridor of $C_j^\prime$, and vice--versa. Therefore $\bar
C_i\cap \bar C_j^\prime$ contains a unique edge. Hence the
intersection of $E$ and $E^\prime$ contains one vertex and one class
of lateral arcs (corresponding to the edge of  $\bar C_i\cap \bar
C_j^\prime$). Therefore the union of the images of $E$ and $E^\prime$
form a configuration of the form shown in Figure \ref{keel}: in  which
case we say $E$ and $E^\prime$ form 
a {\em keel}.

We collect together our conclusions on the intersection of corridors
in the following Lemma.
\begin{lemma}
If $C$ and $C^\prime$ are maximal designated $j$--corridors, $j=0,1$
or 2, then either $C$ and $C^\prime$ do not overlap or $j=2$ and
$C$ and $C^\prime$ form a keel.
\end{lemma}
From now on all corridors are assumed to be designated.
\subsection{$Z$--cancellation}
Let $C=(D,\G_D,\mu)$ be a $j$--corridor, where $j=0,1$ or 2 and
$D=[0,c]\+[0,1]$, for some integer $c$. Recall that $\G$ has boundary
partition $\mbf b=(b_1,\ldots ,b_{W_1})$ and suppose that $C$ binds
boundary intervals $P_0=b_{p_0}$ and $P_1=b_{p_1}$, where $1\le p_1\le
W_1$, that $P_i$ has prime label
$(h_i,f_i)$ and that $\al(f_i)=m_i$: so  the label of $P_i$ is
$h_i\^m_i$. Taking the point $\mu((0,i))$ to precede $\mu((c,i))$
gives an orientation $\zeta_i$ of $P_i$. If $\zeta_i$ coincides  with the fixed
orientation $\zeta(P_i)$ of $P_i$ set $\e_i=1$ and if not set
$\e_i=-1$. 
Fix an integer $s$ with $0\le s\le c-1$ and assume that $\G_D(s)$ has 
base $D_s[x_0,x_1,y_0,y_1]$. 
The label $w_i$ of the sub--interval of $P_i$, between the point 
$\mu((s+\frac{1}{x_i+1},i))$ and the end of the interval $P_i$, read  
in the direction of $\zeta_i$, therefore begins with a terminal
segment of the word $h_i^{\e_i}$ of uniquely determined length $r_i$
such that $0\le r_i<l(h_i)$. We define the {\em left boundary marking
  of rank} $s$ of $C$ to be the sequence 
\[p_0,p_1,\e_0,\e_1,r_0,r_1.\]
The {\em left boundary marking
  of rank} $s$ of an extended 2--corridor $E$ is defined to be the  
left boundary marking of rank $s$ of its foundation.

\begin{figure}
\psfrag{(p+1/(x0+1),0)}{$\left(p+\frac{1}{x_0+1},0\right)$}
\psfrag{(p+1/(x1+1),1)}{$\left(p+\frac{1}{x_1+1},1\right)$}
\psfrag{(q+1/(x1+1),1)}{$\left(q+\frac{1}{x^\prime_1+1},1\right)$}
\psfrag{(q+1/(x0+1),0)}{$\left(q+\frac{1}{x_0^\prime+1},0\right)$}
\psfrag{Dp}{$D_p$}
\psfrag{Dq}{$D_q$}
\psfrag{q+1}{$q+1$}
\psfrag{I0}{$I_0$}
\psfrag{I1}{$I_1$}
\psfrag{p}{$p$}
\psfrag{q}{$q$}
\begin{center}
{ \includegraphics[scale=0.6]{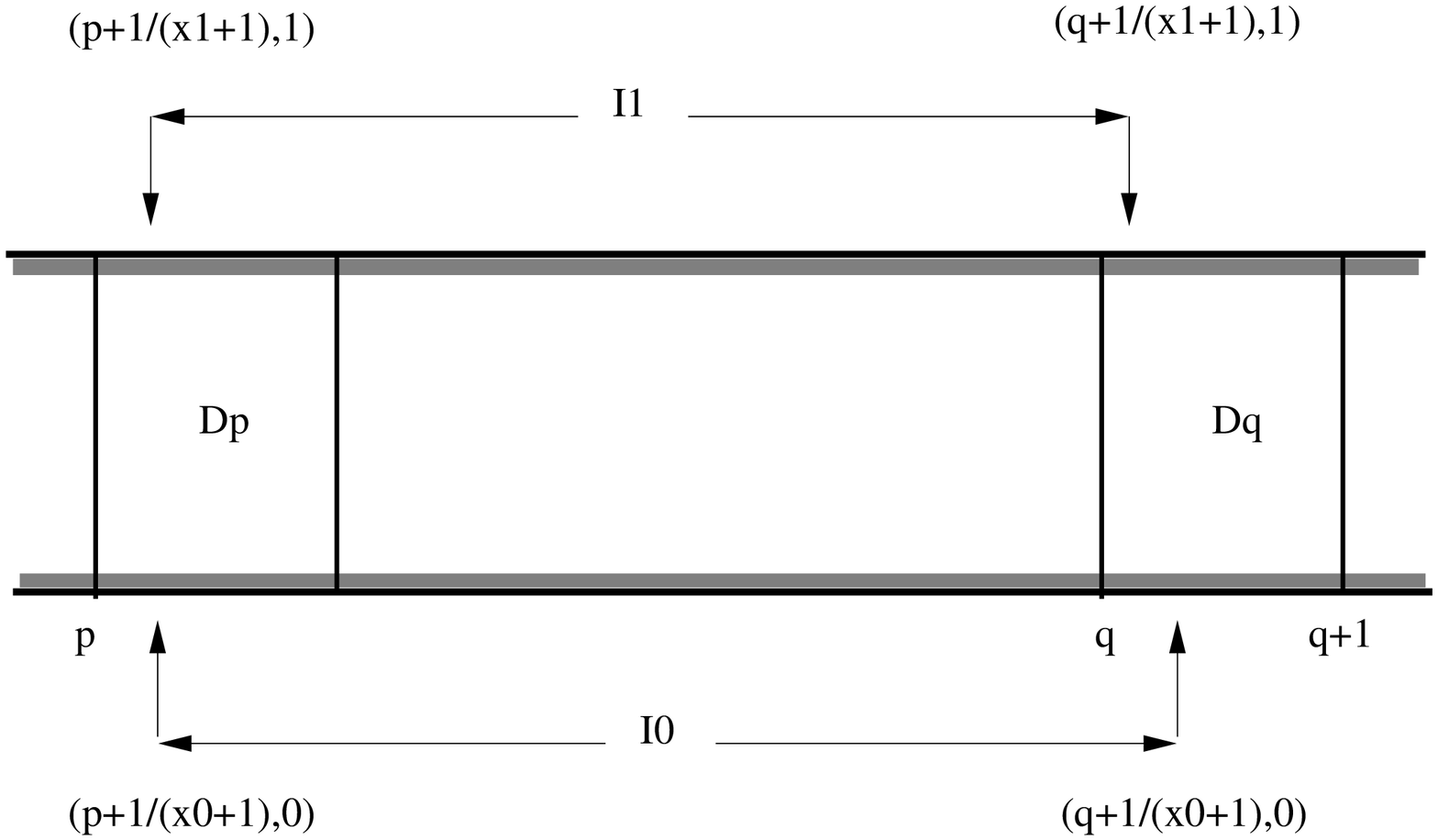}
   }\\
\caption{}\label{Z-canc}
\end{center}
\end{figure}

Now assume that $C$ above is a maximal simple (designated) $j$--corridor. Let
$0<p<q< c$ and 
consider the sub--corridor $C[p,q]$ 
of
$C$
of width $p-q\le c-1$.
Suppose that the left boundary marking $p_0,p_1,\e_0,\e_1,r_0,r_1$ of
rank $p$ of $C$ is 
equal to the  left boundary marking of rank $q$ of $C$. Let
$\G_{D}(p)$ have base $D_p(x_0,x_1,y_0,y_1)$ and 
$\G_{D}(q)$ have base
$D_{q}(x^\prime_0,x^\prime_1,y^\prime_0,y^\prime_1)$
and let $I_i$ be the sub--interval
\[I_i=\left[p+\frac{1}{x_i+1},q+\frac{1}{x^\prime_i+1}\right]\+\{i\}\]
of $[0,c]\+\{i\}$, for $i=0,1$ (see Figure \ref{Z-canc}).
Then the label of  
$\mu(I_i)$ read in the direction $\zeta_i$ is a cyclic permutation of
$(h^{\e_i})^{\al_i}$,
for some positive integer $\al_i$. 
Assume further that the left
interior markings of $\G_{D}(p)$ and $\G_{D}(q)$ are
equal (so in particular $y_0=y_0^\prime$). Under these circumstances we
form a new picture $\G^\prime$ from $\G$ as follows. 
Let $T:\RR^2\maps\RR^2$ be the translation defined by
\begin{equation}\label{T-translation}
T(x,y)=(x-(q-p),y), \textrm{ for }(x,y)\in \RR^2. 
\end{equation}
Let $S$ be a homeomorphism
$S:D_p\maps [p,q+1]\+[0,1]$ such that $S$ fixes \[\{p\}\+[0,1]\cup
[p,p+\frac{1}{x_0+1}]\+\{0\}
\cup [p,p+\frac{1}{x_1+1}]\+\{1\}\] and $S((p+1,t))=(q+1,t)$, for all
$t\in [0,1]$. Define $D^l=[0,p]\+[0,1]$ and
$D^r=[p+1,c-(q-p)]\+[0,1]$. 
Consider the triple $C^\prime=(D^\prime,\G_{D^\prime},\mu^\prime)$
where $D^\prime=[0,c-(q-p)]\+[0,1]$,
\[\G_{D^\prime}(i)=\left\{
\begin{array}{ll}
\G_D(i) & \textrm{ if } 0\le i\le p-1\\
T(\G_D(i+(q-p)))  & \textrm{ if } p\le i\le c-(q-p)-1
\end{array}
\right.
,
\] 
$\mu^\prime|_{D^l}=\mu|_{D^l}$, $\mu^\prime|_{D^r}=(\mu\circ
T^{-1})|_{D^r}$
and $\mu^\prime|_{D^\prime_p}=(\mu\circ S)|_{D^\prime_p}$.
%
Then $\G_{D^\prime}$ is a sum of square--pictures with left--trivial
labelling on the disk $D^\prime$. 
Remove $\mu(\G_D)$ from $\G$ and replace it with
$\mu^\prime(\G_{D^\prime})$. Under the given 
hypotheses this results in  a picture $\G^\prime$ over $G$ on
$\S$. The triple $C^\prime$ is, by construction, a maximal simple corridor
of $\G^\prime$. Furthermore
$\mu([p,q+1]\+[0,1])=\mu^\prime([p,p+1]\+[0,1])$ and $\G\cap (\S\bl \mu([p,q+1]\+[0,1]))=
\G^\prime\cap (\S\bl \mu^\prime ([p,p+1]\+[0,1]))$. Hence all
corridors of $\G$, except $C$ are also corridors of $\G^\prime$.
We 
designate $C^\prime$ as the representative for its equivalence class of
maximal corridors in $G^\prime$. 

\begin{figure}
\psfrag{u}{$u$}
\psfrag{v}{$v$}
\psfrag{e1}{$e_1$}
\psfrag{en}{$e_n$}
\psfrag{d1}{$d_1$}
\psfrag{dn}{$d_n$}
\psfrag{a}{$a$}
\psfrag{b}{$b$}
\psfrag{a'}{$a^\prime$}
\psfrag{b'}{$b^\prime$}
\psfrag{p}{$p$}
\psfrag{q}{$q$}
\psfrag{p'}{$p^\prime$}
\psfrag{q'}{$q^\prime$}
\psfrag{...}{$\cdots$}
\begin{center}
{ \includegraphics[scale=0.34]{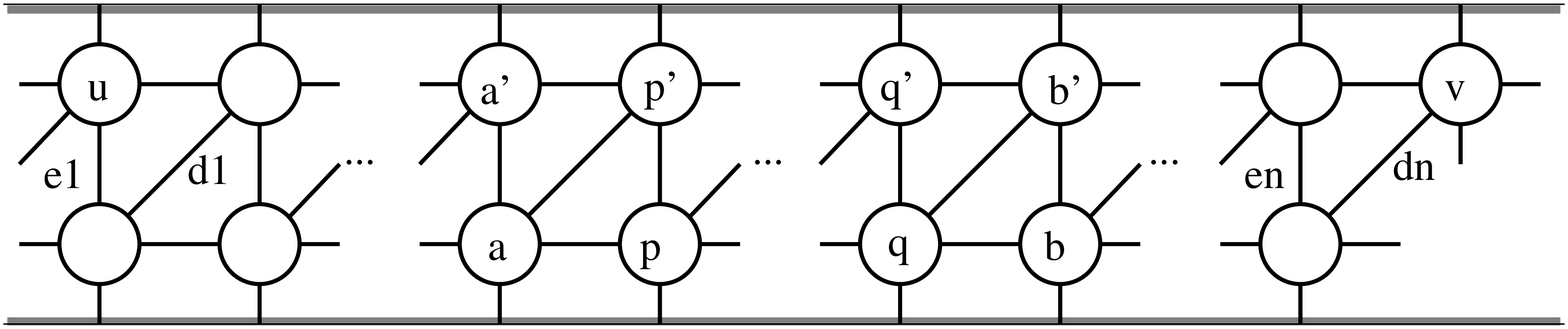}
   }\\
\caption{}\label{diag_cancel_1}
\end{center}
\end{figure}

Now assume that $C$ is the foundation of a designated maximal extended $2$--corridor $E$ 
with extension $C_e$. For example the graphs $\bar C$ and $\bar C_e$ may be as shown in
Figure \ref{diag_cancel_1}, where the vertical edges $e_1,\ldots ,e_n$ are composed of internal--lateral
arcs of compartments of $C$ and the diagonal edges $d_1,\ldots ,d_n$ are composed of internal--lateral
arcs of compartments of $C_e$. Thus $u$ is the vertex of $C$ not in $C_e$ and $v$ is the vertex of
$C_e$ not in $C$.

\begin{figure}
\psfrag{u}{$u$}
\psfrag{v}{$v$}
\psfrag{e1}{$e_1$}
\psfrag{en}{$e_n$}
\psfrag{d1}{$d_1$}
\psfrag{dn}{$d_n$}
\psfrag{a}{$a$}
\psfrag{b}{$b$}
\psfrag{a'}{$a^\prime$}
\psfrag{b'}{$b^\prime$}
\psfrag{p}{$p$}
\psfrag{q}{$q$}
\psfrag{p'}{$p^\prime$}
\psfrag{q'}{$q^\prime$}
\psfrag{...}{$\cdots$}
\begin{center}
{ \includegraphics[scale=0.34]{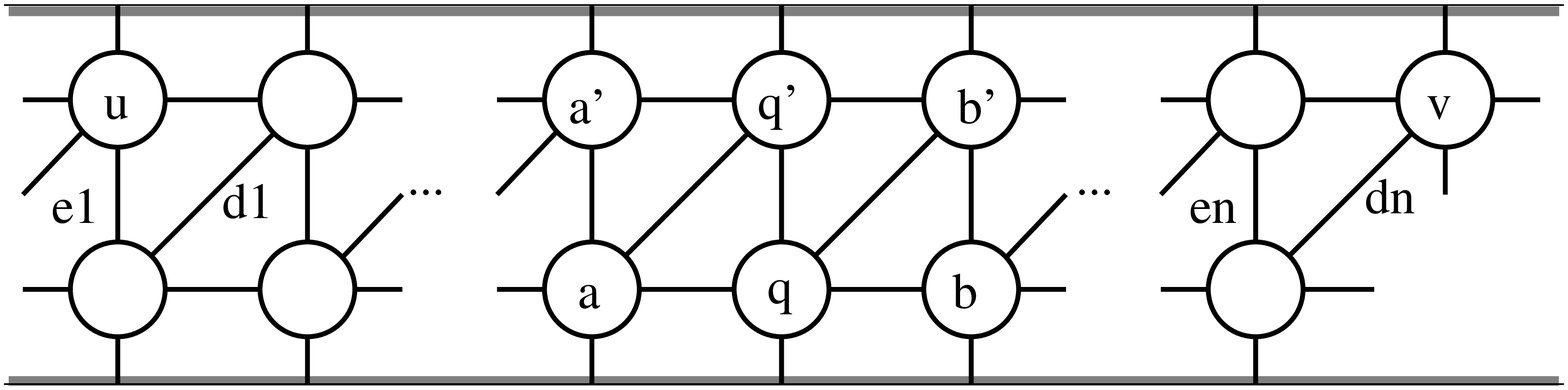}
   }\\
\caption{}\label{diag_cancel_2}
\end{center}
\end{figure}

Cancellation from $p$ to $q$ along $C$ results in corridors $C^\prime$ and $C^\prime_e$
 which have graphs $\bar C^\prime$ and $\bar C^\prime_e$, as shown in Figure \ref{diag_cancel_2}.
The compartments of $C_e$ containing $a^\prime$ and all those to its left become compartments of
$C_e^\prime$, as does the compartment of $C_e$ containing $b^\prime$ and all those to its right.
In addition $C^\prime_e$ has one compartment containing $a$ and $q^\prime$ which replaces all the
compartments of $C_e$ between those containing $a^\prime$ and $b^\prime$. Therefore the only compartment
of $C^\prime_e$ which is not a compartment of $C_e$ is that containing vertices $a$ and $q^\prime$. Note
that in $C$  the left interior markings of the compartments containing $p$ and $q$ are equal. It follows
that the left interior markings of the compartments of $C_e$ and $C_e^\prime$ containing $a$ are equal.
From the graph it is clear that $C^\prime$ and $C^\prime_e$ are diagonal and, since $C$ and $C_e$ are maximal
so are $C^\prime$ and $C^\prime_e$. Thus $E^\prime=C^\prime+C^\prime_e$ is a maximal extended $2$--corridor. 
We choose $E^\prime$ as the designated representative of its class in $\G^\prime$ and $C^\prime$ and $C^\prime_e$
as its foundation and extension, respectively. 

The formation of $\G^\prime$ from $\G$ along a corridor $C$ in this way is called $Z$--{\em
cancellation along} $C$ {\em from} $p$ to $q$. If $C$ is  of type $j$
we refer to $Z(j)$--cancellation along $C$, $j=0,1,2$.
We remark that if  $\G^\prime$ is formed from $\G$ by
$Z$--cancellation then the left interior and left boundary markings of
compartments of corridors of $\G^\prime$ are
left interior and left boundary markings of
compartments of corridors of $\G$.
Also, the boundary partition $\mbf b$ is unaffected by
$Z$--cancellation, the prime labels of boundary intervals are the same
in $\G^\prime$ as in $\G$
and for $j\neq p_i$ the boundary interval $b_j$ has the
same label in $\G^\prime$ as in $\G$. With the above notation, the
label of $P_i=b_{p_i}$ in 
$\G^\prime$ is $h_i\^(m_i-l(h_i)\al_i)$, for $i=0,1$.

\begin{theorem}\label{z-canc-poss}
If $C$ is a $j$--corridor, for $j=0,1$ or 2 and $C$ has more than
$M(j)$ arcs then $Z$--cancellation is possible along some
sub--corridor of $C$. ($M(j)$ is defined in (\ref{M_j}).)
\end{theorem}
\medskip

\noindent{\em Proof.}
Assume first that $C$ is  a simple corridor of width $c$,
$C=(D,\G_D,\mu)$. Then $Z$--cancellation is possible if 
there are integers $q_0,q_1$,  with $0< q_0<q_1<c$, such that
the left interior marking of $C$ of rank $q_0$ equals the left interior
marking of $C$ of rank $q_1$ and 
the boundary marking of $C$ of rank $q_0$ equals the boundary
marking of $C$ of rank $q_1$.

Let $p_0,p_1,\e_0,\e_1,r_0,r_1$ be the left boundary marking of rank
$q_0$  and let $(h_{p_i},m_{p_i})$ be the prime label of the boundary
interval $b_{p_i}$, for $i=0,1$. Since $p_i$ and $\e_i$ are determined
by $C$ then, for all $q_1$ with $0<q_1<c$, the left boundary
marking of rank $q_1$ is  $p_0,p_1,\e_0,\e_1,s_0,s_1$, for some
integers $s_0$ and $s_1$. By definition $0\le r_i,s_i<l(h_i)\le
W_2(\mbf z)$. Hence there are at most $W_2^2$ distinct left boundary
markings of $C$. If $j=0$  the left interior
markings of all ranks of $C$ are identical and the number of arcs of
$C$ is $c$. If $c>W_2^2+1=M(0)$ then from the above it follows that
$Z$-cancellation along $C$ is possible, as claimed.

If $j=1$ then the left interior marking of rank $s$ of $C$ is
$\de,o,\rho,\s,0,0,0,0$, where
$\de\in\{-1,0,1\}$, $o,\s\in\{0,\ldots ,l(r)-1\}$,  $0\le \rho<m$ and $\de=o=\rho=\s=0$ if and
only if $\G_D(s)$ is of type square--$0$. Therefore there are at most
$3ml(r)^2$ distinct left interior markings of $C$. It follows that if 
$C$ has width $c>W_2^2\cdot 3ml(r)^2+1$ then $Z$--cancellation along $C$
is possible. A square--picture of type 0 has 1 arc while a
square--picture of type 1 has $ml(r)$ arcs. Hence if $C$ has more that 
$(W_2^2\cdot 3ml(r)^2+1)ml(r)\le M(1)$ arcs then $Z$--cancellation along $C$
is possible.

If $j=2$ then  the left interior marking of rank $s$ of $C$ is
$\de_0,o_0,\rho_0,\s_0,\de_1,o_1,\rho_1,\s_1$, where
$\de_i=\pm 1$, $o_i,\s_i\in\{0,\ldots ,l(r)-1\}$ and $0\le \rho_i<m$.
Therefore there are at most
$4m^2l(r)^4$ distinct left interior markings of $C$.
It follows that if 
$C$ has width $c>W_2^2\cdot 4m^2l(r)^4+1$ then $Z$--cancellation along $C$
is possible. A square--picture of type 2 has less than $2ml(r)$ arcs
so if  $C$ has more that 
$(W_2^2\cdot 4m^2l(r)^4+1)ml(r)$ arcs then $Z$--cancellation along $C$
is possible.

If $C$ is not simple then $j=2$ and $C$ has foundation a simple
2--corridor $C^f$. Then $C$ has at most $ml(r)$ more arcs than
$C^f$. Thus if $C$ has more than $(W_2^2\cdot 4m^2l(r)^4+2)ml(r)=M(2)$
arcs then, from the above $Z$--cancellation is possible along $C^f$
and therefore along $C$, as claimed.
\medskip

From this Theorem it follows that a $Z$--cancellation along a
$j$--corridor may be achieved
by a sequence of $Z$--cancellations along sub--corridors of width at
most $M(j)$.

\begin{defn}\label{z_reduced}
A picture on which no $Z$--cancellation is possible is called $Z${\em --reduced}.
\end{defn} 

\begin{corol}\label{corridor-size} 
If $C$ is  a $j$--corridor of a $Z$--reduced picture then $C$ has at
most $M(j)$ arcs.
\end{corol}

\begin{lemma}\label{canc-effic}
Let $\G$ be a picture on a surface $\S$ of
type $(\mbf n,\mbf t,\mbf p)$ with  prime
labels $\mbf z$, labelled by $\hat\al(\mbf z)$ and with boundary partition 
$\mbf b=(b_1,\ldots ,b_{W_1})$. Assume that $\G^\prime$ is obtained from $\G$
by $Z$--cancellation.  Then $\G^\prime$ is reduced if and only if $\G$ is reduced. Furthermore 
$\G^\prime$ is reduced and minimalistic if and only if $\G$ is reduced and minimalistic. 
\end{lemma}
\medskip

\noindent{\em Proof.} Suppose first that $\G$ is reduced.
If $\G^\prime$ has a pair of cancelling vertices
$u,v$ then these must be adjacent in $\G^\prime$ but not in
$\G$. Hence $u,v$ belong to consecutive compartments of some corridor
$C=(D,\G_D,\mu)$ of $\G^\prime$ 
which results from 
$Z$--cancellation in $\G$. It follows
that $C$ must have type $j=1$ or $2$.  As $u$ and $v$ cancel in
$\G^\prime$ they are the images of vertices which 
cancel in $\G_D$; so we restrict attention to $\G_D$ and assume $u$
and $v$ belong to $\G_D$. Suppose $u$ belongs to $\G_D(p)$, $v$
belongs to $\G_D(p+1)$ and that $D_p$ and $D_{p+1}$ have bases
$D_p(x_0,x_1,y_0,y_1)$ and
$D_p(x_0^\prime,x_1^\prime,y_1,y_1^\prime)$, respectively. If $j=1$
then it must be possible to perform a bridge move on arcs
joining $u$ and $v$ to $[p,p+1]\+\{i\}$, for $i=0$ or 1. Assume that
this occurs for $i=0$ and let $a_0$ and $a_1$ be the arcs joining $u$
to $q_0=(p+1-\frac{x_0}{1+x_0},0)$ and $v$ to
$q_1=(p+1+\frac{1}{1+x^\prime_0},0)$, respectively: then a bridge move
is possible transforming $a_0$ and $a_1$ into arcs $b_0$, joining $u$
to $v$,  and  $b_1$ with end points $q_0$ and $q_1$. However this
means that the interval $[q_0,q_1]$ of $\pd D$ has trivial label,
contrary to the standing hypothesis that $\mbf z$ is special. Hence if
$\G_D$ is of type 1 no such cancellation is possible. 

Now assume that $j=2$. If $u$ and $v$ are both incident
to $[p,p+2]\+\{i\}$ then, as in the  case $j=1$, no bridge move is
possible on arcs joining $u$ and $v$ to $\pd D$.  Hence $u$ and $v$ do
not cancel. We assume therefore that $u$ is incident to
$[p,p+1]\+\{0\}$ and that $v$ is incident to
$[p+1,p+2]\+\{1\}$. However no bridge moves are possible now because
both regions incident to the edge of $\ovr{\G}_D$ joining $u$ to $v$
are triangular: a bridge move on the sides of a triangle would result
in a trivial label on the corner of the triangle opposite the edge joining $u$
to $v$. Thus no cancellation is possible if $j=2$.  
It follows that if $\G$ is reduced then $\G^\prime$ is reduced.
A similar argument shows that the converse also holds.

Now suppose that $\G$ is reduced and minimalistic.
Then $\G^\prime$ is reduced.
Regions involved in $Z$--cancellation all have Euler characteristic $1$
and both $\G$ and $\G^\prime$ are pictures on the surface $\S$, so the 
condition on Euler characteristic of regions is satisfied by $\G^\prime$.
No arc of $\G^\prime$ is a closed curve since no such arcs  
are created in the process of $Z$--cancellation.
Suppose that  $\G^\prime$ has interior edge of width greater than $l-2$.
Since $\G^\prime$ is reduced it follows from Lemma \ref{bound_edge_width}
that both ends of this edge are incident to the same vertex. No such edge 
results from the process of $Z$--cancellation so this edge also belongs
to $\G$, which cannot therefore be minimalistic. It follows that
$\G^\prime $ is reduced and minimalistic.
The converse follows similarly.

\subsection{Corridor--sections}\label{corridor-section}

Let $\G_D$ be a sum of square--pictures of type $j$ on the disk $D=[0,c]\+[0,1]$
and let $w_i$ denote the label of $[0,c]\+\{i\}$, for $i=0,1$, $j=0,1$
or $2$.
Assume that the following conditions hold.
\be
\item There exist integers $p_0,p_1, \e_0,\e_1,t_0$ and $t_1$ such that $1\le p_i\le
  W_1$, the $p_i$th  letter $(h_i,f_i)$  of $\mbf z$ is not 
  a degenerate letter or
  a minor letter of $H^\Lm$, $\e_i=\pm 1$, $0\le 
  t_i<l(h_i)$ and $w_i=g_i^{\al_i}$, for some $\al_i>0$, and 
  $g_i=\t(h_i^{\e_i},t_i)\io(h^{\e_i},l(h_i)-t_i)$, for
  $i=0,1$.
\item $\G_D$ has at most $M(j)$ arcs.
\item Let $D^e$
  be the disk $[c,c+1]\+[0,1]$.
  Let $T_c$ be the translation $T(x,y)=(x+c,y)$, for $(x,y)\in \RR^2$, and
let $\G_{D^e}=T_c(\G_D(0))$, a picture of type square--$j$ on $D_e$. Then
$\G_D\cup \G_{D^e}$ is a sum of square--pictures
  on the disk $[0,c+1]\+[0,1]$.
\ee
Then $(\G_D,p_0,p_1,\e_0,\e_1,t_0,t_1)$ 
is called  a {\em corridor--section} on $\G_D$, of width $c$, binding
the $p_0$th and $p_1$th letters of $\mbf z$, and of {\em type} $j$ if
$\G_D$ is a sum of square--pictures of type $j$. 

The corridor--section $(\G_D,p_0,p_1,\e_0,\e_1,t_0,t_1)$ is said to have 
{\em left interior marking} equal to the left interior marking of $\G_D$,
{\em left boundary marking} 
the sequence 
\[p_0,p_1,\e_0,\e_1,t_0,t_1\] and  
{\em extent} 
the pair$(\al_0,\al_1)$. 

Let $C$ be a corridor section of width $c$.
The conditions on a corridor--sections imply that if $0< p<c-1$ then  a new
corridor--section $C_p$ may be formed by cutting $C$ along $\{p\}\+[0,1]$ and then
identifying $(0,t)$ with $(c,t)$, for all $t\in [0,1]$. We set $C=C_0$. 
The left interior and left
boundary markings of $C_p$ are called the {\em left interior} and {\em left boundary}
markings, respectively, of $C$ of {\em rank} $p$. Thus the left interior marking and 
left boundary marking of $C$ coincide with the left interior and boundary markings 
of rank 0 of $C$.

Let $L$ denote the set of all corridor sections. As the number of arcs
of a corridor section is bounded $L$ is finite. Let $L$ have $|L|$
elements indexed by the integers $1,\ldots, |L|$ and denote these
elements $L(1),\ldots,L(|L|)$. For $i=1,\ldots, |L|$, define the quintuple
\[L_5(i)=(L(i),s(i),t(i),\al(i,0),\al(i,1))\] where $s(i)$ and $t(i)$ are
the left boundary marking and left interior marking, respectively, of
$L(i)$ and $(\al(i,0),\al(i,1))$ is the extent of $L(i)$.
Furthermore write 
\[s(i)=p(i,0),p(i,1),\e(i,0),\e(i,1),r(i,0),r(i,1)\]
for the left boundary marking of $L(i)$.
For fixed $i$ with $1\le i\le |L|$ and each integer $p$ with $1\le p\le W_1$ define
\[
\zeta(p,i)=\left\{
\begin{array}{ll}
0,& \textrm{ if } p\neq p(i,0) \textrm{  and } p\neq p(i,1)\\
\al(i,0), & \textrm{ if } p=p(i,0) \textrm{  and } p\neq p(i,1)\\
\al(i,1), & \textrm{ if }  p\neq p(i,0) \textrm{  and } p=p(i,1)\\
\al(i,0)+\al(i,1), & \textrm{ if }  p=p(i,0)=p(i,1)
\end{array}
\right.
.
\]

The following lemma follows directly from the
definitions above and  the remarks preceding Theorem
\ref{z-canc-poss}.  
\begin{lemma}\label{corr-sec-remove}
Let $C$ be a $j$--corridor of $\G$
 and let $C[p,q]=(D,\G_D,\mu)$ be a sub--corridor (of
width $q-p$) with at most $M(j)$ arcs.
If $Z$--cancellation takes place along $C$ from $p$ to
$q$ and the left boundary marking of $C$ of rank $p$ is $t$ 
then $(\G_D,t)$ is a corridor--section of width $q-p$ and
type $j$. Let $(\G_D,t)=L(i)$, for some $i$ with $1\le
i\le |L|$, let $\G^\prime $ be the picture formed by this
$Z$--cancellation and let 
the boundary interval $b_p$ of $\G$ have label
$h_p\^m_p$, for $p=1,\ldots,W_1$. Then  the boundary interval of $b_p$
in $\G^\prime$ has label 
\[h_p\^(m_p-l(h_p)\zeta(p,i)),\] 
for $p=1,\ldots,W_1$.
\end{lemma}

The $Z$--cancellation of Lemma \ref{corr-sec-remove}, along $C$ from $p$ to
$q$, where the left boundary marking, of rank $p$, of $C$ is $t$ is
called $Z$--cancellation of corridor--section $(\G_D,t)$. As remarked
(above  Definition \ref{z_reduced}) any $Z$--cancellation may be effected by a sequence of
$Z$--cancellations of corridor--sections.

\subsection{$Z$--insertion}
Let $C=(D,\G_D,\mu)$ be a maximal simple $j$--corridor, with
$D=[0,c]\+[0,1]$, binding $b_{p_0}$ and $b_{p_1}$, for $p_i$ with
$1\le p_i\le W_1$. Let the left interior marking and left boundary marking of
of $C$, of rank $p$, be $s$ and $t$, respectively.
Let $(\G_\D,t)$ be a corridor section with left interior marking $s$, where $\D=[0,q]\+[0,1]$.
Let $T_p:\RR^2\maps\RR^2$ and $T_q:\RR^2\maps\RR^2$ be the translations defined by
$T_p(x,y)=(x+p,y)$ and $T_q(x,y)=(x+q,y)$, for $(x,y)\in \RR^2$ and positive integers $p$ and $q$. 
Let $D_p$ have base $D_p(x_0,x_1,y_0,y_1)$ and $\D_0$ have base
$\D_0(x_0^\prime,x_1^\prime,y_0^\prime,y_1^\prime)$. Let   $S$ be a homeomorphism
$S:[p,p+q+1]\+[0,1]\maps [p,p+1]\+[0,1]$ such that $S$ fixes 
$\{p\}\+[0,1]$, 
\[S([p,p+\frac{1}{x^\prime_i+1}]\+\{i\})=[p,p+\frac{1}{x_i+1}]\+\{i\},\]
for $i=0,1$,
and $S((p+q+1,t))=(p+1,t)$, for all
$t\in [0,1]$. Define $D^l=[0,p]\+[0,1]$, $D^c=[p,p+q+1]\+[0,1]$ and
$D^r=[p+q+1,c+q]\+[0,1]$. 

Consider the triple $C_I=(D_I,\G_{D_I},\mu_I)$, where
$D_I=[0,c+q]\+[0,1]$,
\[
\G_{D_I}(i)=\left\{
\begin{array}{ll}
\G_D(i) & \textrm{ for  } i=0,\ldots , p-1\\
T_p(\G_{\D}(i-p)) &\textrm{ for } i=p,\ldots ,p+q-1\\
T_q(\G_{D}(i-q)) &\textrm{ for } i=p+q,\ldots, c+q-1 
\end{array}
\right.,
\]
$\mu_I|_{D^l}=\mu|_{D^l}$, 
$\mu_I|_{D^c}=(\mu\circ S)|_{D^c}$ and
$\mu_I|_{D^r}=(\mu\circ T^{-1}_q)|_{D^r}$.

Then $\G_{D_I}$ is a sum of square--pictures with left--trivial
labelling on the disk $D_I$.
Remove $\mu(\G_D)$ from $\G$ and replace it with
$\mu_I(\G_{D_I})$. This results on  a picture $\G_I$ over $G$ on
$\S$. The triple $C_I$ is, by construction, a maximal corridor
of $\G_I$. Furthermore
$\mu([p,p+1]\+[0,1])=\mu_I([p,p+q+1]\+[0,1])$ and $\G\cap (\S\bl \mu([p,p+1]\+[0,1]))=
\G_I\cap (\S\bl \mu_I([p,p+q+1]\+[0,1]))$. Hence all
corridors of $\G$, except $C$, are also corridors of $\G_I$.
We 
designate $C_I$ as the representative for its equivalence class of
maximal corridors in $\G_I$. Moreover if $C$ is the foundation
of a designated maximal extended 2--corridor of $\G$ then $C_I$
is a maximal simple sub--corridor of a maximal extended $2$--corridor of 
$\G_I$ (as can be seen by reversing the argument for cancellation along
extended $2$--corridors). In this case $C_I$ 
is chosen as the foundation of the corresponding maximal extended
2--corridor of $\G_I$, which becomes a designated $2$--corridor of $\G_I$.   

The formation of $\G_I$ from $\G$ in this way is called $Z$--{\em
insertion of} $(\G_\D,t)$ {\em into} $C$ {\em at position} $p$. If $C$ is  of type $j$
we refer to $Z(j)${\em --insertion}, $j=0,1,2$.

\begin{lemma}\label{one-insertion}
If  $\G_I$ is formed from $\G$ by
$Z$--insertion of $(\G_\D,t)$ then the left interior and left boundary markings of
compartments of corridors of $\G_I$ are
left interior and left boundary markings of square--pictures of $\G_\D$ or of
compartments of corridors of $\G$.
Also, the boundary partition $\mbf b$ is unaffected by $Z$--insertion
and, given that
the boundary interval $b_p$ of $\G$ has label
$h_p\^m_p$, for $p=1,\ldots,W_1$, then  the boundary interval of $b_p$
in $\G_I$ has label 
\[h_p\^(m_p+l(h_p)\zeta(p,i)),\] 
for $p=1,\ldots,W_1$.
\end{lemma}

The operations of $Z$--cancellation and $Z$--insertion are inverse to
each other in an obvious way. To be more precise, let $\G_I$ be formed
from $\G$
by $Z$--insertion of $(\G_\D,t)$ into $C$ at position $p\ge 1$ and let the
resulting corridor of $\G_I$ be $C_I$ as above. Then $Z$--cancellation
of corridor--section $(\G_\D,t)$ along $C_I$ from $p$ to $q$ is
possible and results in $\G$. A similar statement holds starting with
$Z$--cancellation. 
\subsection{The $Z$-graph}

Let $L$ denote the set of corridor--sections of 
$\mbf z$ 
(indexed by the integers $1,\ldots ,|L|$ as before)
and 
define the {\em left marking}, of {\em rank} $p$, of $L(i)\in L$ to be  
$m(i,p)=s(i,p),t(i,p)$, where $s(i,p)$ and $t(i,p)$ are the left interior and boundary
markings, respectively, of rank $p$ of  $L(i)$. Define 
the {\em left marking} of $L(i)$ to be $m(i)=m(i,0)$,
let \[m_{L(i)}=\{m(i,p): 0\le p< \textrm {width of } L(i)\}\]
 and let 
\[m_{\mbf z}=\bigcup_{i=1}^{|L|}m_{L(i)}.\]

Let $C$ be a corridor of the picture $\G$ and define the {\em left marking}, 
of rank $p$,
of $C$ to be $s,t$, where $s$ and $t$ are the left interior and
boundary markings, respectively, of rank $p$, of $C$.
Define the {\em active} left markings of $\G$ to be those left markings of 
corridors of $C$ which are elements of $m_{\mbf z}$.

The $Z${\em --graph} $G(\mbf z)$ of $\mbf z$ is defined as follows.
The vertices of $G(\mbf z)$ are 
the subsets of $m_{\mbf z}$. There is one directed edge labelled $L(i)$ 
from vertex $R$ to vertex $T$ corresponding to each corridor--section $L(i)$ such that 
$L(i)$ has left marking in $R$ and $T=R\cup m_{L(i)}$. Corresponding to the picture
$\G$ is the vertex 
\[v(\G)=\{\textrm{active left markings of } \G\}\] 
of $G(\mbf z)$ which we call the {\em initial} vertex of $\G$.

Now, for $i=0,\ldots n$, let $\G_i$ be a picture, over $G$ on $\S$, such
that $\G_{i+1}$ is obtained from $\G_i$ by $Z$--cancellation of a
corridor--section $C_i$ along a
corridor of $\G_i$. Assume that $\G=\G_n$. 
Then this sequence of $Z$--cancellations uniquely
determines a path $S_n,\ldots,S_0$ in $G(\mbf z)$, where $S_i\in
V(G(\mbf z))$, $S_n=v(\G)$, $S_{i}=S_{i+1}\cup \{\textrm{left markings
  of } C_i\}$ and the edge of the path
joining $S_{i+1}$ to $S_{i}$ is labelled $C_i$. Then $\G_i$ is
obtained from $\G_{i+1}$ by $Z$--insertion of $C_i$. Hence the set of active left
markings of $\G_{i}$ is the union of the left markings of $C_i$ with
the active left markings of $\G_{i+1}$. It follows, by induction, that
the set of active left markings of $\G_0$ is $S_0$.  

Conversely, let $S_n,\ldots,S_0$ be a directed path in $G(\mbf z)$,
where $S_n=v(\G)$ and the edge of the path joining $S_{i+1}$ to $S_{i}$
is labelled $C_i$. Then $S_{i}=S_{i+1}\cup \{\textrm{left markings
  of } C_i\}$. Let $\G_n=\G$ and let $m$ be an integer with $n\ge
m>0$. Suppose inductively that we have
constructed pictures $\G_i$ with active left markings $S_i$, for
$i=n,\ldots , m$ such that $\G_{i+1}$ is obtained from $\G_{i}$ by
$Z$--cancellation of the corridor--section $C_i$, for $i=n-1,\ldots m$.
As $C_{m-1}$ is the label of the edge joining $S_{m}$ to $S_{m-1}$ it
follows that $S_{m-1}=S_{m}\cup\{\textrm{left markings of }C_{m-1}\}$ and that
the left marking of $C_{m-1}$ is in $S_{m}=\{\textrm{active left
  markings of } \G_m\}$. Hence we may form a
picture $\G_{m-1}$ from $\G_m$ by $Z$--insertion of $C_{m-1}$. The
active left markings of $\G_{m-1}$ are then the left markings of
$C_{m-1}$ together with the active left markings of $\G_m$: that is
$S_{m-1}$. Moreover as $\G_{m-1}$ is formed from $\G_m$ by
$Z$--insertion of $C_{m-1}$ it follows that $\G_m$ is formed from
$\G_{m-1}$ by $Z$--cancellation of $C_{m-1}$. Hence the path
$S_n,\ldots,S_0$ determines a sequence of $Z$--cancellations 
\[\G_0\maps\cdots\maps\G_n=\G,\] 
with $\G_{i+1}$ obtained from $\G_{i}$ by $Z$--cancellation of $C_i$ and 
the set of active left markings of $\G_0$ equal to $S_0$. 

We have therefore proved the following.
\begin{lemma}
Let $\G$ and $\G^\prime$ be pictures over $G$ on $\S$. Then there
exists a sequence 
\[\G^\prime=\G_0\maps\cdots\maps\G_n=\G,\] 
with $\G_{i+1}$ obtained from $\G_{i}$ by $Z$--cancellation of  the
corridor--section $C_i$ if and only if there exists a directed path
\[v(\G)=S_n,\ldots,S_0=v(\G^\prime)\] in $G(\mbf z)$ such that the
edge joining $S_{i+1}$ to $S_i$ is the corridor section $C_i$. 
\end{lemma}

\begin{defn}
Let $p=S_n,\ldots,S_0$ be a directed path in $G(\mbf z)$ such that the
edge joining $S_{i+1}$ to $S_i$ is the corridor section $C_i$. Let 
\[\G_0\maps\cdots\maps\G_n,\] be a sequence of pictures
over $G$ on $\S$ with $\G_{i+1}$ obtained from $\G_{i}$ by $Z$--cancellation of  the
corridor--section $C_i$. If $v(\G_n)=S_n$ then we say that $\G_0$ is
the picture obtained from $\G_n$ by $Z$--insertion {\em following}  $p$. 
If $v(\G_0)=S_0$  then we say that $\G_n$ is the picture obtained from
$\G_0$ by $Z$--cancellation {\em following} $p$. 
\end{defn}

Let $Q$ denote the quadratic exponential equation $Q(\mbf z,\cL,\mbf n,\mbf t,\mbf
  p)$ and $Q_\cH=Q(\ovr{\mbf z},\cH,\mbf n,\mbf t,\mbf
  p)$ the homogeneous equation corresponding to $\mbf z$.
Let the $q$th letter of $\mbf z$ be $(h_q,f_q)$, for $1\le q\le W_1$.

\begin{lemma}\label{follow-path}
Let $p=S_n,\ldots,S_0$ be a directed path in $G(\mbf z)$ with initial
vertex $v(\G_n)$. Let
$e_1,\ldots,e_t$ be the edges of $G(\mbf z)$ occuring in $p$,
 let $L(s_i)$ be the label of $e_i$ and let
$n_i$ be the number of occurrences of $e_i$ in $p$, for $i=1,\ldots, t$.
Let $\G^\prime$ be the picture obtained from $\G_n$ by 
$Z$--insertion following $p$. Then the
following hold.
\be
\item\label{follow-path-1} Let $\al_n$ be a solution to $\cH$ such that $(\al_n,\G_n)$
  corresponds to a solution to $Q_\cH$ and write $m_q=\al_n(\lm_q)$, for $1\le q\le W_1$.
The boundary interval $b_q$ of
$\G^\prime$  has label
\[h_q\^\left(m_q+l(h_q)\left(\sum_{i=1}^{t}n_i\z(q,s_i)\right)\right),\]
for $1\le q\le W_1$.
\item\label{follow-path-2} Let $\al$ be a solution to $\cL$ such that $(\al,\G^\prime)$
  corresponds to  a solution to $Q$. Let $\al_n:M\maps \ZZ$ be the retraction defined 
by setting 
\[\al_n(\lm_q)= \al(f_q)-l(h_q)\left(\sum_{i=1}^{t}n_i\z(q,s_i)\right),\]
for all $q\in \{1,\ldots ,W_1\}$ such that the $q$th letter of $\mbf z$ is
a proper exponential $H^\Lm$--letter, and $\al_n(\lm)=0$, for all other $\lm\in\Lm$.
Then $(\al_n,\G_n)$ corresponds to a solution of $Q_\cH$. 
\ee
\end{lemma}
\medskip

\noindent{\em Proof.}
The first part of the Lemma follows from Lemma \ref{one-insertion}.
For the second part note that the boundary label of $b_q$ in $\G_n$ is
$h_q\^\al_n(\lm_q)$, for all $q$ (using  Lemma \ref{one-insertion} again).
By definition of $Z$--cancellation therefore $\al_n(\lm_q)>0$. Also
$\al_n(\lm_q)\equiv \al(f_q)\equiv k_q (\mod l(h_q))$. 
It follows that $(\al_n,\G_n)$
corresponds to a solution to $Q_\cH$.
\medskip

A picture $\G$ such that there exists a solution $\al$ to $\cH$ and 
$(\al,\G)$ corresponds to a solution to $Q_\cH$ will be called a 
{\em basic} picture for $(\mbf z,\cL)$.

We call a path in a graph with no repeated vertices a {\em simple
  path} and a path with initial and terminal vertices identical but no
  other pair of vertices identical a {\em cycle}. A cycle of length 1
  is called a {\em loop}. The {\em path--subgraph} of a path  is the
  subgraph consisting of all edges and vertices of the path.
The directed path $p$ of $G(\mbf z)$ corresponding to a sequence of
$Z$--cancellations contains no cycle of length
greater than 1 but may contain loops. 
We call a 
directed path $p$ in $G(\mbf z)$ a $Z${\em --path} if $p$ 
has no cycle of length greater than 1. 
The path--graph of a $Z$--path
consists of the
path--subgraph of  a directed simple path together with a (possibly
  empty) collection of loops. Given a directed path $p$
assign weight $w_p(e)$ to each edge of $G(\mbf z)$, where $w_p(e)$ is
the number of times the edge appears in $p$. 
We call the edge weighting $w_p$ the {\em weighting} of the path $p$.
Then there are finitely
many directed paths in $G(\mbf z)$ which have 
the same weighting as $p$. 
We call a weighting of the edges of $G(\mbf z)$ a
$Z${\em --weighting} if it equals the weighting of a $Z$--path. 
Given a $Z$--weighting $w$ of $G(\mbf z)$ the subgraph $P$ 
consisting of all edges with positive weight
(and their incident vertices) is the path--subgraph of a $Z$--path $p$,  
every edge of $P$ has positive weight
and if $e$ has weight greater than 1 then $e$ is a loop. Hence any path $p^\prime$ 
with weighting $w$ is  a $Z$--path with the same initial and terminal vertices as
$p$. 
Given a sequence of $Z$--cancellations the corresponding path $p$
gives rise to  a $Z$--weighting of $G(\mbf z)$. 

\begin{lemma}\label{label-subgraph}
Let $p_0$ and $p_1$ be  $Z$--paths in $G(\mbf z)$ with weightings $w_{p_0}$ and $w_{p_1}$,
respectively, and initial vertex $v=v(\G_s)$, for some basic picture $\G_s$. Let
$\G_i$ be the picture obtained from $\G_s$ by
$Z$--insertion following $p_i$, $i=0,1$. If $w_{p_0}=w_{p_1}$ then the boundary
interval $b_q$ has the same label in $\G_0$ and $\G_1$, for $1\le q\le W_1$.   
\end{lemma} 
\medskip

\noindent{\em Proof.} From Lemma \ref{follow-path} the boundary labels
of $\G_i$ depend only on the number of occurrences of edges of $P$ in 
$p_0$ and $p_1$. As  $w_{p_0}=w_{p_1}$ every edge of $G(\mbf z)$
occurs the same number of times in $p_0$ as in $p_1$.
\medskip

Let $P$ be the path--subgraph of a $Z$--path in $G(\mbf z)$, let 
$ e_1,\ldots,e_t$ be the edges of $P$ and let $L(s_i)$ be the label of
$e_i$. Let $\G_s$ be a basic picture for $(\mbf z,\cL)$ with prime labels
${\mbf{\ovr z}}$, labelled by $\hat\al_s({\mbf{\ovr z}})$, for some solution
$\al_s$ to $\cH$. Let $\lm_1^\prime,\ldots ,\lm_t^\prime$ be elements
of $\Lm$ not occuring in ${\mbf z}$, $\mbf{\ovr z}$ or $\cL$. Let $(h_q,f_q)$ be the
$q$th letter of $\mbf z$: so the label of the boundary interval $b_q$
of $\G_s$ is $h_q\^\al_s(\lm_q)$, for $q\in\{1,\ldots, W_1\}$.
Define the system of parameters $\cL(\G_s,P)$ as follows.
\[f_q=\al_s(\lm_q)+l(h_q)\left(\sum_{i=1}^t \lm^\prime_i\z(q,s_i)
\right),\textrm{ for }1\le q\le W_1,\]
\[\lm_i^\prime=1, \textrm{ if } e_i\textrm{ is not a loop},\]
\[\lm_i^\prime\ge 1, \textrm{ for } i=1,\ldots ,t.\]

\begin{theorem}\label{basic-reduce}
Let $P$ be the path--subgraph of a $Z$--path with initial vertex 
$v(\G_s)$ and let 
$\al_s$ be a solution to $\cH$ such that
$(\al_s,\G_s)$ corresponds to a solution of $Q_\cH$. Then the following
are equivalent.
\be
\item\label{basic-reduce-1}
There exists a solution $\al$ to $\cL$ and a picture $\G$ such
that $(\al,\G)$ corresponds to a solution of $Q$ and such that $\G_s$
is obtained from $\G$ by $Z$--cancellation, along some $Z$--path with path--subgraph
$P$ and initial vertex $v(\G_s)$. 
\item\label{basic-reduce-2} There exists a solution $\al_0$ to $\cL\cup \cL(\G_s,P)$.
\ee
Furthermore $\al$ and $\al_0$ may be chosen so that $\al_0(\lm)=\al(\lm)$, for all $\lm$ 
occuring in $\mbf z$ or $\cL$.
\end{theorem}
\medskip

\noindent{\em Proof.}
Assume (\ref{basic-reduce-1}) holds and let $\G_s$ be obtained from
$\G$ by $Z$--cancellation following the $Z$--path $p$, with path--subgraph $P$. Assume $p$ has
edges $e_1,\ldots, e_t$, let $L(s_i)$ be the label of $e_i$  
and let $n_i$ be the number of occurrences of
edge $e_i$ in $p$, $i=1,\ldots,t$.
Then, using Theorem \ref{follow-path}(\ref{follow-path-1}),  the label of the boundary interval $b_q$ of $\G$ is  
\begin{equation}\label{exp-lab}
h_q\^\al(f_q)=h_q\^\left(m_q+l(h_q)\left(\sum_{i=1}^{t}n_i\z(q,s_i)\right)\right),
\end{equation}
where $m_q=\al_s(\lm_q)$.
Define a retraction $\al_0:M\maps \ZZ$ by $\al_0(\lm)=\al(\lm)$, for 
$\lm\in\Lm\bl\{\lm_1^\prime,\ldots,\lm_t^\prime\}$ and
$\al_0(\lm_i^\prime)=n_i$, for $i=1,\ldots,t$. Then $\al_0(f_q)=\al(f_q)$,
for $1\le q\le W_1$, since $\lm_i^\prime$ does not occur in $\mbf
z$. Hence, from (\ref{exp-lab}), 
\[\al_0(f_q)=m_q+l(h_q)\left(\sum_{i=1}^{t}\al_0(\lm^\prime_i)\z(q,s_i)\right)\]
and, in the light of the above remarks about weightings of $Z$--paths,
\[\al_0(\lm_i^\prime)\ge 1,\] for all $i$, with equality if $e_i$ is not
a loop. Hence $\al_0$ is a solution of $\cL(\G_s,P)$. As
$\al_0(\lm)=\al(\lm)$, when $\lm\notin\{\lm_1^\prime,\ldots,\lm_t^\prime\}$,
it follows that $\al_0$ is a solution to $\cL$ and that the final statement
of the Theorem holds.

Conversely, suppose there exists a solution $\al_0$ to $\cL\cup \cL(\G_s,P)$, where
  $P$ is the path--subgraph of some $Z$--path in $G(\mbf z)$ with
  initial vertex $v(\G_s)$. Let $P$ have edges $e_1,\ldots ,e_t$ and
  let
$L(s_i)$ be the corridor--section labelling $e_i$. 
Then assigning weight
$n_i=\al_0(\lm^\prime_i)$ to edge $e_i$ of $P$ (and 0 to all other edges of $G(\mbf
z)$) gives a $Z$--weighting $w$ of $G(\mbf z)$. Let $p$ be a $Z$--path
with path--subgraph $P$ and weighting $w$: so edge $e_i$ occurs
$w(e_i)=\al_0(\lm^\prime_i)$ times in $p$. Then $p$ has initial vertex
$\G_s$. Let $\G$ be the picture obtained from $\G_s$ by $Z$--insertion
following $p$ and let $m_q=\al_s(\lm_q)$, for $q=1,\ldots ,W_1$. 
Then the label of the boundary interval $b_q$ in $\G$
is
\begin{align*}
& h_q\^\left(m_q+l(h_q)\left(\sum_{i=1}^t n_i\z(q,s_i)\right)\right)\\
&=h_q\^\al_0(f_q),
\end{align*}
as $\al_0$ is a solution to $\cL(\G_s,P)$. Thus $\al_0$ is a solution to
$\cL$ and the boundary label $b_q$ of $\G$ has label
$\hat\al_0(h_q,f_q)$: that is $(\al_0,\G)$ corresponds to a solution
of $Q$.

%% file: in3.tex
\section{Angle assignment} \label{angle}
In the following sections we consider figures showing various sub--configurations of a
picture $\G$. By convention 
every region shown in these figures is planar. 
Solid lines incident to vertices represent edges of the graph $\bar \G$.
Furthermore a dotted line joining two vertices, $u$ and $v$ of a figure,
represents either a single edge of $\bar \G$ joining $u$ and $v$ or 
a sequence $a_0,c_1,a_1,\ldots,c_n,a_n$, where $a_i$ is an
edge of $\bar \G$, $c_i$ is a boundary corner of $\G$, $a_i$ is incident to
$c_{i}$ and $c_{i+1}$, for $i=1,\ldots, n-1$, $a_0$ is incident to $u$
and $c_1$, $a_n$ is incident to $c_n$ and $v$ and $n\ge 1$.
If a vertex (or region) 
of a picture $\G$ appears as vertex (or region) 
$x$ in some
configuration $Y$ we call it a vertex (or region) 
of type $Y(x)$ or just $Y$, if the meaning is clear. Not all arcs or edges incident
to vertices are necessarily shown in the diagrams.
In 
these configurations 
the orientation $\z$ on intervals of $\pd \S$ is shown by an
accompanying arrow: 
for example the orientation
of $\bt$ in configuration $C0^-$ of Figure \ref{C0} is from left to right.

Let $(\mbf h,\mbf n,\mbf t,\mbf p)$ be a positive 4--partition of length $k$, $\cL$
a consistent system of parameters, $\mbf z$ a special element of
$(H^\Lm,\cL,\mbf n,\mbf h)$ and $\al$ a solution to $\cL$.  
Let $\G$ be a picture over $G$, on a compact surface 
$\S $  
of type $(\mbf n,\mbf t,\mbf p)$ with  
prime labels $\mbf z$, 
labelled by $\hat\al(\mbf z)$ and with boundary partition 
$\mbf b=\mbf b_1,\ldots ,\mbf b_{k}$.
We assign 
angles to the corners of regions of $\G $ as follows.
Let 
$\bt$ be a connected component of $\pd\S$ with boundary partition 
$\mbf b_u=(b_{l+1},\ldots , b_{l+\eta})$, for some integers $l$ and $\eta$,  such that 
 $[b_j,b_{j+1}]$ has
prime label $(h_j,m_j)$, with $l(h_j)=d_j$, $j=l+1,\ldots,l+ \eta$
(where $l+\eta+1\equiv l+1$).
Define
$$
\phi_j(\bt) = \frac{d_j(d_j-1)(\al(m_j) -1)\pi}{\al(m_j) ^2}, 
\textrm{ for } j=l+1,\ldots l+\eta,
$$
noting that under the given hypotheses on $\G$ the integer
$\al(m_j)>0$, for all $j$. For notational convenience we also define 
$\phi_i(\bt)=\phi_{l+i}$, $i=1,\ldots ,\eta$,  given $\bt$, $l$ and $\eta$ as above. 
Repeat this for all boundary components of
$\S$ so that $\phi_j$ is defined for $j=1,\ldots, W_1$.
If $[b_j,b_{j+1}]$ is a partisan boundary interval then we have $ \phi_j = 0$.

\begin{figure}
\psfrag{c}{$c$}
\psfrag{bi}{$b_i$}
\psfrag{bi+1}{$b_{i+1}$}
\psfrag{bi+k-1}{$b_{i+k-1}$}
\psfrag{bi+k}{$b_{i+k}$}
\psfrag{...}{$\ldots$}
\begin{center}
{\includegraphics[scale=0.4,clip]{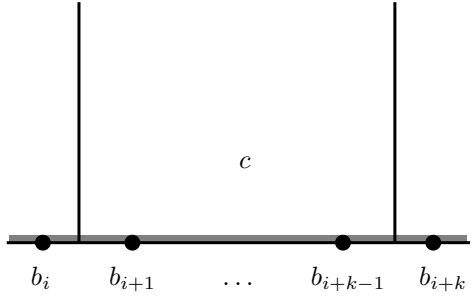}} 
\caption{The angle at $c$ is $-\phi_i-\phi_{i+k-1}$}\label{multi_int_corner}
\end{center}
\end{figure}

We first assign angles to boundary corners of $\G $. Let $\D$ be a region of $\G$ and 
let $c$ be a boundary corner of $\D$ contained in a boundary component
$\bt$ of $\pd\S$.
Assume that $\bt$ has boundary partition $(b_{l+1},\ldots , b_{l+n})$ and that $c$ meets exactly
$k$ boundary intervals in this partition. Then there exists a unique integer 
$i$ such that $c\subseteq [b_i,b_{i+1}]\cup \cdots \cup [b_{i+k-1},b_{i+k}]$
(where subscripts are taken modulo $n$  and lie between $l+1$ and $l+n$).
If $k=1$ then assign angle $\phi(c)=-\phi_i$ to $c$. If $k>1$ then assign angle 
$\phi(c)=-\phi_i-\phi_{i+k-1}$ to $c$.
(See Figure \ref{multi_int_corner}). 
 
We now assign angles to vertex corners
Let $\D $ be a region of $\G $ and $v$ a vertex meeting 
$\pd \D $ on a corner $c$.  If both arcs meeting $v$ at $c$ 
are incident only
to vertices we say that $c$ is a corner of {\em incidence}  0.  If one arc meets 
the boundary $\pd \S $ and the other arc meets only vertices we say $c$ is 
an {\em incidence} 
1 corner and if both arcs at $c$ meet $\pd \S $ we 
say $c$ is an {\em incidence} 
2 corner. 

\begin{figure}
\psfrag{a}{$a$}
\psfrag{b}{$b$}
\psfrag{'b}{$\bt$}
\psfrag{v}{$v$}
\psfrag{'D}{$\D$}
\begin{center}
{\includegraphics[scale=0.4,clip]{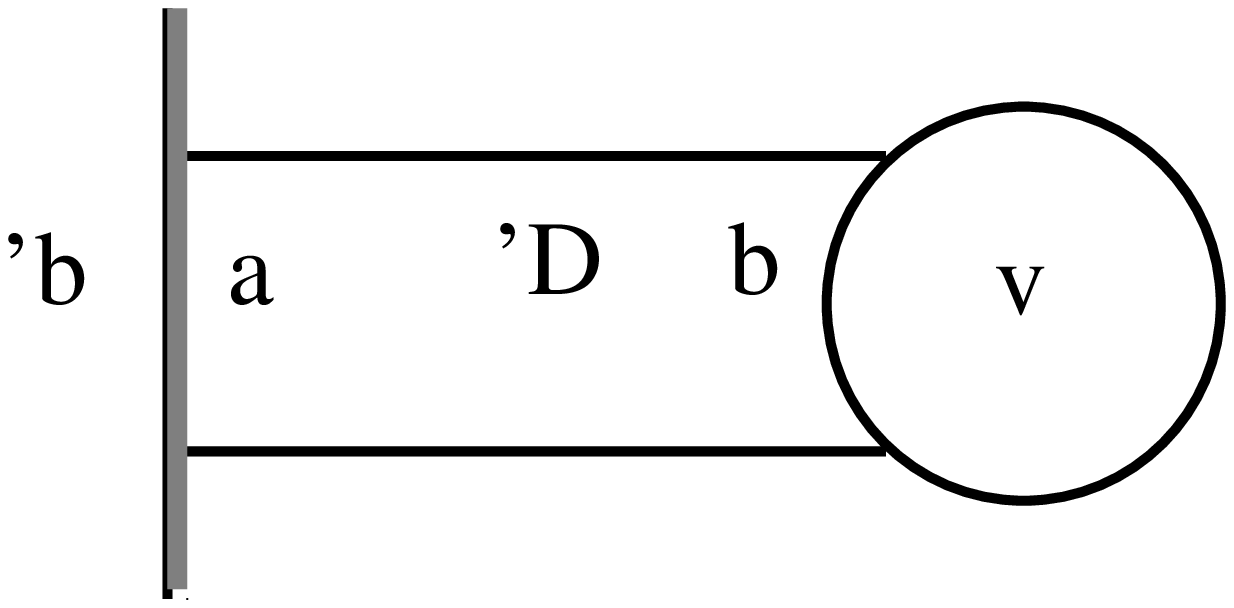}} 
\caption{}\label{bdpara}
\end{center}
\end{figure}

We consider two cases.  
\begin{enumerate}
\item\label{collapsible_corner} First suppose that $\x (\D )=1$, 
$\rho (\D )=1$ and $\bt (\D )=1$, (so $t(\D )=2$ and $\D$ is collapsible), as shown in 
Figure \ref{bdpara}.
Assume that $\D $ meets the boundary component $\bt \subseteq \pd \S $,
where $\bt$ has boundary partition $(b_{l+1},\ldots , b_{l+n})$ and that $\pd\D$ meets
$k$ boundary intervals in this partition. Then there exists a unique integer 
$i$ such that $\pd\D\cap\pd\S\subseteq [b_i,b_{i+1}]\cup \cdots \cup [b_{i+k-1},b_{i+k}]$
(where subscripts are taken modulo $n$ and lie between $l+1$ and $l+n$).
If $k=1$ then assign angle $\phi(c)=\phi_i$ to $c$. If $k>1$ then assign angle 
$\phi(c)=\phi_i+\phi_{i+k-1}$ to $c$. Note that in the latter case 
$[b_j,b_{j+1}]$ is a partisan boundary interval, for $j=i+1,\ldots, i+k-2$, and $v$ is an
$H^\Lm$--vertex.

\item Now consider all those vertex corners $c$ which do not come into the first case 
above.  If $c$ is an incidence $i$ corner we assign an angle of
%
$$
\phi=\phi(c) = \left(  \frac{\rho (\D )-2\x (\D )}{\rho (\D )} + 
\frac{E(\D)i}{2}\right)\pi 
$$
to $c$, where 
$$E(\D) =
\left\{ 
\begin{array}{ll}
1,& \mbox{ if }
\bt(\D)=0;\\
\bt(\D )/(\bt(\D )-\e(\D )),& \mbox{ if }
\rho(\D )\ge \e(\D );\\
(2\bt(\D )-\e(\D ) +\rho(\D )) /2(\bt(\D )-\e(\D )),& \mbox{ if } 
\rho(\D )< \e(\D ).
\end{array}
\right.
$$
\end{enumerate}
\noindent

\begin{figure}
\psfrag{u}{$u$}
\psfrag{v^1}{$v^1$}
\psfrag{v^2}{$v^2$}
\psfrag{x^1}{$x^1$}
\psfrag{x^2}{$x^2$}
\psfrag{'D4^1}{$\D_4^1$}
\psfrag{'D4^2}{$\D_4^2$}
\begin{center}
{\includegraphics[scale=0.6]{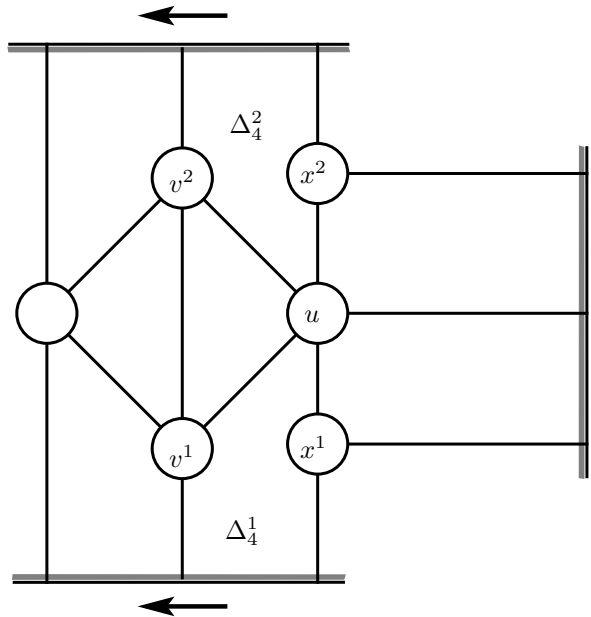}} 
\caption{Configuration $AA$}\label{AA}
\end{center}
\end{figure}

Now we adjust angles on corners of vertices in configuration $AA$ of Figure \ref{AA}.
(The arrows indicate the orientation of boundary 
components.) 
We adjust the angle on the corners of $x^1$ and $x^2$ in 
$\Delta ^1_4$ and $\Delta ^2_4$, respectively, as well as those on $u$ and $v^i$.  
We make the angle on $x^i$ in $\Delta ^i_4$ equal to $4{\pi / 6}$ and
add ${\pi / 9}$ to the corner of $v^i$ in $\Delta ^i_4$ and ${\pi / 18}$ 
to both the corner of $u$ in $\Delta ^1_4$ and in $\Delta ^2_4$.  
\noindent

\begin{figure}
\psfrag{D}{$\D$}
\psfrag{u1}{$u_1$}
\psfrag{u2}{$u_2$}
\psfrag{c1}{$c_1$}
\psfrag{c2}{$c_2$}
\psfrag{c}{$c$}
\psfrag{v}{$v$}
\begin{center}
{\includegraphics[scale=0.6,clip]{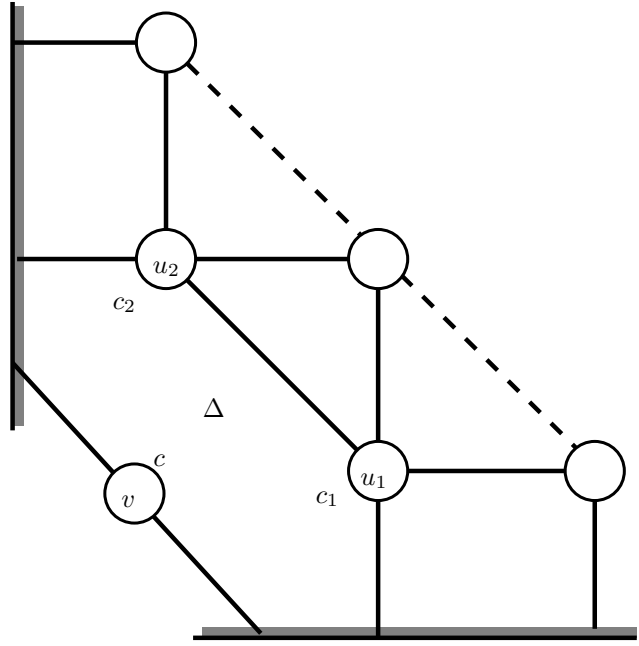}} 
\caption{Configuration $AB$}\label{AB}
\end{center}
\end{figure}

Next angles on corners of vertices in configuration $AB$ of
Figure \ref{AB} are adjusted. 
We adjust the angle $\phi(c)$ from $8\pi /6$ to $13\pi/12$ and angles
$\phi(c_i)$ from $5\pi/6$ to $11\pi/12$, for $i=1,2$. 
\noindent

\begin{figure}
\psfrag{'D1}{$\D_1$}
\psfrag{'D2}{$\D_2$}
\psfrag{'D3}{$\D_3$}
\psfrag{'D4}{$\D_4$}
\psfrag{'D5}{$\D_5$}
\psfrag{'D6}{$\D_6$}
\psfrag{'m1}{$\mu_1$}
\psfrag{'m2}{$\mu_2$}
\psfrag{'n}{$\nu$}
\psfrag{u1}{$u_1$}
\psfrag{u2}{$u_2$}
\psfrag{v}{$v$}
\psfrag{x}{$x$}
\psfrag{y}{$y$}
\begin{center}
\mbox{
\subfigure[Configuration $AC^{+}$\label{AC+}]{
 \includegraphics[scale=0.43,clip]{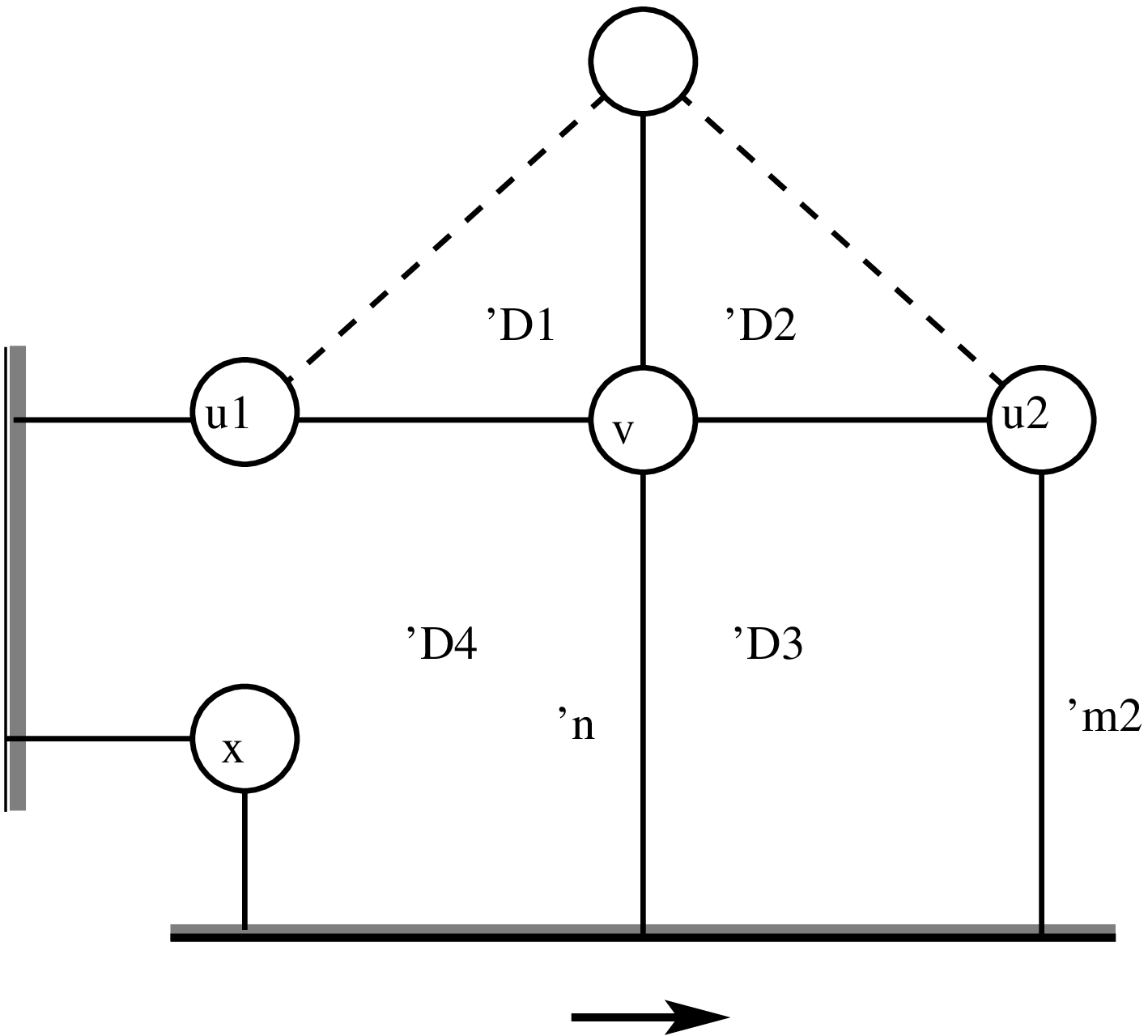}} \qquad
\subfigure[Configuration $AC^{-}$\label{AC-}]{ 
 \includegraphics[scale=0.43,clip]{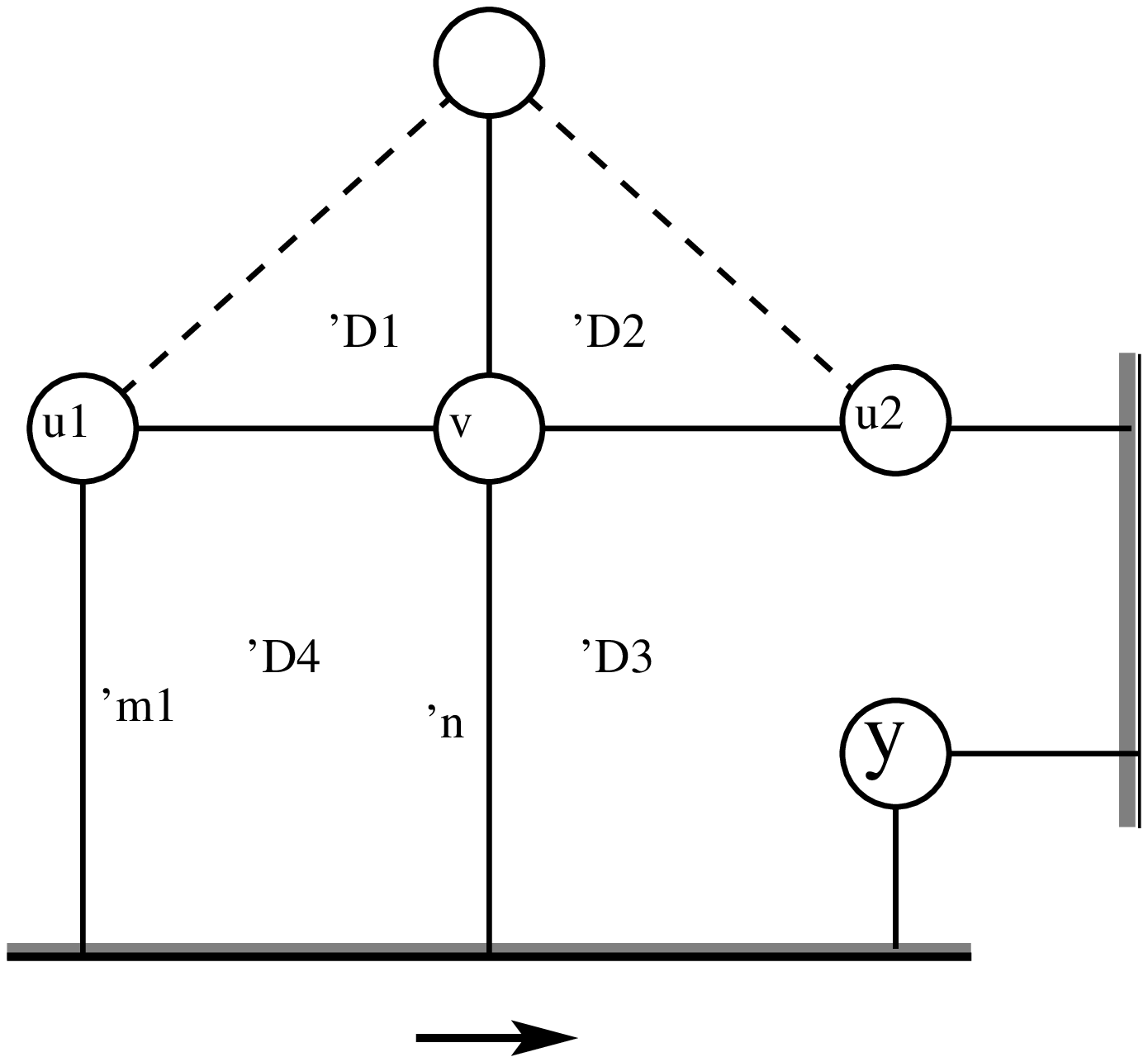} }
} 
\caption{}\label{AC}
\end{center}
\end{figure}

\begin{figure}
\psfrag{'D1}{$\D_1$}
\psfrag{'D2}{$\D_2$}
\psfrag{'D3}{$\D_3$}
\psfrag{'D4}{$\D_4$}
\psfrag{'D5}{$\D_5$}
\psfrag{'D6}{$\D_6$}
\psfrag{'m2}{$\mu_2$}
\psfrag{'n}{$\nu$}
\psfrag{u1}{$u_1$}
\psfrag{u2}{$u_2$}
\psfrag{v}{$v$}
\psfrag{x}{$x$}
\begin{center}
{ \includegraphics[scale=0.4,clip]{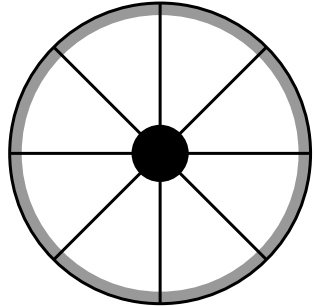} } 
\caption{A single vertex disk}\label{single_vertex_disk}
\end{center}
\end{figure}

Finally we make an 
adjustment of angles assigned to vertices appearing in configurations $AC^\pm$ of 
Figure \ref{AC}.  In configuration $AC^+$ of Figure \ref{AC+} the angle on the corner of vertex $x $ in 
$\D _4$ is decreased from $4\pi/3$ 
to ${7\pi / 6}$ and the angle on the corner of $v$
in $\D _4$ is increased from $5\pi/6$ to $\pi $.  Similarly in configuration 
$AC^-$ of Figure \ref{AC-} the angle on the corner of vertex $y $ in 
$\D _3$ is decreased from $4\pi/3$ 
to ${7\pi / 6}$ and the angle on the corner of $v$
in $\D _3$ is increased from $5\pi/6$ to $\pi $.
\begin{note}
\label{Notes}
\be
\item If $\rho (\D ) \ge  {\e}(\D )$ then $E(\D)\ge 1$.  If 
$\rho (\D )<{\e}(\D )$ then $$E(\D)=1 + 
\left( 
\frac{{\e}(\D )+\rho (\D )}
{2(\bt (\D )-{\e}(\D ))} 
\right)  
>  1,$$ 
as $\bt \le \rho + \e$ and $E(\D)$ is only defined for regions with 
$\rho (\D )\ge 1.$
\item Suppose $c$ is of incidence 1.  Then if $\x (\D )=1$ it follows that 
$\rho (\D )\ge 2$.  Hence, if $c$ is of incidence 1, the angle at $c$ is at 
least $\pi /2$ and is at least $3\pi/2 $ if $\x (\D )\le 0.$
\item Suppose $c$ is of incidence $2$ and not in case (i) of the 
definition of angle above.  If 
$\x (\D )\le 0$ then $\phi \ge 2\pi $.  Assume now that $\x (\D )=1$, 
so $\rho (\D )\ge 1$.  If $\rho (\D )=1$ then $\bt (\D )\ge 2$ 
and ${\e}(\D )=\bt (\D )-1$ so $\phi \ge \pi $.  If 
$\rho (\D )>1$ then again $\phi \ge \pi $.  If $\rho (\D )>2$ then 
$\phi >\pi .$
\item If $\G$ has exactly one vertex and $\S$ is a disk then every vertex
corner falls into case (\ref{collapsible_corner}) of the definition above. We call
such a picture a {\em single vertex disk} (see Figure \ref{single_vertex_disk}).
As minor boundary intervals contribute angles of zero, the sum of angles on the 
vertex  of a single vertex disk may be small: if all the boundary intervals are 
minor the angle sum will be zero. 
\ee
\end{note}
\begin{lemma}
\label{AA0}     
Let $c$ be a vertex--corner of a 
non--collapsible region $\D $ and let $\phi $ be the angle assigned to $c$. 
If $\x (\D )\le 0$ then either 

i) $c$ is incidence $0$ and $\phi  \ge  \pi $ or 

ii) $c$ is incidence $1$ and $\phi  \ge  3{\pi / 2}$ or 

iii) $c$ is incidence $2$ and $\phi  \ge  2\pi $.
\end{lemma}

\noindent
{\em Proof.~}  The result is immediate, from the definition of angles and
Note \ref{Notes}.

\begin{figure}
\psfrag{'f}{$c$}
\begin{center}
  \mbox{
\subfigure[$\bt(\D)=0$\label{P3a}]
{ \includegraphics[scale= 0.4,clip]{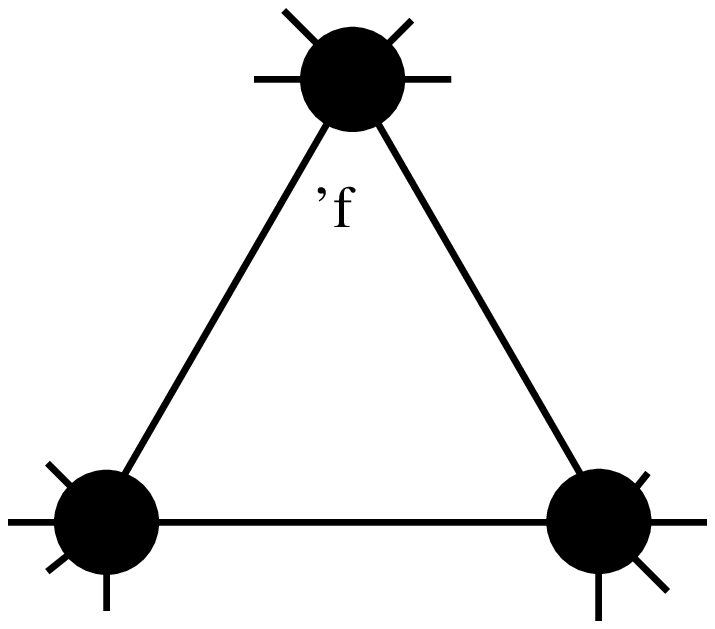} }\qquad
\subfigure[$\bt(\D)=1$\label{P3b}]
{ \includegraphics[scale= 0.4,clip]{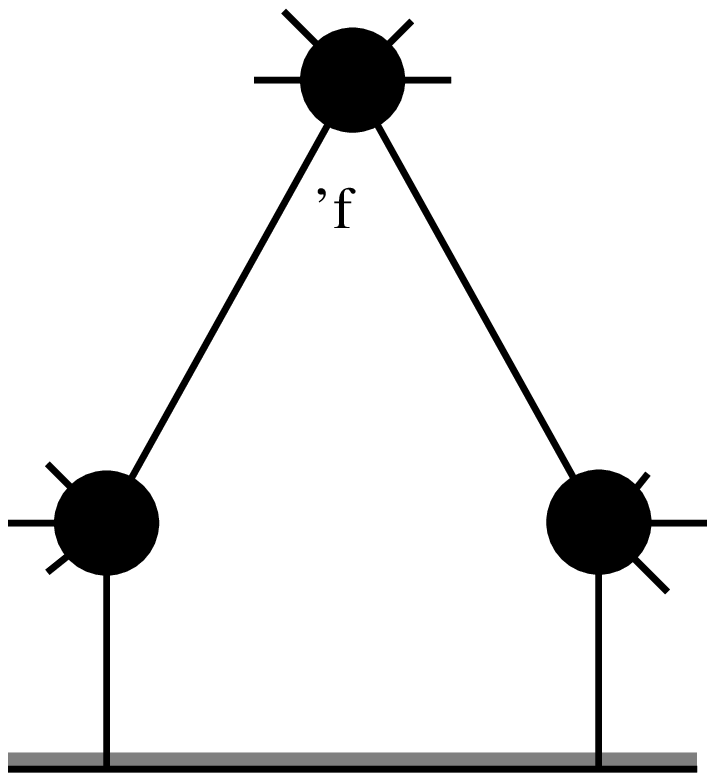} }\qquad
\subfigure[$\bt(\D)>1$, $\e(\D)>0$\label{P3c} ]  
{ \includegraphics[scale=0.4,clip]{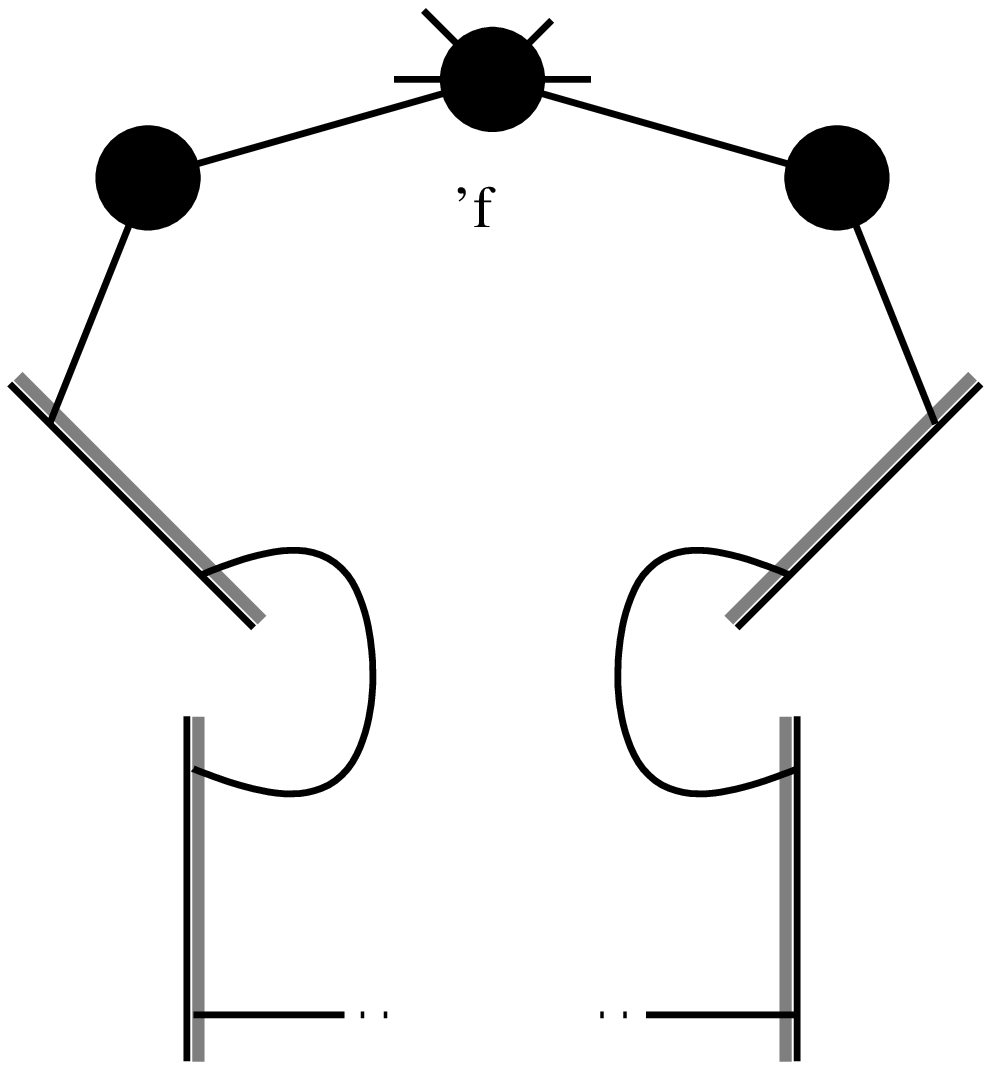} } 
}
\caption{$\rho(\D)=3$ and $\phi(c)=\pi/3$}\label{P3abc}
\end{center}
\end{figure}

\medskip
\noindent
\begin{lemma}
\label{AA10}     
Let $c$ be a vertex--corner 
of a 
non--collapsible region $\D $ and let $\phi=\phi(c)$. 
Suppose $\x (\D ) = 1$ 
and $c$ is of incidence $0$.  Then  
$\rho(\D )\ge 3$ 
and one of the following holds.
\be
\item\label{AA10_1}   $\D$ is of type $AA(\D_4^i)$, $c$ is a corner of a vertex of type $AA(u)$ 
and $\phi=7\pi/18$ (see Figure \ref{AA}).
\item
$$\phi =\left(\frac{\rho(\D )-2}{\rho(\D )}\right)\pi$$ 
and either $\bt(\D )=\e(\D )=0$ or
$0\le \e(\D )\le \bt(\D ) -1$.
In particular
if $\rho(\D )= 3,4,5$ or $6$ then $\phi  = \pi /3,\ \pi /2,\ 3\pi /5$ or
$2\pi /3$, respectively  and $\phi\ge 5\pi /7$ when $\rho(\D ) \ge 7$.
(See Figures \ref{P3a}, \ref{P3b} and \ref{P3c}.)
\ee
\end{lemma}

\noindent
{\em Proof.~}  
We have $\rho (\D )\ge 3$, for otherwise the label 
on $c$ is trivial or $\D $ is collapsible.  
Since (\ref{AA10_1})  follows directly from the definition
of angles we may assume that $\D$ is not of type $AA(\D_4^i)$. 
Hence 
$\phi \ge {\pi /3}$. 
  Clearly ${\e}(\D )=0$ if $\bt (\D )=0$. Otherwise  as $c$ is 
incidence $0$ and $\x(\D )=1$ it follows that 
${\e}(\D )\le\bt (\D )-1$.
The remaining statements follow directly from the definition of angles.
\medskip

\noindent
\begin{figure}
\psfrag{c}{$c$}
\psfrag{u}{$u$}
\psfrag{v}{$v$}
\psfrag{'i}{$\io$}
\psfrag{'m}{$\mu$}
\psfrag{'n}{$\nu$}
\psfrag{pi/2}{$\pi/2$}
\psfrag{5pi/6}{$5\pi/6$}
\psfrag{'D}{$\D$}
\begin{center}
  \mbox{
\subfigure[$\phi(c)=\pi/2$\label{P2}]
{ \includegraphics[scale= 0.4,clip]{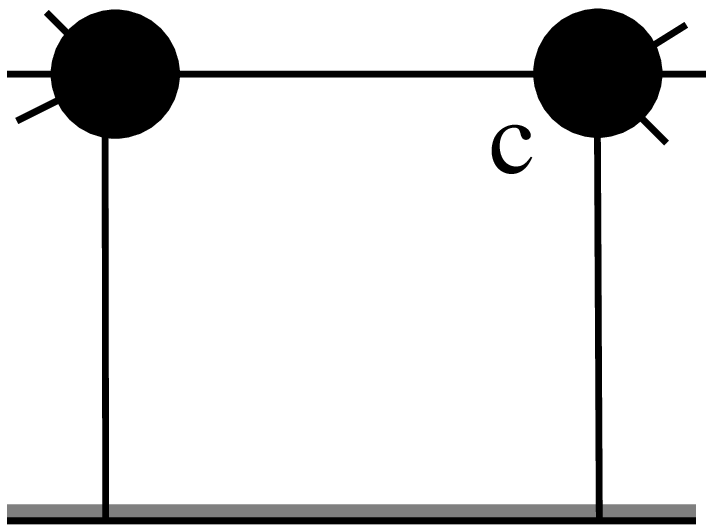} }\qquad
\subfigure[$\phi(c)=5\pi/6$\label{P56a}]
{ \includegraphics[scale= 0.4,clip]{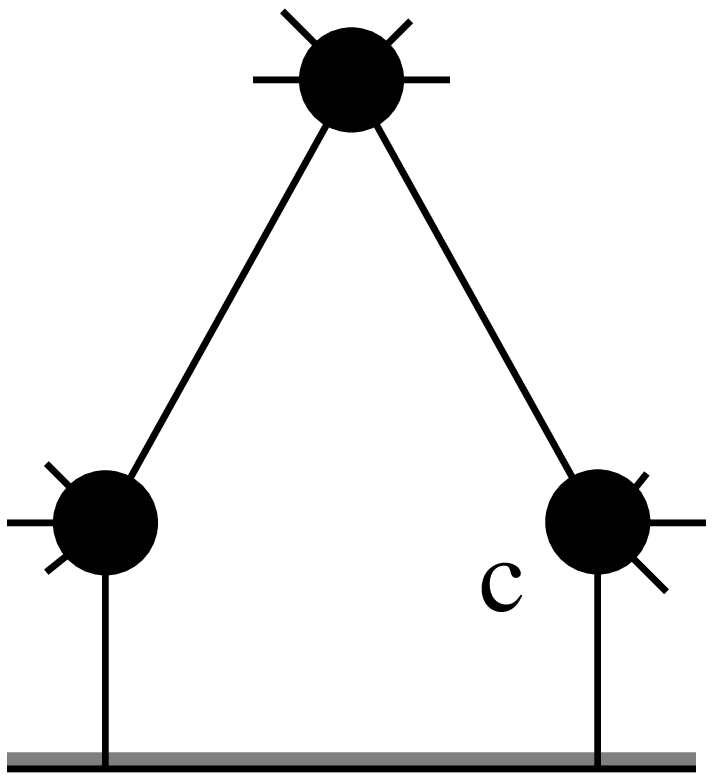} }\qquad
\subfigure[$\phi(c)=5\pi/6$\label{P56b} ]  
{ \includegraphics[scale=0.4,clip]{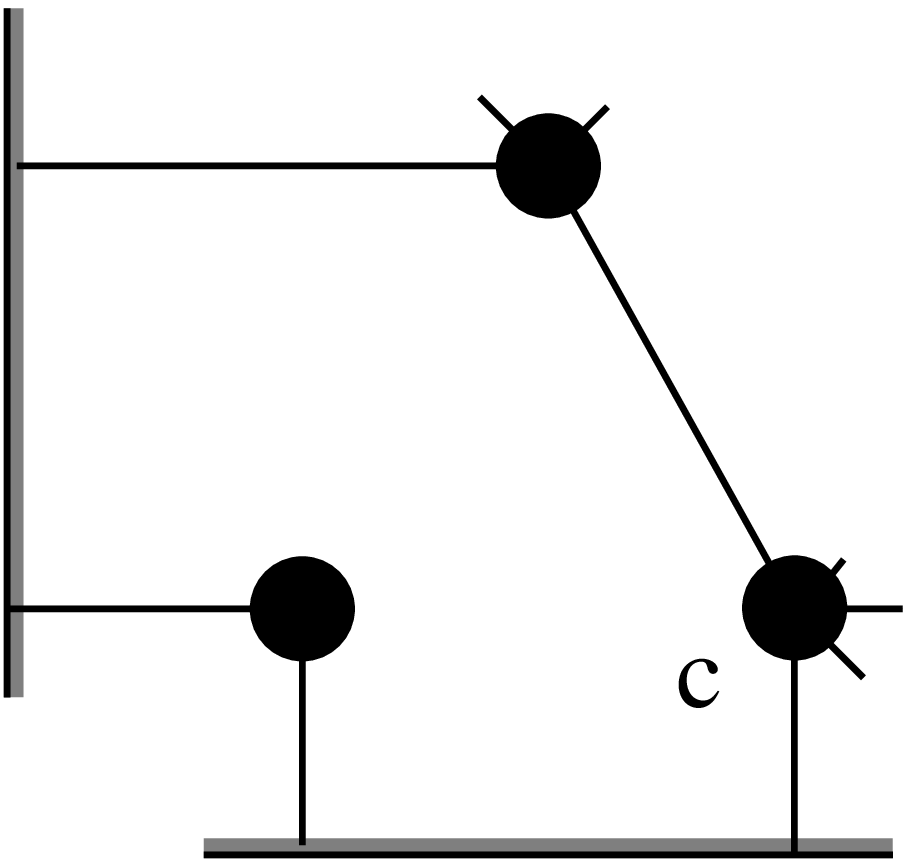} } 
}
\caption{}\label{P2_56}
\end{center}
\end{figure}

\begin{lemma}
\label{AA11}     
Let a vertex $v$ have a corner $c$ 
in a 
non--collapsible region $\D $ and let $\phi $ be the angle assigned to $c$. 
If $\x (\D ) = 1$ and $c$ is incidence $1$ then one of the following holds 
\be
\item\label{AA11_1} $\phi  = {\pi / 2}$, $\rho (\D )=2$,
$\bt (\D )=1$ and ${\e}(\D )=0$ (see Figure \ref{P2}).
\item\label{AA11_2} $\phi=2\pi/3$, $\D$ is of type $AA(\D_4^i)$ and 
$v$ is 
of type $AA(x^i)$ (see Figure \ref{AA}).
\item\label{AA11_3} ${\phi} = 5{\pi / 6}$, $\rho (\D )=3$,
$\bt (\D )=1$ or $2$, 
${\e}(\D )=0$, $\D$ is not of type $AB(\D)$ 
and if $\D$ is of type $AC^{+}(\D_4)$ then 
$v$ is of 
type $AC^+(u_1)$ and if $\D$ is of type $AC^-(\D_3)$ then  $v$ is of 
type $AC^-(u_2)$
(see Figures \ref{P56a}, \ref{P56b}, \ref{AB} and \ref{AC}).
\item\label{AA11_4} $\phi =11\pi/12$, $\D$ is of type $AB(\D)$ and $c$ is a corner of 
type
$AB(c_i)$, $i=1$ or $2$ (see Figure \ref{AB}) or
\item\label{AA11_5} $\phi =17\pi/18$, $\D$ is of type $AA(\D_4^i)$ and 
$v$ is of 
type $AA(v^i)$, $i=1$ or $2$ (see Figure \ref{AA}).
\item\label{AA11_6} $\phi\ge \pi$. 
\ee
\end{lemma}

\noindent
{\em Proof.~} 

 As $\x (\D )=1$ we assume 
$\g (\D )=0$ (or $\G $ is an empty picture on the disk).  As $c$ is
incidence $1$, $\bt (\D )>0$ and $\rho (\D )\ge 2$.  Therefore 
$\phi \ge E(\D){\pi / 2}\ge\pi/2$. If $\phi ={\pi / 2}$ then 
$\rho (\D )$ must equal $2$ and $E(\D)$ must equal 1.  Since 
$\bt (\D )>{\e}(\D )$ this occurs only if 
$\rho (\D )\ge {\e}(\D )$ and 
$\bt (\D )={\e}(\D )+1=1$ as required.  

If $\phi >{\pi / 2}$ 
then $\rho (\D )>2$ or $E(\D)>1$.  
If $\rho (\D )=2$ then 
${\e}(\D ) = \bt(\D )-1$. Hence if
$2 = \rho (\D )\ge \e(\D )$ then $E(\D)>1$ implies that
$E(\D)=\bt (\D ) = 2$ or  $3$.
On the other hand if ${\e}(\D )>\rho (\D )$ 
then\ $E(\D)=(\bt (\D )+3)/ 2$ and $\bt (\D )\ge 4$.  In either 
case $\phi \ge {E(\D)\pi / 2}\ge \pi$.

Now suppose that $\rho (\D )>2$. If $\D$ is of type $AA(\D_4^i)$  
then either (\ref{AA11_2}) or (\ref{AA11_5}) holds. 
If $\D$ is of type $AB(\D)$ then (\ref{AA11_4}) holds. 
If $\D$ is not of type $AB(\D)$ but $\D$ is of type $AC^+(\D_4)$ or type $AC^-(\D_3)$ then
either (\ref{AA11_3}) or (\ref{AA11_6}) holds. 
Hence we may assume that $\D$ is not of type $AA(\D_4^i)$, $AC^+(\D_4)$, $AC^-(\D_3)$ or $AB(\D)$.
Thus
$\phi \ge ({1/ 3} +
{E(\D)/ 2})\pi \ge 5{\pi / 6}$.  
If $\rho(\D ) \ge 4$ then $\phi \ge (1/2 +E(\D) /2)\pi \ge \pi$.
Hence we assume $\rho(\D )=3$.
Now $c$ is incidence $1$, $\x(\D )=1$
and $\rho(\D )=3$, so $\e(\D )=\bt(\D)-1$ or $\bt(\D)
-2$.
If $3=\rho(\D) < \e(\D)$ then 
$$E(\D)
=1+ 
\left(\frac{\e(\D)+3}
{2(\bt(\D)-\e(\D))} \right) \ge 1 +7/4 > 2,$$
so $\phi > \pi$.
Hence we may assume that $\e(\D) \le \rho(\D)$
and $\rho(\D) = 3$.

Now $\phi\ge 5\pi /6$ with equality only if $E(\D)=1$. 
We have $\bt (\D )={\e}(\D )+1$ or 
$\bt (\D )={\e}(\D )+2$ and $E(\D) =\bt (\D )$
or $E(\D)=\bt (\D )/2$, respectively. Hence $E(\D)=1$ implies
$\bt (\D )=1$ and 
${\e}(\D )=0$
or 
$\bt (\D )=2$ and 
${\e}(\D )=0$.  These 
two possibilities give $\phi =5{\pi / 6}$ and configurations of Figures 
\ref{P56a} and \ref{P56b}, respectively. 
Furthermore if $E(\D)> 1$ then $E(\D)\ge 3/2$ so $\phi \ge 13\pi /12 > \pi$.
\medskip

\noindent
\noindent
\begin{figure}
\psfrag{c}{$c$}
\begin{center}
  \mbox{
\subfigure[\label{P1a}]
{ \includegraphics[scale= 0.5,clip]{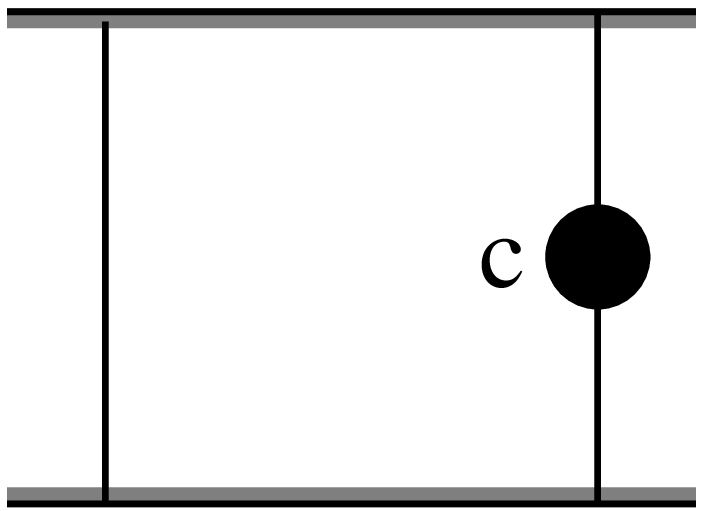} }\qquad
\subfigure[\label{P1b}]
{ \includegraphics[scale= 0.5,clip]{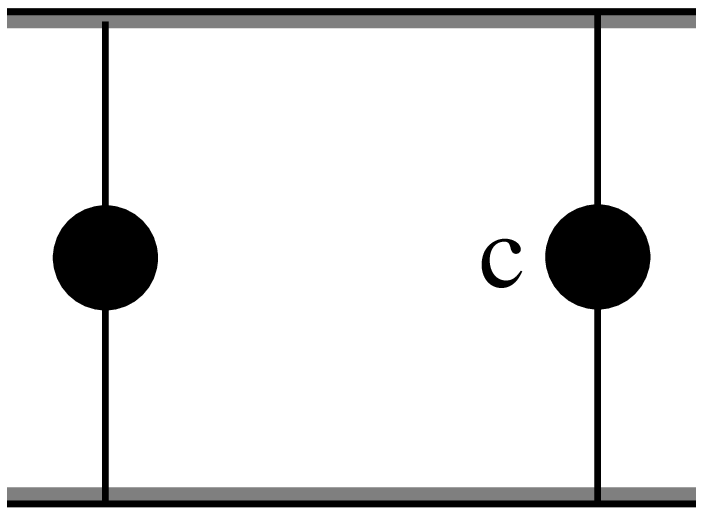} }
}
\caption{$\phi(c)=\pi$}\label{P1}
\end{center}
\end{figure}

\begin{lemma}
\label{AA12}
Let a vertex $v$ have a corner $c$ 
in a   
non--collapsible region $\D $ and let $\phi $ be the angle assigned to $c$. 
If $\x (\D ) = 1$ and $c$ is incidence $2$  then one of the following holds. 
\be
\item $\phi = \pi $,
$\bt (\D )=2$, $\rho (\D )=1$ or $2$ and 
${\e}(\D )=\bt (\D ) - \rho (\D )$ (see Figures 
\ref{P1a} and \ref{P1b}).
\item\label{AA12_2} $\phi=13\pi/12$ and $c$ is a corner of type $AB(c)$ (see Figure \ref{AB}).
\item \label{ugly} $\phi=7\pi/6$, $c$ is not a corner of type $AB(c)$ and either
 $\D$ is of type $AC^+(\D_4)$ and $v$ is of type $AC^+(x)$
or  $\D$ is of type $AC^-(\D_3)$ and $v$ is of type $AC^-(y)$
(see Figures \ref{AC} and \ref{AB}).
\item $\phi \ge 4{\pi / 3}.$ 
\ee
\end{lemma}

\noindent
{\em Proof.~}   If $\rho (\D )=1$ then as 
$\x (\D )=1$, $\bt (\D )\ge 2$ and 
${\e}(\D )=\bt (\D )-1$. Hence if
$\bt (\D )=2$ then 
$\e (\D)=1 $ and $\phi=\pi$, as shown in Figure \ref{P1a}. 
Otherwise $\bt(\D)\ge 3$ and $\rho(\D)=1<\e(\D)$
so
$$\phi =\left(-1+\frac{\bt(\D )+2}{2}\right)\pi \ge {3\pi / 2}.$$

Suppose then that $\rho (\D )\ge 2$. If $\rho(\D)>2$ then 
$\phi\ge 4\pi /3$ unless either case (\ref{AA12_2}) or (\ref{ugly}) holds. Hence we may assume $\rho(\D)=2$.  In this 
case $\bt (\D )={\e}(\D )+2$.  If 
$\rho (\D )\ge {\e}(\D )$ then 
$$E(\D)=\frac{\bt (\D )}{ 2}\ge 1,$$
with equality only if $\bt (\D )=2$ and ${\e}(\D )=0$.  If 
$2=\rho (\D )<{\e}(\D )$ then 
$$E(\D)=\frac{\bt (\D )+2+\rho (\D )}{4} =
\frac{\bt (\D )}{4} + 1
\ge  2.$$ 

Hence $\phi = E(\D)\pi \ge \pi $ with equality only if $\rho (\D )=2$, 
$\bt (\D )=2$ and ${\e}(\D )=0$, as shown in Figure 
\ref{P1b}.  If $\phi >\pi $  and $\rho(\D)>2$
then, from the above, 
$$\phi \ge\frac{3\pi}{2}>\frac{4\pi}{3}.$$  
\medskip

\noindent




\section{Curvature} \label{curvature}

Let $(\mbf h,\mbf n,\mbf t,\mbf p)$ be a positive 4--partition, $\cL$
a consistent system of parameters, $\mbf z$ a special element of
$(H^\Lm,\cL,\mbf n,\mbf h)$ and $\al$ a solution to $\cL$.  
Throughout this section  $\G$ is assumed to be a reduced, minimalistic picture over $G$
on a compact surface $\S$ 
of type $(\mbf n,\mbf t,\mbf p)$ with 
 boundary partition 
$\mbf b=(b_1,\ldots ,b_{W_1})$,
prime labels $\mbf z$ and 
labelled by $\hat\al(\mbf z)$. Moreover we assume that $\G$ is not a single vertex disk.

Given a vertex $v$ we denote the sum of angles on all corners of $v$ by 
$\s (v)$.  Similarly the sum of angles over all vertex and boundary
corners of a region 
$\D $ is denoted $\s (\D )$ and over all corners of a boundary 
component $\bt $ of $\pd \S$ by $\s (\bt ).$ 

We define the {\em curvature} of the vertex $v$ of $\G $ to be
$$
\ka (v) = 2\pi -\s (v)\hbox{ ,}
$$

\noindent
of the region $\D $ to be
$$
\ka (\D ) = \s (\D ) - t(\D )\pi + 2\pi \x (\D )
$$

\noindent
and of the boundary component $\bt $ of $\pd \S $ to be 
$$
\ka (\bt ) = -\s (\bt )\hbox{ .}
$$
\medskip
\noindent

\begin{lemma}
\label{bd-curve}
Let $\bt $ be a boundary component of $\pd \S $ with prime label
 $\prod_{i=1}^n(h_i,m_i)$.  Then 
$\ka (\bt )\le \pi\sum_{i=1}^n l (h_i)(l (h_i)-1) .$
\end{lemma}
\medskip

\noindent
{\em Proof.~}  Let $l (h_i)=d_i$.  We have  
$l(h_i\^\al(m_i))= \al(m_i)$, since $\mbf z$ is special.
Suppose that $\bt$ has boundary partition $b_1,\ldots ,b_n$. A boundary
interval $I=[b_i,b_{i+1}]$ meets $|[b_i,b_{i+1}]|$ boundary corners.
If $I$ meets a boundary corner $c$ then $I$ contributes $-\phi_i(\bt)$ to the 
angle of $c$. Summing over all boundary intervals of $\bt$ we see that 
\[\s(\bt)=\sum_{i=1}^n -\phi_i(\bt)|[b_i,b_{i+1}]|.\]
As $\phi_i=0$ whenever $[b_i,b_{i+1}]$ is a minor boundary interval we have
\[\sum_{i=1}^n -\phi_i(\bt)|[b_i,b_{i+1}]|=\sum_{i=1}^n -\phi_i(\bt)\al(m_i),\]
so 
\[
\ka (\bt ) = -\s (\bt ) =\sum_{i=1}^n \al(m_i)\phi_i(\bt)\le
\pi\sum_{i=1}^n \frac{d_i(d_i-1)}{\al(m_i)}\le \pi\sum_{i=1}^n d_i(d_i-1).
\] 

\medskip

\noindent
\begin{lemma}
\label{r-curve}
If $\D $ is a collapsible region of $\G $ then $\ka(\D)=0$ when $\rho(\D) =1$ or $2$.
When $\rho(\D) =0$ then $\ka(\D)\le 0$,
with equality only if all boundary intervals meeting $\D$ are minor.
If $\D $ is a non--collapsible region of $\G $ then 
\begin{enumerate}
\item if $\bt(\D) = 0$ then $\ka (\D )\le 0$;
\item if $\bt(\D) > 0$ and $\rho(\D) =0$ then $\ka(\D) \le -\pi$;
\item if $\bt(\D) > 0$ and $\rho(\D) > 0$ then 
$$\ka(\D) \le \left\{ \begin {array}{ll}
                       0 & \mbox{ if $\rho(\D) \ge \e(\D)$} \\
                       -\pi/2 & \mbox{ if $\rho(\D) < \e(\D)$}
                      \end{array}
              \right. .$$
\end{enumerate}
\end{lemma}
\medskip
\noindent
{\em Proof.~}   The result follows immediately 
from the definitions for collapsible
regions so we can assume that $\D$ is non--collapsible. 
First suppose that $\D $ is an interior region; that is
$\bt(\D)=0$.  Then 
$\rho (\D )=t(\D )$ so, as $\G$ is minimalistic, 
\[\ka (\D ) = \rho (\D ) 
\left(\frac{\rho (\D )-2\x (\D )}{\rho (\D )}\right) \pi  - 
\rho (\D )\pi  + 2\pi \x (\D ) = 0.\]
Now suppose $\D $ is a boundary region ($\bt(\D)>0$).  
Note that the angle adjustments made to corners of vertices of 
types $AC^\pm(v)$, $AC^+(x)$, $AC^-(y)$, $AB(v)$, $AB(u_i)$, $AA(x^i)$, $AA(v^i)$ and 
$AA(u)$  have no effect on curvature of regions so
may be ignored in this context.
If $\rho(\D )=0$ then
$\s(\D )\le 0$ and so $\ka(\D )\le (-t(\D )+ 2\x(\D ))\pi$.
If $\rho(\D )=0$ and  $\x(\D)\le 0$ then $\ka(\D)<-\pi$. 
If $\rho(\D )=0$ and  $\x(\D)=1$ then, as $\G$ has no trivial
boundary labels, $t(\D)\ge 2$. As $\D$ is non--collapsible we have  
$t(\D)>2$ and so $\ka(\D) \le -\pi$. 
Hence we may assume that $\bt(\D)>0$, $\rho(\D)>0$ and
$t(\D )=\rho (\D )+\bt (\D )$.  Let $\rho_i=|\{$corners of $\D $
of incidence $i\}|$, for $i=0,1,2$.  We have $\rho (\D )=\rho_0+\rho_1+\rho_2$ and 
$$
\s (\D ) \le  \rho (\D ) 
\left( \frac{\rho (\D )-2\x (\D )}{\rho (\D )}\right) \pi  + 
\left( \rho_1/ 2 + \rho_2\right)E(\D) \pi \hbox{ ,}
$$
since all angles on boundary corners are non--positive.  Now ${\rho_1/ 2} + \rho_2
= \bt (\D )-{\e}(\D )$ so
$$
\s (\D) \le  (\rho(\D)-2\x (\D))\pi  + 
(\bt (\D )-{\e}(\D ))E(\D)\pi .
$$
Thus $\ka (\D )=\s (\D ) - (\rho (\D )+\bt (\D ))\pi  +
2\x (\D )\pi $
$$
\le  (\bt (\D ) - {\e}(\D ))\pi E(\D) - \bt (\D )\pi .
$$
If $\rho (\D ) \ge  {\e}(\D )$\ then\ $$E(\D) = 
\frac{\bt (\D )}{\bt (\D )-{\e}(\D )}$$\ and so 
$\ka (\D )\le 0.$
If $\rho (\D ) < {\e}(\D )$ then $$E(\D) = 
\frac{2\bt (\D )-{\e}(\D )+\rho (\D )}{2(\bt (\D )-%
{\e}(\D ))}$$
so that $$\ka (\D ) \le
\frac{(\rho (\D )-{\e}(\D ))\pi}{2}
\le  - \pi /2.
$$
This completes the proof of Lemma \ref{r-curve}.
\medskip
\noindent

\begin{figure}
\psfrag{s1}{$\s_1$}
\psfrag{s2}{$\s_2$}
\psfrag{v}{$v$}
\begin{center}
  \mbox{
\subfigure[$\s_1\neq\s_2$\label{angle_bdry_classes_a}]
{ \includegraphics[scale=0.55,clip]{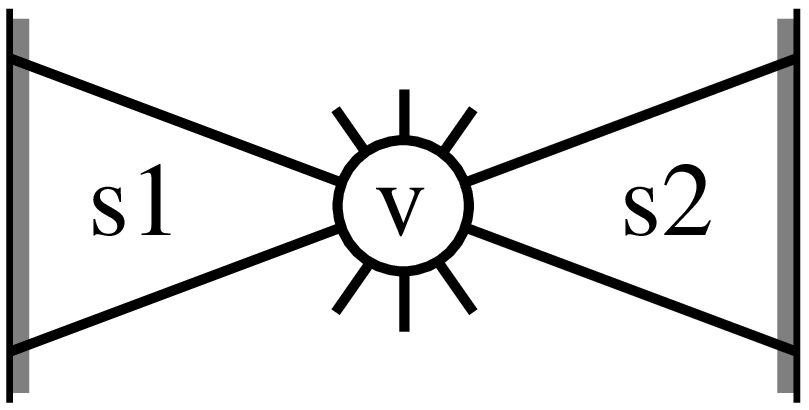} }\qquad
\subfigure[$\s_1=\s_2$\label{angle_bdry_classes_b}]
{ \includegraphics[scale= 0.55,clip]{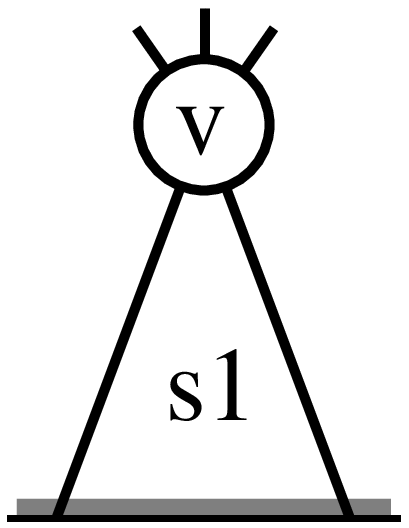} }
}
\caption{}\label{angle_bdry_classes}
\end{center}
\end{figure}

%
%
\begin{lemma}
\label{bd-angle}
Let $v$ be a vertex incident to two (not necessarily distinct)
boundary classes of arc 
$\s _1$ and $\s _2$ (see Figure \ref{angle_bdry_classes}).  Then the sum of angles on corners of $v$ between
$\s _1$ and $\s _2$ is at least $\pi .$
\end{lemma}
\medskip

\noindent
\noindent
{\em Proof.~}  Suppose that $\s _1$ and $\s _2$ are not adjacent classes of arc
at $v_i$: that is, there is at least one class of arcs separating them.  The 
angle on $v$ between a boundary class and any other class is at least 
${\pi / 2}$, from Lemmas \ref{AA11} and \ref{AA12}.  Hence, in this case, $\s _1$ and $\s _2$ are 
separated by an angle of at least $\pi $.  On the other hand if $\s _1$ and
$\s _2$ are adjacent classes then, using Lemmas \ref{AA0}
and \ref{AA12}, the angle on $v$ 
between them is less than $\pi $ only if they meet $v$ at the corner of a 
collapsible region. 
In this
case $\s _1$ and $\s _2$ form 1 class and $\G$ is a single vertex disk, a contradiction.  Hence in all 
cases the angle between $\s _1$ and $\s _2$ is at least $\pi .$
\begin{corol}
\label{bd-c-curve}   
Let $v$ be  a boundary vertex. Then $\ka(v)\le \pi$. If 
$v$ meets $\pd \S $ in more than one 
class of arcs then $\ka (v)\le 0$.  If $v$ meets $\pd \S $ in more 
than two classes of arcs then $\ka (v) \le  -\pi .$
\end{corol}

Let $\t$ be a boundary class of arcs and suppose that one end
of $\t$ is contained in a boundary interval $[b,c]$ of some boundary
component $\bt$ of $\pd \S$.
In this case we say that $\t$ {\em meets} $\pd\S$ in $[b,c]$.

\begin{lemma}
\label{wv}
Let $v$ be a 
boundary vertex of $\G $ incident to a class of arcs $\t$ which meets
 $\pd\S$ in a boundary interval $[b,c]$ with prime label $(w,f)$.
If the width of $\tau $ is more than ${ml / 2}+1$ then $w$ is a cyclic subword
of $r^{\pm m}$, 
$(w,f)$ is a proper exponential $H$--letter, the width of $\tau $ is at most ${ml / 2}+l-1$,  
and 
${ml / 2} - l  + 3 \le  
l (w) \le  {ml / 2}$.
\end{lemma}
\medskip
\noindent
{\em Proof.~}  The class $\tau $ identifies a
subword $s$ of $r^{\pm m}$ which is also a subword of $w\^\al(f)$.  To 
see this fix a transverse orientation of $\tau $.  Starting at the first edge 
of $\tau $ in this orientation we can read round $v$ or along $\bt $ to the 
last edge of $\tau $.  In the former case we obtain a subword $s$ of 
$r^{\pm m}$ and in the latter of $w\^\al(f)$.  Note that $l (s) =  
|\tau | - 1$, where $|\tau |$ denotes the width of $\tau $.
Hence $\al(f)\ge l (s)
\ge  {ml / 2} + 1$.  Let $d = l (w)$.  If $d \ge  {ml / 2} +
1$ then it follows that $w$ has a cyclic subword $s^\prime $, which is a cyclic subword
of $s$, such that $l (s^\prime ) \ge  {ml / 2} + 1$ and $s^\prime $ is a 
cyclic subword of $r^{\pm m}$.  Therefore $(w,f)$ is not relator--reduced, contradicting 
the hypothesis that $\mbf z$ is special. 
Hence $l(w)=d \le  {ml / 2}$ and so $\al(f)>d$ 
and in particular $(w,f)$ is a proper
exponential $H$--letter. 

Let $\delta  = gcd(l ,d).$ 
First suppose that $d + l  - \delta  \le  l (s)$.  Then since $s$ has 
periods $l $ and $d=l (w)$, it follows from \cite[Section 5, Proposition 1]{Howie89},  
that $s$ has period $\delta $.  This implies that 
$\delta =l $ so that $l | l (w)$ and $w$ is a cyclic permutation of 
$r^{\pm t}$ for some $t> 0$, $t\in \ZZ$.  As $(w,f)$ is  $H^\Lm$--irredundant 
it follows that $t=1.$ 
However this means that $(w,f)$ is relator--unconstrained, 
a contradiction. 

\begin{figure}
\psfrag{'t1}{$\t_1$}
\psfrag{'t2}{$\t_2$}
\psfrag{a}{$a$}
\psfrag{'D}{$\D$}
\psfrag{'D1}{$\D_1$}
\psfrag{'D2}{$\D_2$}
\begin{center}
  \mbox{
\subfigure[Configuration $I1$\label{IV1}]
{ \includegraphics[scale= 0.5,clip]{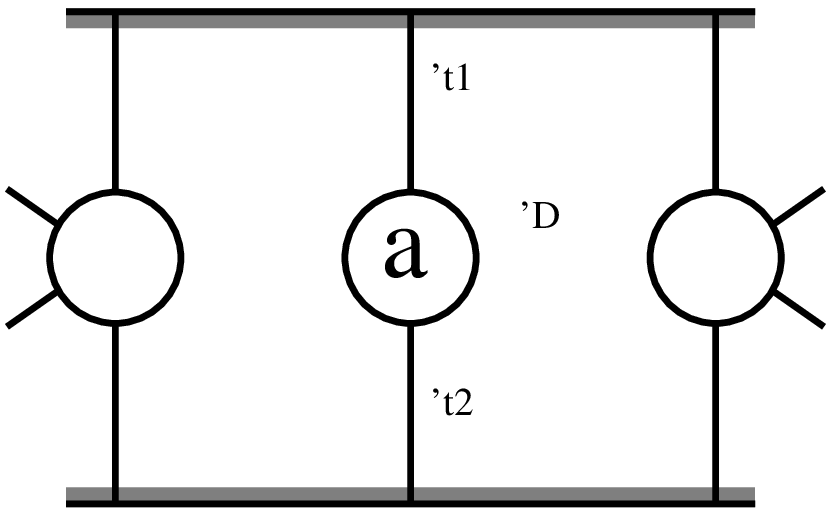} }
\subfigure[Configuration $I2$\label{IV2}]
{ \includegraphics[scale= 0.5,clip]{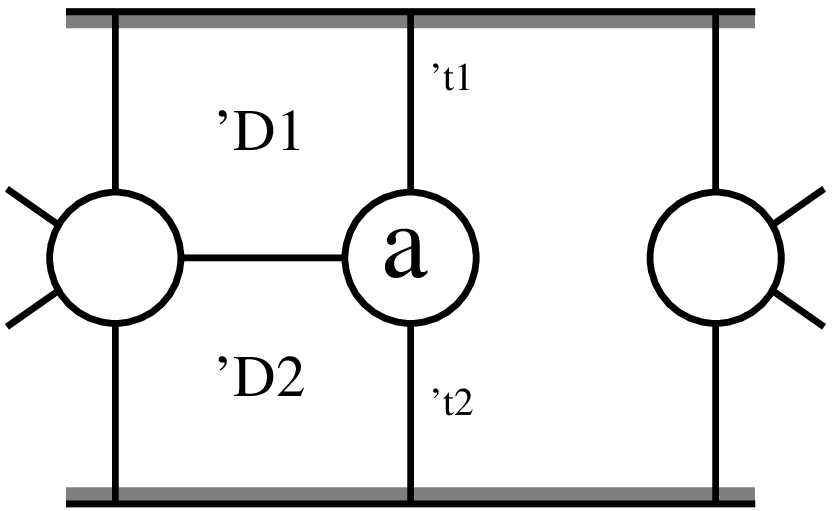} }
\subfigure[Configuration $I3$\label{IV3} ]  
{ \includegraphics[scale=0.5,clip]{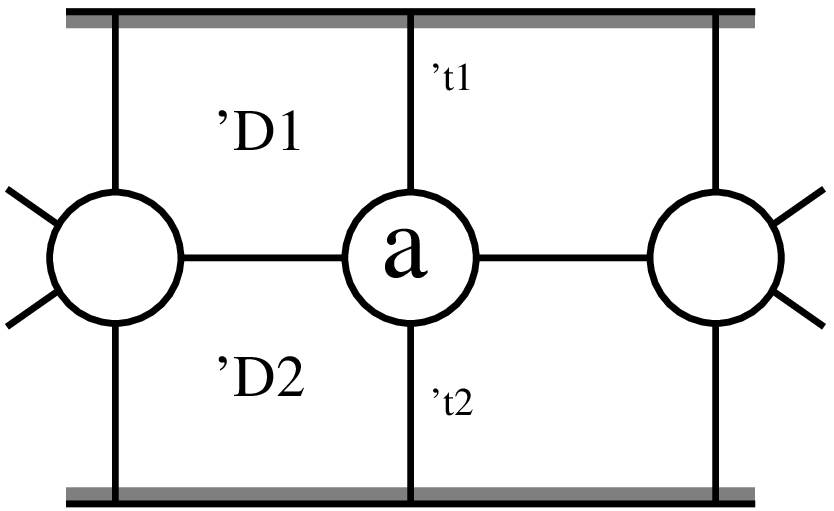} } 
}
~\\
 \mbox{
\subfigure[Configuration $I4$\label{IV4}]
{ \includegraphics[scale= 0.5,clip]{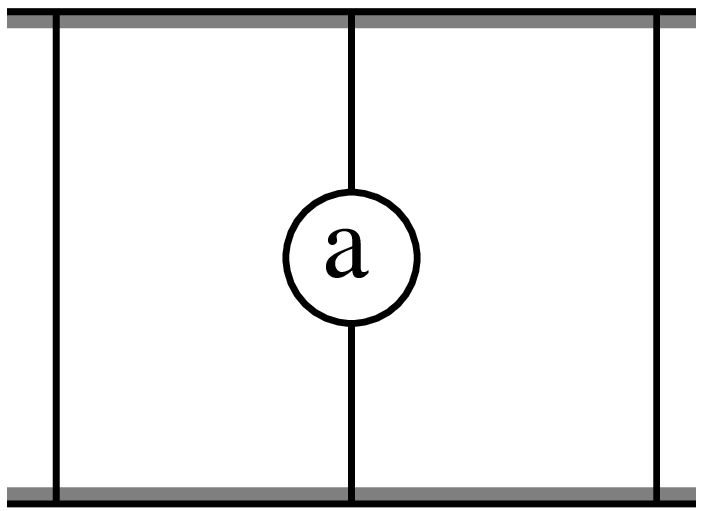} }\quad
\subfigure[Configuration $I5$\label{IV5}]
{ \includegraphics[scale= 0.5,clip]{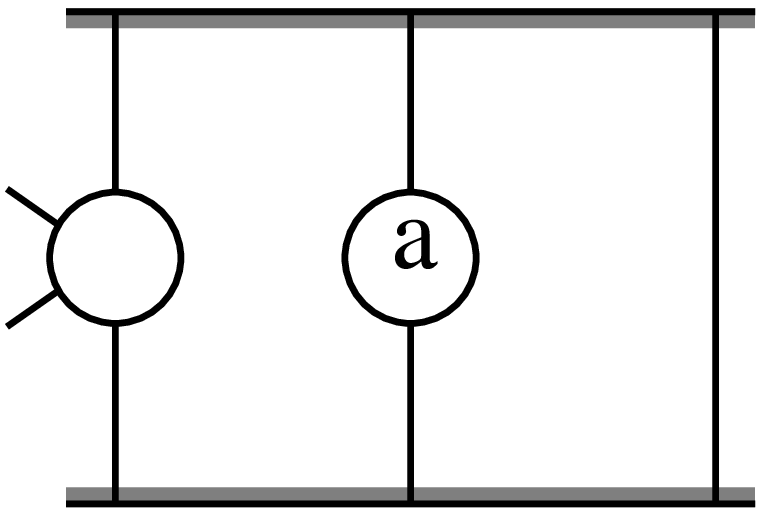} }\quad
\subfigure[Configuration $I6$\label{IV6} ]  
{ \includegraphics[scale=0.5,clip]{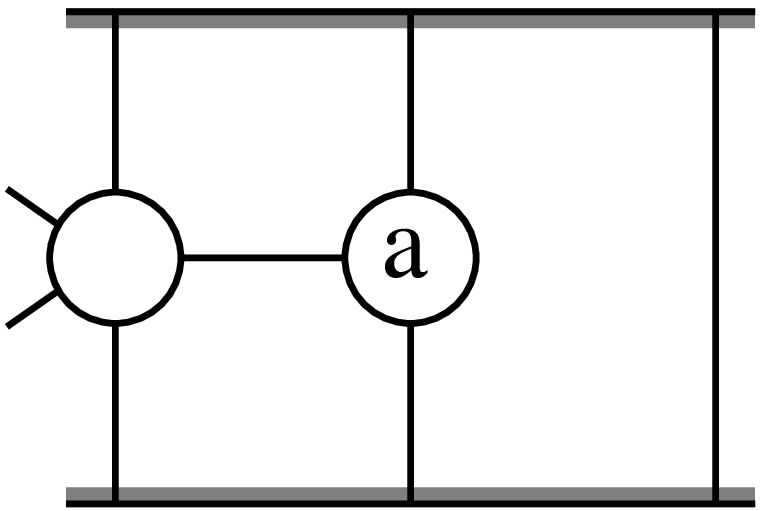} } 
}
\caption{}\label{IV123456}
\end{center}
\end{figure}

Hence we may assume that $d+l -\delta >l (s)$.
As $\delta \ge 1$  and $l(s) \ge ml/2 +1$ this implies that 
$d+l - 1>l (s)\ge {ml / 2}+1$, so
$d\ge {ml / 2}-l +3$. Furthermore 
$l (s)<{ml / 2}+l -\delta $
so $|\t| \le {ml / 2}+l-1.$ 
This completes the proof of Lemma \ref{wv}.
\begin{prop}
\label{prop-I}   
If $v$ is a vertex of $\G $ meeting $\pd \S $ in 
exactly 2 classes of arc then either 
$\ka (v)\le -{\pi / 12}$ or 
$v$ appears as the vertex $a$ in a configuration of 
type I1, I2, I3, I4, I5 or I6 as shown in Figures
\ref{IV1}, \ref{IV2}, \ref{IV3}, 
\ref{IV4}, \ref{IV5} or \ref{IV6} 
respectively, or their mirror images.
\end{prop}
\medskip
\noindent
{\em Proof.~}  The condition that 
$\ka (v)>-{\pi / 12}$ implies $\s (v)<25{\pi / 12}$ which 
in turn implies that the boundary classes incident at $v$ are separated by 
angles less than 
$13\pi/12$, using Lemma \ref{bd-angle}.  
Let the boundary classes at $v$ be 
$\tau _1$ and $\tau _2$ and consider the angle $\phi $ between $\tau _1$
and $\tau _2$ (in the direction of orientation of $\pd v )$.  
Suppose $v$ is incident to some class of arcs $\tau $ 
between $\tau _1$ and $\tau _2$.  Then it follows,
from Lemmas \ref{AA0} to \ref{AA12},
as $\phi<13\pi/12$
that $\t$ is the unique class of arcs between $\t_1$ and $\t_2$
and that the regions $\D _i$ incident to $v$, $\tau $ and $\tau _i$, 
for $i=1$ and 2, are discs with $\rho (\D _i)=2$, $\bt (\D _i)=1$ and
${\e}(\D _i)=0$.  To see that $v$ belongs to one of the given 
configurations it remains to check the case in which there is no class of arcs 
incident at $v$ between $\tau _1$ and $\tau _2$.  In this case it follows, from
Lemmas \ref{AA0} to \ref{AA12}, that 
as $\phi<13\pi/12$
the region 
$\D $ incident to $v,\tau _1$ and $\tau _2$ is a disc 
with $\bt (\D )=2$ and 
$\rho (\D )=1$ or 2.  Hence $v$ must lie in one of the given 
configurations.
\medskip 
\noindent

\noindent 
If $v$ is a vertex with $\ka (v)>-{\pi / 12}$ and $v$ is the vertex $a$ in 
one of the configurations $Ij$ we say $v$ is an $Ij${\em --vertex} or $v$ is of 
{\em type} $Ij(a).$ 

From now on we make the following assumption.
\begin{assum} \label{A6}   The integer $m$ is a least $6$  and $G$ 
is not of type 
$E(2,*,m).$
\end{assum}
\begin{theorem} \label{int-v}
Let $v$ be an interior vertex of $\G$. Then 
$\ka(v)\le -\pi/3$.
\end{theorem}
\medskip

\noindent
{\em Proof.~} 
Since $\G$ is minimalistic and $G$ is not of type $E(2,*,m)$, no interior class of arcs
can have width more than $l-2$. 
It follows that $v$ is
incident to at least $7$ classes of arcs. Therefore $\ka(v)\le -\pi/3$.
\medskip

We now consider vertices meeting $\pd \S $ in precisely one class of 
boundary arcs.  Suppose $v$ is such a vertex.  Then we shall show that if $v$ 
has curvature $k(v) > - {\pi / 24}$ and $v$ is not an 
$H^\Lm$--vertex
then $v$ belongs to 
one of a finite number of configurations.

Let $v$ be a vertex incident to exactly one boundary class of arcs.  Then $v$ 
is incident to $\deg (v)-1$ interior classes of arc, (where $\deg (v)$ denotes 
the degree of $v$ in the graph $\bar{\G}$ associated to $\G$).  Each 
interior class of arcs incident to $v$ is separated from its neighbour by an 
angle of at least ${\pi / 3}$ so if $\deg (v)>5$ we have $\ka (v) < - 
{\pi / 3}$. \ Hence we may assume $\deg (v) \le  5.$ 

\begin{figure}
\psfrag{'D1}{$\D_1$}
\psfrag{'D2}{$\D_2$}
\psfrag{'D3}{$\D_3$}
\psfrag{'D4}{$\D_4$}
\psfrag{'D5}{$\D_5$}
\psfrag{a1}{$a_1$}
\psfrag{a2}{$a_2$}
\psfrag{a3}{$a_3$}
\psfrag{a4}{$a_4$}
\psfrag{a5}{$a_5$}
\psfrag{'pi/2}{$\pi/2$}
\psfrag{'pi/3}{$\pi/3$}
\psfrag{a}{$a$}
\begin{center}
{ \includegraphics[scale=0.6]{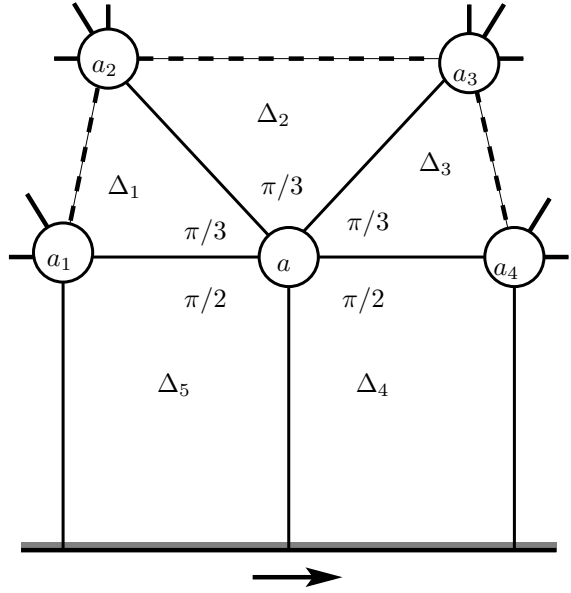} } 
\caption{Configuration $B4$}\label{B4}
\end{center}
\end{figure}

\begin{lemma}
\label{lb4}     
Let $v$ be a vertex of $\G $ 
incident to at least 4 interior and precisely 1 boundary class of arcs.  Then 
either  
$\ka (v) \le  - {\pi / 9}$ or $v$ is
of type $B4(a)$ as shown in Figure \ref{B4}.
\end{lemma}

\noindent
{\em Proof.~}  
If $\deg (v)
 \ge  6$ then  $\ka (v) \le  - 
{\pi / 3}$ so
 we may assume $\deg (v) \le  5$.  Then, from Lemmas  \ref{AA10} and
 \ref{AA11}, 
$\ka (v) \le  0$ with equality only if $v$ is of type $B4(a)$, as shown in Figure 
\ref{B4}, and not of type $AA(u)$.  If $v$ is of type $AA(u)$ then $k(v)=-\pi/9$. 
If $v$ is not of the type $B4(a)$
then $\ka (v)<0$  and,  from Lemmas  \ref{AA10} and \ref{AA11}
 again, either
one of the incidence 0 corners at $v$ has angle $\phi \ge {\pi / 2}$ or one of 
the incidence 1 corners has angle $\phi \ge 5{\pi / 6}$.  In either case 
$\ka (v) \le  - {\pi / 6}.$

\noindent

\noindent
\begin{corol}
\label{lc4}     
Let $v$ be a vertex of 
$\G $ meeting $\pd \S $ in precisely one class of arc of width at 
most ${ml / 2}+1$.  Then either $\ka (v) 
\le  - {\pi / 9}$ or $v$ is of type $B4(a)$ as shown in Figure \ref{B4}.
\end{corol}
\medskip

\noindent
{\em Proof.~}  There are at least ${ml / 2} - 1$ non--boundary arcs incident 
at $v$.  Since $\G$ is minimalistic, ${ml / 2} - 1 \ge  3l  - 1$ and $G$ is not of type 
$E(2,*,m)$ there are at least $4$ non--boundary classes of arc incident at $v$.  
The result follows from Lemma \ref{lb4}.

\begin{figure}
\psfrag{s}{$s$}
\psfrag{u}{$u$}
\psfrag{v}{$v$}
\psfrag{t}{$t$}
\psfrag{'m}{$\mu$}
\psfrag{'n}{$\nu$}
\psfrag{'D}{$\D$}
\begin{center}
\includegraphics[scale= 0.6,clip]{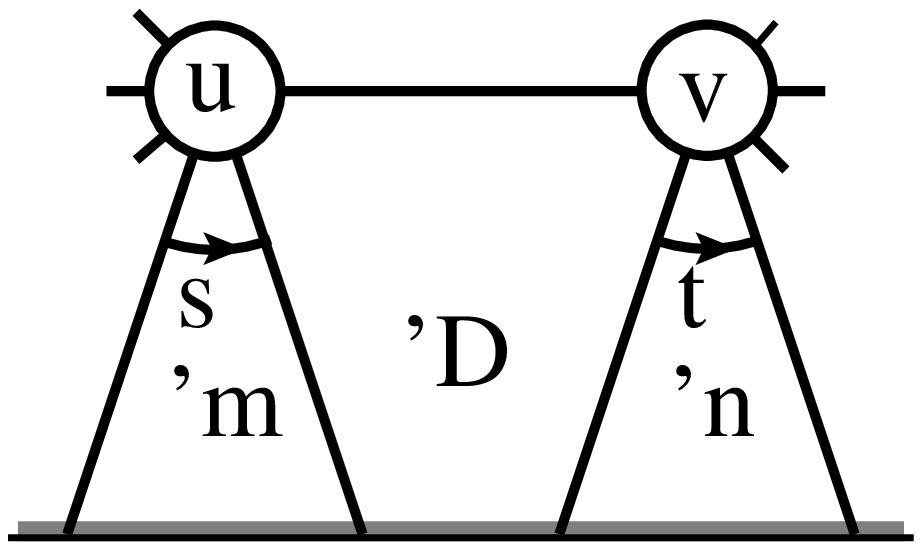} 
\caption{}\label{P2_big}
\end{center}
\end{figure}

\begin{lemma}
\label{uv}     
Let $u$ and $v$ be adjacent 
vertices of $\G $ with incidence $1$ corners on a region $\D $ such that 
$\rho (\D )=2$, $\bt (\D )=1$, $\e(\D )=0$ and 
$\x (\D )=1$, as shown in Figure \ref{P2_big}.  Let  $\mu $ and $\nu $ 
be the boundary classes of arcs incident to $u$ and $v$, respectively, and 
meeting $\D $ and assume that $\pd\S\cap \mu$, $\pd\S\cap \nu$ and the boundary
corner of $\D$ are all contained in a boundary interval $\bt=[b,c]$ with prime
label $(w,f)$.  Assume that 
$\ka (v)> - {\pi / 2}$ and that $ |\nu | \ge  {ml / 2} + 2$.
 Then $|\mu | \le  l$ and $w$ is a cyclic subword of 
$r^{\pm m}$.  
\end{lemma}
\medskip

\noindent
{\em Proof.~}  The class $\mu $ meets $u$ in $|\mu |$ arcs.  
Reading round $u$ in 
the direction determined by the transverse orientation of $\mu$ induced
from the chosen orientation $\z$ of $\pd \S $ these arcs identify a 
subword $s$ of $r^{\e m}$, where $\e=\pm 1$. 
Similarly the class $\nu $ meets $v$ and 
identifies a subword $t$ of $r^{\e m}$, where $\e=\pm 1$.  
For the duration of this proof we
call $\mu$ ($\nu$) {\em positive}  if $\e=1$ and {\em negative} otherwise. 
We have $l (s)=|\mu |-1$ and 
$l (t)=|\nu |-1$.  Assume that $|\nu|\ge ml/2+2$ and that 
$|\mu |\ge l +1.$ 

We may assume, without loss of generality, that the orientation $\z$ of $\bt$ is
such that $\D \cap \bt $ is 
oriented from $\mu $ to $v$ and that an occurrence of $w$ in the label of 
$[b,c]$ begins in the boundary corner following the first arc of $\nu $ (in this 
orientation of $\bt )$.  Then as above $\mu $ identifies a 
cyclic subword $x$ of $w\^\al(f)$ and $\nu$ identifies an initial subword 
$y$ of $w\^\al(f)$.  Furthermore in $H$ we have $x=s$ and $y=t$.  From Lemma 
\ref{wv} it follows that $w$ is a cyclic subword of $r ^{\pm m}$ and $2l  + 3 \le  {ml / 2} - l 
+ 3 \le  l (w) \le  {ml / 2}$, so $w$ is an initial subword of $y$.
Let $x ^\prime $ be the terminal segment of $x $ of length $l $ and let $y^\prime $ 
be the segment of $y$ beginning in position $l (w)-l $ and of length 
$l $.  Then both $x ^\prime $ and $y^\prime $ are subsegments of $w$ beginning in 
position $l (w)-l $ and ending in position $l (w)-1$.  Hence $x^\prime =y^\prime $, 
in $H$.  As $x=s$ and $y=t$ both $x^\prime $ and $y^\prime $ are cyclic permutations of 
$r^{\pm 1}$. 

We consider separately the cases in which $\mu$ is positive and 
negative, under the assumption that $\nu$ is positive. The case $\nu$ negative
follows by symmetry. 
First suppose $\mu$ is positive so $x^\prime $ and $y^\prime $ are both cyclic permutations of 
$r$.  As cyclic subwords of $r$ we may assume the position of the first letters
of $x^\prime $ and $y^\prime $ are $\lambda \le l /2$ apart. (Otherwise interchange $x^\prime $
and $y^\prime $.) The union of $x ^\prime $ and $y^\prime $
then forms a cyclic subword $z$ of $r^m$ of length 
$2l -\lambda $, and of periods $l $ and $\lambda $.  Let 
$\delta =\gcd(l ,\lambda ).$
Then $l  + \lambda  - \delta  \le  {3l / 2} - 1 < 2l  - 
\lambda  = l (z)$. \ Hence from \cite[Section 3, Proposition\ 1]{Howie89}
the word $z$ has period $\delta $.  As $r$ is not a proper power it follows 
that $\delta =l $ and since $\delta|\lambda $ and 
$\lambda \le l /2$ we have $\lambda =0$.  Thus $x^\prime $ and $y^\prime $ are 
identically positioned as cyclic subwords of $r.$ 

Now suppose that the last letter of $w$, as a word in $A\ast B$, is $\w $. 
Then the letter following $y^\prime $ in $y$, and therefore in $r^m$, is $\w $.  
Hence the letter following $x ^\prime $ in $r^m$ is $\w $ and so the label on 
the corner of $u$ in $\D $ is $\w$.  Therefore the label on the corner
of $\D$ meeting $u$, read with orientation induced by that on $\bt$,  is $\w^{-1}$. As $w$ ends in the letter 
$\w $ the boundary component of $\D $ is labelled $\w $.  Hence the
corner of $v$ in $\D $ has trivial label, a contradiction.  Hence $\mu$ 
cannot be positive. 

If $\mu$ is negative then $x ^\prime $ is a cyclic permutation of $r^{-1}$ and $y^\prime $ is a 
cyclic permutation of $r$.  As $x^\prime =y^\prime $ in $A\ast B$ this means that $G$ is of 
type $E(2,2,m)$, contrary to the standing assumption \ref{A6}.
Therefore $|\mu |\le l .$ 

\section{Configurations C}\label{config_C}

Let $(\mbf h,\mbf n,\mbf t,\mbf p)$ be a positive 4--partition, $\cL$
a consistent system of parameters, $\mbf z$ a special element of
$(H^\Lm,\cL,\mbf n,\mbf h)$ and $\al$ a solution to $\cL$.  
As in Section \ref{curvature}, $\G$ is assumed, throughout this section, 
to be a reduced, minimalistic picture over $G$
on a compact surface $\S$ 
of type $(\mbf n,\mbf t,\mbf p)$ with  
boundary partition 
$\mbf b=(b_1,\ldots ,b_{W_1})$,
prime labels $\mbf z$ and 
labelled by $\hat\al(\mbf z)$.
We are now in a position to prove the following.

\begin{figure}
\psfrag{'D1}{$\D_1$}
\psfrag{'D2}{$\D_2$}
\psfrag{'D3}{$\D_3$}
\psfrag{'D4}{$\D_4$}
\psfrag{'a4}{$\al_4$}
\psfrag{'a5}{$\al_5$}
\psfrag{'b}{$\bt$}
\psfrag{'n}{$\nu$}
\psfrag{'m1}{$\mu_1$}
\psfrag{'m2}{$\mu_2$}
\psfrag{u1}{$u_1$}
\psfrag{u2}{$u_2$}
\psfrag{v}{$v$}
\begin{center}
\includegraphics[scale= 0.5,clip]{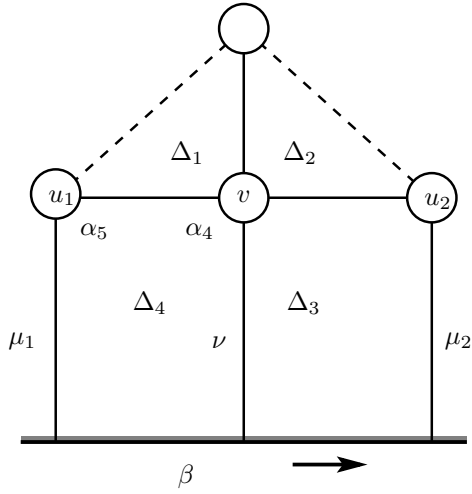} 
\caption{Configuration $C0$}\label{C0}
\end{center}
\end{figure}

\begin{figure}
\psfrag{'D1}{$\D_1$}
\psfrag{'D2}{$\D_2$}
\psfrag{'D3}{$\D_3$}
\psfrag{'D4}{$\D_4$}
\psfrag{'a4}{$\al_4$}
\psfrag{'a5}{$\al_5$}
\psfrag{'b}{$\bt$}
\psfrag{'n}{$\nu$}
\psfrag{'m1}{$\mu_1$}
\psfrag{'m2}{$\mu_2$}
\psfrag{u1}{$u_1$}
\psfrag{u2}{$u_2$}
\psfrag{v}{$v$}
\begin{center}
\mbox{
\subfigure[Configuration $C1^+$\label{C1+}]{
\includegraphics[scale= 0.5,clip]{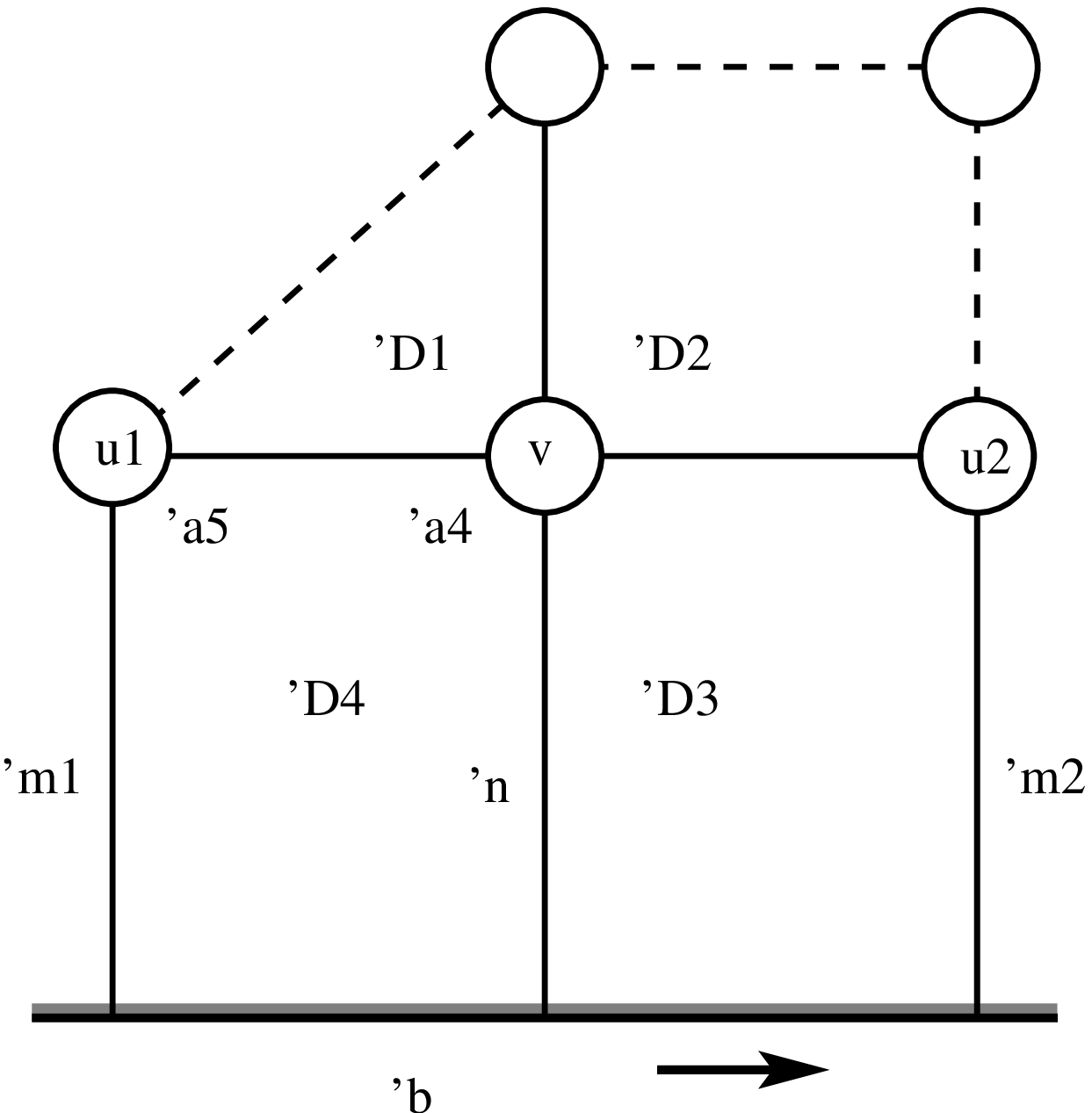}
}\qquad
\subfigure[Configuration $C1^-$\label{C1-}]{
\includegraphics[scale=0.5,clip]{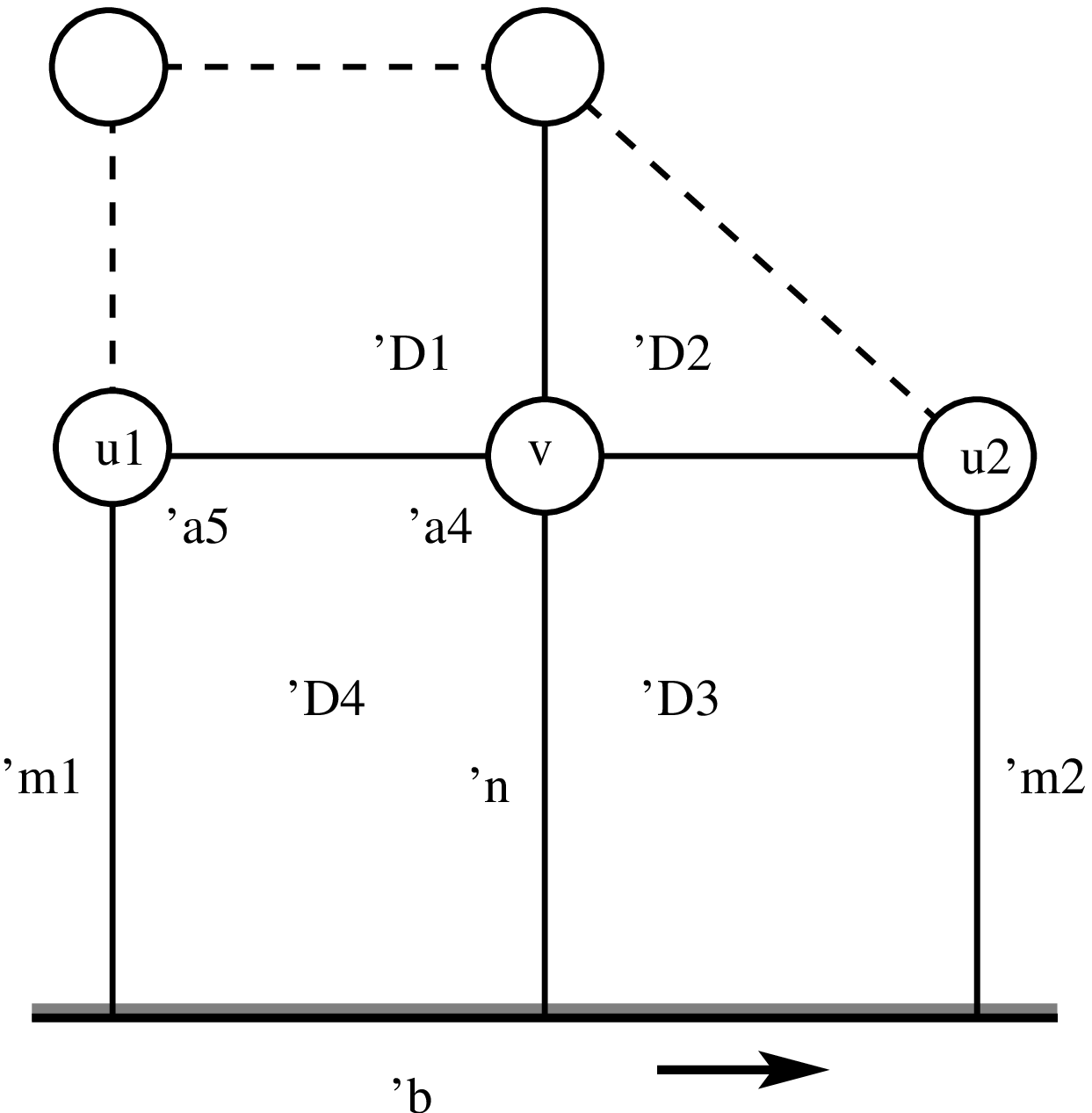}
}
}
\caption{}\label{C1}
\end{center}
\end{figure}

\begin{figure}
\psfrag{'D1}{$\D_1$}
\psfrag{'D2}{$\D_2$}
\psfrag{'D3}{$\D_3$}
\psfrag{'D4}{$\D_4$}
\psfrag{'a4}{$\al_4$}
\psfrag{'a5}{$\al_5$}
\psfrag{'b}{$\bt$}
\psfrag{'n}{$\nu$}
\psfrag{'m1}{$\mu_1$}
\psfrag{'m2}{$\mu_2$}
\psfrag{u1}{$u_1$}
\psfrag{u2}{$u_2$}
\psfrag{v}{$v$}
\begin{center}
  \mbox{
\subfigure[Configuration $C2^+$\label{C2+}]
{ \includegraphics[scale= 0.5,clip]{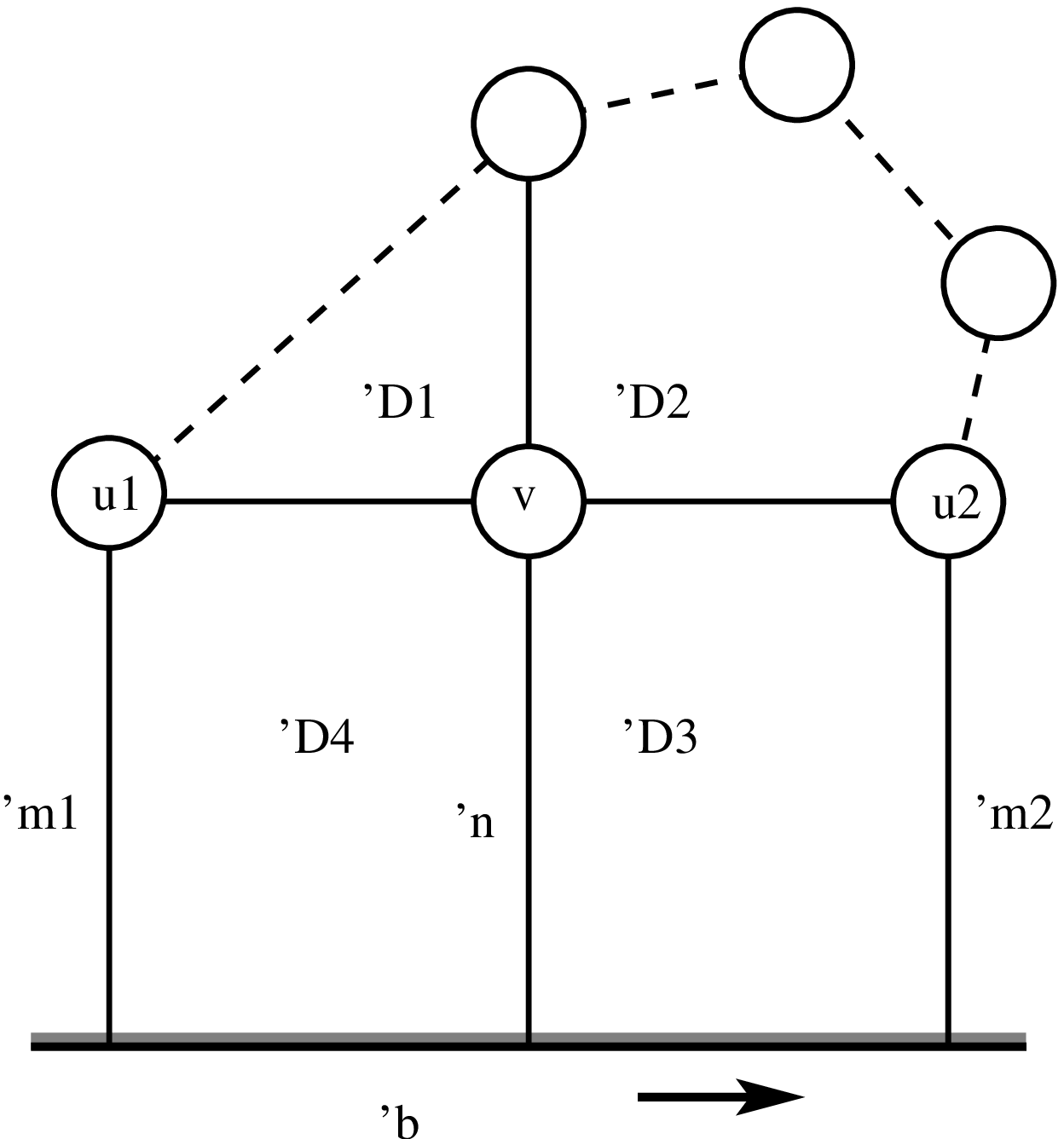} }\qquad
\subfigure[Configuration $C2^-$\label{C2-} ]  
{ \includegraphics[scale=0.5,clip]{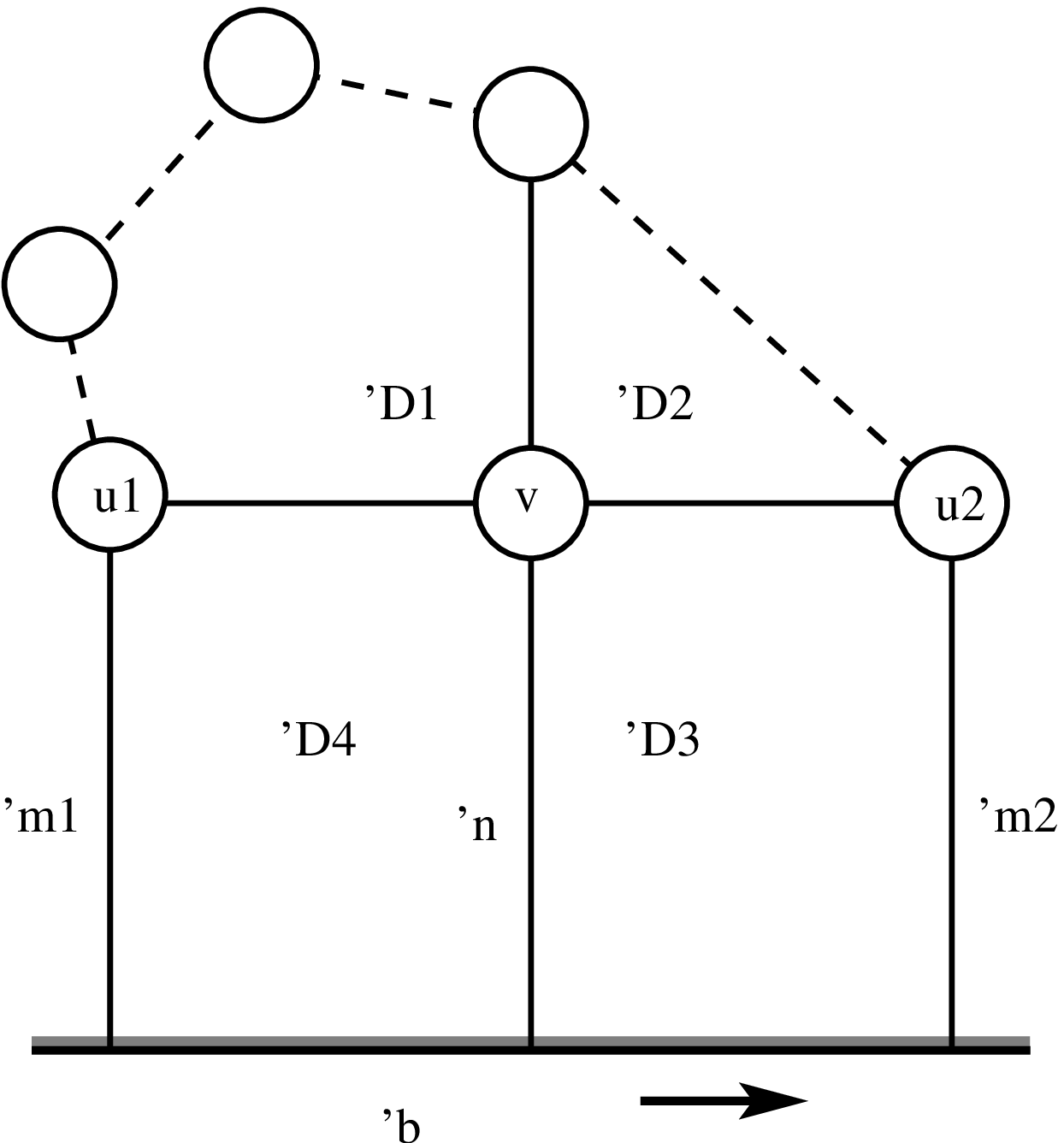} } 
}
\caption{}\label{C2}
\end{center}
\end{figure}

\begin{figure}
\psfrag{'D1}{$\D_1$}
\psfrag{'D2}{$\D_2$}
\psfrag{'D3}{$\D_3$}
\psfrag{'D4}{$\D_4$}
\psfrag{'a4}{$\al_4$}
\psfrag{'a5}{$\al_5$}
\psfrag{'b}{$\bt$}
\psfrag{'n}{$\nu$}
\psfrag{'m1}{$\mu_1$}
\psfrag{'m2}{$\mu_2$}
\psfrag{u1}{$u_1$}
\psfrag{u2}{$u_2$}
\psfrag{v}{$v$}
\begin{center}
  \mbox{
\subfigure[Configuration $C3^+$\label{C3+}]
{ \includegraphics[scale= 0.5,clip]{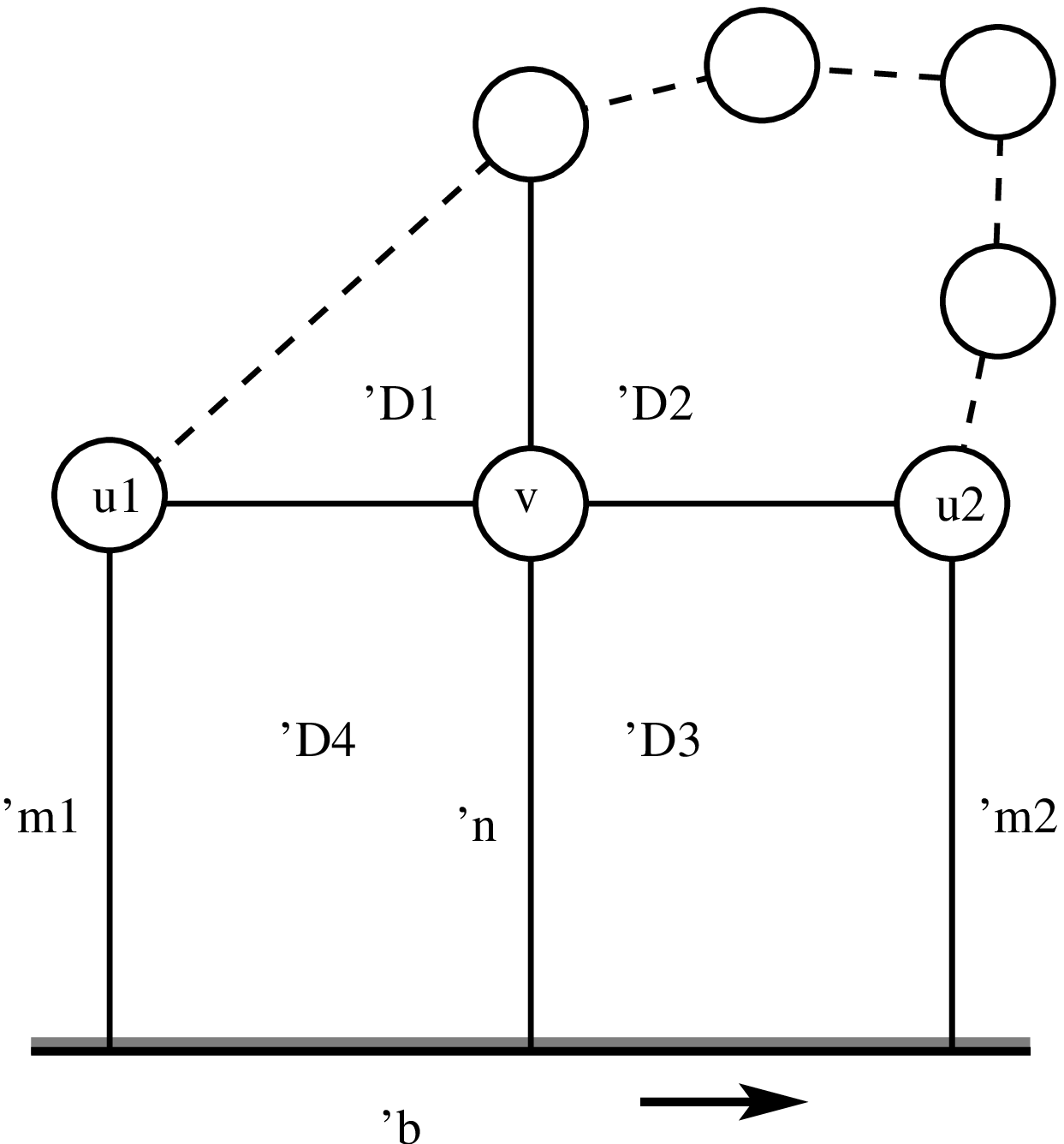} }\qquad
\subfigure[Configuration $C3^-$\label{C3-}]  
{ \includegraphics[scale=0.5,clip]{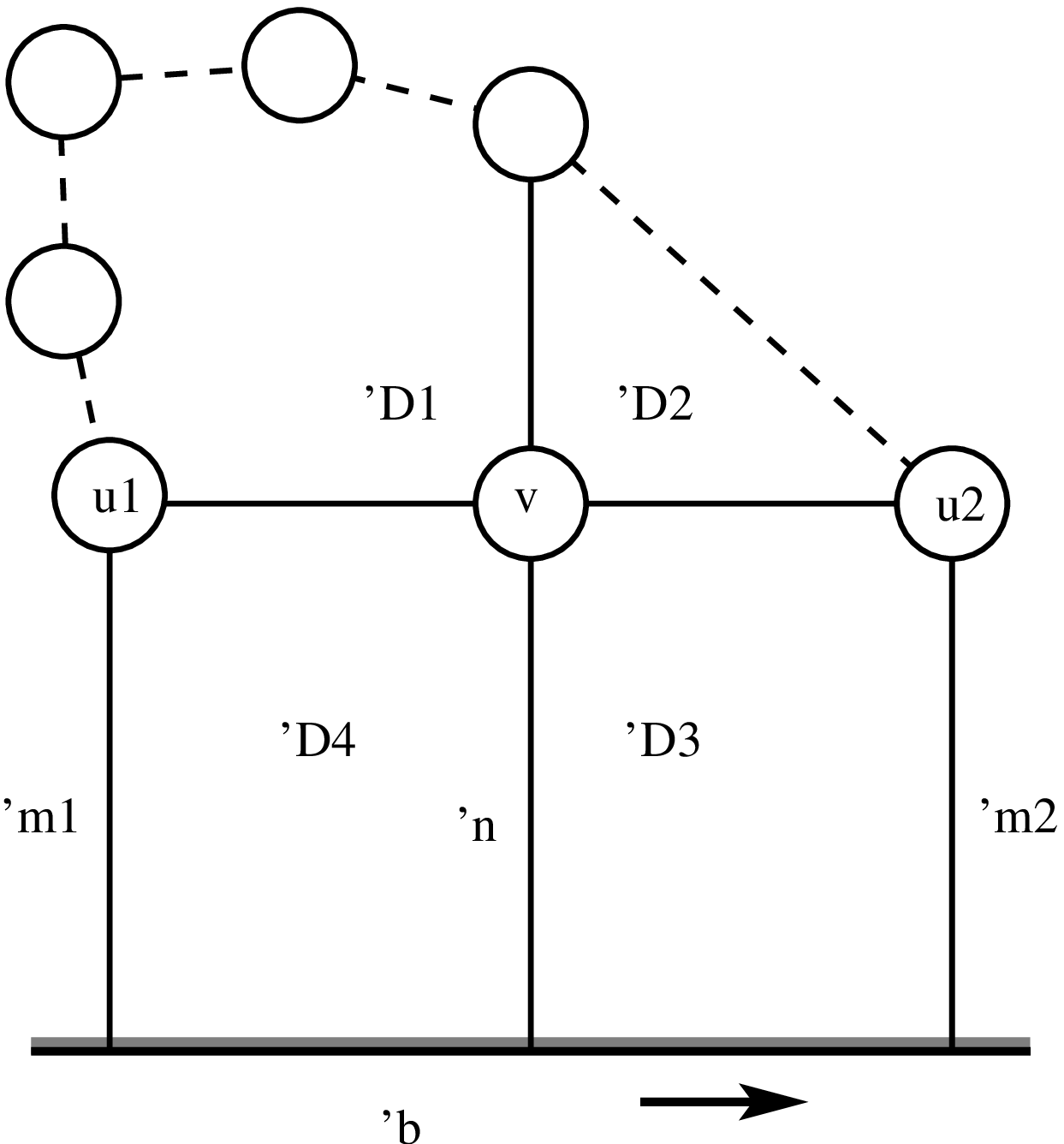} }  
}
\caption{}\label{C3}
\end{center}
\end{figure}

\begin{figure}
\psfrag{'D1}{$\D_1$}
\psfrag{'D2}{$\D_2$}
\psfrag{'D3}{$\D_3$}
\psfrag{'D4}{$\D_4$}
\psfrag{'a4}{$\al_4$}
\psfrag{'a5}{$\al_5$}
\psfrag{'b}{$\bt$}
\psfrag{'n}{$\nu$}
\psfrag{'m1}{$\mu_1$}
\psfrag{'m2}{$\mu_2$}
\psfrag{u1}{$u_1$}
\psfrag{u2}{$u_2$}
\psfrag{v}{$v$}
\begin{center}
\includegraphics[scale= 0.5,clip]{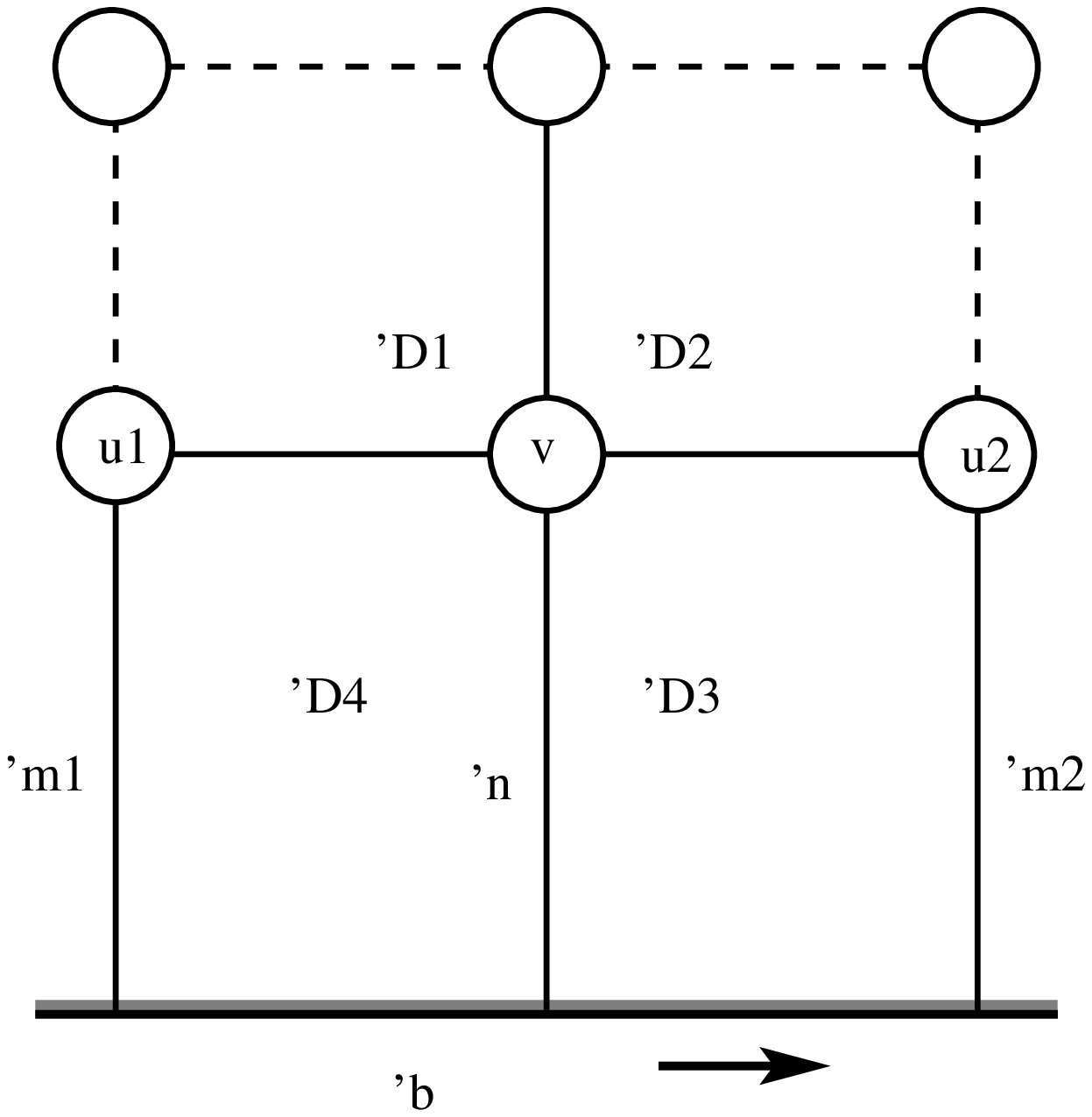}
\caption{Configuration $C4$}\label{C4}
\end{center}
\end{figure}

\begin{figure}
\psfrag{'D1}{$\D_1$}
\psfrag{'D2}{$\D_2$}
\psfrag{'D3}{$\D_3$}
\psfrag{'D4}{$\D_4$}
\psfrag{'a4}{$\al_4$}
\psfrag{'a5}{$\al_5$}
\psfrag{'b}{$\bt$}
\psfrag{'n}{$\nu$}
\psfrag{'m1}{$\mu_1$}
\psfrag{'m2}{$\mu_2$}
\psfrag{u1}{$u_1$}
\psfrag{u2}{$u_2$}
\psfrag{v}{$v$}
\begin{center}
  \mbox{
\subfigure[Configuration $CC^+$\label{CC+}]
{ \includegraphics[scale= 0.5,clip]{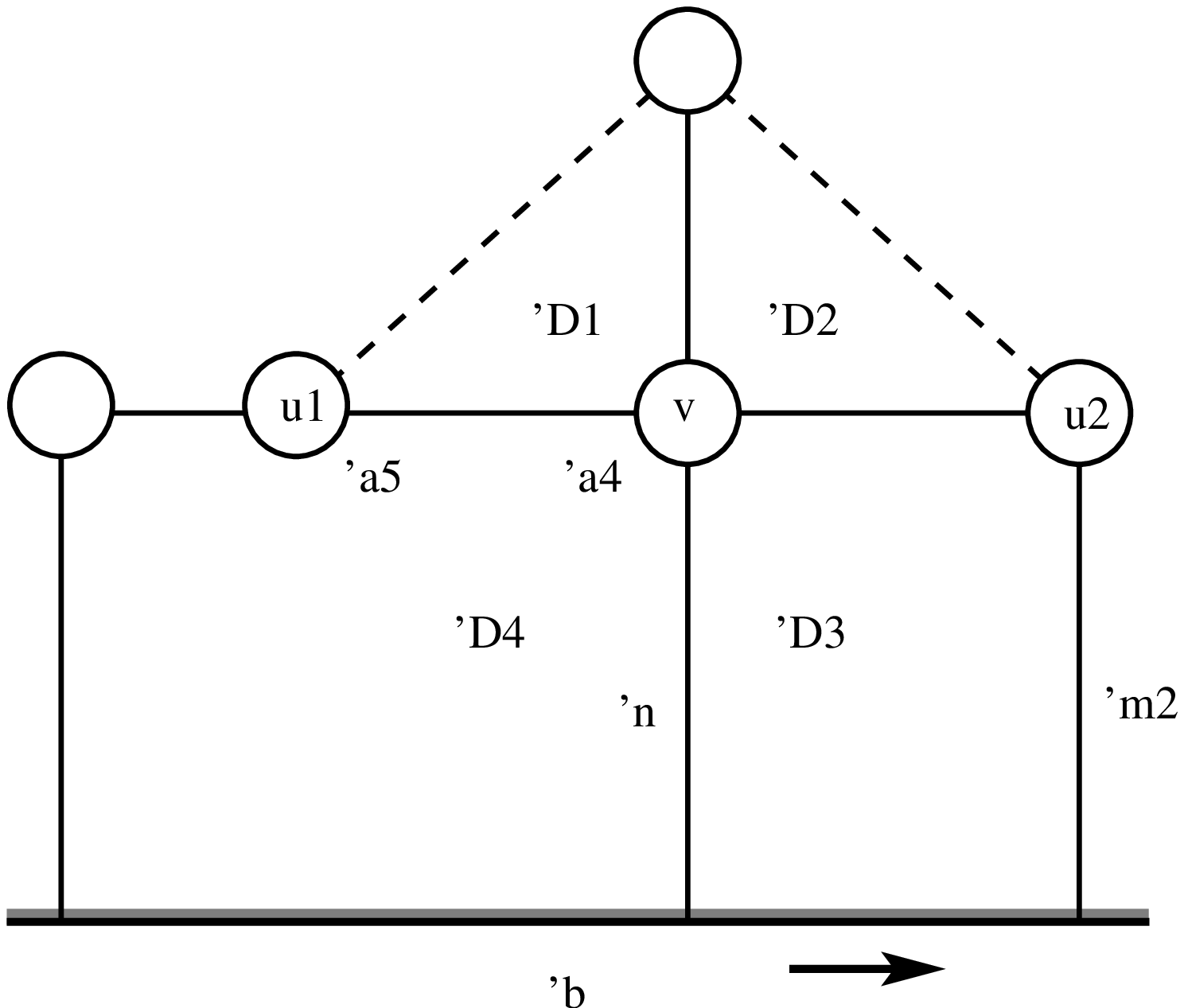} }\qquad
\subfigure[Configuration $CC^-$\label{CC-} ]  
{ \includegraphics[scale=0.5,clip]{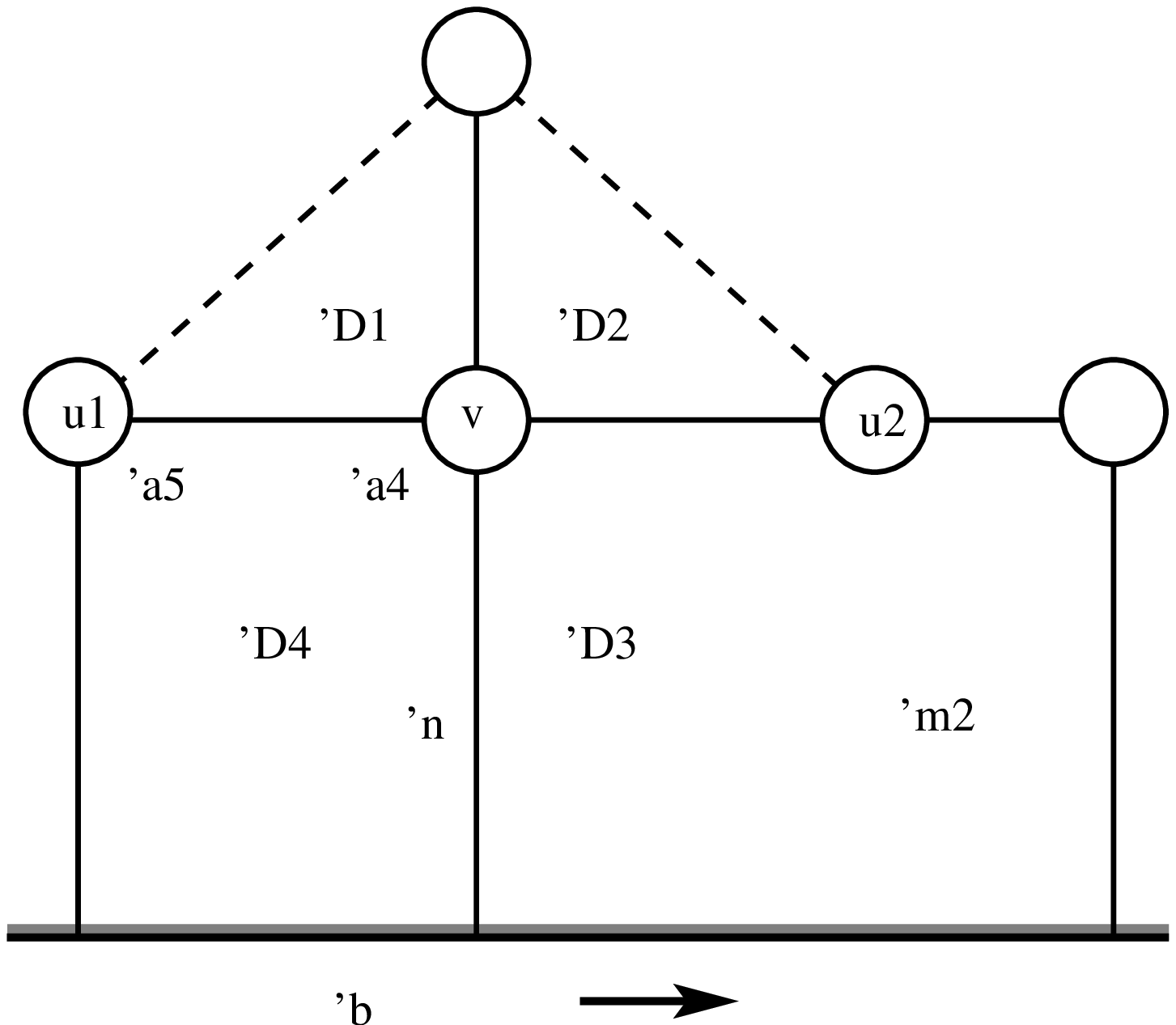} } 
}
\caption{}\label{CC}
\end{center}
\end{figure}

\begin{theorem}
\label{c-config}
Let $v$ be a vertex of $\G $ which meets $\pd \S $ in precisely one 
boundary class of arcs $\nu $.  Then either 
\be
\item\label{c-config-h} $v$ is distance at most $1$ from an $H^\Lm$--vertex or
\item\label{c-config-24} $\ka (v) \le  -\pi /24$ or 
\item\label{c-config-b} $v$ is of type $B4(a)$ (see Figure 
\ref{B4}) or 
\item\label{c-config-c} $v$ appears as the vertex $v$ in configuration $C0$,
$C1^{\pm}$, $C2^{\pm}$, $C3^{\pm}$, $C4$ or $CC^{\pm}$, 
$|\nu |> {ml / 2}+1$, and in each of 
these configurations $|\mu _i|\le l $ (see Figures \ref{C0}, \ref{C1}, \ref{C2}, \ref{C3}, \ref{C4} and \ref{CC}).
\ee
\end{theorem}

{\em Proof.~}  Assume that $v$ satisfies none of the conditions 
{\em(\ref{c-config-h})--(\ref{c-config-b})}. As $v$ is not an $H^\Lm$--vertex it follows that $\nu$ meets
$\pd\S$ in a boundary interval $[b,c]$ which we shall assume has prime label
$(w,f)$. 
Note that failure of condition {\em (\ref{c-config-24})} implies that 
$v$ is not of type $AC^{\pm}(v)$, $AB(u_i)$, $AA(u)$ or $AA(v^i)$.  
From Lemma \ref{lb4} it follows that $v$ is incident to at most 3 
interior classes of arcs.  From Corollary \ref{lc4} it follows that 
$|\nu |>{ml / 2}+1$ and from\ Lemma \ref{wv} it follows 
that $(w,f)$ is a proper exponential $H$--letter, $w$ is a subword of $r^{\pm m}$, ${ml / 2}-  l  + 3 \le  
l (w) \le  {ml / 2}$ and $|\nu | \le  {ml / 2} + l  - 1$.
\ As $m\ge 6$ this implies that there are at least $2l +1$ non--boundary 
arcs incident at $v$.  Hence $v$ is incident to exactly 3 interior classes of 
arcs and 1 boundary class of arcs $\nu .$ 

\begin{figure}
\psfrag{'D1}{$\D_1$}
\psfrag{'D2}{$\D_2$}
\psfrag{'D3}{$\D_3$}
\psfrag{'D4}{$\D_4$}
\psfrag{'a1}{$\al_1$}
\psfrag{'a2}{$\al_2$}
\psfrag{'a3}{$\al_3$}
\psfrag{'a4}{$\al_4$}
\psfrag{'b}{$\bt$}
\psfrag{'n}{$\nu$}
\psfrag{v}{$v$}
\begin{center}
{ \includegraphics[scale=0.4]{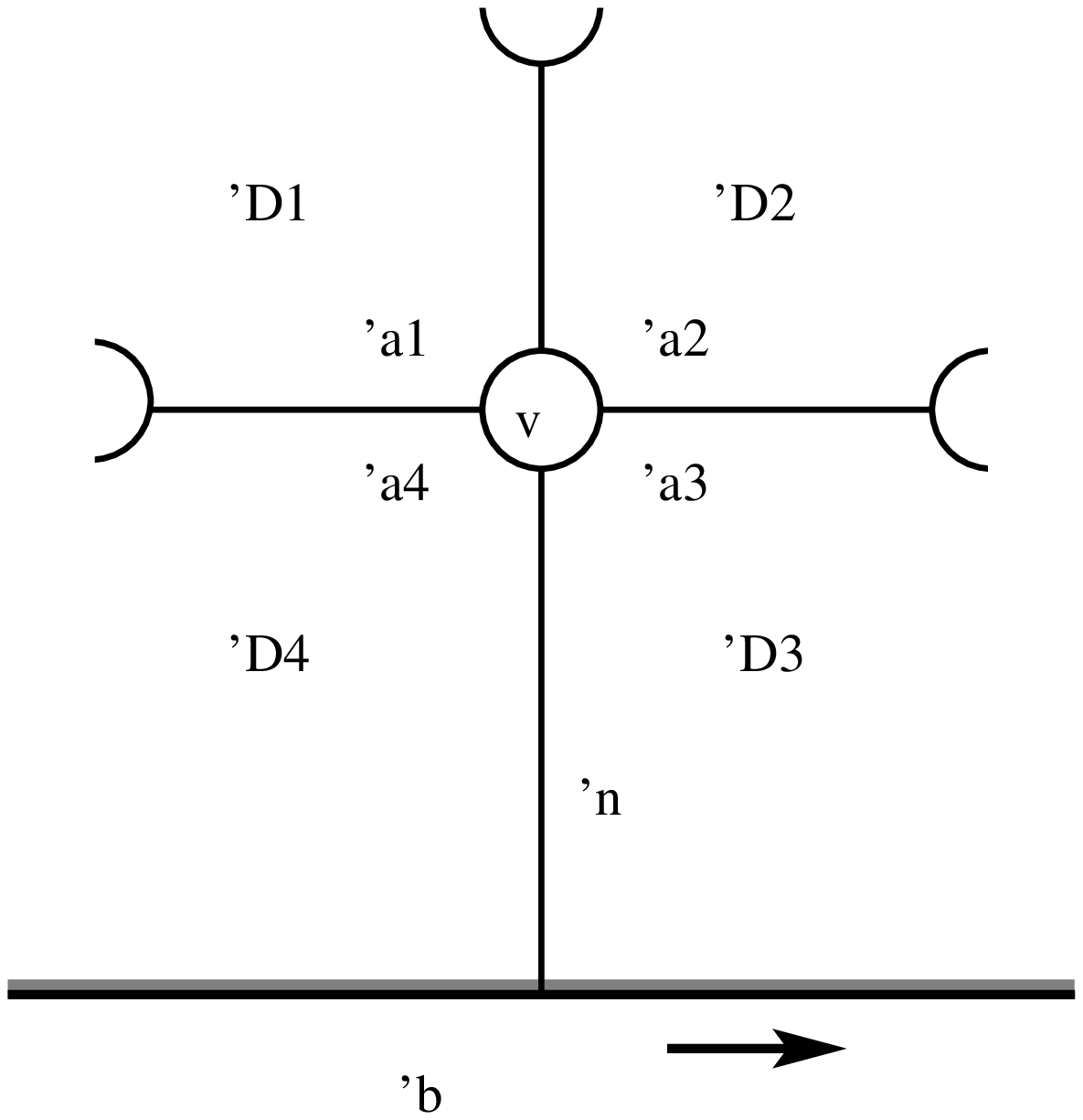} } 
\caption{\label{V3}}
\end{center}
\end{figure}

We denote the regions incident to $v$ by $\D _1$, $\D _2$, $\D _3$ 
and $\D _4$, as in Figure \ref{V3}.  
\noindent
(Note that these $\D$'s are not {\em a priori} 
distinct.)  Let $\al _i$ denote the angle on the corner of $v$ in 
$\D _i$ for all $i.$ 

We have $\al _i \ge  \pi /3$, for $i=1,2$, and $\al _i\ge \pi /2$ for 
$i=3,4$ (using Lemmas \ref{AA0}, \ref{AA10} and \ref{AA11} and 
the fact that $v$ is incident to exactly one boundary class of arcs). 
If $\al _i \ge  5\pi /7$ for $i=1$ or $2$ then $\ka (v) \le  -\pi /21.$ 
If $\al _i \ge  17\pi/18 $ for $i=3$ or $4$ then $\ka (v) \le  -\pi /9.$ 
If $\al_i=11\pi/12$, for $i=3$ or $4$, then $v$ is of type $AB(u_i)$, a
possibility that has already been excluded.

Hence we may assume $\al _i \le  2\pi /3$, for $i=1,2$ and $\al _i \le 
5\pi /6$ for $i=3,4$.  (Using Lemmas \ref{AA0}, \ref{AA10} and \ref{AA11}
again). In particular
$\x (\D _i)=1$, for $i=1,2,3,4$, which implies that all the $\D$ 's are 
distinct. 
Also $3 \le  \rho (\D _i) \le  6$ for $i=1,2$ and $2 \le  \rho (\D _i) 
\le  3$, for $i=3,4$. For $i=3$ 
or 4, write 
$$T(\D _i)=\left\{ 
\begin{array}{ll}
a & \mbox{ if } \rho (\D _i)=3 \mbox{ and }
\bt (\D _i)=1\\ 
b & \mbox{ if } \rho (\D _i)=3 \mbox{ and }
\bt (\D _i)=2 \\
c & \mbox{ if } \rho (\D _i)=2
\end{array}
\right.
.
$$ 

If $(T(\D _3), T(\D _4))=(a,a)$, $(a,b)$,  $(b,a)$ or $(b,b)$ 
then $\ka (v) 
\le -\pi /3$ so we must have $(T(\D _3),T(\D _4)) = 
(a,c),(b,c),(c,c),(c,a)$ or $(c,b).$ 
We consider the possible values of $\rho (\D _i)$, with $i=1,2$ in each of 
these cases. 

Consider first the case $(T(\D _3),T(\D _4))=(c,a)$.  If 
$\rho (\D _i) \ge  4$, for $i=1$ or $2$, then $\ka (v) \le  -\pi /6$.  
Hence $\rho (\D _i)=3$ in this case, for $i=1$ and $i=2$.  This gives 
configuration $CC^-$.  The case $(T(\D _3),T(\D _4))=(a,c)$ gives
$CC^+$.  A similar argument shows that in the cases 
$(T(\D _3),T(\D _4))=(c,b)$ and $(b,c)$ the vertex $v$ is of type $AC^{\pm}(v)$ or 
$AB(u_i)$ so these cases are already covered in (\ref{c-config-24}). 

Now consider the case $(T(\D _3),T(\D _4))=(c,c)$.  We have $3 \le  
\rho (\D _i) \le  6$, for $i=1,2$ and it is easy to check that $\ka (v)
> -\pi /24$ implies $(\rho (\D _1),\rho (\D _2)) = (3,3), (3,4), (3,5),
(3,6), (4,4), (6,3), (5,3)$ or $(4,3)$. 

These values of $\rho (\D _i)$ give rise to configurations 
$C0$, $C1^{\pm}$, $C2^{\pm}$, 
$C3^{\pm}$ and  $C4$. As $v$ is distance more than $1$ from all $H^\Lm$--vertices it 
follows, from Lemma \ref{uv} that the boundary classes $\mu_i$ in these configurations
all have width at most $l$. 
\medskip

\noindent
\noindent
\noindent
\noindent
\noindent
\noindent
\begin{defn}  \label{c-vertex}.  A configuration of type $C0$, 
$C1^{\pm}$, $C2^{\pm}$, $C3^{\pm}$, $C4$ or 
$CC^{\pm}$ is called a $C${\em --configuration} 
(see Figures \ref{C0}, \ref{C1}, \ref{C2}, \ref{C3}, \ref{C4} and \ref{CC}). 
A vertex $u$ of $\G $,
with $\ka(u)>-\pi/24$, which is distance at least $2$ from all  
$H^\Lm$--vertices
and which
appears as vertex $v$ in a $C$--configuration  is called
a $C${\em --vertex}.
\end{defn}

We shall for convenience refer to configurations $Cj^\pm$, $j=0,\ldots ,4$, taking
$Cj^+=Cj^-=Cj$ when $j=0$ or $4$.
\section{Configurations $D$}\label{config_D}
\begin{prop}  \label{D-config}  If $v$ is a $C$--vertex then 
either 
\be
\item\label{D-config-a} we can adjust the angle assignment to corners of $u_i$ and $v$
 so that $\ka (v),
\ka(u_i) 
\le  -\pi /30$, for $i=1,2$, and $\ka(\D)$ is unchanged for all regions $\D$; or 
\item\label{D-config-d} $v$ is the vertex $a$  and $u_1$ is the vertex $b$ 
in one of configurations $D0$, $D1$, $D2$, $D3$,
$D4$, $D5$, $D6$, $D7$, $D8$ or $D9$.
\ee
Such angle adjustments may be carried out over all $C$--configurations simultaneously.
\end{prop}
\noindent
\begin{figure}
\psfrag{'phi}{$\phi$}
\psfrag{'th}{$\theta$}
\psfrag{'D1}{$\D_1$}
\psfrag{'D2}{$\D_2$}
\psfrag{'D3}{$\D_3$}
\psfrag{'D4}{$\D_4$}
\psfrag{'D5}{$\D_5$}
\psfrag{'D6}{$\D_6$}
\psfrag{'D6}{$\D_7$}
\psfrag{'b}{$\bt$}
\psfrag{'n}{$\nu$}
\psfrag{'m}{$\mu$}
\psfrag{a}{$a$}
\psfrag{w}{$w$}
\psfrag{b}{$b$}
\psfrag{f}{$f$}
\psfrag{c}{$c$}
\psfrag{d}{$d$}
\begin{center}
{ \includegraphics[scale=0.5]{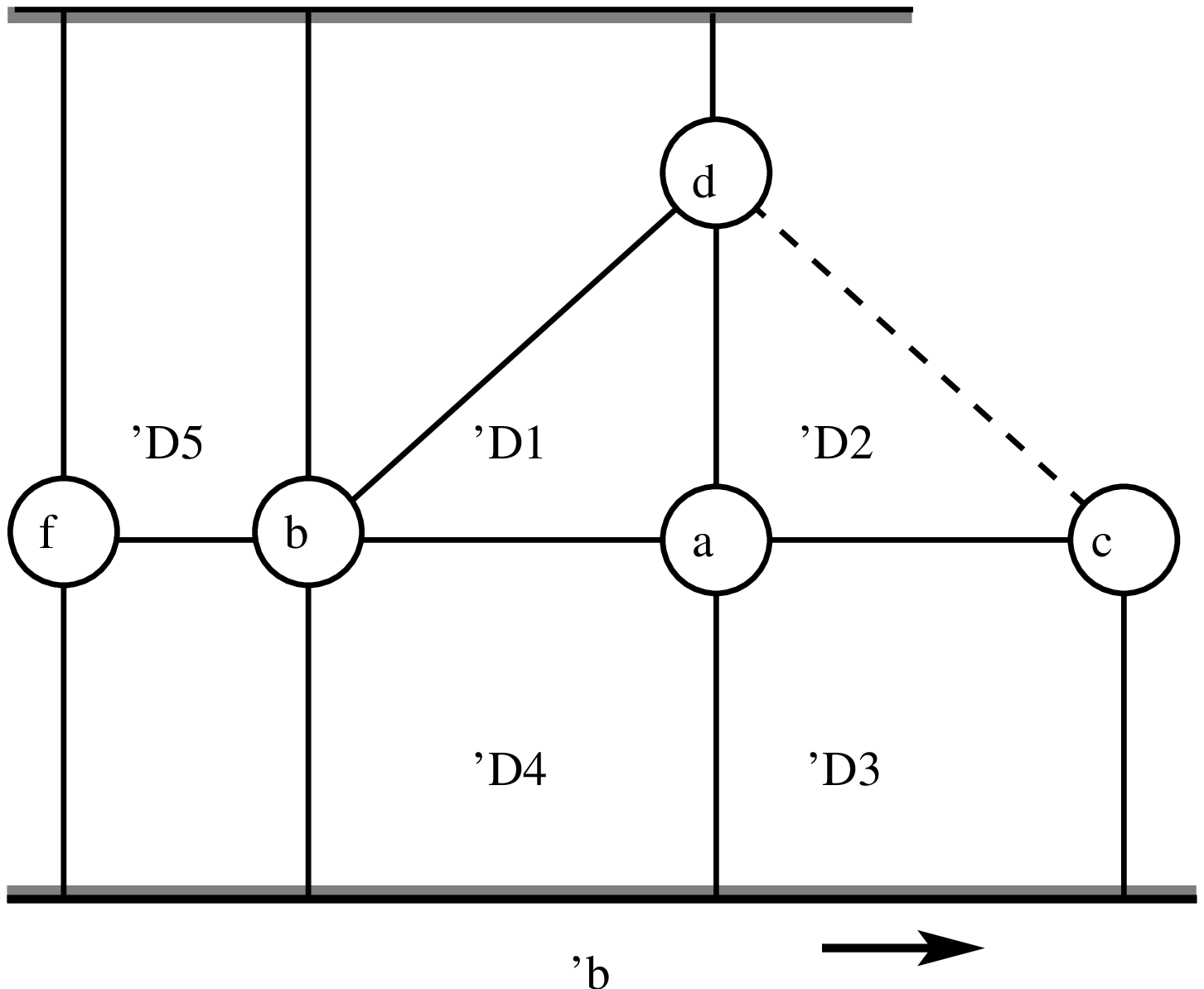} } 
\caption{Configuration $D0$\label{D0}}
\end{center}
\end{figure}
\noindent
\begin{figure}
\psfrag{'phi}{$\phi$}
\psfrag{'th}{$\theta$}
\psfrag{'D1}{$\D_1$}
\psfrag{'D2}{$\D_2$}
\psfrag{'D3}{$\D_3$}
\psfrag{'D4}{$\D_4$}
\psfrag{'D5}{$\D_5$}
\psfrag{'D6}{$\D_6$}
\psfrag{'D6}{$\D_7$}
\psfrag{'b}{$\bt$}
\psfrag{'n}{$\nu$}
\psfrag{'m}{$\mu$}
\psfrag{a}{$a$}
\psfrag{b}{$b$}
\psfrag{f}{$f$}
\psfrag{c}{$c$}
\psfrag{d}{$d$}
\begin{center}
{ \includegraphics[scale=0.5]{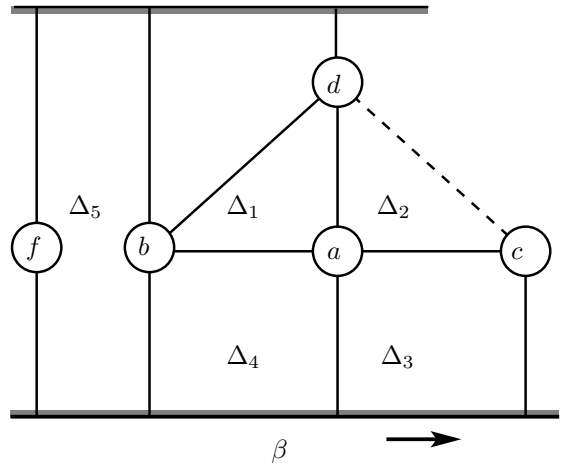} } 
\caption{Configuration $D1$\label{D1}}
\end{center}
\end{figure}
\noindent
\begin{figure}
\psfrag{'phi}{$\phi$}
\psfrag{'th}{$\theta$}
\psfrag{'D1}{$\D_1$}
\psfrag{'D2}{$\D_2$}
\psfrag{'D3}{$\D_3$}
\psfrag{'D4}{$\D_4$}
\psfrag{'D5}{$\D_5$}
\psfrag{'D6}{$\D_6$}
\psfrag{'D6}{$\D_7$}
\psfrag{'b}{$\bt$}
\psfrag{'n}{$\nu$}
\psfrag{'m}{$\mu$}
\psfrag{a}{$a$}
\psfrag{f}{$f$}
\psfrag{c}{$c$}
\psfrag{d}{$d$}
\psfrag{b}{$b$}
\begin{center}
{ \includegraphics[scale=0.5]{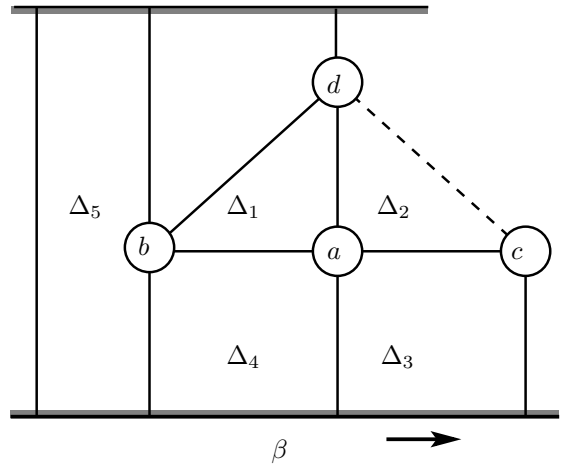} } 
\caption{Configuration $D2$\label{D2}}
\end{center}
\end{figure}
\noindent
\begin{figure}
\psfrag{'phi}{$\phi$}
\psfrag{'th}{$\theta$}
\psfrag{'D1}{$\D_1$}
\psfrag{'D2}{$\D_2$}
\psfrag{'D3}{$\D_3$}
\psfrag{'D4}{$\D_4$}
\psfrag{'D5}{$\D_5$}
\psfrag{'D6}{$\D_6$}
\psfrag{'b}{$\bt$}
\psfrag{a}{$a$}
\psfrag{b}{$b$}
\psfrag{c}{$c$}
\psfrag{d}{$d$}
\psfrag{f}{$f$}
\begin{center}
{ \includegraphics[scale=0.5]{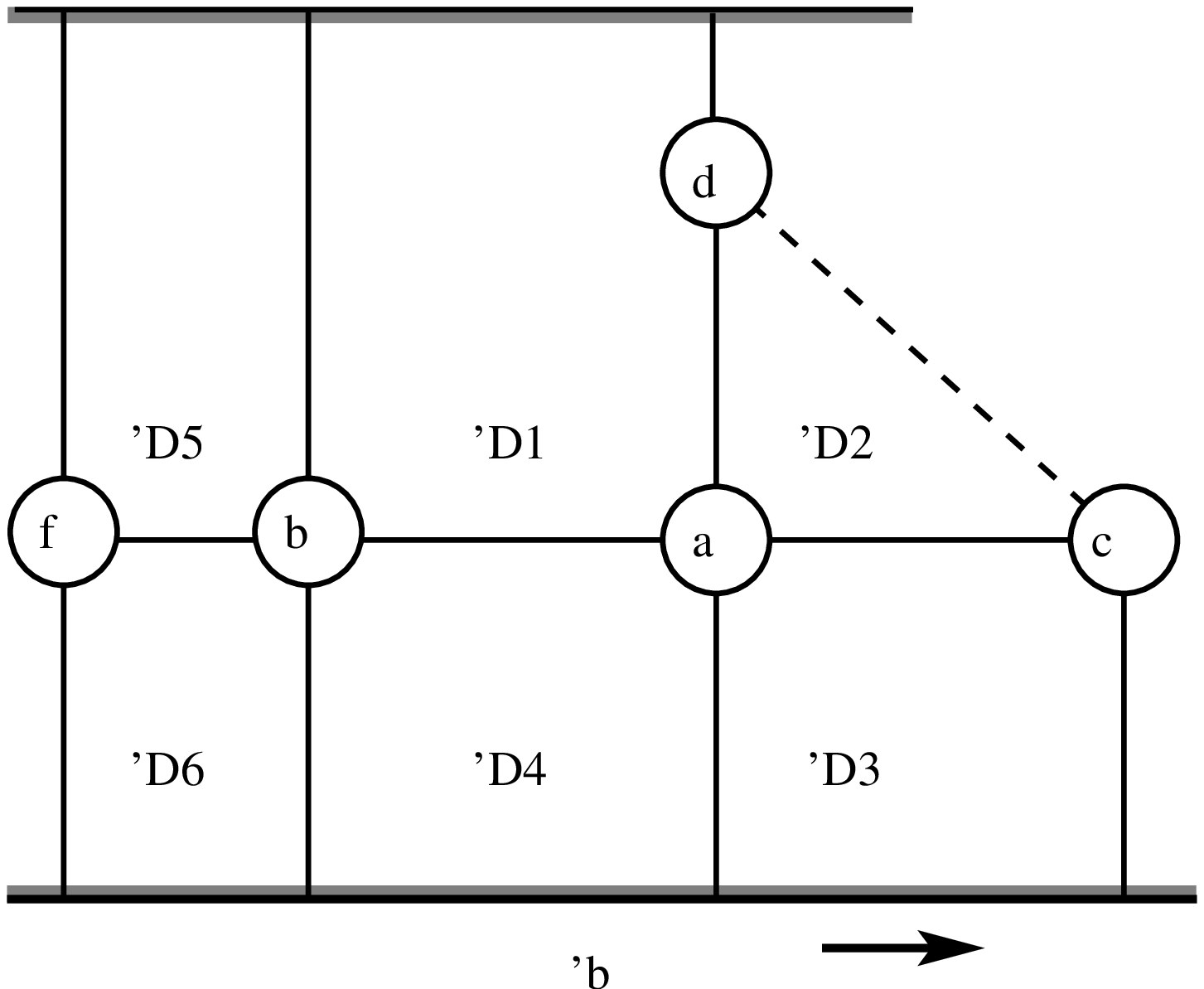} } 
\caption{Configuration $D3$\label{D3}}
\end{center}
\end{figure}
\noindent
\begin{figure}
\psfrag{'D1}{$\D_1$}
\psfrag{'D2}{$\D_2$}
\psfrag{'D3}{$\D_3$}
\psfrag{'D4}{$\D_4$}
\psfrag{'b}{$\bt$}
\psfrag{'n}{$\nu$}
\psfrag{e}{$e$}
\psfrag{f}{$f$}
\psfrag{a}{$a$}
\psfrag{b}{$b$}
\psfrag{c}{$c$}
\psfrag{d}{$d$}
\begin{center}
{ \includegraphics[scale=0.5]{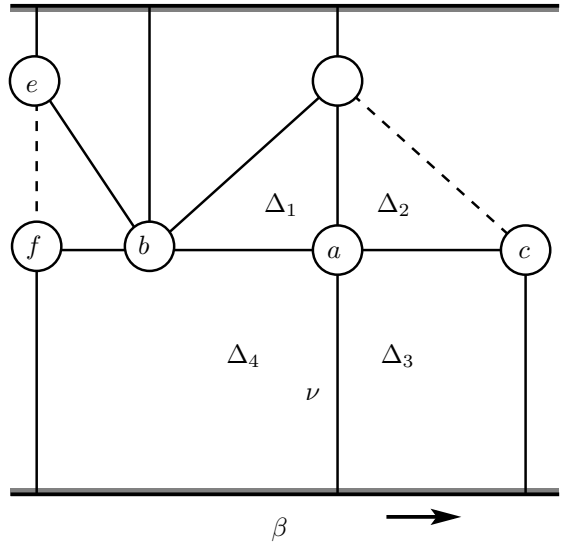} } 
\caption{Configuration $D4$\label{D4}}
\end{center}
\end{figure}
~\\
\noindent
\begin{figure}
\psfrag{'D1}{$\D_1$}
\psfrag{'D2}{$\D_2$}
\psfrag{'D3}{$\D_3$}
\psfrag{'D4}{$\D_4$}
\psfrag{'b}{$\bt$}
\psfrag{'n}{$\nu$}
\psfrag{e}{$e$}
\psfrag{f}{$f$}
\psfrag{a}{$a$}
\psfrag{b}{$b$}
\psfrag{c}{$c$}
\psfrag{d}{$d$}
\begin{center}
{ \includegraphics[scale=0.5]{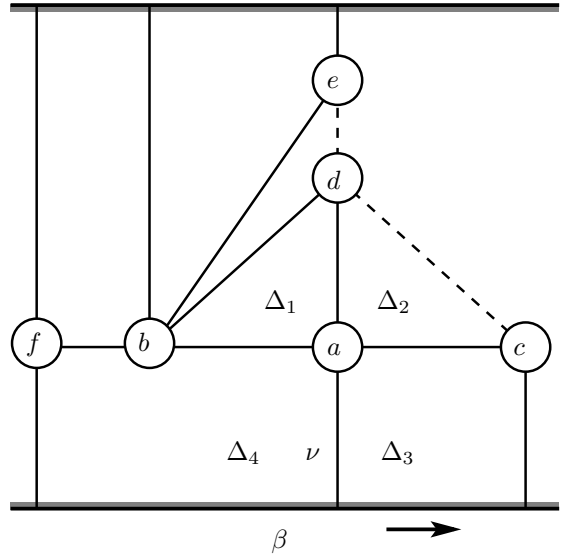} } 
\caption{Configuration $D5$\label{D5}}
\end{center}
\end{figure}
\noindent
\begin{figure}
\psfrag{'phi}{$\phi$}
\psfrag{'th}{$\theta$}
\psfrag{'D1}{$\D_1$}
\psfrag{'D2}{$\D_2$}
\psfrag{'D3}{$\D_3$}
\psfrag{'D4}{$\D_4$}
\psfrag{'D5}{$\D_5$}
\psfrag{'D6}{$\D_6$}
\psfrag{'D7}{$\D_7$}
\psfrag{'b}{$\bt$}
\psfrag{'n}{$\nu$}
\psfrag{'m}{$\mu$}
\psfrag{e}{$e$}
\psfrag{f}{$f$}
\psfrag{a}{$a$}
\psfrag{b}{$b$}
\psfrag{c}{$c$}
\psfrag{d}{$d$}
\psfrag{w}{$w$}
\begin{center}
{ \includegraphics[scale=0.5]{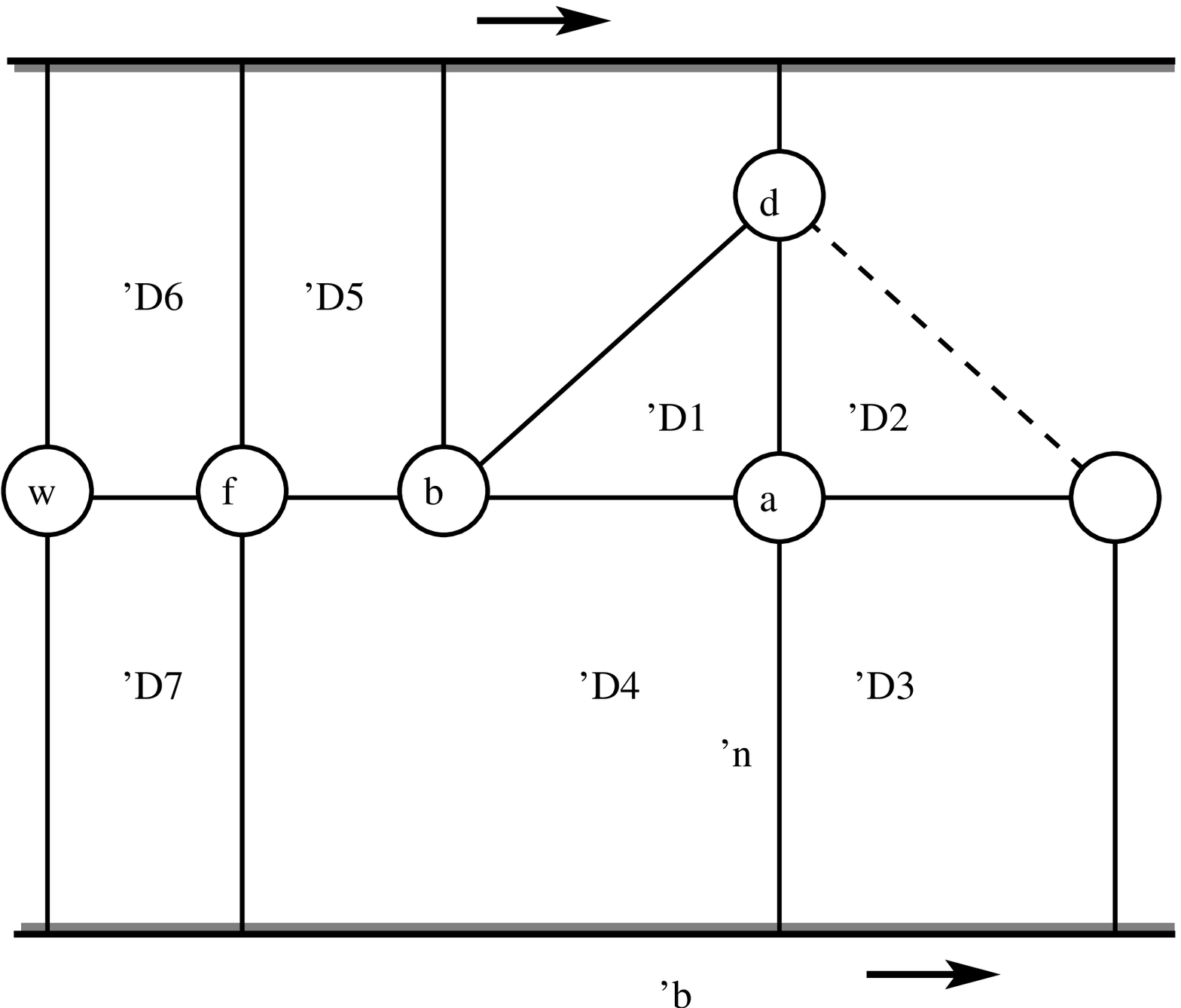} } 
\caption{Configuration $D6$\label{D6}}
\end{center}
\end{figure}
\begin{figure}
\psfrag{'phi}{$\phi$}
\psfrag{'th}{$\theta$}
\psfrag{'D1}{$\D_1$}
\psfrag{'D2}{$\D_2$}
\psfrag{'D3}{$\D_3$}
\psfrag{'D4}{$\D_4$}
\psfrag{'D5}{$\D_5$}
\psfrag{'D6}{$\D_6$}
\psfrag{'D6}{$\D_7$}
\psfrag{'b}{$\bt$}
\psfrag{'n}{$\nu$}
\psfrag{'m}{$\mu$}
\psfrag{e}{$e$}
\psfrag{f}{$f$}
\psfrag{a}{$a$}
\psfrag{b}{$b$}
\psfrag{c}{$c$}
\psfrag{d}{$d$}
\psfrag{w}{$w$}
\begin{center}
{ \includegraphics[scale=0.5]{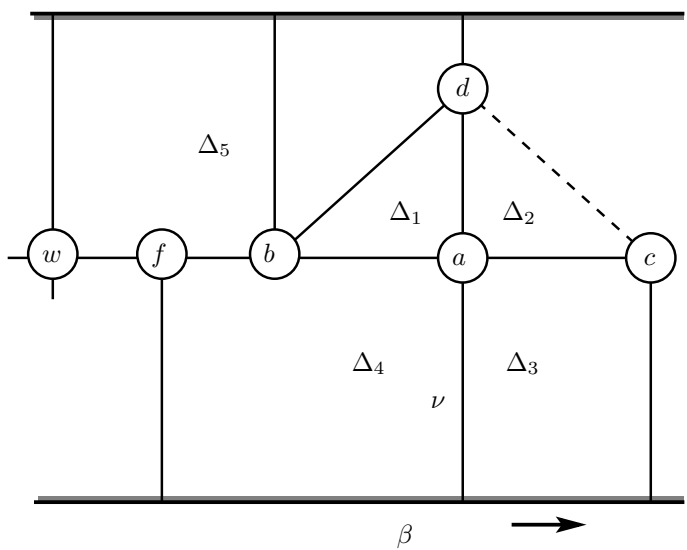} } 
\caption{Configuration $D7$\label{D7}}
\end{center}
\end{figure}
\begin{figure}
\psfrag{'phi}{$\phi$}
\psfrag{'th}{$\theta$}
\psfrag{'D1}{$\D_1$}
\psfrag{'D2}{$\D_2$}
\psfrag{'D3}{$\D_3$}
\psfrag{'D4}{$\D_4$}
\psfrag{'D5}{$\D_5$}
\psfrag{'D6}{$\D_6$}
\psfrag{'D6}{$\D_7$}
\psfrag{'b}{$\bt$}
\psfrag{'n}{$\nu$}
\psfrag{'m}{$\mu$}
\psfrag{e}{$e$}
\psfrag{f}{$f$}
\psfrag{a}{$a$}
\psfrag{b}{$b$}
\psfrag{c}{$c$}
\psfrag{d}{$d$}
\begin{center}
{ \includegraphics[scale=0.5]{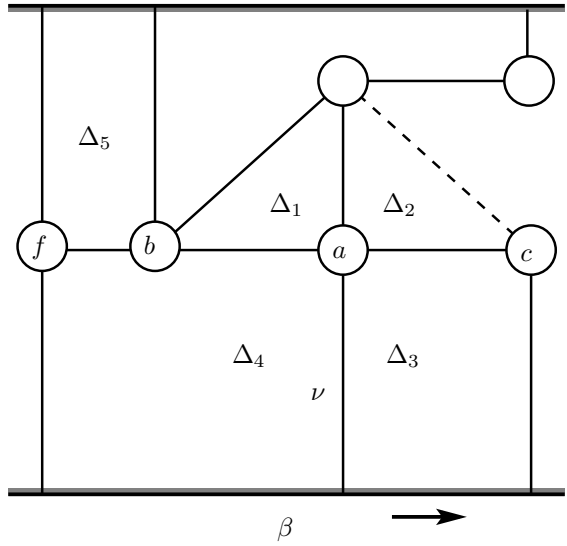} } 
\caption{Configuration $D8$\label{D8}}
\end{center}
\end{figure}
\begin{figure}
\psfrag{'phi}{$\phi$}
\psfrag{'th}{$\theta$}
\psfrag{'D1}{$\D_1$}
\psfrag{'D2}{$\D_2$}
\psfrag{'D3}{$\D_3$}
\psfrag{'D4}{$\D_4$}
\psfrag{'D5}{$\D_5$}
\psfrag{'D6}{$\D_6$}
\psfrag{'D6}{$\D_7$}
\psfrag{'b}{$\bt$}
\psfrag{'n}{$\nu$}
\psfrag{'m}{$\mu$}
\psfrag{a}{$a$}
\psfrag{e}{$e$}
\psfrag{f}{$f$}
\psfrag{b}{$b$}
\psfrag{c}{$c$}
\psfrag{d}{$d$}
\begin{center}
{ \includegraphics[scale=0.5]{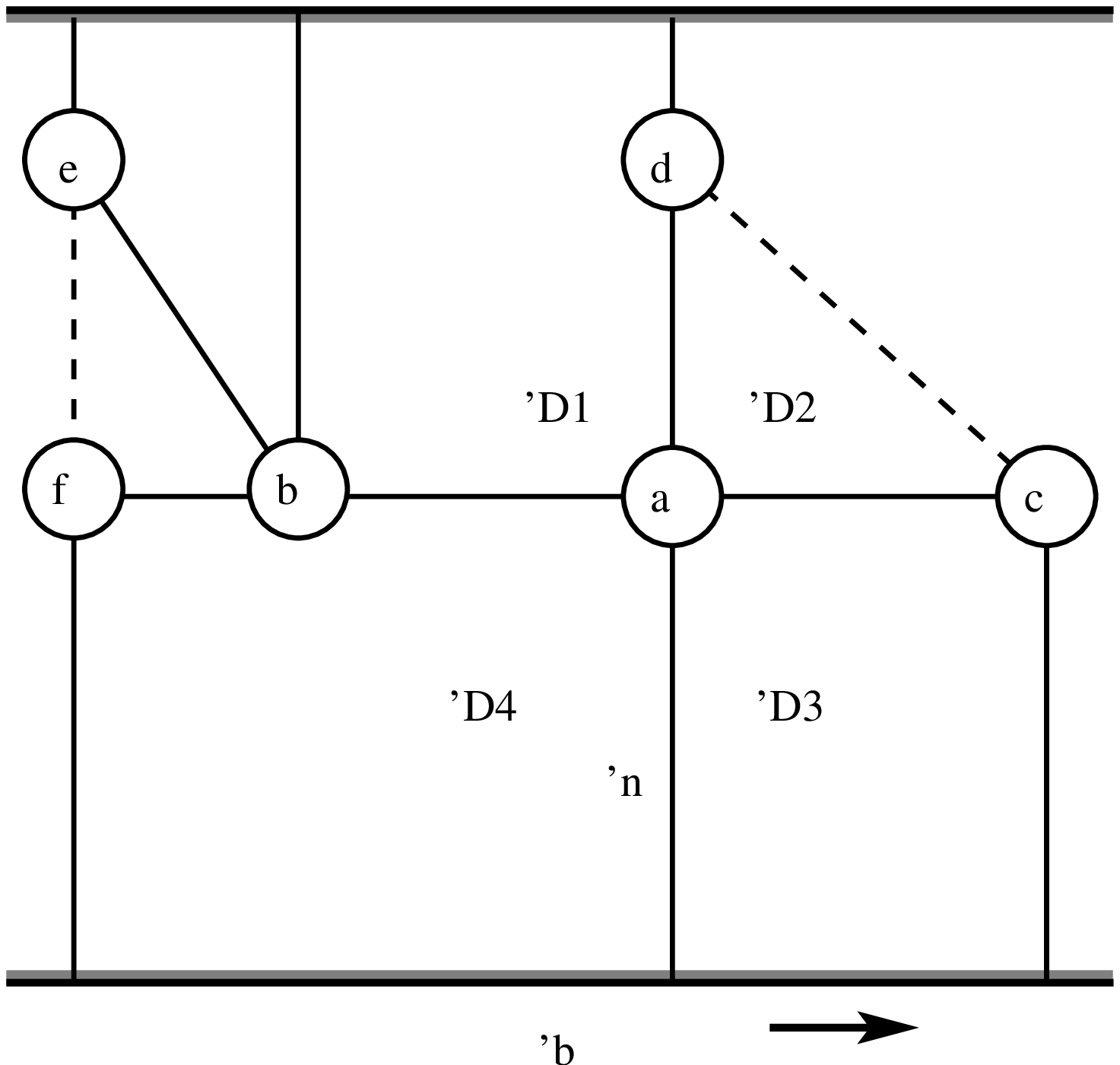} } 
\caption{Configuration $D9$\label{D9}}
\end{center}
\end{figure}

\noindent
{\em Proof.~}
We shall define angle adjustments to satisfy (\ref{D-config-a}) and (\ref{D-config-d}) for all
$C$--configurations of $\G$ such that,  
considering each connected boundary component $\bt$ of $\pd \S$ in turn, these angle
adjustments move angle from corner to corner of regions incident to $\bt$,
only in the direction of orientation of $\bt$. Furthermore these angle 
adjustments leave $\ka(\D)$ fixed, for all regions $\D$ of $\G$.

Given a $C$--configuration with $C$--vertex $v$ set $u=u_1$ and $\mu =\mu _1$ 
(if we're not in case $CC^+$). 
Denote the class of arcs joining $u$ and $v$ by $\iota$.  
Denote the angle on the 
corners of $v$ and $u$ in $\D _4$ by $\al_4$ and $\al _5$, respectively. 
Orient $\pd u$ in the direction from 
$\mu$ to $\iota$ inside $\D_4$. 
Note that $v$ is not of type $AA(v^i)$ these
vertices have curvature $\ka(v^i)\le -\pi/9$. By definition $v$ is not of type
$AA(x^i)$, $AC^\pm(v)$, $AC^+(x)$, $AC^-(y)$, $AB(u_i)$ or $AB(v)$. Therefore
no corners of $v$ have had angles altered by previous angle adjustments. 

To ensure that we can make angle adjustments simultaneously over all $C$--configurations
we need to consider whether or not $u$ and $v$ may appear in more than one $C$--configuration.
As $v$ is incident to only one boundary class of arcs it cannot occur as the $C$--vertex
of any other $C$--configuration. As $|\mu|\le l$ and $|\nu|\ge ml/2+1$, $u$ cannot be the 
$C$--vertex of a $C$--configuration in which $v$ occurs as the vertex of type $u_1$ or $u_2$.
Suppose that $u$ occurs as the $C$--vertex in some $C$--configuration. Let $C_0$ be the original
$C$--configuration, in which $u$ is of type $u_1$, and let $C_1$ be the configuration in which
$u$ occurs as the $C$--vertex. Since $u$ is incident to a boundary class in $C_1$ it follows
that $C_0$ is of type $CC^+$, so $u$ being the vertex of type $CC^+(u_1)$ in $C_0$ is not incident
to a further boundary class of arcs. If $C_1$ is a configuration of type $Cj$ then, since  $\rho (\D_1)=3$, it 
follows that $C_1$ is of type $C0$, with $u$ of type $C0(v)$. This gives configuration $D^\prime 6$ of 
Figure \ref{D'6}. Otherwise $C_1$ is of type $CC^\pm$ and we obtain configuration $D8$ of Figure \ref{D8}
or $D9$ of Figure \ref{D9}, with $u$ of type $Dj(b)$. In all three cases it can be easily verified that
$u$ is of type $u_1$ in $C_0$ alone.    

In general $u$ may occur as the vertex of type $u_1$ in more than one $C$--configuration. 
Suppose then that $u$ is the vertex of type $u_1$ in $p$ distinct 
$C$--configurations. 
We show that {\em (\ref{D-config-d})}
holds only if $p=1$ and  that  if $p=1$, {\em (\ref{D-config-d})} does not hold and $u$ and $v$ 
are not of type $D^\prime 6(u)$ and $D^\prime 6(v)$, respectively, then 
we can decrease the angle $\al_5$ and increase
$\al_4$ to obtain $\ka(u),\ka(v)\le -\pi/12$. 
If $p>1$ we show that
we can decrease the angle $\al_5$ and increase
$\al_4$ to obtain $\ka(u),\ka(v)\le -\pi/15$ over
all $p$ configurations, in which $u$ occurs as the vertex of type $u_1$, simultaneously.
Finally suppose that $u$ and $v$ 
are of type $D^\prime 6(u)$ and $D^\prime 6(v)$, respectively, and let $\bt^\prime$ 
be the boundary component to which $u$ is incident in Figure \ref{D'6}. 
If $\bt^\prime$ is oriented from left to 
right in Figure \ref{D'6} 
we show that we can adjust angles by increasing  $\al_4$ and $\al_7$ and decreasing $\al_5$ and  $\al_6$ so
that $\ka(u),\ka(v),\ka(x)\le -\pi/30$. 
Similarly if $\bt^\prime$ is oriented from right to left 
in Figure \ref{D'6} we show that we can adjust angles  by increasing $\al_4$ and  $\al_8$ and 
decreasing  $\al_5$ and $\al_9$ so
that $\ka(u),\ka(v),\ka(y)\le -\pi/30$. These adjustments all leave curvature of regions unaltered. 
This is sufficient, given the remark at the beginning of the proof,  to prove the theorem.
We divide the proof into Case A, where $p=1$, and Case B, where $p>1$.

\subsection*{{Case A: $p=1$.}}\label{caseA}~\\

In this case $u$ occurs as the vertex of type $u_1$
in no $C$--configuration except that in which $v$ occurs. Under this
assumption we have two further cases to consider: Case A1 in which $u$ is adjacent to 
no boundary classes (in configuration $CC^+$) 
except possibly $\mu $\ (in configurations
$Cj^{\pm},j=0,\ldots ,4$ or $CC^{-}$) 
and Case A2 in which $u$ is incident to boundary classes other than $\mu$. 
\subsection*{\normalsize Case A1.}\label{caseA1}
As $G$ is not of type $E(2,*,m)$ no interior class of arcs contains more than 
$l-2$ arcs and $\mu$ contains no more than $l$ arcs, from Theorem \ref{c-config}.
It follows that in Case A1, configuration $CC^+$, $u$ is adjacent to at least 
$7$ interior classes of arc and in Case A1, configuration $Cj^{\pm}$ or $CC^-$, $u$ is adjacent 
to at least $6$ interior classes of arc.

In configuration $Cj^\pm$, with $0\le j\le 4$, or $CC^-$ with $u$  
incident to 6 or more interior classes of arc we have $\s (u)\ge \pi +5\pi /3$, so 
we can adjust $\al _5$, by $-\pi/2$, to 0 and $\al _4$, by $\pi/2$, to $\pi $ giving 
$\ka (u),  \ka (v)\le -\pi /6$.  In 
configuration $CC^+$ with $u$ is 
incident to $7$ or more classes of arc we have $\s (u)\ge 7\pi /3$,
$\al _5=\pi /3$ and 
$\al _4=5\pi /6$ so we may adjust to $\al _5=\pi /6$ and $\al _4=\pi $, to give $\ka(u),\ka(v)\le -\pi/12$. 

Hence in Case A1 it is always possible to adjust angles 
as required. 
\subsection*{\normalsize Case A2.}\label{caseA2}
\begin{figure}
\psfrag{'phi}{$\phi$}
\psfrag{'th}{$\theta$}
\psfrag{'D1}{$\D_1$}
\psfrag{'D4}{$\D_4$}
\psfrag{a4}{$\al_4$}
\psfrag{a5}{$\al_5$}
\psfrag{mu}{$\mu$}
\psfrag{v}{$v$}
\psfrag{xi}{$\xi$}
\psfrag{u}{$u$}
\psfrag{'b}{$\bt$}
\psfrag{pi/3}{$\pi/3$}
\psfrag{i}{$\io$}
\begin{center}
\includegraphics[scale=0.5]{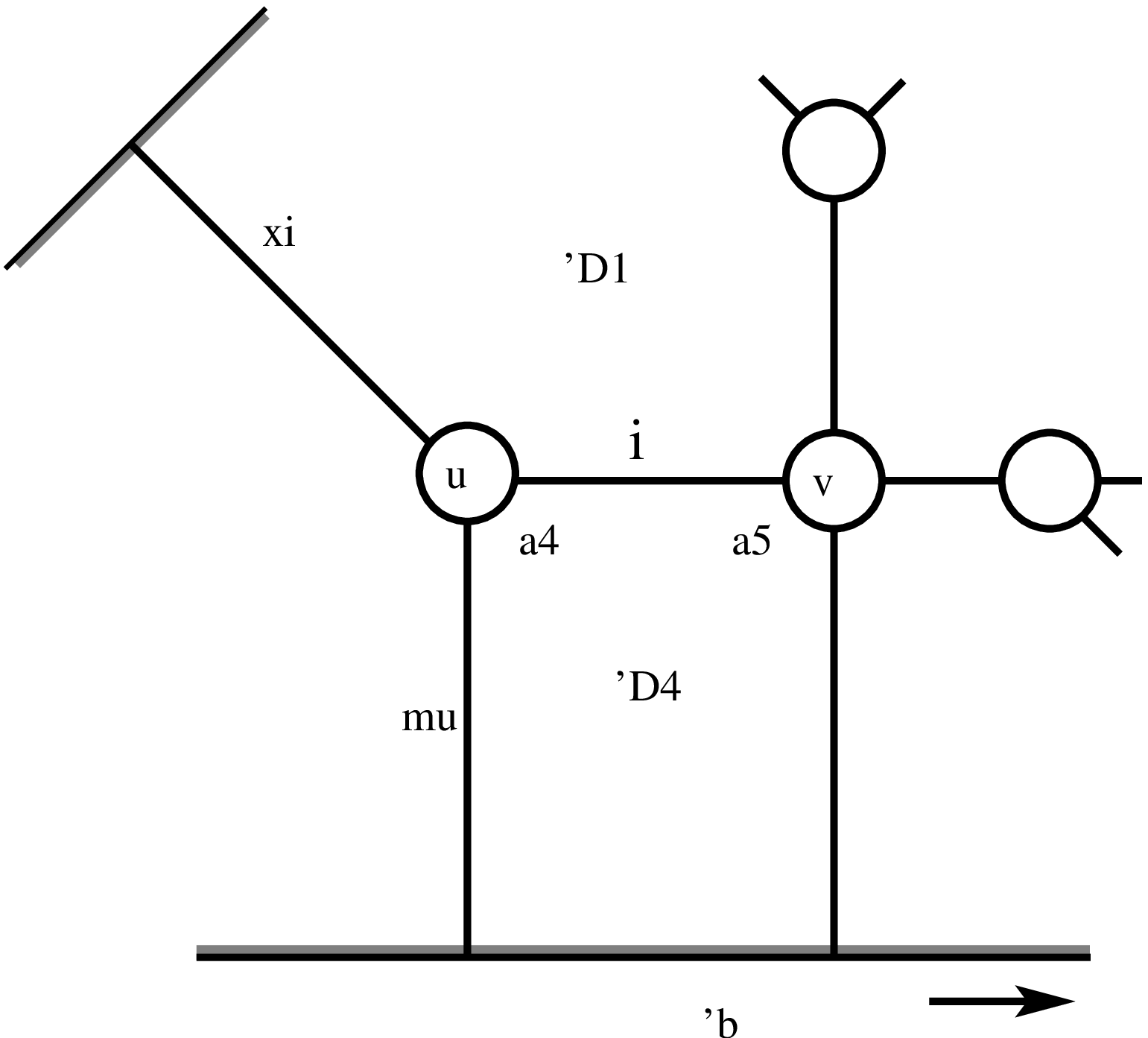}
\caption{}\label{A_C_vertex}
\end{center}
\end{figure}

Now suppose that $u$ is incident to some boundary class of arcs $\xi $
not meeting $\pd \D _4$.  Suppose first we have one of configurations 
$C1^{\pm}$, $C2^{\pm}$, $C3^{\pm}$, $C4$ or $CC^-$ (see Figure \ref{A_C_vertex}).  
In each case $\ka (v)\le \pi /6$.  Let
$\theta $ be the sum of angles on $u$ between $\iota$ and $\xi $ (that is,
 taken in the 
direction of $\D _1$).  If (in the same direction) there is 1 or more 
class of arcs incident at $u$ between $\iota$ 
and $\xi $ then it follows from 
Lemmas \ref{AA0}, \ref{AA10} and 
\ref{AA11} that $\theta \ge 5\pi /6$.  If there are no arcs 
between $\iota$ and $\xi $ then since $\rho (\D _1)\ge 3$ and, as is easily
verified, $\D_1$ is not a region of type $AA(\D_4^i)$, it follows that 
$\theta \ge 5\pi /6$, again.  Hence 
$\s (u)\ge \pi +5\pi /6+\pi /2=7\pi /3$.  We may therefore increase 
$\al _4$ by $\pi /4$ and decrease $\al _5$ by $\pi /4$ to obtain 
$\ka (u)$, $\ka (v) \le  -\pi /12.$ 

Hence we need only consider configurations $C0$ and $CC^+$. 

\begin{figure}
\psfrag{xi}{$\xi$}
\psfrag{'D1}{$\D_1$}
\psfrag{u}{$u$}
\psfrag{v}{$v$}
\psfrag{'b}{$\bt$}
\psfrag{pi/3}{$\pi/3$}
\psfrag{i}{$\io$}
\begin{center}
{ \includegraphics[scale=0.4]{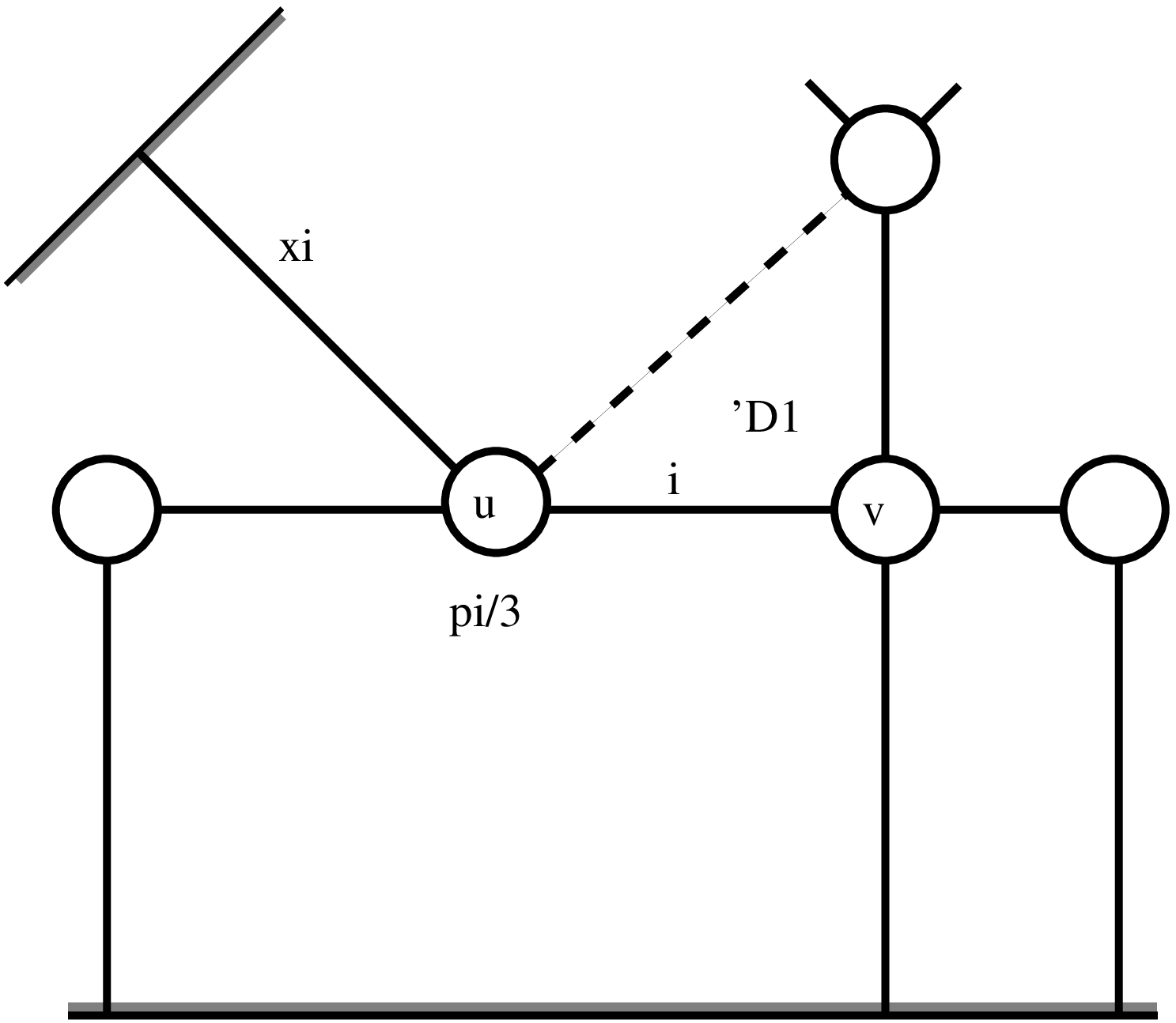} } 
\caption{\label{A_A_vertex}}
\end{center}
\end{figure}

Assume we have configuration $CC^+$, so $\ka (v)=0$ (see Figure \ref{A_A_vertex}).
Suppose first $\D_1$ is an 
interior region, so triangular with all angles $\pi /3$.  In this case 
if any 
region meeting $u$ has Euler characteristic less than $1$ then 
$\s (u)\ge 2\pi /3+\pi +\pi /2=13\pi /6$, so we may decrease $\al _5$ 
to $\pi /4$ and increase $\al _4$ to $11\pi /12$ to obtain 
$\ka (u),\ka (v)\le -\pi /12$.  Hence we may assume all regions 
$\D $ incident to $u$ have $\x (\D )=1$.  Now under the current assumptions 
$u$ is incident to $3$
 or 
more interior classes of arc and one or more boundary classes of arc.  We 
consider four subcases. Subcase A2.1 where $u$ is 
incident to more than one boundary 
class; Subcase A2.2 where $u$ is incident to one boundary class and more than four 
interior classes; Subcase A2.3 where $u$ is incident to 1 boundary class and four 
interior classes and Subcase A2.4 where $u$ is incident to 1 boundary class and 
three interior classes.

\bd
\item[Subcase A2.1] $u$ is incident to more than one boundary class of arcs. 

In this case $\s (u)\ge 2\pi /3+\pi +\pi /2+\pi /2$ so we may adjust 
angles $\al _4$ and $\al _5$ to give $\ka (u), 
\ka (v) \le -\pi /3.$

\item[Subcase A2.2] $u$ is incident to one boundary class and more than four interior 
classes of arcs. 

In this case $\s (u)\ge 4\pi /3+\pi /2+\pi /2$ so we may adjust angles 
$\al _4$ and $\al _5$ to obtain $\ka (u),\ka (v)\le -\pi /6.$

\item[Subcase A2.3] $u$ is incident to one boundary class and four interior classes of 
arcs. 

In this case either $\ka (u)\le -\pi /6$ or $u$ is 
a vertex of type $B4(a)$ (see Figure \ref{B4}) or $u$ is of type $AA(u)$ (see
Figure \ref{AA}).  If $\ka (u)\le -\pi /6$ then we can adjust 
$\al _4$ and $\al _5$ by $\pi /12$ to obtain 
$\ka (u),\ka (v)\le -\pi /12$. 
If $u$ is of type $AA(u)$ then $v$ is of type $AA(v^i)$ or $AA(x_i)$, so this does not
occur.
Therefore we may assume $u$ is a vertex of 
type $B4(a)$. In this case every incidence 0 corner on $u$ must have angle $\pi /3$ 
and both incidence 1 corners must have angle $\pi /2$.  There are therefore two 
possibilities as shown in Figures \ref{D4} and \ref{D5}, 
configurations $D4$ and $D5$ respectively.

\item[Subcase A2.4] $u$ is incident to one boundary class and 3 interior 
classes of arcs. 

\begin{figure}
\psfrag{'phi}{$\phi$}
\psfrag{'th}{$\theta$}
\psfrag{'D1}{$\D_1$}
\psfrag{'D2}{$\D_2$}
\psfrag{'D3}{$\D_3$}
\psfrag{'D4}{$\D_4$}
\psfrag{'b}{$\bt$}
\psfrag{v}{$v$}
\psfrag{u}{$u$}
\begin{center}
{ \includegraphics[scale=0.5]{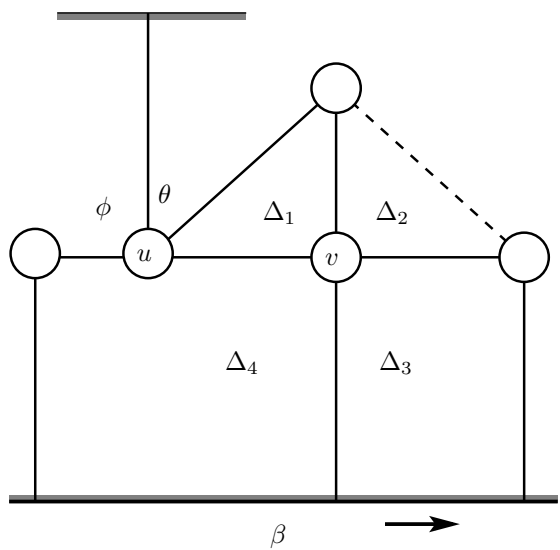} } 
\caption{Configuration $AC4$
\label{AC4}}
\end{center}
\end{figure}

\begin{figure}
\psfrag{'b}{$\bt$}
\psfrag{'n}{$\nu$}
\psfrag{'m}{$\mu$}
\psfrag{v}{$v$}
\psfrag{y}{$y$}
\psfrag{x}{$x$}
\psfrag{u}{$u$}
\psfrag{z}{$z$}
\psfrag{c}{$c$}
\psfrag{d}{$d$}
\psfrag{'D1}{$\D_1$}
\psfrag{'D2}{$\D_2$}
\psfrag{'D3}{$\D_3$}
\psfrag{'D4}{$\D_4$}
\psfrag{'D5}{$\D_5$}
\psfrag{'D6}{$\D_6$}
\psfrag{a1}{$\al_1$}
\psfrag{a2}{$\al_2$}
\psfrag{a3}{$\al_3$}
\psfrag{a4}{$\al_4$}
\psfrag{a5}{$\al_5$}
\psfrag{a6}{$\al_6$}
\psfrag{a7}{$\al_7$}
\psfrag{a8}{$\al_8$}
\psfrag{a9}{$\al_9$}
\psfrag{a10}{$\al_{10}$}
\begin{center}
{ \includegraphics[scale=0.5]{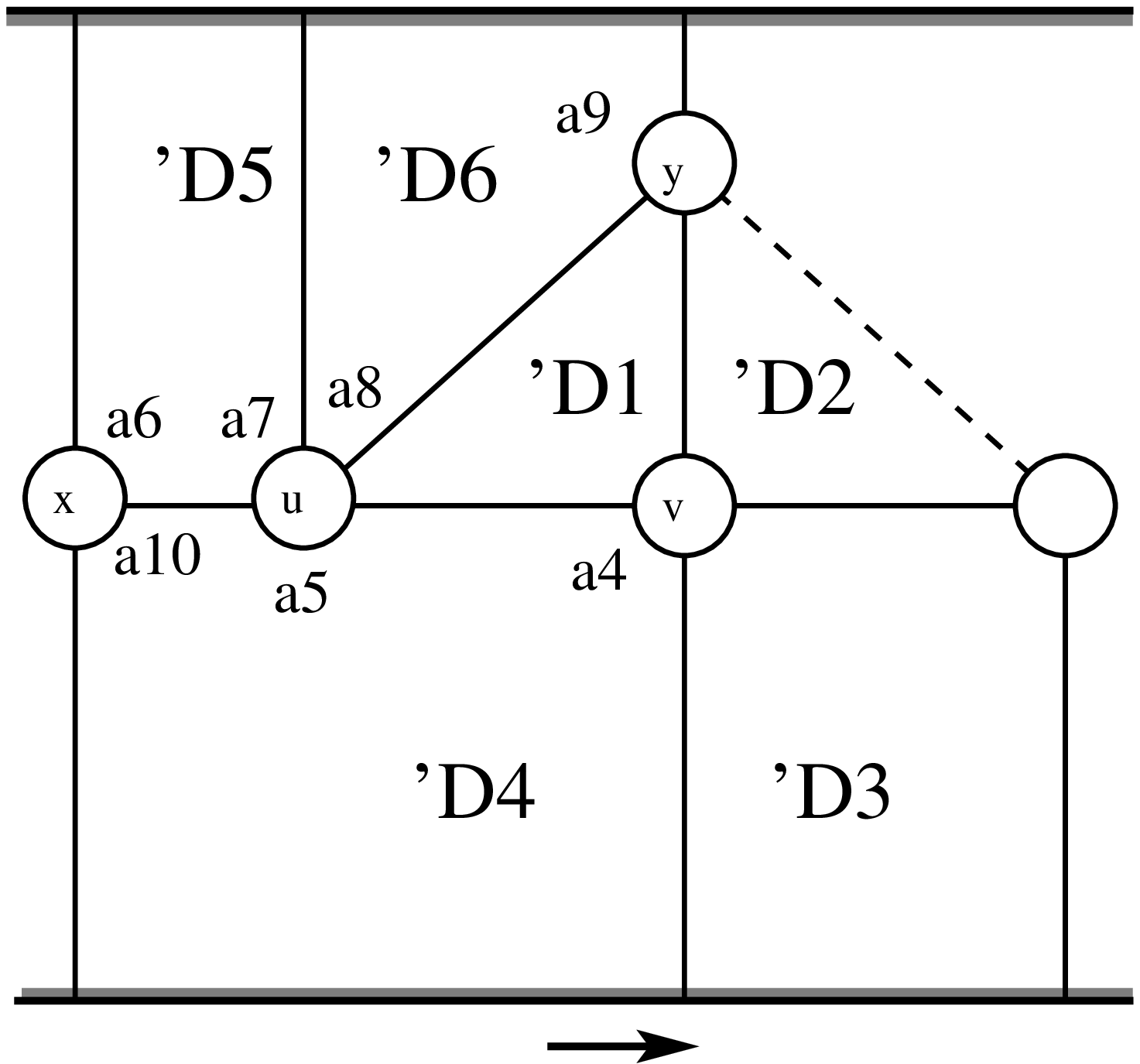} } 
\caption{Configuration $D^\prime 6$\label{D'6}}
\end{center}
\end{figure}

In this case we have the situation of Figure \ref{AC4}.  
\noindent
The 
angles $\theta $ and $\phi $ in Figure \ref{AC4} are both either 
equal to $\pi/2$ or at 
least $5\pi /6$.  If $\phi =\theta =\pi /2$ then $u$ is of type $C0(v)$ and we 
have configuration $D^{\prime } 6$ of Figure \ref{D'6}. 
\noindent
We defer further consideration of $D^{\prime } 6$
until the end of the proof. 
If $\phi =5\pi /6$ and $\theta =\pi /2$ we have 
configuration $D7$, of Figure \ref{D7}, with $u$ of type $CC^{\pm}(v)$.  If $\phi =\pi /2$ and 
$\theta =5\pi /6$ then we have configuration $D8$, of Figure \ref{D8},  
with $u$ of type $CC(v)^{\pm}$.  If 
$\phi \ge 5\pi /6$ and $\theta \ge 5\pi /6$ then $\s (u)\ge  
{2\pi / 3} + {5\pi / 3} = {7\pi / 3}$ so we may adjust $\al _4$ 
and $\al _5$ by $\pi /6$ to give $\ka (u),\ka (v)\le -\pi /6.$ 
\ed

\begin{figure}
\psfrag{'xi}{$\xi$}
\psfrag{'th}{$\theta$}
\psfrag{'D1}{$\D_1$}
\psfrag{'D2}{$\D_2$}
\psfrag{'D3}{$\D_3$}
\psfrag{'D4}{$\D_4$}
\psfrag{'b}{$\bt$}
\psfrag{v}{$v$}
\psfrag{u}{$u$}
\begin{center}
{ \includegraphics[scale=0.5]{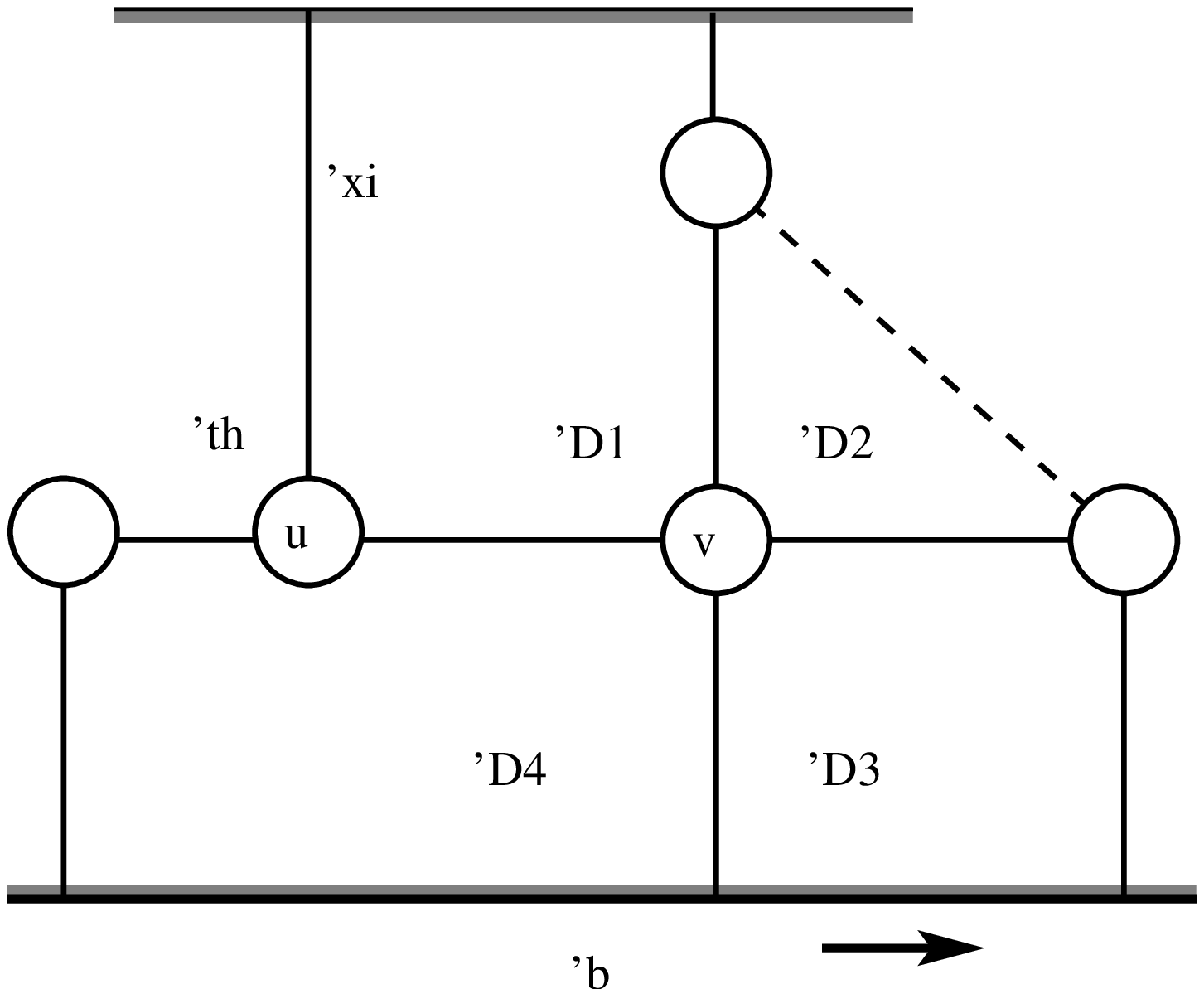} } 
\caption{
\label{ACB}}
\end{center}
\end{figure}

This completes consideration of Case A2, 
configuration $CC^+$ with $\D _1$ an interior 
region (except for configuration $D^{\prime } 6$).  
Now suppose that $\D _1$ is a boundary region.  As $v$ is of type 
$CC^+$ we have $\rho (\D _1)=3$ and $\bt (\D _1)=1$. \ Hence the angle on
the corner of $u$ in $\D _1$ is $5\pi /6$.  If $u$ is incident to a further
boundary class of arcs (other than that in $\pd \D _1)$ then 
$\s (u)\ge 5\pi /6+\pi /3+\pi /2+\pi =8\pi /3$ so we may adjust 
$\al _4$ and $\al _5$ by $\pi /3$ to give 
$\ka (u),\ka (v)\le -\pi /3$.  Hence we may assume $u$ is incident to no
further boundary classes of arcs. \ We have the situation of Figure 
\ref{ACB},  with $\xi $ be the boundary class incident at $u$.  
Now, from Lemma \ref{wv} and the fact that $v$ is not joined to an
$H^\Lm$--vertex, it follows that $|\xi |\le ml /2+l-1$ and so 
$u$ is adjacent to at least one more vertex.  Hence $\theta \ge 5\pi /6$.  If 
$\theta >5\pi /6$ then $\theta \ge \pi $, using Lemmas \ref{AA10}
and \ref{AA11}, so $\s (u)\ge 13\pi /6$ and we may adjust 
$\al _4$ and $\al _5$ by $\pi /12$ to obtain 
$\ka (u),\ka (v)\le -\pi /12$.  Hence we may assume $\theta =5\pi /6$ in
which case we have configuration $D9$, with $u$ of type $CC^{\pm}(v).$ 
\noindent

\begin{figure}
\psfrag{'i}{$\io$}
\psfrag{'phi}{$\phi$}
\psfrag{'th}{$\theta$}
\psfrag{'b}{$\bt$}
\psfrag{'n}{$\nu$}
\psfrag{'m}{$\mu$}
\psfrag{v}{$v$}
\psfrag{u}{$u$}
\psfrag{pi/2}{$\pi/2$}
\begin{center}
{ \includegraphics[scale=0.5]{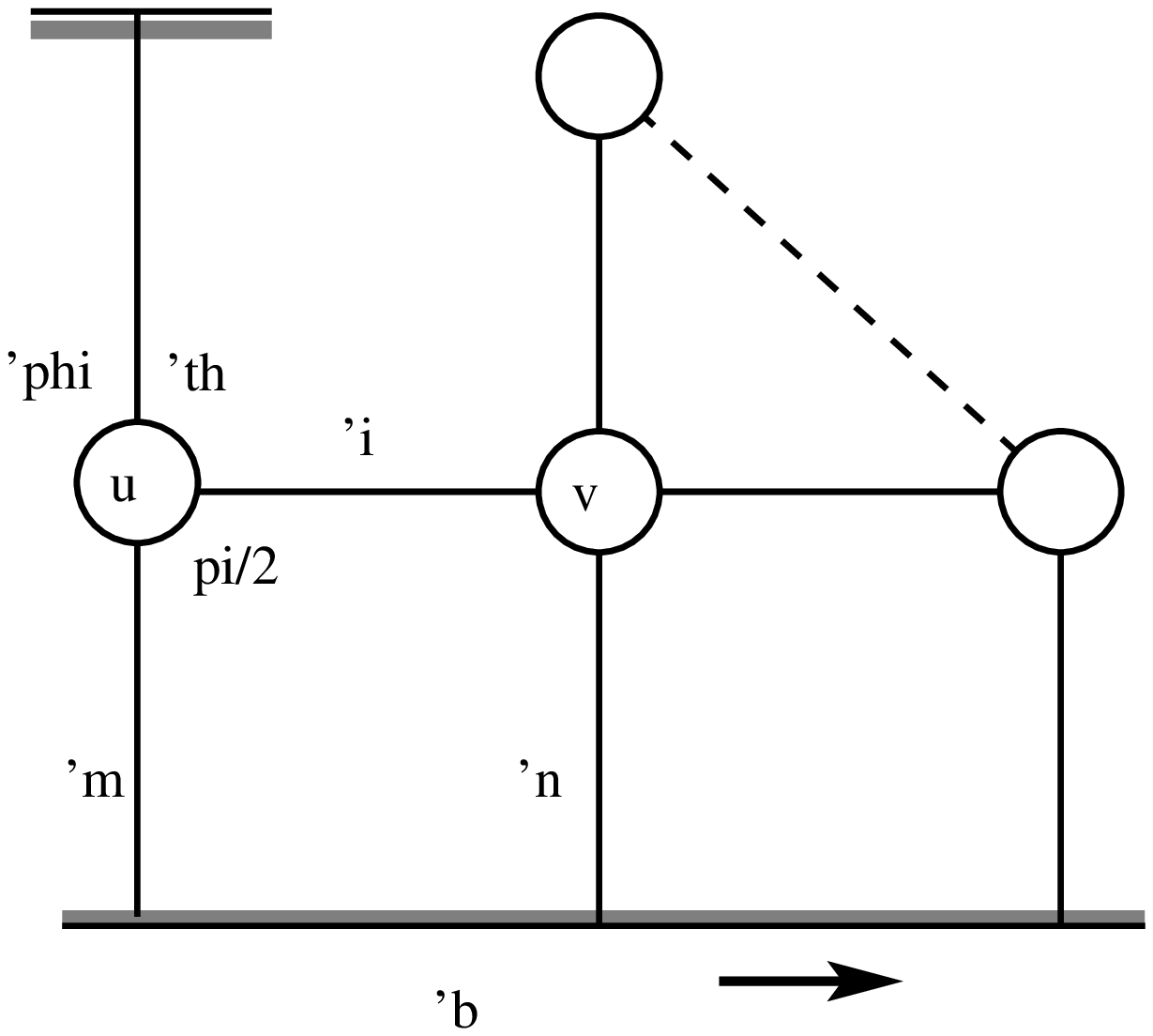} } 
\caption{
\label{COB}}
\end{center}
\end{figure}

It remains to consider configurations $C0$ to complete Case A2 (except for
$D^\prime 6$ which is treated later). 
In Case A2, configuration $C0$, we 
have the situation of Figure \ref{COB} with 
$\s (u)=\theta +\phi +\pi /2$ and $\ka(v)\le -\pi/6$. \ If $u$ is incident to any region $\D $
with $\x (\D )\le 0$ then, since $\x (\D _1)=1$, Lemma \ref{AA0} implies that
either $\theta \ge 11\pi /6$ or $\phi \ge 2\pi $ so 
$\s (u)\ge 20\pi /6$.  Hence we may adjust $\al _4$ and $\al _5$ by 
$\pi /2$ to obtain $\ka (u)$, $\ka (v)\le -\pi /6$.  Hence we may 
assume $\x (\D )=1$, for all regions $\D $ meeting $u$.
\noindent

We consider first the case that $\D _1$ is an interior region, that is 
$\D _1$ is triangular.  In this case $\theta =\pi /3+\psi $, for some angle
$\psi $ with $\psi \ge \pi /2$.  As $\phi \ge \pi $ we have 
$\phi +\psi \ge 3\pi /2$.  If $\phi +\psi =3\pi /2$ then $\psi =\pi /2$ and 
$\phi =\pi $ and it is easy to check that we have one of configurations $D0$, $D1$ 
or $D2$. \ If $\psi>\pi/2$ then $\psi \ge 5\pi/6$, because if $\psi=2\pi/3$ then $u$ is
of type $AA(x^i)$ and $v$ is forced to be of type $AA(v^i)$.
Also  $\phi \ge 7\pi /6$ 
(in fact $\phi \ge 4\pi /3$ unless $u$ is of type $AC^+(x)$, $AC^-(y)$ or $AA(x^i)$).  In both cases we may 
adjust $\al _4$ and $\al _5$ by $5\pi /12$ to obtain 
$\ka (u),\ka (v)\le -\pi /12.$ 

\begin{figure}
\psfrag{'phi}{$\phi$}
\psfrag{'th}{$\theta$}
\psfrag{'D1}{$\D_1$}
\psfrag{'D2}{$\D_2$}
\psfrag{'D3}{$\D_3$}
\psfrag{'D4}{$\D_4$}
\psfrag{'D5}{$\D_5$}
\psfrag{'D6}{$\D_6$}
\psfrag{'b}{$\bt$}
\psfrag{v}{$v$}
\psfrag{u}{$u$}
\psfrag{z}{$z$}
\psfrag{c}{$c$}
\psfrag{d}{$d$}
\begin{center}
{ \includegraphics[scale=0.5]{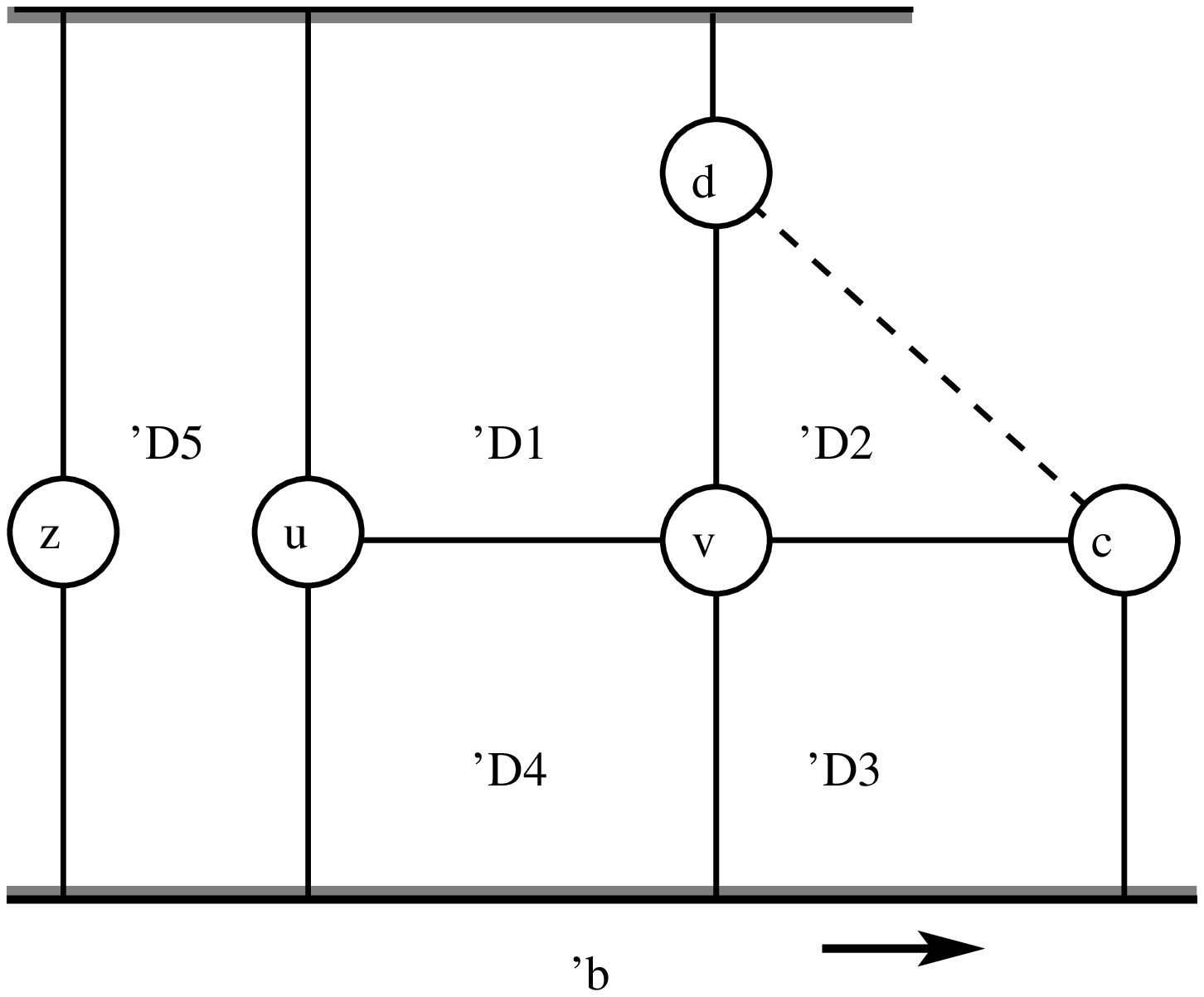} } 
\caption{Configuration $D3^\prime$\label{D3'}}
\end{center}
\end{figure}

\begin{figure}
\psfrag{'phi}{$\phi$}
\psfrag{'th}{$\theta$}
\psfrag{'D1}{$\D_1$}
\psfrag{'D2}{$\D_2$}
\psfrag{'D3}{$\D_3$}
\psfrag{'D4}{$\D_4$}
\psfrag{'D5}{$\D_5$}
\psfrag{'D6}{$\D_6$}
\psfrag{'b}{$\bt$}
\psfrag{v}{$v$}
\psfrag{u}{$u$}
\psfrag{z}{$z$}
\psfrag{c}{$c$}
\psfrag{d}{$d$}
\begin{center}
{ \includegraphics[scale=0.5]{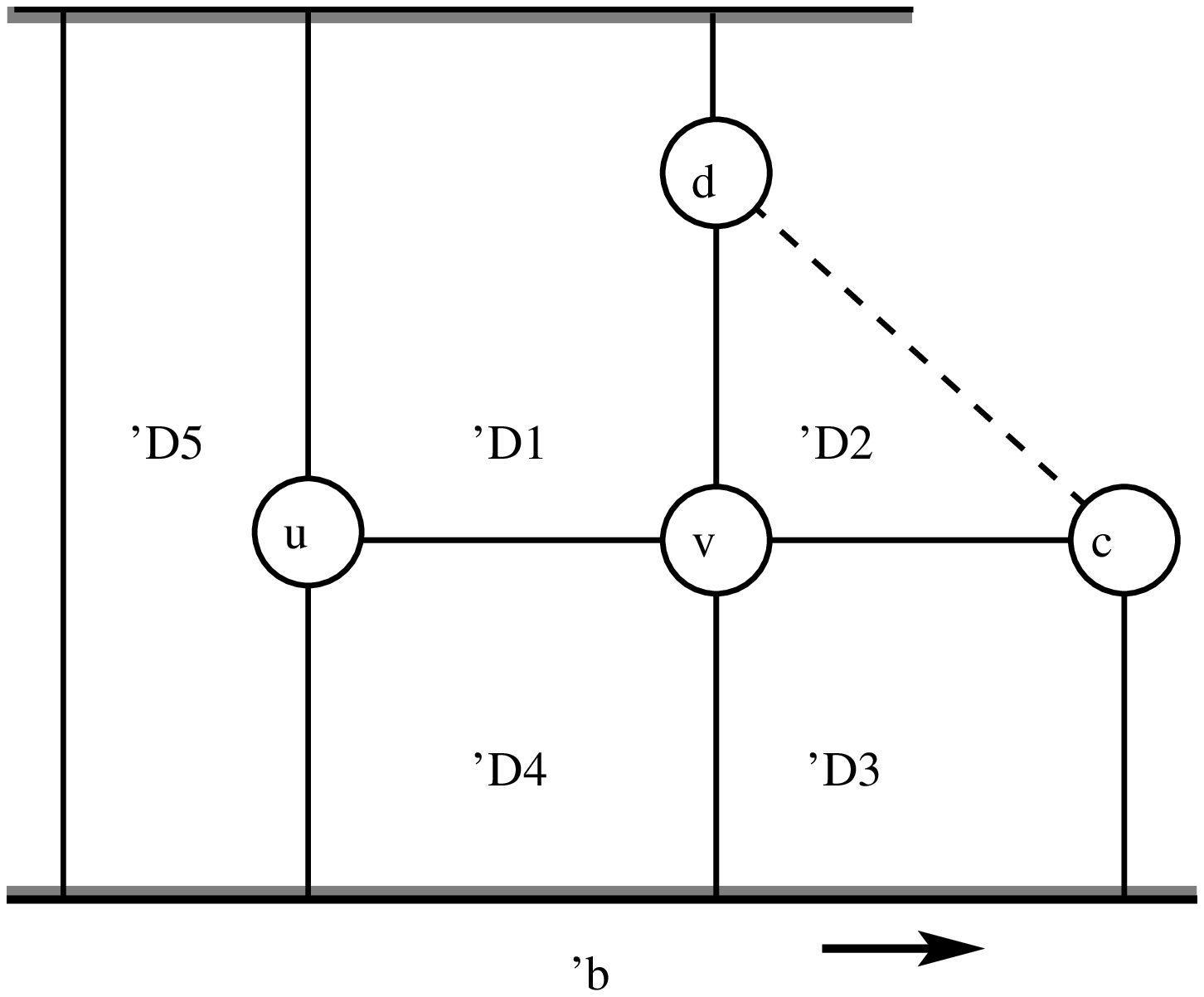} } 
\caption{Configuration $D3^\pprime$\label{D3''}}
\end{center}
\end{figure}

Now suppose that $\D _1$ is a boundary region.  Then, in Figure \ref{COB},
$\theta=5\pi/6$ and $\phi \ge \pi$. If $\phi =\pi$ then  
we have configurations $D3$, $D3^\prime $ or $D3^{\prime \prime} $ of
Figures \ref{D3}, \ref{D3'} and \ref{D3''}, respectively. 
\noindent
\noindent
Theorem \ref{uv} implies that,
in configurations $D3^\prime $ and
$D3^{\pprime} $, the boundary class incident to $u$ in $\pd \D _4$ has size at most 
$l$ so the boundary class incident to $u$ in $\pd \D _1$ has size 
at least $(m-2)l>(ml/2)+l-1$.
Since 
$u$ is not 
an $H^\Lm$--vertex this contradicts Lemma \ref{wv}. Hence the
configurations $D3^\prime $ and $D3^{\prime \prime} $ cannot occur.
If $\phi \ge \pi$ then $\phi \ge 7\pi/6$ and 
we can adjust $\al _4$ and 
$\al _5$ by $5\pi /12$ to give $\ka (u), \ka (v)\le -\pi /12$. 
 This completes Case A2 and so proves 
Case A: $p=1$,  of the theorem, with the exception of configuration $D^\prime 6$.
\subsection*{{Case B: $p>1$.}}
\begin{figure}
\psfrag{'i}{$\io$}
\psfrag{'b}{$\bt$}
\psfrag{'n}{$\nu$}
\psfrag{'m}{$\mu$}
\psfrag{v}{$v$}
\psfrag{u}{$u$}
\psfrag{'i^p}{$\io^p$}
\psfrag{'b^p}{$\bt^p$}
\psfrag{'n^p}{$\nu^p$}
\psfrag{'m^p}{$\mu^p$}
\psfrag{v^p}{$v^p$}
\psfrag{'i^2}{$\io^2$}
\psfrag{'b^2}{$\bt^2$}
\psfrag{'n^2}{$\nu^2$}
\psfrag{'m^2}{$\mu^2$}
\psfrag{v^2}{$v^2$}
\psfrag{blue}{blue}
\psfrag{green}{green}
\begin{center}
{ \includegraphics[scale=0.7]{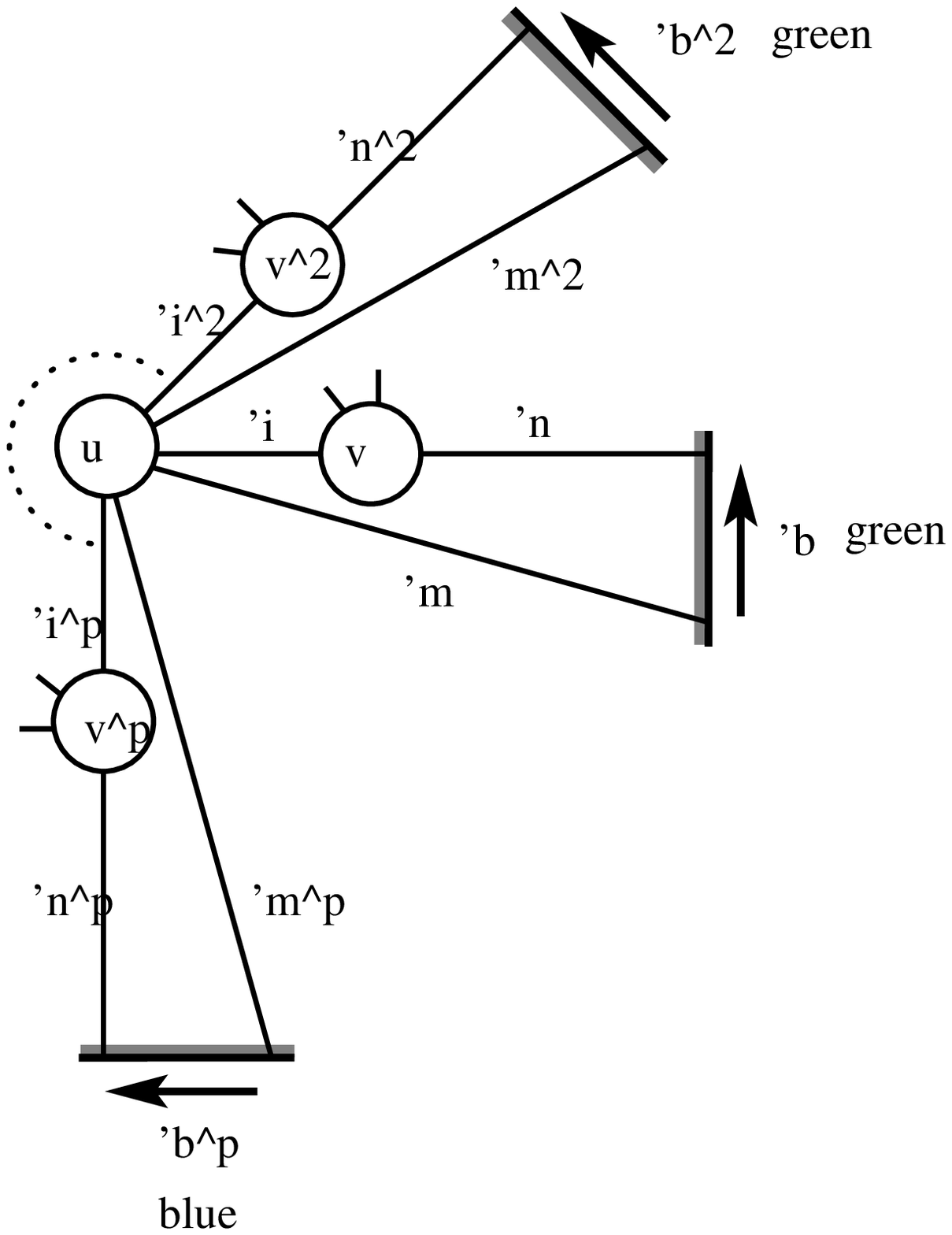} } 
\caption{
\label{uc}}
\end{center}
\end{figure}

~\\
Suppose $u$ is the vertex of type $u_1$ adjacent to 
$C$--vertices $v^1,\cdots ,v^p$, occuring in cyclic order around $u$,
 with $v=v^1$ and $p>1$.  
Note that it is possible that $u$ is a vertex of type $Cj^\pm(u_2)$ or 
$CC^\pm(u_2)$ but then the corresponding vertex $Cj^\pm(v)$ or $CC^\pm(v)$ is
not among $v^1,\cdots ,v^p$. 
If $v^k$ is the vertex $v$ 
in one of configurations $C0$, $C1^\pm$, $C2^\pm$, $C3^\pm$, 
$C4$ or $CC^\pm$ 
we denote the corresponding 
vertices $u_i$, classes $\mu _i$, $\iota_i$ and regions
$\D _i$ by $u^k_i$, $\mu ^k_i$, 
$\iota ^k_i$ and $\D^k_i$, for $i=1,2$ and $k=1,\cdots ,p.$ 
Hence $u=u_1^k$, for $k=1,\cdots ,p$. Write $\mu^k$ and $\io^k$ for $\mu^k_1$ and 
$\io^k_1$, respectively and recall that $\pd u$ is oriented from $\mu$ to $\io$ inside $\D_4$.
If $\mu^k$ precedes $\iota^k$
in the cyclic order around $u$ then we say that $v^k$ is {\em green} 
and otherwise that 
$v^k$ is {\em blue}. (See Figure \ref{uc}.)

\noindent
 
%
We define 
$$T(v^k)=
\left\{
\begin{array}{ll}
X, & \mbox{ if } v ^k \mbox{  is of type } CC(v)^+;\\
Y, & \mbox{ if } v ^k \mbox{  is of type } CC(v)^-;\\
Z, & \mbox{ if } v ^k \mbox{  is of type } Cj(v), \mbox{ with } 0\le j\le 4:
\end{array}
\right. 
$$
and say that $v^k$ is of {\em type} $X,Y$ or $Z$ as appropriate.

We now compute a lower bound for the sum of  angles 
$\theta ^k$ on $u$ between $i^k$ and
$i^{k+1}$, read in the direction of orientation of $u$
(superscripts modulo $p)$.  We consider separately each of the four 
possible colourings, in  green and blue, of $v^k$ and $v^{k+1}$.
\subsection*{\normalsize Case B1:}
$v^k$ and $v^{k+1}$ are green. 

We treat the different combinations of types of $v^k$ and $v^{k+1}$ as distinct
subcases.
\bd
\item[Subcase B1.1]  $T(v^k) = X$, $Y$\ or $Z$, $T(v^{k+1})=X$, see Figure 
\ref{rrXX}. In the following  Figures vertices with broken  outlines may or may
not be present. If they're absent their incident arcs are replaced by a single arc.
In Figure \ref{rrXX} either one or none of the dotted vertices may be present.
\noindent
\begin{figure}
\psfrag{'i^k}{$\io^k$}
\psfrag{'m^k}{$\mu^k$}
\psfrag{v^k}{$v^k$}
\psfrag{'i^k+1}{$\io^{k+1}$}
\psfrag{'m^k+1}{$\mu^{k+1}$}
\psfrag{v^k+1}{$v^{k+1}$}
\psfrag{u}{$u$}
\psfrag{'phi}{$\phi$}
\psfrag{pi/3}{$\pi/3$}
\begin{center}
{ \includegraphics[scale=0.7]{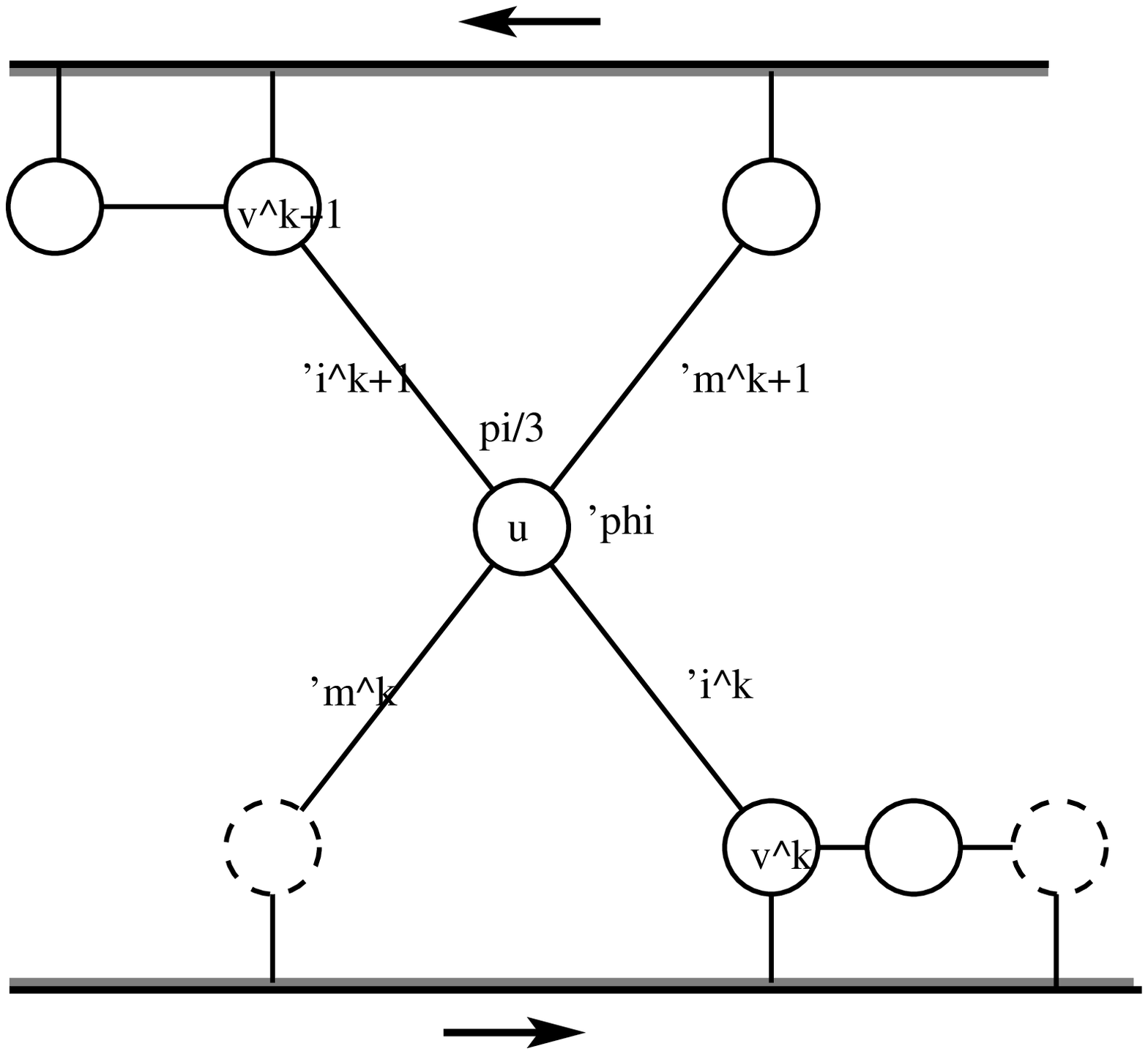} } 
\caption{
\label{rrXX}}
\end{center}
\end{figure}
 
In this subcase as $v^k$ is incident to only one boundary class the angle 
$\phi $ in Figure \ref{rrXX} is at least $\pi /3$, and so 
$\theta ^k\ge 2\pi /3.$

\begin{figure}
\psfrag{'i^k}{$\io^k$}
\psfrag{'m^k}{$\mu^k$}
\psfrag{v^k}{$v^k$}
\psfrag{'i^k+1}{$\io^{k+1}$}
\psfrag{'m^k+1}{$\mu^{k+1}$}
\psfrag{v^k+1}{$v^{k+1}$}
\psfrag{u}{$u$}
\psfrag{'phi}{$\phi$}
\psfrag{pi/2}{$\pi/2$}
\begin{center}
{ \includegraphics[scale=0.7]{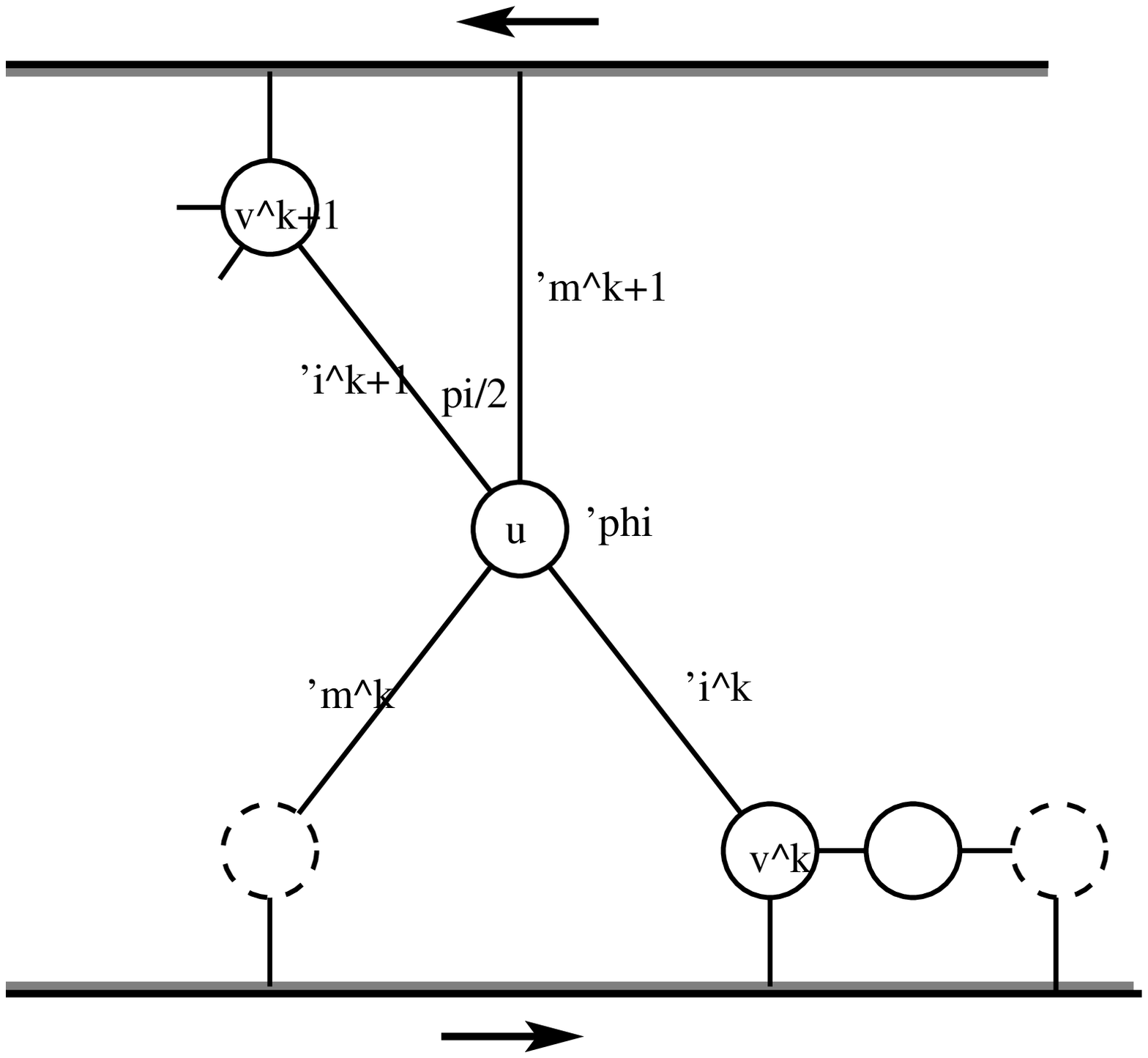} } 
\caption{
\label{rrXY}}
\end{center}
\end{figure}

\item[Subcase B1.2]  $T(v^k) = X$, $Y$ or $Z$, $T(v^{k+1}) = Y$ or $Z$,
(see Figure \ref{rrXY}, where at most one dotted vertex may be present). 
 
\noindent
In this case, as $v^k$ is incident to only one boundary class, 
$\phi \ge 5\pi /6$ so that $\theta ^k\ge 4\pi /3.$
\ed
\subsection*{\normalsize Case B2:}
\begin{figure}
\psfrag{'i^k}{$\io^k$}
\psfrag{'m^k}{$\mu^k$}
\psfrag{v^k}{$v^k$}
\psfrag{'i^k+1}{$\io^{k+1}$}
\psfrag{'m^k+1}{$\mu^{k+1}$}
\psfrag{v^k+1}{$v^{k+1}$}
\psfrag{u}{$u$}
\psfrag{'phi}{$\phi$}
\psfrag{pi/3}{$\pi/3$}
\begin{center}
{ \includegraphics[scale=0.7]{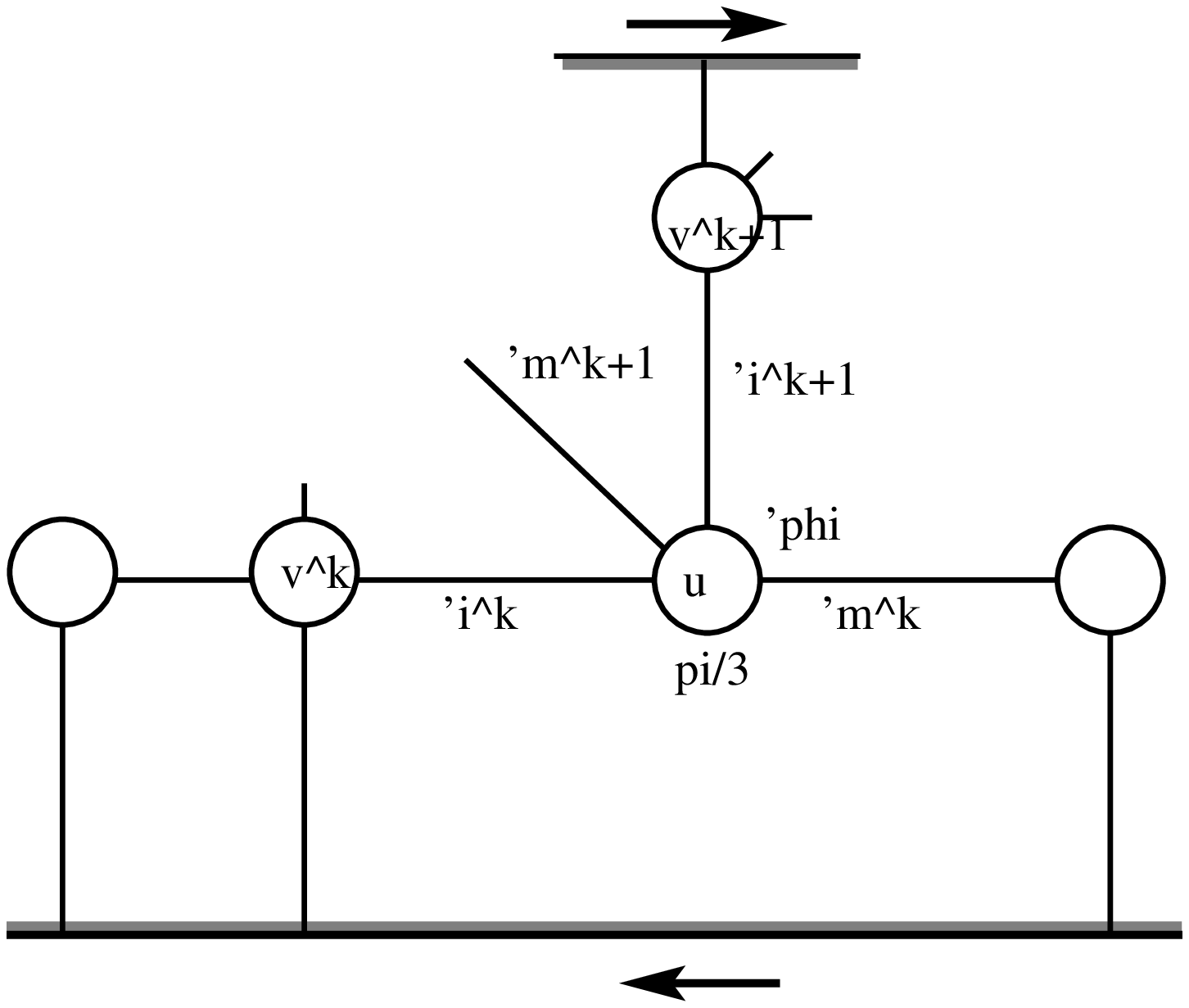} } 
\caption{
\label{bbX}}
\end{center}
\end{figure}

$v^k$ and $v^{k+1}$ are both blue.
\bd
\item[Subcase B2.1]  $T(v^k) = X$, (see Figure \ref{bbX}). 
\noindent
In this case $\phi \ge \pi /3$ so $\theta ^k\ge 2\pi /3.$

\begin{figure}
\psfrag{'i^k}{$\io^k$}
\psfrag{'m^k}{$\mu^k$}
\psfrag{v^k}{$v^k$}
\psfrag{'i^k+1}{$\io^{k+1}$}
\psfrag{'m^k+1}{$\mu^{k+1}$}
\psfrag{v^k+1}{$v^{k+1}$}
\psfrag{u}{$u$}
\psfrag{'phi}{$\phi$}
\psfrag{pi/2}{$\pi/2$}
\begin{center}
{ \includegraphics[scale=0.7]{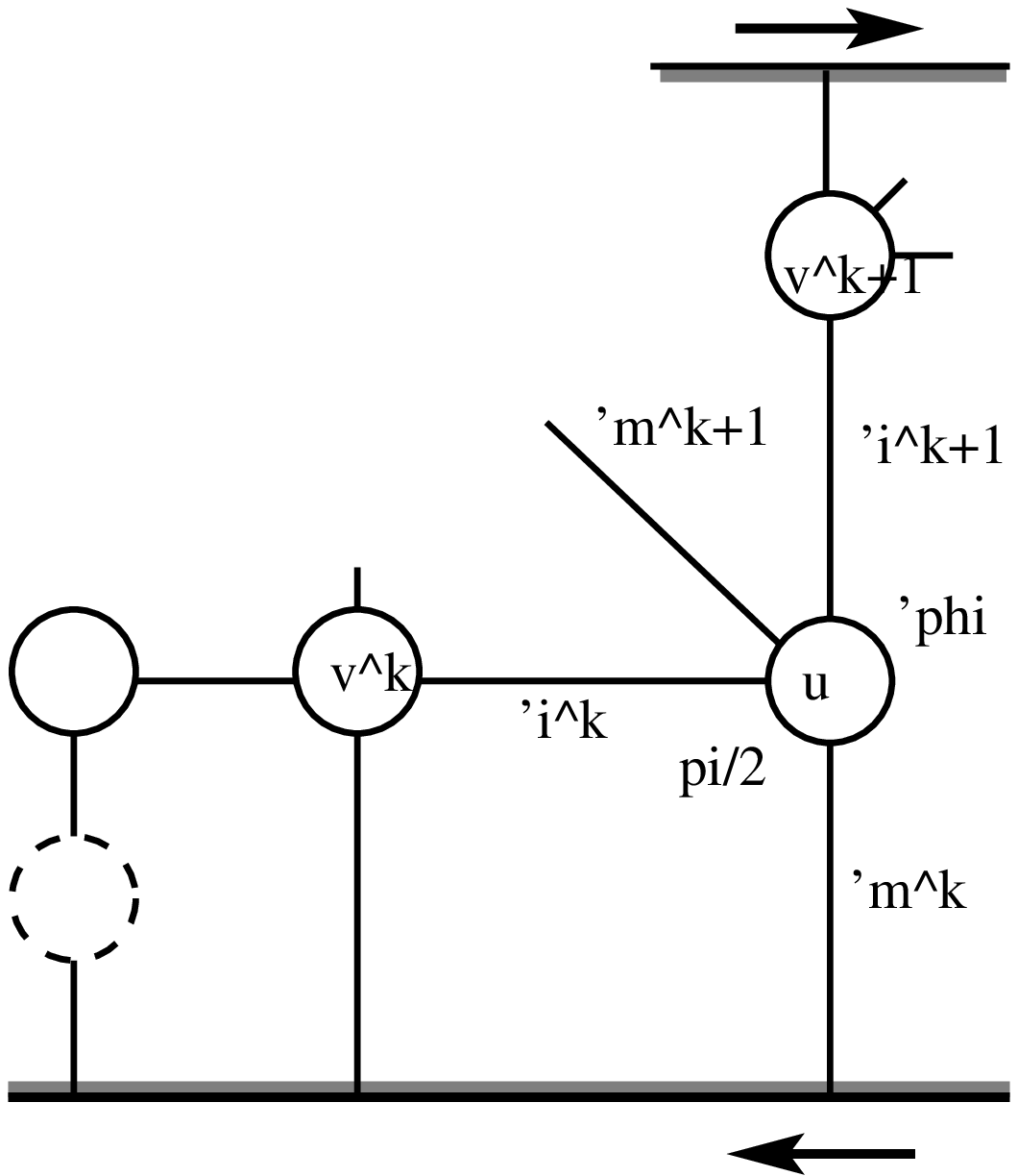} } 
\caption{
\label{bbY}}
\end{center}
\end{figure}

\item[Subcase  B2.2]  $T(v^k) = Y$ or $Z$, (see Figure \ref{bbY}). 
\noindent
As $v^{k+1}$ is incident to only one boundary class we have 
$\phi \ge 5\pi /6$.  Hence $\theta ^k\ge 4\pi /3.$
\ed
\subsection*{\normalsize Case B3:}
$v^k$ is blue, $v^{k+1}$ is green.

\begin{figure}
\psfrag{'i^k}{$\io^k$}
\psfrag{'m^k}{$\mu^k$}
\psfrag{v^k}{$v^k$}
\psfrag{'i^k+1}{$\io^{k+1}$}
\psfrag{'m^k+1}{$\mu^{k+1}$}
\psfrag{v^k+1}{$v^{k+1}$}
\psfrag{u}{$u$}
\psfrag{x}{$x$}
\psfrag{'phi}{$\phi$}
\psfrag{pi/2}{$\pi/2$}
\psfrag{pi/3}{$\pi/3$}
\begin{center}
{ \includegraphics[scale=0.7]{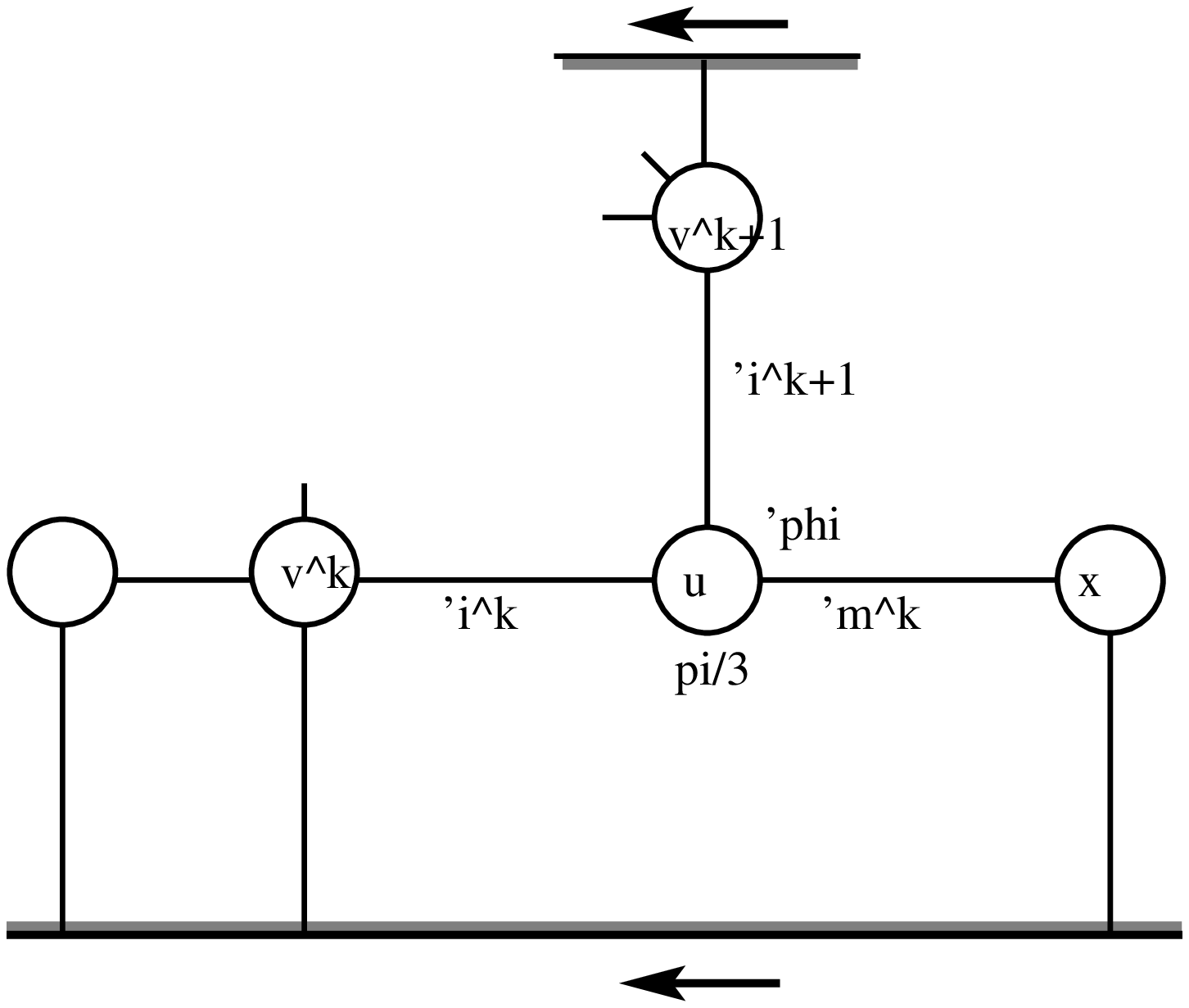} } 
\caption{
\label{brXX}}
\end{center}
\end{figure}

\bd
\item[Subcase B3.1]  $T(v^k) = X$, $T(v^{k+1}) = X$, (see Figure 
\ref{brXX}). 
\noindent
Note that the vertex $x $ in Figure \ref{brXX} cannot be 
$v^{k+1}$, given that we have defined $u=u^{k+1}_1$ 
with respect to the orientation of 
the boundary component $\bt^{k+1} $ of $\pd \S $ incident to $v^{k+1}$. \ We 
have $\phi \ge \pi /3$, so $\theta ^k\ge 2\pi /3.$

\begin{figure}
\psfrag{'i^k}{$\io^k$}
\psfrag{'m^k}{$\mu^k$}
\psfrag{v^k}{$v^k$}
\psfrag{'i^k+1}{$\io^{k+1}$}
\psfrag{'m^k+1}{$\mu^{k+1}$}
\psfrag{v^k+1}{$v^{k+1}$}
\psfrag{u}{$u$}
\psfrag{x}{$x$}
\psfrag{'phi}{$\phi$}
\psfrag{pi/2}{$\pi/2$}
\psfrag{pi/3}{$\pi/3$}
\begin{center}
{ \includegraphics[scale=0.7]{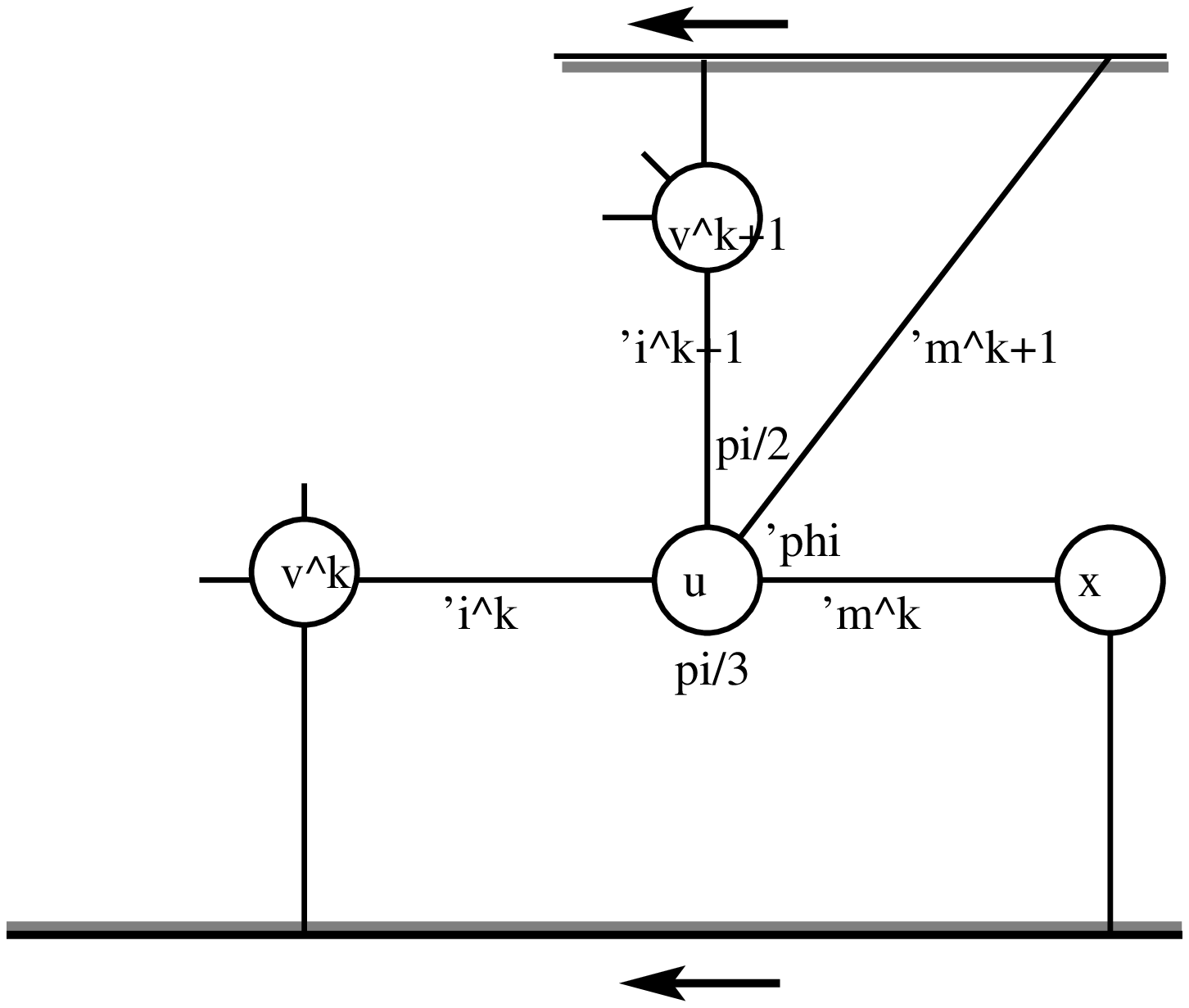} } 
\caption{
\label{brXY}}
\end{center}
\end{figure}

\item[Subcase B3.2]  $T(v^k) = X$, $T(v^{k+1}) = Y$ or $Z$, (see Figure 
\ref{brXY}). 
\noindent
As $\phi \ge \pi /2$ we have $\theta ^k\ge 4\pi /3.$

\begin{figure}
\psfrag{'i^k}{$\io^k$}
\psfrag{'m^k}{$\mu^k$}
\psfrag{v^k}{$v^k$}
\psfrag{'i^k+1}{$\io^{k+1}$}
\psfrag{'m^k+1}{$\mu^{k+1}$}
\psfrag{v^k+1}{$v^{k+1}$}
\psfrag{u}{$u$}
\psfrag{x}{$x$}
\psfrag{'phi}{$\phi$}
\psfrag{pi/2}{$\pi/2$}
\psfrag{pi/3}{$\pi/3$}
\begin{center}
{ \includegraphics[scale=0.7]{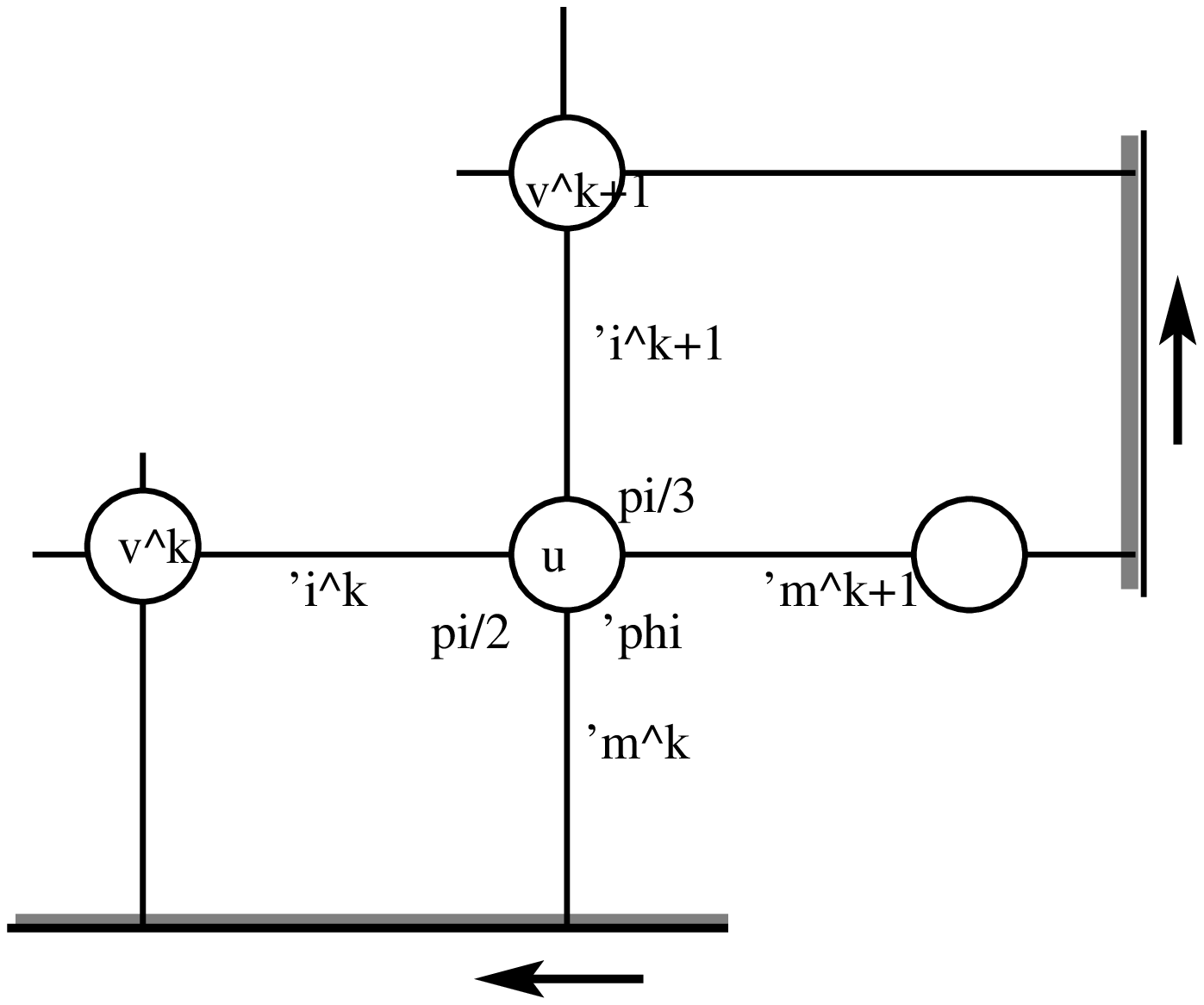} } 
\caption{
\label{brYX}}
\end{center}
\end{figure}

\item[Subcase B3.3]  $T(v^k) = Y$ or $Z$, $T(v^{k+1}) = X$, (see Figure 
\ref{brYX}) 
\noindent
Here $\phi \ge \pi /2$ so $\theta ^k\ge 4\pi /3.$

\begin{figure}
\psfrag{'i^k}{$\io^k$}
\psfrag{'m^k}{$\mu^k$}
\psfrag{v^k}{$v^k$}
\psfrag{'i^k+1}{$\io^{k+1}$}
\psfrag{'m^k+1}{$\mu^{k+1}$}
\psfrag{v^k+1}{$v^{k+1}$}
\psfrag{u}{$u$}
\psfrag{x}{$x$}
\psfrag{'phi}{$\phi$}
\psfrag{pi/2}{$\pi/2$}
\psfrag{pi/3}{$\pi/3$}
\begin{center}
{ \includegraphics[scale=0.7]{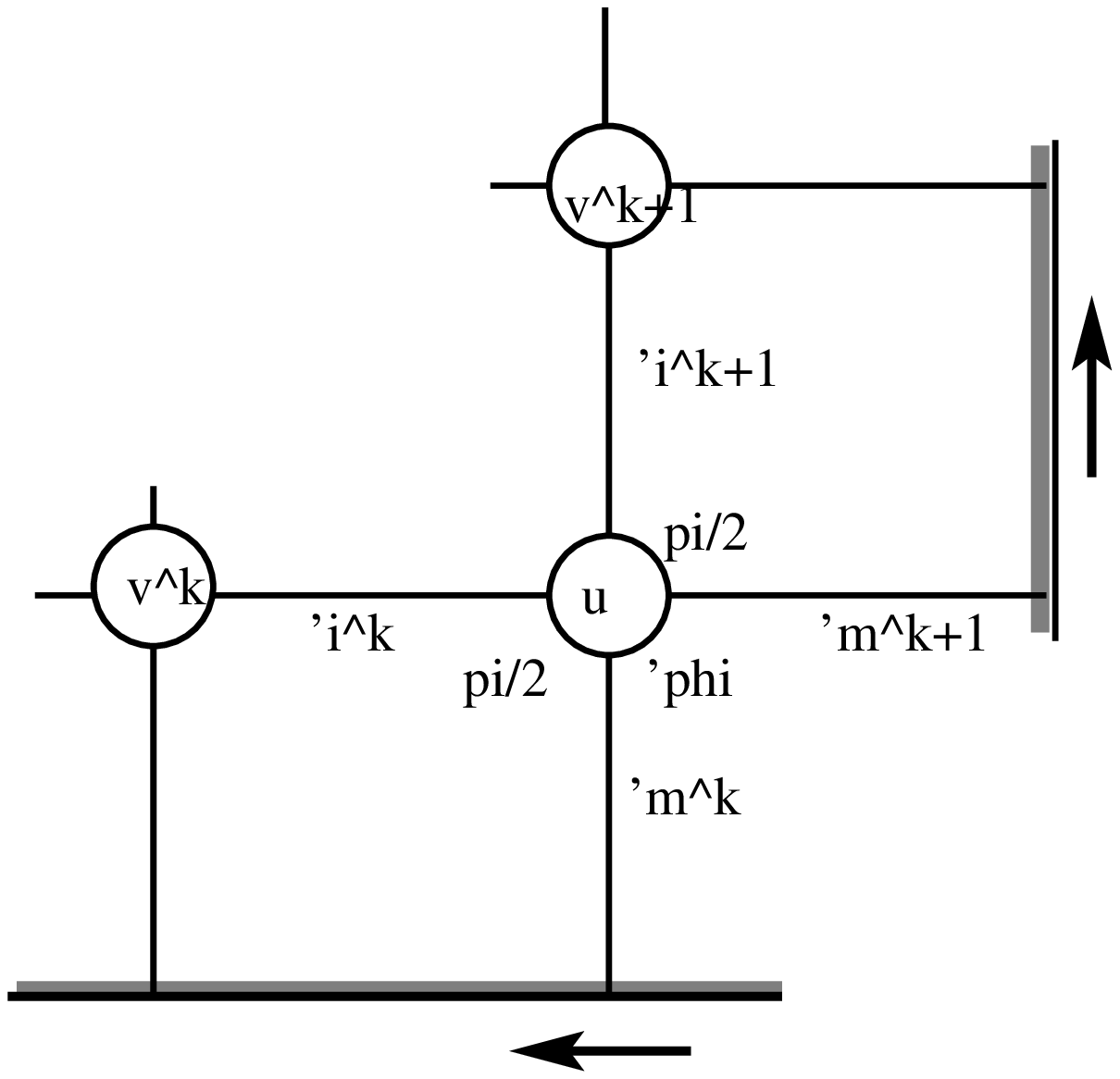} } 
\caption{
\label{brYY}}
\end{center}
\end{figure}

\item[Subcase B3.4]  $T(v^k) = Y$ or $Z$, $T(v^{k+1}) = Y$ or $Z$, (see Figure 
\ref{brYY}) 
\noindent
We have $\phi \ge \pi $ so $\theta ^k\ge 2\pi .$
\ed
\subsection*{\normalsize Case B4:}
\begin{figure}
\psfrag{'i^k}{$\io^k$}
\psfrag{'m^k}{$\mu^k$}
\psfrag{v^k}{$v^k$}
\psfrag{'i^k+1}{$\io^{k+1}$}
\psfrag{'m^k+1}{$\mu^{k+1}$}
\psfrag{v^k+1}{$v^{k+1}$}
\psfrag{u}{$u$}
\psfrag{x}{$x$}
\psfrag{'th^k}{$\theta^k$}
\begin{center}
{ \includegraphics[scale=0.7]{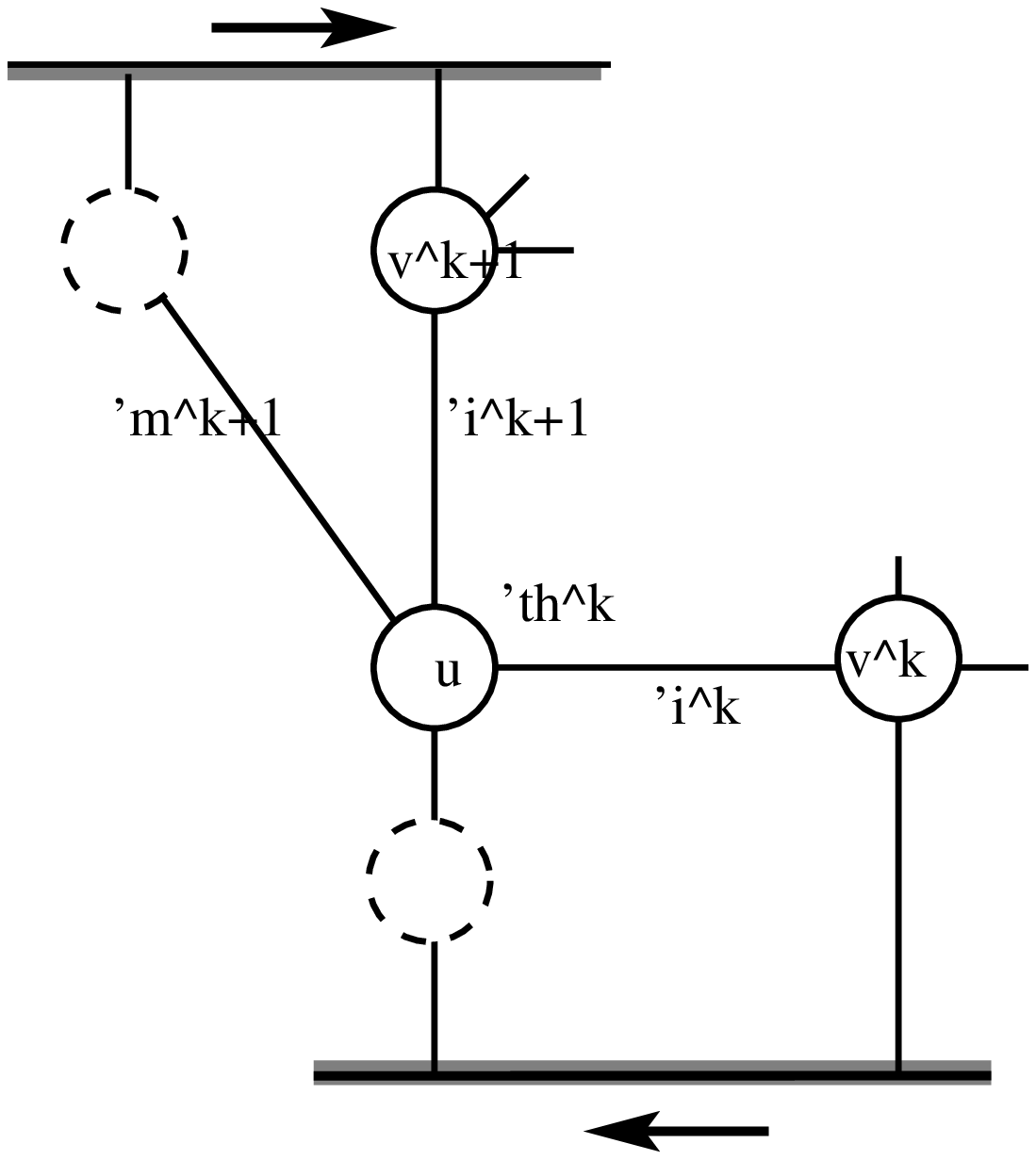} } 
\caption{
\label{rb}}
\end{center}
\end{figure}

$v^k$ is green, $v^{k+1}$ is blue, (see Figure \ref{rb}). 
\noindent
In this case $\theta ^k\ge \pi /3.$
We have therefore a lower bound for $\theta ^k$ in all cases. 

If $v^k$ is green for all $k$ with $c\le k\le d$ and $v^{c-1}$ and $v^{d+1}$ 
are blue (or\ $d+1-c\equiv 0({\rm mod}\ p)$) then $v^c,v^{c+1},\cdots ,v^d$ 
is called a {\em green
component} (of $u$).  {\em Blue components} are defined analogously.

%% file: in4.tex
\subsection*{\normalsize Case B5: {\rm $v^1,\cdots v^p$ are all green.}}
Suppose that $v^1,\cdots v^p$ are all green.  
Let $t=|\{k:T(v^k)=Y\ {\rm or}\ Z\}|$.  
Then $\s (u) \geq  (2\pi /3)p + (2\pi /3)t$.  If $p\geq 4$ then 
$\s (u) \geq  2\pi  + (p-3)(2\pi /3) + (2\pi /3)t$
so $$\sigma (u) \geq  2\pi  + (p+1) {2\pi / 15} + 
(2\pi  / 3)t.
$$
 
\noindent
Noting that $\ka(v^k) \leq  0$ if $T(v^k)=X$  or $Y$ and that 
$\ka(v^k) \leq  {\pi / 3}$ if $T(v^k)=Z$ we see that we can 
adjust angles on $u$ and $v^k$ inside $\D ^{k}_4$ so that 
$\ka (u) \leq  - 
{2\pi / 15}$  and $\ka(v^k) \leq  - {2\pi / 15}$, for all $k$.  (This
follows  because $t \geq  |\{k:T(v^k) = Z\}|$.) 
 
\begin{figure}
\psfrag{x^1}{$x^1$}
\psfrag{'a^1}{$\al^1$}
\psfrag{v^1}{$v^1$}
\psfrag{'phi^1}{$\phi^{1}$}
\psfrag{x^2}{$x^2$}
\psfrag{'a^2}{$\al^2$}
\psfrag{v^2}{$v^2$}
\psfrag{'phi^2}{$\phi^{2}$}
\psfrag{x^3}{$x^3$}
\psfrag{'a^3}{$\al^3$}
\psfrag{v^3}{$v^3$}
\psfrag{'phi^3}{$\phi^{3}$}
\psfrag{u}{$u$}
\begin{center}
{ \includegraphics[scale=0.7]{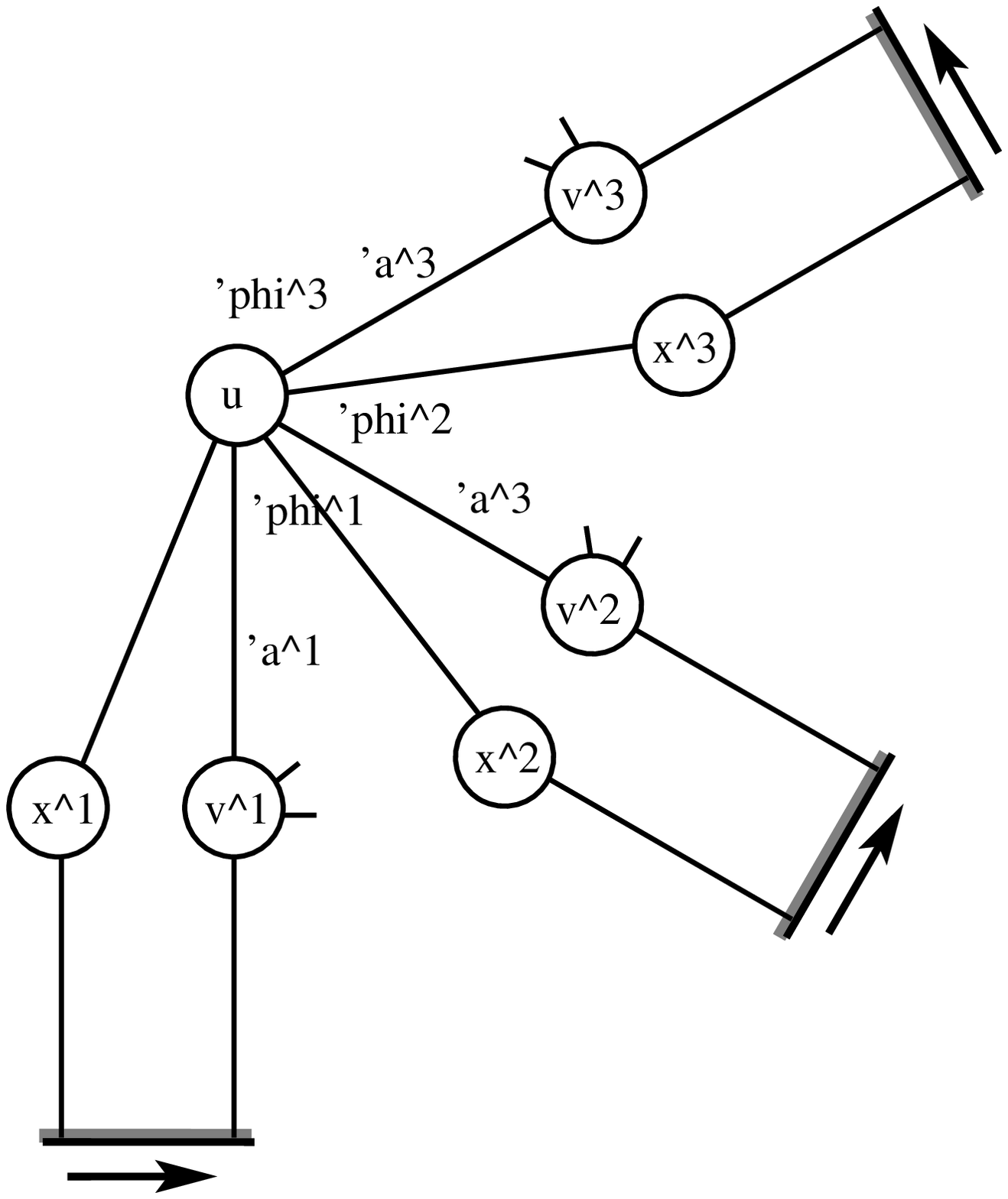} } 
\caption{
\label{r30}}
\end{center}
\end{figure}

\begin{figure}
\psfrag{w^1}{$w^1$}
\psfrag{'a^1_1}{$\al^1_1$}
\psfrag{v^1}{$v^1$}
\psfrag{'phi^1}{$\phi^{1}$}
\psfrag{'i^1}{$\io^1$}
\psfrag{'m^1}{$\mu^1$}
\psfrag{w^2}{$w^2$}
\psfrag{'a^2_1}{$\al^2_1$}
\psfrag{v^2}{$v^2$}
\psfrag{'phi^2}{$\phi^{2}$}
\psfrag{'i^2}{$\io^2$}
\psfrag{'m^2}{$\mu^2$}
\psfrag{pi/2}{$\pi/2$}
\psfrag{u}{$u$}
\begin{center}
{ \includegraphics[scale=0.7]{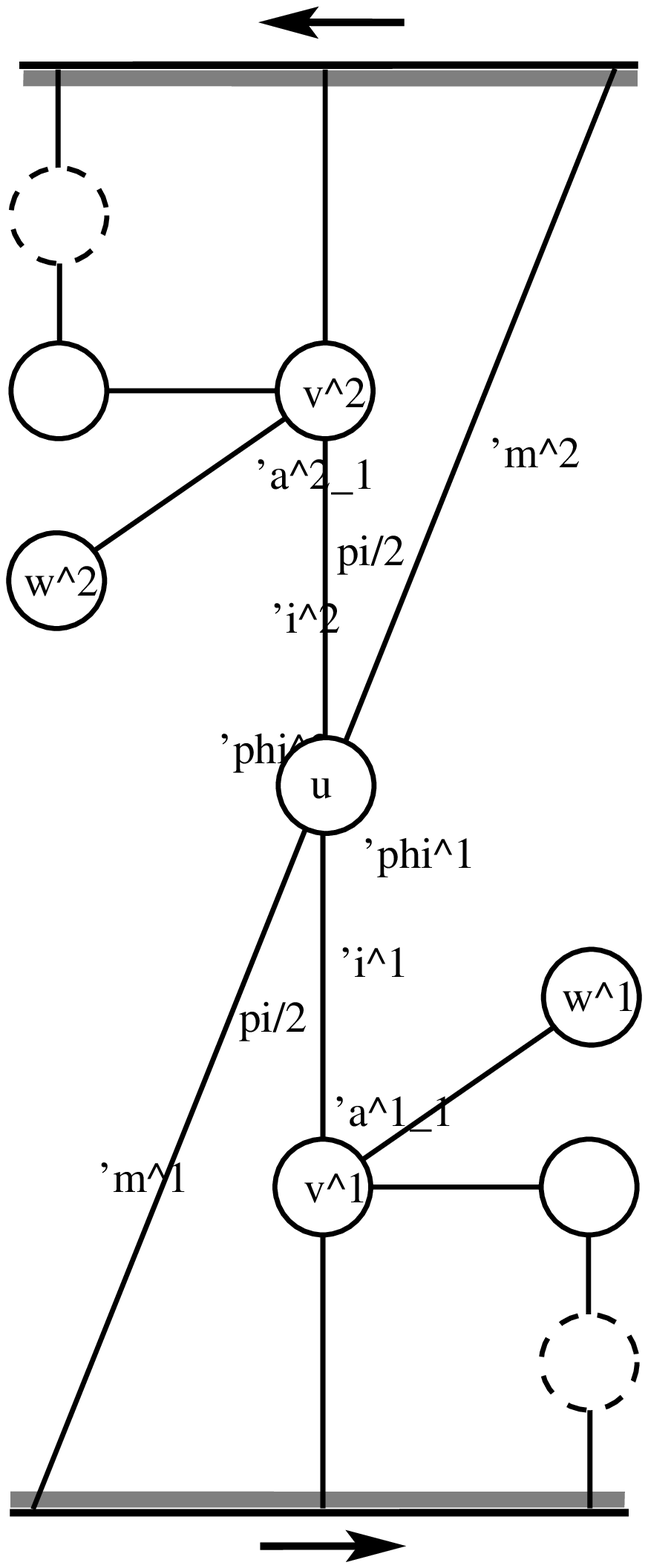} } 
\caption{
\label{r22}}
\end{center}
\end{figure}

Thus given that  all $v^k$'s are green it remains to consider the cases
$2\leq p\leq 3$. 
\bd
\item[Subcase B5.1: $p=3$.] 
In this case if $t\geq 1$ then
$\s(u) \geq  2\pi  + (2\pi /3)t \geq 2\pi  + ({\pi / 3}) + ({\pi / 3})t$.  
It follows that we can adjust angles $\alpha ^k_4$ and $\alpha ^k_5$ for 
$k=1,2,3$ so that $\ka(u)$ and $\ka(v^k) \leq  - {\pi / 12}$.  
We may therefore
suppose that $t=0$.  In this case we have the configuration of figure 
\ref{r30}, 
and $\sigma (u) = 2\pi $ if an only if $\phi ^1 = \phi ^2 = \phi ^3 = 
{\pi / 3}$ (in figure \ref{r30}).  
\noindent
As $v^i$ is incident to only one 
boundary class this occurs only if $v^i$ is incident to $x^{i+1}$, for 
$i=1,2,3$ (superscripts modulo 3).  However the resulting configuration 
cannot occur unless $G$ is of type
$E(2,*,m)$, using Lemma \ref{bound_edge_width}.  Hence $\phi ^i > 
{\pi / 3}$, for some $i$.  As $\rho(\D^i_1)=3$ it follows that $\phi ^i \geq  {2\pi / 3}$, 
for such $i$, so $\s(u) \geq  7{\pi / 3}$ and we may adjust angles so that
$\ka(u)$, $\ka(v^i) \leq  - {\pi / 12}.$ 
\item[Subcase B5.2: $p=2$.]
First suppose $t=2$.  Then we have the 
configuration of Figure \ref{r22} and $\s(u) =  \pi+\phi^1 +\phi^2$, with 
$\phi ^i \geq  5{\pi / 6}$, for $i=1,2$ (as $u$ is not of type $AA(x^i)$, Figure \ref{AA}). 
\noindent
 As $|\mu ^i| \leq  l $ we have 
$\phi ^i = 5{\pi / 6}$ only if $u$ is 
adjacent to $w^i$, for $i=1,2$.  (Otherwise $u$ is incident to fewer than $ml$ 
arcs.)  However if $u$ is adjacent to $w^i$ and incident to no further classes of 
arcs then we have a contradiction to Lemma \ref{bound_edge_width}.  Hence 
$\phi ^i \geq  7\pi/6 $, for $i=1$ or 2 and $\s(u) \geq  3\pi$.  
Therefore we may adjust angles so that $\ka(u)$, $\ka(v^i) \leq  - {\pi / 9}.$

\begin{figure}
\psfrag{w^1}{$w^1$}
\psfrag{'a^1_1}{$\al^1_1$}
\psfrag{v^1}{$v^1$}
\psfrag{'phi^1}{$\phi^{1}$}
\psfrag{'i^1}{$\io^1$}
\psfrag{'m^1}{$\mu^1$}
\psfrag{w^2}{$w^2$}
\psfrag{'a^2_1}{$\al^2_1$}
\psfrag{'D^2_1}{$\D^2_1$}
\psfrag{v^2}{$v^2$}
\psfrag{x^2}{$x^2$}
\psfrag{'phi^2}{$\phi^{2}$}
\psfrag{'i^2}{$\io^2$}
\psfrag{'m^2}{$\mu^2$}
\psfrag{pi/2}{$\pi/2$}
\psfrag{pi/3}{$\pi/3$}
\psfrag{u}{$u$}
\psfrag{b'}{$\bt^\prime$}
\begin{center}
{ \includegraphics[scale=0.7]{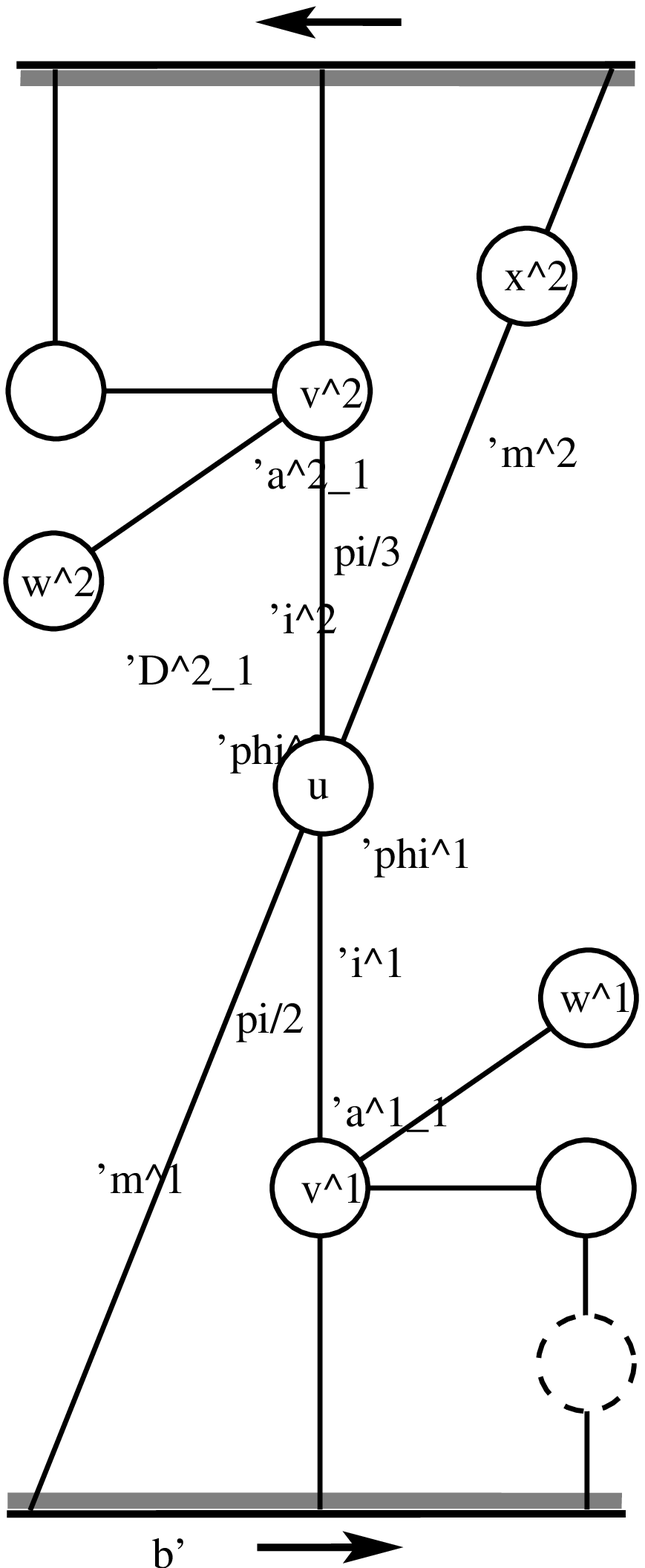} } 
\caption{
\label{r21}}
\end{center}
\end{figure}

Next suppose $t=1$.  Then we have the configuration of Figure \ref{r21},
where $\phi ^1 \geq  {\pi / 3}$ and $\phi ^2 \geq  
5{\pi / 6}$ and $\s(u) =5\pi/6 +\phi^1+\phi^2 $ (as $u$ is not of type $AA(x^i)$, Figure \ref{AA}).  
\noindent
If $\phi^2=5\pi/6$
either $w^2$ is adjacent
to $\bt^1$ and there are no arcs incident to $u$ between $\iota ^2$ and 
$\mu ^1$ or $u$ is adjacent to $w^2$ via the unique class of arcs incident to $u$ 
between $\iota ^2$ and $\mu ^1$. In both of these cases, since $|\mu ^i|\leq l$
there must be arcs incident to $u$ between $\iota ^1$ and $\mu ^2$. 
If none of these arcs are boundary arcs then
there must be  at least two classes of arcs incident
to $u$ between $\iota ^1$ and $\mu ^2$ 
from
which it follows that $\phi^1\geq \pi$.
On the other hand if there is a boundary arc incident to $u$ between 
$\iota ^1$ and $\mu^2$ then $\phi^1\geq 4\pi/3$. Hence $\phi^2=5\pi/6$
implies that $\s(u)\geq 8\pi/3$.

If $\phi^2 > 5\pi/6$ then, as $\al_1^2 =\pi/3$, it follows that 
$\phi^2\geq 7\pi/6$.  If $\phi^2=7\pi/6$ then $\D_1^2$ must be triangular and
there can be at most 2 classes of arc incident to $u$ between $\iota^2$ and $\mu ^1$,
all of which must be interior. In this case there must be  a class of arcs incident
to $u$ between $\iota ^1$ and  $\mu^2$ so $\phi^1\geq 2\pi/3$.
Hence $\phi^2=7\pi/6$ implies that $\s(u)\geq 8\pi/3$. 
If $\phi^2 > 7\pi/6$
then $\phi^2\ge 5\pi/4$, with equality only if $\D^2_1$ is triangular,
there is only one class of arcs incident to $u$ between $\io^2$ and $\mu^1$, and 
$u$ is a vertex of type $AB(u_1)$. In this case there must be at least $2$ classes 
of arc incident to $u$ between $\io^1$ and $\mu^2$, so $\phi^1\ge \pi$ and 
$\s(u)\ge 37\pi/12$. If $\phi^2>5\pi/4$
then $\phi^2\geq 3\pi/2$ so, as $\phi^1\geq \pi/3$, $\s(u)\geq 8\pi/3$.
Therefore $t=1$ implies $\s(u)\geq 8\pi/3$ and we may adjust angles so that
$\ka(u),\ka(v^i)\leq -\pi/9$.

\begin{figure}
\psfrag{w^1}{$w^1$}
\psfrag{'a^1_1}{$\al^1_1$}
\psfrag{v^1}{$v^1$}
\psfrag{'phi^1}{$\phi^{1}$}
\psfrag{'i^1}{$\io^1$}
\psfrag{'m^1}{$\mu^1$}
\psfrag{w^2}{$w^2$}
\psfrag{'a^2_1}{$\al^2_1$}
\psfrag{'D^2_1}{$\D^2_1$}
\psfrag{v^2}{$v^2$}
\psfrag{x^2}{$x^2$}
\psfrag{'phi^2}{$\phi^{2}$}
\psfrag{'i^2}{$\io^2$}
\psfrag{'m^2}{$\mu^2$}
\psfrag{pi/2}{$\pi/2$}
\psfrag{pi/3}{$\pi/3$}
\psfrag{u}{$u$}
\begin{center}
{ \includegraphics[scale=0.7]{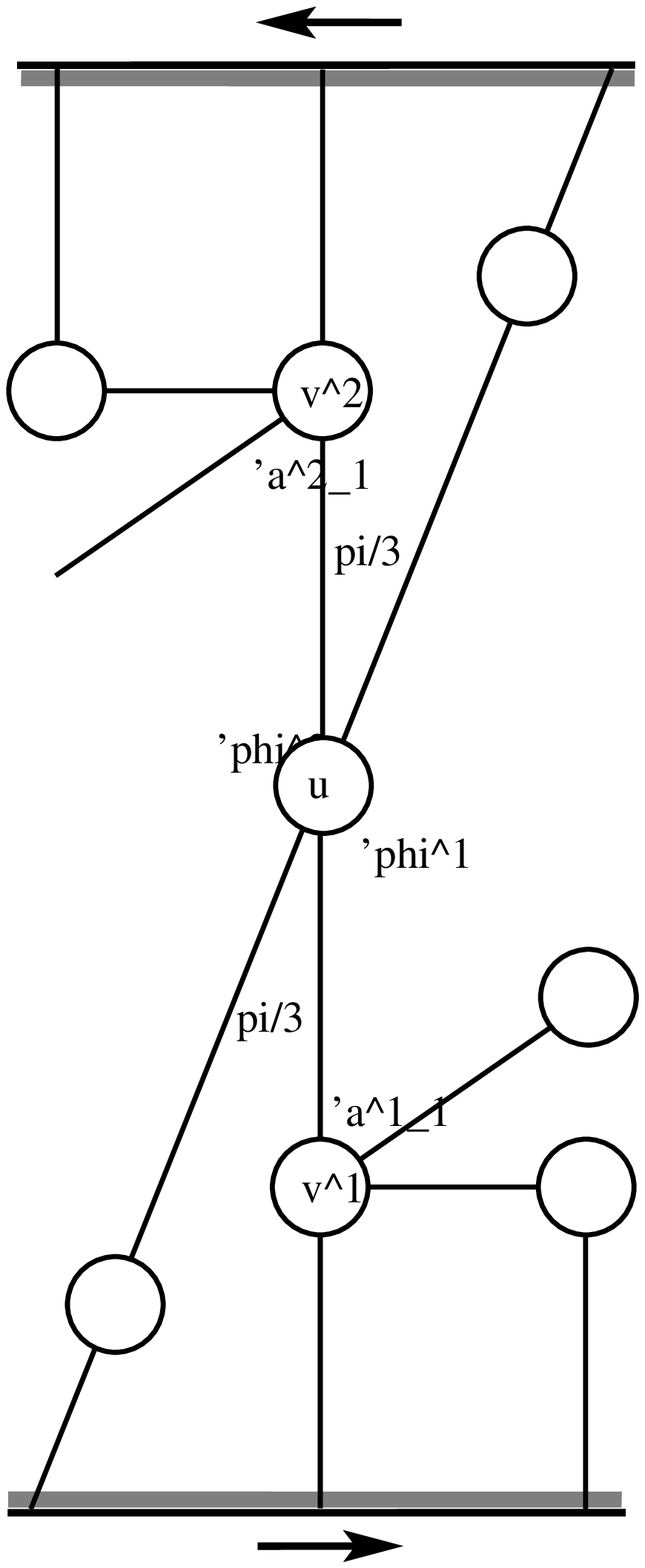} } 
\caption{
\label{r20}}
\end{center}
\end{figure}

Finally suppose $t=0$, in which case we have the configuration of figure 
\ref{r20}.  
\noindent
If $u$ is incident to 3 more classes of arc (than those shown in
Figure \ref{r20}) then $\s(u) \geq  7{\pi / 3}$.  If $u$ is incident to a 
boundary class of arcs then $\phi ^i \geq  4{\pi / 3}$, since $\alpha ^i_1 
= {\pi / 3}$, for $i = 1$ or 2.  In this case $\s(u) \geq  7{\pi / 3}$ 
again.  
Hence, using Lemma \ref{bound_edge_width}, we may 
assume $\s(u) \geq  7{\pi / 3}$ and adjust angles so that $\ka(u)$, $\ka(v^i) 
\leq  - {\pi / 9}$. 

Hence in Case B5 where all $v^k$'s are green we may adjust angles suitably.  It 
follows, by symmetry, that this is also possible if all $v^k$'s are blue. 
\ed

\subsection*{\normalsize Case B6: {\rm there are $q $ green and $q $ blue components.}}

Suppose that there are $q \geq  1$ green and $q \geq  1$ blue components, 
$G_1,..,G_q$ and $B_1,\cdots ,B_q$ respectively, in cyclic order round $u$.
Let $g_j(A)$ be the number of vertices of type $A$ in $G_j$ and 
$b_j(A)$ be the 
number of vertices of type $A$ in $B_j$, where $A$ is $X$, $Y$ or $Z$.
Let $g_j=g_j(X)+g_j(Y)+g_j(Z)$ and $b_j=b_j(X)+b_j(Y)+b_j(Z)$; so 
$$
p = \sum ^q_{j=1} g_j + b_j.
$$

Furthermore let 
$$
g(A) = \sum ^q_{j=1}g_j(A)\hbox{, }\, b(A) = \sum ^q_{j=1}b_j(A)\hbox{ and}
$$
$$
p(A) = g(A) + b(A)\hbox{, where }\ A = X\hbox{, } Y\hbox{ or }\ Z.
$$

\noindent
To complete the proof of case B of the theorem we aim to show that
\begin{equation}\label{su}
\s(u) \geq  2\pi  + p(Z){\pi / 3} + (p+1)/24.
\end{equation}

\noindent
This is sufficient to allow a suitable angle adjustment to be made. In fact
transferring angle of $a$ from the corner of $u$ in $\D_4^i$ to the corner
of $v^i$ in $\D_4^i$, where 
$$a= \left\{
\begin{array}{ll}
9\pi / 24,& \hbox{ if } T(v^i)=Z\\
\pi / 24,& \hbox{ if } T(v^i)=X \hbox{ or } Y
\end{array}
\right. ,
$$
we have $\s (v^i)\geq 2\pi +\pi/24$, for all $i$, and
$$\s (u) \geq 2\pi + p(Z)\pi/3 +(p+1)\pi/24-p(Z)9\pi/24-(p(X)+p(Y))\pi/24$$ 
$$=2\pi+\pi/24.$$
Hence we may adjust angles as described to obtain
$\ka(u),\ka(v^i)\leq-\pi/24$, for all $i$.

A green or blue 
component $v^c,\cdots ,v^{c+d}$ 
is said to {\em begin} with $v^c$ and {\em end} with $v^{c+d}$. 
We examine the effect of the type of the extremal vertices of a component 
on the sum of angles between its $C$--vertices.

If a component $C$ begins with $v^c$ and ends with $v^{c+d}$ then the sum of 
angles around $u$ between $\io^c$ and $\io^{c+1}$, $\io^{c+1}$ and $\io^{c+2},\cdots $, $\io^{c+d-1}$ 
and $\io^{c+d}$ is called the {\em angle} of the component and is denoted $\sigma (C).$ 

Consider a green component $G_j$.  If $G_j$ begins with a $Y$ or $Z$ vertex 
then
$$
\s(G_j) \geq  (4{\pi / 3})(g_j(Y) + g_j(Z) - 1) + (2{\pi / 3})g_j(X).
$$

\noindent
If $G_j$ begins with an $X$ vertex then 
$$
\sigma (G_j) \geq  (4{\pi / 3})(g_j(Y) + g_j(Z)) + 
(2{\pi / 3})(g_j(X)-1).
$$

Now consider a blue component $B_j$.  If $B_j$ ends with an $Y$ or $Z$ vertex
then
$$
\sigma (B_j) \geq  (4{\pi / 3})(b_j(Y) + b_j(Z) - 1) + 
(2{\pi / 3})b_j(X).
$$

\noindent
If $B_j$ ends with an $X$ vertex then
$$
\sigma (B_j) \geq  (4{\pi / 3})(b_j(Y) + b_j(Z)) + 
(2{\pi / 3})(b_j(X)-1).
$$


Let $e(S,T)$ denote the number of blue components which end with a vertex of 
type $S$ and such that the following green component begins with a vertex of 
type $T$, where $S$, $T \in  \{X,Y,Z\}.$ 

Then $q =$ $e(X,X) + e(X,Y) + e(X,Z) +$   
$e(Y,X) + e(Y,Y) + e(Y,Z) +$   
$e(Z,X) + e(Z,Y) + e(Z,Z),$

\noindent
the number of blue components ending with a vertex of type $A$ is
$$
e(A,X) + e(A,Y) + e(A,Z)
$$

\noindent
and the number of green components beginning with a vertex of type $A$ is
$$
e(X,A) + e(Y,A) + e(Z,A).
$$

\noindent
Now $\s(u) = $ [sum of angles  of green components]  

$+$ [sum of angles  of blue components]  

$+$ [sum of angles between green and blue components]  

$+$ [sum of angles between blue and green components]  

$\geq  \left[ \left( 4{\pi / 3}\right)\sum ^q_{j=1}\left( g_j(Y) + g_j(Z) - 
1\right)   + 
\left( 2{\pi / 3}\right) \sum ^q_{j=1}\left( g_j(X)-1\right) \right. $  

$+ \left( 2{\pi / 3}\right) \left( e(X,Y) + e(Y,Y) + e(Z,Y) + e(X,Z) + 
e(Y,Z) + e(Z,Z)\right) $  

$+ \left( 4{\pi / 3}\right) \left( e(X,X) + e(Y,X) + 
e(Z,X)\right) \left. \right] $  

$+ \left[ \left( 4{\pi / 3}\right) \sum ^q_{j=1}\left( b_j(Y) + b_j(Z) - 
1\right)  + \left( 2{\pi / 3}\right) \sum ^q_{j=1}\left( b_j(X) - 
1\right) \right. $  

$+ \left( 2{\pi / 3}\right) \left( e(Y,X) + e(Y,Y) + e(Y,Z) + e(Z,X) + 
e(Z,Y) + e(Z,Z)\right) $  

$+ \left( 4{\pi / 3}\right) \left( e(X,X) + e(X,Y) + 
e(X,Z)\right) \left. \right] $  

$+ \left[ q {\pi / 3}\right] $  

$+ \left[ \left( 2{\pi / 3}\right)  e(X,X) + 
\left( 4{\pi / 3}\right) \left( e(X,Y) + e(X,Z) + e(Y,X) + 
e(Z,X)\right) \right. $  

$+ 2\pi  \left( e(Y,Y) + e(Y,Z) + e(Z,Y) + e(Z,Z)\right) \left. \right] $  

$= \left( 2{\pi / 3}\right)  \sum ^q_{j=1}\left( g_j(Y) + g_j(Z) + g_j(X) +
b_j(Y) + b_j(Z) + b_j(X)\right) $  

$+ \left( 2{\pi / 3}\right)  \sum ^q_{j=1}\left( g_j(Y) + g_j(Z) + b_j(Y) +
b_j(Z)\right) $  

$- 4\pi q + ({\pi / 3})q$  

$+ \left( 10{\pi / 3}\right) \left( e(X,X) + e(X,Y) + e(X,Z) +\right. $   

$e(Y,X) + e(Y,Y) + e(Y,Z) +$   

$e(Z,X) + e(Z,Y) + e(Z,Z)\left. \right) $  

$= \left( 2{\pi / 3}\right) p - \left( {\pi / 3}\right) q + 
\left( 2{\pi / 3}\right) (p(Z)+p(Y)).$  

\noindent
Note that $p\geq 2q$ and if $p$ is odd $p \geq 2q+1$.  Hence if $p\geq 5$ then
$p-3\geq q$ and so  

$\s(u) \geq  2\pi  + (p-3) \left( 2{\pi / 3}\right)  - 
\left( {\pi / 3}\right) q + \left( 2{\pi / 3}\right) p(Z)$  

$\geq  2\pi  + \left( {\pi / 3}\right) (p-3) + 
\left( 2{\pi / 3}\right) p(Z)$  

$\geq 2\pi  + \left( {\pi / 9}\right) (p+1) + 
\left( {\pi / 3}\right) p(Z),$

\noindent
since $3(p-3) \geq  (p+1)$, when $p\geq 5.$

\noindent
Hence (\ref{su}) holds when $p\geq 5$ and 
it remains to check the cases $2\leq p\leq 4.$
\bd
\item[Subcase B6.1] $p = 4$. 

If $p=4$ then $q=1$ or 2 and

$$
\s(u) \geq  \left( 8{\pi / 3}\right)  - \left( {\pi / 3}\right) q + 
\left( 2{\pi / 3}\right) (p(Z) + p(Y)).
$$
We shall show that, when $p=4$, we have 
$$
\s(u) \geq 7\pi/3  + p(Z)\left( {\pi / 3}\right)
$$

\noindent
so  that
$$
\s(u) \geq 2\pi  + p(Z)\left( {\pi / 3}\right)  + (p+1){\pi / 15}
$$

\noindent
as required.
\bd
\item[Subcase B6.1.1] If $p(Y)$  is non--zero  then we have
$$
\s(u) \geq  {8\pi / 3} + \left( 2{\pi / 3}\right) p(Z) > 7\pi/3+(\pi/3)p(Z).
$$
\item[Subcase B6.1.2] 
Suppose $p(Y)= 0$.  Then, if $p(Z)
\geq  1$, we have
$$
\s(u) \geq 2\pi  + \left( 2{\pi / 3}\right) p(Z) \geq  7\pi/3  + 
\left( {\pi / 3}\right) p(Z).
$$

\noindent
If $p(Y) = p(Z) = 0$ and $q = 1$ then
$$
\s(u) \geq  7{\pi / 3} = 7\pi/3  + \left( {\pi / 3}\right) p(Z).
$$

\begin{figure}
\psfrag{x^1}{$x^1$}
\psfrag{'a^1_1}{$\al^1_1$}
\psfrag{v^1}{$v^1$}
\psfrag{'phi^1}{$\phi^{1}$}
\psfrag{x^2}{$x^2$}
\psfrag{'a^2_1}{$\al^2_1$}
\psfrag{v^2}{$v^2$}
\psfrag{'phi^2}{$\phi^{2}$}
\psfrag{x^3}{$x^3$}
\psfrag{'a^3_1}{$\al^3_1$}
\psfrag{v^3}{$v^3$}
\psfrag{'phi^3}{$\phi^{3}$}
\psfrag{x^4}{$x^4$}
\psfrag{'a^4_1}{$\al^4_1$}
\psfrag{v^4}{$v^4$}
\psfrag{'phi^4}{$\phi^{4}$}
\psfrag{u}{$u$}
\begin{center}
{ \includegraphics[scale=0.7]{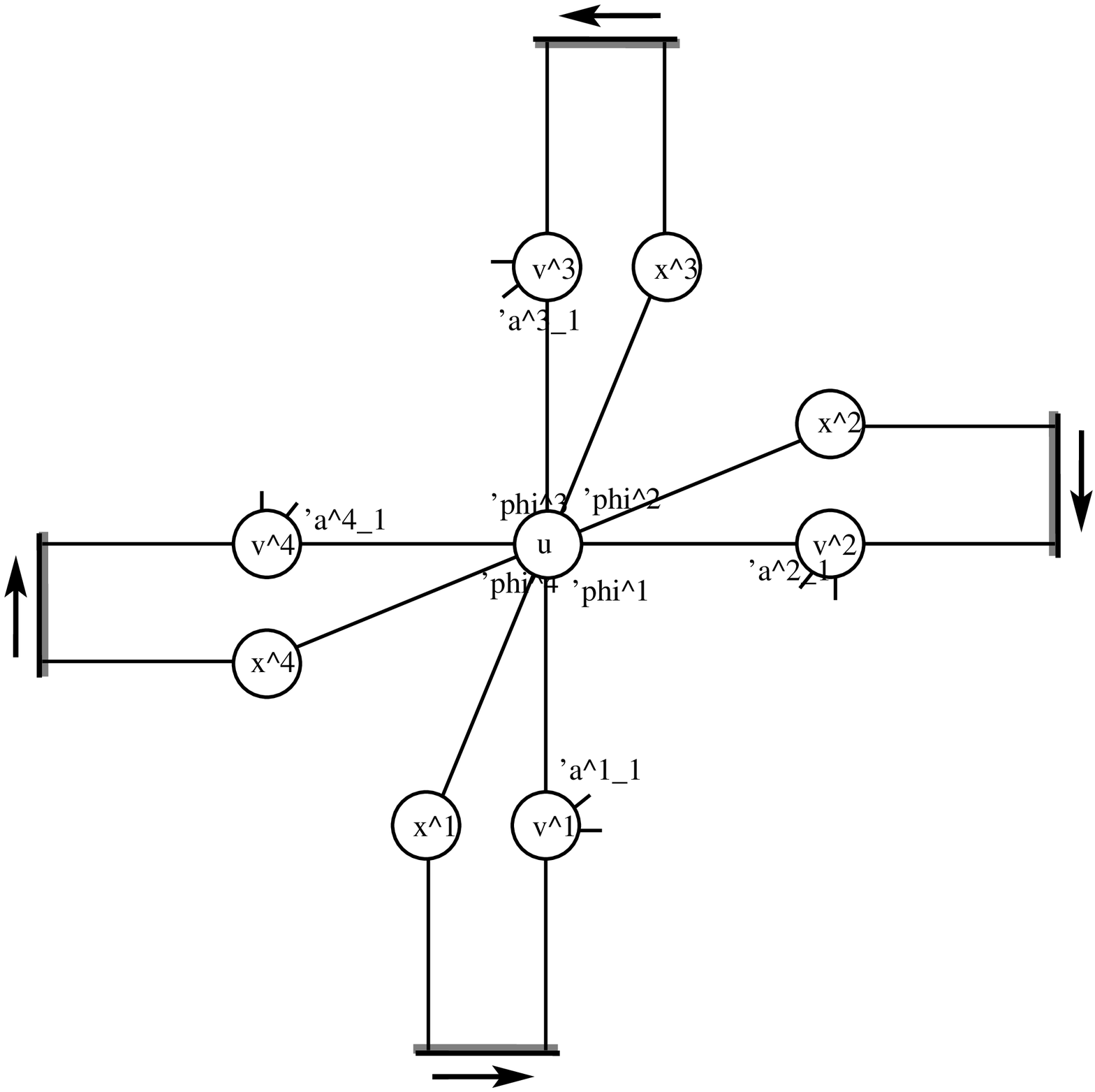} } 
\caption{
\label{p4Z0}}
\end{center}
\end{figure}

\noindent
If $p(Y) = p(Z) = 0$ and $q = 2$ then we have the configuration of figure 
\ref{p4Z0} where $v^1$ and $v^3$ are green, of type $X$, and $v^2$
and $v^4$ are blue of type $X$.  
\noindent
(Note that, in Figure \ref{p4Z0}, 
$\phi ^2$ may equal zero if $x^2=x^3$ and $\phi ^4$ may equal zero if 
$x^1=x^4.)$  Suppose first $\phi ^2=\phi ^4=0$.  As $\alpha ^1_1= 
\alpha ^2_1= \alpha ^3_1=\alpha ^4_1={\pi / 3}$ 
we have $\phi ^1 \geq  2{\pi / 3}$ and $\phi ^3 \geq  {\pi / 3}$, or 
vice--versa, otherwise $G$ is of type $E(2,3,6)$.  Hence $\sum ^4_{i=1} \phi ^i 
\geq  \pi $.  Now suppose $\phi ^2 \neq  0$ or $\phi ^4 \neq  0$, so $\phi ^2 +
\phi ^4 \geq {\pi / 3}$, and again $\sum ^4_{i=1} \phi ^i \geq  \pi .$ 

We have then
$$
\s(u) \geq 7{\pi / 3} = 7\pi/3  + p(Z){\pi / 3}.
$$
\ed
This completes Subcase B6.1.
\item[Subcase B6.2] $p=3$. 

In this case $q=1$ and 
$$
\s(u) \geq  \left( 5{\pi / 3}\right)  + 
\left( 2{\pi / 3}\right) \left( p(Z) + p(Y)\right) .
$$

\noindent
If $p(Y) \neq  0$ then
$$
\s(u) \geq  \left( 7{\pi / 3}\right)  + \left( 2{\pi / 3}\right) p(Z)
\geq 2\pi+(\pi/3)p(Z)+\pi/3,
$$

\noindent
as required.

Assume then that $p(Y) = 0$.  If $p(Z) \geq  2$ then
$$
\s(u) \geq  \left( 7{\pi / 3}\right)  + \left( {\pi / 3}\right) p(Z)
\geq 2\pi+(\pi/3)p(Z)+\pi/3.
$$

\noindent
Next consider the case $p(Z)=1$.  In this case either we have two blue and one 
green $C$--vertices adjacent to $u$ or two green and one blue.
\bd
\item[Subcase B6.2.1] $p(Y)=0$, $p(Z)=1$ and $u$ is adjacent to two blue and 
one green $C$--vertices. 

\begin{figure}
\psfrag{x^1}{$x^1$}
\psfrag{'a^1_1}{$\al^1_1$}
\psfrag{v^1}{$v^1$}
\psfrag{'phi^1}{$\phi^{1}$}

\psfrag{'a^2_1}{$\al^2_1$}
\psfrag{v^2}{$v^2$}
\psfrag{'phi^2}{$\phi^{2}$}

\psfrag{'a^3_1}{$\al^3_1$}
\psfrag{v^3}{$v^3$}
\psfrag{'phi^3}{$\phi^{3}$}
\psfrag{pi/2}{$\frac{\pi}{2}$}
\psfrag{pi/3}{$\pi/3$}
\psfrag{u}{$u$}
\begin{center}
{ \includegraphics[scale=0.7]{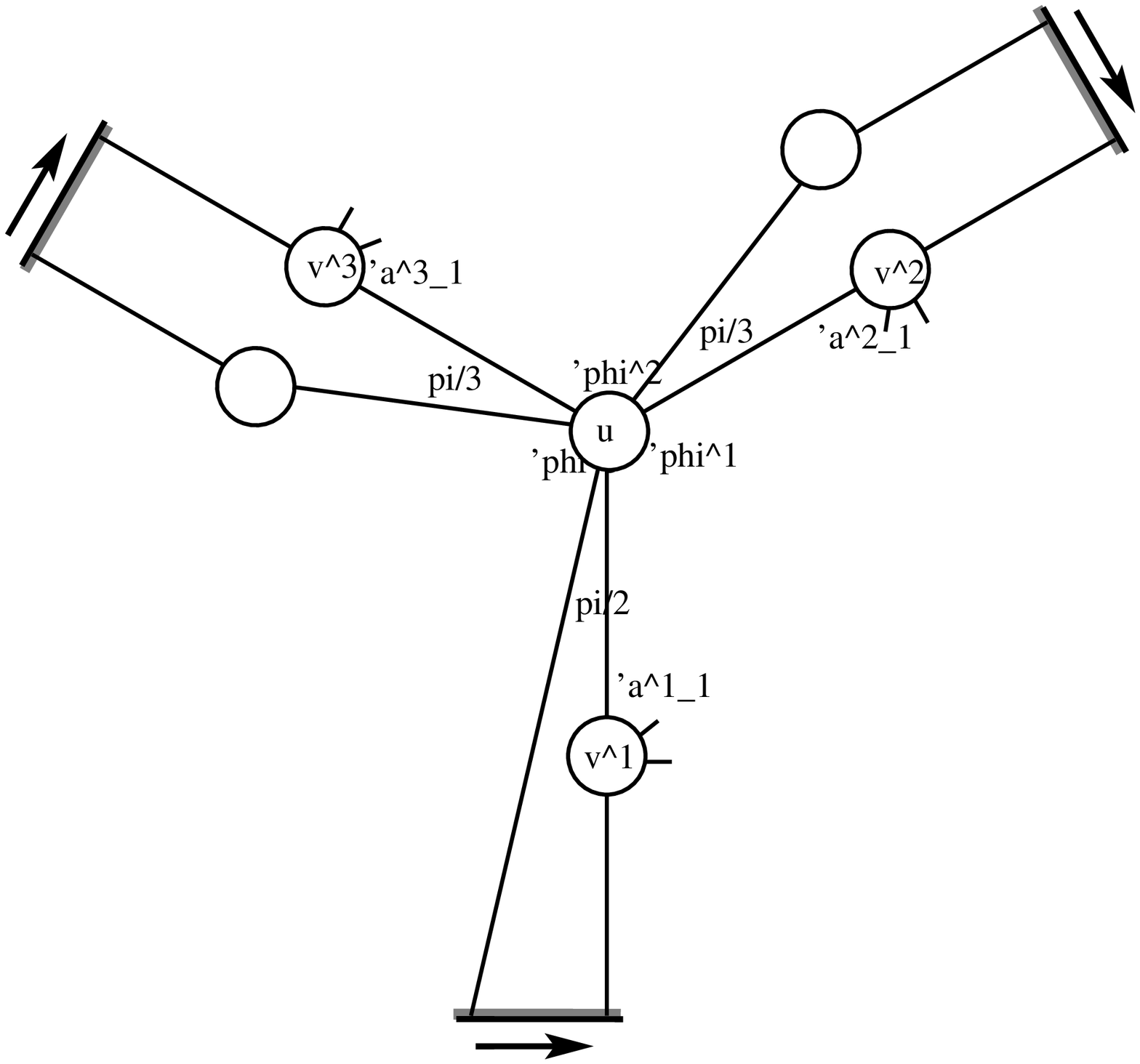} } 
\caption{
\label{p3Z11}}
\end{center}
\end{figure}

Suppose first that $g(Z)=1$ in which case $b(Z)=0$ and $b(X)=2$.  Then we 
have the configuration of Figure \ref{p3Z11}.  
\noindent
From Lemma \ref{wv} and Assumption \ref{A6}, 
$u$ must be incident to some 
class of arcs other than those shown in Figure \ref{p3Z11}.  It follows 
that $\sum ^3_{i=1}\phi ^i \geq  3{\pi / 2}$ so $\s(u) \geq  
8{\pi / 3}.$ 

\begin{figure}
\psfrag{x^1}{$x^1$}
\psfrag{'a^1_1}{$\al^1_1$}
\psfrag{v^1}{$v^1$}
\psfrag{'phi^1}{$\phi^{1}$}
\psfrag{x^2}{$x^2$}
\psfrag{'a^2_1}{$\al^2_1$}
\psfrag{v^2}{$v^2$}
\psfrag{'phi^2}{$\phi^{2}$}
\psfrag{x^3}{$x^3$}
\psfrag{'a^3_1}{$\al^3_1$}
\psfrag{v^3}{$v^3$}
\psfrag{'phi^3}{$\phi^{3}$}
\psfrag{'pi/3}{$\frac{\pi}{3}$}
\psfrag{pi/3}{$\pi/3$}
\psfrag{pi/2}{$\pi/2$}
\psfrag{u}{$u$}
\begin{center}
{ \includegraphics[scale=0.7]{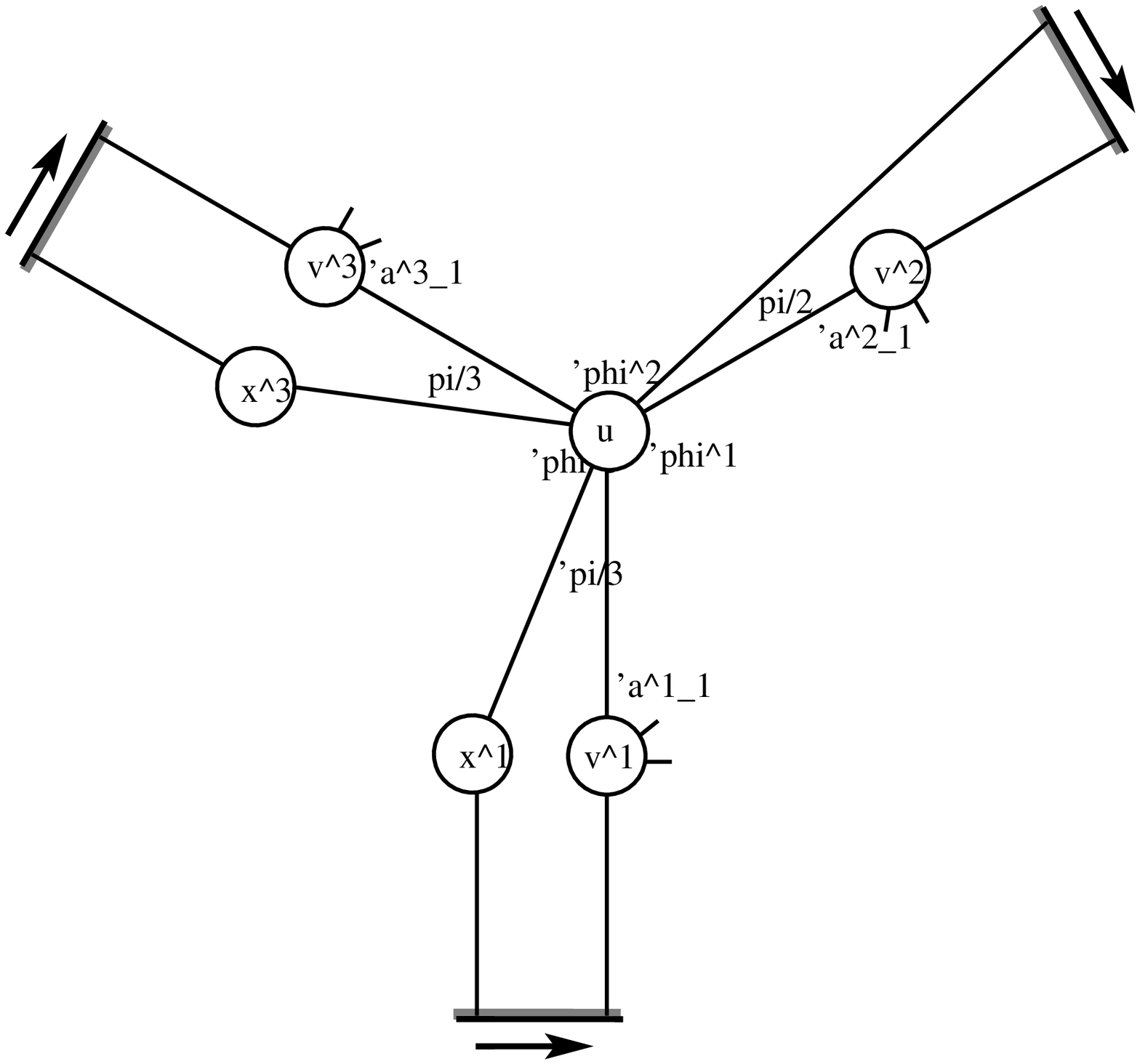} } 
\caption{
\label{p3Z12}}
\end{center}
\end{figure}

\begin{figure}
\psfrag{x^1}{$x^1$}
\psfrag{'a^1_1}{$\al^1_1$}
\psfrag{v^1}{$v^1$}
\psfrag{'phi^1}{$\phi^{1}$}
\psfrag{x^2}{$x^2$}
\psfrag{'a^2_1}{$\al^2_1$}
\psfrag{v^2}{$v^2$}
\psfrag{'phi^2}{$\phi^{2}$}
\psfrag{x^3}{$x^3$}
\psfrag{'a^3_1}{$\al^3_1$}
\psfrag{v^3}{$v^3$}
\psfrag{'phi^3}{$\phi^{3}$}
\psfrag{'pi/3}{$\frac{\pi}{3}$}
\psfrag{pi/3}{$\pi/3$}
\psfrag{pi/2}{$\pi/2$}
\psfrag{u}{$u$}
\begin{center}
{ \includegraphics[scale=0.7]{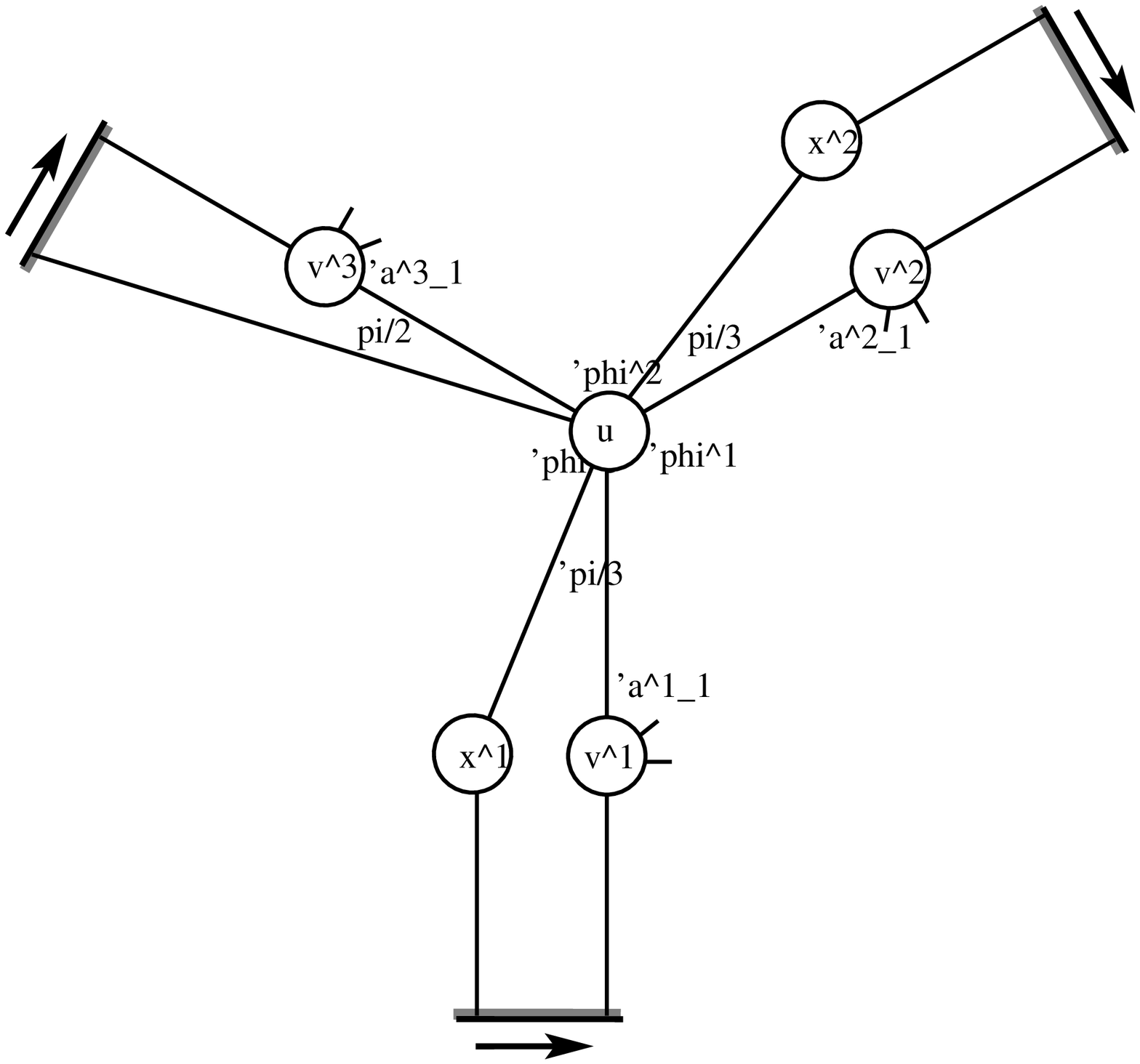} } 
\caption{
\label{p3Z13}}
\end{center}
\end{figure}

Now suppose $b(Z)=1$ so that $g(X)=b(X)=1$.  This gives rise to 
two possible configurations depending on the cyclic order of the blue types $X$
and $Z$ around $u$.  These two configurations correspond to Figures 
\ref{p3Z12} and \ref{p3Z13}. 
\noindent
\noindent
In the configuration of Figure \ref{p3Z12} it is possible that $\phi ^3=0$ 
if $x^1=x^3$.  In this case $u$ is incident to at least $2$ classes of
arcs not shown in Figure 
\ref{p3Z12}, using Lemma \ref{wv} and Assumption \ref{A6}, as before.
Therefore if 
$\phi ^3=0$ we have $\sum \phi ^i \geq  3{\pi / 2}$.  If $\phi ^3>0$ then 
as $\phi ^2 \geq  5{\pi / 6}$ we have again $\sum \phi ^i \geq  3{\pi / 2}$. 
Hence
$
\s(u) \geq  8{\pi / 3}.
$

In the configuration of Figure \ref{p3Z13} again $u$ must be incident to at 
least one more class than shown and so as before
$
\s(u) \geq  8{\pi / 3}.
$
Hence in Subcase B6.2.1 we have
$$
\s(u) \geq  8{\pi / 3}= 2\pi+p(Z)\pi/3+\pi/3.
$$
\item[Subcase B6.2.2] $p(Z)=1$ and $u$ is adjacent to two green 
and one blue $C$--vertices.

\begin{figure}
\psfrag{x^1}{$x^1$}
\psfrag{'a^1_1}{$\al^1_1$}
\psfrag{v^1}{$v^1$}
\psfrag{'phi^1}{$\phi^{1}$}
\psfrag{x^2}{$x^2$}
\psfrag{'a^2_1}{$\al^2_1$}
\psfrag{v^2}{$v^2$}
\psfrag{'phi^2}{$\phi^{2}$}
\psfrag{x^3}{$x^3$}
\psfrag{'a^3_1}{$\al^3_1$}
\psfrag{v^3}{$v^3$}
\psfrag{'phi^3}{$\phi^{3}$}
\psfrag{'pi/3}{$\frac{\pi}{3}$}
\psfrag{pi/3}{$\pi/3$}
\psfrag{pi/2}{$\pi/3$}
\psfrag{u}{$u$}
\begin{center}
{ \includegraphics[scale=0.7]{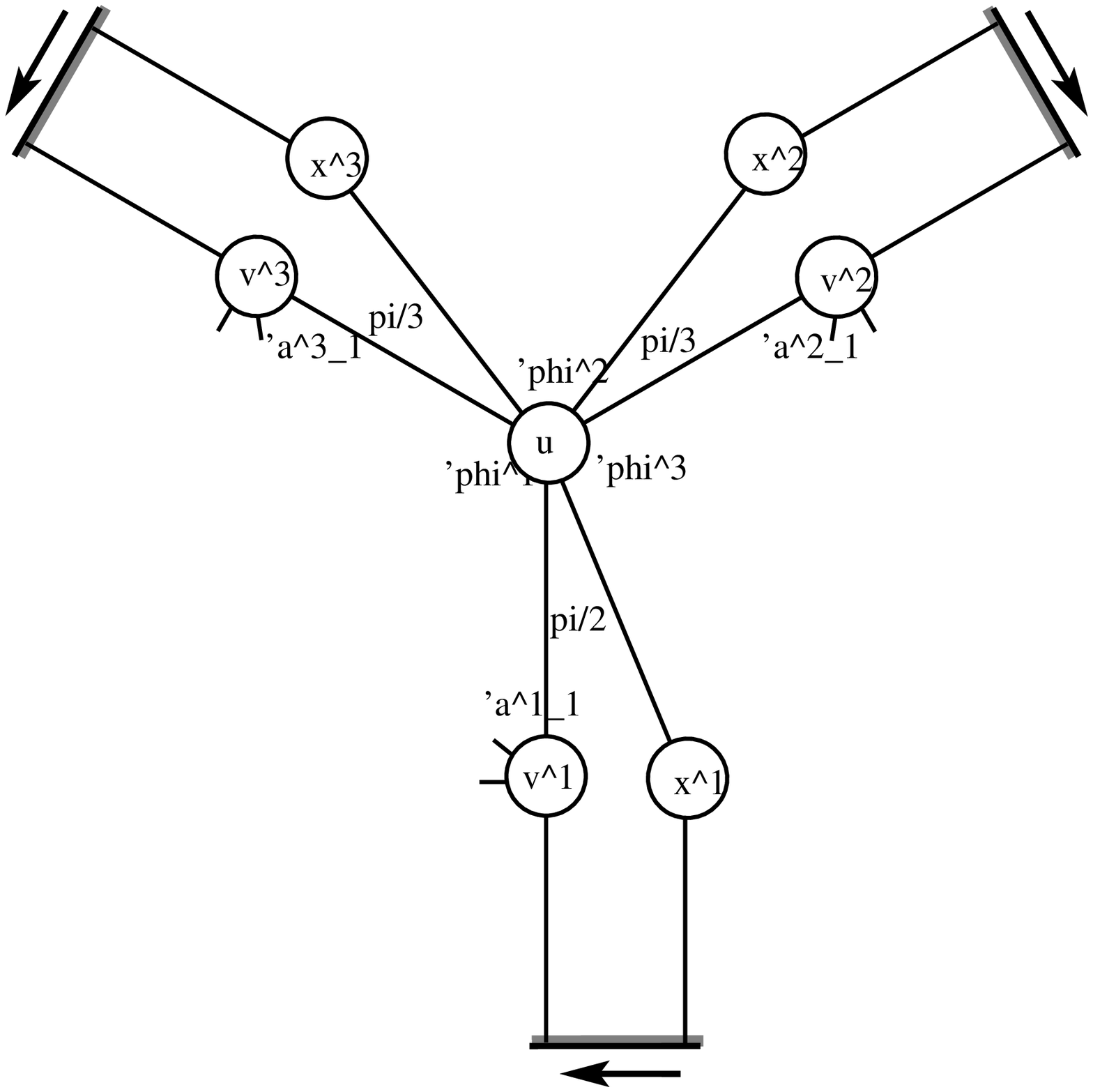} } 
\caption{
\label{p3Z0}}
\end{center}
\end{figure}

\noindent
This subcase follows by symmetry from Subcase B6.2.1.
\item[Subcase B6.2.3]  
The remaining Subcase of Subcase B6.2 is that in which
$p(Y)=p(Z)=0$. We consider only the case in 
which $u$ is adjacent to 
two blue and one green $C$--vertices.  The case in which $u$ is adjacent to 
one blue and two green $C$--vertices follows by symmetry.  We have 
$b_1(X)=2$ and $g_1(X)=1$ giving the configuration of Figure 
\ref{p3Z0}.  
\noindent
If $\phi ^2=0$ (that is $x^2=x^3)$ then 
$u$ must be incident to at least $2$ classes of arc which are
not shown, so 
$\sum \phi ^i \geq  4{\pi / 3}$.  
If
$\phi ^2 > 0$ then $\phi ^2 \geq  {\pi / 3}$.  As $u$ is incident to at least
$7$ classes of arcs 
$\sum \phi ^i \geq  4{\pi / 3}$ again.
Hence, when
$p(Z)=0$, we have $\s(u) \geq  7{\pi / 3}=2\pi+p(Z)\pi/3+\pi/3$.
\ed
This completes Subcase B6.2. 
\item[Subcase B6.3: $p=2$.] ~~\\
In this case $q = 1$ and
$$
\s(u) \geq  \pi  + 2{\pi / 3}\left( p(Z) + p(Y)\right),
$$

\noindent
with $p(Z) = 0$, 1 or 2. We shall show that  
\begin{equation}\label{su2}
\s(u) \geq  2\pi  +  p(Z)(\pi /3) + \pi/3\geq 2\pi  +  p(Z)(\pi /3) +
(p+1)\pi/9,  
\end{equation}
as required.

\begin{figure}
\psfrag{'i^1}{$\io^1$}
\psfrag{'m^1}{$\mu^1$}
\psfrag{v^1}{$v^1$}
\psfrag{'phi^1}{$\phi^{1}$}
\psfrag{'i^2}{$\io^2$}
\psfrag{'m^2}{$\mu^2$}
\psfrag{v^2}{$v^2$}
\psfrag{'phi^2}{$\phi^{2}$}
\psfrag{pi/2}{$\pi/2$}
\psfrag{u}{$u$}
\begin{center}
{ \includegraphics[scale=0.7]{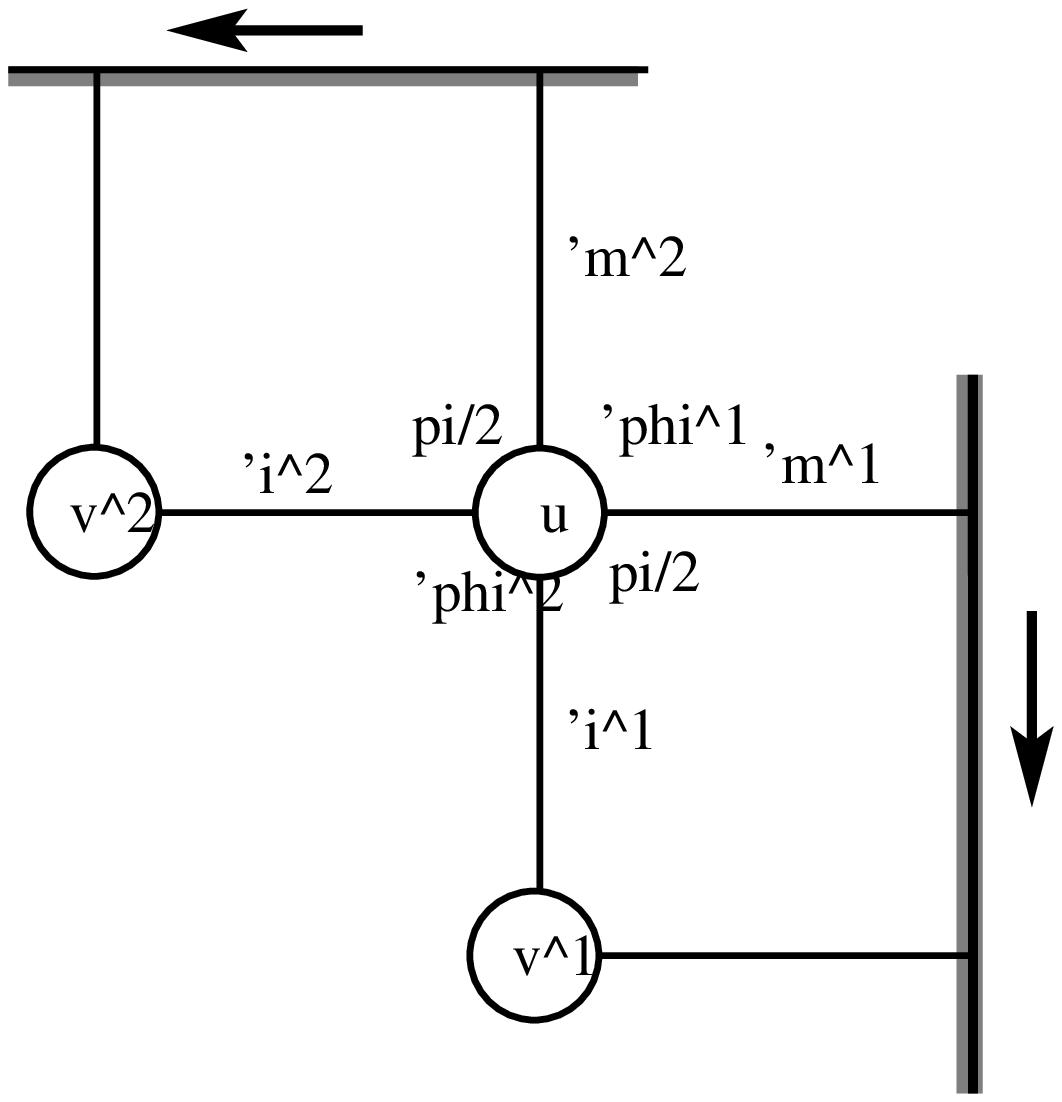} } 
\caption{
\label{p2Z2}}
\end{center}
\end{figure}

First suppose $p(Z) = 2$.  Then (up to reflection) 
we have the configuration of Figure \ref{p2Z2}.  
\noindent
As $|\mu ^i| 
\leq  l $, for $i=1,2$, it follows that there must be additional 
classes of arcs 
incident at $u$.  If $u$ is incident to no boundary classes of arcs other than 
$\mu ^1$ and $\mu ^2$ then $u$ is incident to at least $7$ classes of arcs.  If 
$u$ 
is incident to 7 or more classes of arcs including exactly 2 boundary classes 
then either $\phi ^2 \geq  {\pi / 3}$ and $\phi ^1 \geq  5{\pi / 3}$ or
$\phi ^2 \geq  2{\pi / 3}$ and $\phi ^1 \geq  4{\pi / 3}$ or $\phi ^2 
\geq  \pi $ and $\phi ^1 \geq  \pi $.  Hence $\phi ^1 + \phi ^2 \geq  2\pi $.  On the other 
hand if $u$ is incident to a boundary class of arcs other than $\mu ^1$ and 
$\mu ^2$ then either $\phi ^2 \geq  {\pi / 3}$ and $\phi ^1 \geq  2\pi $ or
$\phi ^2 \geq  \pi $ and $\phi ^1 \geq  \pi $.  Therefore, in  Figure 
\ref{p2Z2}, $\phi ^1 + \phi ^2 \geq  2\pi$ so $\s(u) \geq  3\pi
=2\pi+p(Z)\pi /3+\pi /3$. 

Now suppose $p=2$ and $p(Z)=1$.  Assume that $v^1$ is blue and $v^2$ is green. 
If $b(Y) = g(Z) = 1$ or $b(Z) = g(Y) = 1$ then the argument above, for 
the case $p=p(Z)=2$, can be used to show that $\s(u) \geq 2\pi+p(Z)\pi /3+2\pi /3$.  Hence we may assume that either
$b(X) = g(Z) = 1$ or $b(Z) = g(X) = 1.$
\bd
\item[Subcase B6.3.1] $p = 2$, $p(Z) = 1$ and $b(X) = g(Z) = 1.$

\begin{figure}
\psfrag{'i^1}{$\io^1$}
\psfrag{'m^1}{$\mu^1$}
\psfrag{v^1}{$v^1$}
\psfrag{'phi^1}{$\phi^{1}$}
\psfrag{'i^2}{$\io^2$}
\psfrag{'m^2}{$\mu^2$}
\psfrag{v^2}{$v^2$}
\psfrag{'phi^2}{$\phi^{2}$}
\psfrag{pi/2}{$\pi/2$}
\psfrag{u}{$u$}
\psfrag{'a^1_1}{$a^1_1$}
\begin{center}
{ \includegraphics[scale=0.7]{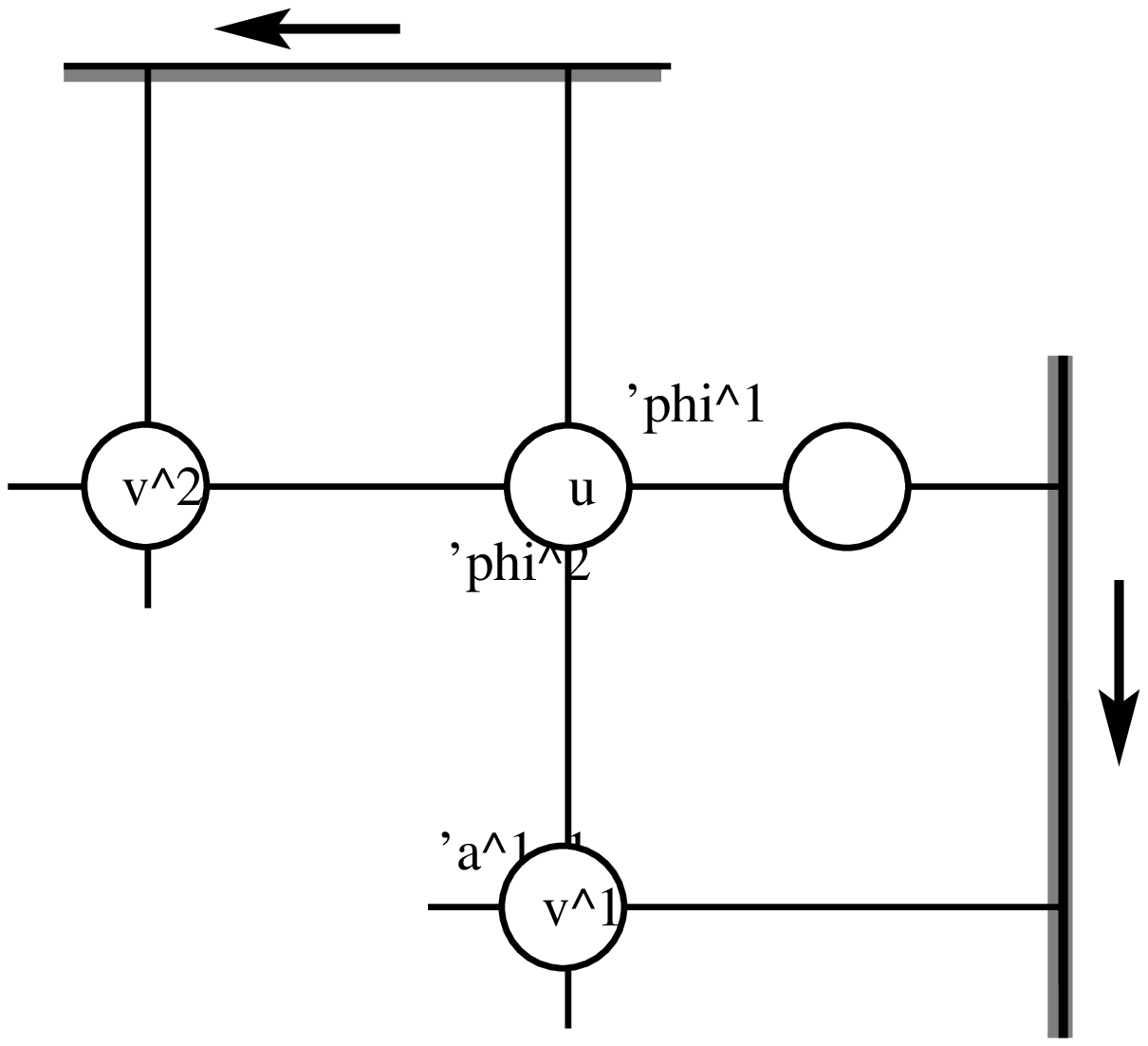} } 
\caption{
\label{p2Z1}}
\end{center}
\end{figure}

\noindent
We have the configuration of Figure \ref{p2Z1}.  
\noindent
Again $u$ must 
be incident to at least $3$ further non--boundary classes of arc or 1 further 
boundary class.  

If $u$ is incident to three or more further classes of arcs then $\phi ^1 \geq 
{\pi / 2}$ and $\phi ^2 \geq  4{\pi / 3}$ or $\phi ^1 \geq  
5{\pi / 6}$ and $\phi ^2 \geq  \pi $ or $\phi ^1 \geq  7{\pi / 6}$ and 
$\phi ^2 \geq  2{\pi / 3}$ or $\phi ^1 \geq  3{\pi / 2}$ and $\phi ^2 
\geq  {\pi / 3}$.  Thus $\phi ^1 + \phi ^2 \geq  11{\pi / 6}.$ 

If $u$ is incident to more than one boundary class then either $\phi ^1 \geq  
{\pi / 2}$ and $\phi ^2 \geq  4{\pi / 3}$ (since $\alpha ^1_1 = 
{\pi / 3})$ or $\phi ^1 \geq  3{\pi / 2}$ and $\phi ^2 \geq  
{\pi / 3}$.  Thus $\phi ^1 + \phi ^2 \geq  11{\pi / 6}$ in Subcase B6.3.1.  
Therefore $\s(u) > 8{\pi / 3}\geq 2\pi +p(Z)\pi /3 +\pi /3$.

\noindent
\item[Subcase B6.3.2] $p = 2$, $p(Z) = 1$ and $b(Z) = g(X) = 1.$

\noindent
This follows by symmetry from Subcase B6.3.1.
\ed

\begin{figure}
\psfrag{x^1}{$x^1$}
\psfrag{'i^1}{$\io^1$}
\psfrag{'m^1}{$\mu^1$}
\psfrag{v^1}{$v^1$}
\psfrag{'phi^1}{$\phi^{1}$}
\psfrag{x^2}{$x^2$}
\psfrag{'i^2}{$\io^2$}
\psfrag{'m^2}{$\mu^2$}
\psfrag{v^2}{$v^2$}
\psfrag{'phi^2}{$\phi^{2}$}
\psfrag{pi/2}{$\pi/2$}
\psfrag{u}{$u$}
\psfrag{'a^1_1}{$a^1_1$}
\psfrag{'a^1_2}{$a^1_2$}
\begin{center}
{ \includegraphics[scale=0.7]{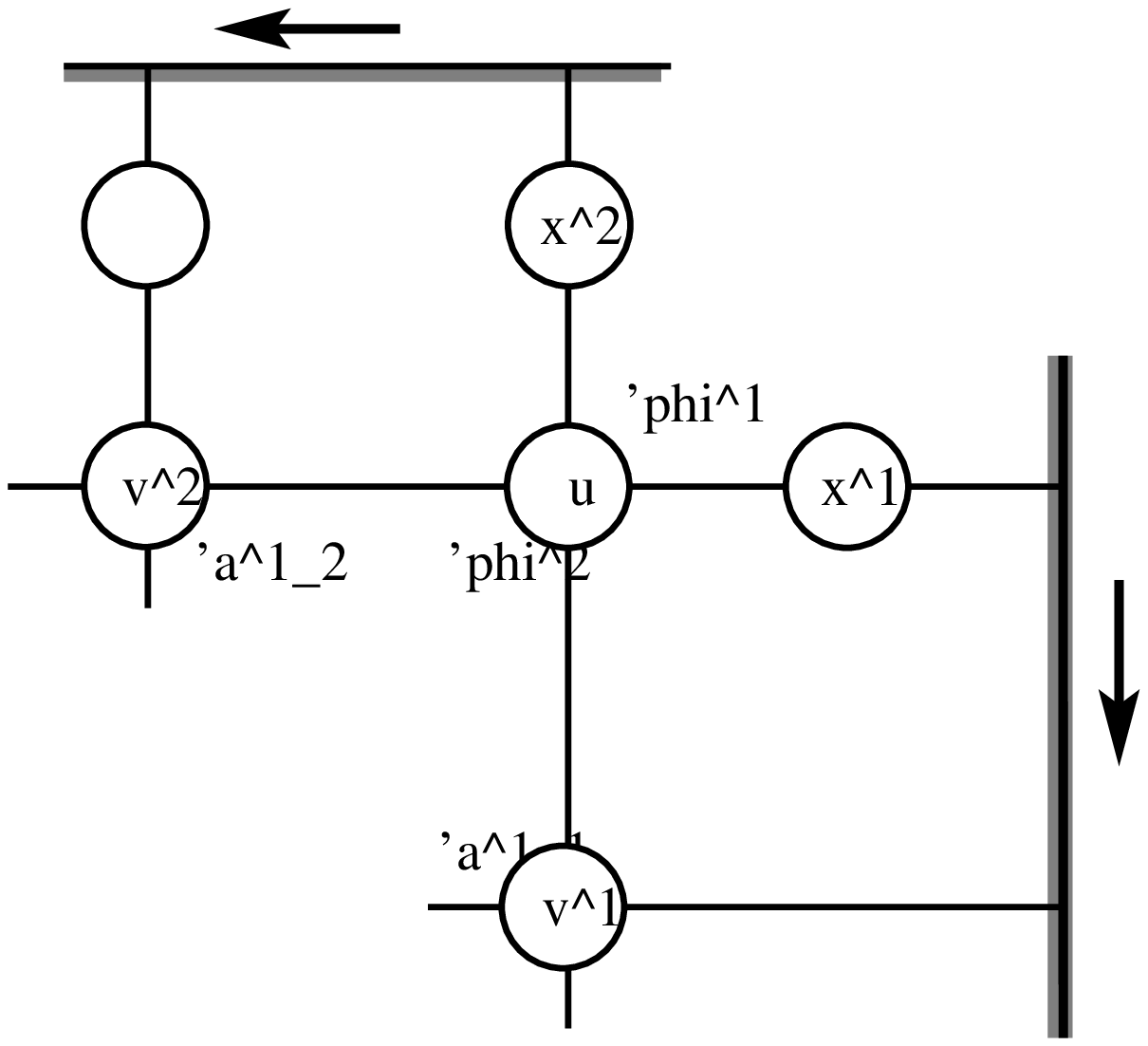} } 
\caption{
\label{p2Z0}}
\end{center}
\end{figure}

Finally we consider the case $p=2$ and $p(Z) = 0$.  
Assume that
$v^1$ is blue and $v^2$ is green.  Then in the cases 
$g(Y) = b(Y) = 1$, $g(Y) = b(X) = 1$ and 
$g(X) = b(Y) = 1$ the arguments used above in the cases $p(Z) = 1$ 
and 2 go through to give (\ref{su2}). Suppose then 
that $g(X) = b(X) = 1$, giving the configuration of Figure 
\ref{p2Z0}.  
\noindent
Note that $\phi ^1 = 0$ if and only if $x^1 = x^2.$ 

If $u$ is incident to more than 6 classes of arc then $\phi ^1 + 
\phi ^2 \geq  5{\pi / 3}$.  If $u$ is incident to a boundary class then 
either  $\phi ^2 = {\pi / 3}$ and $\phi ^1 = \pi $
or  $\phi ^2 =  {\pi / 3}$ and $\phi ^1 \geq  4\pi /3$
or  $\phi ^2 \geq  2{\pi / 3}$ and $\phi ^1 \geq  \pi $
or $\phi ^2 \geq  5{\pi / 3}$ and $\phi ^1 \geq  0$.  
In all except the first case $\phi ^1 + \phi ^2 \geq  
5{\pi / 3}$, so $\s(u) \geq  7{\pi / 3}\geq 2\pi +p(Z)\pi /3 +\pi /3$.
Also the first case, with $\phi ^2 = 
{\pi / 3}$ and $\phi ^1 = \pi $, gives (up to reflection) 
configuration $AA$ of Figure  \ref{AA}, in which the $v^i$ are not
$C$--vertices as angles have previously been adjusted, so this does not occur. 
\ed
This completes consideration of Subcase B6.3 and shows that if $u$ is 
incident to more than one vertex of type $Cj^\pm(v)$ then an angle adjustment may be
made to give all vertices involved negative curvature of at most $- 
{\pi / 15}$, hence completing Case B of the proof.

To complete the proof of the Proposition we consider 
configuration $D^{\prime } 6$, where $u$
is of type $C0(v)$, $v$ of type 
$CC^+(v)$ and $p=1$ as in case A2.4. 
Let $\bt^\prime$ be the boundary component of $\S$  to which $u$ is incident.
Suppose first that $\bt^\prime$ is oriented, in Figure \ref{D'6}, from left to
right. Assume first that
$x$ occurs as vertex $u_1$ in only $1$ $C$--configuration (namely that
appearing in Figure \ref{D'6}). Then in Case A2 with $u$ of type $C0(v)$ either we
have adjusted angles in $\D_5$, increasing $\al_7$ and decreasing $\al_6$, so
that $\ka(x), \ka(u)\le -\pi/12$, or we found that $u$ was of type $Dj(a)$, with $j=0,1,2$ or $3$.
If we adjusted angles $\al_6$ and $\al_7$ in Case A2 then we make further angle adjustments,
increasing $\al_4$ by $\pi/18$ and decreasing $\al_5$ and $\al_{10}$ by $\pi/36$, 
so that $\ka(x),\ka(u),\ka(v)\le -\pi/18$. On the other hand if $u$ is of type $Dj(a)$ then, since 
$\bt^\prime$ is oriented from left to right, $x$ is of type $Dj(b)$ and the only possibility
is that $j=3$. This gives rise to configuration $D6$, with $v$ of type $D6(a)$ and $u$ of 
type $D6(b)$. 

\begin{figure}
\psfrag{'b}{$\bt$}
\psfrag{b'}{$\bt^\prime$}
\psfrag{'n}{$\nu$}
\psfrag{'m}{$\mu$}
\psfrag{a}{$a$}
\psfrag{y}{$d$}
\psfrag{f}{$f$}
\psfrag{b}{$b$}
\psfrag{z}{$z$}
\psfrag{c}{$c$}
\psfrag{d}{$d$}
\psfrag{v2}{$v_2$}
\psfrag{vq}{$v_q$}
\psfrag{'D1}{$\D_1$}
\psfrag{'D2}{$\D_2$}
\psfrag{'D3}{$\D_3$}
\psfrag{'D4}{$\D_4$}
\psfrag{'D5}{$\D_5$}
\psfrag{'D6}{$\D_6$}
\psfrag{a1}{$\al_1$}
\psfrag{a2}{$\al_2$}
\psfrag{a3}{$\al_3$}
\psfrag{a4}{$\al_4$}
\psfrag{a5}{$\al_5$}
\psfrag{a6}{$\al_6$}
\psfrag{a7}{$\al_7$}
\psfrag{a8}{$\al_8$}
\psfrag{a9}{$\al_9$}
\psfrag{a10}{$\al_{10}$}
\begin{center}
{ \includegraphics[scale=0.5]{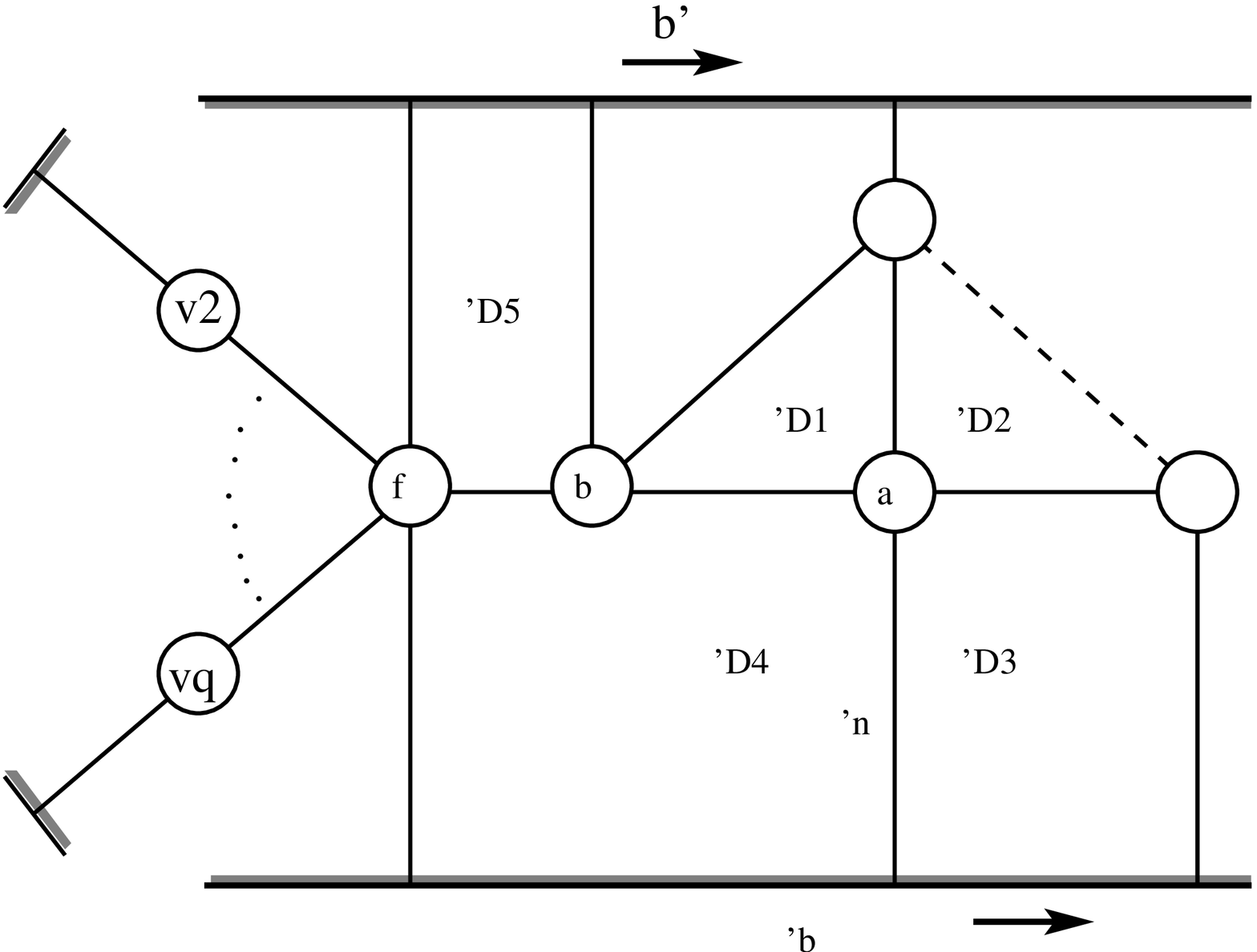} } 
\caption{}\label{D'6L}
\end{center}
\end{figure}

Now with $\bt^\prime$ still oriented from left to right assume that $x$ occurs as $u_1$ in $q>1$
$C$--configurations. If the $C$--vertices of these $C$--configurations are $u,v_2,\ldots ,v_q$ 
then we have the configuration of Figure \ref{D'6L}.
In Case B above angles on $u,x,v_2,\ldots, v_q$ were adjusted so that each of these vertices had 
curvature at most $-\pi/15$. Hence we may further adjust angles $\al_4$ and $\al_5$ so that 
$\ka(u),\ka(v)\le -\pi/30$. If say $v_i$ is also a vertex of type $D^\prime 6(u)$ and $x$ is the 
corresponding vertex of type $D^\prime 6(x)$ then similar adjustments are made to the angles of $v_i$
and the corresponding vertex of type $D^\prime 6(v)$.

Now suppose that $\bt^\prime$ is oriented from right to left in Figure \ref{D'6}. If $y$
occurs as vertex $u_1$ in only one $C$--configuration then, from Case A2 with $u$ of type
$C0(v)$, either  we may adjust angles $\al_8$ and $\al_9$ in Figure \ref{D'6} so that
$\ka(u),\ka(y)\le -\pi/12$ or $u$ is of type $Dj(a)$, with $j=0,1,2$ or $3$. The latter cannot occur
as $\bt^\prime$ is oriented from right to  left. Therefore we may make further adjustments to
angles $\al_4$, $\al_5$, $\al_8$ and $\al_9$ so that $\ka(u),\ka(v),\ka(y)\le -\pi/18$.
If on the other hand $y$
occurs as vertex $u_1$ in more than one $C$--configuration then, as in the case of $\bt^\prime$ 
oriented from left to right, we may adjust angles $\al_4$ and $\al_5$ so that $\ka(u),\ka(v)\le -\pi/30$. 
Furthermore we may make such adjustments simultaneously 
over all configurations of type $D^\prime 6$ in which $y$
occurs as the vertex of type $D^\prime 6(y)$.  
\noindent
This completes 
the proof of Proposition \ref{D-config}.
\section{Final angle adjustment}\label{angle_final}
We now make final adjustments to the assignment of angles to corners of 
vertices.  These adjustments are carried out on those configurations
of type $Dj$ where the angle of type $Dj(a)$ has curvature greater than $-\pi/30$:
that is, those arising in the proof of Proposition \ref{D-config}
in which angles were not adjusted.
In all such configurations the vertex $b$ is incident only to one $C$--vertex 
(namely the vertex $Dj(a)$).
In configurations $Dj$, with $0\leq j\leq 3$ we transfer angle of 
${\pi / 3}$ from the corner of $b$ in $\D _1$ to the corner of $a$ in 
$\D _1$.  
Note that since configuration $D6$ contains configuration $D3$ this results in a
transfer angle of ${\pi / 3}$ from the 
corner of $f$ in $\D _4$ of $D6$ to the corner of $b$ in $\D _4$ of $D6$.  
%
Having made all the above adjustments to the angle assignations we make the 
following definitions.  
\begin{defn}\label{vdiv}
Let 
$\cV_N = \left\lbrace v\in \cV:\ka(v) \leq - {\pi / 30}\right\rbrace $,
let $\cV_G=\cV_N\cup \cV_B$, let 
$$\cV_{k,B}=\{v\in\cV: v \textrm{ is distance at most } k \textrm{
  from } \cV_B\},$$ 
and let $\cV_Z = \cV\backslash \cV_G$ (where $\cV_B$ is defined in Definition \ref{h-vertex}).
\end{defn}

\begin{lemma} \label{zero-curve}
After final angle adjustment  $\ka(v) \leq 0$, for all
$v \in  \cV\backslash \cV_{1,B}$ and $\ka(v)\le \pi$, for all $v\in\cV_{1,B}$. Furthermore
if $v\in\cV_Z\backslash \cV_{1,B}$ then either
\be
\item\label{zero-curve-1}  $v$ is of type $B4(a)$ or 
\item\label{zero-curve-2} $v$ is of type $Dj(a)$, $j=0,\ldots, 9$ or
\item\label{zero-curve-3}  $v$ is of type $Dj(b)$, $j=0,\ldots, 9$ or 
\item\label{zero-curve-4}  $v$ is of type $Ij(a)$ or its mirror image, $j=0,\ldots,6$.
\ee
\end{lemma}
\medskip

\noindent{\em Proof.}
It follows from Theorem \ref{int-v} and  Corollary \ref{bd-c-curve} that 
before the angle adjustments made in the proof of Proposition \ref{D-config} either
$v$ is interior and $\ka(v)\le 0$ or $v$ is a boundary vertex and $\ka(v)\le \pi$,
for all $v\in\cV$.
Thus $\ka(v)\le \pi$, for all $v\in\cV_{1,B}$. From Theorem \ref{int-v}, 
if $v$ is an interior vertex not in $\cV_{1,B}$ then
$\ka(v)\le -\pi/3$, so $v\in\cV_N$. If $v$ is a boundary vertex not in $\cV_{1,B}$ and 
$v$ meets $\pd \S$ in $2$ or more classes of arcs then Corollary \ref{bd-c-curve} and 
Proposition \ref{prop-I} imply that either $v\in\cV_N$ or (\ref{zero-curve-4}) of the Lemma holds. 
If $v$ is a boundary vertex not in $\cV_{1,B}$ and 
$v$ meets $\pd \S$ in exactly $1$ class of arcs then, from Theorem \ref{c-config}, either
$v\in\cV_N$, (\ref{zero-curve-1}) of the Lemma holds or $v$ is a $C$--vertex.
Moreover $\ka(v)\le 0$, when (\ref{zero-curve-1}) or  (\ref{zero-curve-4}) of the Lemma holds.

The angle adjustments made in the proof of Proposition \ref{D-config} affect only
$C$--vertices and their adjacent vertices of type $u_1$ and result in curvature
being at most $-\pi/30$ for all vertices affected. Hence all of the above holds after
these adjustments have been made with the following two exceptions. Firstly, 
if $v$ is a boundary vertex not in
$\cV_{1,B}$ and $v$ meets $\pd \S$ in $2$ or more classes of arcs then
either $v\in\cV_N$; or (\ref{zero-curve-4}) of the Lemma holds; or $v$ is of type
$Dj(b)$, $0\le j\le 3$. Secondly, if $v$ is a boundary vertex not in
$\cV_{1,B}$ and $v$ meets $\pd \S$ in exactly one class of arcs then
either $v\in \cV_N$ or (\ref{zero-curve-1}), (\ref{zero-curve-2}) or (\ref{zero-curve-3}) of the Lemma holds.
 
Final angle adjustment affects only vertices in (\ref{zero-curve-2}) or (\ref{zero-curve-3}) of the Lemma,
so all of the above still holds and in addition $\ka(v)\le 0$ 
when (\ref{zero-curve-2}) or (\ref{zero-curve-3}) of the Lemma
holds. 

\begin{figure}
\psfrag{a}{$a$}
\psfrag{b}{$b$}
\psfrag{c}{$c$}
\psfrag{d}{$d$}
\psfrag{e}{$e$}
\psfrag{'D1}{$\D_1$}
\psfrag{'D2}{$\D_2$}
\psfrag{'D3}{$\D_3$}
\psfrag{'D4}{$\D_4$}
\psfrag{'D5}{$\D_5$}
\psfrag{'D6}{$\D_6$}
\psfrag{'D7}{$\D_7$}
\psfrag{'D8}{$\D_8$}
\begin{center}
{ \includegraphics[scale=0.4]{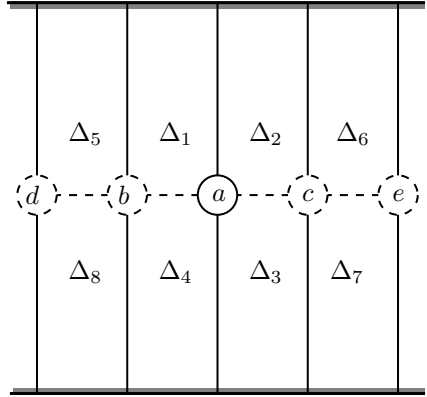} } 
\caption{Configuration $F1$
\label{F1}}
\end{center}
\end{figure}

\begin{figure}
\psfrag{a}{$a$}
\psfrag{b}{$b$}
\psfrag{c}{$c$}
\psfrag{d}{$d$}
\psfrag{e}{$e$}
\psfrag{'D1}{$\D_1$}
\psfrag{'D2}{$\D_2$}
\psfrag{'D3}{$\D_3$}
\psfrag{'D4}{$\D_4$}
\psfrag{'D5}{$\D_5$}
\psfrag{'D6}{$\D_6$}
\psfrag{'D7}{$\D_7$}
\psfrag{'D8}{$\D_8$}
\begin{center}
{ \includegraphics[scale=0.4]{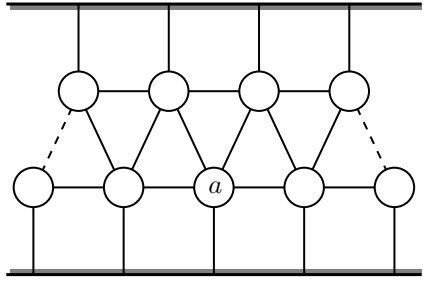} } 
\caption{Configuration $F2$
\label{F2}}
\end{center}
\end{figure}

\begin{defn}\label{F_vertex}
We say that a vertex $v$ is of type $Fj(a)$ if 
\be
\item $v$ occurs as vertex $a$ in a configuration $F1$ or $F2$ of Figures \ref{F1} or \ref{F2}, respectively, and 
\item the configuration $Fj$ in which $v$ occurs contains no $H^\Lm$--arcs. 
\ee
\end{defn}

\begin{theorem} \label{bad-vertex}
If $v\in \cV$ then either 
\be
\item\label{bad-vertex_1} $v$ is  of type $F1(a)$ or $F2(a)$ or 
\item\label{bad-vertex_2} $v$ is distance at most $3$ from $\cV_G$.
\ee
\end{theorem}
\noindent
\noindent

\noindent
{\em Proof.} We may assume that $v\in\cV_Z$ otherwise (\ref{bad-vertex_2}) holds trivially.
Let $v\in \cV_Z$ be distance more than 1 from $\cV_G$. From Theorem \ref{zero-curve}
it follows that 
$v$ and all its 
adjacent vertices are 
of type $Ij(a)$ or its mirror image, with  $j=1,\cdots , 6$, type $B4(a)$, 
type $Dj(a)$ or  type $Dj(b)$ with $0\leq j \leq 9$. 
Note that vertices of type $Dj(b)$, with $j=4$ or $5$ are also of type $B4(a)$ and that 
vertices of type $D6(b)$ are of type $D3(a)$. Thus once we've shown that a particular
vertex is not of type $B4(a)$ or $D3(a)$ then this vertex cannot be of type $Dj(b)$, $j=4,5,6$.

\begin{figure}
\psfrag{v}{$v$}
\psfrag{y}{$y$}
\psfrag{'a}{$\al$}
\psfrag{v_1}{$v_1$}
\psfrag{v_2}{$v_2$}
\psfrag{v_3}{$v_3$}
\psfrag{v_4}{$v_4$}
\psfrag{'D1}{$\D_1$}
\psfrag{'D2}{$\D_2$}
\psfrag{'D3}{$\D_3$}
\psfrag{'D4}{$\D_4$}
\psfrag{'D5}{$\D_5$}
\psfrag{'D6}{$\D_6$}
\psfrag{'D7}{$\D_7$}
\psfrag{'D8}{$\D_8$}
\begin{center}
  \mbox{
\subfigure[\label{BD1}]
{ \includegraphics[scale= 0.5,clip]{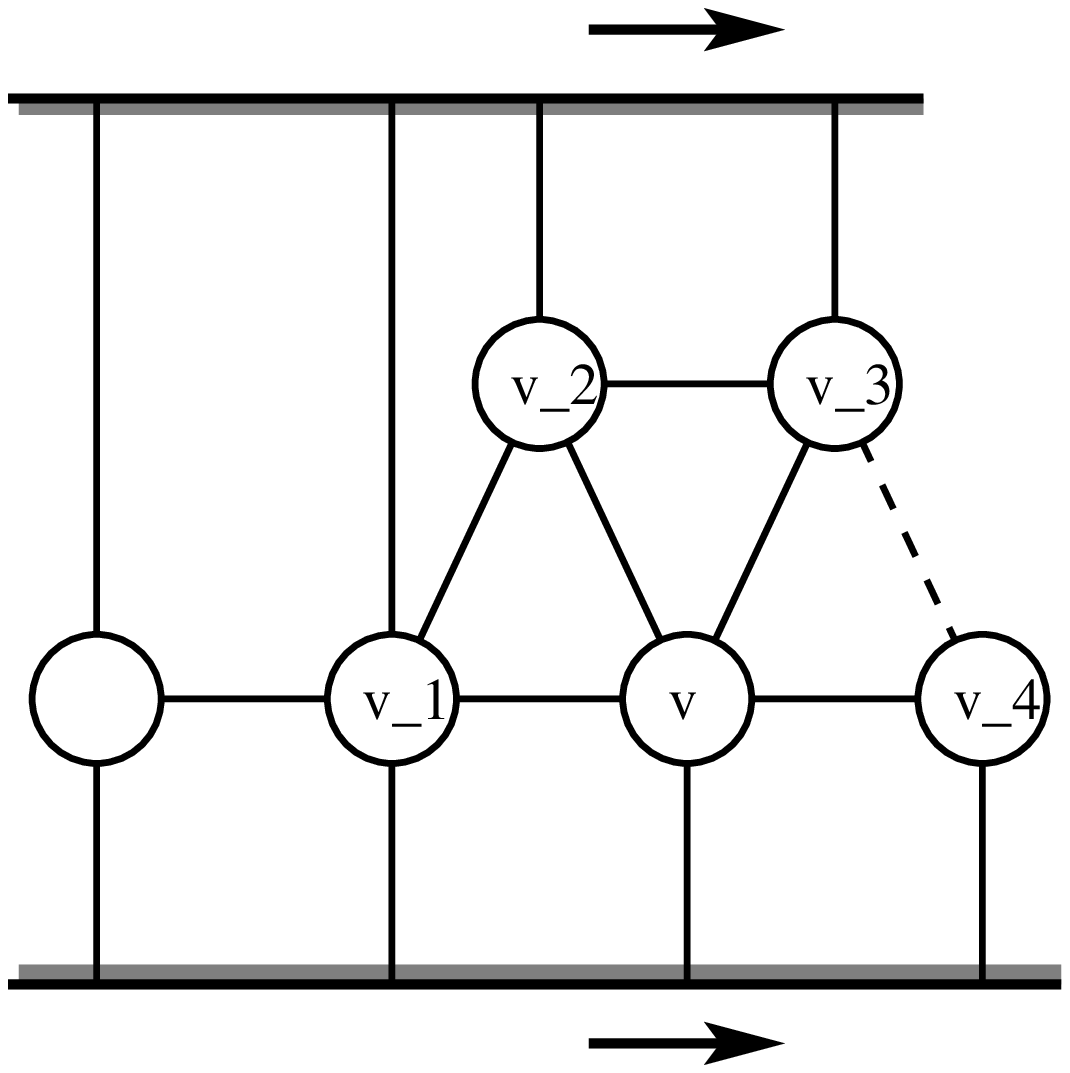} }\qquad
\subfigure[\label{BD2} ]  
{ \includegraphics[scale=0.5,clip]{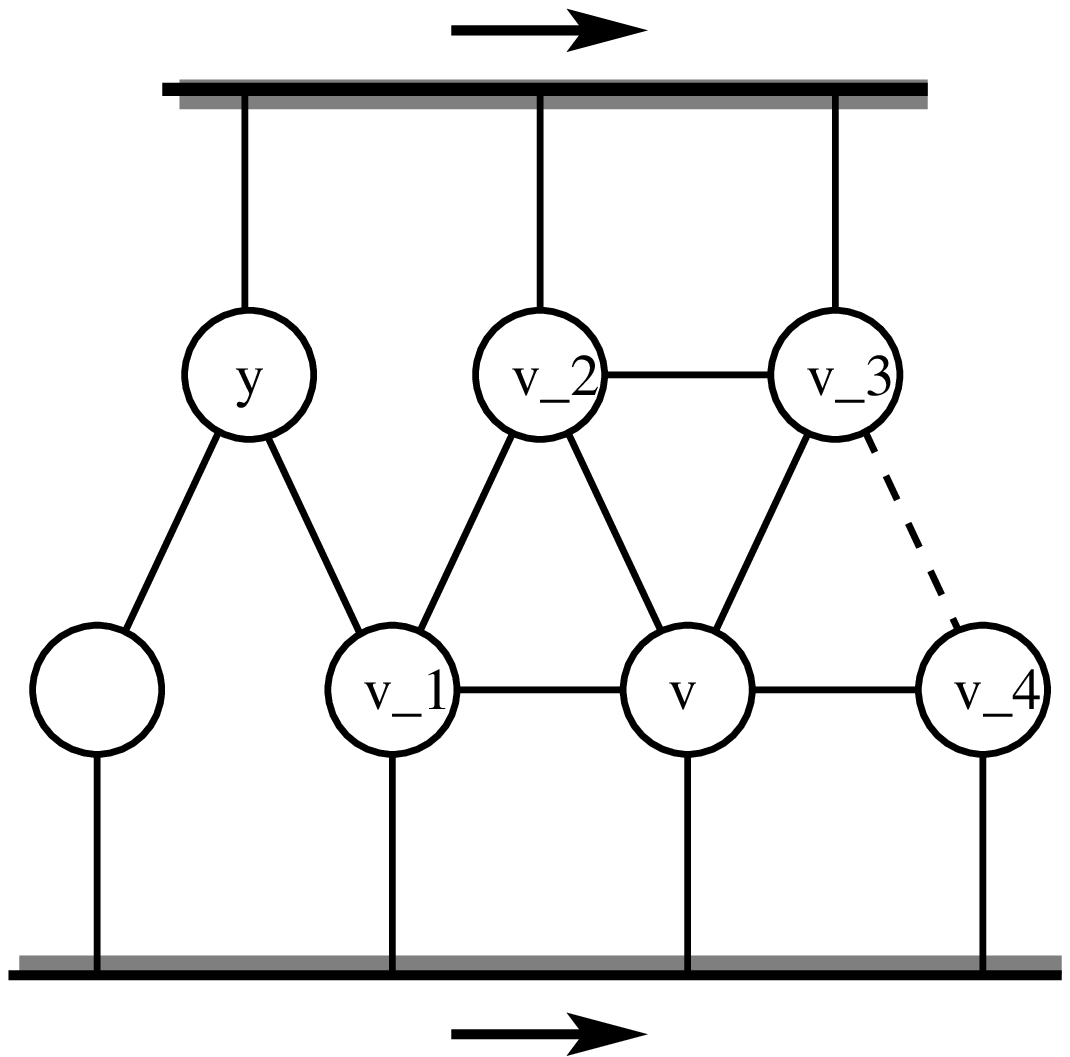} } 
}
\caption{}\label{BD1_BD2}
\end{center}
\end{figure}

\bd
\item[Case : $v$ of type $B4(a)$] ~\\
We show that in this case $v$ is either of type $F2(a)$ or is distance 
at most 2 from $\cV_G$. Assume then that $v$ is distance at least 3 from $\cV_G$.
Let the vertices incident to $v$ be $v_i$ of type 
$B4(a_i)$, $i=1,\cdots ,4$, as in Figure \ref{B4}.  Consider first 
$v_1$.  If $v_1$ where of type $Dj(a)$ then $v$ would necessarily be of type 
$Dj(c)$ so, from Theorem \ref{c-config}, the boundary class incident at $v$ would have width at 
most $l $.  As this boundary class has width more than $2l $, $v_1$ is not 
of type $Dj(a)$.  Furthermore $v_1$ is incident to a region $\D _1$ with 
$\rho (\D _1) = 3$, so $v_1$ cannot be of type $Ij(a)$.  Therefore $v_1$ 
can only be of type $Dj(b)$, $j=0,\ldots, 9$ or $B4(a)$.  Moreover 
$v_1$ cannot be of type $D3(b)$, since such vertices are 
not incident to any vertex of type $B4(a)$.  
Suppose that $v_1 = Dj(b)$ with $0\leq j \leq 2$.  
Since $v$ is neither of type $D0(f)$ or $Dj(a)$  
it follows that 
$v$ is of  type $Dj(d)$, forcing the configuration of 
Figure \ref{BD1}.  
Since we consider the case where $v_1$ is of type $B4(a)$ below and 
$v_1$ cannot be of type $D3(a)$, we need not consider the possibilities that $v_1$ is 
of type $Dj(b)$, with $j=4,5$ or $6$. 
If $v_1$ is of type $D7(b)$ then $v$ must be of type $D7(d)$ and 
we obtain the configuration of Figure \ref{BD2}.
If $v_1$ is of type $D8(b)$ then,  as $v$ is of type $B4(a)$, the only possibility
is that $v$ is 
the corresponding vertex of type $D8(d)$. However $v_1$ cannot then be of type $B4(a_1)$
so this cannot occur. If $v$ is of type 
$D9(b)$ then $v$ must be of type $D9(e)$. However in this case $v_1$ is a $C$--vertex and the 
class of boundary arcs incident at $v$ can have width at most $l$. As $v$ is of type $B4(a)$ this
class is in fact of width at least $2l$, so $v_1$ cannot be of type $D9(b)$.

\begin{figure}
\psfrag{v}{$v$}
\psfrag{y}{$y$}
\psfrag{'a}{$\al$}
\psfrag{'phi}{$\phi$}
\psfrag{v_1}{$v_1$}
\psfrag{v_2}{$v_2$}
\psfrag{v_3}{$v_3$}
\psfrag{v_4}{$v_4$}
\psfrag{'D1}{$\D_1$}
\psfrag{'D2}{$\D_2$}
\psfrag{'D3}{$\D_3$}
\psfrag{'D4}{$\D_4$}
\psfrag{'D5}{$\D_5$}
\psfrag{'D6}{$\D_6$}
\psfrag{'D7}{$\D_7$}
\psfrag{'D8}{$\D_8$}
\begin{center}
  \mbox{
\subfigure[\label{BB1}]
{ \includegraphics[scale= 0.5,clip]{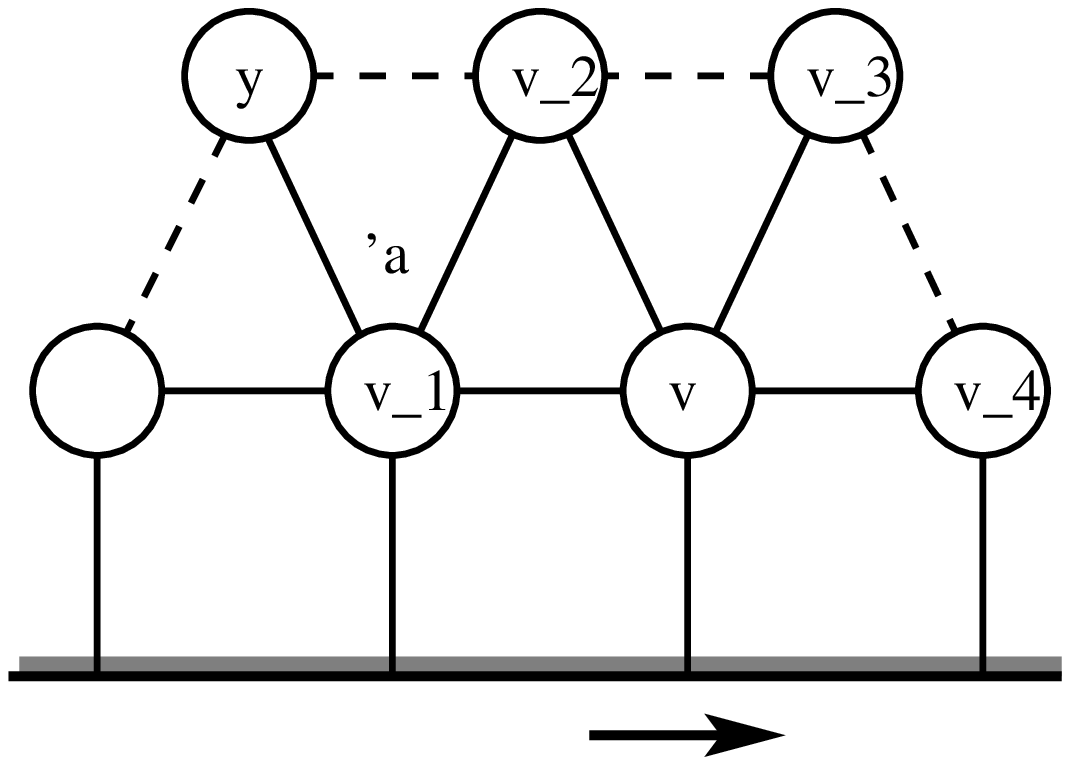} }\qquad
\subfigure[\label{BB2} ]  
{ \includegraphics[scale=0.5,clip]{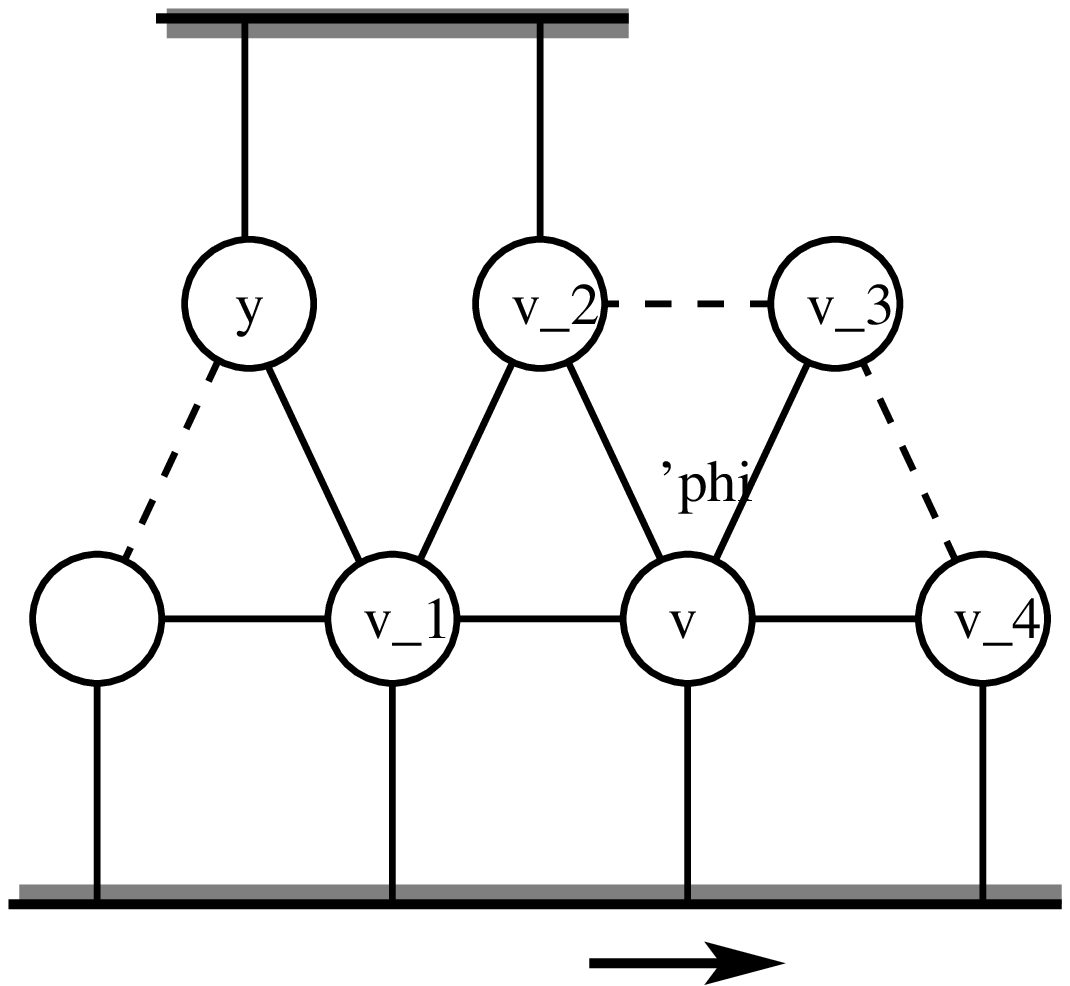} } 
}
\caption{}\label{BB1_BB2}
\end{center}
\end{figure}

\begin{figure}
\psfrag{v}{$v$}
\psfrag{y}{$y$}
\psfrag{v_1}{$v_1$}
\psfrag{v_2}{$v_2$}
\psfrag{v_3}{$v_3$}
\psfrag{v_4}{$v_4$}
\psfrag{'phi}{$\phi$}
\psfrag{'D1}{$\D_1$}
\psfrag{'D2}{$\D_2$}
\psfrag{'D3}{$\D_3$}
\psfrag{'D4}{$\D_4$}
\psfrag{'D5}{$\D_5$}
\psfrag{'D6}{$\D_6$}
\psfrag{'D7}{$\D_7$}
\psfrag{'D8}{$\D_8$}
\begin{center}
\mbox{
\subfigure[\label{BB5}]
{\includegraphics[scale=0.5,clip]{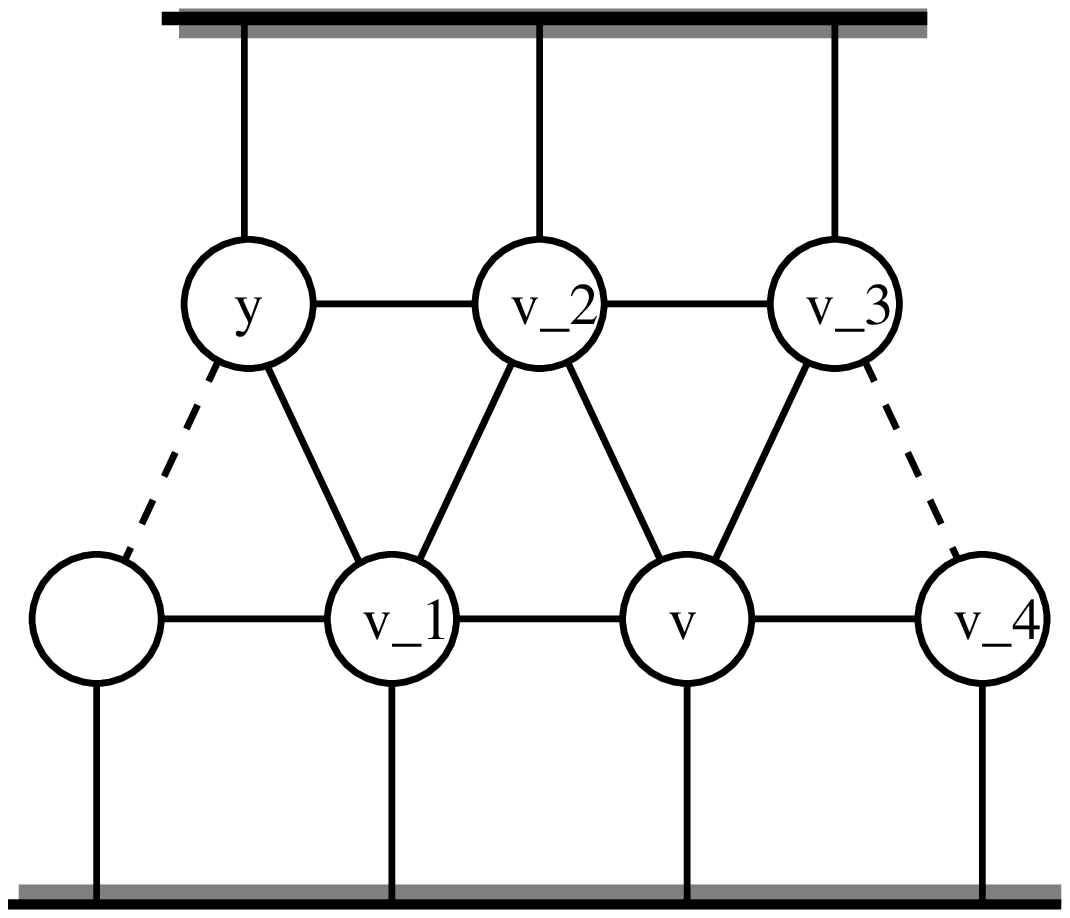} }\qquad
\subfigure[\label{DD1}]
{\includegraphics[scale=0.5,clip]{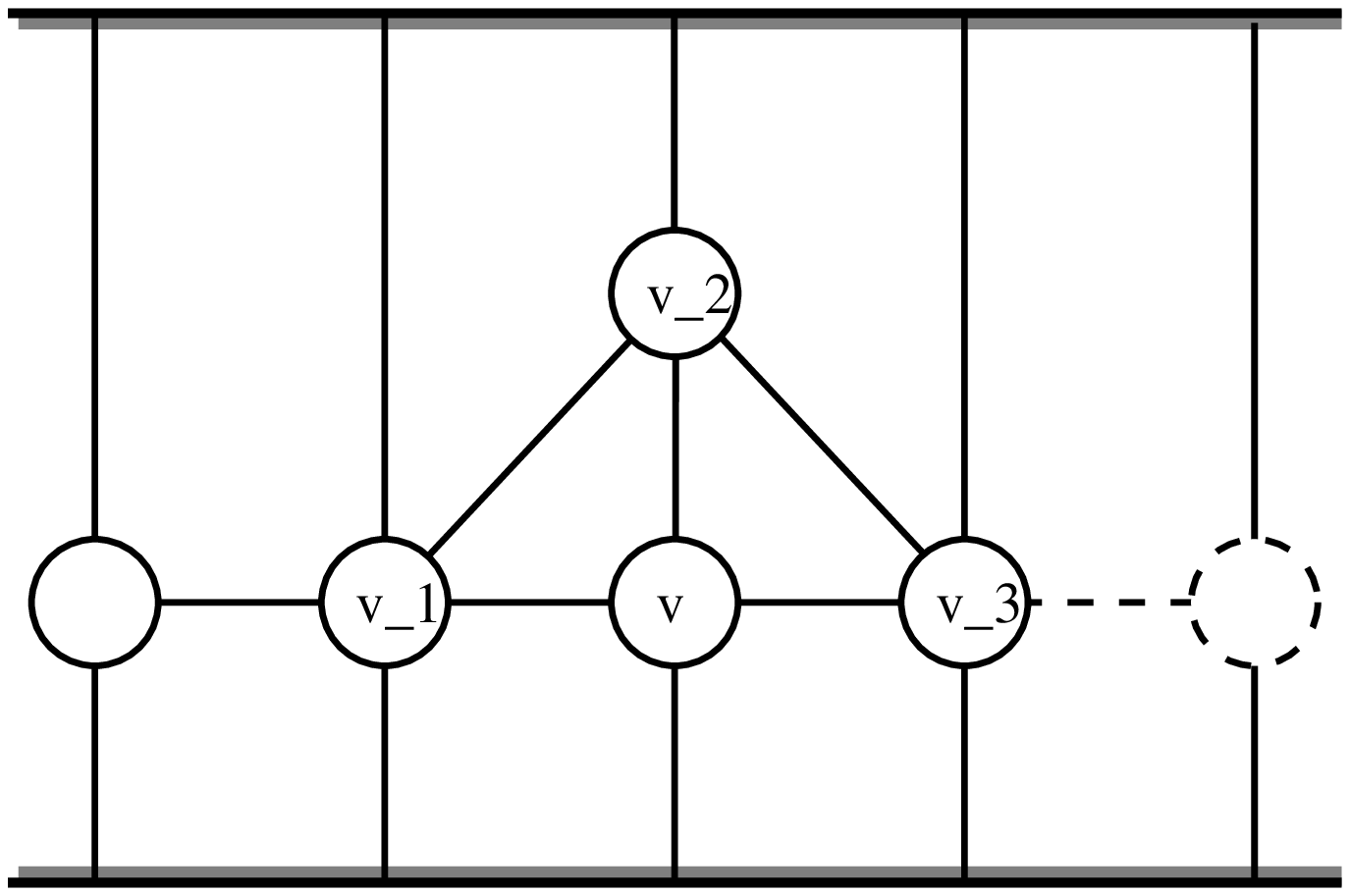} }
}
\caption{}\label{BB5_DD1}
\end{center}
\end{figure}

If
$v_1$ is of type $B4(a)$ then we have the configuration of
Figure \ref{BB1}, and 
we consider vertex $v_2$. 
 Clearly $v_2$ is not of type $D3(b)$.  If $v_2$ is 
not adjacent to $y$ then, as $v_1$ is of type $B4(a)$, $y$ and $v_2$ must both be
incident to the same boundary component, as shown in Figure 
\ref{BB2}.  
\noindent
\begin{figure}
\psfrag{v}{$v$}
\psfrag{y}{$y$}
\psfrag{v_1}{$v_1$}
\psfrag{v_2}{$v_2$}
\psfrag{v_3}{$v_3$}
\psfrag{v_4}{$v_4$}
\psfrag{'phi}{$\phi$}
\psfrag{'D1}{$\D_1$}
\psfrag{'D2}{$\D_2$}
\psfrag{'D3}{$\D_3$}
\psfrag{'D4}{$\D_4$}
\psfrag{'D5}{$\D_5$}
\psfrag{'D6}{$\D_6$}
\psfrag{'D7}{$\D_7$}
\psfrag{'D8}{$\D_8$}
\begin{center}
  \mbox{
\subfigure[\label{BB3}]
{ \includegraphics[scale=0.5,clip]{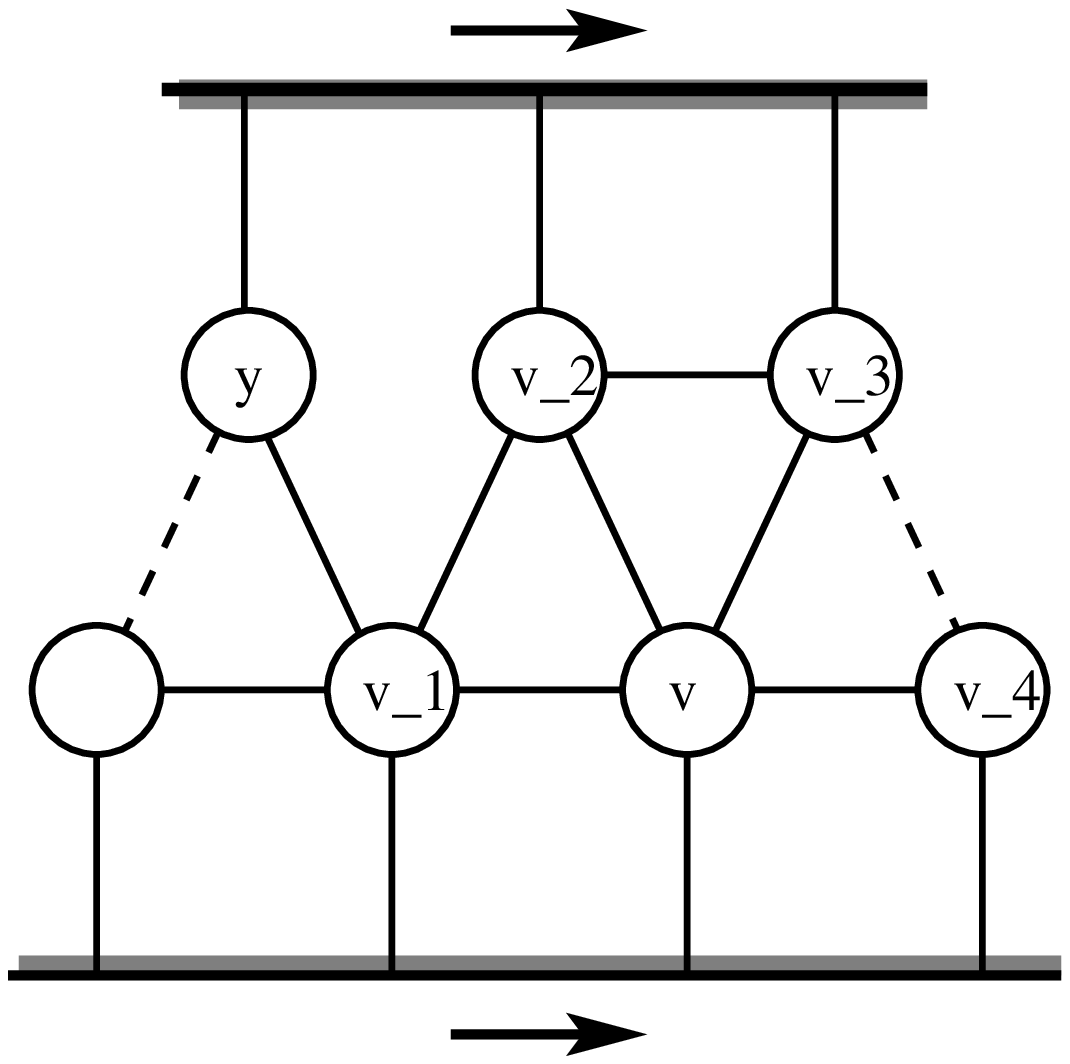} }\qquad
\subfigure[\label{BB4} ]  
{ \includegraphics[scale=0.5,clip]{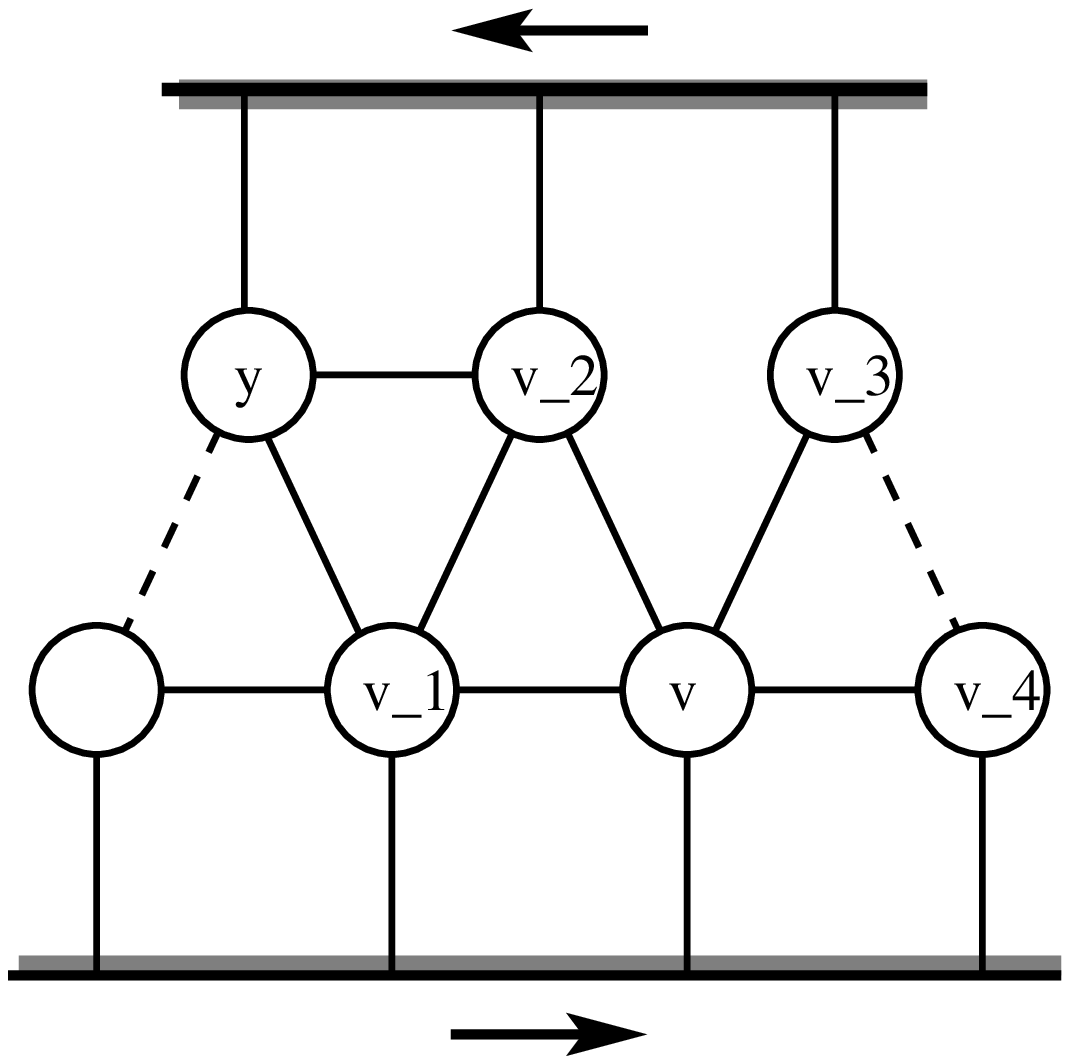} } 
}
\caption{}\label{BB3_BB4}
\end{center}
\end{figure}
As $v_2$ must be incident to some further arcs and 
the angle $\phi $ in Figure \ref{BB2} must be ${\pi / 3}$ it 
follows that $v_2$ and $v_3$ are adjacent.  By inspection of the regions 
incident at $v_2$ we see that $v_2$ cannot be of type $B4(a)$, $Ij(a)$
or $Dj(b)$.  Therefore 
$v_2$ is of type $Dj(a)$ and we obtain the configuration of Figure 
\ref{BB3}.  Suppose, on the other hand that, in Figure \ref{BB1}, 
$v_2$ is adjacent to $y$. Then 
if $v_2$ is of type $B4(a)$ we obtain the configuration of Figure 
\ref{BB5}.  If $v_2$ is of type $Dj(a)$ then, as $v_2$ is 
already adjacent to $v$, $v_1$ and $y$, it follows that $v_2$ must be of type 
$Dj(a)$ with $j\geq 4$ and $v_3$ must be incident to $\partial \S $,
as shown in 
Figure \ref{BB4}. It follows that $v_2$ must have type 
$D4(a)$ and $y$ must have type $D4(c)$. 
Inspection of regions incident to $v_2$, in Figure \ref{BB1}, show that $v_2$ cannot be
 adjacent to $y$ and of type $Dj(b)$,
with $0\leq j\leq 9$.
\noindent

We conclude that if $v$ is of type $B4(a)$ then we have one of the configurations 
of Figures \ref{BD1}, \ref{BD2}, \ref{BB3}, 
\ref{BB4} or \ref{BB5}.  In Figures 
\ref{BD1}, \ref{BD2} and \ref{BB3} the vertex $v_2$ is a $C$--vertex and
this forces the boundary 
class incident at $v_3$ in Figures \ref{BD1}, \ref{BD2}  and 
\ref{BB3} to have width at most $l$. Hence $v_3$ can only be of 
type $Dj(b)$, with $0\le j\le 3$. Consideration of its incident regions shows
that  $v_3$ cannot be any of these types so 
Figures 
\ref{BD1}, \ref{BD2} and \ref{BB3} do not occur.
In Figure \ref{BB4} as, by assumption,  $v$ is distance more than 
2 from $\cV_G$ we have 
$y\in \cV_Z$. However, as $y$ is of type $D4(c)$ the boundary class incident at $y$
in Figure \ref{BB4} has width at most $l$. Hence $y$ must be of type $Dj(b)$, with 
$0\le j\le 3$ or $7\le j\le 9$. As this is clearly not the case
the configuration of Figure \ref{BB4} 
does not arise. 

Hence  the only possible 
configuration is that of Figure 
\ref{BB5}.  Clearly $v_3$ is not then of type $Ij(a)$. 
If $v_3$ were of type $Dj(a)$ then the boundary class incident
at $v_2$ could have width at most $l$, whereas it has, in this case, 
width greater than $2l$. As above $v_3$ is not of type $Dj(b)$ so $v_3$
must be of type $B4(a)$. A similar argument shows that $v_4$ is of type
$B4(a)$.
This implies that $v$ is of type $F2(a)$.
\item[Case : $v$ of type $Dj(a)$, $0\le j\le 9$]~\\
Now consider the case where $v$ is of type $Dj(a)$, for $0\leq j\leq 9$, and is distance more than
$1$ from $\cV_G$.  
Let the vertices adjacent to $v$ 
be $v_1$, $v_2$ and $v_3$ of type $Dj(b)$, $Dj(d)$ and $Dj(c)$ respectively.  
Consideration of  regions to which $v_i$ is incident, shows that 
$v_i$ cannot be of type $Ik(a)$, 
for all $i$, $j$ and $k$.  Furthermore $v_3$ cannot be of type $B4(a)$ or a $C$--vertex, 
as $v_3$ meets $\partial \S $ in a class of width at most $l $.  Hence 
$v_3$ must be of type $Dj(b)$, with $0\leq j\leq 3$.
\bd
\item[Subcase : $v$ of type $D0(a)$] ~\\ 
Suppose $v_3$ is of type $D0(b)$.  Then consideration of regions incident to 
$v_3$ shows that $v_2$ and $v_3$ are adjacent and that $v_2$ is the vertex of 
type $D0(a)$ corresponding to the type $D0(b)$ vertex $v_3$.  We therefore have 
the configuration of Figure \ref{DD1}, including the dashed vertex and arc.  
However in this figure both boundary classes incident to 
$v_1$ are of width less than $l$, 
so  $v_1$ is incident to fewer than $ml $ arcs, a contradiction.  A 
similar argument  shows that $v_3$ cannot be of type $Dj(b)$ with $j =1$ or 2.  
If $v_3$ is of type $D3(b)$ then $v$ must be of type $D3(a)$, which 
is not the case.  
Therefore
$v$ cannot be of type $D0(a)$.
\item[Subcase: $v$ of type $Dj(a)$, $j = 1$ or $2$] ~\\
The argument of the previous case applies to show that $v$ cannot be of these 
types.
\item[Subcase: $v$ of type $D3(a)$] ~\\
If $v_3$ is of type 
$Dj(b)$ with $j\leq 2$ then $v_2$ must be of type $Dj(a)$, with $j\leq 2$, and 
$v_2$ and $v_3$ must be adjacent.  
Here $v_3$ is incident to fewer than $ml$ arcs, so this cannot happen.
If $v_3$ 
is of type $D3(b)$ then $v$ must be the corresponding vertex of type $D3(a)$, but
this is impossible, given the orientation of boundary components and the existence of 
a vertex of type $D3(f)$.
Therefore $v$ cannot be of type $D3(a)$.
\item[Subcase : $v$ of type $D4(a)$ or $D5(a)$] ~\\
In this case $v_1$ is of type $B4(a)$ and hence, using the Case: $v$
of type $B4(a)$,  
is a distance at most $2$ from $\cV_G$. Therefore $v$ is distance at most
$3$ from $\cV_G$.
\item[Subcase : $v$ of type $D6(a)$] ~\\
As $v_1$ is of type $D3(a)$ the distance of $v_1$
from $\cV_G$ is at most $1$ so that $v$ is distance at most $2$ from $\cV_G$.
\item[Subcase : $v$ of type $D7(a)$] ~\\
Since 
$v_1$ is of type $CC^{\pm}(v)$ it follows that  $v_2$
is of type $Dj(b)$ and $v_2$ cannot be a vertex of type $B4(a)$ or a $C$--vertex, so $0\le j\le 3$.  
It is easy to check that $v_2$ cannot be of 
these types, so $v$ cannot be of type $D7(a)$.

\item[Subcase : $v$ of type $D8(a)$] ~\\
If $v_3$ is of type $Dj(b)$, with $0\leq j\leq 2$, then $v_3$ is adjacent
to $v_2$ which must therefore be of type $B4(a)$. In this case $v$ is distance
at most $3$ from $\cV_G$. If $v_3$ is of type $D3(b)$ then $v$ is of type
$D3(a)$, which it cannot be. 
Therefore 
if $v$  is of type $D8(a)$ then $v$ is distance at most $3$ from $\cV_G$.
\item[Subcase : $v$ of type $D9(a)$]  ~\\
If $v_3$ is of type $Dj(b)$, with $j\leq 3$, then
$v_2$ is incident to fewer than $ml$ arcs.
Therefore $v$ is not of type $D9(a)$.
\ed
\item[Case : $v$ of type $Dj(b)$, with $0\leq j\leq 3$] ~\\
The vertex $v$ is adjacent to a vertex of type $Dj(a)$, with $0\leq j\leq 3$, 
which, from the above, is distance at most 1 from $\cV_G$.
Hence $v$ is distance at most $2$ from $\cV_G$. 
\item[Case : $v$ of type $Dj(b)$, with $4\leq j\leq 6$] ~\\
In this case $v$ is of type $B4(a)$ or $D3(a)$ so, 
from the above it follows that 
$v$ is distance at most $2$ from $\cV_G$.
\item[Case : $v$ of type $Dj(b)$, with $7\leq j\leq 9$] ~\\
If $j=7$ or $9$ then $Dj(a)$ is distance at most $1$ from $\cV_G$ so
$v$ is distance at most $2$ from $\cV_G$.
If $j=8$ then, from the case $v$ of type $D8(a)$, if $v$ is distance more than $2$ from
$\cV_G$ then the corresponding vertex of type $D8(d)$ is also of type $B4(a)$, so is distance
at most $2$ from $\cV_G$. Hence $v$ is distance at most $3$ from $\cV_G$.
\item[Case : $v$ of type $Ij(a)$] ~\\
Assume $v$ is distance more than 3 from $\cV_G$. 
As vertices of type $B4(a)$ and $Dj(a)$ are
neither adjacent nor paired to a 
vertex of type $Ij(a)$ it follows that $v$ can be adjacent 
or paired only to vertices of type $Ij(a)$ or  
$Dj(b)$ with $j = 0,\cdots , 9$.  If $v$ is adjacent or
paired to a vertex of 
type $Dj(b)$ then, by inspection, $0\le j\le 3$ and  from the above $v$ is distance at most 
$3$ from $\cV_G$, a contradiction. 
Hence all vertices to which $v$ is adjacent or
paired are of type $Ij(a)$. It follows that $v$ is of type $F1(a)$, Figure
\ref{F1}, and the result follows.
\ed
This completes the proof of the theorem.
\medskip

\noindent

We partition ${\cV}_Z$ into ${\cV}_F=\{v\in {\cV}_Z : v$ is of type $F1(a)$ or $F2(a)\,\}$
and ${\cV}_D={\cV}_Z\backslash {\cV}_F$.

\begin{lemma}   \label{D-number}
The number of vertices in ${\cV}_D$ is bounded above by
$$
|{\cV}_D| \leq  (ml )^3 |{\cV}_G|.
$$
\end{lemma}

\noindent
{\em Proof.}   If $v\in {\cV}_D$, then $v\notin {\cV}_G$ but the
distance of $v$ from $\cV_G$ is at most $3$ (Theorem \ref{bad-vertex}).
\medskip

\section{Isoperimetry}\label{isoperimetry}

Let $(\mbf h,\mbf n,\mbf t,\mbf p)$ be a positive 4--partition of 
$(h,n,t,p)\in \ZZ^4$, $\cL$
a consistent system of parameters and  $\mbf z$ a special element of
$(H^\Lm,\cL,\mbf n,\mbf h)$ and let $\cH$ be the homogeneous system
of parameters associated with $\mbf z$ and $\cL$. Then $\ovr{\mbf z}$ is a special
element of $(H^\Lm,\cH,\mbf n,\mbf h)$
Throughout this section  $\G$ is assumed to be a reduced, minimalistic picture over $G$
on a compact surface $\S$ 
of type $(\mbf n,\mbf t,\mbf p)$ with 
 boundary partition 
$\mbf b=(b_1,\ldots ,b_{W_1})$,
prime labels $\ovr{\mbf z}$ and 
labelled by $\hat\al(\ovr{\mbf z})$, for some solution $\al$ to $\cH$.
Let ${\cal D}$ denote the set of boundary components of $\S$. Then
\begin{equation*}
\sum _{\cal V}\ka (v) + \sum _{\cal R} \ka (\D ) + \sum _{\cal D} \ka (\bt )
\end{equation*}
\begin{equation*}
= \sum _{\cal V}\left( 2\pi -\s(v)\right)  + \sum _{\cal R} \left( \s(\D ) - 
t(\D )\pi  + 2\pi \chi (\D )\right) 
\end{equation*}
\begin{equation*}
+ \sum _{\cal D}  - \s(\bt )
\end{equation*}
\begin{equation*}
= \sum _{\cal V}2\pi  + \sum _{\cal R} \left( -t(\D )\pi  + 
2\pi \chi (\D )\right) 
\end{equation*}
\begin{equation*}
= 2\pi \left( |\cV| + \sum _{\cR} \left( - (1/ 2)t(\D ) + 
\chi (\D )\right) \right) 
\end{equation*}
\begin{equation*}
= 2\pi \left( |\cV| - |\cA| + \sum _{\cR} \chi (\D )\right)  = 2\pi \chi (\S ).
\end{equation*}
That is 
\begin{equation} \label{total-curve}
\sum _{v\in {\cal V}}\ka (v) + \sum _{\D \in {\cal R}} \ka (\D ) + \sum _{\bt \in {\cal D}}\ka (\bt )
= 2\pi \chi (\S )
\end{equation}

\begin{defn}
The $Z${\em --subspace} of $\G$ is the union of all maximal designated
corridors of $\G$.
\end{defn}

\begin{defn}\label{I-corridor-boundary}
Let 
$C=(D,\G_D,\mu)$ be a maximal designated simple $j$--corridor, where $j=0,1$ or $2$.
Assume 
that the image of $C$ under $\mu$ is a disk, let $c$ be the width of
$C$ and let $\G_\bt=\G_D(0)$ or $\G_D(c-1)$. If $\G_\bt$ is a picture of
type square--j then $\mu(\G_\bt)$ is called an $j${\em --corridor--boundary}.
Let $E$ be a  maximal designated $2$--corridor with 
foundation $C_0=(D_0,\G_{D_0},\mu_0)$ and extension $C_1=(D_1,\G_{D_1},\mu_1)$ and with
image a disk. Then $E$ has one compartment $\G_0$ which is not a compartment of $C_1$ and 
one compartment $\G_1$ which is not a compartment of $C_0$. The image $\mu(\G_i)$ is 
called 
a $2${\em --corridor--boundary}, for $i=1,2$. A $j$--corridor-boundary with $j=0,1$ or $2$ is
called a {\em corridor--boundary}.
\end{defn}

The image of every corridor is an annulus a M\"obius band or a disk. Every corridor, of width at least $2$,
which has image a disk has two corridor--boundaries or consists of a unique arc of type $I$.

\begin{prop}   \label{I0-bd} 
The number of $0$--corridor--boundaries is bounded above by 
$$3(ml(ml+1)^3|{\cal V}_G| + W_0 - 2\chi (\S ))+(6ml+10)W_1.$$
\end{prop}

\noindent
{\em Proof.}   Let  ${\cR}_1=\{\D\in {\cal R} : \D$ is not contained in the $Z$--subspace
of $\G$ $\}$.  For $\D \in {\cR}_1$ define 
$I(\D ) =$ number of $0$--corridor--boundaries meeting $\partial \D $.  Let ${\cR}_2 = 
\{\D \in {\cR}_1 : 0 < I(\D ) \leq  \rho (\D )\}$.
Then
\begin{equation} \label{st2}
\sum _{\D \in {\cR}_2}I(\D ) \leq  \sum _{\D \in {\cR}_2}\rho (\D )
.
\end{equation}

Let 
$$\cR_3=\{\D\in\cR_1 : I(\D)>\rho(\D) \textrm { and } \D \textrm{ is
  not collapsible}\}$$ and 
$$\cR_4=\{\D\in\cR_1 : I(\D)>\rho(\D) \textrm { and } \D \textrm{ is
  collapsible}\}.$$
Then $\D\in \cR_3\cup\cR_4$ implies that $\e(\D)\ge I(\D)>\rho(\D).$
Let $\D\in \cR_3$. If $\rho(\D)\ge 1$ then
$$
\ka (\D ) \leq  
\left( {\rho (\D )-\e(\D )}\right) \pi/ 2  \leq  
\left( {\rho (\D )-I(\D )/ 3}\right) \pi ,
$$
using the proof of Lemma \ref{r-curve}.
\noindent
If $\rho (\D )=0$, then $$\ka (\D ) \leq  (-\e(\D ) + 
2\chi (\D ))\pi .$$
\noindent
If $\rho (\D )=0$ and $\chi (\D )\leq 0$ then
$$
\ka (\D ) \leq  -\e (\D )\pi  \leq  -I(\D )\pi  < - 
{I(\D )\pi / 3}
$$
$$
= \left( {\rho (\D )-I(\D )/ 3}\right) \pi .
$$

\noindent
If $\rho (\D )=0$ and $\chi (\D )=1$, then 
$t(\D )=\e(\D )\geq 2$, since $\G $ has no trivial 
boundary labels.  
If $\e(\D )=2$ then $\D $ is collapsible so $\D\notin \cR_3$. 
Hence $\e(\D )\geq 3$ and
$$
\ka (\D ) \leq  (-\e(\D )+2)\pi  \leq  - 
{\e(\D )\pi / 3} \leq  
\left( {\rho (\D )-I(\D )/ 3}\right) \pi .
$$

\noindent
Therefore if $\D\in\cR_3$ then  
%
$$
\ka (\D ) \leq  \left( {\rho (\D )-I(\D )/ 3}\right) \pi .
$$
\noindent
As $\ka(\D)\le 0$, for all regions $\D$ (Lemma \ref{r-curve}), we have
$$
\sum _{\cal R} \ka (\D ) \leq  \sum _{\D \in {\cR}_3}\ka (\D ) \leq 
\sum _{\D \in {\cR}_3}\left( {\rho (\D )-I(\D )/ 3}\right) \pi .
$$

\noindent
Also, from 
Lemma \ref{zero-curve} and Lemma \ref{HLbound}
\[\sum _{\cal V}\ka (v) \leq  |{\cal V}_{1,B}|\pi\le 2(ml+1)W_1\pi ,\] 
as the number of vertices at distance at most $1$ from an $H^\Lm$--vertex
is at most $ml$.  Hence, using (\ref{total-curve}) and Lemma 
\ref{bd-curve},
$$
2\pi \chi (\S ) \leq  2(ml+1)W_1\pi  + 
\sum _{\D \in \cR_3}\left( {\rho (\D )-I(\D )/ 3}\right) \pi  + 
W_0\pi ,
$$
so
\begin{equation}\label{st3}
\sum _{\D\in \cR_3}I(\D)\leq  6(ml+1)W_1+\sum_{\D\in \cR_3}3\rho(\D)+3W_0-6\chi(\S)
\end{equation}
Now let $\D\in \cR_4$. 
Then $\D\in \cR_1$ implies that $\pd \D$ contains
an $H^\Lm$--arc so there are at most $2W_1$ regions in $\cR_4$. 
As $\D$ is collapsible we have $\e(\D)=t(\D)=2$ and so $I(\D)\le 2$. Therefore
\begin{equation} \label{st4}
\sum_{\D\in \cR_4}I(\D)\le 4W_1.
\end{equation}

Now consider a  vertex $v$ which meets $\pd \D$, 
for some region $\D\in \cR_2\cup \cR_3$. Then $v$ cannot be of 
type $F2(a)$ as no region incident to $a$ in $F2$ meets a $0$--boundary.
If $v$ is of type $F1(a)$ then $\D$ and $v$ both belong to the
same $j$--corridor, contrary to choice of $\D\in\cR_1$, unless $v$ is an $H^\Lm$--vertex (that
is $v\in\cV_B\subseteq \cV_G$). Hence, 
from Theorem \ref{bad-vertex},
$v$ is distance at most $3$ from $\cV_G$. Every vertex has degree $ml$, 
so if a vertex is paired to $p$ vertices and joined to $q$ vertices then $p+q\le ml$.
Hence the number of vertices which are distance at most 3 from $|\cV_G|$ is less than
or equal to $(ml+1)^3|\cV_G|$. A vertex meets at most $ml$ regions so 
\begin{equation}\label{st5}
\sum_{\D\in\cR_2\cup\cR_3} \rho(\D)\le ml(ml+1)^3|\cV_G|.
\end{equation}
From (\ref{st2}), (\ref{st3}) and (\ref{st4}), it follows that
$$\sum_{\D\in\cR_2\cup\cR_3\cup\cR_4} I(\D)\le
\sum_{\D\in\cR_2\cup\cR_3} 3\rho(\D) +6(ml+1)W_1 +3W_0-6\chi(\S) +4W_1.$$
As the number of $0$--corridor--boundaries is at most $\sum_{\D\in\cR_2\cup\cR_3\cup\cR_4} I(\D)$
the result follows from  (\ref{st5}).

\begin{lemma}\label{corridor-number}
There are fewer than 
$$
n+
3((ml+1)^4|\cV_G| + W_0 - 2\chi (\S ))+(6ml+10)W_1
$$
maximal designated corridors of $\G$.
\end{lemma}
\medskip

\noindent{\em Proof.}
Suppose $C$ is a maximal designated $j$--corridor and let $C_D$ be the
image of the $C$. If $C_D$ is an annulus or
M\"obius band then $C_D$ is a connected component of
$\S$. As $\S$ has $n$ boundary components there are
at most $n$  maximal designated corridors
of this form. Assume then that $C_D$ is not an annulus or
M\"obius band, so $C_D$ is a disk. If $j=0$ then $C$ has a
$0$--corridor--boundary. From Proposition \ref{I0-bd} there are at
most 
$$c_0=3(ml(ml+1)^3|{\cal V}_G| + W_0 - 2\chi (\S ))+(6ml+10)W_1$$
corridors with a $0$--corridor-boundary.
If $j\ge 1$ then either $C$ has 2
$0$--corridor boundaries or  a $j$--corridor--boundary. If $C$ has a
$j$--corridor--boundary $\G_\bt$ with vertex $v$ then $v$ cannot be of
type $Fj(a)$. Hence $v$ is distance at most 3 from $\cV_G$,
from Theorem \ref{bad-vertex}. The number of vertices distance at most
3 from $\cV_G$ is less than or equal to 
$c_1=(ml+1)^3|\cV_G|$ so there are at most
$c_1$ corridors with a $j$--corridor--boundary, $j\ge 1$. The result follows.  
\medskip

\begin{lemma}   \label{N-number}
The number of vertices in ${\cal V}_N$ is bounded above by
$$
|\cV_N| \leq  30 (2(ml+1)W_1 + W_0 -2\chi (\S )).
$$
\end{lemma}

\noindent
{\em Proof.}   
From (\ref{total-curve}) and Lemma \ref{zero-curve} we have
$$
2\pi \chi (\S ) \leq  \sum _{v\in \cV_N} \ka (v) + \sum _{v\in \cV_{1,B}}\ka (v)
+ \sum_{\D\in \cR} \ka (\D ) + \sum_{\bt\in \cD} \ka (\bt )
$$
and, using Lemma \ref{HLbound}, Lemma \ref{bd-curve} and  Lemma \ref{r-curve}, this gives
$$2\pi\chi(\S)\le - (\pi / 30) |\cV_N|+ 2(ml+1)W_1\pi +W_0\pi,$$
from which the result follows.
\begin{defn}\label{Bbound}
We define
\begin{align*}
\cB_1= \cB_1(\mbf z) & = 2W_1+30(2(ml+1)W_1+W_0-2\chi(\S)) \textrm{ and }\\
\cB = \cB(\mbf z)  &=2W_1+ ml((ml)^3+1)\cB_1+M(\mbf z)(n+3((ml+1)^4\cB_1+W_0-2\chi(\S))\\
 &+(6ml+10)W_1).
\end{align*}
\end{defn}

\begin{theorem}   \label{arc-bound}
Let $\G $ be a picture over $G$ as described at the beginning of this
section.
Then the 
number of arcs $|{\cal A}|$ of $\G$ satisfies $|{\cal A}|\leq {\cal B}$.
\end{theorem}

\noindent
{\em Proof.}   Let $y$ be an arc of $\G$ and assume first that $y$ is
not contained in the $Z$--subspace.  
If $y$ is incident to no vertex then the conditions on $\G$ imply
that $y$ is an $H^\Lm$--arc. 
If $y$ is incident to a vertex $v$ then $v\in \cV_G\cup \cV_D$, as if
$v$ is of type $Fj(a)$ then $y$ is contained in the $Z$--subspace. 
Hence the number of arcs of $\G$ which are not contained in the
$Z$--subspace of $\G$ is at most 
$$2W_1 +ml|\cV_G\cup\cV_D|.$$
From Corollary \ref{corridor-size} and Lemma \ref{corridor-number}
there are at most $$M(\mbf z)(n+
3((ml+1)^4|\cV_G| + W_0 - 2\chi (\S ))+(6ml+4)W_1)$$
arcs in the $Z$--subspace (note that $M(\mbf z)=M({\mbf {\ovr z}})$: see Section \ref{qpe_and_genus}).
The result now follows from Lemma \ref{D-number} and Lemma \ref{N-number}. 

\section{Proof of Theorem \ref{mainc}}\label{proof_mainc}

Let $\cN=\cN(\mbf h,\mbf n,\mbf t,\mbf p)$ and $l=l(r)$.
Consider first the $\cPG(\mbf h,\mbf n,\mbf t,\mbf p)$--problem when $\mbf h$ is 
a partition of $h=0$. In this case the problem reduces to 
the $\cQ$--problem and, as in Example \ref{exx_classic_genus}, under the 
current hypotheses is solvable in 
$(G,A,B,H_i:i\in I)$.

We may therefore assume that $\mbf h$ is a partition of a positive integer $h$, that
$\mbf n$ is a partition of a positive integer and that there is some $i$ such that 
$h_in_i>0$,
since otherwise the
system reduces to the $\cQ$--problem again. It follows that the 
positive $4$--partition, of length $1$, $(1,1,0,0)\in \cN$. The hypothesis on 
$(A,B,H_i : i\in I)$ now implies that $A$, $B$, $H_i$ all have soluble word problem.
It follows from Lemma \ref{specialg} that, in
$(G, A,B,H_i:i\in I)$, the 
$\cPG(\mbf h,\mbf n,\mbf t,\mbf p)$--problem is solvable if the 
$\cSPG(\mbf h^{\pprime},\mbf n^{\pprime},\mbf t^{\pprime},\mbf p^{\pprime})$--problem
is  solvable for all positive $4$--partitions 
$(\mbf h^{\pprime},\mbf n^{\pprime},\mbf t^{\pprime},\mbf p^{\pprime})$
of $(h,n^\pprime,t,p)$ 
such that $n^\pprime\in \ZZ$ and
$(\mbf h^{\pprime},\mbf n^{\pprime},\mbf t^{\pprime},\mbf p^{\pprime})
\le (\mbf h,\mbf n,\mbf t,\mbf p)$. Such partitions 
$(\mbf h^{\pprime},\mbf n^{\pprime},\mbf t^{\pprime},\mbf p^{\pprime})$ are
in $\cN$. It therefore suffices, using the argument at the beginning of the 
proof of Theorem \ref{main} in Section \ref{dec_in_op}, to show that the 
$\cSPG(\mbf h,\mbf n,\mbf t,\mbf p)$--problem is solvable 
in $(G, A, B, H_i:i\in I)$, where $\mbf n$ is a partition of
an integer $n\ge 1$. 

Let $i_A$, $i_B$ and $i_*$ be symbols not in $I$ and let $J=I\cup\{i_A,i_B\}$. Let
$K=J\cup\{i_*\}$, $X_{i_*}=\{i_A,i_B\}$ and $X_j=\{j\}$, for all $j\in J$. Set
$H_{i_A}=A$ and  $H_{i_B}=B$ so $H_{*,i_*}=H$ and, with $\mbf s_{i_*}=\{s\}$ we have
$G_{i_*}=G$. Set $\mbf s_j=\emptyset$, for all $j\in J$. Then 
$H_{*,i_A}=A$, $H_{*,i_B}=B$, $H_{*,i}=H_i$, for all $i\in I$ and $G_j=H_{*,j}$, for
all $j\in J$.
If $\cL$ is a consistent system of parameters and
\[\mbf z \in \left(\bdcup_{k\in K} H^\Lm_{*,k},\cL,\mbf n,\mbf h
              \right)\]
then the equation $Q$ associated to $\mbf z$ has environment 
\[((X_k)_{k\in K};(H_j)_{j\in J}; (G_k)_{k\in K})^\Lm.\] 
The $\cSPG(\mbf h,\mbf n,\mbf t,\mbf p)$--problem is solvable if,
given such $\mbf z$, which is in addition special, there exists an algorithm to decide whether
or not $Q$ has a solution and if so to output one.

Assume then that $\mbf z$ is a special element of 
$\left(\bdcup_{k\in K} H^\Lm_{*,k},\cL,\mbf n,\mbf h
              \right)$. Reordering the elements of $\mbf z$ if necessary we may assume that
we have a positive integer $u$ such that
$\mbf z$ is supported over $\mbf n$ with basis $(a_1,\ldots ,a_k)$, where $a_i=i_*$, for $i=1,\ldots ,u$ and 
$a_i\in J$, for $i>u$.
Then  $\mbf z= \mbf z_1,\mbf z_2$,
$\mbf h=\mbf h_1,\mbf h_2$ and $\mbf n=\mbf n_1,\mbf n_2$
where 
\begin{align*}
\mbf z_1 &=(z_1,\ldots ,z_{\mu_{u+1}})\in (H^\Lm, \cL,\mbf n_1,\mbf h_1),\\
\mbf z_2 &= (z_{\mu_{u+1}+1},\ldots ,z_n)\in \left(A^\Lm \dcup B^\Lm \dcup \bdcup_{i\in I} H^\Lm_i,\cL,\mbf n_2,\mbf h_2
              \right),\\
\mbf h_1 &= (h_1,\ldots ,h_u),\\
\mbf h_2 &= (h_{u+1},\ldots ,h_k),\\
\mbf n_1 &= (n_1,\ldots ,n_u)\textrm{ and }\\
\mbf n_2 &= (n_{u+1},\ldots ,n_k),
\end{align*}
where
$\mu_i=\mu_i(\mbf n)$ as defined in
(\ref{muj}).
Define 
\begin{align*}
\mbf t_1 &=(t_1,\ldots, t_u),\\
\mbf t_2 &=(t_{u+1},\ldots , t_k),\\
\mbf p_1 &=(p_1,\ldots, p_u)\textrm{ and}\\
\mbf p_2 &=(p_{u+1},\ldots , p_k),
\end{align*}
so that $\mbf t=\mbf t_1,\mbf t_2$ and $\mbf p=\mbf p_1,\mbf p_2$.
Then 
$$(\mbf h_i,\mbf n_i,\mbf t_i,\mbf p_i)\in \cN,\textrm{ for } i=1,2.$$
Let 
$$Q=Q(\mbf z,\cL,\mbf n,\mbf t,\mbf p)=(\mbf q=1,\bt,\cL),$$ 
let
$$\mbf q_1=(q(\xi_j,n_j,t_j,p_j): j=1,\ldots u)\textrm{ and } 
\mbf q_2=(q(\xi_j,n_j,t_j,p_j): j=u+1,\ldots k),$$
$$\bt_i=\bt|_{L_D(\mbf q_i)}\textrm{ and }Q_i=(\mbf q_i=1, \bt_i,\cL), \textrm{ for } i=1,2.$$

Let $Q_\cH$ be the homogeneous equation associated to $\mbf z$ and $\cL$ and suppose
that $Q_\cH=(\mbf q=1,\bt_\cH,\cH)$. Define $\cH_i$ to be the set of elements of $\cH$
containing parameters which occur in $\ovr{\mbf z_i}$ and 
$$Q_{\cH_i}=(\mbf q_i=1,\bt_{\cH_i},\cH_i),$$
where $\bt_{\cH_i}=\bt_\cH |_{L_D(\mbf q_i)},$ for $i=1,2$. 
Then $Q_{\cH_1}$ is the homogeneous equation associated to $\mbf z_1$.
Up to relabelling of parameters it is also the case that $Q_{\cH_2}$ is the 
homogeneous equation associated to $\mbf z_2$.
As $\mbf z$ is special
so is $\mbf z_i$ and hence also, using Lemma \ref{homog-special}, is $Q_{\cH_i}$. 
The equations $Q_1$ and $Q_{\cH_1}$ have environment
\[(X_{i_*};H_{i_A},H_{i_B};G_{i_*})=(X_{i_*};A,B;G)^\Lm,\]
while $Q_2$ and 
$Q_{\cH_2}$ have environment 
\[((X_j)_{j\in J};(H_j)_{j\in J};(H_j)_{j\in J})^\Lm.\]

If $(\phi ,\al)$ is a solution to the equation $Q$ 
then define 
\begin{equation*}
\phi_i=\phi|_{L(\mbf q_i)},\textrm{ for } i=1,2.
\end{equation*}
Then $(\phi_1,\al)$ is a solution to $Q_1$ and $(\phi_2,\al)$
is a solution to $Q_2$.
Conversely suppose we are 
given solutions $(\phi_i,\al)$ to $Q_i$, for $i=1,2$.
Define
\begin{gather}
\phi:F(L(\mbf q))\maps H*A*B**_{i\in I}H_i\nonumber\\ 
\textrm{by }\phi(x)=\phi_i(x),\textrm{ for }x\in L(\mbf q_i).\label{phi1_2}
\end{gather}
 Then $(\phi ,\al)$ is a solution
to $Q$.

Assuming that a solution $(\phi, \al)$ to $Q$ exists
then there is a picture $\G$ on a surface $\S$ of type $(\mbf n_1,\mbf t_1,\mbf p_1)$ 
such that$(\al,\G)$ corresponds
to the solution $(\phi_1,\al)$ to $Q_1$.
A minimal picture on $\S$ with the same boundary labels as $\G$ then exists, and corresponds
to a, possibly different, solution to $Q_1$. Together with $\phi_2$ this solution to $Q_1$ can be
used, as in (\ref{phi1_2}), to construct a new solution to $Q$. Hence we may assume that $\G$ is a minimal picture.
(See Definition \ref{pic_soln}, Definition
\ref{pic_min} 
and Corollary \ref{pic_qpe}.) Let $\G_s$ be a $Z$--reduced picture obtained from 
$\G$ by $Z$--cancellation. Then $\G_s$ is reduced and minimalistic and, using
Lemma \ref{follow-path}(\ref{follow-path-2}), 
determines a solution $\al_s$ to $\cH_1$ such that
$(\al_s,\G_s)$ corresponds to some solution $(\phi_s,\al_s)$ to 
$Q_{\cH_1}$: 
that is $\G_s$ is a basic picture for $(\mbf z_1,\cL)$. 
From Theorem \ref{basic-reduce} there exists a solution $\al_0$ to 
$\cL\cup \cL(\G_s,P)$, where $P$ is the subgraph of some path in the $Z$--graph
$G(\mbf z_1)$ of $\mbf z_1$, with initial vertex $v(\G_s)$. Furthermore $\al_0$ may
be chosen so that $\al_0(\lm)=\al(\lm)$, for all $\lm\in\Lm$ occuring in $Q$.

Since partisan regions are not involved in $Z$--cancellation, the partisan regions
of $\G$ and $\G_s$ are the same. Let $\D$ be the closure of a partisan region
of $\G$. Assume that $\D$ is a $V$--region, where $V=A$ or $B$ and that $\D$ is the
sum of $t^\prime=t^\prime(\D)$ torii with $p^\prime=p^\prime(\D)$ 
projective planes and has
$k^\prime=k^\prime(\D)$ boundary components $\bt_1^\prime,\ldots, \bt_{k^\prime}^\prime$.
Assume further that $\bt_1^\prime,\ldots, \bt_{d^\prime}^\prime$ contain partisan boundary
intervals of $\G$, whilst $\bt_{d^\prime+1}^\prime,\ldots, \bt_{k^\prime}^\prime$ do not.
Fix as base point in $\bt_j^\prime$ an end point of an arc of $\G$ (or 
a point of a boundary partition of a component of $\pd \S$, if $\bt^\prime_j$ does
not meet an arc of $\G$).  
Starting at the base point and 
reading prime labels of boundary corners and labels of vertex corners as they occur
in order around $\bt_j^\prime$ gives the {\em prime label} 
$$v_j=v_j(\D)=\prod_{i=1}^{n^\prime_j}(h_{j_i},f_{j_i})\in V^\Lm$$
of $\bt^\prime_j$.
Reading labels of both boundary corners and vertex labels we obtain the label of
$\bt^\prime_j$ which, by construction, is equal to 
$\hat{\al}(v_j)$. If we identify $(v,1)\in V^\Lm$ with $v\in V$ then,
for $j=d^\prime+1,\ldots ,k^\prime$, the prime
label is the same as the label of $\bt^\prime_j$, since $\bt^\prime_j$ contains 
no partisan boundary intervals for such $j$. Thus, for $j=d^\prime+1,\ldots ,k^\prime$,
the label (and prime label) of $\bt^\prime_j$ is a fixed element $v_j$ of $V$.
The prime labels of boundary components of $\D$ are called the {\em prime}
boundary labels of $\D$. 

Setting 
$$\mbf v(\D)= v_1,\ldots ,v_{k^\prime}$$ and 
$h^\prime(\D) = \sum_{j=1}^{k^\prime}|v_j|_\Lm$
we have $\mbf v(\D)\in (V^\Lm,\cL,(k^\prime), (h^\prime))$ and, using property
$P$\ref{pic_cond2} of pictures, $\D$ determines a solution 
to $Q(\mbf v(\D), \cL, (h^\prime), (k^\prime))$ in $V$.

We have 
\[\G=\G_1\dcup\cdots \dcup\G_u,\]
where $\G_i$ is a picture on a surface $\Sigma_i$ of type $(n_i,t_i,p_i)$.
Let $\G_0$ be the connected component of $\G$ containing $\D$. Then $\G_0$ is 
a component of a sub--picture $\G_i$ of $\G$, 
for some $i$ with $1\le i \le u$. Thus $\G_0$ is
a picture on a connected surface $\S_0$ of type $(n_0,t_0,p_0)$, where $n_0\le n_i$,
$t_0\le t_i$  and $p_0\le p_i$.  As $\D$ is a region of a picture on $\S_0$ we 
have 
\begin{equation}\label{degn_reg_cond_1}
2t^\prime + p^\prime \le 2t_0+p_0\le 2t_i+p_i. 
\end{equation}
As $\G$ is minimal so is $\G_0$ and so, from Proposition \ref{genreg}, we have
$\chi(\D)\ge \chi(\S_0)$. That is 
\begin{align}
2-2t^\prime -p^\prime -k^\prime &\ge 2- 2t_0-p_0-n_0\ge 2-2t_i-p_i-n_i \nonumber \\
\textrm{ so }k^\prime &\le n_i+2(t_i-t^\prime)+p_i-p^\prime \label{degn_reg_cond_2}.
\end{align}

Let $\G$ have a total of $n^\prime$ partisan regions $\D_1,\ldots , \D_{n^\prime}$.
Assume that $\D_j$ is a  $V_j$--region, where $V_j=A$ or $B$, and that $\D_j$ is the
sum of $t^\prime(\D_j)$ torii with $p^\prime(\D_j)$ 
projective planes and has
$k^\prime(\D_j)$ boundary components of which $d^\prime(\D_j)$ are partisan 
whilst the others are
non--partisan.
Let $l_i$ be the number of partisan regions of $\G_i$, let 
\begin{equation*}\label{parti_count}
\w_0=0\textrm{ and let } \w_i=\sum_{j=1}^i l_i, \textrm{ for } i=1,\ldots ,u,
\end{equation*}
so $\w_u=n^\prime$. We may assume that 
$\D_j$ is a region of $\G_i$ for all $j$ such that
$\w_{i-1}+1\le j\le \w_i$. 
Set
\begin{align*} 
\mbf v_j &=\mbf v(\D_j), \textrm{ for } j=1,\ldots ,n^\prime,\\
\mbf v &= \mbf v_1,\ldots ,\mbf v_{n^\prime},\\
\mbf h^\prime & = h^\prime(\D_1),\ldots , h^\prime(\D_{n^\prime}),\\
\mbf k^\prime & = k^\prime(\D_1),\ldots , k^\prime(\D_{n^\prime}),\\
\mbf t^\prime & = t^\prime(\D_1),\ldots , t^\prime(\D_{n^\prime})\textrm{ and }\\
\mbf p^\prime & = p^\prime(\D_1),\ldots , p^\prime(\D_{n^\prime}).
\end{align*}
Then 
\begin{equation}\label{degn_reg_cond_3}
\mbf v\in (A^\Lm\dcup B^\Lm,\cL,\mbf h^\prime, \mbf k^\prime).
\end{equation}
Also $0\le l_i\le h_i$, as the boundary of every partisan region contains
at least one proper exponential $H$--letter. Moreover 
\[\sum_{j=\w_{i-1}+1}^{\w_i} h^\prime(D_j)\le h_i\]
and 
(\ref{degn_reg_cond_1}) and 
(\ref{degn_reg_cond_2}) hold, with $k^\prime(D_j),t^\prime(D_j),p^\prime(D_j)$ in place of 
$k^\prime,t^\prime,p^\prime$, for 
$\w_{i-1}+1\le j\le \w_i$ and 
$i=1,\ldots ,u$. Hence 
\begin{equation}\label{degn_reg_cond_4}
(\mbf h^\prime, \mbf k^\prime, \mbf t^\prime, \mbf p^\prime)
\in \cN(\mbf h_1,\mbf k_1, \mbf t_1, \mbf p_1)\subseteq \cN.
\end{equation}
The equation associated to $\mbf v$ is $Q_\D=Q(\mbf v,\cL,\mbf h^\prime,\mbf k^\prime)$,
an equation with environment $(X_{i_A},X_{i_B};A,B;A,B)^\Lm$.
Set
$\phi_{\D}=\phi |_{L(\mbf v)}$. Then, since $\al_0(\lm)=\al(\lm)$ for all $\lm$
occuring in $Q$,  the pair $(\phi_\D,\al_0)$ is a solution to 
$Q_\D=Q(\mbf v,\cL,\mbf h^\prime,\mbf k^\prime)$.
Let $\mbf q_\D$ be the system of quadratic words in standard form and $\bt_\D$ the map
from $L_D(\mbf q_\D)$ to $A^\Lm\dcup B^\Lm$ such that
$Q_\D=(\mbf q_\D=1,\bt_\D,\cL)$. 

Let $(h_1,f_1),\ldots, (h_{W_1},f_{W_1})$ be the letters of $\mbf z_1$, where
$W_1=W_1(\mbf z_1)$. Form a system of parameters $\cL^-(\G_s,P)$ from $\cL(\G_s,P)$
by deleting those equations $f_q=\al_s(\lm_q)$ from $\cL(\G_s,P)$, for all $q$ such
that $(h_q,f_q)$ is a proper exponential $H$--letter with $h_q\in A\dcup B$.
Now, since elements of $\Lm$ occurring in $\cL(\G_s,P)$ do not occur in 
$Q_1$ and we may choose then so that they do not occur in $Q$, we have
 $$\mbf z^\pprime= 
\mbf v,\mbf z_2\in 
(A^\Lm\dcup B^\Lm\dcup\bdcup_{i\in I}H^\Lm_{i},\cL\cup\cL^-(\G_s,P),\mbf h^\pprime,\mbf n^\pprime),$$
where $\mbf h^\pprime=\mbf h^\prime, \mbf h_2$, and 
$\mbf n^\pprime = \mbf n^\prime, \mbf n_2$. Define 
$ \mbf t^\pprime = \mbf t^\prime, \mbf t_2$ and 
$\mbf p^\pprime = \mbf p^\prime, \mbf p_2$  and 
\begin{align*}
Q^\pprime &= Q(\mbf z^\pprime,\cL\cup\cL^-(\G_s,P),\mbf n^\pprime,\mbf t^\pprime,
\mbf p^\pprime)\\
&=(\mbf q^\pprime=1, \bt^\pprime, \cL\cup \cL^-(\G_s,P)),
\end{align*}
a system of equations with environment 
$((X_j)_{j\in J}; (H_j)_{j\in J};(H_j)_{j\in J})^\Lm$.
Define 
$$\phi_0:F(L(\mbf q^\pprime))\maps A*B**_{i\in I}H_i$$
by $\phi_0|_{L(\mbf q_\D)}=\phi_\D$ and $\phi_0|_{L(\mbf q_2)} = \phi|_{L(\mbf q_2)}
= \phi_2$. Then $(\phi_0, \al_0)$ is a solution to $Q^\pprime$ and 
$(\mbf h^\pprime, \mbf n^\pprime, \mbf t^\pprime, \mbf p^\pprime)\in \cN$.

Now $\G_s$ is a picture over $G$ on a surface of type $(\mbf n_1,\mbf t_1,\mbf p_1)$ 
with partisan regions $\D_1,\ldots ,\D_{n^\prime}$. Removing boundary labels from
partisan boundary intervals of $\D_j$ we replace these labels with the prime labels
of the corresponding boundary intervals of $\G_s$. We call the result a picture with
{\em partisan boundary labels stripped} and denote it $\G_s^-$. Note that given
$(\mbf h, \mbf n,\mbf t,\mbf p)$, $\mbf z$ and $\cL$ the graph $G(\mbf z_1)$ is uniquely
determined.
Given $\G_s$ the vertex $v(\G_s)$ of $G(\mbf z_1)$ is unique and, since 
$\{$active left markings of $\G_s \} = \{$active left markings of $\G_s^-\}$ we
may determine $v(\G_s)$ from $\G_s^-$. We may therefore refer to this
vertex as $v(\G_s^-)$. Furthermore, from the above and 
Theorem \ref{basic-reduce}, if there
exists a solution $(\phi, \al)$ to $Q(\mbf z, \cL, \mbf h,\mbf n)$ then there exists
a $Z$--path in $G(\mbf z_1)$ with path--subgraph $P$ and  with 
initial vertex $v(\G_s^-)$ and a solution $\al_0$
to $\cL\cup \cL(\G_s,P)$. Again $\cL^-(\G_s,P)$ may be constructed using only $\G_s^-$,
since the equation of $\cL(\G_s,P)$ corresponding to a partisan 
boundary interval with prime label $(h_q,f_q)$, the $q$th letter of $\mbf z_1$,
 on partisan
boundary component of $\G$ has the form $f_q=\al_s(\lm_q)$, and so has been deleted. 
Hence we may refer to $\cL^-(\G_s,P)$ as $\cL^-(\G_s^-,P)$. Also, given
$\mbf z$, $\cL$, $\G_s^-$ and $\cL^-(\G_s^-,P)$ we may construct the system $Q^\pprime$,
which is again uniquely determined by these data. The existence of the solution 
$(\phi,\al)$ then guarantees the existence of the solution $(\phi_0,\al_0)$ to 
$Q^\pprime$. Hence given a solution $(\phi,\al)$ to $Q$ we obtain the following.
\be[(a)]
\item\label{geta} 
A picture $\G_s^-$ with partisan boundary labels stripped which is obtained
(by stripping partisan boundary labels) from  a $Z$--reduced, reduced, minimalistic,
picture $\G_s$ such
that $(\al_s,\G_s)$ corresponds to a solution to $Q_{\cH_1}$ for some retraction 
$\al_s:M\maps \ZZ$.
\item\label{getb} The path--subgraph $P$ of some $Z$--path 
in $G(\mbf z_1)$ with initial vertex $v(\G_s^-)$ 
and a corresponding system 
of parameters $\cL^-(\G_s^-,P)$.
\item\label{getc} A system of equations 
$Q^\pprime$ as above and solution $(\phi_0,\al_0)$ to $Q^\pprime$.
\ee 

Conversely, given the data of (\ref{geta}),  (\ref{getb}) and   (\ref{getc}) 
we may form a picture $\G_t$ from $\G_s^-$ as follows. 
Let $(h_1,f_1),\ldots ,(h_{W_1},f_{W_1})$ be the letters of $\mbf z_1$ and let
$b_q$ be the boundary interval of $\G_s^-$ with prime label $(h_q,\lm_q)$.  Form $\G_t$ by
assigning to the partisan boundary
interval $b_q$ of $\G_s^-$ the label 
$\hat\al_0(h,f)$. 
Define a retraction
$\al_t:M\maps \ZZ$ by $\al_t(\lm_q)=\al_0(f_q)$, if $b_q$ is a partisan boundary interval
and $\al_t(\lm_q)=|b_q|$ otherwise. Then $\al_t$ is a solution to $\cH_1$ and the boundary interval
$b_q$ of $\G_t$ has label $\hat \al_t(h_q,\lm_q)$, for $q=1,\ldots, W_1=W_1(\mbf z_1)$. Hence 
$(\al_t,\G_t)$ corresponds to a solution $(\phi_t,\al_t)$ to $Q_{\cH_1}$. Note that $v(\G_t)=
v(\G_s^-)=v(\G_s)$ in $G(\mbf z_1)$ and that $\cL^-(\G_s,P)=\cL^-(\G_t,P)$. Therefore
$\al_0$ is a solution to $\cL^-(\G_t,P)$ and so also to $\cL(\G_t,P)$. It follows from 
Theorem \ref{basic-reduce} that there exists a solution $\al_c$ to $\cL$ and a picture
$\G_c$ obtained from $\G_t$ by $Z$--insertion along a path with path--subgraph $P$ such that
$(\al_c,\G_c)$ corresponds to a solution $(\phi_c,\al_c)$ to $Q_1$. In fact, since $\al_0$ is a 
solution to $\cL\cup\cL(\G_t,P)$ and $\al_0(\lm)=\al_c(\lm)$, for all $\lm$ occuring in $\mbf z_1$ or $\cL$
we may take $\al_c=\al_0$. From $(\al_0,\G_0)$ we obtain a solution $(\phi_c,\al_0)$ to $Q_1$.
If $\phi_2=\phi_0|_{L(\mbf q_2)}$ then $(\phi_2,\al_0)$ is a solution to $Q_2$. Using
(\ref{phi1_2}) we obtain a solution $(\phi,\al_0)$ to $Q$.

 It therefore suffices to 
determine whether the data of  (\ref{geta}),  (\ref{getb}) and   (\ref{getc})
exist and if so to find it. We describe an algorithm, similar to 
those described in \cite{\DHa}, \cite{\DHb} and \cite{\DHc}, which solves this problem.
First note that given $\mbf z$ (which we are assuming is special)
we may construct $\mbf z_1$ and calculate $\cB=\cB(\mbf z_1)$ of Definition \ref{Bbound},
since the word problem is solvable in $A$ and $B$. The number of arcs of $\G_s$,
if it exists, is at most $\cB$, using Theorem \ref{arc-bound}. The same is true
of $\G_s^-$. Assume that $\G_s$ exists and that 
the set of arcs of $\G_s^-$ is $\cA$, so 
$0\le |\cA|\le \cB$. There are, up to homeomorphism, finitely many surfaces of 
type $(\mbf n_1,\mbf t_1,\mbf p_1)$ which have no closed components. We fix one such surface, call it $\S$ and
assume that $\G_s$ is a picture on $\S$.
Then $\S$ is a disjoint union of surfaces $\S_1,\ldots ,\S_u$, where $\S_j$ is
of type $(n_j,t_j,p_j)$, for $j=1,\ldots, u$. Note that $\S_j$ need not be connected
and each connected component of $\S$ is a connected component of $\S_j$,  for some 
$j$ with $1\le j\le u$.
Assume then that $\S$ has $c$ connected components $\W_1,\ldots ,\W_c$, for some $c$ such that
$u\le c\le \mu_{u+1}$. 
If $\W$ is  a connected component of $\S$ we denote by
$\G^-_{s}(\W)$ the restriction of $\G_s^-$ to $\W$.  
There are finitely many partitions $(A_1,\ldots, A_c)$
of $|\cA|$ and (if $\G_s^-$ exists) one such partition satisfies $A_j=$ 
the number of arcs of 
of $\G_s^-(\W_j)$, $j=1,\ldots ,c$.  
Fix some $j$, with $1\le j\le c$ and write $\W=\W_j$. Let $\W$ be the connected
sum of $t_\W$ torii and $p_\W$ projective planes and have $n_\W$ boundary components
with prime boundary labels $z_i$, for $i=d+1,\ldots , d+n_\W$ and an appropriate
integer $d$.  

Let $\W_0$ be the surface
$$\W_0=\W\backslash\{\textrm{int}(v): v \textrm{ is a vertex of } \G_s^-(\W)\}.$$
Let $V_j$ be the number of vertices of $\G_s^-(\W)$. As $\G_s^-(\W)$ has $A_j$ 
arcs and its vertices all have index $ml$ it follows that
$V_j\le A_j/ml$. 

Let $\W_{00}$ be the surface formed by cutting $\W_0$ along arcs of $\G_s^-(\W)$.
Each connected component of $\W_{00}$ is the closure of a region of $\G_s^-(\W)$.
Since $\G_s$ is minimalistic so is $\G_s(\W)$, so for each region $\D$ of $\G_s^-(\W)$ we
have $\chi(\D)\ge \chi(\W)$.
As $\W_0$ is connected $\W_{00}$ has at most $A_j+1$ connected components and,
by construction, 
\begin{align*}
\chi(\W_{00}) & = \chi(\W)-V_j+A_j\\
 & = 2 - t_\W -p_\W/2 -n_\W -V_j+A_j\\
 & \ge 2 - (t_\W +p_\W/2) -n_\W +A_j(1-1/ml).
\end{align*}

Set $\chi_j= 2 - (t_\W +p_\W/2) -n_\W -A_j(1-1/ml)$. There are finitely many
surfaces of Euler characteristic $\chi_j$ with at most $A_j+1$ connected 
components each of Euler characteristic at least $\chi(\W)$. 
The algorithm takes each such surface and attempts to construct 
$\G^-_s(\W)$ from it. Given a such a surface which we take to be $\W_{00}$ 
suppose that $\pd \W_{00}$ has connected components $c_1,\ldots ,c_D$. There
are finitely many ways of assigning a positive integer $\rho_i$ to $c_i$, for $i=1,\ldots,D$,
so that $\sum_{i=1}^D \rho_i=2(A_j+1)$. Choose one such assignment and divide $c_i$ into
$\rho_i$ intervals, $i=1,\ldots,\rho_i$. 
There are now finitely
many ways of identifying intervals in pairs to obtain a surface of type 
$(n_\W+V_0,t_\W,p_\W)$, for each choice of $V_0$ such that $0\le V_0\le A_j/ml$. 
If $\G_s^-$ exists then one
such gives rise to a surface homeomorphic to  
$\W_0$. Fix such a choice of $V_0$ and of pairings and assume the resulting
surface is $\W_0$. The paired intervals become properly embedded arcs on $\W_0$ which
we refer to as arcs, as they become arcs of $\G_s^-$ if the process is successful.
To obtain $\W$ all but $n_\W$ of the boundary components of 
$\W_0$ are capped off with disks, which we refer to as vertices (they'll become
the vertices of $\G_s^-(\W)$). 
This must be done in such a way that each vertex meets exactly $ml$ arcs. There are
finitely many ways of doing this, one of which results in a surface homeomorphic to 
$\W$. Fix one such choice and assume the resulting surface is $\W$ and that 
the vertices and arcs are those of $\G^-_s(\W)$.
Images of connected components of $\W_{00}$ in $\W$  are called
regions of $\W$. Choose one of the finitely many possible orientations of $\G_s^-(\W)$
consistent with the definition of orientation of a sketch in Section \ref{pictures}.
The union of $\W$ with these arcs, vertices, regions and orientation is called
$\G_b$.

To recover $\G_s^-(\W)$ the components of 
$\{$vertices of $\G_b\}\backslash \{$arcs of $\G_b\}$ and of
$\{$boundary components of $\W\}\backslash \{$arcs of $\G_b\}$ must be labelled to give a 
picture with partisan boundary labels stripped. There are 
finitely many ways of labelling corners of vertices with elements of $A\dcup B$ 
so that each vertex has label
$r^{\pm m}$.  
There are finitely many ways of assigning the $z_j$'s to 
the boundary components of $\W$. 
This must be done in such a way that if $z_{d+j}\in A^\Lm\dcup B^\Lm$ then
$z_{d+j}$ is assigned to a boundary component which meets
no arc of $\G_b$. Suppose $z$ is assigned to $\bt$. 
There are then finitely 
many ways of assigning a boundary
partition to $\bt$ which $\cH_1$ admits $z$, since the arcs of $\G_b$ are already fixed.
Suppose that a boundary partition of $\bt$, which $\cL$--admits $z$ has been chosen. There
are finitely many ways of assigning the letters of $z$ as prime labels of the boundary
intervals of this boundary partition. Assume then that $\bt$ has prime label $z$ with some 
boundary partition. 
If $[b,c]$ is a  non--partisan boundary interval 
of this boundary partition then the label of $[b,c]$ and of all components 
of $[b,c]\backslash \G_b$ is fixed by the assignations made so far. Choosing boundary
partitions and prime labels  in this way for all boundary components of $\pd \W$ we obtain 
labels for all corners of $\G_b$ except those containing partisan boundary intervals.

Since $A$ and $B$ have solvable word problem, 
the conditions $P$\ref{pic_cond0}, $P$\ref{pic_cond1} and $P$\ref{pic_cond3}, for
 $\G_b$ to be a picture with partisan boundary labels stripped
can now be checked, as prime labels are enough to verify
these conditions. Furthermore  $P$\ref{pic_cond2}
can be verified for every non--partisan region of $\G_b$. 
To see this note that, since $\G_s$ is obtained from 
$\G$ by $Z$--cancellation and $\G_s$ is minimalistic, by construction
$\x(\D)\ge \x(\W)$, for all regions $\D$ of $\G_b$. 

Suppose $\D$ is a region of
$\G_b$ 
which has $s^\prime$ boundary components 
all of which are non--partisan, with labels $(u_1,\cdots ,u_{s^\prime})$ in
$V$, for $V=A$ or $B$.
Suppose that $\D$ is a union of $t^\prime$ torii and $p^\prime$
projective planes. Since  $\x(\D)\ge \x(\W)$, 
(\ref{degn_reg_cond_1}) and (\ref{degn_reg_cond_2}) hold,
with $j$ in place of $i$. 
Setting 
$$\mbf u^\prime= u_1,\ldots ,u_{s^\prime}$$ 
we have $\mbf u\in (V^\Lm,\cH_1,s^\prime, t^\prime, p^\prime)$.
Furthermore
$(0, s^\prime,t^\prime,p^\prime)\in \cN$ and so, 
by hypothesis, 
the ${\cPG}(0, s^\prime,t^\prime,p^\prime)$--problem
is solvable in $(A,B)$. Therefore we may check whether or not 
$(s^\prime, t^\prime,p^\prime)\in \cH_1$--genus$(\mbf u^\prime)$ in $V$: that is 
whether condition $P$\ref{pic_cond2} holds for $\D$. 
If all these conditions are satisfied then $\G_b$ is a candidate
for $\G_s^-(\W)$. We now repeat this process for all connected components
$\W=\W_1,\ldots ,\W_c$ of our candidate $\S$. If successful, this gives rise to 
a candidate for $\G^-_s$.

Assume we have a candidate for $\G^-_s$, which we continue to call $\G_b$.
Next we consider the partisan regions of $\G_b$. Using the same
notation for partisan regions of $\G_b$ as used for $\G_s$ above, 
the prime labels boundary components of partisan
boundary regions of $\G_b$ give rise to an element  
$$
\mbf v\in (A^\Lm\dcup B^\Lm,\cL,\mbf h^\prime, \mbf k^\prime)
$$
as in (\ref{degn_reg_cond_3}) and furthermore (\ref{degn_reg_cond_4}) 
holds, as before.

Since the word problem is solvable in $A$ and $B$ the graph $G(\mbf z_1)$ may
be constructed and the vertex $v(\G_s^-)$ located, using $\G_s^-$. 
There are finitely many path--subgraphs of $Z$--paths in
$G(\mbf z_1)$, since $G(\mbf z_1)$ is a finite graph, and we fix one such and call
it $P$. This allows construction of the system of parameters $\cL^-(\G_s^-,P)$.
We may now construct the system of equations 
$Q^\pprime=Q(\mbf z^\pprime,\cL\cup\cL^-(\G_s^-,P),\mbf h^\pprime,\mbf n^\pprime)$
described above. As 
$(\mbf h^\pprime,\mbf n^\pprime,\mbf t^\pprime,\mbf p^\pprime)\in \cN$
and $Q^\pprime$ has  environment 
$((X_j)_{j\in J}; (H_j)_{j\in J};(H_j)_{j\in J})^\Lm$,
there is an algorithm to determine whether or not $Q^\pprime$ has a solution and if
so to find one.

If a solution $(\phi_0,\al_0)$ to $Q^\pprime$ is found then, as above, a solution to $Q$ may
be constructed. If no such solution exists then either there is no solution to $Q$ or
one of the choices above was not the right one. At each stage there are finitely
many choices so once they've all been tried either we have a solution or there is 
none. This completes the description of the algorithm.

%% file: qpe.bbl
\begin{thebibliography}{10}

\bibitem{Baumslag60}
G.\ Baumslag.
\newblock Some aspects of groups with unique roots.
\newblock {\em Acta Math.}, 104:217--303, 1960.

\bibitem{Baumslag65}
G.\ Baumslag.
\newblock On free {Q}--group.
\newblock {\em Commun.\ on Pure and Appl.\ Math.}, 18:25--30, 1965.

\bibitem{BogleyPride92}
W.A.\ Bogley and S.J.\ Pride.
\newblock Aspherical relative presentations.
\newblock {\em Proc.\ Edinburgh Math.\ Soc.}, 35:1--39, 1992.

\bibitem{ComerfordEdmunds81}
Leo~P.\ Comerford and Charles~C.\ Edmunds.
\newblock Quadratic equations over free groups and free products.
\newblock {\em J.\ Algebra}, pages 276--297, 1981.

\bibitem{ComerfordEdmunds84}
Leo~P.\ Comerford and Charles~C.\ Edmunds.
\newblock Quadratic parametric equations over free groups.
\newblock {\em Cont.\ Math.}, 33:159--196, 1984.

\bibitem{Culler81}
M.\ Culler.
\newblock Using surfaces to solve equations in free groups.
\newblock {\em Topology}, pages 237--300, 1981.

\bibitem{DuncanHowie91a}
Andrew~J.\ Duncan and J.\ Howie.
\newblock The genus problem for one-relator products of locally indicable
  groups.
\newblock {\em Math.\ Z.}, 208:225--237, 1991.

\bibitem{DuncanHowie91b}
Andrew~J.\ Duncan and J.\ Howie.
\newblock The nonorientable genus problem for one-relator products.
\newblock {\em Commun.\ in Alg.}, 19:2547--2556, 1991.

\bibitem{DuncanHowie93}
Andrew~J.\ Duncan and J.\ Howie.
\newblock One relator products with high--powered relators.
\newblock In Graham~A.\ Niblo and Martin~A.\ Roller, editors, {\em Geometric
  Group Theory, Volume 1 of Proceedings of the Symposium held in Sussex 1991},
  volume 181 of {\em London Mathematical Society Lecture Note Series}, pages
  48--74. Cambridge University Press, 1994.

\bibitem{EdjvetThomas97}
M.\ Edjvet and R.M.\ Thomas.
\newblock The groups $(l,m|n,k)$.
\newblock {\em J. Pure Applied Algebra}, 114:175--208, 1997.

\bibitem{Fenn83}
R.A.\ Fenn.
\newblock {\em Techniques of Geometric Topology}.
\newblock Number~57 in London Math. Soc. Lecture Notes. Cambridge University
  Press, Cambridge, 1983.

\bibitem{BaumslagMyasnikovRemeslennikov99}
A.\~Myasnikov G.\~Baumslag and V.\ Remeslennikov.
\newblock Algebraic geometry over groups {I}. {A}lgebraic sets and ideal
  theory.
\newblock {\em J. of Algebra.}, 219:16--79, 1999.

\bibitem{Gersten87}
S.M.\ Gertsten.
\newblock Reducible diagrams and equations over groups.
\newblock In S.M.\ Gersten, editor, {\em Essays in group theory}, volume~8 of
  {\em MSRI Publications}, pages 15--73. Springer, 1987.

\bibitem{GoldsteinTurner79}
R.Z.\ Goldstein and E.C.\ Turner.
\newblock Applications of topological graph theory to group theory.
\newblock {\em Math.\ Z.}, 165:1--10, 1979.

\bibitem{GrigorchukLysionok92}
R.~I.~\ Grigorchuk and I.~G.\ Lysionok.
\newblock A description of solutions of quadratic equations in hyperbolic
  groups.
\newblock {\em Internat.\ J.\ Algebra Comput.}, 2(3):237--274, 1992.

\bibitem{Howie89}
J~Howie.
\newblock The quotient of a free product of groups by a single high--powered
  relator. {I}. {P}ictures. {F}ifth and higher powers.
\newblock {\em Proc.\ London Math.\ Soc. (3)}, 59:507--540, 1989.

\bibitem{Howie90}
J~Howie.
\newblock The quotient of a free product of groups by a single high--powered
  relator. {II}. {F}ourth powers.
\newblock {\em Proc.\ London Math.\ Soc. (3)}, 61:33--62, 1989.

\bibitem{Howie91}
J~Howie.
\newblock The quotient of a free product of groups by a single high--powered
  relator. {III}.the word problem.
\newblock {\em Proc.\ London Math.\ Soc. (3)}, 62:590--606, 1991.

\bibitem{HowieThomas93}
J.\ Howie and R.M.\ Thomas.
\newblock The groups $(2,3,p,q)$; asphericity and a conjecture of {C}oxeter.
\newblock {\em J.\ Algebra}, 154:289--309, 1993.

\bibitem{KharlampovichMyasnikov984}
O.\ Kharlampovich and A.\ Myasnikov.
\newblock Equations in a free {Q}--group.
\newblock {\em Trans.\ Amer.\ Math.\ Soc.}, 350(3):947--974, 1998.

\bibitem{KharlampovichMyasnikov982}
O.\ Kharlampovich and A.\ Myasnikov.
\newblock Irreducible affine varieties over a free group. {I}: {I}rreducibility
  of quadratic equations and nullstellensatz.
\newblock {\em J. of Algebra}, 200(2):472--516, 1998.

\bibitem{KharlampovichMyasnikov983}
O.\ Kharlampovich and A.\ Myasnikov.
\newblock Irreducible affine varieties over a free group. {II}: {S}ystems in
  triangular quasi--quadratic form and description of residually free groups.
\newblock {\em J. of Algebra}, 200(2):517--570, 1998.

\bibitem{LipschutzMiller71}
Seymour Lipschutz and Charles F.~Miller III.
\newblock Groups with certain solvable and unsolvable decision problems.
\newblock {\em Comm.\ Pure Appl.\ Math.}, XXIV:7--15, 1971.

\bibitem{Lyndon60a}
R.C.\ Lyndon.
\newblock Equations in free groups.
\newblock {\em Trans.\ Amer.\ Math.\ Soc.}, 96:445--457, 1960.

\bibitem{Lyndon60b}
R.C.\ Lyndon.
\newblock Groups with parametric exponents.
\newblock {\em Trans.\ Amer.\ Math.\ Soc.}, 96:518--553, 1960.

\bibitem{LyndonSchutzenberger62}
R.C.\ Lyndon and M.P.\ Sch{\"{u}}tzenberger.
\newblock The equation $a^m=b^nc^p$ in a free group.
\newblock {\em Michigan Math.\ J.}, 9:289--298, 1962.

\bibitem{Makanin82}
G.S.\ Makanin.
\newblock Equations in a free group. (russian).
\newblock {\em Izv.\ Akad.\ Nauk SSSR Ser.\ Mat.}, 46(6):1199--1273, 1982.
\newblock English translation in: Math. USSR--Izv. 21 (1983), no. 3, 546--582.

\bibitem{Matiyasevic93}
I.U.\ Matiyasevi{\u{c}}.
\newblock {\em Hilbert's tenth problem}.
\newblock Foundations of Computing. MIT Press, Cambridge, Mass., 1993.

\bibitem{MyasnikovRemeslennikov96}
A.G.\ Myasnikov and V.N.\ Remeslennikov.
\newblock Exponential groups {II}: extensions of centralisers and tensor
  completion of {CSA}--groups.
\newblock {\em Int. J. of Algebra and Comp.}, 6(6):687--711, 1996.

\bibitem{Pride91}
S.J.\ Pride.
\newblock Identities among relations of group presentations.
\newblock In (E.\ Ghys, A.\ Haefliger, and A.\ Verjovsky, editors, {\em Group
  Theory from a Geometric Viewpoint, Conference Proceedings, ICTP}, pages
  687--716. World Scientific, 1991.

\bibitem{RipsSela95}
E.\ Rips and Z.\ Sela.
\newblock Canonical representatives and equations in hyperbolic groups.
\newblock {\em Invent. Math.}, 120(3):489--512, 1995.

\bibitem{Rourke79}
C.P.\ Rourke.
\newblock Presentations and the trivial group.
\newblock In R.\ Fenn, editor, {\em Topology of Low Dimensional Manifolds},
  number 722 in Lecture Notes in Mathematics, pages 134--143. Springer, 1979.

\bibitem{Sela02}
Zlil Sela.
\newblock Diophantine geometry over groups. {II}. {C}ompletions, closures and
  formal solutions.
\newblock to appear.

\bibitem{Sela01}
Zlil Sela.
\newblock Diophantine geometry over groups. {I}. {M}akanin--{R}azborov
  diagrams.
\newblock {\em Publ. Math. Inst. Hautes \'{E}tudes Sci.}, 93:31--105, 2001.

\bibitem{Short84}
Hamish Short.
\newblock {\em Topological methods in group theory: the adjunction problem}.
\newblock PhD thesis, University of Warwick, 1984.

\end{thebibliography}
